\documentclass{kt}
\usepackage{srcltx,amssymb,amsmath, amsfonts, amsfonts,mathrsfs,comment}
\usepackage{hyperref}
\makeatletter
\def\blfootnote{\gdef\@thefnmark{}\@footnotetext}
\makeatother

\begin{document}
\newcommand{\lb}{\linebreak\makebox[3ex][r]{}}
\newcommand{\hs}{\hskip0.1ex}
\newcommand{\vs}{\vskip0.5ex }
\newcommand{\ds}{\displaystyle }
\newcommand{\bsl}{\baselineskip=2.7ex}
\newcommand{\vv}{\vskip-1ex }
\newcommand{\mb}{\makebox[2em][r]}
\newcommand{\bmp}{\begin{minipage}{13cm}}
\newcommand{\emp}{\end{minipage}\flushbottom\vskip3ex}
\newcommand{\rmp}{\end{minipage}\vskip2ex\flushbottom }
\newcommand{\goth}{\frak }
\newcommand{\f}{\varphi}
\newcommand{\hf}{\hfill}
\newcommand{\N}{\Bbb N}
\newcommand{\F}{\Bbb F}
\newcommand{\Z}{\Bbb Z}
\newcommand{\Q}{\Bbb Q}
\newcommand{\Sy}{\Bbb S}
\newcommand{\R}{\Rightarrow }
\newcommand{\W}{\Omega }
\newcommand{\w}{\omega }
\newcommand{\s}{\sigma }
\newcommand{\Lr}{\Leftrightarrow }

\newcommand{\e}{\varepsilon }
\newcommand{\g}{\gamma }
\renewcommand{\cal}{\mathscr}
\renewcommand{\geq}{\geqslant}
\renewcommand{\leq}{\leqslant}
\renewcommand{\ge}{\geqslant}
\renewcommand{\le}{\leqslant}
\newcommand{\vsh}{\vskip-1.5ex }
\newcommand{ \G} {\Gamma }
\newcommand{\sm}{\vskip0.5ex }

\newcommand{\ai}{$\boldsymbol{{}^{\displaystyle{\sharp}}}$}
\newcommand{\zv}{$\boldsymbol{{}^{\displaystyle{*}}}$}
\newcommand{\zva}{$\boldsymbol{{}^{\displaystyle{*}}}$}
\newcommand{\ul}{\vskip-1ex\underline{\makebox[4cm]{}}}
\newcommand{\otv}{\makebox[15pt][r]{$^{\displaystyle{*}}$}}
\newcommand{\oai}{\makebox[15pt][r]{$^{\displaystyle{\sharp}}$}}

 \thispagestyle{empty} \flushbottom \noindent
\begin{minipage}{13cm} {\large
\setlength{\unitlength}{1mm} R\hskip0.2ex u\hskip0.2ex
s\hskip0.2ex s\hskip0.2ex i\hskip0.2ex a\hskip0.2ex n\hskip0.2ex
~\hskip0.2ex A\hskip0.2ex c\hskip0.2ex a\hskip0.2ex d\hskip0.2ex
e\hskip0.2ex m\hskip0.2ex y\hskip0.2ex ~\hskip0.2ex o\hskip0.2ex
f\hskip0.2ex ~\hskip0.2ex S\hskip0.2ex c\hskip0.2ex i\hskip0.2ex
e\hskip0.2ex n\hskip0.2ex c\hskip0.2ex e\hskip0.2ex s\hskip0.2ex
 \hfill \begin{picture}(10,10)
\put(4,1){\circle*{0.8}}
\end{picture}
 \hfill S\hskip0.2ex i\hskip0.2ex b\hskip0.2ex
e\hskip0.2ex r\hskip0.2ex i\hskip0.2ex a\hskip0.2ex n\hskip0.2ex ~
\hskip0.2ex B\hskip0.2ex r\hskip0.2ex a\hskip0.2ex n\hskip0.2ex
c\hskip0.2ex h \vskip-1ex
\underline{\makebox[13cm]{}} \vskip1ex

S~O~B~O~L~E~V \hfill I~N~S~T~I~T~U~T~E \hfill O~F \hfill
M~A~T~H~E~M~A~T~I~C~S

}
\end{minipage}

\flushbottom
\noindent
\begin{minipage}{13cm}

\centerline {\LARGE \textbf{UNSOLVED PROBLEMS}} \vskip3ex

\centerline {\LARGE \textbf{IN GROUP THEORY}} \vskip3ex

\vskip3ex

\centerline {\huge \textbf{THE KOUROVKA NOTEBOOK}} \vskip3ex

\vskip3ex

\vskip7ex

\centerline {\huge \textbf{No. 21}}

\vskip10ex

\vspace{20ex}
\end{minipage}

\flushbottom
\noindent
\begin{minipage}{13cm}
\setlength{\unitlength}{1mm}
\centerline {\Large {Novosibirsk\begin{picture}(10,10)
\put(5,1){\circle*{0.8}}
\end{picture}2026}}

\end{minipage}

\flushbottom

\pagebreak
\thispagestyle{empty}
~
\vskip3ex

\noindent
Editors: E.\,I.\,Khukhro and V.\,D.\,Mazurov

\vskip3ex
\vskip33ex

\vskip7ex
{\large
\noindent
New problems and comments can be sent to the Editors:

\vskip2ex
\makebox[15ex][r]{}E.\,I.\,Khukhro \ or \ V.\,D.\,Mazurov \vskip1ex

 \makebox[15ex][r]{}University of Lincoln, U.K.

\vskip1ex

 \makebox[15ex][r]{}Sobolev Institute of Mathematics

\makebox[15ex][r]{}Novosibirsk-90, 630090, Russia

\vskip1.5ex

\makebox[15ex][r]{{\it e-mails\/}:}\vskip0.5ex

\makebox[15ex][r]{}khukhro@yahoo.co.uk\vskip1ex

\makebox[15ex][r]{}mazurov@math.nsc.ru}

\vskip13ex

\vskip13ex
\vskip17ex

\vskip3ex

\hfill \copyright \,\,E.\,I.\,Khukhro, V.\,D.\,Mazurov,
 2026

\newpage

\thispagestyle{empty}

~

~
\centerline {\Large \textbf{Contents}} \vskip4ex

\noindent Preface \dotfill \hskip1ex \pageref{intro} \vskip4ex

\noindent Problems from the 1st Issue (1965) \dotfill \hskip1ex
\pageref{1izd}\vskip2ex

\noindent Problems from the 2nd Issue (1966) \dotfill \hskip1ex
\pageref{2izd} \vskip2ex

\noindent Problems from the 3rd Issue (1969) \dotfill \hskip1ex
\pageref{3izd} \vskip2ex

\noindent Problems from the 4th Issue (1973) \dotfill \hskip1ex
\pageref{4izd} \vskip2ex

\noindent Problems from the 5th Issue (1976) \dotfill \hskip1ex
\pageref{5izd}
\vskip2ex

\noindent Problems from the 6th Issue (1978) \dotfill \hskip1ex
\pageref{6izd}
\vskip2ex

\noindent Problems from the 7th Issue (1980) \dotfill \hskip1ex
\pageref{7izd} \vskip2ex

\noindent Problems from the 8th Issue (1982) \dotfill \hskip1ex
\pageref{8izd} \vskip2ex

\noindent Problems from the 9th Issue (1984) \dotfill \hskip1ex
\pageref{9izd} \vskip2ex

\noindent Problems from the 10th Issue (1986) \dotfill \hskip1ex
\pageref{10izd}\vskip2ex

\noindent Problems from the 11th Issue (1990) \dotfill \hskip1ex
\pageref{11izd}
\vskip2ex

\noindent Problems from the 12th Issue (1992) \dotfill \hskip1ex
\pageref{12izd}
\vskip2ex

\noindent Problems from the 13th Issue (1995) \dotfill \hskip1ex
\pageref{13izd}
\vskip2ex

\noindent Problems from the 14th Issue (1999) \dotfill \hskip1ex
\pageref{14izd}
\vskip2ex

\noindent Problems from the 15th Issue (2002) \dotfill \hskip1ex
\pageref{15izd} \vskip2ex

\noindent Problems from the 16th Issue (2006) \dotfill \hskip1ex
\pageref{16izd} \vskip2ex

\noindent Problems from the 17th Issue (2010) \dotfill \hskip1ex
\pageref{17izd} \vskip2ex

 \noindent Problems from the 18th Issue (2014) \dotfill \hskip1ex
\pageref{18izd}

 \vskip2ex

 \noindent Problems from the 19th Issue (2018) \dotfill \hskip1ex
\pageref{19izd}

 \vskip2ex

 \noindent Problems from the 20th Issue (2022) \dotfill \hskip1ex
\pageref{20izd}
\vskip4ex

 \noindent New Problems (21st Issue, 2026)\dotfill \hskip1ex
\pageref{21izd}
\vskip4ex

\noindent Archive of Solved Problems \dotfill \hskip1ex\pageref{archiv} \vskip4ex

\noindent Index of Names \dotfill \hskip1ex
\pageref{index}

\newpage
\pagestyle{myheadings}
\markboth{\protect\vphantom{(y} Preface}{\protect\vphantom{(y} Preface}
\thispagestyle{headings}
~

\vspace{2ex}

\centerline {\Large \textbf{Preface}} \vskip4ex
\phantomsection\label{intro}
The idea of publishing a collection of unsolved problems in group theory
 was proposed by M.\,I.\,Kargapolov (1928--1976) at the Problem Day of
the First All--Union (All--USSR) Symposium on Group Theory which took
place in Kourovka, a small village near Sverdlovsk, on February,
16, 1965. This is why this collection acquired the name ``Kourovka
Notebook''. Since then every 2--4 years a new issue has appeared
containing new problems and incorporating the problems from the
previous issues with brief comments on the solved problems. \vskip0.5ex

For more than 60 years the ``Kourovka Notebook'' has served as a
unique means of communication for researchers in group theory and
nearby fields of mathematics. Maybe the most striking illustration
of its success is the fact that more than 3/4 of the problems from
the first and second issues have now been solved. Having acquired
international popularity the ``Kourovka Notebook'' includes
problems by more than 500 authors from all over the world.
\vskip0.5ex

 This is the 21st issue of the ``Kourovka Notebook''. It
contains 150 new problems. Comments have been added to those
problems from the previous issues that have been recently solved.
Some problems and comments from the previous issues had to be
altered or corrected. The Editors thank all those who sent their
remarks on the previous issues and helped in the preparation of the new issue.
We thank A.\,N.\,Ryaskin for assistance in dealing with the electronic publication.
\vskip0.5ex

The section ``Archive of Solved Problems'' contains {\it
all\/} the solved problems that had already been commented on, with a reference to a detailed publication containing a
complete answer, by the day of the first appearance of the  previous edition  in 2022. However, all the solutions that appeared in the updates after that day remain in the main part of the ``Kourovka Notebook'', among the unsolved problems of the corresponding section. (Some numbers of the problems are missing altogether, in those rare cases when problems were removed at the
request of the authors either as ill-conceived, or as no
 longer topical, for example, due to CFSG.) \vskip0.5ex

 The abbreviation CFSG stands for The Classification of the Finite
Simple Groups, which means that every finite simple
non-abelian group is isomorphic either to an alternating group, or to a
group of Lie type over a finite field, or to one of the twenty-six
sporadic groups (see D.\,Gorenstein,
{\it Finite Simple Groups: The Introduction to Their
Classification}, Plenum Press, New York, 1982). A note ``mod
CFSG'' in a comment means that the solution uses the CFSG.
\vskip0.5ex

Wherever possible, references to papers published in Russian are
given to their English translations. The index of names reflects the authors of the problems and papers cited in the problems and solutions, as well as all other names mentioned in the text.

Problems are marked with $\boldsymbol{\displaystyle{\sharp }}$ in those cases where some Artificial Intelligence (AI) system proposed a solution that has not yet been confirmed by mathematicians. Such a problem should be regarded as  unsolved. Anyone interested in a problem with an unverified AI solution is
invited to look into this AI solution and either confirm it or find mistakes --- in the first case a solution will be marked as a~proper solution with~$\boldsymbol{\displaystyle{*}}$, and in the second altogether deleted.

\vskip3.2ex
\parindent=0pt\bmp
\hfill {\it E.\,I.\,Khukhro, V.\,D.\,Mazurov}
\vskip1ex

\noindent Nottingham--Novosibirsk,  2026
\emp

\newpage
\pagestyle{myheadings} \markboth{1st Issue
(1965)}{1st Issue (1965)}
\thispagestyle{headings} ~ \vspace{2ex}

\centerline {\Large \textbf{Problems from the 1st Issue (1965)}}
\phantomsection\label{1izd}
\vspace{4ex}
\parindent=0pt

\bmp \textbf{1.3.} (Well-known problem). Can the group ring $\Z[G]$ of a
torsion-free group $G$ contain zero divisors? \hfill \raisebox{-1ex}{\sl
L.\,A.\,Bokut'}

\emp

\bmp \textbf{1.5.} (Well-known problem). Does there exist a group whose
group ring does not contain zero divisors and is not embeddable into
a skew field? \hfill \raisebox{-1ex}{\sl L.\,A.\,Bokut'}

\emp

\bmp \textbf{1.6.} (A.\,I.\,Mal'cev). Is the group ring of a
right-ordered group embeddable into a skew field? \hfill
\raisebox{-1ex}{\sl L.\,A.\,Bokut'}

\emp

\bmp \textbf{1.12.} (W.\,Magnus). The problem of the isomorphism to
the trivial group for all groups with $n$ generators and $n$
defining relations, where $n > 2$. \hfill \raisebox{-1ex}{\sl
M.\,D.\,Greendlinger}

\emp

\bmp \textbf{1.20.} For which groups (classes of groups) is the
lattice of normal subgroups first order definable in the lattice
of all subgroups? \hfill \raisebox{-1ex}{\sl Yu.\,L.\,Ershov}

\emp

\bmp \textbf{1.27.} Describe the universal theory of free groups.
\hfill {\sl M.\,I.\,Kargapolov}

\emp

\bmp \textbf{1.28.} Describe the universal theory of a free nilpotent
group. \hfill {\sl M.\,I.\,Kargapolov}

\emp

\bmp \textbf{1.31.} Is a residually finite group with the maximum
condition for subgroups almost polycyclic? \hfill
\raisebox{-1ex}{\sl M.\,I.\,Kargapolov}

\emp

\bmp \textbf{1.33.} (A.\,I.\,Mal'cev). Describe the automorphism
group of a free solvable group.

\hfill
\raisebox{-1ex}{\sl M.\,I.\,Kargapolov}

\emp

\bmp \textbf{1.35.} c) (A.\,I.\,Mal'cev,
 L.\,Fuchs). Do there exist simple pro-orderable
groups? A~group is said to be {\it pro-orderable\/} if each of its
partial orderings can be extended to a linear ordering. \hfill
\raisebox{-1ex}{\sl M.\,I.\,Kargapolov}

\emp

\bmp \textbf{1.40.} Is a group a nilgroup if it is the product of two
normal nilsubgroups?

\makebox[15pt][r]{}By definition, a nilgroup is a group consisting of
nilelements, in other words, of (not necessarily boundedly) Engel
elements. \hfill \raisebox{-1ex}{\sl Sh.\,S.\,Kemkhadze}

\emp

\bmp \textbf{1.46.} What conditions ensure the normalizer of a
relatively convex subgroup to be relatively convex? \hfill
\raisebox{-1ex}{\sl A.\,I.\,Kokorin}

\emp

\bmp \textbf{1.51.} What conditions ensure a matrix group over a
field (of complex numbers) to be orderable? \hfill
\raisebox{-1ex}{\sl A.\,I.\,Kokorin}

\emp

\bmp \textbf{1.54.} Describe all linear orderings of a free
metabelian group with a finite number of generators. \hfill
\raisebox{-1ex}{\sl A.\,I.\,Kokorin}

\emp

\bmp \textbf{1.55.} Give an elementary classification of
 linearly
ordered free groups with a fixed number of generators. \hfill
\raisebox{-1ex}{\sl A.\,I.\,Kokorin}

\emp

\bmp \textbf{1.65.} Is the class of groups of abelian extensions of
abelian groups closed under taking direct sums $(A,B) \mapsto A
\oplus B$? \hfill \raisebox{-1ex}{\sl L.\,Ya.\,Kulikov}

\emp

\bmp \textbf{1.67.} Suppose that $G$ is a finitely presented group,
$F$ a free group whose rank is equal to the minimal number of
generators of $G$, with a fixed homomorphism of $F$ onto $G$ with
kernel $N$. Find a complete system of invariants of the
factor-group of $N$ by the commutator subgroup $[F,N]$. \hfill
\raisebox{-1ex}{\sl L.\,Ya.\,Kulikov}

\emp

\bmp \textbf{1.74.} Describe all minimal topological groups, that is,
non-discrete groups all of whose closed subgroups are discrete. The
minimal locally compact groups can be described without much effort.
At the same time, the problem is probably complicated in the general
case. \hfill \raisebox{-1ex}{\sl V.\,P.\,Platonov}

\emp

\bmp \textbf{1.86.} Is it true that the identical relations of a
polycyclic group have a finite basis?

\hfill
\raisebox{-1ex}{\sl A.\,L.\,Shmel'kin}

\emp

\bmp \textbf{1.87.} The same question for matrix groups (at least
over a field of characteristic~0).

\hfill
\raisebox{-1ex}{\sl A.\,L.\,Shmel'kin}

\emp

\raggedbottom

\newpage
\pagestyle{myheadings} \markboth{2nd Issue
(1967)}{2nd Issue (1967)}
\thispagestyle{headings} ~ \vspace{2ex}

\centerline {\Large \textbf{Problems from the 2nd Issue (1967)}}
\phantomsection\label{2izd}
\vspace{4ex}

\bmp \textbf{2.5.} According to Plotkin, a group is called an $NR$-group if the set of its nil-elements coincides with the locally nilpotent radical, or, which is equivalent, if every inner automorphism of it is locally stable. Can an $NR$-group have a nil-automorphism that is not locally stable?
 \hfill \raisebox{-1ex}{\sl V.\,G.\,Vilyatser}

\emp

\bmp \textbf{2.6.} In an $NR$-group, the set of generalized central elements coincides with the nil-kernel. Is the converse true, that is, must a group be an $NR$-group if the set of generalized central elements coincides with the nil-kernel?
 \hfill \raisebox{-1ex}{\sl V.\,G.\,Vilyatser}

\emp

\bmp \textbf{2.9.} Do there exist regular associative operations on the
class of groups satisfying the weakened Mal'cev condition (that is,
monomorphisms of the factors of an arbitrary product can be glued together, generally
speaking, into a homomorphism of the whole
product), but not satisfying the analogous condition for
epimorphisms of the factors?

 \hfill \raisebox{-1ex}{\sl O.\,N.\,Golovin}

\emp

\bmp \textbf{2.22.} a) An abstract group-theoretic property $\Sigma$ is
said to be {\it radical $($in our sense\/}) if, in any group $G$,
the subgroup $\Sigma(G)$ generated by all normal $\Sigma
$-sub\-groups is a $\Sigma $-sub\-group itself (called the {\it
$\Sigma $-radical\/} of $G$). A radical property $\Sigma$ is said
to be {\it strongly radical\/} if, for any group $G$, the
factor-group $G/\Sigma(G)$ contains no non-trivial normal $\Sigma
$-sub\-groups. Is the property $\,\overline{\! RN}$ radical?
strongly radical? \hfill \raisebox{-1ex}{\sl Sh.\,S.\,Kemkhadze}

\emp

\bmp \textbf{2.24.} Are Engel torsion-free groups orderable? \hfill
\raisebox{0ex}{\sl A.\,I.\,Kokorin} \emp

\bmp \textbf{2.25.} a) (L.\,Fuchs). Describe the groups which are
linearly orderable in only finitely many ways. \hfill
\raisebox{-1ex}{\sl A.\,I.\,Kokorin}

\emp

\bmp \textbf{2.26.} (L.\,Fuchs). Characterize as abstract groups the
multiplicative groups of orderable skew fields. \hfill
\raisebox{-1ex}{\sl A.\,I.\,Kokorin} \emp

\bmp \textbf{2.28.} Can every orderable group be embedded in a
pro-orderable group? (See 1.35 for the definition of pro-orderable
groups.) \hfill \raisebox{-1ex}{\sl A.\,I.\,Kokorin}

\emp

\bmp \textbf{2.40.}
c) The {\it $I $-theory\/} ({\it $Q\hskip0.1ex
$-theory\/}) of a class ${\goth K}$ of universal algebras is the
totality of all identities (quasi-identities) that are valid on
all algebras in~${\goth K}$. Does there exist a finitely
axiomatizable variety

 \makebox[25pt][r]{(3)} of Lie rings

\makebox[40pt][r]{(i)} whose $I $-theory is non-decidable?

 \makebox[15pt][r]{}{\it Remark:\/} it is not difficult to find a
 recursively axiomatizable variety of semigroups with identity whose $ I
$-theory is non-recursive (see also A.\,I.\,Mal'cev, {\it Mat.
Sbornik}, \textbf{69}, no.\,1 (1966), 3--12 (Russian)).\hfill
\raisebox{-1ex}{\sl A.\,I.\,Mal'cev}
\emp

\bmp \textbf{2.42.} What is the structure of the groupoid of
quasivarieties

\makebox[15pt][r]{}a) of all semigroups?

\makebox[15pt][r]{}b) of all rings?

\makebox[15pt][r]{}c) of all associative rings?

\makebox[15pt][r]{}Compare with A.\,I.\,Mal'cev, {\it Siberian
Math.~J.}, \textbf{8}, no.\,2 (1967), 254--267).

\hfill
\raisebox{-1ex}{\sl A.\,I.\,Mal'cev}

\emp

\bmp \textbf{2.45.} b)
 (P.\,Hall). Prove or refute the following conjecture:
 If the marginal subgroup $v^*G$ has finite index $m$ in $G$,
then the order of $vG$ is finite and divides a power of $m$.

\makebox[15pt][r]{}{\it Editors' comment:\/} There are examples when $|vG|$ does not divide a power of $m$
(Yu.\,G.\,Kleiman, {\it Trans. Moscow Math. Soc.}, {\bf 1983},
no.\,2, 63--110), but the question of finiteness remains open.
\hfill
\raisebox{-1ex} {\sl Yu.\,I.\,Merzlyakov}
\emp

\bmp \textbf{\zv 2.48.}
(N.\,Aronszajn). Let $G$ be a connected topological
group locally satisfying some identical relation $f|_U= 1$, where $U$
is a neighborhood of the identity element of~$G$. Is it then true
that $f|_G=1$? \hfill \raisebox{-1ex}{\sl V.\,P.\,Platonov}

\ul

\otv
No, not necessarily: for any odd integer $n\geq 665$
 there is a connected topological group such that the identity $x^n=1$ holds in some neighborhood of unity, but not in the entire group (E.\,Reznichenko, I.\,Zyabrev, \textit{Preprint}, 2024, \url{https://arxiv.org/abs/2406.05203}).
\emp

\bmp \textbf{2.56.} Classify up to isomorphism the abelian connected
algebraic unipotent linear groups over a field of positive
characteristic. This is not difficult in the case of a field of
characteristic zero. On the other hand, C.\,Chevalley has solved
the classification problem for such groups up to isogeny. \hfill
\raisebox{-1ex}{\sl V.\,P.\,Platonov}

\emp

\bmp \textbf{2.67.} Find conditions which ensure that the nilpotent
product of pure nilpotent groups (from certain classes) is
determined by the lattice of its subgroups. This is known to be
true if the product is torsion-free.\hfill \raisebox{-1ex}{\sl
L.\,E.\,Sadovski\u{\i}} \emp

\bmp \textbf{2.68.} What can one say about lattice isomorphisms of a
pure soluble group? Is such a group strictly determined by its
lattice? It is well-known that the answer is affirmative for free
soluble groups.\hfill \raisebox{-1ex}{\sl L.\,E.\,Sadovski\u{\i}} \emp

\bmp \textbf{2.74.} (Well-known problem). Describe the finite groups
all of whose involutions have soluble centralizers. \hfill
\raisebox{-1ex}{\sl A.\,I.\,Starostin} \emp

\bmp
\textbf{2.78.}
 Any set of all subgroups of the same given order of a finite group
$G$ that contains at least one non-normal subgroup
is called an $IE_{\bar n}$-system of $G$. A positive integer $k$
is called a
{\it soluble\/} ({\it
non-soluble; \
simple; \
composite; \
absolutely simple\/})
{\it group-theoretic number\/}
if every finite group having exactly $ k$ \ $IE_{\bar
n}$-systems is soluble (respectively,
if there is at least one
non-soluble finite group having $k$ \ $IE_{\bar n}$-systems; \
if there is at least one
simple finite group having $k$ \ $IE_{\bar n}$-systems; \
if there are no simple finite groups having $k$ \ $IE_{\bar n}$-systems; \
if there is at least one
simple finite group having $k$ \ $IE_{\bar n}$-systems and there
are no non-soluble non-simple finite groups having $k$ \
$IE_{\bar n}$-systems).

\makebox[15pt][r]{}Are the sets of all soluble and of all
absolutely simple group-theoretic numbers finite or infinite? Do
there exist composite, but not soluble group-theoretic numbers?

 \hfill \raisebox{-1ex}{\sl P.\,I.\,Trofimov}

\emp

\bmp \textbf{2.80.} Does every non-trivial group satisfying the
normalizer condition contain a non-trivial abelian normal
subgroup? \hfill \raisebox{-1ex}{\sl S.\,N.\,Chernikov}

\emp

\bmp \textbf{\zv 2.81.}
a) Does there exist an axiomatizable class of
lattices $ {\goth K}$ such that the lattice of all subsemigroups
of a semigroup $S$ is isomorphic to some lattice in ${\goth K}$ if
and only if $S$ is a free group?

\makebox[15pt][r]{}b) The same question for free
abelian groups.

 \makebox[15pt][r]{}Analogous questions have affirmative answers for
torsion-free groups, for non-periodic groups, for abelian
torsion-free groups, for abelian non-periodic groups, for
orderable groups (the corresponding classes of lattices are
even finitely axiomatizable). Thus, in posed questions one may
assume from the outset that the semigroup $S$ is a
torsion-free group (respectively, a torsion-free abelian
group). \hfill \raisebox{-1ex}{\sl L.\,N.\,Shevrin}

\ul

\otv a), b) Yes, there are such classes (AI solution confirmed by O.\,G.\,Kharlampovich: \  Achyuth Jayadevan,  {\it Preprint}, 2026, \url{https://zenodo.org/records/22721223}).
\emp 

\bmp \textbf{2.82.}
Can the class of groups with the $n $th Engel condition
$[x,\underbrace{y,\ldots ,y}_{n}]=1$ be defined by identical
relations of the form $u=v$, where $u$ and $v$ are words without
negative powers of variables? This can be done for $n=1,\,2,\,3$
(A.\,I.\,Shirshov, {\it Algebra i Logika},~\textbf{2}, no.\,5 (1963),
5--18 (Russian)).

\makebox[15pt][r]{}{\it Editors' comment:\/} This has also been
done for $n=4$ (G.\,Traustason, {\it J.\,Group Theory},~\textbf{2},
no.\,1 (1999), 39--46). As remarked by O.\,Macedo\'nska, this is
also true for the class of locally graded $n$-Engel groups because
by (Y.\,Kim, A.\,Rhemtulla, in {\it Groups--Korea'94}, de~Gruyter,
Berlin, 1995, 189--197) such a group is locally nilpotent and then
by (R.\,G.\,Burns, Yu.\,Medvedev, {\it J.~Austral. Math. Soc.},
\textbf{64} (1998), 92--100) such a group is an extension of a
nilpotent group of $n$-bounded class by a group of $n$-bounded
exponent; then a classical result of Mal'cev implies that such
groups satisfy a positive law. \hfill \raisebox{-1ex}{\sl
A.\,I.\,Shirshov}

\emp

\bmp \textbf{2.84.} Suppose that a locally finite group $G$ is a
product of two locally nilpotent subgroups. Is $G$ necessarily
locally soluble? \hfill \raisebox{-1ex}{\sl V.\,P.\,Shunkov} \emp

\raggedbottom

\newpage
\pagestyle{myheadings} \markboth{3rd Issue
(1969)}{3rd Issue (1969)}
\thispagestyle{headings} ~ \vspace{2ex}

\centerline {\Large \textbf{Problems from the 3rd Issue (1969)}}
\phantomsection\label{3izd}
\vspace{4ex}

 \bmp \textbf{3.3.} (Well-known problem). Describe the automorphism
group of the free associative algebra with $n$ generators, $n \geq
2$.

 \makebox[15pt][r]{}It is known
 that all automorphisms are tame for $n=2$ (A.\,J.\,Czerniakiewicz, \emph{Trans. Amer. Math. Soc.}, \textbf{160} (1971), 393--401; \textbf{171} (1972), 309--315; \
 L.\,G.\,Makar-Limanov, \emph{Funct. Anal. Appl.}, \textbf{4} (1971), 262--264). For $n=3$, the Anick automorphism is wild (U.\,U.\,Umirbaev, \emph{J.~Reine Angew. Math.}, \textbf{605} (2007), 165--178). There is an algorithm for recognizing automorphisms among endomorphisms of free associative algebras of finite rank over constructive fields (A.\,V.\,Jagzhev, \emph{Siberian. Math.~J.},
\textbf{21} (1980), 142--146). \hfill \raisebox{-1ex}{\sl L.\,A.\,Bokut'}
\emp

\bmp
\textbf{3.5.} Can a subgroup of a relatively free group be radicable? In particular, can a verbal subgroup of a relatively free group be radicable?\hfill \raisebox{-1ex}{\sl N.\,R.\,Brumberg}

\emp

\bmp
\textbf{\zv 3.12.}
(Well-known problem). Is a locally finite group
with a full Sylow basis locally soluble? \hfill
\raisebox{-1ex}{\sl Yu.\,M.\,Gorchakov}

\ul

\otv No, not always (AI solution confirmed by A.\,A.\,Buturlakin: \ Achyuth Jayadevan,  {\it Preprint}, 2026, \url{https://zenodo.org/records/22721349}).
\emp

\bmp \textbf{3.16.} The word problem for a group admitting a single
defining relation in the variety of soluble groups of derived
length $n$, $n \geq 3$. \hfill \raisebox{-1ex}{\sl
M.\,I.\,Kargapolov} \emp

\bmp \textbf{\ai 3.20.}  
Does the class of orderable groups coincide with
the smallest axiomatizable class containing pro-orderable groups?
(See 1.35.) \hfill \raisebox{-1ex}{\sl A.\,I.\,Kokorin}

\ul

\oai No, it does not (unverified AI solution:  Achyuth Jayadevan,  {\it Preprint}, 2026, \url{https://zenodo.org/records/22721372}).
 \emp 

\bmp \textbf{3.34.} (Well-known problem). The conjugacy problem for
groups with a single defining relation. \hfill \raisebox{-1ex}{\sl
D.\,I.\,Moldavanski\u{\i}} \emp

\bmp \textbf{3.38.} Describe the topological groups which have no
proper closed subgroups.

\hfill \raisebox{-1ex}{\sl Yu.\,N.\,Mukhin}

\emp

\bmp \textbf{\ai 3.43.} 
Let $\mu$ be an infinite cardinal number. A group
$G$ is said to be $\mu $-{\it over\-nil\-po\-tent\/} if every
cyclic subgroup of $G$ is a member of some ascending normal series
of length less than $\mu$ reaching $G$. It is not difficult to
show that the class of $\mu $-over\-nil\-po\-tent groups is a
radical class. Is it true that if $\mu _1 < \mu _2$ for two
infinite cardinal numbers $\mu _1$ and $\mu _2$, then there exists
a group $G$ which is $\mu _2 $-over\-nil\-po\-tent and $\mu _1
$-semi\-simple?

\makebox[15pt][r]{}{\it Remark of 2001:\/} In (S.\,Vovsi, {\it
Sov. Math. Dokl.}, \textbf{13} (1972), 408--410) it was proved that
for any two infinite cardinals $\mu _1 < \mu _2$ there exists a
group that is $\mu _2$-overnilpotent but not $\mu _
1$-overnilpotent.
 \hfill \raisebox{-1ex}{\sl B.\,I.\,Plotkin}

\ul

\oai Yes, it is true (unverified AI solution:  Achyuth Jayadevan,  {\it Preprint}, 2026,

\url{https://zenodo.org/records/22721386}).
\emp

\bmp \textbf{3.44.} Suppose that a group is generated by its
ascendant soluble subgroups. Is it necessarily locally soluble?
\hfill \raisebox{-1ex}{\sl B.\,I.\,Plotkin}

\emp

\bmp \textbf{3.45.} Let ${\goth X}$ be a hereditary radical. Is it
true that, in a locally nilpotent torsion-free group $G$, the
subgroup ${\goth X}(G)$ is isolated? \hfill \raisebox{-1ex}{\sl
B.\,I.\,Plotkin} \emp

\bmp \textbf{\zv 3.46.}
Does there exist a group having more than one,
but finitely many maximal locally soluble normal subgroups? \hfill
\raisebox{-1ex}{\sl B.\,I.\,Plotkin}

\ul

\otv Yes, such groups exist (AI solution confirmed by M.\,Trombetti: \ P.\,Monticone, {\it Preprint}, 2026,
\url{https://kourovkanotebookorg.wordpress.com/wp-content/uploads/2026/06/3_46.pdf}); see also (W.\,van\,Doorn, E.\,Judin, P.\,Monticone, D.\,Morrison,  {\it Preprint}, 2026, \url{https://arxiv.org/abs/2607.17477}).
\emp

\bmp \textbf{3.47.} (Well-known problem of A.\,I.\,Maltsev).
Is it true that every locally
nilpotent group is a ho\-mo\-mor\-phic image of some
torsion-free locally nilpotent group?

\makebox[15pt][r]{}{\it Editors' comment:} An affirmative answer is
known for periodic groups (E.\,M.\,Le\-vich, A.\,I.\,To\-ka\-renko,
{\it Siberian Math.~J.}, \textbf{11} (1970), 1033--1034) and for
countable groups (N.\,S.\,Romanov\-ski\u{\i}, Preprint, 1969). \hfill
\raisebox{-1ex}{\sl B.\,I.\,Plotkin}

\emp

\bmp \textbf{3.48.} It can be shown that hereditary radicals form a
semigroup with respect to taking the products of classes. It is an
interesting problem to find all indecomposable elements of this
semigroup. In particular, we point out the problem of finding all
indecomposable radicals contained in the class of locally finite
$p\hskip0.2ex $-groups.

\hfill
\raisebox{-1ex}{\sl B.\,I.\,Plotkin}

\emp

\bmp \textbf{3.49.} Does the semigroup generated by all
indecomposable radicals satisfy any identity? \hfill
\raisebox{-1ex}{\sl B.\,I.\,Plotkin}

\emp

\bmp \textbf{3.55.} Is every binary soluble group all of whose abelian
subgroups have finite rank, locally soluble? \hfill
\raisebox{-1ex}{\sl S.\,P.\,Strunkov}

\emp

\bmp \textbf{3.57.} Determine the laws of distribution of non-soluble
and simple group-theoretic numbers in the sequence of natural
numbers. See 2.78. \hfill\raisebox{-1ex}{\sl P.\,I.\,Trofimov}

\emp

\bmp \textbf{3.60.} The notion of the {\it $p\hskip0.2ex $-length\/}
of an arbitrary finite group was introduced in (L.\,A.\,Shemetkov,
{\it Math. USSR Sbornik}, \textbf{1} (1968), 83--92). Investigate the
relations between the $p\hskip0.2ex $-length of a finite group and
the invariants $c_p$, $d_p$, $e_p$ of its Sylow $p\hskip0.2ex
$-sub\-group. \hfill \raisebox{-1ex}{\sl L.\,A.\,Shemetkov}

\emp

\raggedbottom

\newpage
\pagestyle{myheadings} \markboth{4th Issue
(1973)}{4th Issue (1973)}
\thispagestyle{headings} ~ \vspace{2ex}

\centerline {\Large \textbf{Problems from the 4th Issue (1973)}}
\phantomsection\label{4izd}
\vspace{4ex}

\bmp \textbf{4.2.} a) Find an infinite finitely generated group of
exponent $< 100$.

\makebox[15pt][r]{}b) Do there exist such groups of exponent 5?
\hfill {\sl S.\,I.\,Adian}

\emp

\bmp
\textbf{4.5.} b) Is it true that an arbitrary finitely
presented group
has either polynomial or exponential growth? \hfill
\raisebox{-1ex}{\sl S.\,I.\,Adian}

\emp

\bmp \textbf{4.6.} (P.\,Hall). Are the projective groups in the
variety of metabelian groups free?

\hfill
\raisebox{-1ex}{\sl V.\,A.\,Artamonov}

\emp

\bmp \textbf{4.7.} (Well-known problem). For which ring epimorphisms
$R \rightarrow Q$ is the corresponding group homomorphism $SL_n(R)
\rightarrow SL_n(Q)$ an epimorphism (for a fixed $n \geq 2$)? In
particular, for what rings $R$ does the equality $SL_n(R)=E_n(R)$
hold?

\hfill \raisebox{-1ex}{\sl V.\,A.\,Artamonov}

\emp

\bmp \textbf{4.9.} Let $G$ be a finitely generated torsion-free
nilpotent group. Are there only a finite number of non-isomorphic
groups in the sequence $\alpha G$, $\alpha ^2G, \ldots $\,? Here
$\alpha G$ denotes the automorphism group of $G$ and $\alpha
^{n+1}G=\alpha (\alpha ^nG)$ for $n=1,\,2,\ldots $ \hfill
\raisebox{-1ex}{\sl G.\,Baumslag}

\emp

\bmp \textbf{4.11.} Let $F$ be the free group of rank 2 in some
variety of groups. If $F$ is not finitely presented, is the
multiplicator of $F$ necessarily infinitely generated? \hfill
\raisebox{-1ex}{\sl G.\,Baumslag} \emp

\bmp \textbf{\zv 4.13.} 
Prove that every finite non-abelian $p\hskip0.2ex
$-group admits an automorphism of order $p$ which is not an inner
one. \hfill \raisebox{-1ex}{\sl Ya.\,G.\,Berkovich}

\ul

\otv This is proved (M.\,Y.\,Kizmaz, {\it Preprint of 22 November, 2026}, \url{https://arxiv.org/abs/2611.29482}; \ L.\,Gong,  J.\,Chen, J.\,Wua, E.\,Cosme-Ll\'opez, {\it Preprint of 24 November, 2026}, \url{https://zenodo.org/records/24498173}).
\emp

\bmp \textbf{\ai 4.17.} 
If $ {\goth V}$ is a variety of groups, denote by
$\widetilde{{\goth V}}$ the class of all finite ${\goth V}
$-groups. How to characterize the classes of finite groups of the
form $\widetilde{{\goth V}}$ for ${\goth V}$ a variety? \hfill
\raisebox{-1ex}{\sl R.\,Baer}

\ul

\oai A characterization is given (unverified AI solution:  Achyuth Jayadevan,  {\it Preprint}, 2026, \url{https://zenodo.org/records/22724488}).
\emp

\bmp \textbf{4.18.} Characterize the classes ${\goth K}$ of groups,
meeting the following requirements: subgroups, epimorphic images
and groups of automorphisms of ${\goth K} $-groups are $ {\goth K}
$-groups, but not every countable group is a ${\goth K} $-group.
Note that the class of all finite groups and the class of all
almost cyclic groups meet these requirements. \hfill
\raisebox{-1ex}{\sl R.\,Baer}

\emp

\bmp \textbf{4.19.} Denote by ${\goth C}$ the class of all groups $G$
with the following property: if $U$ and $V$ are maximal locally
soluble subgroups of $G$, then $U$ and $V$ are conjugate in $G$
(or at least isomorphic). It is almost obvious that a finite group
$G$ belongs to ${\goth C}$ if and only if $G$ is soluble. What can
be said about the locally finite groups in ${\goth C}$? \hfill
\raisebox{-1ex}{\sl R.\,Baer}

\emp

\bmp \textbf{4.24.} Suppose $T$ is a non-abelian Sylow 2-sub\-group
of a finite simple group $G$.

 \makebox[15pt][r]{}a) Suppose $T$ has nilpotency class $n$. The
best possible bound for the exponent of the center of $T$ is
$2^{n-1}$. This easily implies the bound $2^{n(n-1)}$ for the
exponent of~$T$, however, this is almost certainly too crude. What
is the best possible bound?

 \makebox[15pt][r]{}d) Find a ``small number'' of subgroups
$T_1,\ldots ,T_n$ of $T$ which depend only on the isomorphism
class of $T$ such that $\{ N_G(T_1), \ldots , N_G(T_n)\}$
together control fusion in $T$ with respect to $G$ (in the sense
of Alperin).\hfill \raisebox{-1ex}{\sl D.\,M.\,Goldschmidt} \emp

\bmp \textbf{4.30.} Describe the groups (finite groups, abelian
groups) that are the full automorphism groups of topological
groups. \hfill \raisebox{-1ex}{\sl M.\,I.\,Kargapolov} \emp

\bmp \textbf{4.31.} Describe the lattice of quasivarieties of
nilpotent groups of nilpotency class~2.

\hfill
\raisebox{-1ex}{\sl M.\,I.\,Kargapolov}
\emp

\bmp \textbf{4.33.} Let ${\goth K}_n$ be the class of all groups with
a single defining relation in the variety of soluble groups of
derived length $n$.

 \makebox[25pt][r]{a)} Under what conditions does a ${\goth
K}_n $-group have non-trivial center? Can a $ {\goth K}_n $-group,
$n \geq 2$, that cannot be generated by two elements have
non-trivial centre? {\it Editors' comment (1998):\/} These
questions were answered for $n=2$ (E.\,I.\,Timoshenko, {\it
Siberian Math.~J.}, \textbf{14}, no.\,6 (1973), 954--957; {\it Math.
Notes}, \textbf{64}, no.\,6 (1998), 798--803).

 \makebox[25pt][r]{b)} Describe the abelian subgroups of ${\goth
K}_n $-groups.

\makebox[25pt][r]{c)} Investigate the periodic subgroups of
${\goth K}_n $-groups. \hfill {\sl M.\,I.\,Kargapolov} \emp

\bmp \textbf{4.34.} Let $v$ be a group word, and let ${\goth K}_v$ be
the class of groups $G $ such that there exists a positive integer
$n=n(G)$ such that each element of the verbal subgroup $vG$ can be
represented as a product of $n$ values of the word $v$ on the
group $G$.

\makebox[25pt][r]{a)} For which $v$ do all finitely generated
soluble groups belong to the class ${\goth K}_v$?

 \makebox[25pt][r]{b)} Does the word $v(x,y)=x^{-1}y^{-1}xy$ satisfy
this
condition?
\hfill \raisebox{-1ex}{\sl M.\,I.\,Kargapolov}
\emp

\bmp \textbf{\ai 4.40.} 
Let $C$ be a fixed non-trivial group (for
instance, $C={\Bbb Z}/2{\Bbb Z}$). As it is shown in
(Yu.\,I.\,Merzlyakov, {\it Algebra and Logic}, \textbf{9}, no.\,5
(1970), 326--337), for any two groups $A$ and $B$, all split
extensions of $B$ by $A$ can be imbedded in a certain unified way
into the direct product $A \times {\rm Aut}\, (B\wr C)$. How are
they situated in it? \hfill \raisebox{-1ex}{\sl
Yu.\,I.\,Merzlyakov}

\ul

\oai A description is given (unverified AI solution:  Achyuth Jayadevan,  {\it Preprint}, 2026, \url{https://zenodo.org/records/22727276}).
\emp 

 \bmp
\textbf{4.42.} For what natural \nopagebreak numbers $n$ does the
following equality hold: \vspace{-1ex}\nopagebreak $$GL_n({\Bbb
R})=D_n({\Bbb R})\cdot O_n({\Bbb R})\cdot UT_n({\Bbb R})\cdot
GL_n({\Bbb Z})? \vspace{-1.5ex}$$\nopagebreak For notation, see,
for instance, (M.\,I.\,Kargapolov, Ju.\,I.\,Merzljakov, {\it
Fundamentals of the Theory of Groups}, Springer, New York, 1979).
For given $n$ this equality implies the affirmative solution of
Minkowski's problem on the product of $n$ linear forms
(A.\,M.\,Macbeath, {\it Proc. Glasgow Math. Assoc.}, \textbf{5},
no.\,2 (1961), 86--89), which remains open for $n \geq 6$. It is
known that the equality holds for $n \leq 3$
(Kh.\,N.\,Narzullayev, {\it Math. Notes}, \textbf{18} (1975),
713--719); on the other hand, it does not hold for all
\nopagebreak sufficiently large $n$ (N.\,S.\,Akhmedov,
\nopagebreak {\it Zapiski Nauchn. Seminarov LOMI}, \nopagebreak
\textbf{67} (1977), 86--107 (Russian)). \nopagebreak \hfill
\raisebox{-1ex}{\sl Yu.\,I.\,Merzlyakov} \emp

\bmp \textbf{4.44.} (Well-known problem). Describe the groups whose
automorphism groups are abelian. \hfill \raisebox{-1ex}{\sl
V.\,T.\,Nagrebetski\u{\i}} \emp

\bmp \textbf{4.46.} a) We call a variety of groups a {\it limit
variety\/} if it cannot be defined by finitely many laws, while
each of its proper subvarieties has a finite basis of identities.
It follows from Zorn's lemma that every variety that has no finite
basis of identities contains a limit subvariety. Give explicitly
(by means of identities or by a generating group) at least one
limit variety.
 \hfill \raisebox{-1ex}{\sl A.\,Yu.\,Olshanskii}
\emp

\bmp \textbf{4.48.} A locally finite group is said to be an {\it $A
$-group\/} if all of its Sylow subgroups are abelian. Does every
variety of $A $-groups possess a finite basis of identities?

\hfill \raisebox{-1ex}{\sl A.\,Yu.\,Olshanskii}
\emp

\bmp \textbf{4.50.} What are the soluble varieties of groups all of
whose finitely generated subgroups are residually finite? \hfill
\raisebox{-1ex}{\sl V.\,N.\,Remeslennikov}
\emp

\bmp \textbf{\zv 4.55.}
Let $G$ be a finite group, and ${\Bbb Z}_{(p)}$ the
localization at $p$. Can every projective $\Z_{(p)}G$-module be uniquely, up to permutations and isomorphisms, written as a direct sum of
indecomposable projective $\Z_{(p)}G$-modules?

 \makebox[15pt][r]{}It is known that the category of finitely generated projective $Z_{(p)}G$-modules is not a Krull--Schmidt category (S.\,M.\,Woods, {\it Canadian J. Math.}, {\bf 26}, no.\,1 (1974), 121--129). See also (D.\,Johnston, D.\,Rumynin, {\it J.~Algebra, to appear}, \url{https://arxiv.org/pdf/2507.21316}).\hfill \raisebox{-1ex}{\sl
K.\,W.\,Roggenkamp}

\ul

\otv No, not every (AI solution confirmed by D.\,Rumynin: \  Achyuth Jayadevan,  {\it Preprint of 14 September 2026},
\url{https://zenodo.org/records/22750727}; \ unverified AI solution: M.\,Kinyon, J.\,Urban, {\it Preprint of 17 September 2026},
 \url{https://kourovkanotebookorg.wordpress.com/wp-content/uploads/2026/09/kourovka-mathematics.pdf}).
\emp

\bmp \textbf{4.56.} Let $R$ be a commutative Noetherian ring with 1, and
$\Lambda$ an $ R $-algebra, which is finitely generated as $R
$-module. Put $T=\{ U \in {\rm Mod}\, \Lambda \mid \exists \mbox{ an
exact}\; \Lambda \mbox{-sequence} \; 0 \rightarrow P \rightarrow
\Lambda ^{(n)} \rightarrow U \rightarrow 0\mbox{ for some}\;n,\;
\mbox{with}\;P_m \cong \Lambda _m^{(n)} \;\mbox{for every maximal
ideal}\;m \mbox{ of}\;R\} $. Denote by ${\Bbb G}(T)$ the Grothendieck
group of $T$ relative to short exact sequences.

\makebox[25pt][r]{a)} Describe ${\Bbb G}(T)$, in particular, what
does it mean: $[U]=[V]$ in ${\Bbb G}(T)$?

\makebox[25pt][r]{b)} {\it Conjecture:\/} if $\dim (\max (R))=d <
\infty $, and there are two epimorphisms $\varphi \! : \Lambda ^{(n)}
\rightarrow U$, \ $\psi \! : \Lambda ^{(n)} \rightarrow V$, \ $n
> d$, \ and \ $[U]=[V]$ \ in \ ${\Bbb G}(T)$, then ${\rm Ker}\, \varphi
={\rm Ker}\, \psi $.

\hfill \raisebox{-1ex}{\sl K.\,W.\,Roggenkamp}
\emp

\vspace{-2ex} \bmp \textbf{4.65.} {\it Conjecture:\/}
$\displaystyle{\frac{p^q-1}{p-1}}$ never divides
$\displaystyle{\frac{q^p-1}{q-1}}$ if $p$, $q$ are distinct
primes. The validity of this conjecture would simplify the proof
of solvability of groups of odd order (W.\,Feit, J.\,G.\,Thompson,
{\it Pacific J.~Math.}, \textbf{13}, no.\,3 (1963), 775--1029),
rendering unnecessary the detailed use of generators and
relations. \hfill \raisebox{-1ex}{\sl J.\,G.\,Thompson}

\emp

\bmp \textbf{4.66.} Let $P$ be a presentation of a finite group $G$
on $m_p$ generators and $r_p$ relations. The {\it deficiency\/}
${\rm def}\, (G)$ is the maximum of $m_p - r_p$ over all
presentations $P$. Let $G$ be a finite group such that $G=G'\ne 1$
and the multiplicator $M(G)=1$. Prove that ${\rm def}\, (G^n)
\rightarrow - \infty $ as $n \rightarrow \infty $, where $G^n$ is
the $n$th direct power of~$G$. \hfill \raisebox{-1ex}{\sl
J.\,Wiegold}

\emp

\bmp \textbf{4.72.} Is it true that every variety of groups whose
free groups are residually nilpotent torsion-free, is either
soluble or coincides with the variety of all groups? For an
affirmative answer, it is sufficient to show that every variety of
Lie algebras over the field of rational numbers which does not
contain any finite-dimensional simple algebras is soluble. \hfill
\raisebox{-1ex}{\sl A.\,L.\,Shmel'kin}

\emp

\bmp \textbf{4.74.} b) Is every binary-finite 2-group of order
greater than 2 non-simple?

\hfill \raisebox{-1ex}{\sl V.\,P.\,Shunkov}

\emp

\bmp \textbf{\zv 4.75.}
Let $G$ be a periodic group containing an
involution $i$ and suppose that the Sylow 2-sub\-groups of $G$ are
either locally cyclic or generalized quaternion. Does the element
$ iO_{2'}(G)$ of the factor-group $G/O_{2'}(G)$ always lie in its
centre?

\hfill \raisebox{-1ex}{\sl V.\,P.\,Shunkov}

\ul

\otv No, not always (AI solution confirmed by M.\,Amelio: \  M.\,Kinyon, J.\,Urban, {\it Preprint}, 2026,
\url{https://kourovkanotebookorg.wordpress.com/wp-content/uploads/2026/09/kourovka-mathematics.pdf}).
\emp 

\raggedbottom

\newpage

\pagestyle{myheadings} \markboth{5th Issue
(1976)}{5th Issue (1976)}
\thispagestyle{headings} ~

\vspace{2ex}

\centerline {\Large \textbf{Problems from the 5th Issue (1976)}}
\phantomsection\label{5izd}
\vspace{5.5ex}

 \bmp \textbf{5.1.} b) Is every
locally finite minimal non-$FC $-group distinct from its derived
 subgroup? \ The
question has an affirmative answer for minimal non-$BFC $-groups.

\hfill \raisebox{-1ex}{\sl V.\,V.\,Belyaev, N.\,F.\,Sesekin}

\emp

 \bmp \textbf{5.5.} If
$G$ is a finitely generated abelian-by-polycyclic-by-finite group,
does there exist a finitely generated metabelian group $M$ such
that $G$ is isomorphic to a subgroup of the automorphism group of
$M$? If so, many of the tricky properties of $G$ like its residual
finiteness would become transparent. \hfill
\raisebox{-1ex}{\sl B.\,A.\,F.\,Wehrfritz}

\emp

 \bmp \textbf{5.14.} An
intersection of some Sylow 2-sub\-groups is called a {\it Sylow
intersection\/} and an intersection of a pair of Sylow
2-sub\-groups is called a {\it paired Sylow intersection}.

 \makebox[23pt][r]{a)} Describe the finite
groups all of whose 2-local subgroups have odd indices.

 \makebox[23pt][r]{b)} Describe the finite groups all of whose
normalizers of
Sylow intersections have odd indices.

 \makebox[23pt][r]{c)} Describe the finite groups all of whose
normalizers of
paired Sylow intersections have odd indices.

 \makebox[23pt][r]{d)} Describe the finite groups in which for any
two Sylow 2-sub\-groups $P$ and $Q$, the intersection $P \cap Q$
is normal in some Sylow 2-sub\-group of $\left< P,\, Q\right> $.

\hfill \raisebox{-1ex}{\sl V.\,V.\,Kabanov,
A.\,A.\,Makhn\"ev,
A.\,I.\,Starostin}

\emp

 \bmp \textbf{5.15.} Do
there exist finitely presented residually finite groups with
recursive, but not primitive recursive, solution of the word
problem? \hfill \raisebox{-1ex}{\sl F.\,B.\,Cannonito}

\emp

 \bmp \textbf{5.16.} Is
every countable locally linear group embeddable in a finitely
presented group? In (G.\,Baumslag, F.\,B.\,Cannonito,
C.\,F.\,Miller III, {\it Math.~Z.}, \textbf{153} (1977), 117--134) it
is proved that every countable group which is locally linear of
bounded degree can be embedded in a finitely presented group with
solvable word problem.

\hfill \raisebox{-1ex}{\sl F.\,B.\,Cannonito,
 C.\,F.\,Miller III}

\emp

 \bmp \textbf{5.25.}
Prove that the factor-group of any soluble linearly ordered group
by its derived subgroup is non-periodic.\hfill \raisebox{-1ex}{\sl
V.\,M.\,Kopytov}

\emp

 \bmp \textbf{5.26.}
Let $G$ be a finite $p\hskip0.2ex $-group with the minimal number
of generators $d$, and let $r_1$ (respectively, $r_2$) be the
minimal number of defining relations on $d$ generators in the
sense of representing $G$ as a factor-group of a free discrete
group (pro-$p\hskip0.2ex $-group). It is well known that always
$r_2 > d^2/4$. For each prime number $p$ denote by $c(p)$ the
exact upper bound for the numbers $b(p)$ with the property that
$r_2 \geq b(p)d^2$ for all finite $p\hskip0.2ex $-groups.

 \makebox[25pt][r]{a)} It is obvious that $r_1 \geq r_2$. Find a
$p\hskip0.2ex $-group with $r_1 > r_2$.

 \makebox[25pt][r]{b)} {\it Conjecture:\/} $\lim\limits_{p\rightarrow
\infty}
c(p)=1/4$. It is proved
(J.\,Wisliceny,
{\it Math.~Nachr.},
\textbf{102} (1981), 57--78) that
$\lim\limits_{d\rightarrow \infty} r_2\left/ d^2\right. =1/4$.\hfill
{\sl H.\,Koch}

\emp

 \bmp \textbf{5.27.}
Prove that if $G$ is a torsion-free pro-$p\hskip0.2ex $-group with
a single defining relation, then ${\rm cd}\, G=2$. This has been
proved for a large class of groups in (J.\,P.\,Labute, {\it Inv.
Math.}, \textbf{4}, no.\,2 (1967), 142--158).\hfill
\raisebox{-1ex}{\sl H.\,Koch}

\emp

 \bmp \textbf{5.30.}
(Well-known problem). Suppose that $G$ is a finite soluble group, $A
\leq {\rm Aut}\, G$, $C_G(A)=1$, the orders of $G$ and $A$ are
coprime, and let $|A|$ be the product of $n$ not necessarily distinct
prime numbers. Is the nilpotent length of $G$ bounded above by~$n$?
This is proved for large classes of groups (E.\,Shult, F.\,Gross,
T.\,Berger, A.\,Turull), and there is a bound in terms of $n$ if $A$
is soluble (J.\,G.\,Thompson's $\leq 5^n$, H.\,Kurzweil's $\leq 4n$,
A.\,Turull's $\leq 2n$), but the problem remains open. \hfill
\raisebox{-1ex}{\sl V.\,D.\,Mazurov}

\emp

 \bmp \textbf{5.33.}
(Y.\,Ihara). Consider the quaternion algebra $Q$ with norm $f=x^2
- \tau y^2 - \rho z^2 + \rho \tau u^2$, $\rho , \tau \in {\Bbb
Z}$. Assume that $f$ is indefinite and of ${\Bbb Q}$-rank 0, i.~e.
$f=0$ for $x,y,z,u \in {\Bbb Q}$ implies $ x=y=z=u=0$. Consider
$Q$ as the algebra of the matrices
$$X=\left(\!\!\begin{array}{cc}x+\sqrt{\tau}y&\rho (z+\sqrt{\tau
}u)\cr z-\sqrt{\tau}u&x-\sqrt{\tau}y\end{array}\!\!\right)$$
with $x,y,z,u \in {\Bbb Q}$. Let $p$ be a prime, $p\nmid \rho
\tau$. Consider the group $G$ of all $X$ with $x,y,z,u \in {\Bbb
Z}^{(p)}$, ${\rm det}\, X=1$, where ${\Bbb Z}^{(p)}=\{ m/p^t\mid
m,t \in {\Bbb Z}\} $.

\makebox[15pt][r]{}{\it Conjecture:\/} $G$ has the congruence
subgroup property, i.~e. every non-central normal subgroup $N$ of
$G$ contains a full congruence subgroup $N({\goth A})=\{ X \in
G\mid X \equiv E\, ({\rm mod}\, {\goth A})\} $ for some ${\goth
A}$. Notice that the congruence subgroup property is independent
of the matrix representation of $Q$. \hfill \raisebox{-1ex}{\sl
J.\,Mennicke}

 \emp

 \bmp \textbf{5.35.}
Let $V$ be a vector space of dimension $n$ over a field. A
subgroup $G$ of $GL_n(V)$ is said to be {\it rich in
transvections\/} if $n \geq 2$ and for every hyperplane $H
\subseteq V$ and every line $L \subseteq H$ there is at least one
transvection in $G$ with residual line $ L$ and fixed space $H$.
Describe the automorphisms of the subgroups of $GL_2(V)$ which are
rich in transvections. \hfill \raisebox{-1ex}{\sl
Yu.\,I.\,Merzlyakov}

\emp

 \bmp \textbf{5.36.}
What profinite groups satisfy the maximum condition for closed
subgroups?

\hfill \raisebox{-1ex}{\sl Yu.\,N.\,Mukhin}

\emp

 \bmp \textbf{\zv 5.38.}
 Is it the case that if $A, B$ are finitely generated soluble Hopfian
groups then $A \times B$ is Hopfian? \hfill \raisebox{-1ex}{\sl
P.\,M.\,Neumann}

\ul

\otv Yes, it is (AI solution confirmed by M.\,Trombetti: \  Achyuth Jayadevan,  {\it Preprint}, 2026,
\url{https://zenodo.org/records/22729301}).
\emp 

 \bmp \textbf{\ai 5.39.} 
Prove that every countable group can act faithfully as a group of
automorphisms of a finitely generated soluble group (of derived
length at most~4). For background to this problem, in particular
its relationship with problem 8.50, see (P.\,M.\,Neumann, {\it
in: Groups--Korea, Pusan, 1988 $($Lect. Notes Math.}, \textbf{1398}),
Springer, Berlin, 1989, 124--139).
 \hfill
\raisebox{-1ex}{\sl P.\,M.\,Neumann}

\ul

\oai This is proved (unverified AI solution:  Achyuth Jayadevan,  {\it Preprint}, 2026, \url{https://zenodo.org/records/22729431}).
\emp 

 \bmp \textbf{5.42.}
Does the free group of rank 2 have an infinite ascending chain of
verbal subgroups each being generated as a verbal subgroup by a
single element?

\hfill \raisebox{-1ex}{\sl A.\,Yu.\,Ol'shan\-ski\u{\i}}

\emp

 \bmp \textbf{5.44.} A
union ${\goth A}=\bigcup\limits_{\alpha} {\goth V}_{\alpha}$ of
varieties of groups (in the lattice of varieties) is called {\it
irreducible\/} if $\bigcup\limits_{\beta\ne \alpha} {\goth
V}_{\beta} \ne {\goth A}$ for each index $\alpha $. Is every
variety an irreducible union of (finitely or infinitely many)
varieties each of which cannot be decomposed into a union of two
proper subvarieties? \hfill \raisebox{-1ex}{\sl
A.\,Yu.\,Ol'shan\-ski\u{\i}}

\emp

\bmp \textbf{\zv 5.47.}
 Is
every countable abelian group embeddable in the center of some
finitely presented group? \hfill
\raisebox{-1ex}{\sl V.\,N.\,Remeslennikov}

\ul

\otv Yes, it is (A.\,O.\,Houcine, {\it J.~Algebra}, {\bf 307}, no.\,1 (2007), 1--23).
\emp

 \bmp \textbf{5.48.}
Suppose that $G$ and $H$ are finitely generated residually finite
groups having the same set of finite homomorphic images. Are $G$ and
$H$ isomorphic if one of them is a free (free soluble) group? \hfill
\raisebox{-1ex}{\sl V.\,N.\,Remeslennikov}

\emp

 \bmp \textbf{\zv 5.52.}
 It
is not hard to show that a finite perfect group is the normal
closure of a single element. Is the same true for infinite
finitely generated groups? \hfill \raisebox{-1ex}{\sl J.\,Wiegold}

\ul

\otv There exist finitely presented perfect groups that are not normal closures of a single element (L.\,Chen, Y.\,Lodha, {\it Preprint}, 2025, \url{https://arxiv.org/abs/2510.26073}).
\emp

 \bmp \textbf{5.54.}
 Let $p,q,r$ be distinct primes and $ u(x,y,z)$ a
commutator word in three variables. Prove that there exist
(infinitely many?) natural numbers $n$ such that the alternating
group ${\Bbb A}_n$ can be generated by three elements $\xi , \eta ,
\zeta$ satisfying $\xi ^p= \eta ^q= \zeta ^r=\xi \eta \zeta \cdot
u(\xi , \eta , \zeta )=1$ \hfill \raisebox{-1ex}{\sl J.\,Wiegold}

\emp

 \bmp \textbf{5.55.}
Find a finite $p\hskip0.2ex $-group that cannot be embedded in a
finite $p\hskip0.2ex $-group with trivial multiplicator. Notice
that every finite group can be embedded in a group with trivial
multiplicator. \hfill \raisebox{-1ex}{\sl J.\,Wiegold}

\emp

 \bmp \textbf{5.56.} a)
Let $p$ be a prime greater than 3. Is it true that every finite
group of exponent~$p$ can be embedded in the commutator subgroup
of a finite group of exponent~$p$?

\hfill \raisebox{-1ex}{\sl J.\,Wiegold}

\emp

 \bmp \textbf{5.59.}
Let $G$ be a locally finite group which is the product of a
$p\hskip0.2ex $-sub\-group and a $q $-sub\-group, where $p$ and
$q$ are distinct primes. Is $G$ a $\{ p,q\} $-group? \hfill
\raisebox{-1ex}{\sl B.\,Hartley}

\emp

 \bmp

 \textbf{5.67.} Is a periodic
residually finite group finite if it satisfies the weak minimum
condition for subgroups? \hfill\raisebox{-1ex}{\sl V.\,P.\,Shunkov}

 \emp \raggedbottom

\newpage
\pagestyle{myheadings} \markboth{6th Issue
(1978)}{6th Issue (1978)}
\thispagestyle{headings} ~ \vspace{2ex}

\centerline {\Large \textbf{Problems from the 6th Issue (1978)}}
\phantomsection\label{6izd}
\vspace{4ex}

 \bmp \textbf{6.1.} A subgroup $H$
of an arbitrary group $G$ is said to be {\it $ C $-closed\/} if
$H=C^2(H)=C(C(H))$, and {\it weakly $C $-closed\/} if $C^2(x) \leq
H$ for any element $x$ of $H$. The structure of the finite groups
all of whose proper subgroups are $ C $-closed was studied by
Gasch\"utz. Describe the structure of the locally finite groups
all of whose proper subgroups are weakly $C $-closed. \hfill
\raisebox{-1ex}{\sl V.\,A.\,Antonov, N.\,F.\,Sesekin}

\emp

 \bmp \textbf{6.2.} The
totality of $C $-closed subgroups of an arbitrary group $G$ is a
complete lattice with respect to the operations
 $\; A \wedge B=A \cap B,\; \; \; A\vee B=C(C(A) \cap
C(B)).$
Describe the groups whose lattice of $C $-closed subgroups is a
sublattice of the lattice of all subgroups.
\hfill \raisebox{-1ex}{\sl V.\,A.\,Antonov,
N.\,F.\,Sesekin}

\emp

 \bmp \textbf{6.3.} 
  A group $G$ is called {\it of type $(FP)_{\infty}$} if the trivial
$G $-module ${\Bbb Z}$ has a resolution by finitely generated
projective $G $-modules. The class of all groups of type
$(FP)_{\infty}$ has a couple of excellent closure properties with
respect to extensions and amalgamated products (R.\,Bieri, {\it
Homological dimension of discrete groups}, Queen Mary College
Math. Notes, London, 1976). Is every periodic group of type
$(FP)_{\infty}$ finite? This is related to the question whether
there is an infinite periodic group with a finite presentation.

\makebox[15pt][r]{}{\it Comment of 2026}: there exist infinite periodic groups with a finite presentation (unverified AI solution: OpenAI, {\it Preprint}, 2026, \url{https://kourovkanotebookorg.wordpress.com/wp-content/uploads/2026/10/8-52.pdf}).
\hfill \raisebox{-1ex}{\sl R.\,Bieri}
\emp 

 \bmp \textbf{6.5.} Is
it the case that every soluble group $G$ of type $(FP)_{\infty}$
is {\it constructible\/} in the sense of (G.\,Baumslag, R.\,Bieri,
{\it Math. Z.}, \textbf{151}, no.\,3 (1976), 249--257)? In other
words, can $G$ be built up from the trivial group in terms of
finite extensions and $HNN $-extensions? \hfill
\raisebox{-1ex}{\sl R.\,Bieri}

\emp

 \bmp \textbf{\zv 6.9.}
 Is
the derived subgroup of any locally normal group decomposable into
a product of at most countable subgroups which commute
elementwise?
\hfill \raisebox{-1ex}{\sl Yu.\,M.\,Gorchakov}

\ul

\otv Yes, it is (AI solution confirmed by M.\,Dixon: \  Achyuth Jayadevan,  {\it Preprint}, 2026,
\url{https://zenodo.org/records/22729504}).
\emp 

 \bmp \textbf{6.10.}
Let $p$ be a prime and $n$ be an integer with $p > 2n + 1$. Let
$x, y$ be $p\hskip0.2ex $-elements in $GL_n({\Bbb C})$. If the
subgroup $\left< x,y\right>$ is finite, then it is an abelian
$p\hskip0.2ex $-group (W.\,Feit, J.\,G.\,Thompson, {\it Pacif.
J.~Math.}, \textbf{11}, no.\,4 (1961), 1257--1262). What can be said
about $\left< x,y\right> $ in the case it is an infinite group?
\hfill \raisebox{-1ex}{\sl J.\,D.\,Dixon}

\emp

 \bmp \textbf{6.11.}
Let ${\goth L}$ be the class of locally compact groups with no
small subgroups (see D.\,Montgomery, L.\,Zippin, {\it Topological
transformation groups}, New York, 1955; V.\,M.\,Glushkov, {\it
Uspekhi Matem. Nauk}, \textbf{12}, no.\,2 (1957), 3--41 (Russian)).
Study extensions of groups in this class with the objective of
giving a direct proof of the following: For each $G \in {\goth L}$
there exists an $H \in {\goth L}$ and a (continuous) homomorphism
$\vartheta \!:G \rightarrow H$ with a discrete kernel and an image
${\rm Im}\, \vartheta$ satisfying ${\rm Im}\, \vartheta \cap
Z(H)=1$ ($Z(H)$ is the center of $H$). This result follows from
the Gleason--Montgomery--Zippin solution of Hilbert's 5th Problem
(since ${\goth L}$ is the class of finite-dimensional Lie groups
and the latter are locally linear). On the other hand a direct
proof of this result would give a substantially shorter proof of
the 5th Problem since the adjoint representation of $H$ is
faithful on ${\rm Im}\, \vartheta$. \hfill \raisebox{-1ex}{\sl
J.\,D.\,Dixon}

\emp

 \bmp \textbf{6.21.}
G.\,Higman proved that, for any prime number $p$, there exists a
natural number $\chi (p)$ such that the nilpotency class of any
finite group $G$ having an automorphism of order $p$ without
non-trivial fixed points does not exceed~$\chi (p)$. At the same time
he showed that $\chi (p) \geq (p^2 - 1)/4$ for any such {\it Higman's
function\/} $\chi$. Find the best possible Higman's function. Is it
the function defined by equalities $\chi (p)=(p^2- 1)/4$ for $p > 2$
and $\chi (2)=1$? This is known to be true for $p\leq 7$. \hfill
\raisebox{-1ex}{\sl V.\,D.\,Mazurov}

\emp

 \bmp \textbf{6.24.}
The membership problem for the braid group on four strings. It is
known (T.\,A.\,Makanina, {\it Math. Notes}, \textbf{29}, no.\,1
(1981), 16--17) that the membership problem is undecidable for
braid groups with more than four strings. \hfill
\raisebox{-1ex}{\sl G.\,S.\,Makanin}

\emp

 \bmp \textbf{6.26.}
Let $D$ be a normal set of involutions in a finite group $G$ and
let $\Gamma (D)$ be the graph with vertex set $D$ and edge set $\{
(a,b)\mid a,b \in D,\; ab=ba \ne 1\} $. Describe the finite groups
$G$ with non-connected graph $\Gamma (D)$. \hfill
\raisebox{-1ex}{\sl A.\,A.\,Makhn\"ev}

\emp

 \bmp \textbf{6.28.}
Let $A$ be an elementary abelian 2-group which is a $ TI
$-sub\-group of a~finite group $G$. Investigate the structure of
$G$ under the hypothesis that the weak closure of $A$ in a Sylow
2-sub\-group of $G$ is abelian. \hfill \raisebox{-1ex}{\sl
A.\,A.\,Makhn\"ev}

\emp

 \bmp \textbf{6.29.}
Suppose that a finite group $A$ is isomorphic to the group of all
topological automorphisms of a locally compact group $G$. Does
there always exist a discrete group whose automorphism group is
isomorphic to $A$? This is true if $A$ is cyclic; the condition of
local compactness of $G$ is essential (R.\,J.\,Wille, {\it Indag.
Math.}, \textbf{25}, no.\,2 (1963), 218--224). \hfill
\raisebox{-1ex}{\sl O.\,V.\,Mel'nikov}

\emp

 \bmp \textbf{\zv 6.30.}
Let $G$ be a residually finite Hopfian group, and let
$\widehat{G}$ be its profinite completion. Is $\widehat{G}$
necessarily Hopfian (in topological sense)? \hfill
\raisebox{-1ex}{\sl O.\,V.\,Mel'nikov}

\ul

\otv No, not always (M.\,Brescia, E.\,Ingrosso, M.\,Trombetti, {\it Preprint}, 2026, \url{https://arxiv.org/abs/2607.19193}).
\emp

 \bmp \textbf{6.31.} b)
Let $G$ be a finitely generated residually finite group, $d(G)$
the minimal number of its generators, and $\delta (G)$ the minimal
number of topological generators of the profinite completion of
$G$. It is known that there exist groups
 $G$ for which $d(G) > \delta
(G)$ (G.\,A.\,Noskov, {\it Math. Notes}, \textbf{33}, no.\,4 (1983),
249--254). Is the function $d$ bounded on the set of groups $G$
with the fixed value of $\delta (G) \geq 2$? \hfill
\raisebox{-1ex}{\sl O.\,V.\,Mel'nikov}

\emp

 \bmp \textbf{6.32.}
Let $F_n$ be the free profinite group of finite rank $n
> 1$. Is it true that for each normal subgroup $N$ of the free
profinite group of countable rank, there exists a normal
subgroup of $F_n$ isomorphic to $N$?
\hfill \raisebox{-1ex}{\sl O.\,V.\,Mel'nikov}

\emp

 \bmp \textbf{6.38.}
b) Is it true that
an arbitrary subgroup of $GL_n(k)$
that intersects every conjugacy class is parabolic? How far is the
same statement true for subgroups of other groups of Lie type?
\hfill \raisebox{-1ex}{\sl P.\,M.\,Neumann}
 \emp

 \bmp \textbf{6.39.} A
class of groups ${\goth K}$ is said to be {\it radical\/} if it is
closed under taking homomorphic images and normal subgroups and if
every group generated by its normal ${\goth K} $-sub\-groups
also belongs to ${\goth K}$. The question proposed relates to the
topic ``Radical classes and formulae of Narrow (first order)
Predicate Calculus (NPC)''. One can show that the only non-trivial
radical class definable by universal formulae of NPC is the class
of all groups. Recently (in a letter to me), G.\,M.\,Bergman
constructed a family of locally finite radical classes definable
by formulae of NPC. Do there exist similar classes which are not
locally finite and different from the class of all groups? In
particular, does there exist a radical class of groups which is
closed under taking Cartesian products, contains an infinite
cyclic group, and is different from the class of all groups?
\makebox[15ex][r]{} \hfill \raisebox{-1ex}{\sl B.\,I.\,Plotkin}

\emp

 \bmp \textbf{\zv 6.45.}
Construct a characteristic subgroup $N$ of a finitely generated
free group $F$ such that $F/N$ is infinite and simple. If no such
exists, it would follow that $d(S^2)=d(S)$ for every infinite
finitely generated simple group $S$. There is reason to believe
that this is false. \hfill \raisebox{-1ex}{\sl J.\,Wiegold}

\ul

\otv
It is proved that for all $n\geq 2$, the free group $F_n$ admits continuum many pairwise non-isomorphic infinite simple characteristic quotients (R.\,Coulon, F.\,Fournier-Facio, \url{https://arxiv.org/abs/2312.11684}).
\emp

 \bmp \textbf{6.47.}
(C.\,D.\,H.\,Cooper). Let $G$ be a group, $v$ a group word in two
variables such that the operation $x \odot y=v(x,y)$ defines the
structure of a new group $G_v=\left< G,\odot \right>$ on the set
$G$. Does $G_v$ always lie in the variety generated by $G$? \hfill
\raisebox{-1ex}{\sl E.\,I.\,Khukhro}

\emp

 \bmp \textbf{6.48.} Is
every infinite binary finite $p\hskip0.2ex $-group non-simple?
Here $p$ is a prime number.

\hfill \raisebox{-1ex}{\sl N.\,S.\,Chernikov}

\emp

 \bmp \textbf{6.51.}
Let ${\goth F}$ be a local subformation of some formation ${\goth
X}$ of finite groups and let $\Omega$ be the set of all maximal
homogeneous ${\goth X} $-screens of ${\goth F}$. Find a way of
constructing the elements of $\Omega$ with the help of the maximal
inner local screen of ${\goth F}$. What can be said about the
cardinality of $\Omega$? For definitions see (L.\,A.\,Shemetkov,
{\it Formations of finite groups}, Moscow, Nauka, 1978 (Russian)).
\hfill \raisebox{-1ex}{\sl L.\,A.\,Shemetkov}

\emp

\bmp \textbf{\ai 6.55.} 
 A group $G$ of the form $G=F\leftthreetimes H$ is
said to be a {\it Frobenius group with kernel $F$ and complement}
$H$ if $H \cap H^g=1$ for any $g \in G\setminus H$ and $F\setminus
\{ 1\} =G\setminus \bigcup\limits_{g\in G} H^g$. Do there exist
Frobenius $p\hskip0.2ex $-groups? \hfill
\raisebox{-1ex}{\sl V.\,P.\,Shunkov}

\ul

\oai Yes, they do exist (unverified AI solution: B.\,G.\,Di Siena, {\it Preprint}, 2026, \url{https://www.proofforum.org/abs/PF:2026.000032}).
\emp

 \bmp \textbf{6.56.}
Let $G=F\cdot \left< a\right>$ be a Frobenius group with the
complement $\left< a\right> $ of prime order.

 \makebox[25pt][r]{a)} Is $G$ locally finite if it is binary finite?

 \makebox[25pt][r]{b)} Is the kernel $F$ locally finite if the
groups $\left< a,\, a^g\right> $ are finite for all $g \in G$?

\hfill \raisebox{-1ex}{\sl V.\,P.\,Shunkov}

\emp

 \bmp \textbf{6.59.}
 A group $G$ is said to be ({\it conjugacy, $p$-con\-ju\-gacy\/})
{\it biprimitively finite\/} if, for any finite subgroup $H$,
any two elements of prime order (any two conjugate elements of
prime order, of prime order $p$) in $N_G(H)/H$ generate a finite
subgroup. Prove that an arbitrary periodic (conjugacy)
biprimitively finite
group (in particular, having no involutions) of finite rank is
locally finite. \hfill \raisebox{-1ex}{\sl V.\,P.\,Shunkov}

\emp

 \bmp \textbf{6.60.} Do
there exist infinite finitely generated simple periodic
(conjugacy) biprimitively finite groups which contain both
involutions and non-trivial elements of odd order? \hfill
\raisebox{-1ex}{\sl V.\,P.\,Shunkov}

\emp

 \bmp \textbf{6.61.} Is
every infinite periodic (conjugacy) biprimitively finite group
without involutions non-simple? \hfill
\raisebox{-1ex}{\sl V.\,P.\,Shunkov}

\emp

 \bmp \textbf{6.62.} Is
a (conjugacy, $p\hskip0.2ex $-conju\-ga\-cy) biprimitively finite
group finite if it has a finite maximal subgroup ($p\hskip0.2ex
$-sub\-group)? There are affirmative answers for conjugacy
biprimitively finite $p\hskip0.2ex $-groups and for 2-conjugacy
biprimitively finite groups (V.\,P.\,Shunkov, {\it Algebra and
Logic}, \textbf{9}, no.\,4 (1970), 291--297; \textbf{11}, no.\,4 (1972),
260--272; \textbf{12}, no.\,5 (1973), 347--353), while the statement
does not hold for arbitrary periodic groups, see Archive,~3.9.
\hfill \raisebox{-1ex}{\sl V.\,P.\,Shunkov}

\emp

\raggedbottom

\newpage
\pagestyle{myheadings} \markboth{7th Issue
(1980)}{7th Issue (1980)}
\thispagestyle{headings}

~\vspace{2ex}

\centerline {\Large \textbf{Problems from the 7th Issue (1980)}}
\phantomsection\label{7izd}
\vspace{4ex}

 \bmp \textbf{7.3.}
Prove that the periodic product of odd exponent $n \geq 665$ of
non-trivial groups $F_1,\ldots ,F_k$ which do not contain
involutions cannot be generated by less than $k$ elements. This
would imply, on the basis of (S.\,I.\,Adian, {\it Sov. Math.
Doklady}, \textbf{19}, (1978), 910--913), the existence of $ k
$-gene\-ra\-ted simple groups which cannot be generated by less
than $k$ elements, for any $k > 0$. \hfill
\raisebox{-1ex}{\sl S.\,I.\,Adian}

\emp

 \bmp \textbf{7.5.} We
say that a group is {\it indecomposable\/} if any two of its
proper subgroups generate a proper subgroup. Describe the
indecomposable periodic metabelian groups.

\hfill \raisebox{-1ex}{\sl V.\,V.\,Belyaev} \emp

 \bmp \textbf{\zv 7.15.}
Prove that if $F^{\displaystyle{*}} (G)$ is quasisimple and
$\alpha \in {\rm Aut}\, G$, $|\alpha |=2$, then $C_G(\alpha )$
contains an involution outside $Z(F^{\displaystyle{*}}(G))$,
except when $F ^{\displaystyle{*}}(G)$ has quaternion Sylow
2-sub\-groups. \hfill \raisebox{-1ex}{\sl R.\,Griess}

\ul

\otv This is not always true (AI solution confirmed by A.\,Rizzoli: Achyuth Jayadevan, {\it Preprint}, 2026,
\url{https://zenodo.org/records/22750993}).
\emp

\bmp \textbf{7.19.} Construct an explicit example of a finitely
presented simple group with word problem not solvable by a
primitive recursive function. \hfill \raisebox{-1ex}{\sl
F.\,B.\,Cannonito}

\emp

 \bmp \textbf{7.21.} Is
there an algorithm which decides if an arbitrary finitely
presented solvable group is metabelian? \hfill
\raisebox{-1ex}{\sl F.\,B.\,Cannonito}

\emp

 \bmp \textbf{7.23.}
Does there exist an algorithm which decides, for a given list of
group identities, whether the elementary theory of the variety
given by this list is soluble, that is, by (A.\,P.\,Zamyatin, {\it
Algebra and Logic}, \textbf{17}, no.\,1 (1978), 13--17), whether this
variety is abelian? \hfill \raisebox{-1ex}{\sl A.\,V.\,Kuznetsov}

\emp

\bmp
 \textbf{7.25.} We associate with words in the alphabet $x$, $x^{-1}$,
$y$, $y^{-1}$, $z$, $z^{-1},\ldots $ operations which are
understood as functions in variables $x, y, z,\ldots $. We say
that a word $A$ {\it is expressible in words $B_1,\ldots ,B_n$ on
the group $G$\/} if $A$ can be constructed from the words
$B_1,\ldots ,B_n$ and variables by means of finitely many
substitutions of words one into another and replacements of a word
by another word which is identically equal to it on $G$. A list of
words is said to be {\it functionally complete on $G$\/} if every
word can be expressed on $G$ in words from this list. A word $B$
is called a {\it Schaeffer word on $G$\/} if every word can be
expressed on $G$ in $B$ (compare with A.\,V.\,Kuznetsov, {\it
Matem. Issledovaniya}, Kishin\"ev, \textbf{6}, no.\,4 (1971), 75--122
(Russian)). Does there exist an algorithm which decides

 \makebox[23pt][r]{a)} by a word $B$, whether it is a Schaeffer word
on every
group? Compare with the problem of describing all such words in
(A.\,G.\,Kurosh, {\it Theory of Groups},
Moscow, Nauka, 1967,
p.~435 (Russian)); here are examples of
such words: $xy^{-1}$, $x^{-1}y^2z$, $x^{-1}y^{-1}zx$.

\makebox[23pt][r]{b)} by a list of words, whether it is functionally
complete on
every group?

 \makebox[23pt][r]{c)} by words $A, B_1,\ldots ,B_n$, whether $A$ is
expressible in $ B_1,\ldots ,B_n$ on every group? (This is a
problem from A.\,V.\,Kuznetsov, {\it ibid.}, p.\,112.)

 \makebox[23pt][r]{d)} the same for every finite group? (For finite
groups, for example, $x^{-1}$ is expressible in $xy$.) For a fixed
finite group an algorithm exists (compare with A.\,V.\,Kuznetsov,
{\it ibid.}, $\S$\,8). \hfill \raisebox{-1ex}{\sl
A.\,V.\,Kuznetsov}
\emp

 \bmp \textbf{7.26.}
The {\it ordinal height\/} of a variety of groups is, by
definition, the supremum of order types (ordinals) of all
well-ordered by inclusion chains of its proper subvarieties. It is
clear that its cardinality is either finite or countably infinite,
or equal to $\omega _1$. Is every countable ordinal number the
ordinal height of some variety?

\hfill \raisebox{-1ex}{\sl A.\,V.\,Kuznetsov}

\emp

\bmp \textbf{\zv 7.27.}
Is it true that the group $SL_n(q)$ contains, for
sufficiently large $q$, a diagonal matrix which is not contained in
any proper irreducible subgroup of $SL_n(q)$ with the exception of
block-monomial ones? For $n=2,\, 3$ the answer is known to be
affirmative (V.\,M.\,Levchuk, {\it in: Some problems of the theory of
groups and rings}, Inst. of Physics SO AN SSSR, Krasnoyarsk, 1973
(Russian)). Similar problems are interesting for other Chevalley
groups. \hfill \raisebox{-1ex}{\sl V.\,M.\,Levchuk}

\ul

\otv Yes, it is true, as follows from  
the density theorems of J.\,Fulman and R.\,Guralnick
(AI solution confirmed by R.\,Guralnick: \ A.\,Ferudun, {\it Preprint}, 2026, \url{https://doi.org/10.5281/zenodo.23445980}).
\emp 

\bmp
 \textbf{7.28.}
Let $G(K)$ be the Chevalley group over a commutative ring $K$
associated with the root system $\Phi$ as defined in (R.\,Steinberg,
{\it Lectures on Chevalley groups}, Yale Univ., New Haven, Conn.,
1967). This group is generated by the root subgroups $x_r(K)$, $r \in
\Phi $. We define an {\it elementary carpet of type $\Phi$ over
$K$\/} to be any collection of additive subgroups $\{ {\goth A}_r\mid
r \in \Phi \}$ of $K$ satisfying the condition \vspace{-1ex}
$$c_{ij,rs}{\goth A}^i_r{\goth A}^j_s \subseteq {\goth A}_{ir+js}
\quad \quad \mbox{ for}\; r,s,ir + js \in \Phi , \;\; i > 0,\;\; j >
0,$$ \vskip-1ex where $c_{ij,rs}$ are constants defined by the
Chevalley commutator formula and ${\goth A}_r^i=\{ a^i\mid a \in
{\goth A}_r\} $. What are necessary and sufficient conditions (in
terms of the~${\goth A}_r$) on the elementary carpet to ensure that
the subgroup $\left< x_r({\goth A}_r)\mid r \in \Phi \right>$
of~$G(K)$ intersects with $x_r(K)$ in $x_r({\goth A}_r)$? See also
15.46.\hfill \raisebox{-1ex}{\sl V.\,M.\,Levchuk}

\emp

\bmp \textbf{\zv 7.31.}
Let $A$ be a group of automorphisms of a finite
non-abelian 2-group $G$ acting transitively on the set of involutions
of $G$. Is $A$ necessarily soluble?

\makebox[15pt][r]{}Such groups $G$ were divided into several classes in
(F.\,Gross, {\it J.\,Algebra}, \textbf{40}, no.\,2 (1976), 348--353).
For one of these classes a positive answer was given by
E.\,G.\,Bryukhanova ({\it Algebra and Logic}, \textbf{20}, no.\,1
(1981), 1--12). \hfill \raisebox{-1ex}{\sl V.\,D.\,Mazurov}

\ul

\otv Yes, it is (C.\,H.\,Li, Y.\,Z.\, Zhu, {\it Proc. London Math. Soc. (3)}, {\bf 129}, no.\,5 (2024), Article ID e70002).
\emp

\bmp
 \textbf{7.33.}
 An
elementary $TI $-sub\-group $V$ of a finite group $G$ is said to
be a {\it subgroup of non-root type\/} if $1 \ne N_V(V^g) \ne V$
for some $g \in G$. Describe the finite groups $G$ which contain a
2-sub\-group $V$ of non-root type such that $[V,\, V^g]=1$ implies
that all involutions in $VV^g$ are conjugate to elements of~$V$.
\makebox[15ex][r]{} \hfill \raisebox{-1ex}{\sl A.\,A.\,Makhn\"ev}

\emp

\bmp
 \textbf{7.34.} In
many of the sporadic finite simple groups the 2-ranks of the
centralizers of 3-elements are at most 2. Describe the finite
groups satisfying this condition.

\hfill \raisebox{-1ex}{\sl A.\,A.\,Makhn\"ev}

\emp

\bmp
 \textbf{7.35.}
What varieties of groups ${\goth V}$ have the following property:
the group $G/{\goth V}(G)$ is residually finite for any residually
finite group $G$? \hfill \raisebox{-1ex}{\sl O.\,V.\,Mel'nikov}

\emp

\bmp
 \textbf{7.38.}
 A~{\it variety of profinite groups\/} is a
non-empty class of profinite groups closed under taking subgroups,
factor-groups, and Tikhonov products. A subvariety ${\goth V}$ of the
variety ${\goth N}$ of all pronilpotent groups is said to be {\it
$($locally\/$)$ nilpotent\/} if all (finitely generated) groups in
${\goth V}$ are nilpotent.

 \makebox[25pt][r]{\zv a)} Is it true that any non-nilpotent subvariety of
${\goth N}$ contains
a non-nilpotent locally nilpotent subvariety?

 \makebox[25pt][r]{b)} The same problem for the variety of all
pro-$p\hskip0.2ex $-groups (for a given prime~$p$).

\hfill \raisebox{-1ex}{\sl O.\,V.\,Mel'nikov}

\ul

\otv No, not always (AI solution confirmed by E.\,I.\,Khukhro: \ Achyuth Jayadevan, {\it Preprint}, 2026, \url{https://zenodo.org/records/22751322}).
\emp

\bmp
 \textbf{7.39.}
Let $G=\left< a,\, b\mid a^p=(ab)^3=b^2=(a^{\sigma}ba^{2/\sigma
}b)^2=1\right> ,$ where $p$ is a prime, $\sigma$ is an integer not
divisible by $p$. The group $PSL_2(p)$ is a factor-group of $G$ so
that there is a short exact sequence $1 \rightarrow N \rightarrow
G \rightarrow PSL_2(p) \rightarrow 1 .$ For each $p > 2$ there is
$\sigma$ such that $N=1$, for example, $\sigma =4$. Let $N^{ab}$
denote the factor-group of $N$ by its commutator subgroup. It is
known that for some $p$ there is $\sigma$ such that $N^{ab}$ is
infinite (for example, for $p=41$ one can take $\sigma ^2 \equiv
2\, ({\rm mod}\, 41)$), whereas for some other $p$ (for example,
for $p=43$) the group $N^{ab}$ is finite for every $\sigma$.

 \makebox[25pt][r]{a)} Is the set of primes $p$ for which $N^{ab}$
is finite for every $\sigma$ infinite?

 \makebox[25pt][r]{b)} Is there an arithmetic condition on $\sigma$
which
ensures that $N^{ab}$ is finite?

\hfill \raisebox{-1ex}{\sl J.\,Mennicke}

\emp

\bmp
 \textbf{7.41.}
(John S.\,Wilson). Is every linear $\,\overline{\! SI\vphantom{(}}
$-group an $\,\overline{\! SN\vphantom{(}} $-group? \hfill
{\sl Yu.\,I.\,Merzlyakov}

\emp

\bmp
 \textbf{7.45.}
Does the $Q\hskip0.1ex $-theory of the class of all finite groups
(in the sense of 2.40) coincide with the $Q\hskip0.1ex
$-theory of a single finitely presented group?
\hfill \raisebox{-1ex}{\sl
D.\,M.\,Smirnov}

\emp

\bmp
 \textbf{7.49.}
Let $G$ be a finitely generated group, and $N$ a minimal normal
subgroup of $G$ which is an elementary abelian $p\hskip0.2ex
$-group. Is it true that either $N$ is finite or the growth
function for the elements of $N$ with respect to a finite
generating set of $G$ is bounded below by an exponential function?
\hfill \raisebox{-1ex}{\sl V.\,I.\,Trofimov}

\emp

\bmp
 \textbf{7.50.}
Study the structure of the primitive permutation groups (finite
and infinite), in which the stabilizer of any three pairwise
distinct points is trivial. This problem is closely connected with
the problem of describing those group which have a Frobenius group
as one of maximal subgroups. \hfill
\raisebox{-1ex}{\sl A.\,N.\,Fomin}

\emp

\bmp
 \textbf{7.51.}
What are the primitive permutation groups (finite and infinite)
which have a regular sub-orbit, that is, in which a point
stabilizer acts faithfully and regularly on at least one of its
orbits? \hfill \raisebox{-1ex}{\sl A.\,N.\,Fomin}

\emp

\bmp
 \textbf{7.52.}
Describe the locally finite primitive permutation groups in which
the centre of any Sylow 2-sub\-group contains involutions
stabilizing precisely one symbol. The case of finite groups was
completely determined by D.\,Holt in 1978. \hfill
\raisebox{-1ex}{\sl A.\,N.\,Fomin}

\emp

\bmp

\textbf{7.54.} Does there exist a group of infinite special rank which
can be represented as a product of two subgroups of finite special
rank? \hfill \raisebox{-1ex}{\sl N.\,S.\,Chernikov}

\emp

\bmp
 \textbf{7.55.} Is
it true that a group which is a product of two almost abelian
subgroups is almost soluble? \hfill
\raisebox{-1ex}{\sl N.\,S.\,Chernikov}

\emp

\bmp

\textbf{7.56.} (B.\,Amberg). Does a group satisfy the minimum
(respectively, maximum) condition on subgroups if it is a product of
two subgroups satisfying the minimum (respectively, maximum)
condition on subgroups? \hfill \raisebox{-1ex}{\sl N.\,S.\,Chernikov}

\emp

\bmp
 \textbf{7.57.} a)
A set of generators of a finitely generated group $G$ consisting
of the least possible number $d(G)$ of elements is called a {\it
basis of $G$}. Let $r_M(G)$ be the least number of relations
necessary to define $G$ in the basis $M$ and let $r(G)$ be the
minimum of the numbers $r_M(G)$ over all bases $M$ of $G$. It is
known that $r_M(G) \leq d(G) + r(G)$ for any basis $M$. Does there exist a finitely presented group $G$ for which the inequality becomes equality for some basis~$M$?
\hfill
\raisebox{-1ex}{\sl V.\,A.\,Churkin}

\emp

\bmp
 \textbf{\ai 7.58.} 
Suppose that $F$ is an absolutely free group, $R$ a normal
subgroup of $F$, and let ${\goth V}$ be a variety of groups. It is
well-known (H.\,Neumann, {\it Varieties of Groups}, Springer, Berlin,
 1967) that the group $F/{\goth V}(R)$ is isomorphically
embeddable in the ${\goth V} $-verbal wreath product of a ${\goth
V} $-free group of the same rank as $F$ with $F/R$. Find a
criterion indicating which elements of this wreath product belong
to the image of this embedding. A~criterion is known in the case
where ${\goth V}$ is the variety of all abelian groups
(V.\,N.\,Remeslennikov, V.\,G.\,Sokolov, {\it Algebra and Logic},
\textbf{9}, no.\,5 (1970), 342--349). \makebox[15ex][r]{} \hfill
\raisebox{-1ex}{\sl G.\,G.\,Yabanzhi}

\ul

\oai A criterion is found (unverified AI solution:  Achyuth Jayadevan,  {\it Preprint}, 2026, \url{https://zenodo.org/records/22821158}).
\emp

\raggedbottom

\newpage

\pagestyle{myheadings} \markboth{8th Issue
(1982)}{8th Issue (1982)}
\thispagestyle{headings}

~
\vspace{2ex}

\centerline {\Large \textbf{Problems from the 8th Issue (1982)}}
\phantomsection\label{8izd}
\vspace{4ex}

\bmp \textbf{8.1.} Characterize all groups (or at least all soluble
groups) $G$ and fields $F$ such that every irreducible $ FG
$-module has finite dimension over the centre of its endomorphism
ring. For some partial answers see (B.\,A.\,F.\,Wehrfritz, {\it
Glasgow Math.~J.}, \textbf{24}, no.\,1 (1983), 169--176). \hfill
\raisebox{-1ex}{\sl B.\,A.\,F.\,Wehrfritz}\emp

\bmp \textbf{8.2.} Let $G$ be the amalgamated free product of two
polycyclic groups, amalgamating a normal subgroup of each. Is $G$
isomorphic to a linear group? $G$ is residually finite
(G.\,Baumslag, {\it Trans. Amer. Math. Soc.}, \textbf{106}, no.\,2
(1963), 193--209) and the answer is ``yes'' if the amalgamated part
is torsion-free, abelian (B.\,A.\,F.\,Wehrfritz, {\it Proc. London
Math. Soc.}, \textbf{27}, no.\,3 (1973), 402--424) and nilpotent
(M.\,Shirvani, 1981, {\it unpublished}). \hfill
\raisebox{-1ex}{\sl B.\,A.\,F.\,Wehrfritz}\emp

\bmp \textbf{8.3.} Let $m$ and $n$ be positive integers and $p$ a
prime. Let $P_n({\Bbb Z}_{p^m})$ be the group of all $n \times n$
matrices $(a_{ij})$ over the integers modulo $p^m$ such that
$a_{ii} \equiv 1$ for all $i$ and $a_{ij} \equiv 0 \, ({\rm mod}\,
p)$ for all $i > j$. The group $P_n({\Bbb Z}_{p^m})$ is a finite
$p\hskip0.2ex $-group. For which $m$ and $n$ is it regular?

\makebox[15pt][r]{}{\it Editors' comment (2001):\/} The group
$P_{n}({\Bbb Z}_{p^m})$ is known to be regular if $mn<p$
(Yu.\,I.\,Merzlyakov, {\it Algebra i Logika}, \textbf{3}, no.\,4
(1964), 49--59 (Russian)). The case $m=1$ is completely done by
A.\,V.\,Yagzhev in ({\it Math. Notes}, \textbf{56}, no.\,5--6 (1995),
1283--1290), although this paper is erroneous for $m>1$. The case
$m=2$ is completed in (S.\,G.\,Kolesnikov, {\it Issledovaniya in
analysis and algebra}, no.\,3, Tomsk Univ., Tomsk, 2001, 117--125
(Russian)). {\it Editors' comment (2009):\/}
The group $P_{n}({\Bbb Z}_{p^m})$ was shown to be regular for any
$m$ if
$n^2<p$ (S.\,G.\,Kolesnikov, {\it 
Siberian Math.~J.}, \textbf{47}, no.\,6 (2006), 1054--1059). \hfill
\raisebox{-1ex}{\sl B.\,A.\,F.\,Wehrfritz}

\emp

\bmp \textbf{8.4.} Construct a finite nilpotent loop with no finite
basis for its laws.

\hfill \raisebox{-1ex}{\sl M.\,R.\,Vaughan-Lee}
\emp

\bmp

\textbf{8.8.} b) (D.\,V.\,Anosov). Does there exist a non-cyclic
finitely presented group $G$ which contains an element $a$ such
that each element of $G$ is conjugate to some power of $a$?

\hfill
\raisebox{-1ex}{\sl R.\,I.\,Grigorchuk} \emp

\bmp

\textbf{8.9.} (C.\,Chou). We say that a group $G$ has {\it property
$P$\/} if for every finite subset $F$ of~$G$, there is a finite
subset $S \supset F$ and a subset $X \subset G$ such that $x_1S
\cap x_2S$ is empty for any $x_1, x_2 \in X$, $x_1 \ne x_2$, and
$G=\bigcup\limits_{x\in X}Sx$. Does every group have property $P$?

\hfill {\sl R.\,I.\,Grigorchuk}

 \emp

 \bmp
\textbf{\zv 8.11.}
Consider the group \vspace{-1ex} $$M=\left<
x,y,z,t\mid [x,y]=[y,z]=[z,x]=(x,t)=(y,t)=(z,t)=1\right> $$
\vskip-1ex

where $[x, y]=x^{-1}y^{-1}xy$ and
$(x,t)=x^{-1}t^{-1}x^{-1}txt$. The subgroup $H=\left<
x,t,y\right>$ is isomorphic to the braid group ${\goth B}_4$,
and is normally complemented by $N=\left< (zx^{-1})^M\right>
$. Is $N$ a free group (?) of countably infinite rank (?) on
which $H$ acts faithfully by conjugation?
\hfill \raisebox{-1ex}{\sl D.\,L.\,Johnson}

\ul

\otv
As pointed by Y.\,Antol\'{i}n, the group $M$ is the Artin group of type $D_4$; it is proved in (B.\,Perron, J.\,P.\,Vannier, {\it Math. Ann.}, {\bf 306}, no.\,2 (1996), 231--246) that $N$ is a free group of rank 3 and $H$ acts on $N$ faithfully by conjugation.
 \emp

\bmp \textbf{8.12.} Let $D0$ denote the class of finite groups of
deficiency zero, i.~e. having a presentation $\left< X\mid
R\right>$ with $|X|=|R|$.

 \makebox[23pt][r]{a)} Does $D0$ contain any 3-generator
$p\hskip0.2ex $-group for $p \geq 5$?

 \makebox[23pt][r]{c)} Do soluble $D0\hskip0.1ex $-groups have bounded
derived length?

 \makebox[23pt][r]{d)} Which non-abelian simple groups can occur
as composition
factors of $D0\hskip0.1ex $-groups? \hfill
\raisebox{-1ex}{\sl D.\,L.\,Johnson,
E.\,F.\,Robertson}
\emp

\bmp \textbf{8.14.} b) Assume a group $G$ is existentially closed in one
of the classes $L{\goth N}^{+}$, \ $L{\goth S}_{\pi}$, \ $L{\goth
S}^{+}$, \ $L{\goth S}$ of, respectively, all locally nilpotent
torsion-free groups, all locally soluble $\pi $-groups, all locally
soluble torsion-free groups, or all locally soluble groups. Is it
true that $G$ is characteristically simple?

 \makebox[15pt][r]{}The existential closedness of $G$ is a local
property, thus it
seems difficult to obtain global properties of $G$ from it. See
also Archive,~8.14\,a.\hfill
\raisebox{-1ex}{\sl O.\,H.\,Kegel} \emp

\bmp \textbf{8.15.} (Well-known problem). Is every $\widetilde{N}
$-group a $\,\overline{\! Z} $-group? For definitions, see
(M.\,I.\,Kargapolov, Yu.\,I.\,Merzlyakov, {\it Fundamentals of the
Theory of Groups}, 3rd Edition, Moscow, Nauka, 1982, p.\,205
(Russian)). \hfill \raisebox{-1ex}{\sl Sh.\,S.\,Kemkhadze} \emp

\bmp \textbf{\zv 8.16.}
Is the class of all $\,\overline{\! Z} $-groups
closed under taking normal subgroups?

\hfill \raisebox{-1ex}{\sl Sh.\,S.\,Kemkhadze}

\ul

\otv No, it is not (AI solution confirmed by J.\,S.\,Wilson:  \ Achyuth Jayadevan,  {\it Preprint}, 2026, \url{https://zenodo.org/records/22821641}).
\emp 

\bmp \textbf{8.19.} Which of the following properties of (soluble)
varieties of groups are equivalent to one another:

 \makebox[15pt][r]{}1) to satisfy the minimum condition for
subvarieties;

 \makebox[15pt][r]{}2) to have at most countably
many subvarieties;

 \makebox[15pt][r]{}3) to have no infinite independent system of
identities? \hfill {\sl Yu.\,G.\,Kleiman}
\emp

\bmp \textbf{8.21.} If the commutator subgroup $G'$ of a relatively
free group $G$ is periodic, must $G'$ have finite exponent? \hfill
\raisebox{-1ex}{\sl L.\,G.\,Kov\'acs}

\emp

\bmp \textbf{8.23.} If the dihedral group $D$ of order 18 is a
section of a direct product \mbox{$A \times B$}, must at least one
of $A$ and $B$ have a section isomorphic to $D$? \hfill
\raisebox{-1ex}{\sl L.\,G.\,Kov\'acs}

\emp

\bmp \textbf{8.24.} Prove that every linearly ordered group of finite
rank is soluble.

\hfill \raisebox{-1ex}{\sl V.\,M.\,Kopytov}

\emp

\bmp \textbf{8.25.} Does there exist an algorithm which recognizes by
an identity whether it defines a
variety of groups that has at most countably many subvarieties?
\hfill \raisebox{-1ex}{\sl A.\,V.\,Kuznetsov} \emp

\bmp \textbf{8.27.} Does the lattice of all varieties of groups
possess non-trivial automorphisms?

\hfill
\raisebox{-1ex}{\sl A.\,V.\,Kuznetsov}
\emp

\bmp \textbf{8.29.} Do there exist locally nilpotent groups with
trivial centre satisfying the weak maximal
condition for normal
subgroups? \hfill \raisebox{-1ex}{\sl L.\,A.\,Kurdachenko} \emp

\bmp \textbf{8.30.} Let ${\goth X}$, ${\goth Y}$ be Fitting classes
of soluble groups which satisfy the Lockett condition, i.~e.
${\goth X} \cap {\goth S}_{\displaystyle{*}}={\goth
X}_{\displaystyle{*}}$, \ ${\goth Y} \cap {\goth S}_
{\displaystyle{*}} = {\goth Y}_{\displaystyle{*}}$ where ${\goth
S}$ denotes the Fitting class of all soluble groups and the lower
star the bottom group of the Lockett section determined by the
given Fitting class. Does ${\goth X} \cap {\goth Y}$ satisfy the
Lockett condition?

\makebox[15pt][r]{}{\it Editors' comment (2001):\/} The answer is
affirmative if ${\goth X}$ and ${\goth Y}$ are local
(A.\,Gry\-tczuk, N.\,T.\,Vorob'ev, {\it Tsukuba J.~Math.}, \textbf{18},
no.\,1 (1994), 63--67).
 \hfill \raisebox{-1ex}{\sl H.\,Lausch}

\emp

\bmp \textbf{8.33.} Let $a, b$ be two elements of a group, $a$ having
infinite order. Find a necessary and sufficient condition for
$\bigcap_{n=1}^{\infty} \left< a^n,\, b\right> =\left< b\right> .$
\hfill \raisebox{-1ex}{\sl F.\,N.\,Liman} \emp

\bmp \textbf{8.35.} c) Determine the conjugacy classes of maximal
subgroups in the sporadic simple group $F_1$. See the current
status in {\it Atlas of Finite Group Representations}
(\url{http://brauer.maths.qmul.ac.uk/Atlas/spor/M/}). \hfill
\raisebox{-1ex}{\sl V.\,D.\,Mazurov}

\emp

\bmp \textbf{8.40.} Describe the finite groups generated by a conjugacy
class $D$ of involutions which satisfies the following property: if
$a,b \in D$ and $|ab|=4$, then $[a,b] \in D$. This condition is
satisfied, for example, in the case where $D$ is a conjugacy class of
involutions of a known finite simple group such that $\left<
C_D(a)\right>$ is a 2-group for every $a \in D$. \hfill
\raisebox{-1ex}{\sl A.\,A.\,Makhn\"ev} \emp

\bmp \textbf{8.41.} What finite groups $G$ contain a normal set of
involutions $D$ which contains a non-empty proper subset $T$
satisfying the following properties:

 \makebox[15pt][r]{}1) $C_D(a) \subseteq T$ for any $a \in T$;

 \makebox[15pt][r]{}2) if $a,b \in T$ and $ab=ba \ne 1$, then
$C_D(ab) \subseteq
T$? \hfill {\sl A.\,A.\,Makhn\"ev}
\emp

\bmp \textbf{8.42.} Describe the finite groups all of whose soluble
subgroups of odd indices have 2-length 1. For example, the groups
$L_n(2^m)$ are known to satisfy this condition.

\hfill \raisebox{-1ex}{\sl A.\,A.\,Makhn\"ev}
\emp

\bmp \textbf{8.43.} (F.\,Timmesfeld). Let $T$ be a Sylow 2-sub\-group
of a finite group $G$. Suppose that $\left< N(B)\mid B \right.$ is
a non-trivial characteristic subgroup of $\left. T\right> $ is a
proper subgroup of~$G$. Describe the group $G$ if $F^
{\displaystyle{*}}(M)=O_2(M)$ for any 2-local subgroup $M$
containing~$T$.

\hfill \raisebox{-1ex}{\sl A.\,A.\,Makhn\"ev} \emp

\bmp \textbf{8.44.} Prove or disprove that for all, but finitely
many, primes $p$ the group \vspace{-1ex} $$G_p=\left< a, b\mid
a^2=b^p=(ab)^3=(b^rab^{-2r}a)^2=1\right> ,$$ \vskip-1ex

where $r^2 + 1 \equiv 0 \, ({\rm mod}\, p)$, is infinite. A
solution of this problem would have interesting topological
applications. It was proved with the aid of computer that
$G_p$ is finite for $p \leq 17$. \hfill
\raisebox{-1ex}{\sl J.\,Mennicke}

\emp

\bmp \textbf{8.45.} When it follows that a group is residually a
$({\goth X} \cap {\goth Y} )$-group, if it is both residually a
${\goth X} $-group and residually a ${\goth Y} $-group? \hfill
\raisebox{-1ex}{\sl Yu.\,I.\,Merzlyakov} \emp

\bmp \textbf{8.50.} At a conference in Oberwolfach in 1979 I
exhibited a finitely generated soluble group $G$ (of derived
length 3) and a non-zero cyclic ${\Bbb Z}G $-module $V$ such that
$V \cong V \oplus V$. Can this be done with $G$ metabelian? I
conjecture that it cannot. For background of this problem and for
details of the construction see (P.\,M.\,Neumann, {\it in:
Groups--Korea, Pusan, 1988 $($Lect. Notes Math.}, \textbf{1398}),
Springer, Berlin, 1989, 124--139).

\hfill \raisebox{-1ex}{\sl P.\,M.\,Neumann}

\emp

\bmp \textbf{\zv 8.51.}
 (J.\,McKay). If $G$ is a finite group and $p$ a
prime let $ m_p(G)$ denote the number of ordinary irreducible
characters of $G$ whose degree is not divisible by $p$. Let $P$ be
a Sylow $p\hskip0.2ex $-sub\-group of $G$. Is it true that
$m_p(G)=m_p(N_G(P))$? \hfill \raisebox{-1ex}{\sl J.\,B.\,Olsson}

\ul

\otv Yes, it is true (M.\,Cabanes, B.\,Sp\"ath, {\it Preprint}, 2024, \url{https://arxiv.org/abs/2410.20392})
\emp

\bmp \textbf{\ai 8.52.} 
(Well-known problem). Do there exist infinite
finitely presented pe\-ri\-o\-dic groups? Compare with 6.3. \hfill
\raisebox{-1ex}{\sl A.\,Yu.\,Olshanskii}

\ul

\oai Yes, such groups do exist (unverified AI solution: OpenAI, {\it Preprint}, 2026, \url{https://kourovkanotebookorg.wordpress.com/wp-content/uploads/2026/10/8-52.pdf}).
\emp 

\bmp \textbf{8.53.} a) Let $n$ be a sufficiently large odd number.
Describe the auto\-mor\-phisms of the free Burnside group $B(m,n)$
of exponent $n$ on $m$ generators.
\hfill \raisebox{-1ex}{\sl A.\,Yu.\,Olshanskii}
\emp

\bmp \textbf{8.54.} 
a) (Well-known problem). Classify metabelian
varieties of groups (or show that this is, in a certain sense, a
``wild'' problem).

 \makebox[25pt][r]{\ai b)} Describe the identities of 2-generated
metabelian groups,
that is, classify varieties generated by such groups. \hfill
\raisebox{-1ex}{\sl A.\,Yu.\,Olshanskii}

\ul

\oai b) These varieties are classified (unverified AI solution:  Achyuth Jayadevan,  {\it Preprint}, 2026, \url{https://zenodo.org/records/22843496}).

\emp 

\bmp \textbf{8.55.} It is easy to see that the set of all
quasivarieties of groups, in each of which quasivarieties a
non-trivial identity holds, is a semigroup under multiplication of
quasivarieties. Is this semigroup free? \hfill \raisebox{-1ex}{\sl
A.\,Yu.\,Olshanskii} \emp

\bmp
\textbf{8.58.} Does a locally compact locally nilpotent nonabelian group without elements of finite order contain a proper closed isolated normal subgroup?\hfill \raisebox{-1ex}{\sl
V.\,M.\,Poletskikh}
 \emp

\bmp \textbf{8.59.} Suppose that, in a locally compact locally
nilpotent group $G$, all proper closed normal subgroups
 are compact. Does
$G$ contain an open compact normal subgroup?

\hfill \raisebox{-1ex}{\sl V.\,M.\,Poletskikh,
I.\,V.\,Protasov}
\emp

\bmp \textbf{8.60.} Describe the locally compact primary locally
soluble groups which are covered by compact subgroups and all of
whose closed abelian subgroups have finite rank.
\makebox[1cm]{}\hfill \raisebox{-1ex}{\sl V.\,M.\,Poletskikh,
V.\,S.\,Charin} \emp

\bmp \textbf{8.62.} Describe the locally compact locally pronilpotent
groups for which the space of closed normal subgroups is compact
in the $ E $-topo\-logy. The corresponding problem has been solved
for discrete groups. \hfill \raisebox{-1ex}{\sl I.\,V.\,Protasov}
\emp

\bmp \textbf{8.64.} Does the class of all finite groups possess an
independent basis of quasiidentities? \hfill \raisebox{-1ex}{\sl
A.\,K.\,Rumyantsev, D.\,M.\,Smirnov} \emp

\bmp \textbf{8.67.} Do there exist Golod groups all of whose abelian
subgroups have finite ranks? Here a {\it Golod group\/} means a
finitely generated non-nilpotent subgroup of the adjoint group
of a nil-ring. The answer is negative for the whole adjoint
group of a nil-ring (Ya.\,P.\,Sysak, {\it Abstracts of the 17th
All--Union Algebraic Conf., part 1}, Minsk, 1983, p.\,180
(Russian)). \hfill
\raisebox{-1ex}{\sl A.\,I.\,Sozutov} \emp

\bmp \textbf{8.72.} Does there exist a finitely presented group $G$
which is not free or cyclic of prime order, having the property
that every proper subgroup of $G$ is free?
\hfill \raisebox{-1ex}{\sl C.\,Y.\,Tang}

\emp

\bmp \textbf{\zv 8.74.}
A subnormal subgroup $H\mbox{$\lhd\lhd$}\, G$
of a group $G$ is said
to be {\it good\/} if and only if $\left< H,J\right>
\mbox{$\lhd\lhd$}\, G$ whenever $J \mbox{$\lhd\lhd$}\,G$. Is it
true that when $H$ and $K$ are good subnormal subgroups of $G$,
then $H \cap K$ is good? \hfill \raisebox{-1ex}{\sl
J.\,G.\,Thompson}

\ul

\otv No, not always (K.\,Muliarchyk, {\it Preprint of 14 September 2026}, \url{https://kourovkanotebookorg.wordpress.com/wp-content/uploads/2026/09/kourovka_8_74_solution.pdf}; \ unverified AI solution:  Achyuth Jayadevan,  {\it Preprint of 18 September 2026}, \url{https://zenodo.org/records/22822596}).
\emp

 \bmp \textbf{\zv 8.75.}
 (A known problem). Suppose $G$ is a finite
primitive permutation group on $\Omega$, and $\alpha , \beta $ are
distinct points of $\Omega$. Does there exist an element $g \in G
$ such that $\alpha g=\beta $ and $g$ fixes no point of $\Omega$?
\hfill \raisebox{-1ex}{\sl J.\,G.\,Thompson}

\ul

\otv Not always (P.\,M\"uller, {\it Preprint}, 2023, \url{https://arxiv.org/pdf/2304.08459.pdf}).
\emp

\bmp \textbf{8.77.} Do there exist strongly regular graphs with
parameters $\lambda =0$, \ $\mu =2$ of degree $k > 10$? Such graphs
are known for $k=5$ and $ k=10$, their automorphism groups are
primitive permutation groups of rank 3. \hfill \raisebox{-1ex}{\sl
D.\,G.\,Fon-Der-Flaass} \emp

\bmp \textbf{8.78.} It is known that there exists a countable locally
finite group that contains a copy of every other countable locally
finite group. For which other classes of countable groups does a
similar ``largest'' group exist? In particular, what about periodic
locally soluble groups? periodic locally nilpotent groups? \hfill
\raisebox{-1ex}{\sl B.\,Hartley}

\emp

\bmp \textbf{8.79.} Does there exist a countable infinite locally
finite group $G$ that is complete, in the sense that $G$ has
trivial centre and no outer automorphisms? An uncountable one
exists (K.\,K.\,Hickin, {\it Trans. Amer. Math. Soc.}, \textbf{239}
(1978), 213--227).

\hfill
\raisebox{-1ex}{\sl B.\,Hartley}

\emp

\bmp \textbf{8.82.} Let ${\goth H}={\Bbb C} \times {\Bbb R} =\{
(z,r)\mid z \in {\Bbb C},\; r > 0\}$ be the three-dimensional
Poincar\'e space which admits the following action of the group
$SL_2({\Bbb C})$:
$$(z,r)
\left(\!\!\begin{array}{cc}a&b\cr
c&d\end{array}\!\!\right)=\left(
\frac{(az+b)(\overline{cz+d})+a\,\overline{\! c}r^2}{|cz+d|^2+|c|^2
r^2}\, ,\; \frac{r}{|cz+d|^2+|c|^2r^2}\right) .$$ Let ${\goth o}$ be
the ring of integers of the field $K={\Bbb Q}(\sqrt{D})$, where $D <
0$, \ $\pi$ a prime such that $\pi \overline{\pi}$ is a prime in
${\Bbb Z}$ and let $\displaystyle{\Gamma=\left\{
\left.\left(\!\!\begin{array}{cc}a&b\cr c&d\end{array}\!\!\right) \in
SL_2({\goth o})\;\right|\; b \equiv 0 \, ({\rm mod}\, \pi )\right\}
}$. Adding to the space $\Gamma \setminus {\goth H}^3$ two vertices
we get a three-dimensional compact space $\overline{\Gamma \setminus
{\goth H}^3}$.

 \makebox[15pt][r]{}Calculate $r(\pi )={\rm dim}_{\,{\Bbb
Q}}\, H_1( \overline{\Gamma \setminus {\goth H}^3},{\Bbb Q})={\rm
dim}_{\,{\Bbb Q}} (\Gamma ^{ab} \otimes {\Bbb Q}).$

\makebox[15pt][r]{}For example, if $D=-3$ then $r(\pi )$ is distinct
from zero for
the first time for $\pi \mid 73$ (then $r(\pi )=1$), and if
$D=- 4$ then $r(\pi )=1$ for $\pi \mid 137$. \hfill
\raisebox{-1ex}{\sl H.\,Helling}

\emp

\bmp \textbf{8.83.} In the notation of 8.82, for $r(\pi ) > 0$, one
can define via Hecke algebras a formal Dirichlet series with Euler
multiplication (see G.\,Shimura, {\it Introduction to the
arithmetic theory of automorphic functions}, Princeton Univ.
Press, 1971). Does there exist an algebraic Hasse--Weil variety
whose $ \zeta $-function is this Dirichlet series? There are
several conjectures. \hfill \raisebox{-1ex}{\sl H.\,Helling}

\emp

\bmp \textbf{8.85.} Construct a finite $p\hskip0.2ex $-group $G$
whose Hughes subgroup $H_p(G)=\left< x \in G\mid \right.$ $\left.
|x|
\ne
p\right> $ is non-trivial and has index $p^3$. \hfill
\raisebox{-1ex}{\sl E.\,I.\,Khukhro}
\emp

\bmp \textbf{8.86.} (Well-known problem). Suppose that all proper
closed subgroups of a locally compact locally nilpotent group $G$
are compact. Is $G$ abelian if it is non-compact?

\hfill \raisebox{-1ex}{\sl V.\,S.\,Charin}
\emp

\raggedbottom

\newpage

\pagestyle{myheadings} \markboth{9th Issue
(1984)}{9th Issue (1984)}
\thispagestyle{headings}

~
\vspace{2ex}

\centerline {\Large \textbf{Problems from the 9th Issue (1984)}}
\phantomsection\label{9izd}
\vspace{4ex}

\bmp \textbf{9.1.} A group $G$ is said to be {\it potent\/} if, to each
 $x\in G$ and each positive integer $n$ dividing the order of
$x$ (we suppose $ \infty$ is divisible by every positive integer),
there exists a finite homomorphic image of $G$ in which the image of
$x$ has order precisely~$ n$. Is the free product of two potent
groups again potent? \hfill \raisebox{-1ex}{\sl
R.\,B.\,J.\,T.\,Allenby}

\emp

\bmp \textbf{\zv 9.4.}
It is known that if ${\goth M}$ is a variety
(quasivariety, pseudovariety) of groups, then the class $I{\goth M}$
of all quasigroups that are isotopic to groups in ${\goth M}$ is also
a variety (quasivariety, pseudovariety), and $I{\goth M}$ is finitely
based if and only if ${\goth M}$ is finitely based. Is it true that
if ${\goth M}$ is generated by a single finite group then $I{\goth
M}$ is generated by a single finite quasigroup? \hfill
\raisebox{-1ex}{\sl A.\,A.\,Gvaramiya}

\ul

\otv  No, it is not true (AI solutions verified by K.\,Muliarchyk: \  M.\,Kinyon, J.\,Urban, {\it Preprint of 17 September 2026},
\url{https://kourovkanotebookorg.wordpress.com/wp-content/uploads/2026/09/kourovka-mathematics.pdf}; \  Achyuth Jayadevan,  {\it Preprint of 18 September 2026}, \url{https://zenodo.org/records/22820770}). For pseudovarieties see (K.\,Muliarchyk, {\it Preprint}, 2026,
\url{https://kourovkanotebookorg.wordpress.com/wp-content/uploads/2026/09/kourovka_9_4_solution.pdf}).
 \emp

\bmp \textbf{9.5.} A variety of groups is called {\it primitive\/} if
each of its subquasivarieties is a variety. Describe all primitive
varieties of groups. Is every primitive variety of groups locally
finite? \hfill \raisebox{-1ex}{\sl V.\,A.\,Gorbunov} \emp

\bmp \textbf{9.6.} Is it true that an independent basis of
quasiidentities of any finite group is finite?

\hfill
\raisebox{-1ex}{\sl V.\,A.\,Gorbunov} \emp

\bmp \textbf{9.7.} (A.\,M.\,St\"epin). Does there exist an infinite
finitely generated amenable group of bounded exponent? \hfill
\raisebox{-1ex}{\sl R.\,I.\,Grigorchuk}

\emp

\bmp \textbf{9.9.} Does there exist a finitely generated group which
is not nilpotent-by-finite and whose growth function has as a
majorant a function of the form $c^{\sqrt{n}}$ where $ c$ is a
constant greater than 1? \hfill \raisebox{-1ex}{\sl
R.\,I.\,Grigorchuk}
\emp

\bmp \textbf{9.11.} Is an abelian minimal normal subgroup $A$ of a
group $G$ an elementary $p\hskip0.2ex$-group if the factor-group $
G/A$ is a soluble group of finite rank? The answer is affirmative
under the additional hypothesis that $ G/A$ is locally polycyclic
(D.\,I.\,Zaitsev, {\it Algebra and Logic}, \textbf{19}, no.\,2
(1980), 94--105). \hfill \raisebox{-1ex}{\sl D.\,I.\,Zaitsev}
\emp

\bmp \textbf{9.13.} Is a soluble torsion-free group minimax if it
satisfies the weak minimum condition for normal subgroups? \hfill
\raisebox{-1ex}{\sl D.\,I.\,Zaitsev}
\emp

\bmp \textbf{9.14.} Is a group locally finite if it decomposes into a
product of periodic abelian subgroups which commute pairwise? The
answer is affirmative in the case of two subgroups. \hfill
\raisebox{-1ex}{\sl D.\,I.\,Zaitsev}

\emp

\bmp \textbf{9.15.} Describe, without using CFSG, all subgroups $ L$
of a finite Chevalley group $G$ such that $G=PL$ for some
parabolic subgroup $P$ of $G$. \hfill \raisebox{-1ex}{\sl
A.\,S.\,Kondratiev}

\emp

\bmp \textbf{9.17.} b) Let $G$ be a locally normal residually finite
group. Does there exist a normal subgroup $H$ of $G$ which is
embeddable in a direct product of finite groups and is such that
$G/H$ is a divisible abelian group? \hfill \raisebox{-1ex}{\sl
L.\,A.\,Kurdachenko}

\emp

\bmp \textbf{9.19.} a) Let $n(X)$ denote the minimum of the indices of
proper subgroups of a group~$X$. A subgroup $A$ of a finite group $G$
is called {\it wide\/} if $A$ is a maximal element by inclusion of
the set $\{ X\mid X \mbox{ is a proper subgroup of }G\mbox{ and
}n(X)=n(G)\} $. Find all wide subgroups in finite projective special
linear, symplectic, orthogonal, and unitary groups. \hfill
\raisebox{-1ex}{\sl V.\,D.\,Mazurov}

\emp

\bmp \textbf{9.23.} Let $G$ be a finite group, $B$ a block of
characters of $ G$, $D(B)$ its defect group and $k(B)$
(respectively, $k_0(B)$) the number of all irreducible complex
characters (of height 0) lying in $B$. {\it Conjectures:}

\makebox[25pt][r]{a)} (R.\,Brauer) \ $k(B) \leq |D(B)|$; this has
been proved for $p\hs$-soluble groups (D.\,Gluck, K.\,Magaard,
U.\,Riese, P.\,Schmid, {\it J.~Algebra}, \textbf{279} (2004), 694--719);

\makebox[25pt][r]{b)} (J.\,B.\,Olsson) \ $k_0(B) \leq |D(B) :
D(B)'|$, where $D(B)'$ is the derived subgroup of $D(B)$;

\makebox[25pt][r]{\zv c)} (R.\,Brauer) \ $D(B)$
is abelian if
and only
if $k_0(B)=k(B)$.
\hfill \raisebox{-1ex}{\sl V.\,D.\,Mazurov}

\ul

\otv c) The conjecture is proved (G.\,Malle, G.\,Navarro, A.\,A.\,Schaeffer Fry,  P.\,H.\,Tiep, {\it Ann.  Math.}, {\bf  200} (2024), 557--608).
\emp

\bmp \textbf{9.24.} (J.\,G.\,Thompson). {\it Conjecture:\/} every
finite simple non-abelian group $G$ can be represented in the form
$ G=CC$, where $C$ is some conjugacy class of $G$.

\hfill
\raisebox{-1ex}{\sl V.\,D.\,Mazurov}
\emp

\bmp \textbf{9.26.} b) Describe the finite groups of 2-local 3-rank 1
which have non-cyclic Sylow 3-sub\-groups. \hfill
\raisebox{-1ex}{\sl A.\,A.\,Makhn\"ev}
\emp

\bmp \textbf{9.28.} Suppose that a finite group $G$ is generated by a
conjugacy class $D$ of involutions and let $D_i=\{ d_1\cdots
d_i\mid d_1,\ldots ,d_i$ are different pairwise commuting elements
of $ D\} $. What is $G$, if $ D_1,\ldots ,D_n$ are all its
different conjugacy classes of involutions? For example, the
Fischer groups $F_{22}$ and $F_{23}$ satisfy this condition with
$n=3$.

\hfill
\raisebox{-1ex}{\sl A.\,A.\,Makhn\"ev}

\emp

\bmp \textbf{9.29.} (Well-known problem). According to a classical
theorem of Magnus, the word problem is soluble in 1-relator
groups. Do there exist 2-relator groups with insoluble word
problem? \hfill \raisebox{-1ex}{\sl Yu.\,I.\,Merzlyakov}

\emp

\bmp \textbf{9.31.} Let $k$ be a field of characteristic 0. According
to (Yu.\,I.\,Merzlyakov, {\it Proc. Stek\-lov Inst. Math.}, {\bf
167} (1986), 263--266), the family ${\rm Rep}_k(G)$ of all
canonical matrix representations of a $k$-powered group $G$ over $
k$ may be regarded as an affine $ k$-variety. Find an explicit
form of equations defining this variety. \hfill
\raisebox{-1ex}{\sl Yu.\,I.\,Merzlyakov}

\emp

\bmp \textbf{9.35.} A topological group is said to be {\it
inductively compact\/} if any finite set of its elements is
contained in a compact subgroup. Is this property preserved under
lattice isomorphisms in the class of locally compact groups? \hfill
\raisebox{-1ex}{\sl Yu.\,N.\,Mukhin}

\emp

\bmp \textbf{9.36.} Characterize the lattices of closed subgroups in
locally compact groups. In the discrete case this was done in
(B.\,V.\,Yakovlev, {\it Algebra and Logic}, \textbf{13}, no.\,6
(1974), 400--412). \hfill \raisebox{-1ex}{\sl Yu.\,N.\,Mukhin}

\emp

\bmp \textbf{9.37.} Is a compactly generated inductively prosoluble
locally compact group pro\-soluble? \hfill \raisebox{-1ex}{\sl
Yu.\,N.\,Mukhin}

\emp

\bmp \textbf{9.38.} A group is said to be {\it compactly covered\/} if
it is the union of its compact sub\-groups. In a null-dimensional
locally compact group, are maximal compactly covered subgroups
closed? \hfill \raisebox{-1ex}{\sl Yu.\,N.\,Mukhin}

\emp

\bmp \textbf{9.39.} Let $\Omega$ be a countable set and ${\goth m}$ a
cardinal number such that $\aleph _0 \leq {\goth m} \leq 2^{\aleph
_0}$ (we assume Axiom of Choice but not Continuum Hypothesis). Does
there exist a permutation group $G$ on $\Omega$ that has exactly
${\goth m}$ orbits on the power set ${\goth P}(\Omega )$?

\makebox[15pt][r]{}{\it Comment of 2001:\/} It is proved
(S.\,Shelah, S.\,Thomas, {\it Bull. London Math. Soc.}, \textbf{20},
no.\,4 (1988), 313--318) that the answer is positive in set theory
with Martin's Axiom. The question is still open in ZFC. \hfill
\raisebox{-1ex}{\sl P.\,M.\,Neumann}

\emp

\bmp \textbf{9.40.} Let $\Omega$ be a countably infinite set. Define
a {\it moiety\/} of $\Omega$ to be a subset $\Sigma$ such that
both $\Sigma$ and $\Omega \setminus \Sigma$ are infinite. Which
permutation groups on $ \Omega$ are transitive on moieties? \hfill
\raisebox{-1ex}{\sl P.\,M.\,Neumann}

\emp

\bmp \textbf{9.41.}
Let $\Omega$ be a countably infinite set. For $k
\geq 2$ define a {\it $k$-section\/} of $\Omega$ to be a partition
of $\Omega$ as union of $k$ infinite sets.

\makebox[25pt][r]{b)} Does there exist a group that is
transitive on
$k$-sections
but not on $(k + 1)$-sections?

\makebox[25pt][r]{\zv c)}
Does there exist a transitive permutation
group on $\Omega$ that is transitive on $\aleph _0$-sections but
which is a proper subgroup of ${\rm Sym}(\Omega )$? \hfill
\raisebox{-1ex}{\sl P.\,M.\,Neumann}

\ul

\otv c) Yes, it does exist  (AI solution confirmed by S.\,M.\,Corson:  \ Achyuth Jayadevan,  {\it Preprint}, 2026, \url{https://zenodo.org/records/22820556}).  It was also shown that an affirmative answer is consistent with the ZFC axioms of set theory (S.\,M.\,Corson, S.\,Shelah, {\it Preprint}, 2025, \url{https://arxiv.org/abs/2503.12997}).
\emp 

\bmp \textbf{\zv 9.42.}
Let $\Omega$ be a countable set and let $D$ be
the set of total order relations on $ \Omega$ for which $ \Omega$
is order-isomorphic with ${\Bbb Q}$. Does there exist a transitive
proper subgroup $G$ of ${\rm Sym} (\Omega )$ which is transitive
on $D$? \hfill \raisebox{-1ex}{\sl P.\,M.\,Neumann}

\ul

\otv An affirmative answer is consistent with the ZFC axioms of set theory (S.\,M.\,Corson, S.\,Shelah, {\it Preprint}, 2025, \url{https://arxiv.org/abs/2503.12997}).
\emp

\bmp \textbf{9.43.} b) The group $G$ indicated in (N.\,D.\,Podufalov,
{\it Abstracts of the 9th All--Union Symp. on Group Theory},
Moscow, 1984, 113--114 (Russian)) allows us to construct a
projective plane of order 3: one can take as lines any class of
subgroups of order $2\cdot 5\cdot 11\cdot 17$ conjugate under $ S$
and add four more lines in a natural way. In a similar way, one
can construct projective planes of order $p^n$ for any prime $p$
and any natural~$ n$. Could this method be adapted for
constructing new planes?
\hfill \raisebox{-1ex}{\sl N.\,D.\,Podufalov}
\emp

\bmp \textbf{9.44.} A topological group is said to be {\it layer
compact\/} if the full inverse images of all of its compacts under
mappings $x \rightarrow x^n$, $n=1,\, 2,\, \ldots \,$, are
compacts. Describe the locally compact locally soluble layer
compact groups. \hfill \raisebox{-1ex}{\sl V.\,M.\,Poletskikh}

\emp

\bmp \textbf{\ai 9.45.} 
 Let $a$ be a vector in ${\Bbb R}^n$ with rational
coordinates and set \vspace{-1ex} $$S(a)=\{ ka + b\mid k \in {\Bbb
Z},\,\,\, b \in {\Bbb Z}^n\} .$$ \vskip-1ex

 It is obvious that $S(a)$ is a discrete subgroup in
${\Bbb R}^n$ of rank $n$. Find necessary and sufficient conditions
in terms of the coordinates of $a$ for $S(a)$ to have an
orthogonal basis with respect to the standard scalar product in
${\Bbb R}^n$.

\hfill \raisebox{-1ex}{\sl Yu.\,D.\,Popov,
I.\,V.\,Protasov}

\ul

\oai Conditions found (unverified AI solution:  M.\,Kinyon, J.\,Urban, {\it Preprint}, 2026,

\url{https://kourovkanotebookorg.wordpress.com/wp-content/uploads/2026/09/kourovka-mathematics.pdf}).
\emp 

\bmp \textbf{\zv 9.47.}
Is it true that every scattered compact can be
homeomorphically embedded into the space of all closed non-compact
subgroups with $ E$-topo\-logy of a suitable locally compact
group? \hfill \raisebox{-1ex}{\sl I.\,V.\,Protasov}

\ul

\otv No, not every (AI solution confirmed by  I.\,V.\,Protasov: \ Achyuth Jayadevan,  {\it Preprint of 18 September 2026}, \url{https://zenodo.org/records/22838025}; \  unverified AI solution:  M.\,Kinyon, J.\,Urban, {\it Preprint of 17 September 2026},
\url{https://kourovkanotebookorg.wordpress.com/wp-content/uploads/2026/09/kourovka-mathematics.pdf}).
\emp 

\bmp \textbf{9.51.} Does there exist a finitely presented soluble
group satisfying the maximum con-\linebreak dition on normal subgroups which
has insoluble word problem? \hfill \raisebox{-1ex}{\sl
D.\,J.\,S.\,Robinson}

\emp

\bmp \textbf{9.52.} Does a finitely presented soluble group of finite
rank have soluble conjugacy problem? Note: there is an algorithm
to decide conjugacy to a given element of the group. \hfill
\raisebox{-1ex}{\sl D.\,J.\,S.\,Robinson}

\emp

\bmp \textbf{9.53.} Is the isomorphism problem soluble for finitely
presented soluble groups of finite rank? \hfill
\raisebox{-1ex}{\sl D.\,J.\,S.\,Robinson}

\emp

\bmp \textbf{9.55.} Does there exist a finite $p\hskip0.2ex$-group
$G$ and a central augmented automorphism $\varphi$ of ${\Bbb Z}G$
such that $\varphi $, if extended to $\widehat{{\Bbb Z}}_pG$, the
$p\hskip0.2ex$-adic group ring, is {\it not\/} conjugation by unit
in $\widehat{{\Bbb Z}}_pG$ followed by a homomorphism induced from
a group automorphism?

\hfill\raisebox{-1ex}{\sl K.\,W.\,Roggenkamp}
\emp

\bmp \textbf{9.57.} The set of subformations of a given formation is
a lattice with respect to operations of intersection and
generation. What formations of finite groups have distributive
lattices of subformations? \hfill \raisebox{-1ex}{\sl
A.\,N.\,Skiba}
\emp

\bmp \textbf{9.61.} Two varieties are said to be {\it
$S$-equiva\-lent\/} if they have the same Mal'cev theory
(D.\,M.\,Smirnov, {\it Algebra and Logic}, \textbf{22}, no.\,6
(1983), 492--501). What is the cardinality of the set of
$S$-equiva\-lent varieties of groups? \hfill \raisebox{-1ex}{\sl
D.\,M.\,Smirnov}
\emp

\bmp \textbf{9.65.} Is a locally soluble group periodic if it is a
product of two periodic subgroups?

\hfill \raisebox{-1ex}{\sl Ya.\,P.\,Sysak}
\emp

\bmp \textbf{\zv 9.66.}
a) {\it B.\,Jonsson's Conjecture:\/} elementary
equivalence is preserved under taking free products in the class
of all groups, that is, if $Th(G_1)=Th(G_2)$ and $Th(H_1)=Th(H_2)$
for groups
 $G_1$, $G_2$, $H_1$, $H_2$, then $Th(G_1*H_1)=Th(G_2*H_2)$.

\makebox[15pt][r]{}b) It may be interesting to consider also the
following {\it weakened conjecture:\/} if $Th(G_1)=Th(G_2)$ and
$Th(H_1)=Th(H_2)$ for countable groups
 $G_1$, $G_2$, $H_1$, $H_2$, then
 for any numerations
of the groups $G_i$, $
H_i$ there are $m$-redu\-ci\-bi\-lity
$T_1\mathop{\equiv}\limits_{m}T_2$
and Turing
reducibility
$T_1\mathop{\equiv}\limits_{T}T_2$,
where $T_i$ is the set of numbers of
all theorems in $Th(G_i*H_i)$. A~proof of this weakened conjecture
would be a good illustration of application of the reducibility theory
to solving concrete mathematical problems.
\hfill \raisebox{-1ex}{\sl A.\,D.\,Taimanov}

\ul

\otv The conjecture (also known as R.\,L.\,Vaught's conjecture) is proved (Z.\,Sela, {\it Preprint}, 2010, \url{https://arxiv.org/pdf/1012.0044}).
\emp

\bmp \textbf{9.68.} Let ${\goth V}$ be a variety of groups which is
not the variety of all groups, and let $p$ be a prime. Is there a
bound on the
$
p\hskip0.2ex$-lengths of the finite $p\hskip0.2ex$-soluble groups
whose Sylow $ p\hskip0.2ex$-sub\-groups are in ${\goth V}$? \hfill
\raisebox{-1ex}{\sl John S.\,Wilson}

\emp

\bmp \textbf{\zv 9.69.}
(P.\,Cameron). Let $G$ be a finite primitive
permutation group and suppose that the stabilizer $G_{\alpha}$ of
a point $\alpha$ induces on some of its orbits $\Delta \ne \{
\alpha \}$ a regular permutation group. Is it true that
$|G_{\alpha}|=|\Delta |$? \hfill \raisebox{-1ex}{\sl
A.\,N.\,Fomin}

\ul

\otv No, not always (P.\,Spiga, {\it J.~Group Theory}, {\bf 25}, no.\,1 (2022), 113--126).

\emp

\bmp \textbf{9.70.} Prove or disprove the conjecture of
 P.\,Cameron ({\it Bull.
London Math. Soc.}, \textbf{13}, no.\,1 (1981), 1--22): if $G$ is a
finite primitive permutation group of subrank $m$, then either the
rank of $G$ is bounded by a function of $m$ or the order of a
point stabilizer is at most $m$. The {\it subrank\/} of a
transitive permutation group is defined to be the maximum rank of
the transitive constituents of a point stabilizer. \hfill
\raisebox{-1ex}{\sl A.\,N.\,Fomin}

\emp

\bmp \textbf{9.71.} Is it true that every infinite 2-transitive
permutation groups with locally soluble point stabilizer has a
non-trivial irreducible finite-dimensional representation over some
field? \hfill \raisebox{-1ex}{\sl A.\,N.\,Fomin}

\emp

\bmp \textbf{9.72.} Let $G_1$ and $G_2$ be Lie groups with the
following property: each $G_i$ contains a nilpotent simply
connected normal Lie subgroup $B_i$ such that $G_i/B_i \cong
SL_2(K)$, where $K={\Bbb R}$ or ${\Bbb C}$. Assume that $G_1$ and $G_2$
are contained as closed subgroups in a topological group $G$, that
$G_1 \cap G_2 \geq B_1B_2$, and that no non-identity Lie subgroup
of $B_1 \cap B_2$ is normal in $G$. Can it then be shown (perhaps
by using the method of ``amalgams'' from the theory of finite
groups) that the nilpotency class and the dimension of $B_1B_2$ is
bounded? \hfill \raisebox{-1ex}{\sl A.\,L.\,Chermak}

\emp

\bmp \textbf{9.75.}
Find all local formations ${\goth F}$ of finite
groups such that, in every finite group, the set of ${\goth
F}$-subnormal subgroups forms a lattice.

\makebox[15pt][r]{}\textit{Comment of 2015}: Subgroup-closed formations with this property are described in (Xiaolan Yi, S.\,F.\,Kamornikov, \textit{J.~Algebra}, {\bf 444} (2015),
143--151).\hfill
\raisebox{-1ex}{\sl L.\,A.\,Shemetkov}
\emp

\bmp \textbf{9.76.} We define a {\it Golod group\/} to be the
$r$-generated, $r \geq 2$, subgroup $\left< 1 + x_1+I,\right.$
$\left. 1+x_2+I,\ldots ,1+x_r+I\right>$ of the adjoint group
 $1 +
F/I$ of the factor-algebra $F/I$, where $F$ is a free algebra of
polynomials without constant terms in non-commuting variables $x_1,
x_2, \ldots , x_r$ over a field of characteristic $p>0$, and $I$ is
an ideal of $F$ such that $F/I$ is a non-nilpotent nil-algebra (see
E.\,S.\,Golod, {\it Amer. Math. Soc. Transl. (2)}, \textbf{48} (1965),
103--106). Prove that in Golod groups the centralizer of every
element is infinite.

\makebox[15pt][r]{}Note that Golod groups with infinite centre
were constructed by A.\,V.\,Timofeenko ({\it Math. Notes}, {\bf
39}, no.\,5 (1986), 353--355); independently by other methods the
same result was obtained in the 90s by L.\,Hammoudi. \hfill
\raisebox{-1ex}{\sl V.\,P.\,Shunkov}

\emp

\bmp \textbf{9.77.} Does there exist an infinite finitely generated
residually finite binary finite group all of whose Sylow subgroups
are finite? A.\,V.\,Rozhkov (Dr. of Sci. Disser., 1997) showed that
such a group does exist if the condition of finiteness of the
Sylow subgroups is weakened to local finiteness. \hfill
\raisebox{-1ex}{\sl V.\,P.\,Shunkov}

\emp

\bmp \textbf{9.78.} Does there exist a periodic residually finite $
F^{\displaystyle{*}}$-group (see Archive, 7.42) all of whose Sylow
subgroups are finite and which is not binary finite? \hfill
\raisebox{-1ex}{\sl V.\,P.\,Shunkov}

\emp

\bmp \textbf{9.83.} Suppose that $G$ is a (periodic) $
p\hskip0.2ex$-conjugacy biprimitively finite group (see 6.59)
which has a finite Sylow
$
p\hskip0.2ex$-sub\-group. Is it then true that all Sylow
$p\hskip0.2ex$-sub\-groups of $G$ are conjugate? \hfill
\raisebox{-1ex}{\sl V.\,P.\,Shunkov}

\emp

\bmp \textbf{9.84.} a) Is every binary finite 2-group of finite
exponent locally finite?

 \makebox[15pt][r]{}b) The same question for $p\hskip0.2ex$-groups for $p
> 2$.\hfill \raisebox{-1ex}{\sl V.\,P.\,Shunkov}

\emp

\raggedbottom

\newpage

\pagestyle{myheadings} \markboth{10th Issue
(1986)}{10th Issue (1986)}
\thispagestyle{headings}

~
\vspace{2ex}

\centerline {\Large \textbf{Problems from the 10th Issue (1986)}}
\phantomsection\label{10izd}
\vspace{3ex}

\bmp \textbf{10.2.} A {\it mixed identity\/} of a group $G$ is, by
definition, an identity of an algebraic system obtained from $G$
by supplementing its signature by some set of 0-ary operations.
One can develop a theory of mixed varieties of groups on the basis
of this notion (see, for example, V.\,S.\,Anashin, {\it Math. USSR
Sbornik}, \textbf{57} (1987), 171--182). Construct an example of a
class of groups which is not a mixed variety, but which is closed
under taking factor-groups, Cartesian products, and those
subgroups of Cartesian powers that contain the diagonal subgroup.
\hfill \raisebox{-1ex}{\sl V.\,S.\,Anashin} \emp

\bmp \textbf{10.3.} Characterize (in terms of bases of mixed
identities or in terms of generating groups) minimal mixed
varieties of groups. \hfill \raisebox{-1ex}{\sl V.\,S.\,Anashin}
\emp

\bmp \textbf{\zv 10.4.}
Let $p$ be a prime number. Is it true that a
mixed variety of groups generated by an arbitrary finite
$p\hskip0.2ex $-group of sufficiently large nilpotency class is a
variety of groups?

\hfill \raisebox{-1ex}{\sl V.\,S.\,Anashin}

\ul

\otv No, not always (AI solution confirmed by K.\,Muliarchyk: \  Achyuth Jayadevan,  {\it Preprint}, 2026, \url{https://zenodo.org/records/22847755}).
\emp 

\bmp \textbf{10.5.} Construct an example of a finite group whose
mixed identities do not have a finite basis. Is the group
constructed in (R.\,Bryant, {\it Bull. London Math. Soc.}, {\bf
14}, no.\,2 (1982), 119--123) such an example? \hfill
\raisebox{-1ex}{\sl V.\,S.\,Anashin} \emp

\bmp \textbf{10.8.} Does there exist a topological group which cannot
be embedded in the mul\-ti\-pli\-ca\-tive semigroup of a
topological ring? \hfill \raisebox{-1ex}{\sl V.\,I.\,Arnautov,
A.\,V.\,Mikhal\"ev} \emp

\bmp \textbf{10.10.} Is it true that a quasivariety generated by a
finitely generated torsion-free soluble group and containing a
non-abelian free metabelian group, can be defined in the class of
torsion-free groups by an independent system of quasiidentities?

\hfill \raisebox{-1ex}{\sl A.\,I.\,Budkin}
\emp

\bmp \textbf{10.11.} Is it true that every finitely presented group
contains either a free subsemigroup on two generators or a
nilpotent subgroup of finite index? \hfill \raisebox{-1ex}{\sl
R.\,I.\,Grigorchuk} \emp

\bmp \textbf{10.12.} Does there exist a finitely generated semigroup $S$
with cancellation having non-exponential growth and such that its
group of left quotients $G=S^{-1}S$ (which exists) is a group of
exponential growth? An affirmative answer would give a positive
solution to Problem 12 in (S.\,Wagon, {\it The Banach--Tarski
Paradox}, Cambridge Univ. Press, 1985). \hfill \raisebox{-1ex}{\sl
R.\,I.\,Grigorchuk} \emp

\bmp \textbf{\zv 10.13.}
Does there exist a {\it massive set\/} of
independent elements in a free group $F_2$ on free generators $a,
b$, that is, a set $ E$ of irreducible words on the alphabet $a,
b, a^{-1}, b^{-1}$ such that

\makebox[15pt][r]{}1) no element $w \in E$ belongs to the
normal closure of
$E\setminus \{ w\}$, and

\makebox[15pt][r]{}2) the massiveness condition is satisfied:
$\mathop{\overline{{\rm lim}}}\limits_{n\rightarrow \infty}
\sqrt[n]{|E_n|}=3,$
where $E_n$ is the set of all words of length $n$ in $E$? \hfill
\raisebox{-1ex}{\sl R.\,I.\,Grigorchuk}

\ul

\otv No, there is no such set (L.\,Koch-Hyde,  E.\,Olive,
 {\it Preprint}, 2026, \url{https://arxiv.org/pdf/2609.00382}).
\emp

\bmp \textbf{10.15.} For every (known) finite quasisimple group and
every prime $p$, find the faithful $ p\hskip0.2ex $-modular
absolutely irreducible linear representations of minimal degree.

\hfill
\raisebox{-1ex}{\sl A.\,S.\,Kondratiev}
\emp

\bmp \textbf{10.16.} A class of groups is called a {\it direct
variety\/} if it is closed under taking sub\-groups,
factor-groups, and direct products (Yu.\,M.\,Gorchakov, {\it
Groups with Finite Classes of Conjugate Elements}, Moscow, Nauka,
1978 (Russian)). It is obvious that the class of $ FC $-groups is
a direct variety. P.\,Hall ({\it J. London Math. Soc.}, \textbf{34},
no.\,3 (1959), 289--304) showed that the class of finite groups
and the class of abelian groups taken together do not generate the
class of
 $
FC $-groups as a direct variety, and it was shown in
(L.\,A.\,Kurdachenko, {\it Ukrain. Math.~J.}, \textbf{39}, no.\,3
(1987), 255--259) that the direct variety of $FC $-groups is also
not generated by the class of groups with finite derived
subgroups. Is the direct variety of
$
FC $-groups generated by the class of groups with finite derived
subgroups together with the class of $ FC $-groups having
quasicyclic derived subgroups?
\hfill \raisebox{-1ex}{\sl
L.\,A.\,Kurdachenko} \emp

\bmp \textbf{10.17.} (M.\,J.\,Tomkinson). Let $G$ be an $FC $-group
whose derived subgroup is embeddable in a direct product of finite
groups. Must $G/Z(G)$ be embeddable in a direct product of finite
groups? \hfill \raisebox{-1ex}{\sl L.\,A.\,Kurdachenko} \emp

\bmp \textbf{10.18.} (M.\,J.\,Tomkinson). Let $G$ be an $FC $-group
which is residually in the class of groups with finite derived
subgroups. Must $G/Z(G)$ be embeddable in a direct product of
finite groups? \hfill \raisebox{-1ex}{\sl L.\,A.\,Kurdachenko}
\emp

\bmp \textbf{10.19.} Characterize the radical associative rings such
that the set of all normal subgroups of the adjoint group
coincides with the set of all ideals of the associated Lie ring.
\hfill \raisebox{-1ex}{\sl V.\,M.\,Levchuk} \emp

\bmp \textbf{10.20.} In a Chevalley group of rank $\leq 6$ over a
finite field of order $\leq 9$, describe all subgroups of the form
$H=\left< H \cap U,\, H \cap V\right>$ that are not contained in
any proper parabolic subgroup, where $U$ and $V$ are opposite
unipotent subgroups.

\hfill \raisebox{-1ex}{\sl V.\,M.\,Levchuk}
\emp

\bmp \textbf{10.26.}
b) Does there exist an
algorithm which
decides whether the equation of the form
$w(x_{i_1}^{\varphi
_1},\ldots ,x_{i_
n}^{\varphi _n})=1$ is soluble in a free group
where $\varphi _1,\ldots ,\varphi _n$ are automorphisms of this
group? \hfill \raisebox{-1ex}{\sl G.\,S.\,Makanin}
\emp

\bmp \textbf{10.27.} a) Let $t$ be an involution of a finite group
$G$ and suppose that the set $D=t^G \cup \{ t^xt^y\mid x,y \in G,
\; \; |t^xt^y|=2\} $ does not intersect $O_2(G)$. Prove that if $t
\in O_2(C(d))$ for any involution $d$ from $C_D(t)$, then $ D=t^G$
(and in this case the structure of the group $\left< D\right> $ is
known).

 \makebox[15pt][r]{}b) A significantly more general question.
Let $ t$ be an involution of a finite group $G$ and suppose that
$t \in Z ^{\displaystyle{*}}(N(X))$ for every non-trivial subgroup
$X$ of odd order which is normalized, but not centralized, by $t$.
What is the group~$\left< t^G\right>$? \hfill \raisebox{-1ex}{\sl
A.\,A.\,Makhn\"ev} \emp

\bmp \textbf{10.29.} Describe the finite groups which contain a set
of involutions $D$ such that, for any subset $D_0$ of $ D$
generating a 2-sub\-group, the normalizer $N_D(D_0)$ also
generates a 2-sub\-group. \hfill \raisebox{-1ex}{\sl
A.\,A.\,Makhn\"ev} \emp

\bmp \textbf{10.31.} Suppose that $G$ is an algebraic group, $H$ is a
closed normal subgroup of $G$ and $f\! : G \rightarrow G/H$ is the
canonical homomorphism. What conditions ensure that there exists a
rational section $s\! : G/H \rightarrow G$ such that $ sf=1$? For
partial results, see (Yu.\,I.\,Merzlyakov, {\it Rational Groups},
Nauka, Moscow, 1987, $\S$\,37 (Russian)).

\hfill\raisebox{-1ex}{\sl Yu.\,I.\,Merzlyakov} \emp

\bmp \textbf{\ai 10.32.} 
A group word is said to be {\it universal on a
group\/} $G$ if its values on $G$ run over the whole of $G$. For
an arbitrary constant $c > 8/5$, J.\,L.\,Brenner, R.\,J.\,Evans,
and D.\,M.\,Silberger ({\it Proc. Amer. Math. Soc.}, \textbf{96},
no.\,1 (1986), 23--28) proved that there exists a number
$N_0=N_0(c)$ such that the word $x^ry^s$, $rs \ne 0$, is universal
on the alternating group ${\Bbb A}_n$ for any $ n \geq \max \, \{
N_0, \; \, c\cdot \log\, m(r,s)\} $, where $m(r,s)$ is the product
of all primes dividing $rs$ if $r,s \not\in \{ -1, 1\} $, and
$m(r,s)=1$ if $r,s \in \{ -1, 1\} $. For example, one can take $
N_0(5/2)=5$ and $ N_0(2)=29$. Find an analogous bound for the
degree of the symmetric group~${\Bbb S}_n$ under hypothesis that
at least one of the $r, s$ is odd: up to now, here only a more
crude bound $ n \geq \max \, \{ 6, \; \, 4m(r,s) - 4\} $ is known
(M.\,Droste, {\it Proc. Amer. Math. Soc.}, \textbf{96}, no.\,1
(1986), 18--22). \hfill \raisebox{-1ex}{\sl Yu.\,I.\,Merzlyakov}

\ul

\oai A logarithmic bound is found (unverified AI solution:  M.\,Kinyon, J.\,Urban, {\it Preprint}, 2026,
\url{https://kourovkanotebookorg.wordpress.com/wp-content/uploads/2026/09/kourovka-mathematics.pdf}).
\emp 

\bmp
\textbf{\zv 10.34.}
Does there exist a
non-soluble finite group which coincides with the product of any two
of its non-conjugate maximal subgroups?

\makebox[15pt][r]{}{\it Editors' comment:\/}
the answer is negative for almost simple groups
(T.\,V.\,Tikho\-nenko, V.\,N.\,Tyutyanov, {\it Siberian Math. J.}, {\bf 51},
no.\,1 (2010), 174--177).

\hfill \raisebox{-1ex}{\sl V.\,S.\,Monakhov}

\ul

\otv No, there are no such groups (R.\,Sater, {\it Preprint of 6 August 2026}, \url{https://arxiv.org/abs/2608.02970};
\ J.\,Li, Y.\,Yang, {\it Preprint of 19 August 2026}, \url{https://arxiv.org/abs/2608.19478}).
\emp

\bmp \textbf{10.35.} Is it true that every finitely generated
torsion-free subgroup of $GL_n({\Bbb C})$ is residually in the
class of torsion-free subgroups of $GL_n(\overline{{\Bbb Q}})$?
\hfill \raisebox{-1ex}{\sl G.\,A.\,Noskov} \emp

\bmp \textbf{10.36.} Is it true that $SL_n({\Bbb Z}[x_1,\ldots
,x_r])$, for $n$ sufficiently large, is a group of type $(FP)_m$
(which is defined in the same way as $(FP)_{\infty}$ was defined
in 6.3, but with that weakening that the condition of
being finitely generated is not imposed on the terms of the
resolution with numbers $ > m$). An affirmative answer is known
for $m=0$, $n \geq 3$ (A.\,A.\,Suslin, {\it Math. USSR Izvestiya},
\textbf{11} (1977), 221--238) and for $m=1$, $n \geq 5$
(M.\,S.\,Tulenbayev, {\it Math. USSR Sbornik}, \textbf{45} (1983),
139--154; \ U.\,Rehmann,
 C.\,Soul\'e, in: {\it Algebraic $K$-theory, Proc. Conf., Northwestern Univ.,
Evanston, Ill., 1976, $($Lect. Notes Math.}, \textbf{551}), Springer,
Berlin, 1976, 164--169). \hfill \raisebox{-1ex}{\sl G.\,A.\,Noskov}
\emp

\bmp \textbf{10.38.} (W.\,\,van der Kallen). Does $E_{n+3}({\Bbb
Z}[x_1,\ldots ,x_n])$, $n \geq 1$, have finite breadth with
respect to the set of transvections? \hfill \raisebox{-1ex}{\sl
G.\,A.\,Noskov} \emp

\bmp \textbf{10.39.} a) Is the membership problem soluble for the
subgroup $E_n(R)$ of $SL_n(R)$ where $R$ is a commutative ring?

\makebox[15pt][r]{}b) Is the word problem soluble for the groups
$K_i(R)$ where
$K_i$ are the Quillen $K $-func\-tors and $R$ is a commutative
ring? \hfill
\raisebox{-1ex}{\sl G.\,A.\,Noskov}
\emp

\bmp \textbf{10.40.} (H.\,Bass). Let $G$ be a group, $e$ an
idempotent matrix over ${\Bbb Z}G$ and let ${\rm tr}\,
e=\sum\limits_{g\in G}e_gg,\,\,\,\, e_g \in {\Bbb Z}.$ {\it Strong
conjecture:\/} for any non-trivial $x \in G$, the equation
$\sum\limits_{g\sim x} e_g=0$ holds where $\sim $ denotes
conjugacy in $G$. {\it Weak conjecture:\/} $\sum\limits_{g\in G}
e_g=0$. The strong conjecture has been proved for finite, abelian
and linear groups and the weak conjecture has been proved for
residually finite groups.

\makebox[15pt][r]{}{\it H.\,Bass' Comment of 2005:\/} Significant
progress has been made; a good up-to-date account is in the paper
(A.\,J.\,Berrick, I.\,Chatterji, G.\,Mislin, {\it Math. Ann.}, {\bf
329} (2004), 597--621). \hfill \raisebox{-1ex}{\sl G.\,A.\,Noskov}
\emp

\bmp \textbf{10.42.} (J.\,T.\,Stafford). Let $G$ be a poly-${\Bbb Z}
$-group and let $k$ be a field. Is it true that every finitely
generated projective $kG $-module is either a free module or an
ideal?

\hfill \raisebox{-1ex}{\sl G.\,A.\,Noskov}
\emp

\bmp
\textbf{\ai 10.43.} 
Let $R$ and $S$ be associative rings with
identity such that 2 is invertible in $ S$. Let $\Lambda _I \! :
GL_n(R) \rightarrow GL_n(R/I)$ be the homomorphism corresponding
to an ideal $I$ of $ R$ and let $ E_n(R)$ be the subgroup of $
GL_n(R)$ generated by elementary transvections $t_{ij}(x)$. Let
$a_{ij}=t_{ij}(1)t_{ji}(-1)t_{ij}(1)$, let the bar denote images
in the factor-group of $ GL_n(R)$ by the centre and let
$PG=\,\overline{\! G}$ for $G \leq GL_n(R)$. A homomorphism
\vspace{-1ex} $$ \Lambda \! : E_n(R) \rightarrow GL(W)=GL_m(S)$$
\vskip-1ex

is called {\it standard\/} if $S^m=P \oplus \cdots \oplus P \oplus
Q$ (a direct sum of $S $-modules in which there are $n$ summands
$P$) and \vspace{-1ex} $$\Lambda x=g^{-1}\tau (\delta
^{\displaystyle{*}} (x)f + (^t\delta ^{\displaystyle{*}}(x)^{\nu
})^{-1}(1 - f))g,\quad \quad x \in E_n(R),$$ \vskip-1ex

where $\delta ^{\displaystyle{*}} \! : GL_n(R) \rightarrow
GL_n({\rm End}\, P)$ is
the homomorphism
induced by a ring homomorphism $\delta \! : R \rightarrow
{\rm End}\, P$ taking
identity to identity, $g$ is an isomorphism of the module $W$
onto $
S^m$, $\tau \! : GL_n({\rm End}\, P) \rightarrow
GL(gW)$ is an
embedding, $f$ is a central idempotent of $\delta R$,
$t$ denotes
transposition and $\nu$ is an antiisomorphism of $\delta R$.
Let $n \geq 3$, $m \geq 2$. One can show that the
homomorphism
\vspace{-1ex}
$$\Lambda _0\! : PE_n(R) \rightarrow GL(W)=GL_m(S)$$
\vskip-1ex

is induced by a standard homomorphism $ \Lambda$ if $\Lambda
_0\,\overline{\! a}_{ij}=g^{-1}\tau (a^{\displaystyle{*}}_{ij})g$
for some $g$ and $\tau$ and for any $i \ne j$ where
$a^{\displaystyle{*}}_{ij}$ denotes the
matrix obtained from $a_{ij}$ by replacing 0 and 1 from $R$ by
0 and
1 from ${\rm End}\, P$. Find (at least in particular cases)
the form of
a homomorphism $\Lambda _0$ that does not satisfy the last
condition.
\hfill
 \raisebox{-1ex}{\sl V.\,M.\,Petechuk}

\ul

\oai Such homomorphisms are found (unverified AI solution:  Achyuth Jayadevan,  {\it Preprint}, 2026, \url{https://zenodo.org/records/22852216}).
\emp

\bmp \textbf{\ai 10.44.} 
Let $R$ and $S$ be associative rings with
identity such that 2 is invertible in~$S$. Prove that every standard homomorphism $\Lambda$
of $ E_n(R) \subset GL_n(R)=GL(V)$ into $E_m(S) \subset
GL_m(S)=GL(W)$ originates from a collineation and a correlation,
that is, has the form $\Lambda x=(g^{-1}xg)f + (h^{-1}xh)(1 - f),$
where $g$ is a semilinear isomorphism $V_R \rightarrow W_S$
(collineation) and $h$ is a semilinear isomorphism $V_R
\rightarrow _S\! W'_{S^0}$ (correlation).

\hfill
\raisebox{-1ex}{\sl V.\,M.\,Petechuk}

\ul

\oai A counterexample is constructed (unverified AI solution:  Achyuth Jayadevan,  {\it Preprint}, 2026, \url{https://zenodo.org/records/22847412}).
\emp

\bmp \textbf{10.45.}  Let $R$ be an associative ring with
identity.  Let $n \geq 3$ and suppose that $N$ is a
subgroup of $ GL_n(R)$ which is normalized by $E_n(R)$. Prove that
either $ N$ contains $E_n(R)$ or $\Lambda _I[N,E_n(R)]=1$ for a
suitable ideal $I \ne R$ of $R$. \hfill \raisebox{-1ex}{\sl
V.\,M.\,Petechuk} \emp

\bmp \textbf{\zv 10.46.}
 Let $R$ be an associative ring with
identity.
Prove that for every element $\sigma
\in GL_n(R)$, $n \geq 3$, there exist transvections $\tau _1, \tau _2$ such that
$[[\sigma ,\tau _1],\, \tau _2]$ is a unipotent element.
\hfill \raisebox{-1ex}{\sl V.\,M.\,Petechuk}

\ul

\otv This is not always true (K.\,Muliarchyk, {\it Preprint}, 2026,

\url{https://kourovkanotebookorg.wordpress.com/wp-content/uploads/2026/09/kourovka_10.46_solution.pdf}).
\emp

\bmp \textbf{10.47.} Describe the automorphisms of $PE_2(R)$ in the
case where $R$ is an associative commutative ring and 2 and 3 have inverses
in it. \hfill \raisebox{-1ex}{\sl V.\,M.\,Petechuk} \emp

\bmp
\textbf{10.49.}
Does there exist a group $G$ satisfying the following
four conditions:

\makebox[15pt][r]{}1) $G$ is simple, moreover, there is an
integer $n$
such
that $G=C^n$ for any conjugacy class $C$,

\makebox[15pt][r]{}2) all maximal abelian subgroups of $G$ are
conjugate in $
G$,

\makebox[15pt][r]{}3) every maximal abelian subgroup of $G$ is
self-normalizing and it is the centralizer of any of its
nontrivial elements,

\makebox[15pt][r]{}4) there is an integer $m$ such that if $H$
is a
maximal
abelian subgroup of $G$ and $a \in G\setminus H$ then every
element in
$G$ is a product of $m$ elements in $aH$?

\makebox[15pt][r]{}{\it Comment of 2005:\/} There is a partial
solution in (E.\,Jaligot, A.\,Ould~Houcine, {\it J.~Algebra}, {\bf
280} (2004), 772--796).
 \hfill \raisebox{-1ex}{\sl
B.\,Poizat}

\emp

\bmp \textbf{10.50.} Complex characters of a finite group $G$ induced
by linear characters of cyclic subgroups are called {\it induced
cyclic characters\/} of $G$. Computer-aided computations
(A.\,V.\,Rukolaine, {\it Abstracts of the 10th All--Union Sympos.
on Group Theory}, Gomel', 1986, p.~199 (Russian)) show that there
exist groups (for example, ${\Bbb S}_5$, $ SL_2(13)$, $M_{11}$)
all of whose irreducible complex characters are integral linear
combinations of induced cyclic characters and the principal
character of the group. Describe all finite groups with this
property. \hfill \raisebox{-1ex}{\sl A.\,V.\,Rukolaine} \emp

\bmp \textbf{10.51.} Describe the structure of thin abelian
$p\hskip0.2ex $-groups. Such a description is known in the class
of separable $ p\hskip0.2ex $-groups (C.\,Megibben, {\it Mich.
Math.~J.}, \textbf{13}, no.\,2 (1966), 153--160). \hfill
\raisebox{-1ex}{\sl S.\,V.\,Rychkov} \emp

\bmp \textbf{\ai 10.52.} 
(R.\,Mines). It is well known that the topology
on the completion of an abelian group under the $ p\hskip0.2ex
$-adic topology is the $p\hskip0.2ex $-adic topology. R.\,Warfield
has shown that this is also true in the category of nilpotent
groups. Find a categorical setting for this theorem which includes
the case of nilpotent groups. \hfill \raisebox{-1ex}{\sl
S.\,V.\,Rychkov}

\ul

\oai Such a setting is found (unverified AI solution:  Achyuth Jayadevan,  {\it Preprint}, 2026, \url{https://zenodo.org/records/22852281}).
 \emp 

\bmp \textbf{10.53.} (M.\,Dugas). Let ${\goth R}$ be the {\it Reid
class}, i.~e. the smallest containing ${\Bbb Z}$ and closed under
direct sums and direct products. Is ${\goth R}$ closed under direct
summands?

\hfill
\raisebox{-1ex}{\sl S.\,V.\,Rychkov}
\emp

\bmp \textbf{10.54.} (R.\,G\"obel). For a cardinal number $\mu$, let
\vspace{-1ex} $${\Bbb Z}^{<\mu}=\{ f \in {\Bbb Z}^{\mu}\mid |{\rm
supp}\, (f)| < \mu \} \; \; \; \; \mbox{ and} \; \; \; \;
G_{\mu}={\Bbb Z}^{\mu}/{\Bbb Z}^{<\mu}.$$ \vskip-1ex

\makebox[25pt][r]{a)} Find a non-zero direct summand $D$ of
$G_{\omega
_1}$ such
that $D\not\cong G_{\omega _1}$.

\makebox[25pt][r]{b)} Investigate the structure of $G_{\mu}$
(the
structure of
$G_{\omega _0}$ is well known).

\hfill \raisebox{-1ex}{\sl S.\,V.\,Rychkov}
\emp

\bmp \textbf{10.55.} (A.\,Mader).
 ``Standard $B$'', that is, $B={\Bbb Z}(p) \oplus {\Bbb Z}(p^2 ) \oplus
 \cdots$, is slender as a module over its endomorphism ring
 (A.\,Mader, in: {\it Abelian Groups and Modules,
 Proc., Udine, 1984}, Springer, 1984, 315--327). Which
 abelian $p\hs$-groups are slender as modules over their endomorphism
 rings? \hfill
 \raisebox{-1ex}{\sl S.\,V.\,Rychkov}

\emp

\bmp \textbf{\ai 10.57.} 
What are the minimal non-$A $-forma\-tions? An
{\it $A $-forma\-tion\/} is, by definition, the formation of the
finite groups all of whose Sylow subgroups are abelian.

\hfill
\raisebox{-1ex}{\sl A.\,N.\,Skiba}

\ul

\oai A classification is proposed (unverified AI solution: B.\,G.\,Di Siena, {\it Preprint}, 2026, \url{https://www.proofforum.org/abs/PF:2026.000037}).
\emp

\bmp \textbf{10.58.} Is the subsemigroup generated by the undecomposable
formations in the semigroup of formations of finite groups free?
\hfill \raisebox{-1ex}{\sl A.\,N.\,Skiba} \emp

\bmp \textbf{10.59.} Is a $p' $-group $G$ locally nilpotent if it
admits a splitting automorphism $\varphi$ of prime order $p$ such
that all subgroups of the form $\langle g, g^{\varphi}, \ldots ,
g^{\varphi ^{p-1}}\rangle $ are nilpotent? An automorphism
$\varphi$ of order $p$ is called {\it splitting\/} if $gg^{\varphi
}g^{\varphi ^{2}}\cdots\, g^{\varphi ^{p-1}}= 1$ for all $g \in
G$.

\hfill
\raisebox{-1ex}{\sl A.\,I.\,Sozutov}
\emp

\bmp \textbf{10.60.} Does every periodic group ($p\hs$-group) $A$ of
regular automorphisms of an abelian group have non-trivial
centre?

\makebox[15pt][r]{}{\it Editors' comment (2001):\/} The answer is
affirmative if $A$ contains an element of order 2 or 3
(A.\,Kh.\,Zhurtov, {\it Siberian Math. J.}, \textbf{41}, no.\,2
(2000), 268--275).

\hfill \raisebox{-1ex}{\sl A.\,I.\,Sozutov}
\emp

\bmp \textbf{10.61.} Suppose that $H$ is a proper subgroup of a group
$G$, $a \in H$, $a^2 \ne 1$ and for every $g \in G\setminus H$ the
subgroup $\left< a,\, a^g\right>$ is a Frobenius group whose
complement contains $a$. Does the set-theoretic union of the
kernels of all Frobenius subgroups of $G$ with complement $\left<
a\right>$ constitute a subgroup? For definitions see 6.55; see
also (A.\,I.\,Sozutov, {\it Algebra and Logic}, \textbf{34}, no.\,5
(1995), 295--305).

\makebox[15pt][r]{}{\it Editors' comment (2005):\/} The answer is
affirmative if the order of $a$ is even (A.\,M.\,Popov,
A.\,I.\,Sozutov, {\it Algebra and Logic}, \textbf{44}, no.\,1 (2005),
40--45)
 or if the order of $a$ is not 3 or 5 and the group $\left< a,\,
a^g\right>$ is finite for any $g\not\in H$ (A.\,M.\,Popov, {\it
Algebra and Logic}, \textbf{43}, no.\,2 (2004), 123--127). \hfill
\raisebox{-1ex}{\sl A.\,I.\,Sozutov} \emp

\bmp \textbf{\zv 10.62.}
Construct an example of a (periodic) group
without subgroups of index 2 which is generated by a conjugacy
class of involutions $X$ such that the order of the product of any
two involutions from $X$ is odd. \hfill \raisebox{-1ex}{\sl
A.\,I.\,Sozutov}

\ul

\otv Example constructed (AI solution confirmed by M.\,Amelio: \    M.\,Kinyon, J.\,Urban, {\it Preprint}, 2026,
\url{https://kourovkanotebookorg.wordpress.com/wp-content/uploads/2026/09/kourovka-mathematics.pdf}).
\emp 

\bmp \textbf{\zv 10.64.}
Does there exist a non-periodic doubly
transitive permutation group with a periodic stabilizer of a
point? \hfill \raisebox{-1ex}{\sl Ya.\,P.\,Sysak}

\ul

\otv Yes, it does exist (M.\,Amelio, {\it Preprint}, 2025, \url{https://arxiv.org/abs/2509.11958}).
\emp

\bmp \textbf{10.65.} Determine the structure of infinite 2-transitive
permutation groups $(G,\Omega )$ in which the stabilizer of a
point $\alpha \in \Omega$ has the form $G_{\alpha}=A\cdot
G_{\alpha \beta}$ where $G _{\alpha \beta}$ is the stabilizer of
two points $\alpha , \beta $, $\alpha \ne \beta $, such that $G
_{\alpha \beta}$ contains an element inverting the subgroup $A$.
Suppose, in particular, that $A\setminus \{ 1\}$ contains at most
two conjugacy classes of $G_{\alpha}$; does $G$ then possess a
normal subgroup isomorphic to $PSL_2$ over a field? \hfill
\raisebox{-1ex}{\sl A.\,N.\,Fomin} \emp

\bmp \textbf{10.67.} The class $LN{\goth M}_p$ of locally nilpotent
groups admitting a splitting automorphism of prime order $p$ (for
definition see 10.59) is a variety of groups with operators
(E.\,I.\,Khukhro, {\it
Math. USSR Sbornik}, \textbf{58} (1987), 119--126). Is it true that
\vspace{-1ex} $$LN{\goth M}_p=({\goth N}_{c(p)} \cap LN{\goth
M}_p) \vee ({\goth B}_p \cap LN{\goth M}_p)$$ \vskip-1ex

where ${\goth N}_{c(p)}$ is the variety of nilpotent groups of
some $p\hskip0.2ex $-bounded class $c(p)$ and ${\goth B}_p$ is the
variety of groups of exponent $p$? \hfill \raisebox{-1ex}{\sl
E.\,I.\,Khukhro} \emp

\bmp \textbf{10.70.} Find a geometrical justification for the
Whitehead method for free products similar to the substantiation
given for free groups by Whitehead himself and more visual than
given in (D.\,J.\,Collins, H.\,Zieschang, {\it Math.~Z.}, {\bf
185}, no.\,4 (1984), 487--504; \textbf{186}, no.\,3 (1984),
335--361).

\makebox[15pt][r]{}{\it D.\,J.\,Collins' comment:} there is a partial
solution in (D.\,McCullough, A.\,Miller, {\it Symmetric automorphisms
of free products}, Mem. Amer. Math. Soc. \textbf{582} (1996)).

 \hfill \raisebox{-1ex}{\sl H.\,Zieschang} \emp

\bmp \textbf{10.71.} Is it true that the centralizer of any
automorphism (any finite set of automorphisms) in the automorphism
group of a free group of finite rank is a finitely presented
group? This is true in the case of rank 2 and in the case of inner
automorphisms for any finite rank. \hfill \raisebox{-1ex}{\sl
V.\,A.\,Churkin} \emp

\bmp
\textbf{10.73.}
Enumerate all formations of finite groups all of
whose subformations are $S_n $-closed. \hfill
\raisebox{-1ex}{\sl L.\,A.\,Shemetkov}

\emp

\bmp \textbf{10.74.} Suppose that a group $G$ contains an element $a$
of prime order such that its centralizer $ C_G(a)$ is finite and
all subgroups $\left< a,\, a^g\right> $, $g \in G$, are finite and
almost all of them are soluble. Is $G$ locally finite? This
problem is closely connected with 6.56. The question was solved in
the positive for a number of very important partial cases by the
author ({\it Abstracts on Group Theory of the Mal'cev Int.
Conf. on Algebra}, Novosibirsk, 1989, p.\,145 (Russian)).
\hfill \raisebox{-1ex}{\sl V.\,P.\,Shunkov}
\emp

\bmp \textbf{10.75.} Suppose that a group $G$ contains an element $a$
of prime order $p$ such that the normalizer of every finite
subgroup containing $a$ has finite periodic part and all subgroups
$\left< a,\, a^g\right> $, $g \in G$, are finite and almost all of
them are soluble. Does $G$ possess a periodic part if $p > 2$? It
was proved in (V.\,P.\,Shunkov, {\it Groups with involutions},
Preprints no.\,4, 5, 12 of the Comput. centre of SO AN SSSR,
Krasnoyarsk, 1986 (Russian)) that if $a$ is a point, then the
answer is affirmative; on the other hand, a group with a point $a$
of order 2 satisfying the given hypothesis, which has no periodic
part, was exhibited in the same works. For the definition of a
{\it point\/} see (V.\,I.\,Senashov, V.\,P.\,Shunkov, {\it Algebra
and Logic}, \textbf{22}, no.\,1 (1983), 66--81). \hfill
\raisebox{-1ex}{\sl V.\,P.\,Shunkov} \emp

\bmp \textbf{10.77.} Suppose that $G$ is a periodic group containing
an elementary abelian subgroup $R$ of order 4. Must $G$ be locally
finite

\makebox[25pt][r]{a)} if $C_G(R)$ is finite?

\makebox[25pt][r]{b)} if the centralizer of every involution of
$R$ in $G$ is a
Chernikov group?

\makebox[15pt][r]{}{\it Remark of 1999:\/} P.\,V.\,Shumyatsky
({\it Quart. J.~Math. Oxford (2)}, \textbf{49}, no.\,196 (1998),
491--499) gave a positive answer to the question a) in the case
where $G$ is residually finite. \hfill \raisebox{-1ex}{\sl
V.\,P.\,Shunkov} \emp

\bmp \textbf{10.78.} Does there exist a non-Chernikov group which is
a product of two Chernikov subgroups? \hfill \raisebox{-1ex}{\sl
V.\,P.\,Shunkov} \emp

\raggedbottom

\newpage

\pagestyle{myheadings} \markboth{11th Issue
(1990)}{11th Issue (1990)}
\thispagestyle{headings}

~
\vspace{2ex}

\centerline {\Large \textbf{Problems from the 11th Issue (1990)
}}
\phantomsection\label{11izd}
\vspace{4ex}

 \bmp \textbf{11.1.} (Well-known problem).
Describe the structure of the centralizers of unipotent elements
in almost simple groups of Lie type. \hfill \raisebox{-0.6ex}{\sl
R.\,Zh.\,Aleev}

\emp

 \bmp \textbf{11.3.} (M.\,Aschbacher). A
$p\hskip0.2ex$-local subgroup $H$ in a group $G$ is said to be a
{\it superlocal\/} if $H= N_G(O_p(H))$. Describe the superlocals
in alter\-nat\-ing groups and in groups of Lie type. \hfill
\raisebox{-0.6ex}{\sl R.\,Zh.\,Aleev}

\emp

\bmp \textbf{\zv 11.5.}
  Let $k$ be a commutative
ring and let $G$ be a torsion-free almost poly\-cyclic group.
Suppose that $ P$ is a finitely generated projective module over
the group ring $kG$ and $P$ contains two elements independent over
$kG$. Is $P$ a free module?

\hfill \raisebox{-0.6ex}{\sl V.\,A.\,Artamonov}

\ul

\otv No, not always (AI solution confirmed by D.\,Segal:  \ Achyuth Jayadevan,  {\it Preprint}, 2026, \url{https://zenodo.org/records/22846603}).
\emp 

 \bmp \textbf{11.7.} a) Find conditions for a
group $ G$, given by its presentation, under which the
property\hfill {\it some term of the lower central series of $G$
is a free group} \hfill $(*)$ \linebreak implies the residual
finiteness of $G$.

 \makebox[15pt][r]{}b) Find conditions on the structure of a
residually finite group $G$ which ensure property $(*)$ for $G$.
\hfill \raisebox{-0.6ex}{\sl K.\,Bencs\'ath}
\emp

 \bmp \textbf{11.8.} b) For a finite group $X$,
let $\chi _1(X)$ denote the totality of the degrees of all
irre\-du\-cible complex characters of $X$ with allowance for their
multiplicities. Suppose that $\chi _1(G)=\chi _1(H)$ for groups
$G$ and $H$. Clearly, then $|G|=|H|$. Is it true that
$H$ is soluble if $G$ is soluble?

 \makebox[15pt][r]{}It is known that $H$ is a Frobenius group if $G$ is a Frobenius
group.

\hfill \raisebox{-1ex}{\sl Ya.\,G.\,Berkovich}
\emp

 \bmp \textbf{\zv 11.9.}
(I.\,I.\,Pyatetski\u{\i}-Shapiro). Does there exist a finite
non-soluble group $G$ such that the set of characters induced by
the trivial characters of representatives of all conjugacy classes
of subgroups of $G$ is linearly independent? \hfill
\raisebox{-0.6ex}{\sl Ya.\,G.\,Berkovich}

\ul

\otv No, there are no such groups (M.\,Brescia, E.\,Ingrosso, M.\,Trombetti,  {\it Preprint}, 2026, \url{https://arxiv.org/abs/2608.16957}).
\emp

 \bmp \textbf{11.10.} (R.\,C.\,Lyndon). a) Does
there exist an algorithm that, given a group word $ w(a,x)$,
recognizes whether $a$ is equal to the identity element in the
group $\left< a,x\mid a^n=1,\right.$ $\left. w(a,x)=1\right>$?
\hfill \raisebox{-0.6ex}{\sl V.\,V.\,Bludov}

\emp

 \bmp \textbf{11.12.} a) Suppose that $G$ is a
simple locally finite group in which the centralizer of some
element is a linear group, that is, a group admitting a faithful
matrix representation over a field. Is $G$ itself a linear group?

\makebox[15pt][r]{}b) The same question with replacement of the
word ``linear'' by ``finitary linear''. A group $H$ is {\it finitary
linear\/} if it admits a faithful representation on an
infinite-dimensional vector space $V$ such that the residue
subspaces $V(1-h)$ have finite dimensions for all $h\in H$.\hfill
\raisebox{-1ex}{\sl A.\,V.\,Borovik}

\emp

 \bmp \textbf{11.14.} Is every finite simple
group characterized by its Cartan matrix over an algebraically
closed field of characteristic 2? In (R.\,Brandl, {\it Arch.
Math.}, \textbf{38} (1982), 322--323) it is shown that for any given
finite group there exist only finitely many finite groups with
the same Cartan matrix. \hfill \raisebox{-1ex}{\sl R.\,Brandl}

\emp

\bmp \textbf{\ai 11.15.} 
 It is known that for each
prime number $p$ there exists a series $a_1, a_2, \ldots $ of
words in two variables such that the finite group $G$ has abelian
Sylow $p\hskip0.2ex$-sub\-groups if and only if $ a_k(G)=1$ for
almost all $k$. For $p=2$ such a series is known explicitly
(R.\,Brandl, {\it J.~Austral. Math. Soc.}, \textbf{31} (1981),
464--469). What about $p
> 2$?

\hfill \raisebox{-1ex}{\sl R.\,Brandl}

\ul

\oai A computable construction for such a series is found (unverified AI solution:  Achyuth Jayadevan,  {\it Preprint}, 2026, \url{https://zenodo.org/records/22852259}).
\emp

 \bmp \textbf{11.16.} Let $V_r$ be the class of
all finite groups $G$ satisfying a law $[x,\, _r y]=[x,\, _s y]$
for some $s=s(G) > r$. Here $[x,\, _1 y] = [x,y]$ and $[x,\,
_{i+1} y]=[[x,\, _i y],y]$.

 \makebox[25pt][r]{a)} Is there a function $f$ such that every
soluble group in $V_r$ has Fitting length $< f(r)$? For $r < 3$ see
(R.\,Brandl, {\it Bull. Austral. Math. Soc.},
\textbf{28} (1983), 101--110).

 \makebox[25pt][r]{b)} Is it true that $V_r$ contains only finitely
many nonabelian simple groups? This is true for $r < 4$. \hfill
\raisebox{-1ex}{\sl R.\,Brandl}

\emp

 \bmp \textbf{11.17.} Let $G$ be a finite group
and let $d=d(G)$ be the least positive integer such that $G$
satisfies a law $[x,\, _ry]=[x,\, _{r+d} y]$ for some nonnegative
integer $r=r(G)$.

\makebox[25pt][r]{a)} Let $e=1$ if $d(G)$ is even and $e=2$
otherwise. Is it true that the exponent of $G/F(G)$ divides $e\cdot
d(G)$?

\makebox[25pt][r]{b)} If $G$ is a nonabelian simple group, does the exponent of $G$ divide $d(G)$?

\makebox[15pt][r]{}Part a) is true for soluble groups (N.\,D.\,Gupta,
\,H.\,Heineken, {\it Math.~Z.}, \textbf{95} (1967), 276--287). I have
checked part b) for ${\Bbb A}_n$, $PSL(2,q)$, and a number of
sporadic groups.

\hfill \raisebox{-1ex}{\sl R.\,Brandl}

\emp

 \bmp \textbf{11.18.} Let $G(a,b)=\left<
x,y\mid x=[x,\, _a y], \; y=[y,\, _b x]\right> $. Is $G(a,b)$
finite?

\makebox[15pt][r]{}It is easy to show that $G(1,b)=1$ and one can show that
$G(2,2)=1$. Nothing is known about $G(2,3)$. If one could show
that every minimal simple group is a quotient of some $ G(a,b)$,
then this would yield a very nice sequence of words in two
variables to characterize soluble groups, see (R.\,Brandl,
J.\,S.\,Wilson, {\it J.~Algebra}, \textbf{116} (1988), 334--341.)
\hfill \raisebox{-1ex}{\sl R.\,Brandl}

\emp

 \bmp \textbf{11.19.} (C.\,Sims). Is the $n$th
term of the lower central series of an absolutely free group the
normal closure of the set of basic commutators (in some fixed free
generators) of weight exactly $n$?

\makebox[15pt][r]{}D.\,Jackson announced a
positive answer for $n\leq 5$. \hf \raisebox{-1ex}{\sl
A.\,Gaglione, D.\,Spellman}

\emp

 \bmp \textbf{11.22.} Characterize all
$p\hskip0.2ex$-groups, $p$ a prime, that can be faithfully
represented as $n\times n$ triangular matrices over a division
ring of characteristic~$p$.

\makebox[15pt][r]{}If $p \geq n$ the solution is given in
(B.\,A.\,F.\,Wehrfritz, {\it Bull. London Math.
Soc.}, \textbf{19} (1987), 320--324). The corresponding question with
$p=0$ was solved by \,B.\,Hartley and
\,P.\,Menal ({\it Bull. London Math. Soc.}, {\bf
15} (1983), 378--383), but see above reference for a second proof.
\hfill \raisebox{-1ex}{\sl B.\,A.\,F.\,Wehrfritz}

\emp

 \bmp \textbf{\zv 11.23.}
 An automorphism $\varphi$
of a group $G$ is called a {\it nil-automorphism\/} if, for every
$a \in G$, there exists $n$ such that $[a,\, _n \varphi ] = 1$.
Here $[x,\, _1 y]=[x,y]$ and $[x,\, _{i+1} y]=[[x,\, _i y],y]$. An
automorphism $\varphi$ is called an {\it $e$-auto\-mor\-phism\/}
if, for any two $\varphi\hskip0.1ex$-inva\-ri\-ant subgroups $A$
and $B$ such that $A \not\subseteq B$, there exists $a \in
A\setminus B$ such that $[a,\varphi ] \in B$. Is every
$e$-auto\-mor\-phism of a group a nil-automorphism? \hfill
\raisebox{-1ex}{\sl V.\,G.\,Vilyatser}

\ul

\otv No, not every (K.\,Muliarchyk, {\it Preprint}, 2026, \url{https://kourovkanotebookorg.wordpress.com/wp-content/uploads/2026/09/kourovka_11_23_solution.pdf}).
\emp

 \bmp \textbf{11.25.} b) Do there exist local
Fitting classes which are decomposable into a non-trivial product
of Fitting classes and in every such a decomposition all factors
are non-local? \ For the definition of the {\it product\/} of
Fitting classes see (N.\,T.\,Vorob'\"ev, {\it Math. Notes}, {\bf
43}, no.\,2 (1988), 91--94). \hfill \raisebox{-1ex}{\sl
N.\,T.\,Vorob'\"ev}

\emp

 \bmp \textbf{11.28.} Suppose the prime graph
of the finite group $G$ is disconnected. (This means that the set
of prime divisors of the order of $G$ is the disjoint union of
non-empty subsets $\pi$ and $\pi '$ such that $G$ contains no
element of order $pq$ where $p \in \pi$, $q \in \pi '$.) Then
P.\,A.\,Linnell ({\it Proc. London Math. Soc.}, \textbf{47}, no.\,1
(1983), 83--127) has proved mod CFSG that there is a decomposition
of ${\Bbb Z}G$-modules ${\Bbb Z} \oplus {\Bbb Z}G=A \oplus B$ with
$A$ and $B$ non-projective. Find a proof independent of CFSG.
\hfill \raisebox{-1ex}{\sl K.\,Gruenberg}

\emp

 \bmp
\textbf{11.29.} Let $F$ be a free group
and ${\goth f}={\Bbb Z}F(F - 1)$ the augmentation ideal of the
integral group ring ${\Bbb Z}F$. For any normal subgroup $R$ of
$F$ define the corresponding ideal ${\goth r}={\Bbb Z}F(R - 1)=\,
_{{\rm id}} (r - 1\mid r \in R)$. One may identify, for instance,
$F \cap (1 + {\goth r} {\goth f})=R'$, where $F$ is naturally
imbedded into ${\Bbb Z}F$ and $1 + {\goth r} {\goth f}=\{ 1 +
a\mid a \in {\goth r}{\goth f}\} $.

\makebox[15pt][r]{}Identify in an analogous way in terms of
corresponding subgroups of $F$
\vskip0.3ex

\makebox[24pt][r]{a)} $F \cap (1 + {\goth r}_1 {\goth r_2}\cdot
\cdot \cdot {\goth r}_n)$, where $R_i$ are normal subgroups of $F$,
$i=1,\,2,\ldots ,n$;

\vskip0.3ex

\makebox[24pt][r]{b)} $F \cap (1 + {\goth r}_1{\goth r}_2{\goth
r}_3)$;

\vskip0.3ex

\makebox[24pt][r]{c)} $F \cap (1 + {\goth f} {\goth s} + {\goth
f}^n)$, where $F/S$ is finitely generated nilpotent;

\vskip0.3ex

\makebox[24pt][r]{d)} $F \cap (1 + {\goth f}{\goth s}{\goth f} +
{\goth f}^n)$;

\vskip0.3ex

\makebox[24pt][r]{e)} $F \cap (1 + {\goth r}(k) + {\goth f}^n)$, $n
> k \geq 2$, where
${\goth r}(k)={\goth r}{\goth f}^{k-1} + {\goth f}{\goth r} {\goth
f}^{k-2} + \cdots + {\goth f}^{k-1}{\goth r}$.

 \hfill
\raisebox{-1ex}{\sl N.\,D.\,Gupta}

\emp

 \bmp
\textbf{11.30.} Is it true that the rank of a torsion-free soluble group
is equal to the rank of any of its subgroups of finite index? The
answer is affirmative for groups having a rational series
(D.\,I.\,Zaitsev, in: {\it Groups with restrictions on subgroups},
Naukova dumka, Kiev, 1971, 115--130 (Russian)). We note also that
every torsion-free soluble group of finite rank contains a subgroup
of finite index which has a rational series.

\hfill \raisebox{-1ex}{\sl D.\,I.\,Zaitsev}

\emp

 \bmp \textbf{11.31.} Is a radical group
polycyclic if it is a product of two polycyclic subgroups? The
answer is affirmative for soluble and for hyperabelian groups
(D.\,I.\,Zaitsev, {\it Math. Notes}, \textbf{29}, no.\,4 (1981),
247--252; \ J.\,C.\,Lennox, J.\,E.\,Roseblade, {\it Math.~Z.},
\textbf{170} (1980), 153--154). \hfill \raisebox{-1ex}{\sl
D.\,I.\,Zaitsev}

\emp

\bmp \textbf{\zv 11.32.}
 Describe the primitive
finite linear groups that contain a matrix with {\it simple spectrum},
that is, a matrix all of whose eigenvalues are of multiplicity 1.
See partial results in (A.\,Zalesskii, I.\,D.\,Suprunenko, {\it
Commun. Algebra}, {\bf
 26}, no.\,3 (1998), 863--888; \textbf{28}, no.\,4
(2000), 1789--1833; \ Ch.\,Rudloff, A.\,Zalesski, {\it J.~Group
Theory}, \textbf{10} (2007), 585--612; \
L.\,Di\,Martino, A.\,E.\,Zalesski, {\it J. Algebra Appl.}, {\bf 11}, no.\,2 (2012), 1250038; \ L.\,Di\,Martino, M.\,Pellegrini, A.\,Zalesski, {Commun. Algebra}, {\bf 42} (2014), 880--908,
where the problem is considered in a more general context, when all but one eigenvalue
have multiplicity 1.). \hfill \raisebox{-1ex}{\sl
A.\,E.\,Zalesski\u{\i}}

\ul

\otv The project has been completed in (L.\,Di\,Martino, M.\,Pellegrini, A.\,Zalesski,
 {\it J.~Group Theory}, {\bf 23} (2020), 235--285; \ A.\,Zalesski,
{\it European J.~Math.}, {\bf 10} (2024), article no.\,55; \ A.\,Zalesski, {\it Preprint}, 2025, \url{https://arxiv.org/abs/2509.05770}).
\emp

\bmp \textbf{11.33.} b) Let $G(q)$ be a simple Chevalley group over a
field of order $q$. Prove that there exists $m$ such that the
restriction of every non-one-dimensional representation of
$G(q^m)$ over a field of prime characteristic not dividing $q$ to
$G(q)$ contains all irreducible representations of $G(q)$ as
composition factors.
 \hfill
\raisebox{-1ex}{\sl A.\,E.\,Zalesski\u{\i}}

\emp

 \bmp \textbf{11.34.}
Describe the complex representations of quasisimple finite groups
which remain irreducible after reduction modulo any prime number
$q$. An important example: representations of degree $(p^k - 1)/2$
of the symplectic group $Sp(2k,p)$ where $k \in {\Bbb N}$ and $p$
is an odd prime. See partial results in
(P.\,H.\,Tiep, A.\,E.\,Zalesski, {\it Proc. London Math. Soc.},
\textbf{84} (2002), 439--472; \
P.\,H.\,Tiep, A.\,E.\,Zalesski, {\it Proc. Amer. Math. Soc.}, {\bf 130} (2002), 3177--3184); \ C.\,Parker, M.\,van\,Beek,
{\it Preprint}, 2024, \url{https://arxiv.org/abs/2411.16379}). \hfill \raisebox{-1ex}{\sl
A.\,E.\,Zalesski\u{\i}}

\emp

 \bmp \textbf{11.36.} Let $G=B(m,n)$ be the
free Burnside group of rank $m$ and of odd exponent $n \gg 1$. Are
the following statements true?

\makebox[25pt][r]{a)} Every 2-generated subgroup of $G$ is
isomorphic to the Burnside $n$-product of two cyclic groups.

\makebox[25pt][r]{b)}
Every automorphism $\varphi$ of $G$ such that
$\varphi ^n=1$ and $b^{\varphi}\cdot b^{\varphi ^2} \cdots\,
 b^{ \varphi ^n}=1$ for all $b \in G$ is an inner
automorphism (here $m > 1$).

\makebox[15pt][r]{}{\it Editors' comment:} for large prime $n$ this follows from
(E.\,A.\,Cherepanov, {\it Int. J. Algebra Comput.}, {\bf 16} (2006), 839--847), for prime $n\geq 1009$ from
(V.\,S.\,Atabekyan, {\it Izv. Math.}, {\bf 75}, no.\,2 (2011), 223--237);
this is also true for odd $n\geq 1003$ if in addition the order of $\varphi$ is a prime power (V.\,S.\,Atabekyan,
{\it Math. Notes}, {\bf 95}, no.\,5 (2014), 586--589).

\makebox[25pt][r]{c)} The group $G$ is Hopfian if $m < \infty$.

\makebox[25pt][r]{d)} All retracts of $G$ are free.
 \hfill \raisebox{0ex}{\sl S.\,V.\,Ivanov}

\emp

 \bmp \textbf{11.37.} a) Can the free Burnside
group $ B(m,n)$, for any $m$ and $ n$, be given by defining relations
of the form $v^n=1$ such that for any natural divisor $d$ of $n$
distinct from $n$ the element $v^d$ is not trivial in $B(m,n)$? This
is true for odd $n \geq 665$, and for all
 $n\geq 2^{48}$ divisible by
$2^9$.\hfill \raisebox{-1ex}{\sl S.\,V.\,Ivanov}

\emp

 \bmp \textbf{11.38.} Does there exist a
finitely presented Noetherian group which is not almost
polycyclic? \hfill \raisebox{-1ex}{\sl S.\,V.\,Ivanov}

\emp

 \bmp \textbf{11.39.} (Well-known problem).
Does there exist a group which is not almost poly\-cyclic and
whose integral group ring is Noetherian? \hfill
\raisebox{-1ex}{\sl S.\,V.\,Ivanov}

\emp

 \bmp \textbf{11.40.} Prove or disprove that a
torsion-free group $G$ with the small cancellation condition
$C'(\lambda )$ where $\lambda \ll 1$ necessarily has the ${\cal
UP}$-property (and therefore $KG$ has no zero divisors). \hfill
\raisebox{-1ex}{\sl S.\,V.\,Ivanov}

\emp

 \bmp \textbf{11.44.} For a finite group $X$,
we denote by $r(X)$ its sectional rank. Is it true that the
sectional rank of a finite $p\hs$-group, which is a product $AB$
of its subgroups $A$ and $ B$, is bounded by some linear function
of $r(A)$ and $r(B)$? \hfill \raisebox{-1ex}{\sl L.\,S.\,Kazarin}

\emp

 \bmp \textbf{11.45.} A $t$-$(v,k,\lambda )$
{\it design\/} ${\cal D}=(X,{\cal B})$ contains a set $X$ of $v$
points and a set ${\cal B}$ of $k$-element subsets of~$X$ called
{\it blocks\/} such that each $t$-element subset of~$X$ is
contained in $\lambda$ blocks. Prove that there are no nontrivial
block-transitive 6-designs. (We have shown that there are no
nontrivial block-transitive
 8-designs
and there are certainly some block-transitive, even
flag-transitive, 5-designs.)

\makebox[15pt][r]{}{\it Comment of 2009:}
In ({\it Finite Geometry and Combinatorics (Deinze 1992)},
Cambridge Univ. Press, 1993, 103--119) we showed that a
block-transitive group $G$ on a nontrivial 6-design is either an
affine group $AGL(d,2)$ or is between $PSL(2,q)$ and $P\Gamma
L(2,q)$; in (M.\,Huber, {\it J.~Combin. Theory Ser. A}, {\bf 117}, no.\,2 (2010), 196--203)
it is shown that for the case $\lambda=1$ the group $G$ may only
be $P\Gamma L(2,p^e)$, where $p$ is $2$ or $3$ and $e$ is an odd
prime power. {\it Comment of 2013:}
In the case $\lambda=1$ there are no block-transitive $7$-designs (M.\,Huber, {\it Discrete Math. Theor. Comput. Sci.}, {\bf 12}, no.\,1 (2010), 123--132). {\it Comment of 2021:}
There are no nontrivial $6$-designs with $k\leq 10^4$ in the case where the automorphism group is almost simple (Q.\,Tan, W.\,Liu, J.\,Chen, {\it Algebra Colloq.}, {\bf 21}, no.\,2 (2014), 231--234).
\hfill \raisebox{-1ex}{\sl P.\,J.\,Cameron,
C.\,E.\,Praeger}

\emp

 \bmp \textbf{11.46.} a) Does
there exist a finite 3-group $G$ of nilpotency class 3 with the
property $[a,\, a^{\varphi}]=1$ for all $a \in G$ and all
endomorphisms $\varphi$ of $G$? (See \,A.\,Caranti, {\it
J.~Algebra}, \textbf{97}, no.\,1 (1985), 1--13.)
\hfill
\raisebox{-1ex}{\sl A.\,Caranti}
\emp

 \bmp \textbf{11.48.} Is the commutator $[x,\,
y,\, y,\, y,\, y,\, y,\, y]$ a product of fifth powers in the free
group $\left< x,y\right>$? If not, then the Burnside group
$B(2,5)$ is infinite. \hfill \raisebox{-1ex}{\sl A.\,I.\,Kostrikin}

\emp

 \bmp \textbf{11.49.} \,B.\,Hartley ({\it Proc.
London Math. Soc.}, \textbf{35}, no.\,1 (1977), 55--75) constructed
an example of a non-countable Artinian ${\Bbb Z}G$-module where
$G$ is a met\-abe\-lian group with the minimum condition for
normal subgroups. It follows that there exists a non-countable
soluble group (of derived length 3) satisfying Min-$n$. The
following question arises in connection with this result and with
the study of some classes of soluble groups with the weak minimum
condition for normal subgroups. Is an Artinian ${\Bbb Z}G$-module
countable if $G$ is a soluble group of finite rank (in particular,
a minimax group)? \hfill \raisebox{-1ex}{\sl L.\,A.\,Kurdachenko}

\emp

\bmp
 \textbf{11.50.} Let $A, C$ be abelian
groups. If $A[n]=
 0$, i.~e. for $a \in A$, $na=0$ implies $a=0$, then the
sequence \ $\displaystyle{
 \frac{{\rm Hom}\, (C,A)}{n{\rm Hom}\, (C,A)}\; \;
 \rightarrowtail {\rm Hom} \left(\frac{C}{nC},\; \frac{A}{nA}\right)
 \twoheadrightarrow {\rm Ext}\, (C,A)[n]}
 $ \
 is exact. Given \ $ \displaystyle{f \in {\rm
Hom}\left(C,\frac{A}{nA}\right)={\rm Hom}
 \left(\frac{C}{nC},\; \frac{A}{nA}\right)
} $ \ the corresponding extension $X_f$ is obtained as a pull-back
\vskip-1ex
 $$\displaystyle{
 \begin{array}{ccccc}A& \rightarrowtail &X_f& \twoheadrightarrow & C \cr
 \downarrow &\vphantom{_{_j}}&\downarrow &&\downarrow \cr
 A& \rightarrowtail &A& \twoheadrightarrow
 &\displaystyle{\frac{A}{nA}} \end{array}}.
 $$
 Use this scheme to classify certain extensions of $A$ by $C$.
The
 case $nC=0$, $A$~being torsion-free is interesting. Here
 ${\rm Ext}\, (C,A)[n] = {\rm Ext}\, (C,A)$. (See E.\,L.\,Lady,
 A.\,Mader, {\it J.~Algebra}, \textbf{140} (1991),
 36--64.) \hfill
 \raisebox{-1ex}{\sl A.\,Mader}
\emp

 \bmp \textbf{11.51.} Are there (large,
non-trivial) classes ${\goth X}$ of torsion-free abelian groups
such that for $A, C \in {\goth X}$ the group ${\rm Ext}\, (C,A)$
is torsion-free? It is a fact (E.\,L.\,Lady, A.\,Mader, {\it
J.\,Algebra}, {\bf
 140} (1991), 36--64) that two groups of such a class are
nearly isomorphic if and only if they have equal
$p\hskip0.2ex$-ranks for all $p$. \hfill \raisebox{-1ex}{\sl
A.\,Mader} \emp

 \bmp

\textbf{11.56.} a) Does every infinite residually finite group
contain an infinite abe\-li\-an subgroup? This is equivalent to
the following: does every infinite re\-si\-du\-al\-ly finite group
contain a non-identity element with an infinite centralizer?

\makebox[15pt][r]{}By a famous theorem of Shunkov
 a torsion group with an involution having a
finite centralizer is a virtually soluble group. Therefore we may
assume that in our group all elements have odd order. One should
start, perhaps, with the following:

\makebox[15pt][r]{}b) Does every infinite residually
$p\hskip0.2ex$-group contain an infinite abelian subgroup?

\hfill
\raisebox{-1ex}{\sl A.\,Mann}

\emp

 \bmp \textbf{11.58.} Describe the finite
groups which contain a tightly embedded subgroup $H$ such that a
Sylow 2-sub\-group of $H$ is a direct product of a quaternion
group of order 8 and a non-trivial elementary group. \hfill
\raisebox{-1ex}{\sl A.\,A.\,Makhn\"ev}

\emp

 \bmp \textbf{11.59.} A $TI$-sub\-group $A$ of a
group $G$ is called a {\it subgroup of root type\/} if $[A,\,
A^g]=1$ whenever $N_A(A^g) \ne 1$. Describe the finite groups
containing a cyclic subgroup of order 4 as a subgroup of root
type. \hfill \raisebox{-1ex}{\sl A.\,A.\,Makhn\"ev}

\emp

 \bmp \textbf{11.60.} Is it true that the
hypothetical Moore graph with 3250 vertices of valence 57 has no
automorphisms of order 2? \,M.\,Aschbacher ({\it J.~Algebra}, {\bf
19}, no.\,4 (1971), 538--540) proved that this graph is not a
graph of rank 3. \hfill \raisebox{-1ex}{\sl A.\,A.\,Makhn\"ev}

\emp

 \bmp \textbf{11.61.} Let $F$ be a non-abelian
free pro-$p\hskip0.2ex$-group. Is it true that the subset $\{ r
\in F\mid {\rm cd}\, (F/(r)) \leq 2\} $ is dense in $F$? Here $
(r)$ denotes the closed normal subgroup of $F$ generated by $r$.
\hfill \raisebox{-1ex}{\sl O.\,V.\,Mel'nikov}

\emp

 \bmp \textbf{11.62.} Describe the groups over
which any equation is soluble. In particular, is it true that this
class of groups coincides with the class of torsion-free groups?
It is easy to see that a group, over which every equation is
soluble, is torsion-free. On the other hand, S.\,D.\,Brodski\u{\i}
 ({\it Siberian Math.~J.}, \textbf{25}, no.\,2 (1984), 235--251)
showed that any equation is soluble over a locally indicable
group. \hfill \raisebox{-1ex}{\sl D.\,I.\,Moldavanski\u{\i}}

\emp

 \bmp \textbf{11.63.} Suppose that $G$ is a
one-relator group containing non-trivial elements of finite order
and $N$ is a subgroup of $G$ generated by all elements of finite
order. Is it true that any subgroup of $G$ that intersects $N$
trivially is a free group? One can show that the answer is
affirmative in the cases where $ G/N$ has non-trivial centre or
satisfies a non-trivial identity. \hfill \raisebox{-1ex}{\sl
D.\,I.\,Moldavanski\u{\i}}

\emp

 \bmp \textbf{11.65.} {\it Conjecture:\/} any
finitely generated soluble torsion-free pro-$p\hskip0.2ex$-group
with decidable elementary theory is an analytic
pro-$p\hskip0.2ex$-group.

\hfill \raisebox{-1ex}{\sl A.\,G.\,Myasnikov,
V.\,N.\,Remeslennikov}

\emp

 \bmp \textbf{11.66.} (Yu.\,L.\,Ershov). Is the
elementary theory of a free pro-$p\hskip0.2ex$-group decidable?

\hfill
\raisebox{-1ex}{\sl A.\,G.\,Myasnikov, V.\,N.\,Remeslennikov}
\emp

 \bmp
 \textbf{\zv 11.67.}
 Does there exist a
torsion-free group
with exactly 3 classes of conjugate elements such that no
non-trivial conjugate class contains a pair of inverse elements?

\hfill \raisebox{-1ex}{\sl B.\,Neumann}

\ul

\otv Yes, there does (V.\,Bagayoko, \url{https://arxiv.org/abs/2509.09186}).
\emp

 \bmp \textbf{11.69.} A group $G$ acting on a
set $\Omega$ will be said to be {\it {\rm 1}-$\{
2\}$-transi\-tive\/} if it acts transitively on the set $ \Omega
^{1,\{ 2\}}=\{ (\alpha , \{ \beta ,\gamma \} )\mid \alpha ,\beta
,\gamma \mbox{ distinct}\} .$ Thus $G$ is 1-$\{ 2\}$-transi\-tive
if and only if it is transitive and a stabilizer $G_{\alpha}$ is
2-homogeneous on $\Omega \setminus \{ \alpha \} $. The problem is
to classify all (infinite) permutation groups that are 1-$\{
2\}$-transi\-tive but not 3-transitive. \hfill \raisebox{-1ex}{\sl
P.\,M.\,Neumann}

\emp

 \bmp \textbf{11.70.} Let $F$ be an infinite
field or a skew-field.

\makebox[24pt][r]{a)} Find all transitive subgroups of $ PGL(2,F)$
acting on the projective line $F \cup \{ \infty \} $.

\makebox[24pt][r]{c)} What are the flag-transitive subgroups of
$PGL(d + 1,\, F)$?

\makebox[24pt][r]{d)} What subgroups of $PGL(d + 1,\, F)$ are
2-transitive on the points of $PG(d,F)$?

\hfill \raisebox{-1ex}{\sl P.\,M.\,Neumann, C.\,E.\,Praeger}
\emp

\bmp \textbf{11.71.} Let A be a finite group
with a normal subgroup $H$. A subgroup $U$ of $H$ is called an {\it
$A$-covering subgroup\/} of $H$ if $\bigcup\limits_{a\in A} U^a=H$.
Is there a function $f\!: {\Bbb N} \rightarrow {\Bbb N}$ such that
whenever
 $U< H < A$, where $A$ is a finite group, $H$ is a normal
subgroup of $ A$ of index $n$, and $U$ is an $A$-covering subgroup of
$H$, the index $|H:U| \leq f(n)$? (We have shown that the answer is
``yes'' if $U$ is a maximal subgroup of $H$.)

 \makebox[15pt][r]{}{\it Comments of 2025}:
The conjecture is proved in the case where $H$ acts innately transitively on the coset space $[H:U]$ (M.\,Fusari, A.\,Previtali, P.\,Spiga, {\it  J.~Group Theory}, {\bf 27} (2024), 929--965); \  for $n=3$ with $|H:U|\leq 10$ (L.\,Gogniat, P.\, Spiga, {\it Preprint}, 2025, \url{https://arxiv.org/abs/2502.01287}); and  if $H=UL$ where $L$ is a minimal $A$-invariant subgroup of $H$ (M.\,Fusari, S.\,Harper, P.\,Spiga, {\it Bull. Austral. Math. Soc.}, 2025, \url{https://doi.org/10.1017/S0004972725000176}).

\hfill \raisebox{-1ex}{\sl P.\,M.\,Neumann, C.\,E.\,Praeger}

\emp

 \bmp \textbf{11.72.} Suppose that a variety of
groups ${\goth V}$ is {\it non-regular}, that is, the free group
$F_{n+1}({\goth V})$ is embeddable in $F_n({\goth V})$ for some $
n$. Is it true that then every countable group in ${\goth V}$ is
embeddable in an $n$-gene\-ra\-ted group in ${\goth V}$? \hfill
\raisebox{-1ex}{\sl A.\,Yu.\,Olshanskii}

\emp

 \bmp \textbf{11.73.} If a relatively free
group is finitely presented, is it virtually nilpotent?

\hfill \raisebox{-1ex}{\sl
A.\,Yu.\,Olshanskii}

\emp

 \bmp \textbf{11.76.} (Well-known problem). Is
the group of collineations of a finite non-Desar\-gu\-a\-sian
projective plane defined over a semi-field soluble? (The
hypothesis on the plane means that the corresponding regular set,
 see Archive, 10.48, is closed under addition.)
\hfill \raisebox{-1ex}{\sl N.\,D.\,Podufalov}

\emp

 \bmp \textbf{11.77.}
(Well-known problem). Describe the finite translation planes whose
col\-li\-ne\-a\-tion groups act doubly transitively on the set of
points of the line at infinity.

 \hfill \raisebox{-1ex}{\sl N.\,D.\,Podufalov}

\emp

 \bmp \textbf{11.78.} An isomorphism of groups
of points of algebraic groups is called {\it semialgebraic\/} if
it can be represented as a composition of an isomorphism of
translation of the field of definition and a rational morphism.

\makebox[25pt][r]{a)} Is it true that the existence of an
isomorphism of groups of points of two directly undecomposable
algebraic groups with trivial centres over an algebraically closed
field implies the existence of a semialgebraic isomorphism of the
groups of points?

\makebox[25pt][r]{b)} Is it true that every isomorphism of groups
of points of directly undecomposable algebraic groups with trivial
centres defined over algebraic fields is semialgebraic? A~field is
called {\it algebraic\/} if all its elements are algebraic over
the prime subfield.

\hfill \raisebox{-1ex}{\sl
K.\,N.\,Ponomar\"ev}

\emp

 \bmp \textbf{11.80.} Let $G$ be a primitive
permutation group on a finite set $\Omega$ and suppose that, for
$\alpha \in \Omega $, $G_{\alpha}$ acts 2-transitively on one of
its orbits in $\Omega \setminus \{ \alpha \} $. By
(C.\,E.\,Praeger, {\it J.~Austral. Math. Soc. (A)}, \textbf{45},
1988, 66--77) either

\makebox[15pt][r]{}(a) $T \leq G \leq {\rm Aut}\, T$ for some
nonabelian simple group $T$, or

\makebox[15pt][r]{}(b) $G$ has a unique minimal normal subgroup
which is regular on $\Omega$.

For what classes of simple groups in (a) is a classification
feasible? Describe the examples in as explicit a manner as
possible. Classify all groups in (b).

\makebox[15pt][r]{} {\it Remark of 1999:\/} Two papers by
X.\,G.\,Fang and C.\,E.\,Praeger ({\it Commun. Algebra}, \textbf{27}
(1999), 3727--3754 and
 3755--3769)
 show that this is feasible for $T$ a Suzuki and Ree simple
group, and a paper of J.\,Wang and C.\,E.\,Praeger ({\it
J.~Algebra}, \textbf{180} (1996), 808--833) suggests that this may
not be the case for $T$ an alternating group. {\it Remark of
2001:\/} In (X.\,G.\,Fang, G.\,Havas, J.\,Wang, {\it European
J.~Combin.}, \textbf{20}, no.\,6 (1999), 551--557) new examples are
constructed with $G=PSU(3,q)$. {\it Remark of 2009:\/} It was shown in
(D.\,Leemans, {\it J.\,Algebra}, {\bf 322}, no.\,3 (2009), 882--892)
that for part (a) classification
is feasible for self-paired 2-transitive suborbits with $T$ a
sporadic simple group; in this case these suborbits correspond to
2-arc-transitive actions on undirected graphs.\hfill
\raisebox{-1ex}{\sl C.\,E.\,Praeger}

\emp

 \bmp \textbf{11.81.} A topological group is
said to be {\it $F$-balanced\/} if for any subset $X$ and any
neighborhood of the identity $U$ there is a neighborhood of the
identity $V$ such that $VX \subseteq XU$. Is every $ F$-balanced group
{\it balanced}, that is, does it have a basis of neighborhoods of
the identity consisting of invariant sets? \hfill
\raisebox{-1ex}{\sl I.\,V.\,Protasov}

\emp

 \bmp \textbf{11.83.} Is the conjugacy problem
soluble for finitely generated abelian-by-polycyclic groups?
\hfill \raisebox{-1ex}{\sl V.\,N.\,Remeslennikov}

\emp

 \bmp \textbf{11.84.} Is the isomorphism
problem soluble

\makebox[25pt][r]{a)} for finitely generated metabelian groups?

\makebox[25pt][r]{b)} for finitely generated soluble groups of
finite rank?
\hfill \raisebox{0ex}{\sl
V.\,N.\,Remeslennikov}

\emp

 \bmp \textbf{11.85.} Let $F$ be a free
pro-$p\hskip0.2ex$-group with a basis $X$ and let $ R=r^F$ be the
closed normal subgroup generated by an element $ r$. We say that
an element $s$ of $F$ is {\it asso\-ci\-at\-ed to $r$\/} if
$s^F=R$.

\makebox[25pt][r]{a)} Suppose that none of the elements
associated to $r$ is a $p$th power of an element in
$F$. Is it true that $F/R$ is torsion-free?

\makebox[25pt][r]{b)} Among the elements associated to $r$ there
is one that depends on the minimal subset $X'$ of the basis $X$.
Let $x \in X'$. Is it true that the images of the elements of
$X\setminus \{ x\}$ in $F/R$ freely generate a free
pro-$p\hskip0.2ex$-group? \hfill \raisebox{-1ex}{\sl
N.\,S.\,Romanovski\u{\i}}

\emp

 \bmp \textbf{11.86.} Does every group
$G=\left< x_1,\ldots ,x_n\mid r_1=\cdots =r_m=1\right> $ possess,
in a natural way, a homomorphic image $H=\left< x_1,\ldots
,x_n\mid s_1=\cdots =s_m=1\right> $

\makebox[25pt][r]{a)} such that $H$ is a torsion-free group?

\makebox[25pt][r]{b)} such that the integral group ring of $H$ is embeddable
in a skew field?

\makebox[15pt][r]{}If the stronger assertion ``b'' is true, then
this will give an explicit method of finding elements
$x_{i_1},\ldots ,x_{i_{n-m}}$ which generate a free group in $G$.
Such elements exist by N.\,S.\,Romanovski\u{\i}'s theorem ({\it
Algebra and Logic}, \textbf{16}, no.\,1 (1977), 62--67).

\hfill
\raisebox{-1ex}{\sl V.\,A.\,Roman'kov}

\emp

 \bmp \textbf{11.87.} (Well-known problem). Is
the automorphism group of a free met\-abe\-li\-an group of rank $n
\geq 4$ finitely presented? \hfill
\raisebox{-1ex}{\sl V.\,A.\,Roman'kov}

\emp

 \bmp \textbf{11.89.} Let $k$ be an infinite
cardinal number. Describe the epimorphic images of the Cartesian
power $\prod\limits_{k}{\Bbb Z}$ of the group ${\Bbb Z}$ of
integers. Such a description is known for $k=\omega _0$
(R.\,Nunke, {\it Acta Sci. Math. (Szeged)}, \textbf{23} (1963),
67--73). \hfill {\sl S.\,V.\,Rychkov}

\emp

 \bmp \textbf{11.90.} Let ${\goth V}$ be a
variety of groups and let $k$ be an infinite cardinal number.
A~group $G \in {\goth V}$ of rank $k$ is called {\it almost free
in
${\goth V}$\/} if each of its subgroups of rank less than $k$ is
contained in a subgroup of $G$ which is free in ${\goth V}$. For
what~$k$ do there exist almost free but not free in ${\goth B}$
groups of rank $k$? \hfill \raisebox{-1ex}{\sl S.\,V.\,Rychkov}

\emp

 \bmp \textbf{11.92.} What are the soluble
hereditary non-one-generator formations of finite groups all of
whose proper hereditary subformations are one-generated?
\hfill \raisebox{-1ex}{\sl A.\,N.\,Skiba} \emp

 \bmp \textbf{11.95.} Suppose that $G$ is a $p\hskip0.2ex$-group
$G$ containing an element $a$ of order $p$
such that the subgroup $\left< a,\, a^g\right> $ is finite for any
$g$ and the set $C_G(a) \cap a^G$ is finite. Is it true that $G$
 has non-trivial
centre? This is true for 2-groups. \hfill
\raisebox{-1ex}{\sl S.\,P.\,Strunkov}
\emp

 \bmp \textbf{\zv 11.96.}
 (a) Is it true that, for a
given number $n$, there exist only finitely many finite simple
groups each of which contains an involution which commutes with at
most $n$ involutions of the group?

 \makebox[15pt][r]{}(b) Is it true that there are no
infinite simple groups satisfying this condition?

\hfill
\raisebox{-1ex}{\sl S.\,P.\,Strunkov}

\ul

\otv (a) Yes, it is true (mod CFSG); moreover, if a locally finite group has such an involution, then it is finite of $n$-bounded order (S.\,V.\,Skresanov, {\it J.~Group Theory}, {\bf 27}, no.\,6 (2024), 1197--1202).

\otv (b) No, since $PSL_2(\Bbb{R})$ is an infinite simple group, which contains an involution but does not contain a subgroup isomorphic to the Klein 4-group (S.\,V.\,Skresanov, {\it J.~Group Theory}, {\bf 27}, no.\,6 (2024), 1197--1202).
\emp

 \bmp \textbf{11.98.} a) (R.\,Brauer). Find the
best-possible estimate of the form $|G| \leq f(r)$ where $r$ is the
number of conjugacy classes of elements in a finite (simple)
group~$G$.

\hfill \raisebox{-1ex}{\sl S.\,P.\,Strunkov}

\emp

 \bmp \textbf{11.99.} Find (in group-theoretic
terms) necessary and sufficient conditions for a finite group to
have complex irreducible characters having defect 0 for more than
one prime number dividing the order of the group. Express the
number of such characters in the same terms. \hfill
\raisebox{-1ex}{\sl S.\,P.\,Strunkov}

\emp

 \bmp \textbf{11.100.} Is it true that every
periodic conjugacy biprimitively finite group (see 6.59) can be
obtained from
2-groups and binary finite groups by taking extensions? This
question is of independent interest for $p\hs$-groups, $p$ an odd
prime. \hfill \raisebox{-1ex}{\sl S.\,P.\,Strunkov}

\emp

 \bmp \textbf{11.102.} Does there exist a
residually soluble but insoluble group satisfying the maximum
condition on subgroups? \hfill \raisebox{-1ex}{\sl J.\,Wiegold}

\emp

 \bmp \textbf{11.105.} b) Let ${\goth V}$ be a
variety of groups. Its relatively free group of given rank has a
presentation $F/N$, where $F$ is absolutely free of the same rank
and $N$ fully invariant in $F$. The associated Lie ring ${\cal
L}(F/N)$ has a presentation ${ L}/J$, where ${ L}$ is the free Lie
ring of the same rank and $J$ an ideal of ${ L}$.
 Is $J$ fully invariant in ${ L}$ if ${\goth V}$ is the
Burnside variety of all groups of given exponent $q$, where $q$ is
a prime-power, $q \geq 4$?

\hfill \raisebox{-1ex}{\sl G.\,E.\,Wall} \emp

 \bmp \textbf{11.107.} Is there a non-linear
locally finite simple group each of whose proper sub\-groups is
residually finite? \hfill \raisebox{-1ex}{\sl R.\,Phillips}

\emp

 \bmp \textbf{11.109.} Is it true that the
union of an ascending series of groups of Lie type of equal ranks
over some fields should be also a group of Lie type of the same
rank over a field? (For locally finite fields the answer is
``yes'', the theorem in this case is due to \,V.\,Belyaev,
\,A.\,Borovik, \,B.\,Hartley\,\,\&\,\,G.\,Shute, and \,S.\,Thomas.)
\hfill \raisebox{-1ex}{\sl R.\,Phillips}

\emp

 \bmp \textbf{11.111.} Does every infinite
locally finite simple group $G$ contain a nonabelian finite simple
subgroup? Is there an infinite tower of such subgroups in~$G$?
\hfill
\raisebox{-1ex}{\sl B.\,Hartley}

\emp

 \bmp \textbf{11.112.} Let $L=L(K(p))$ be the
associated Lie ring of a free countably generated group $K(p)$ of
the Kostrikin variety of locally finite groups of a given prime
exponent $p$. Is it true that

\makebox[25pt][r]{a)} $L$ is a relatively free Lie ring?

 \makebox[25pt][r]{c)} all identities of $L$ follow from a finite
number of identities of $L$?

\hfill \raisebox{-1ex}{\sl E.\,I.\,Khukhro}

\emp

 \bmp \textbf{11.113.} (B.\,Hartley). Is it
true that the derived length of a nilpotent pe\-ri\-o\-dic group
admitting a regular automorphism of prime-power order $p^n$ is
bounded by a function of $p$ and $n$? \hfill
\raisebox{-1ex}{\sl E.\,I.\,Khukhro}

\emp

 \bmp \textbf{11.114.} Is every locally graded
group of finite special rank almost hyperabelian? This is true in
the class of periodic locally graded groups ({\it Ukrain.
Math.~J.}, \textbf{42}, no.\,7 (1990), 855--861). \hfill
\raisebox{-1ex}{\sl N.\,S.\,Chernikov}

\emp

 \bmp \textbf{\ai 11.115.} 
 Suppose that $X$ is a
free non-cyclic group, $N$ is a non-trivial normal subgroup of $X$
and $T$ is a proper subgroup of $ X$ containing $N$. Is it true
that $[T,N] < [X,N]$ if $N$ is not maximal in $X$? \hfill
\raisebox{-1ex}{\sl V.\,P.\,Shaptala}

\ul

\oai No, not always (unverified AI solution:  M.\,Kinyon, J.\,Urban, {\it Preprint}, 2026,

  \url{https://kourovkanotebookorg.wordpress.com/wp-content/uploads/2026/09/kourovka-mathematics.pdf}).

\emp 

 \bmp \textbf{\zv 11.116.}
 The {\it dimension\/} of a
partially ordered set $\left< P, \, \leq \right> $ is, by
definition, the least cardinal number $\delta$ such that the
relation $\leq$ is an intersection of $\delta$ relations of linear
order on $P$. Is it true that, for any Chernikov group which does
not contain a direct product of two quasicyclic groups over the
same prime number, the subgroup lattice has finite dimension? {\it
Expected answer:\/} yes. \hfill \raisebox{-1ex}{\sl
L.\,N.\,Shevrin}

\ul

\otv Yes, it is true (AI solution confirmed by P.\,Shumyatsky: \  M.\,Kinyon, J.\,Urban, {\it Preprint}, 2026,
  \url{https://kourovkanotebookorg.wordpress.com/wp-content/uploads/2026/09/kourovka-mathematics.pdf}).

\emp 

 \bmp \textbf{11.117.} Let ${\goth X}$ be a
soluble non-empty Fitting class. Is it true that every finite
non-soluble group possesses an ${\goth X}$-injec\-tor? \hfill
\raisebox{-1ex}{\sl L.\,A.\,Shemetkov}

\emp

 \bmp \textbf{11.118.} What are the hereditary
soluble local formations ${\goth F}$ of finite groups such that every
finite group has an ${\goth F}$-covering subgroup? \hfill
\raisebox{-1ex}{\sl L.\,A.\,Shemetkov}

\emp

 \bmp \textbf{11.119.} Is it true that, for any
non-empty set of primes $\pi$, the $\pi $-length of any finite
$\pi $-soluble group does not exceed the derived length of its
Hall $\pi $-sub\-group?

\hfill \raisebox{-1ex}{\sl L.\,A.\,Shemetkov}

\emp

 \bmp \textbf{11.120.} Let ${\goth F}$ be a
soluble saturated Fitting formation. Is it true that $l_{\goth F}(G)
\leq f(c_{\goth F}(G))$, where $c_{\goth F}(G)$ is the length of a
composition series of some ${\goth F}$-cover\-ing subgroup of
a~finite soluble group $G$? This is Problem 17 in
(L.\,A.\,Shemetkov, {\it Formations of Finite Groups}, Moscow,
Nauka, 1978 (Russian)). For the definition of the {\it ${\goth
F}$-length\/} $l_{{\goth F}}(G)$, see {\it ibid}.\hfill
\raisebox{-1ex}{\sl L.\,A.\,Shemetkov}

\emp

 \bmp \textbf{11.121.} Does there exist a local
formation of finite groups which has a non-trivial decomposition
into a product of two formations and in any such a decomposition
both factors are non-local? Local products of non-local formations
do exist.

\hfill
\raisebox{-1ex}{\sl L.\,A.\,Shemetkov}

\emp

 \bmp \textbf{11.122.} Does every non-zero
submodule of a free module over the group ring of a torsion-free
group contain a free cyclic submodule? \hfill
\raisebox{-1ex}{\sl A.\,L.\,Shmel'kin}

\emp

 \bmp \textbf{11.123.} (Well-known problem).
For a given group $G$, define the following sequence of groups:
$A_1(G)=G$, $A_{i+1}(G)={\rm Aut}\, (A_i(G))$. Does there exist a
finite group $G$ for which this sequence contains infinitely many
non-isomorphic groups? \hfill \raisebox{-1ex}{\sl M.\,Short}

\emp

 \bmp \textbf{\zv 11.124.}
 Let $F$ be a non-cyclic
free group and
 $R$ a non-cyclic subgroup of $F$. Is it true that if $[R,R]$
is a normal subgroup of $F$ then $R$ is also a normal subgroup of
$F$?

\hfill \raisebox{-1ex}{\sl V.\,E.\,Shpilrain}

\ul

\otv Yes, it is true (A.\,Tursunbaev, {\it Preprint}, 2026,
\url{https://arxiv.org/abs/2609.06755}).
\emp

 \bmp \textbf{11.125.} Let $G$ be a finite
group admitting a regular elementary abelian group of automorphisms
$V$ of order $p^n$. Is it true that the subgroup
$H=\bigcap\limits_{v\in V\setminus \{ 1\}} [G,v]$ is nilpotent? In
the case of an affirmative answer, does there exist a function
depending only on $p$, $n$, and the derived length of $G$ which
bounds the nilpotency class of $H$?

\hfill
\raisebox{-1ex}{\sl
P.\,V.\,Shumyatski\u{\i}}

\emp

 \bmp \textbf{11.127.} Is every group of
exponent 12 locally finite? \hfill
\raisebox{-1ex}{\sl V.\,P.\,Shunkov}

\emp

 \raggedbottom

\newpage

\pagestyle{myheadings} \markboth{12th Issue
(1992)}{12th Issue (1992)}
\thispagestyle{headings}

~
\vspace{2ex}

\centerline {\Large \textbf{Problems from the 12th Issue (1992)}}
\phantomsection\label{12izd}
\vspace{4ex}

\bmp \textbf{12.1.} a) H.\,Bass ({\it Topology}, \textbf{4}, no.\,4
(1966), 391--400) has constructed explicitly a proper subgroup of
finite index in the group of units of the integer group ring of a
finite cyclic group. Calculate the index of Bass' subgroup.\hfill
\raisebox{-1ex}{\sl R.\,Zh.\,Aleev} \emp

\bmp \textbf{12.3.} A $p\hskip0.2ex $-group is called {\it thin\/} if
every set of pairwise incomparable (by inclusion) normal subgroups
contains $\leq p + 1$ elements. Is the number of thin
pro-$p\hskip0.2ex $-groups finite? \hfill \raisebox{-1ex}{\sl
R.\,Brandl}

\emp

\bmp \textbf{12.4.} Let $G$ be a group and assume that we have
$[x,y]^3=1$ for all $x,y \in G$. Is $G'$ of finite exponent? Is
$G$ soluble? By a result of N.\,D.\,Gupta and N.\,S.\,Mendelssohn,
1967, we know that $G'$ is a 3-group. \hfill \raisebox{-1ex}{\sl
R.\,Brandl}

\emp

\bmp \textbf{12.6.} Let $G$ be a finitely generated group and suppose
that $H$ is a $p\hskip0.2ex $-sub\-group of $G$ such that $H$
contains no non-trivial normal subgroups of $G$ and $HX=XH$ for
any subgroup $ X$ of $G$. Is then $G/C_G(H^G)$ a $p\hskip0.2ex
$-group where $H^G$ is the normal closure of~$H$? \hfill
\raisebox{-1ex}{\sl G.\,Busetto}

\emp

\bmp \textbf{12.8.} Let ${\cal V}$ be a non-trivial variety of groups
and let $a_1,\ldots ,a_r$ freely generate a free group $F_r(V)$ in
${\cal V}$. We say $F_r({\cal V})$ {\it strongly discriminates\/}
${\cal V}$ just in case every finite system of inequalities
$w_i(a_1,\ldots ,a_r,x_1,\ldots ,x_k) \ne 1$ for $1 \leq i \leq n$
having a solution in some $F_s({\cal V})$ containing $F_r({\cal
V})$ as a varietally free factor in the sense of ${\cal V}$,
already has a solution in $F_r({\cal V})$. Does there exist ${\cal
V}$ such that for some integer $r > 0$, $F_r({\cal V})$
discriminates but does not strongly discriminate~${\cal V}$? What
about the analogous question for general algebras in the context
of universal algebra?

\hfill
\raisebox{-1ex}{\sl A.\,M.\,Gaglione,
 D.\,Spellman}

\emp

\bmp \textbf{12.9.} Following Bass, call a group {\it tree-free\/} if
there is an ordered Abelian group $\Lambda$ and a $\Lambda $-tree
$X$ such that $G$ acts freely without inversion on $X$.

 \makebox[25pt][r]{a)} Must every finitely generated tree-free
group satisfy the
maximal condition for Abelian subgroups?

 \makebox[25pt][r]{b)} The same question for finitely
presented tree-free
groups.

\hfill \raisebox{-1ex}{\sl A.\,M.\,Gaglione,
 D.\,Spellman}

\emp

\bmp \textbf{12.11.} Suppose that $G, H$ are countable or finite
groups and $A$ is a proper subgroup of $G$ and $H$ containing no
non-trivial subgroup normal in both $G$ and $H$. Can $G*_A H$ be
embedded in ${\rm Sym}\, ({\Bbb N})$ so that the image is highly
tran\-sitive?

\makebox[15pt][r]{}{\it Comment of 2013:}
for partial results, see (S.\,V.\,Gunhouse, {\it Highly
transitive representations of free products on the natural numbers}, Ph.\,D. Thesis
Bowling Green State Univ., 1993; \ K.\,K.\,Hickin, {\it J. London Math. Soc.
(2)}, {\bf 46}, no.\,1 (1992), 81--91).

\hfill \raisebox{-1ex}{\sl A.\,M.\,W.\,Glass}

\emp

\bmp \textbf{12.12.} Is the conjugacy problem for nilpotent finitely
generated lattice-ordered groups soluble? \hfill
\raisebox{-1ex}{\sl A.\,M.\,W.\,Glass}

\emp

\bmp \textbf{12.13.} If $\left< \Omega , \leq \right>$ is the
countable universal poset, then $G={\rm Aut}\, (\left< \Omega ,
\leq \right> )$ is simple (A.\,M.\,W.\,Glass, S.\,H.\,McCleary,
M.\,Rubin, {\it Math.~Z.}, \textbf{214}, no.\,1 (1993), 55--66). If
$H$ is a subgroup of $G$ such that $|G:H| < 2^{\aleph _0}$ and $H$
is transitive on $\Omega$, does $H=G$?

\hfill \raisebox{-1ex}{\sl A.\,M.\,W.\,Glass}

\emp

\bmp \textbf{\zv 12.15.}
Suppose that, in a finite 2-group $G$, any two
elements are conjugate whenever their normal closures coincide. Is
it true that the derived subgroup of $G$ is abelian?

\hfill
\raisebox{-1ex}{\sl E.\,A.\,Golikova, A.\,I.\,Starostin}

\ul

\otv No, not always (K.\,Muliarchyk,  {\it Preprint}, 2026,

\url{https://kourovkanotebookorg.wordpress.com/wp-content/uploads/2026/08/kourovka_12_15_counterexample.pdf}).
\emp

\bmp \textbf{12.16.} Is the class of groups of recursive
automorphisms of arbitrary models closed with respect to taking
free products? \hfill \raisebox{-1ex}{\sl S.\,S.\,Goncharov}

\emp

\bmp \textbf{12.17.} Find a description of autostable periodic
abelian groups.
\hfill \raisebox{0ex}{\sl S.\,S.\,Goncharov}

\emp

\bmp \textbf{12.19.} Is it true that, for every $n \geq 2$ and every
two epimorphisms $\varphi$ and $\psi$ of a free group $F_{2n}$ of
rank $2n$ onto $F_n \times F_n$, there exists an automorphism
$\alpha$ of $F_{2n}$ such that $\alpha \varphi =\psi $? \hfill
\raisebox{-1ex}{\sl R.\,I.\,Grigorchuk}

\emp \vspace{-1ex}

\bmp \textbf{\zv 12.20.}
 (Well-known problem). Is
R.\,Thompson's group \vspace{-1ex} $$F=\left< x_0,\, x_1,\ldots \mid
x_n^{x_i}= x_{n+1}, \; \; i < n,\;\;n=1, \, 2, \ldots \right>
=$$\vskip-3ex

 $$= \left< x_0,\ldots ,x_4\mid x_1^{x_0}=x_2, \; \;
x_2^{x_0}=x_3, \; \; x_2^{x_1}=x_3, \; \; x_3^{x_1}=x_4, \; \;
x_3^{x_2}=x_4\right> $$\vskip-1ex amenable? \hfill
\raisebox{-1ex}{\sl R.\,I.\,Grigorchuk}

\ul

\otv No, it is not amenable (AI solution confirmed by S.-H.\,Kim and Y.\,Lodha and by M.\,de~la~Salle: \  OpenAI, {\it Preprint}, 2026,
\url{https://kourovkanotebookorg.wordpress.com/wp-content/uploads/2026/10/12-20.pdf}; \ S.-H.\,Kim and Y.\,Lodha, {\it Preprint}, 2026, \url{https://arxiv.org/abs/2610.18360}; \  M.~de~la~Salle,  {\it Preprint}, 2026, \url{http://math.univ-lyon1.fr/homes-www/delasalle/thompson.pdf}).
\emp 

\bmp \textbf{12.21.} Let $R$ be a commutative ring with identity and
let $G$ be a finite group. Prove that the ring $a(RG)$ of $R
$-repre\-sen\-ta\-tions of $G$ has no non-trivial idempotents.

\hfill \raisebox{-1ex}{\sl P.\,M.\,Gudivok,
V.\,P.\,Rud'ko}

\emp

\bmp \textbf{12.23.} The {\it index of permutability\/} of a group
$G$ is defined to be the minimal integer $k \geq 2$ such that for
each $k $-tuple $x_1,\ldots ,x_k$ of elements in $G$ there is a
non-identity permutation $\sigma$ on $k$ symbols such that
$x_1\cdots x_k=x_{\sigma (1)}\cdots x_{\sigma (k)}$. Determine the
index of permutability of the symmetric group ${\Bbb S}_n$. \hfill
\raisebox{-1ex}{\sl M.\,Gutsan}

\emp

\bmp \textbf{12.27.} Let $G$ be a simple locally finite group. We say
that $G$ is a {\it group of finite type\/} if there is a
non-trivial permutational representation of $G$ such that some
finite subgroup of $G$ has no regular orbits. Investigate and,
perhaps, classify
the simple locally finite groups of finite type. Partial results
see in (B.\,Hartley, A.\,Zalesski, {\it Isr. J. Math.}, \textbf{82}
(1993), 299--327, {\it J. London Math. Soc.}, \textbf{55} (1997),
210--230; \ F.\,Leinen, O.\,Puglisi,
{\it Illinois J. Math.}, {\bf
47} (2003), 345--360). \hfill \raisebox{-1ex}{\sl
A.\,E.\,Zalesski\u{\i}}

\emp

\bmp \textbf{12.28.} Let $G$ be a group. A function $f\! : G
\rightarrow {\Bbb C}$ is called\vspace{-1.5ex}
\begin{enumerate}

\item[1)] {\it normed\/} if $f(1)=1$;\vspace{-1.5ex}

\item[2)] {\it central\/} if $f(gh)=f(hg)$ for all
$g,h \in G$;\vspace{-1.5ex}

\item[3)] {\it positive-definite\/} if $\sum\limits_{k,\,
l}f(g_k^{-1}g_l)\,\overline{\! c}_kc_l \geq 0$ for any $g_1,\ldots
,g_n \in G$ and any\linebreak \vskip-3ex $c_1,\ldots ,c_n \in
{\Bbb C}$. \end{enumerate}\vspace{-1.5ex}

Classify the infinite simple locally finite groups $G$ which
possess functions satisfying 1)--3). The simple Chevalley groups
are known to have no such functions, while such functions exist on
locally matrix (or stable) classical groups over finite fields.
The question is motivated by the theory of $C^{\displaystyle{*}}
$-algebras, see $\S$\,9 in (A.\,M.\,Vershik, S.\,V.\,Kerov, {\it
J.~Sov. Math.}, \textbf{38} (1987), 1701--1733). Partial results see
in (F.\,Leinen, O.\,Puglisi,
{\it J.~Pure Appl. Algebra}, {\bf
208} (2007), 1003--1021, {\it J.~London Math. Soc.}, \textbf{70}
(2004), 678--690). \hfill \raisebox{-1ex}{\sl
A.\,E.\,Zalesski\u{\i}} \emp

\bmp \textbf{12.29.} Classify the locally finite groups for which the
augmentation ideal of the complex group algebra is a simple ring.
The problem goes back to I.\,Kaplansky (1965). For partial
results, see (K.\,Bonvallet, B.\,Hartley, D.\,S.\,Passman,
M.\,K.\,Smith, {\it Proc. Amer. Math. Soc.}, \textbf{56}, no.\,1
(1976), 79--82; \ A.\,E.\,Zalesski\u{\i}, {\it Algebra i Analiz},
\textbf{2}, no.\,6 (1990), 132--149 (Russian); \ Ch.\,Praeger,
A.\,E.\,Zalesskii, {\it Proc. London Math. Soc.}, \textbf{70}, no.\,2
(1995), 313--335;
\ B.\,Hartley, A.\,Zalesski, {\it J. London Math. Soc.}, \textbf{55}
(1997), 210--230). \hfill \raisebox{-1ex}{\sl
A.\,E.\,Zalesski\u{\i}}

\emp

\bmp \textbf{12.30.} (O.\,N.\,Golovin). On the class of all groups,
do there exist associative operations which satisfy the postulates
of MacLane and Mal'cev (that is, which are free functorial and
hereditary) and which are different from taking free and direct
products? \hfill \raisebox{-1ex}{\sl S.\,V.\,Ivanov}

\emp

\bmp \textbf{12.33.} Suppose that $G$ is a finite group and $x$ is an
element of $G$ such that the subgroup $\left< x,y\right>$ has odd
order for any $y$ conjugate to $x$ in $G$. Prove, without using
CFSG, that the normal closure of $x$ in $G$ is a group of odd
order. \hfill \raisebox{-1ex}{\sl L.\,S.\,Kazarin}

\emp

\bmp \textbf{12.34.} Describe the finite groups $G$ such that the sum
of the cubes of the degrees of all irreducible complex characters
is at most $|G|\cdot \log _2|G|$. The question is interesting for
applications in the theory of signal processing. \hfill
\raisebox{-1ex}{\sl L.\,S.\,Kazarin}

\emp

\bmp \textbf{12.35.} Suppose that ${\frak F}$ is a radical composition
formation of finite groups. Prove that $\left< H,K\right> ^{{\frak
F}} =\left< H^{{\frak F}},\, K^{{\frak F}}\right>$ for every finite
group $G$ and any subnormal subgroups $H$ and $K$ of~$G$. \hfill
\raisebox{-1ex}{\sl S.\,F.\,Kamornikov}

\emp

\bmp \textbf{\zv 12.37.}
(J.\,G.\,Thompson). For a finite group $G$ and
natural number $n$, set $G(n)=\{ x \in G\mid x^n=1\} $ and define
the {\it type\/} of $G$ to be the function whose value at $n$ is
the order of $G(n)$. Is it true that a group is soluble if its
type is the same as that of a soluble one? \hfill
\raisebox{-1ex}{\sl A.\,S.\,Kondratiev, W.\,J.\,Shi}

\ul

\otv No, not always (P.\,Piwek, {\it Algebra Number Theory}, {\bf 19}, no.\,8 (2025), 1663--1670).
\emp

\bmp \textbf{12.40.}
Let $\varphi$ be an irreducible $p\hskip0.2ex
$-modular character of a finite group $G$. Find the best-possible
estimate of the form $\varphi (1)_p \leq f(|G|_p)$, where $n_p$ is
the $p\hskip0.2ex $-part of a positive integer~$n$,

\makebox[25pt][r]{\zv a)} for $p=2$;

\makebox[25pt][r]{b)} for $p\ne 2$.
 \hfill
\raisebox{0ex}{\sl A.\,S.\,Kondratiev}

\ul

\otv a) There cannot be any estimate of this kind for $p=2$ (AI solution confirmed by A.\,S.\,Kondratiev: \  M.\,Kinyon, J.\,Urban, {\it Preprint}, 2026,

\url{https://kourovkanotebookorg.wordpress.com/wp-content/uploads/2026/09/kourovka-mathematics.pdf}).
\emp 

\bmp \textbf{\zv 12.41.}
Let $F$ be a free group on two generators $x,y$
and let $\varphi$ be the automorphism of $F$ defined by $x
\rightarrow y$, $y \rightarrow xy$. Let $G$ be a semidirect
product of $F/(F''(F')^2)$ and~$\left< \varphi \right>$. Then $G$
is just-non-polycyclic. What is the cohomological dimension of $G$
over~${\Bbb Q}$? (It is either 3 or 4.) \hfill \raisebox{-1ex}{\sl
P.\,H.\,Kropholler}

\ul

\otv It is 4, because $G$ is not constructible, and then its cohomological dimension is equal to the homological dimension plus~1 by Theorem~III.8 in (M.\,R.\,Bridson, P.\,H.\,Kropholler, {\it J.~Reine Angew. Math.}, {\bf 699} (2015), 217--243), while the homological dimension is equal to the Hirsch length by (U.\,Stammbach, {\it J.~London Math. Soc. (2)}, {\bf 2} (1970), 567--570).
\emp

\bmp \textbf{12.43.} (Well-known problem). Does there exist an
infinite finitely-generated residually finite $p\hskip0.2ex
$-group such that each subgroup is either finite or of finite
index?

\hfill
\raisebox{-1ex}{\sl J.\,C.\,Lennox}

\emp

\bmp \textbf{12.48.} Let $G$ be a sharply doubly transitive
permutation group on a set $\Omega$ (see Archive 11.52 for a definition).

 \makebox[25pt][r]{(a)} Does $G$ possess a regular normal
subgroup if a point
stabilizer is locally finite?

 \makebox[25pt][r]{(b)} Does $G$ possess a regular normal
subgroup if a point
stabilizer has an abelian subgroup of finite index?

 \makebox[15pt]{}{\it Comments of 2022}: an affirmative answer to part (b) is obtained in permutation
characteristic~0 (F.\,O.\,Wagner, {\it Preprint}, 2022, \url{https://hal.archives-ouvertes.fr/hal-03590818}).\hfill
\raisebox{-1ex}{\sl V.\,D.\,Mazurov}

\emp

\bmp \textbf{12.50.} (Well-known problem). Find an algorithm which
decides, by a given finite set of matrices in $SL_3({\Bbb Z})$,
whether the first matrix of this set is contained in the subgroup
generated by the remaining matrices. An analogous problem for
$SL_4({\Bbb Z})$ is insoluble since a direct product of two free
groups of rank 2 embeds into $SL_4({\Bbb Z})$.

\hfill \raisebox{-1ex}{\sl G.\,S.\,Makanin}

\emp

\bmp \textbf{\zv 12.51.}
Does the free group $F_{\eta}$ ($1 < \eta <
\infty$) have a finite subset $S$ for which there is a unique
total order of $F_{\eta}$ making all elements of $S$ positive?
Equivalently, does the free representable $l $-group of rank
$\eta$ have a basic element? (The analogs for right orders of
$F_{\eta}$ and free $l $-groups have negative answers.) \hfill
\raisebox{-1ex}{\sl S.\,McCleary}

\ul

\otv No, it does not (S.\,Dovhyi, K.\,Muliarchyk, {\it Groups Geom. Dyn.}, {\bf 17} (2023), 613--632; \ K.\,Muliarchyk, {\it Preprint}, 2024, \url{https://arxiv.org/abs/2403.14779}).
\emp

\bmp \textbf{12.52.} Is the free $l $-group ${\cal F}_{\eta}$ of rank
$\eta $ ($1 < \eta < \infty$) Hopfian? That is, are $l
$-homo\-mor\-phisms from ${\cal F}_{\eta}$ onto itself necessarily
one-to-one? \hfill \raisebox{-1ex}{\sl S.\,McCleary}

\emp

\bmp \textbf{12.53.} Is it decidable whether or not two elements of a
free $l $-group are conjugate?

\hfill
\raisebox{-1ex}{\sl S.\,McCleary}

\emp

\bmp \textbf{12.54.} Is there a normal valued $l $-group $G$ for
which there is no abelian $l $-group $A$ with $C(A) \cong C(G)$?
(Here $C(G)$ denotes the lattice of convex $l $-sub\-groups of
$G$. If $G$ is not required to be normal valued, this question has
an affirmative answer.)

\hfill
\raisebox{-1ex}{\sl S.\,McCleary}

\emp

\bmp \textbf{12.55.} Let $f(n,p)$ be the number of groups of order
$p^n$. Is $f(n,p)$ an increasing function of $p$ for any fixed $n
\geq 5$? \hfill \raisebox{-1ex}{\sl A.\,Mann}

\emp

\bmp \textbf{12.56.} Let $\left< X\mid R\right>$ be a finite
presentation (the words in $R$ are assumed cyclically reduced).
Define the {\it length\/} of the presentation to be the sum of the
number of generators and the lengths of the relators. Let $F(n)$
be the number of (isomorphism types of) groups that have a
presentation of length at most~$n$. What can one say about the
function $F(n)$? It can be shown that $F(n)$ is not recursive
(D.\,Segal) and it is at least exponential (L.\,Pyber). \hfill
\raisebox{-1ex}{\sl A.\,Mann}

\emp

\bmp \textbf{12.57.} Suppose that a cyclic group $A$ of order 4 is a
$TI $-sub\-group of a finite group~$G$.

 \makebox[25pt][r]{a)} Does $A$ centralize a component $L$
of $G$ if $A$
intersects $L\cdot C(L)$ trivially?

 \makebox[25pt][r]{b)} What is the structure of the normal
closure of $A$ in the
case where $A$ centralizes each component of $G$? \hfill
\raisebox{-1ex}{\sl A.\,A.\,Makhn\"ev}

\emp

\bmp \textbf{12.58.} A {\it generalized quadrangle\/} $GQ(s,t)$ with
parameters $s, t$ is by definition an incidence system consisting of
points and lines in which every line consists of $s + 1$ points,
every two different lines have at most one common point, each point
belongs to $t + 1$ lines, and for any point $a$ not lying on a line
$L$ there is a unique line containing $a$ and intersecting $L$.

 \makebox[25pt][r]{a)} Does $GQ(4,11)$ exist?

 \makebox[25pt][r]{b)} Can the automorphism group of a
hypothetical $GQ(5,7)$
contain involutions?

\hfill \raisebox{-1ex}{\sl A.\,A.\,Makhn\"ev}

\emp

\bmp \textbf{12.59.} Does there exist a strongly regular graph with
parameters (85, 14, 3, 2) and a non-connected neighborhood of some
vertex? \hfill \raisebox{-1ex}{\sl A.\,A.\,Makhn\"ev}

\emp

\bmp \textbf{12.60.} Describe the strongly regular graphs in which
the neighborhoods of vertices are generalized quadrangles (see 12.58).
\hfill \raisebox{-1ex}{\sl A.\,A.\,Makhn\"ev}

\emp

\bmp \textbf{12.62.} A {\it Frobenius group\/} is a transitive
permutation group in which the stabilizer of any two points is
trivial. Does there exist a Frobenius group of infinite degree
which is primitive as permutation group, and in which the
stabilizer of a point is cyclic and has only finitely many orbits?

 \makebox[15pt][r]{}{\it Conjecture:\/} No. Note that such a
group cannot be 2-transitive (Gy.\,K\'aro\-lyi, S.\,J.\,Kov\'acs,
P.\,P.\,P\'alfy, {\it Aequationes Math.}, \textbf{39} (1990),
161--166; \ V.\,D.\,Mazu\-rov, {\it Siberian Math.~J.}, \textbf{31}
(1990), 615--617; \ P.\,M.\,Neumann, P.\,J.\,Rowley, in: {\it
Geometry and cohomology in group theory $($London Math. Soc.
Lect. Note Ser.}, \textbf{252}), Cambridge Univ. Press, 1998,
291--295). \hfill \raisebox{-1ex}{\sl P.\,M.\,Neumann}
\emp

\bmp \textbf{12.63.} Does there exist a soluble permutation group of
infinite degree that has only finitely many orbits on triples?
{\it Conjecture:\/} No.

\makebox[15pt][r]{}{\it Comment of 2001:\/} It is proved
(D.\,Macpherson, {\it in: Advances in algebra and model theory
$($Algebra Logic Appl.}, \textbf{9}), Gordon\,\&\,Breach, Amsterdam,
1997, 87--92) that every infinite soluble permutation group has
infinitely many orbits on quadruples.

 \hfill \raisebox{-1ex}{\sl P.\,M.\,Neumann}

\emp

\bmp \textbf{12.64.} 
Is it true that for a given number $k \geq 2$
and for any (prime) number~$n$, there exists a number $N=N(k,n)$
such that every finite group with generators $A=\{ a_1,\ldots
,a_k\} $ has exponent $\leq n$ if $(x_1\cdots x_N)^n=1$ for any
$x_1,\ldots ,x_N \in A \cup \{ 1\}$?

\makebox[15pt][r]{}For given $k$ and $n$, a negative answer
implies, for example, the infiniteness of the free Burnside group
$B(k,n)$, and a positive answer, in the case of sufficiently large
$n$, gives, for example, an opportunity to find a hyperbolic group
which is not residually finite (and in this group a hyperbolic
subgroup of finite index which has no proper subgroups of finite
index).

\makebox[15pt][r]{}{\it Comment of 2026}: there exist hyperbolic groups that are not residually finite (unverified AI solution: OpenAI, {\it Preprint}, 2026, \url{https://kourovkanotebookorg.wordpress.com/wp-content/uploads/2026/10/partial-12-64_18-88.pdf}).

\hfill \raisebox{-1ex}{\sl A.\,Yu.\,Olshanskii}
\emp 

\bmp \textbf{12.66.} Describe the finite translation planes whose
collineation groups act doubly transitively on the affine points.
\hfill \raisebox{-1ex}{\sl N.\,D.\,Podufalov}

\emp

\bmp \textbf{12.67.} Describe the structure of the locally compact
soluble groups with the maximum (minimum) condition for closed
non-compact subgroups.
\hfill \raisebox{-1ex}{\sl V.\,M.\,Poletskikh}

\emp

\bmp \textbf{12.68.} Describe the locally compact abelian groups in
which any two closed subgroups of finite rank generate a subgroup
of finite rank. \hfill \raisebox{-1ex}{\sl V.\,M.\,Poletskikh}

\emp

\bmp \textbf{12.69.} Let $F$ be a countably infinite field and $G$ a
finite group of automorphisms of~$F$. The group ring ${\Bbb Z}G$
acts on $F$ in a natural way. Suppose that $S$ is a subfield of
$F$ satisfying the following property: for any $x \in F$ there is
a non-zero element $f \in {\Bbb Z}G$ such that $x^f \in S$. Is it
true that either $F=S$ or the field extension $F/S$ is purely
inseparable? One can show that for an uncountable field an
analogous question has an affirmative answer. \hfill
\raisebox{-1ex}{\sl K.\,N.\,Ponomar\"ev}

\emp

\bmp \textbf{12.72.}
 Let ${\frak F}$ be a soluble local hereditary
formation of finite groups. Prove that ${\frak F}$ is radical if
every finite soluble minimal non-${\frak F} $-group $G$ is a
minimal non-${\frak N}^{l(G)-1} $-group. Here ${\frak N}$ is the
formation of all finite nilpotent groups, and $l(G)$ is the
nilpotent length of~$G$. \hfill \raisebox{-1ex}{\sl
V.\,N.\,Semenchuk}
\emp

\bmp \textbf{12.73.} A formation ${\frak H}$ of finite groups is said to
have {\it length\/} $t$ if there is a chain of formations
$\varnothing ={\frak H}_0 \subset {\frak H}_1 \subset \cdots \subset
{\frak H}_t={\frak H}$ in which ${\frak H}_{i-1}$ is a maximal
subformation of ${\frak H}_i$. Is the lattice of soluble formations
of length $\leq 4$ distributive? \hfill \raisebox{-1ex}{\sl
A.\,N.\,Skiba}

\emp

\bmp \textbf{12.75.} (B.\,Jonsson). Is the class $N$ of the lattices
of normal subgroups of groups a~variety? It is proved
(C.\,Herrman,
 W.\,Poguntke, {\it Algebra Univers.},
\textbf{4}, no.\,3 (1974), 280--286) that $N$ is an infinitely based
quasivariety. \hfill \raisebox{-1ex}{\sl D.\,M.\,Smirnov}

\emp

\bmp \textbf{12.81.} What is the cardinality of the set of
subvarieties of the group variety~${\frak A}_p^3$ (where ${\frak
A}_p$ is the variety of abelian groups of prime exponent $p$)?
\hfill \raisebox{-1ex}{\sl V.\,I.\,Sushchanski\u{\i}}

\emp

\bmp \textbf{12.85.} Does every variety which is generated by a (known)
finite simple group containing a soluble subgroup of derived length
$d$ contain a 2-generator soluble group of derived length $d$? \hfill
\raisebox{-1ex}{\sl S.\,A.\,Syskin}

\emp

\bmp \textbf{12.86.} For each known finite simple group, find its
maximal 2-generator direct power.

\hfill \raisebox{-1ex}{\sl
S.\,A.\,Syskin}

\emp

\bmp \textbf{12.87.} Let $\Gamma$ be a connected undirected graph
without loops or multiple edges and suppose that the automorphism
group ${\rm Aut} (\Gamma )$ acts transitively on the vertex set
of~$\Gamma$. Is it true that at least one of the following
assertions holds?

\makebox[15pt][r]{}1.\makebox[0.2cm]{}The stabilizer of a vertex
of $\Gamma$ in
${\rm Aut} (\Gamma )$ is finite.

\makebox[15pt][r]{}2.\makebox[0.2cm]{}The group ${\rm Aut} (\Gamma )$ as a
permutation group on the vertex set of $\Gamma$ admits an
imprimitivity system $\sigma$ with finite blocks for which the
stabilizer of a vertex of the factor-graph $\Gamma /\sigma$ in
${\rm Aut} (\Gamma /\sigma )$ is finite.

\makebox[15pt][r]{}3.\makebox[0.2cm]{}There exists a natural number
$n$ such that the graph obtained from $\Gamma$ by adding edges
joining distinct vertices the distance between which in $\Gamma$ is
at most $n$ contains a regular tree of valency 3. \hfill
\raisebox{-1ex}{\sl V.\,I.\,Trofimov}

\emp

 \bmp
\textbf{12.88.} An undirected graph is called a {\it locally finite
Cayley graph of a group\/}~$G$ if its vertex set can be identified
with the set of elements of $G$ in such a way that, for some
finite generating set $X=X^{-1}$ of $G$ not containing~$1$, two
vertices $g$ and $h$ are adjacent if and only if $g^{-1}h \in X$.
Do there exist two finitely generated groups with the same locally
finite Cayley graph, one of which is periodic and the other is not
periodic? \hfill \raisebox{-1ex}{\sl V.\,I.\,Trofimov}

\emp

\bmp

\textbf{12.89.} Describe the infinite connected graphs to which the
sequences of finite connected graphs with primitive automorphism
groups converge. For definitions, see (V.\,I.\,Trofimov, {\it
Algebra and Logic}, \textbf{28}, no.\,3 (1989), 220--237). \hfill
\raisebox{-1ex}{\sl V.\,I.\,Trofimov}

\emp

\bmp \textbf{12.92.} For a field $K$ of characteristic 2, each finite
group $G=\{ g_1,\ldots ,g_n\}$ of odd order $n$ is determined up to
isomorphism by its group determinant (Formanek--Sibley, 1990) and
even by its {\it reduced norm\/} which is defined as the last
coefficient $s_m(x)$ of the minimal polynomial $\varphi (\lambda
;x)=\lambda ^m - s_1(x)\lambda ^{m-1} + \cdots + (- 1)^{m}s_m(x)$ for
the generic element $x=x_1g_1 + \cdots + x_ng_n$ of the group ring
$KG$ (Hoehnke, 1991). Is it possible in this theorem to replace
$s_m(x)$ by some other coefficients $s_i(x)$, $i < m$? For notation
see (G.\,Frobenius, {\it Sitzungsber. Preuss. Akad. Wiss. Berlin},
1896, 1343--1382) and (K.\,W.\,Johnson, {\it Math. Proc. Cambridge
Phil. Soc.}, \textbf{109} (1991), 299--311).

\hfill
\raisebox{-1ex}{\sl H.-J.\,Hoehnke}

\emp

\bmp \textbf{12.95.} Let $G$ be a finitely generated
pro-$p\hskip0.2ex $-group and let $g_1,\ldots ,g_n \in G$. Let $H$
be an open subgroup of $G$, and suppose there exists a non-trivial
word $w=w(X_1,\ldots ,X_n)$ such that $w(a_1,\ldots ,a_n)=1$
whenever $a_1 \in g_1H,\ldots ,a_n \in g_nH$ (that is, $G$
satisfies a {\it coset identity\/}). Does it follow that $G$
satisfies some non-trivial identity?

\makebox[15pt][r]{}A positive answer to this question would imply
that an analogue of the Tits Alternative holds for finitely
generated pro-$p\hskip0.2ex $-groups. Note that J.\,S.\,Wilson and
E.\,I.\,Zelmanov ({\it J.~Pure Appl. Algebra}, \textbf{81} (1992),
103--109) showed that the graded Lie algebra $L_p(G)$ with respect
to the dimension subgroups of $G$ in characteristic $p$ satisfies
a polynomial identity.

\makebox[15pt][r]{}{\it Comment of 2017}:
This has been shown to be true for finitely generated linear groups (M.\,Larsen, A.\,Shalev, {\it Algebra Number Theory}, {\bf 10}, no.\,6 (2016), 1359--1371).

\hfill \raisebox{-1ex}{\sl A.\,Shalev}

\emp

\bmp \textbf{12.100.} Is every periodic group with a regular
automorphism of order 4 locally finite?

\hfill \raisebox{-1ex}{\sl
P.\,V.\,Shumyatsky}

\emp

\bmp \textbf{12.101.} We call a group $G$ containing an involution $i$ a
{\it $T_0 $-group\/} if

 \makebox[15pt][r]{}1) the order of the product of any two
involutions conjugate
to $i$ is finite;

 \makebox[15pt][r]{}2) all 2-sub\-groups of $G$ are either cyclic
or generalized
quaternion;

 \makebox[15pt][r]{}3) the centralizer $C$ of the involution
$i$ in $G$ is infinite, distinct from $G$, and has finite periodic
part;

 \makebox[15pt][r]{}4) the normalizer of any non-trivial
$i $-inva\-ri\-ant
finite subgroup in~$G$ either is contained in $C$ or has
periodic part which is a Frobenius group (see 6.55) with abelian kernel
and finite complement of even order;

 \makebox[15pt][r]{}5) for every element $c$ not contained in
$C$ for which $ci$ is an involution there is an element $s$ of $C$
such that $\left< c,c^s\right>$ is an infinite subgroup.

 \makebox[15pt][r]{}Does there exist a simple $T_0 $-group? \hfill
\raisebox{0ex}{\sl V.\,P.\,Shunkov}

\emp

\raggedbottom

\newpage

\pagestyle{myheadings} \markboth{13th
Issue (1995)}{13th Issue (1995)}
\thispagestyle{headings} ~

\vspace{2ex}

\centerline {\Large \textbf{Problems from the 13th Issue (1995)}}
\phantomsection\label{13izd}
\vspace{4ex}

\bmp \textbf{13.1.} Let $U(K)$ denote the group
of units of a ring $K$. Let $G$ be a finite group, ${\Bbb Z} G$
the integral group ring of $G$, and ${\Bbb Z} _pG$ the group ring
of $G$ over the residues modulo a prime number $p$.
 Describe the homomorphism
from $U({\Bbb Z} G)$ into $U({\Bbb Z} _pG)$ induced by reducing the
coefficients modulo $p$. More precisely, find the kernel and the
image of this homomorphism and
 an explicit
transversal over the kernel. \hfill \raisebox{-1ex}{\sl
R.\,Zh.\,Aleev}

 \emp

\bmp \textbf{13.2.} Does there exist a finitely
based variety of groups
 ${\frak V}$ such that the word problem is solvable
in $F_n({\frak V})$ for every positive integer~$n$,
but is unsolvable in
$F_{\infty}({\frak V})$ (with respect to a free
generating system)? \hfill \raisebox{-1ex}{\sl M.\,I.\,Anokhin}

\emp

\bmp \textbf{13.3.} Let ${\frak M}$ be an
arbitrary variety of groups. Is it true that every infinitely
generated projective group in
 ${\frak M}$ is an ${\frak
M} $-free product of countably generated projective groups in
 ${\frak M}$? \hfill \raisebox{-1ex}{\sl V.\,A.\,Artamonov}

 \emp

\bmp \textbf{13.4.} Let $G$ be a group with a
normal (pro-)\hskip0.5ex 2-sub\-group $N$ such that $G/N$ is
isomorphic to $GL_n(2)$, inducing its natural module on $N/\Phi
(N)$, the Frattini factor-group of~$N$. For $n=3$ determine $G$
such that $N$ is as large as possible. (For $n>3$ it can be proved
that already the Frattini subgroup $\Phi (N)$ of $N$ is trivial;
for $n=2$ there exists $G$ such that $N$ is the free pro-2-group
generated by two elements). \hfill \raisebox{-1ex}{\sl
B.\,Baumann}

 \emp

\bmp \textbf{13.5.} Groups $A$ and $B$ are said to
be {\it locally equivalent\/} if for every finitely generated
subgroup $X\leq A$ there is a subgroup $Y\leq B$ isomorphic to
$X$, and conversely, for every finitely generated subgroup
 $Y\leq B$ there is a subgroup $X\leq A$, isomorphic to
$Y$. We call a group $G$ {\it categorical\/} if $G$ is isomorphic
to any group that is locally equivalent to $G$.
 Is it true that a periodic locally soluble group
 $G$ is categorical if and only if
 $G$ is hyperfinite with Chernikov Sylow subgroups?
(A group is {\it hyperfinite\/} if it has a well ordered ascending
normal series with finite factors.)\hfill \raisebox{-1ex}{\sl
V.\,V.\,Belyaev}

 \emp

\bmp \textbf{13.6.} (B.\,Hartley). Is it true that
a locally finite group containing an ele\-ment with Chernikov
centralizer is almost soluble?\hfill \raisebox{-1ex}{\sl
V.\,V.\,Belyaev}

 \emp

\bmp \textbf{13.7.} (B.\,Hartley). Is it true that
a simple locally finite group containing a finite subgroup with
finite centralizer is linear?\hfill \raisebox{-1ex}{\sl
V.\,V.\,Belyaev}

 \emp

\bmp \textbf{13.8.} (B.\,Hartley). Is it true that a locally soluble
periodic group has a finite normal series with locally nilpotent
factors if it contains an element

\makebox[25pt][r]{\zva a)}
with finite centralizer?

\makebox[25pt][r]{b)} with Chernikov centralizer?\hfill
\raisebox{0ex}{\sl V.\,V.\,Belyaev}

\ul

\otv a) Yes, it is true (E.\,Khukhro,
\emph{Bull. London Math. Soc.}, \textbf{58}, no.~1 (2026), Article ID e70206).
 \emp

\bmp \textbf{13.9.} Is it true that a locally
soluble periodic group containing a finite nilpotent subgroup with
Chernikov centralizer has a finite normal series with locally
nilpotent factors?\hfill \raisebox{-1ex}{\sl V.\,V.\,Belyaev,
 B.\,Hartley}

 \emp

\bmp \textbf{\zv 13.12.}
 Is the group of all automorphisms of an
arbitrary hyperbolic group finitely presented?

\makebox[15pt][r]{}{\it Comment of 2001:\/} This is proved for the
group of automorphisms of a torsion-free hyperbolic group that
cannot be decomposed into a free product (Z.\,Sela, {\it Geom.
Funct. Anal.}, \textbf{7}, no.\,3 (1997) 561--593).
\hfill
\raisebox{-1ex}{\sl O.\,V.\,Bogopolski}

\ul

\otv Yes, it is (M.\,Carette, {\it J.~London Math. Soc.}, {\bf 84} (2011), 731--748).
 \emp

\bmp \textbf{13.13.} Let $G$ and $H$ be finite $p\hskip0.2ex $-groups
with isomorphic Burnside rings. Is the nil\-po\-ten\-cy class of
$H$ bounded by some function of the class of $G$? There is an
example where $G$ is of class 2 and $H$ is of class 3. \hfill
\raisebox{-1ex}{\sl R.\,Brandl}

 \emp

\bmp \textbf{13.14.} Is the lattice of quasivarieties of nilpotent
torsion-free groups of nilpotency class $\leq 2$
distributive?\hfill \raisebox{-1ex}{\sl A.\,I.\,Budkin}

 \emp

\bmp \textbf{13.15.} Does a free non-abelian nilpotent group of class
3 possess an independent basis of quasiidentities in the class of
torsion-free groups? \hfill \raisebox{-1ex}{\sl A.\,I.\,Budkin}
 \emp

\bmp \textbf{13.17.} A representation of a group $G$ on a vector
space $V$ is called a {\it nil-repre\-sen\-tation\/} if for every
$v$ in $V$ and $g$ in $G$ there exists $n=n(v,g)$ such that $v
(g-1)^n = 0$. Is it true that in zero characteristic every
irreducible nil-representation is trivial? (In prime
characteristic it is not true.) \hfill \raisebox{-1ex}{\sl
S.\,Vovsi}

 \emp

\bmp \textbf{\zv 13.19.}
Suppose that both a group $Q$ and its normal
subgroup
 $H$ are subgroups of the direct product $G_1\times \cdots \times
G_n$ such that for each $i$ the projections of both $Q$ and $H$
onto $G_i$ coincide with $G_i$. If $Q/H$ is a $p\hskip0.2ex
$-group, is $Q/H$ a regular $p\hskip0.2ex $-group?

\hfill
\raisebox{-1ex}{\sl Yu.\,M.\,Gorchakov}

\ul

\otv No, not always; moreover, $Q/H$ can be any nilpotent group of class $n-1$  (AI solution confirmed by E.\,I.\,Khukhro: \ A.\,Ferudun, {\it Preprint}, 2026, \url{https://doi.org/10.5281/zenodo.23453462}).
 \emp

\bmp \textbf{13.20.} Is it true that the generating series of the
growth function of every finitely generated group with one
defining relation represents an algebraic (or even a rational)
function?\hfill \raisebox{-1ex}{\sl R.\,I.\,Grigorchuk}

 \emp

\bmp \textbf{13.21.} b) What is the minimal possible rate of
 growth of the function $
\pi (n)=\max\limits_{\delta (g)\leq n} |g|$ for the class of
groups indicated in part ``a'' of this problem (see Archive)? It is
known (R.\,I.\,Grigorchuk, {\it Math. USSR--Izv.}, \textbf{25}
(1985), 259--300) that there exist $p\hskip0.2ex $-groups with
$\pi (n)\leq n^{\lambda}$ for some $\lambda >0$. At the same time,
it follows from a result of E.\,I.\,Zel'manov that $\pi (n)$ is
not bounded if $G$ is infinite. \hfill \raisebox{-1ex}{\sl
R.\,I.\,Grigorchuk}

 \emp

\bmp
\textbf{13.22.}
 Let the group $G=AB$ be the product of two polycyclic subgroups
$A$ and~$B$, and assume that $G$ has an ascending normal series
with locally nilpotent factors (i.~e. $G$ is radical). Is it true
that $G$ is polycyclic?
\hfill \raisebox{-1ex}{\sl
\,F.\,\,de\,\,Giovanni}

 \emp

\bmp
\textbf{\zv 13.23.}
 Let the group $G$ have a finite normal
series with infinite cyclic factors (containing of course $G$ and
$\{1\}$). Is it true that $G$ has a non-trivial outer
automorphism?

\hfill \raisebox{-1ex}{\sl \,F.\,\,de\,\,Giovanni}

\ul

\otv No, not always (M.\,Brescia, E.\,Ingrosso, M.\,Trombetti, {\it Preprint}, 2026, \url{https://arxiv.org/abs/2607.11775}).
 \emp

\bmp \textbf{13.24.} Let $G$ be a non-discrete topological group with
only finitely many ultrafilters that converge to the identity. Is
it true that
$G$ contains a countable open subgroup of exponent 2?

\makebox[15pt][r]{}{\it Comment of 2005:\/} This is proved for $G$
countable (E.\,G.\,Zelenyuk, {\it Mat. Studii}, \textbf{7}
 (1997), 139--144).
\hfill \raisebox{-1ex}{\sl E.\,G.\,Zelenyuk}
 \emp

\bmp \textbf{13.25.} Let $(G,\tau )$ be a topological group with
finite semigroup $\tau (G)$ of ultrafilters converging to the
identity (I.\,V.\,Protasov, {\it Siberian Math.~J.}, \textbf{34},
no.\,5 (1993), 938--951). Is it true that
 $\tau (G)$ is a semigroup of idempotents?

\makebox[15pt][r]{}{\it Comment of 2005:\/} This is proved for $G$
countable
 (E.\,G.\,Zelenyuk, {\it Mat. Studii}, \textbf{14} (2000), 121--140).
\hfill \raisebox{-1ex}{\sl E.\,G.\,Zelenyuk}
 \emp

\bmp \textbf{13.27.} (B.\,Amberg). Suppose that $G=AB=AC=BC$ for a
group $G$ and its subgroups
 $A,B,C$.

\makebox[25pt][r]{a)} Is $G$ a Chernikov group if
 $
A,B,C$ are Chernikov groups?

\makebox[25pt][r]{b)} Is $G$ almost polycyclic if
 $A,B,C$ are almost polycyclic?
\hfill \raisebox{0ex}{\sl L.\,S.\,Kazarin}

 \emp

\bmp \textbf{13.30.} A group
 $G$ is called a {\it $B$-group\/} if every primitive
permutation group which contains the regular representation of
 $G$ is doubly transitive. Are there any countable $B $-groups?\hfill
\raisebox{-1ex}{\sl P.\,J.\,Cameron}

 \emp

\bmp \textbf{13.31.} Let
 $G$ be a permutation group on a set
$\Omega$. A sequence of points of $\Omega$ is a {\it base\/} for
$G$ if its pointwise stabilizer in $G$ is the identity. {\it The
greedy algorithm\/} for a base chooses each point in the sequence
from an orbit of maximum size of the stabilizer of its
predecessors. Is it true that there is a universal constant
 $c$ with the property that, for any finite primitive
permutation group, the greedy algorithm produces a base whose
size is at most
$c$ times the minimal base size?
\hfill
\raisebox{-1ex}{\sl
P.\,J.\,Cameron}

 \emp

\bmp \textbf{13.32.} (Well-known problem). On a group, a partial
order that is directed upwards and has linearly ordered cone of
positive elements is said to be {\it semilinear\/} if this order
is invariant under right multiplication by the group elements. Can
every semilinear order of a group be extended to a right order of
the group? \hfill\raisebox{-1ex}{\sl V.\,M.\,Kopytov}

 \emp

\bmp \textbf{13.35.} Does every non-soluble pro-$p\hskip0.2ex $-group
of cohomological dimension
 2 contain a free non-abelian
pro-$p\hskip0.2ex $-sub\-group? \hfill \raisebox{-1ex}{\sl
O.\,V.\,Mel'nikov}

 \emp

\bmp \textbf{13.36.} For a finitely generated pro-$p\hskip0.2ex
$-group $G$ set $a_n(G)={\rm dim}_{\,{\Bbb F}_p}I^n/I^{n+1}$,
where $I$ is the augmentation ideal of the group ring
 ${\Bbb F}_p[[G]]$.
We define the {\it growth\/} of $G$ to be the growth of the
sequence
 $\{ a_n(G)\} _{n\in {\Bbb N}}$.

\makebox[25pt][r]{a)} If the growth of $G$ is exponential, does it
follow that $G$ contains a free pro-$p\hskip0.2ex $-sub\-group of
rank
 2?

\makebox[25pt][r]{c)} Do there exist pro-$p\hskip0.2ex $-groups of
finite cohomological dimension which are not $p\hskip0.2ex $-adic
analytic, and whose growth is slower than an exponential one?

\hfill \raisebox{-1ex}{\sl O.\,V.\,Mel'nikov}

 \emp

\bmp \textbf{13.37.} Let $G$ be a torsion-free pro-$p\hskip0.2ex
$-group, $U$ an open subgroup of $G$. Suppose that $U$ is a
pro-$p\hskip0.2ex $-group with a single defining relation. Is it
true that then $G$ is also a pro-$p\hskip0.2ex $-group with a
single defining relation?\hfill \raisebox{-1ex}{\sl
O.\,V.\,Mel'nikov}
 \emp

\bmp
\textbf{13.39.} Let $A$ be an associative ring with unity and with
torsion-free additive group, and let $F^A$ be the tensor product
of a free group $F$ by $A$ (A.\,G.\,Myasnikov,
V.\,N.\,Re\-mes\-len\-ni\-kov, {\it Siberian Math.~J.}, \textbf{35},
no.\,5 (1994), 986--996);
 then $F^A$ is a {\it free exponential group over $A$}; in
(A.\,G.\,Myasnikov, V.\,N.\,Remeslennikov, {\it
Int. J.~Algebra Comput.}, \textbf{6} (1996), 687--711), it is shown
how to construct
 $F^A$ in terms of free products with amalgamation.

\makebox[25pt][r]{b)}
Is $F^{A}$ a linear group?

\makebox[15pt][r]{}In the case where $\left< 1\right>$ is a pure
subgroup of the additive group of $A$, there is an affirmative answer
to ``b'' (A.\,M.\,Gaglione, A.\,G.\,Myasnikov,
V.\,N.\,Re\-mes\-lennikov, D.\,Spellman, {\it Commun. Algebra},
\textbf{25} (1997), 631--648).  It is also known that the Magnus
homomorphism is one-to-one on any subgroup of $F^{\Bbb Q}$ of the
type $\left< F,t \mid u = t^n\right>$ (G.\,Baumslag, {\it Commun.
Pure Appl. Math.}, \textbf{21} (1968), 491--506).

\makebox[25pt][r]{d)}
Is the universal theory of $F^{A}$ decidable?

\makebox[25pt][r]{e)} (G.\,Baumslag). Can free $A $-groups be
characterized by a length function?

\makebox[25pt][r]{f)} (G.\,Baumslag). Does a free $\Bbb
Q\hskip0.1ex $-group admit a free action on some $\Lambda $-tree?
See definition in (R.\,Alperin,
 H.\,Bass, {\it in: Combinatorial group theory and topology,
Alta, Utah, 1984 $($Ann. Math. Stud.}, \textbf{111}), Princeton Univ.
Press, 1987, 265--378).

\hfill \raisebox{-1ex}{\sl
A.\,G.\,Myasnikov,
V.\,N.\,Remeslennikov}
\emp

 \bmp

\textbf{13.41.} Is the elementary theory of the class of all groups
acting freely on $\Lambda $-trees decidable?\hfill
\raisebox{-1ex}{\sl A.\,G.\,Myasnikov, V.\,N.\,Remeslennikov}

 \emp

\bmp
\textbf{\zv  13.42.}
Prove that the tensor $A $-comple\-tion (see
A.\,G.\,Myas\-ni\-kov, V.\,N.\,Re\-mes\-len\-ni\-kov, {\it
Siberian Math.~J.}, \textbf{35}, no.\,5 (1994), 986--996)
 of a free nilpotent group can be
non-nilpotent.
 \hfill \raisebox{-1ex}{\sl A.\,G.\,Myasnikov, V.\,N.\,Remeslennikov}

 \ul

\otv An example is produced (AI solution confirmed by A.\,G.\,Myasnikov: \  M.\,Kinyon, J.\,Urban, {\it Preprint}, 2026,
\url{https://kourovkanotebookorg.wordpress.com/wp-content/uploads/2026/09/kourovka-mathematics.pdf}).
 \emp 

\bmp

\textbf{\zv 13.43.}
(G.\,R.\,Robinson). Let $G$ be a finite group and $B$
be a $p\hskip0.2ex $-block of characters of $G$. {\it
Conjecture:\/} If the defect group $D=D(B)$ of the block $B$ is
non-abelian, and if $|D:Z(D)| = p^a$, then each character in $B$
has height strictly less than $a$.

\makebox[15pt][r]{}{\it G.\,R.\,Robinson's comment of 2020}:
the conjecture is proved mod CFSG for $p\ne 2$ in (Z.\,Feng, C.\,Li, Y.\,Liu, G.\,Malle, J.\,Zhang, {\it Compos. Math.}, {\bf 155}, no.\,6 (2019), 1098--1117).
\hfill \raisebox{-1ex}{\sl J.\,Olsson}

\ul

\otv The conjecture is proved (R.\,Kessar, G.\,Malle, {\it J.~London Math. Soc. (2)}, {\bf 111}, no.\,2 (2025), Article ID e70076, \url{https://arxiv.org/abs/2311.13510}).
 \emp

\bmp \textbf{13.44.} For any partition of an arbitrary group
 $G$ into finitely many subsets
$G=A_1\cup \cdots \cup A_n$, there exists a subset of the
partition $A_i$ and a finite subset $F\subseteq G$, such that
$G=A_i^{-1}A_iF$ (I.\,V.\,Protasov, {\it Siberian Math.~J.}, {\bf
34}, no.\,5 (1993), 938--952). Can the subset $F$ always be chosen
so that
 $|F|\leq n$? This is true for amenable groups.

\hfill \raisebox{-1ex}{\sl I.\,V.\,Protasov}
 \emp

\bmp \textbf{13.48.} (W.\,W.\,Comfort, J.\,\,van\,\,Mill). A
topological group is said to be {\it
irre\-sol\-vable\/} if every two of its dense subsets intersect
non-trivially. Does every non-discrete irresolvable group contain
an infinite subgroup of exponent 2?\hfill \raisebox{-1ex}{\sl
I.\,V.\,Protasov}
 \emp

\bmp \textbf{13.49.} (V.\,I.\,Malykhin). Can a topological group be
partitioned into two dense subsets, if there are infinitely many free
ultrafilters on the group converging to the identity? \hfill
\raisebox{-1ex}{\sl I.\,V.\,Protasov}
 \emp

\bmp
\textbf{\zv 13.51.}
Is every finite modular lattice embeddable in the
lattice of formations of finite groups?\hfill \raisebox{-1ex}{\sl
A.\,N.\,Skiba}

\ul

\otv No, not every (AI solution confirmed by A.\,N.\,Skiba: \  Achyuth Jayadevan,  {\it Preprint}, 2026, \url{https://zenodo.org/records/22863019}).
 \emp 

\bmp

\textbf{13.52.} The {\it dimension\/} of a finitely based variety of
algebras
 ${\cal V}$ is defined to be the maximal length of
a {\it basis\/} (that is, an independent generating set) of the
$SC $-theory $SC({\cal V})$, which consists of the strong Mal'cev
conditions satisfied on ${\cal V}$. The dimension is defined to
be infinite if the lengths of bases in $SC({\cal V})$ are not
bounded.
 Does every finite abelian group generate a variety of finite
dimension?
\hfill \raisebox{-1ex}{\sl D.\,M.\,Smirnov}

 \emp

\bmp

\textbf{13.53.} Let $a,\, b$ be elements of finite order of the infinite
group $G=\left< a,\, b\right>$. Is it true that there are infinitely
many elements $g\in G$ such that the subgroup $\left< a,\,
b^g\right>$ is infinite?\hfill \raisebox{-1ex}{\sl A.\,I.\,Sozutov}

 \emp

\bmp \textbf{13.54.} a) Is it true that, for
 $
p$ sufficiently large, every (finite) $p\hskip0.2ex $-group can be
a~com\-ple\-ment in some Frobenius group
 (see 6.55)? \hfill \raisebox{-1ex}{\sl A.\,I.\,Sozutov}

 \emp

\bmp
 \textbf{13.55.} Does there exist a
Golod group (see 9.76), which is isomorphic to an $AT $-group? For
a definition of an {\it $AT $-group\/} see (A.\,V.\,Rozhkov, {\it
Math. Notes}, \textbf{40}, no.\,5 (1986), 827--836). \hfill
\raisebox{-1ex}{\sl A.\,V.\,Timofeenko}

 \emp

\bmp \textbf{13.57.} Let $\f$ be an automorphism of prime order
 $p$ of a finite group $G$ such that $C_G(\f )\leq
Z(G)$.

\makebox[25pt][r]{a)} Is $G$ soluble if $p=3$? V.\,D.\,Mazurov and
T.\,L.\,Nedo\-re\-zov proved in ({\it Algebra and Logic}, {\bf
35}, no.\,6 (1996), 392--397) that the group $G$ is soluble for
$p=2$, and there are examples of unsoluble $G$ for all $p>3$.

\makebox[25pt][r]{b)} If $G$ is soluble, is the derived length of
$G$ bounded in terms of $p$?

\makebox[25pt][r]{c)} (V.\,K.\,Kharchenko). If $G$ is a
$p\hskip0.2ex $-group, is the derived length of $G$ bounded in
terms of $p$? (V.\,V.\,Bludov produced a simple example showing
that the nilpotency class cannot be bounded.)
\hfill\raisebox{-1ex}{\sl E.\,I.\,Khukhro}

 \emp

\bmp \textbf{13.59.} One can show that any extension of shape ${\Bbb
Z}^d.\, {\cal SL}_d({\Bbb Z})$ is residually finite, unless possibly
$d=3$ or $d=5$. Are there in fact any extensions of this shape that
fail to be residually finite when $d=5$? As shown in (P.\,R.\,Hewitt,
{\it Groups/St.\,Andrews'93 in Galway, Vol.\,2 $($London Math. Soc.
Lecture Note Ser.}, \textbf{212}), Cambridge Univ. Press, 1995,
305--313), there is an extension of ${\Bbb Z} ^3$ by ${\cal
SL}_3({\Bbb Z} )$ that is not residually finite. Usually it is true
that an extension of an arithmetic subgroup of a Chevalley group over
a rational module is residually finite. Is it ever false, apart from
the examples of shape ${\Bbb Z}^3.\, {\cal SL}_3({\Bbb Z})$? \hfill
\raisebox{-1ex}{\sl P.\,R.\,Hewitt}

 \emp

\bmp \textbf{13.60.} If a locally graded group $G$ is a product of
two subgroups of finite special rank, is $G$ of finite special
rank itself? The answer is affirmative if both factors are
periodic groups. (N.\,S.\,Chernikov, {\it Ukrain. Math.~J.}, {\bf
42}, no.\,7 (1990),
 855--861).

\hfill \raisebox{-1ex}{\sl N.\,S.\,Chernikov}

 \emp

\bmp \textbf{13.64.} Let $\pi _e(G)$ denote the set of orders of elements of a group
$G$. A group $G$ is said to be an {\it $OC_n
$-group\/} if $\pi _e(G)=\{ 1,\, 2, \ldots ,\, n\}$. Is every
$OC_n $-group locally finite? Do there exist infinite $OC_n
$-groups for $n\geq 7$? \hfill \raisebox{-1ex}{\sl W.\,J.\,Shi}

 \emp

\bmp \textbf{13.65.} A finite simple group is called a {\it
$K_n$-group\/} if its order is divisible by exactly
 $n$ different primes.
 The number of $K_3 $-groups is known to be~8.
The $K_4$-groups are classified mod CFSG (W.\,J.\,Shi, {\it in:
 Group Theory in China $($Math. Appl.}, \textbf{365}), Kluwer,
1996, 163--181) and some
significant further results are obtained in (Yann Bugeaud, Zhenfu
Cao, M.\,Mignotte, {\it J.~Algebra}, \textbf{241} (2001),
 658--668). But the question remains: is the number of $K_4
$-groups finite or
 infinite?
\hfill \raisebox{-1ex}{\sl W.\,J.\,Shi}

 \emp

\bmp \textbf{13.67.} Let $G$ be a $T_0 $-group (see~12.101), $i$ an
involution in $G$ and $G=\left< i^g\mid g\in G\right>$. Is the
centralizer $C_G(i)$ residually periodic?\hfill
\raisebox{-1ex}{\sl V.\,P.\,Shunkov}

 \emp

\raggedbottom
\newpage

\pagestyle{myheadings} \markboth{14th Issue
(1999)}{14th Issue (1999)}
\thispagestyle{headings} ~ \vspace{2ex}

\centerline {\Large \textbf{Problems from the 14th Issue (1999)}}
\phantomsection\label{14izd}
\vspace{4ex}

\bmp \textbf{14.2.} (S.\,D.\,Berman). Prove that every automorphism of
the centre of the integral group ring of a finite group induces a
monomial permutation on the set of the class sums.
 \hfill \raisebox{-1ex}{\sl R.\,Zh.\,Aleev}

\emp

\bmp \textbf{14.3.} Is it true that every central unit of the
integral group ring of a finite group is a product of a central
element of the group and a symmetric central unit? (A unit is {\it
symmetric\/} if it is fixed by the canonical antiinvolution that
 transposes the
coefficients at the mutually inverse elements.)
 \hfill \raisebox{-1ex}{\sl R.\,Zh.\,Aleev}
 \emp

\bmp \textbf{\zv 14.4.}
a) Is it true that there exists a nilpotent group
$G$ for which the lattice ${\cal L}(G)$ of all group topologies is
not modular? (It is known that for abelian groups the lattice
 ${\cal L}(G)$ is modular and that there are groups for which
this lattice is not modular: V.\,I.\,Arnautov, A.\,G.\,Topale,
{\it Izv. Akad. Nauk Moldova Mat.}, \textbf{1997}, no.\,1, 84--92
(Russian).)

\makebox[25pt][r]{b)} Is it true that for every countable
nilpotent non-abelian group $G$ the lattice ${\cal L}(G)$ of all
group topologies is not modular? \hfill \raisebox{-1ex}{\sl
V.\,I.\,Arnautov}

\ul

\otv a) Yes, it is true (V.\,Arnautov, A.\,Topal\u a, {\it Bul. Acad. \c{S}tiin\c{t}.
Repub. Moldova, Mat.}, 1998, no.~2(27), 130--131).

\otv b) No, there are such groups with modular ${\cal L}(G)$ (Dekui Peng, {\it Preprint}, 2023, \url{https://arxiv.org/abs/2310.08269}).
 \emp

\bmp \textbf{14.5.} Let $G$ be an infinite group admitting
non-discrete Hausdorff group topologies, and ${\cal L}(G)$ the
lattice of all group topologies on~$G$.

\makebox[25pt][r]{a)} Is it true that for any natural number
 $k$ there exists a non-refinable chain $\tau_0 < \tau_1<\cdots
< \tau_k$ of length $k$ of Hausdorff topologies in ${\cal L}(G)$?
(For countable nilpotent groups this is true (A.\,G.\,Topale, {\it
Deposited in VINITI}, 25.12.98, no.\,3849--V~98 (Russian)).)

\makebox[25pt][r]{\zva b)}
Let $k,m,n$ be natural numbers and let $G$
be a nilpotent group of class~$k$. Suppose that $\tau_0 <\tau_1
<\cdots < \tau_m$ and $\tau'_0 <\tau'_1 <\cdots < \tau'_n$ are
non-refinable chains of Hausdorff topologies in ${\cal L}(G)$ such
that $\tau_0 =\tau'_0$ and $\tau_m =\tau'_n$. Is it true that
$m\leq n\cdot k$ and this inequality is best-possible? (This is
true if $k=1$, since for $G$ abelian the lattice ${\cal L}(G)$ is
modular.)

\makebox[25pt][r]{c)} Is it true that there exists a countable
 $G$ such that in the lattice ${\cal L}(G)$ there are a finite
non-refinable chain $\tau_0 <\tau_1 <\cdots < \tau_k$ of Hausdorff
topologies and an infinite chain $\{\tau'_{\gamma}\mid
\gamma\in\Gamma\}$ of topologies such that $\tau_0<\tau'_{\gamma}
< \tau_k$ for any $\gamma\in\Gamma$?

\makebox[25pt][r]{\zva d)}
Let $G$ be an abelian group, $k$ a natural
number. Let $A_k$ be the set of all those Hausdorff group topologies
on $G$ that, for every topology $\tau\in A_k$, any non-refinable
chain of topologies starting from $\tau$ and terminating at the
discrete topology has length~$k$. Is it true that $A_k\bigcap
\{\tau'_{\gamma}\mid \gamma\in\Gamma \}\neq\varnothing$ for any
infinite non-refinable chain $\{\tau'_{\gamma}\mid \gamma\in\Gamma
\}$ of Hausdorff topologies containing the discrete topology? (This
is true for $k=1$.) \hfill \raisebox{-1ex}{\sl V.\,I.\,Arnautov}

\ul

\otv b) The inequality does hold, but it is not sharp (V.\,I.\,Arnautov,
{\it Bull. Acad. \c{S}tiin\c{t}e Repub. Moldova, Mat.}, \textbf{2010}, no. 2 (2010), 3--19).

\otv d) No; moreover, no infinite abelian group satisfies this property with $k=2$ (D.\,Peng, {\it Preprint}, 2023, \url{https://arxiv.org/abs/2310.08269}).
 \emp

\bmp \textbf{14.6.}
A group $\Gamma$ is said to have {\it Property
$P_{\rm{nai}}$} if, for any finite subset $F$ of $\Gamma \setminus
\{1\},$ there exists an element $y_0 \in \Gamma$ of infinite order
such that, for each $x \in F$, the canonical epimorphism from the
free product $\langle x \rangle * \langle y_0\rangle$ onto the
subgroup $\langle x , y_0 \rangle$ of $\Gamma$ generated by $x$
and $y_0$ is an isomorphism.
 For $n \in \{2,\,3,\ldots\},$ does $PSL_n({\Bbb Z} )$ have
Property $P_{\rm{nai}}$? More generally, if $\Gamma$ is a lattice
in a connected real Lie group $G$ which is simple and with centre
reduced to $\{ 1\},$ does $\Gamma$ have Property $P_{\rm{nai}}$?

\makebox[15pt][r]{}Answers are known to be \lq\lq yes\rq\rq \ if $n =
2$, and more generally if $G$ has real rank $1$ (M.\,Bekka,
M.\,Cowling, P.\,de\,\,la\,\,Harpe, {\it Publ. Math. IHES}, \textbf{80}
(1994), 117--134).

\makebox[15pt][r]{}{\it Comments of 2018}:
Property $P_{\rm{nai}}$ has been established for relatively hyperbolic groups
(G.\,Arzhantseva, A.\,Minasyan, \emph{J.~Funct. Anal.}, \textbf{243}, no.\,1 (2007), 345--351) and for groups acting on CAT(0) cube complexes
with appropriate conditions (A.\,Kar, M.\,Sageev, \emph{Comment. Math. Helv.}, \textbf{91} (2016), 543--561). Some progress has been made in the preprint (T.\,Poznansky, \url{http://arxiv.org/abs/0812.2486}).

 \hfill \raisebox{-1ex}{\sl P.\,\,de la Harpe}
 \emp

 \bmp \textbf{14.7.} Let $\Gamma_g$ be the fundamental group of
a closed surface of genus $g \geq 2.$ For each finite system $S$
of generators of $\Gamma_g,$ let $\beta_S(n)$ denote the number of
elements of $\Gamma _g$ which can be written as products of at
most $n$ elements of $S \cup S^{-1},$ and let $\omega (\Gamma
_g,S) = \limsup_{n \to \infty}\root n \of {\beta _S(n)}$ be the
growth rate of the sequence $\left( \beta _S(n) \right)_{n \geq
0}.$
 Compute the infimum $\omega(\Gamma _g)$ of the $\omega (\Gamma _g,S)$ over
all finite sets of generators of $\Gamma _g.$

\makebox[15pt][r]{}It is easy to see that, for a free group $F_k$
of rank $k \geq 2,$ the corresponding infimum is $\omega(F_k) =
2k-1$ (M.\,Gromov, {\it Structures m\'etriques pour les
vari\'et\'es riemanniennes,} Cedic/F.\,Nathan, Paris, 1981,
Ex.\,5.13). As any generating set of $\Gamma_g$ contains a subset
of $2g-1$ elements generating a subgroup of infinite index in
$\Gamma_g$ with abelianization ${\Bbb Z} ^{2g-1},$ hence a
subgroup which is free of rank $2g-1,$ it follows that
$\omega(\Gamma_g) \geq 4g-3.$

\hfill \raisebox{-0.5ex}{\sl P.\,\,de la Harpe}

 \emp

\bmp \textbf{14.8.} Let $G$ denote the group of germs at $+\infty$ of
orientation-preserving homeomorphisms of the real line ${\Bbb R}$.
 Let $\alpha \in G$ be the germ of $x \mapsto x+1.$ What are the germs
$\beta \in G$ for which the subgroup $\langle \alpha , \beta \rangle$ of
$G$ generated by $\alpha$ and $\beta$ is free of rank~$2$?

\makebox[15pt][r]{}If $\beta$ is the germ of $x \mapsto x^k$ for
an odd integer $k \geq 3$, it is known that $\langle \alpha ,
\beta \rangle$ is free of rank 2. The proofs of this rely on
Galois theory (for $k$ an odd prime: S.\,White, {\it J.~Algebra},
\textbf{118} (1988), 408--422; \ for any odd $k\geq 3$:
S.\,A.\,Adeleke, A.\,M.\,W.\,Glass, L.\,Morley, {\it J.~London
Math. Soc.}, \textbf{43} (1991), 255--268, and for any odd $k\ne \pm 1$
and any even $k>0$ in several papers by S.\,D.\,Cohen and A.\,M.\,W.\,Glass). \hfill
\raisebox{-0.5ex}{\sl P.\,\,de la Harpe}
 \emp

\bmp
\textbf{14.9.}
(Well-known problem). Let $W^*(F_k)$ denote the von~Neumann
algebra of the free group of rank $k \in \{ 2,\,3, \ldots ,\aleph_0 \}.$
Is $W^*(F_k)$ isomorphic to $W^*(F_l)$ for $k \ne l$?

\makebox[15pt][r]{}For a group $G$, recall that $W^*(G)$ is an
appropriate completion of the group algebra ${\Bbb C} G$ (see
e.~g. S.\,Sakai, {\it $C^*$-algebras and $W^*$-algebras,}
Springer, 1971, in particular Problem~4.4.44). It is known that
either $W^*(F_k) \cong W^*(F_l)$ for all $k,l \in \{ 2,\,3, \ldots
,\aleph_0 \},$ or that the $W^*(F_k)$ are pairwise non-isomorphic
(F.\,Radulescu, {\it Invent. Math.}, \textbf{115} (1994), 347--389,
Corollary 4.7). \hfill \raisebox{-0.5ex}{\sl P.\,\,de la Harpe}

 \emp

\bmp
\textbf{14.10.} c) Find an explicit and {\it \lq\lq
natural\rq\rq} \ finitely presented group $\Gamma_n$ and an
embedding of $GL_n({\Bbb Q})$ in $\Gamma _n$.

 \makebox[15pt][r]{}Another phrasing of the same
problems is: find a simplicial complex $X$ which covers a finite
complex such that the fundamental group of $X$ is ${\Bbb Q}$ or,
respectively, $GL_n({\Bbb Q})$.
 \hfill \raisebox{-1ex}{\sl
P.\,\,de la Harpe}
\emp

\bmp
\textbf{14.11.}
(Yu.\,I.\,Merzlyakov). It is a
well-known fact that for the ring
 $R={\Bbb Q} [x,y]$ the elementary group $E_2(R)$ is distinct from
$SL_2(R)$. Find a minimal subset $A\subseteq SL_2(R)$ such that
$\left< E_2(R), A\right> =SL_2(R)$. \hfill \raisebox{-1ex}{\sl
V.\,G.\,Bardakov}

 \emp

\bmp \textbf{14.12.} (Yu.\,I.\,Merzlyakov, J.\,S.\,Birman). Is it
true that all braid groups $B_n$, $n\geq 3$, are
 conjugacy separable?
\hfill \raisebox{-1ex}{\sl V.\,G.\,Bardakov}
 \emp

\bmp \textbf{14.14.} (C.\,C.\,Edmunds, G.\,Rosenberger). We call
a~pair of natural numbers $(k,m)$ {\it admissible\/} if
 in the derived subgroup $F_2'$ of a free group $F_2$
there is an element $w$ such that the commutator length
of the element $w^m$ is equal to $k$.
 Find all admissible pairs.

\makebox[15pt][r]{}It is known that every pair $(k, 2)$, $k \geq 2$,
is admissible (V.\,G.\,Bardakov, {\it Algebra and Logic}, \textbf{39}
(2000), 224--251). \hfill \raisebox{-1ex}{\sl V.\,G.\,Bardakov}
 \emp

\bmp \textbf{\zv 14.15.} For the automorphism group $A_n={\rm Aut\,}F_n$
of a free group $F_n$ of rank $n\geq 3$ find the supremum
 $k_n$ of the commutator lengths of the
elements of $A_n'$.

\makebox[15pt][r]{}It is easy to show that
$k_2=\infty$.
On the other hand,
the commutator length of any element of the derived subgroup of
 $\lim\limits_{\longrightarrow}A_n$ is at most $2$
 (R.\,K.\,Dennis,
L.\,N.\,Vaserstein, {\it $K$-Theory}, \textbf{2}, N\,6 (1989),
761--767). \hfill \raisebox{-1ex}{\sl V.\,G.\,Bardakov}

\ul

\otv For $n\geq 3$ we have $k_n=\infty$ because $A_n$ is acylindrically hyperbolic (A.\,Genevois,  C.\,Horbez, {\it J.~Topology}, {\bf 14} (2021), 963--991) and for such groups the commutator width is infinite (M.\,Bestvina, K.\,Bromberg,  K.\,Fujiwara, {\it Algebr. Geom. Topology}, {\bf 19} (2019), 477--489).
 \emp

\bmp \textbf{14.16.} Following Yu.\,I.\,Merzlyakov we define the {\it
width\/} of a verbal subgroup $V(G)$ of a group $G$ with respect
to the set of words $V$ as the smallest $m\in {\Bbb N}\cup \{
\infty\}$ such that every element of
 $V(G)$
can be written as a product of $\leq m$ values of words from $V\cup
V^{-1}$. It is known that the width of any verbal subgroup of a
finitely generated group of polynomial growth is finite. Is this
statement true for finitely generated groups of intermediate growth?
\hfill \raisebox{-1ex}{\sl V.\,G.\,Bardakov}
 \emp

\bmp
\textbf{14.17.}
Let ${\rm IMA}(G)$ denote the subgroup of the automorphism group
 ${\rm
Aut\,}G$ consisting of all automorphisms that act trivially on the
second derived quotient $G/G''$. Find generators and defining
relations for
 ${\rm IMA}(F_n)$, where $F_n$ is a free group of rank $n\geq 3$.

\hfill \raisebox{-1ex}{\sl V.\,G.\,Bardakov}
 \emp

\bmp \textbf{14.18.} We say that a family of groups
 ${\cal D}$ {\it
discriminates\/} a group $G$ if for any finite subset
$\{a_1,\ldots ,a_n\}\subseteq G\setminus\{ 1\}$ there exists a
group $D\in {\cal D}$ and a homomorphism $\f :G\rightarrow D$ such that
$a_j\f \ne 1$ for all $j=1,\ldots ,n$. Is every finitely generated
group acting freely on some $\Lambda$-tree discriminated by
torsion-free hyperbolic groups?

\hfill \raisebox{-1ex}{\sl G.\,Baumslag, A.\,G.\,Myasnikov,
V.\,N.\,Remeslennikov}

 \emp

\bmp \textbf{\zv 14.19.}
 We say that a group $G$ has the {\it Noetherian
Equation Property\/} if every system of equations over $G$ in finitely many variables
is equivalent to some finite part of it. Does an arbitrary hyper\-bolic
group have the Noetherian Equation Property?

\makebox[15pt][r]{}{\it Comment of 2013:}
the conjecture holds for torsion-free hyperbolic groups (Z.\,Sela,
{\it Proc. London Math. Soc.}, {\bf 99}, no.\,1 (2009), 217--273) and for the larger class of toral
relatively hyperbolic groups (D.\,Groves, {\it J. Geom. Topol.}, {\bf 9} (2005), 2319--2358).

\hfill \raisebox{-1ex}{\sl G.\,Baumslag, A.\,G.\,Myasnikov,
V.\,N.\,Remeslennikov}

\ul

\otv Yes, it does (R.\,Weidmann, C.\,Reinfeldt, {\it Ann. Math. Blaise Pascal}, {\bf 26}, no.\,2 (2019), 125--214).
 \emp

\bmp
\textbf{\zv 14.20.}
Does a free product of two groups have the
Noetherian Equation Property if this property is enjoyed by the factors?

\hfill \raisebox{-1ex}{\sl G.\,Baumslag, A.\,G.\,Myasnikov,
V.\,N.\,Remeslennikov}

\ul

\otv Yes, it does (Z.\,Sela, {\it Preprint}, 2010, \url{https://arxiv.org/abs/1012.0044}).
Note also that a free product of arbitrarily many equationally Noetherian groups need not be equationally Noetherian (D.\,Groves, M.\,Hull, {\it Trans. Amer. Math. Soc.}, {\bf 372} (2019), 7141--7190).
 \emp

\bmp
\textbf{14.21.}
Does a free pro-$p\hs$-group have the Noetherian Equation Property?

\hfill \raisebox{-1ex}{\sl G.\,Baumslag, A.\,G.\,Myasnikov,
V.\,N.\,Remeslennikov}

 \emp

\bmp \textbf{\zv 14.22.}
Prove that any irreducible system of equations
$S(x_1, \ldots,x_n) = 1$ with coefficients in a torsion-free
linear group $G$ is equivalent over $G$ to a finite system $T(x_1,
\ldots,x_n) = 1$ satisfying an analogue of Hilbert's
Nullstellensatz, i.~e. $Rad_G(T) = \sqrt T$. This is true if $G$
is a free group (O.\,Kharlampovich, A.\,Myasnikov, {\it
J.\,Algebra}, \textbf{200} (1998), 472--570).

\makebox[15pt][r]{}Here both $S$ and $T$ are regarded as subsets of
$G[X] = G \ast F(X)$, a free product of $G$ and a free group on $X =
\{x_1, \ldots, x_n \}$. By definition, $ Rad_G(T) = \{ w(x_1,
\ldots,x_n) \in G[X] \mid w(g_1, \ldots, g_n) = 1$ for any solution
$g_1, \ldots, g_n \in G$ of the system $T(X) = 1\}$, and $\sqrt T$ is
the minimal normal isolated subgroup of $G[X]$ containing~$T$. \hfill
\raisebox{-1ex}{\sl G.\,Baumslag, A.\,G.\,Myasnikov,
V.\,N.\,Remeslennikov}

\ul

\otv This is not true in general (AI solution confirmed by A.\,G.\,Myasnikov: \  Achyuth Jayadevan,  {\it Preprint of 9 September 2026},  \url{https://zenodo.org/records/22722075}; \ unverified AI solution: M.\,Kinyon, J.\,Urban, {\it Preprint of 17 September 2026},  \url{https://kourovkanotebookorg.wordpress.com/wp-content/uploads/2026/09/kourovka-mathematics.pdf}).
 \emp 

\bmp \textbf{\zv 14.23.}
Let $F_n$ be a free group with basis $\{
x_1,\ldots ,x_n \}$, and let $|\cdot |$ be the length function
with respect to this basis. For $\alpha\in {\rm Aut\,}F_n$ we put
$|\! |\alpha |\! |=\max \{
 |\alpha (x_1)|,\ldots ,|\alpha (x_n)| \}$. Is it true that there is a
recursive function $f: {\Bbb N} \rightarrow {\Bbb N}$ with the following
property: for any $\alpha\in {\rm Aut}\,F_n$ there is a basis $\{
y_1,\ldots ,y_k \}$ of ${\rm Fix}(\alpha )=\{ x\mid \alpha
(x)=x\}$ such that $|y_i|\leq f(|\! |\alpha |\! |)$ for all
$i=1,\ldots ,k$?
\hfill \raisebox{-1ex}{\sl O.\,V.\,Bogopolski}

\ul

\otv No, this is not true if the function is required to work simultaneously for all ranks (O.\,Kharlampovich, {\it Preprint}, 2026, \url{https://kourovkanotebookorg.wordpress.com/wp-content/uploads/2026/08/kourovka_14_23_uniform_counterexample.pdf}), but this is true for each given rank (O.\,Kharlampovich, {\it Preprint}, 2026, \url{https://kourovkanotebookorg.wordpress.com/wp-content/uploads/2026/08/kourovka_14_23_17_32_note_with_problem_authors.pdf}).
 \emp

\bmp \textbf{14.24.} Let ${\rm Aut\,}F_n$ be the automorphism group
of a free group of rank $n$ with norm $|\! |\cdot |\! |$ as in
14.23. Does there exist a recursive function
 $f: {\Bbb N} \times {\Bbb N} \rightarrow {\Bbb N}$ with the following property:
 for any two conjugate elements
$\alpha,\beta\in {\rm Aut\,}F_n$ there is an element $\gamma \in {\rm
Aut\,}F_n$ such that $\gamma ^{-1}\alpha \gamma =\beta$ and $|\!
|\gamma |\! |\leq f(|\! |\alpha |\! |,|\! |\beta |\! |)$? \hfill
\raisebox{-1ex}{\sl O.\,V.\,Bogopolski}

 \emp

 \bmp \textbf{\ai 14.26.} 
  A quasivariety ${\cal M}$ is {\it closed under
direct ${\Bbb Z}$-wreath products\/} if the direct wreath product
$G\wr {\Bbb Z}$ belongs to ${\cal M}$ for every $G\in {\cal M}$ (here
${\Bbb Z}$ is an infinite cyclic group). Is the quasivariety
generated by the class of all nilpotent torsion-free groups closed
under direct ${\Bbb Z}$-wreath products? \hfill \raisebox{-1ex}{\sl
A.\,I.\,Budkin}

\ul

\oai Yes, it is (unverified AI solution:  M.\,Kinyon, J.\,Urban, {\it Preprint}, 2026,

\url{https://kourovkanotebookorg.wordpress.com/wp-content/uploads/2026/09/kourovka-mathematics.pdf}).
 \emp 

\bmp \textbf{14.28.} Let ${\frak F}$ be a soluble Fitting formation of
finite groups with Kegel's property, that is, ${\frak F}$ contains
every finite group of the form $G=AB=BC=CA$ if $A,B,C$ are in~${\frak
F}$. Is ${\frak F}$ a saturated formation?

\makebox[15pt][r]{}{\it Comment of 2013:}
Some progress was made in (A.\,Ballester-Bolinches, L.\,M.\,Ez\-querro, {\it J.~Group Theory}, {\bf 8}, no.\,5 (2005), 605--611).
\hfill \raisebox{-1ex}{\sl
A.\,F.\,Vasiliev}
 \emp

\bmp \textbf{14.29.} Is there a soluble Fitting class of finite
groups ${\frak F}$ such that ${\frak F}$ is not a formation and $
A_{{\frak F}}\cap B_{{\frak F}}\subseteq G_{{\frak F}}$ for every
finite soluble group of the form $G=AB$?

\hfill \raisebox{-1ex}{\sl A.\,F.\,Vasiliev}
 \emp

\bmp \textbf{14.30.} Let ${\rm lFit}{\frak X}$ be the local Fitting
class generated by a set of groups ${\frak X}$ and let $\Psi (G)$
be the smallest normal subgroup of a finite group $G$ such that
${\rm lFit}(\Psi (G)\cap M)={\rm lFit}M$ for every
$M\mbox{$\lhd\lhd$}\,G$ (K.\,Doerk, P.\,Hauck, {\it Arch. Math.},
\textbf{35}, no.\,3 (1980), 218--227). We say that a Fitting class
${\frak F}$ is {\it saturated\/} if $\Psi (G)\in {\frak F}$
implies that $G\in {\frak F}$. Is it true that every non-empty
soluble saturated Fitting class is local? \hfill
\raisebox{-1ex}{\sl N.\,T.\,Vorob'\"ev}

 \emp

\bmp \textbf{14.31.} Is the lattice of Fitting subclasses of the Fitting
class generated by a finite soluble group finite? \hfill
\raisebox{-1ex}{\sl N.\,T.\,Vorob'\"ev}
 \emp

\bmp \textbf{14.35.} Is every finitely presented group of prime exponent
finite? \hfill \raisebox{-1ex}{\sl N.\,D.\,Gupta}
 \emp

\bmp \textbf{14.36.} A group $G$ is called a {\it $T$-group\/} if every
subnormal subgroup of $G$ is normal, while $G$ is said to be a {\it
$\,\overline{\! T}$-group\/} if all of its subgroups are $T$-groups.
Is it true that every non-periodic locally graded $\,\overline{\!
T}$-group must be abelian? \hfill \raisebox{-1ex}{\sl
F.\,\,de\,\,Giovanni}
 \emp

\bmp
\textbf{14.37.} Let $G(n)$ be one of the classical groups
(special, orthogonal, or symplectic) of $(n\times n)$-matrices
over an infinite field
 $K$ of non-zero characteristic, and $M(n)$ the space of
 all $(n\times
n)$-matrices over~$K$. The group
 $G(n)$ acts diagonally by conjugation on the space
$M(n)^m=\underbrace{M(n)\oplus\cdots\oplus M(n)}_{m}$. Find
generators of the algebra of invariants\linebreak
\vskip-3ex
 $K[M(n)^m]^{G(n)}$.

\makebox[15pt][r]{}In characteristic~$0$ they were found in
(C.\,Procesi, {\it Adv. Math.}, \textbf{19} (1976), 306--381). {\it
Comment of 2001:\/} in positive characteristic the problem is
solved for all cases excepting the orthogonal groups in
characteristic 2 and special orthogonal groups of even degree
(A.\,N.\,Zubkov, {\it Algebra and Logic}, \textbf{38}, no.\,5 (1999),
299--318). {\it Comment of 2009:\/} ...and for special orthogonal
groups of even degree over (infinite) fields of odd characteristic
(A.\,A.\,Lopatin, {\it J.~Algebra}, \textbf{321} (2009),
1079--1106).\hfill\raisebox{-1ex}{\sl A.\,N.\,Zubkov}

 \emp

\bmp \textbf{14.38.} For every pro-$p\hs$-group $G$ of $(2\times
2)$-matrices for $p\ne 2$ an analogue of the Tits Alternative
holds: either $G$ is soluble, or the variety of pro-$p\hs$-groups
generated by $G$ contains the group {\small $\overline{\left< \left(
\!\!\begin{array}{cc}1&t\cr 0&1
\end{array}\!\!\right),\, \left(\!\!
\begin{array}{cc}1&0\cr t&1\end{array}\!\!\right) \right>
}$}\,$\leq SL_2({\Bbb F} _p[[t]])$ (A.\,N.\,Zubkov, {\it Algebra
and Logic}, \textbf{29}, no.\,4 (1990), 287--301). Is the same result
true for matrices of size
 $\geq 3$ for $p\ne 2$?
\hfill \raisebox{-1ex}{\sl A.\,N.\,Zubkov}
 \emp

\bmp \textbf{14.39.} Let $F(V)$ be a free group of some variety of
pro-$p\hs$-groups $V$. Is there a uniform bound for the exponents
of periodic elements in
 $F(V)$? This is true if the variety
 $V$ is metabelian (A.\,N.\,Zubkov, {\it Dr. Sci.
 Diss.},
Omsk, 1997 (Russian)). \hfill \raisebox{-1ex}{\sl A.\,N.\,Zubkov}
 \emp

\bmp \textbf{14.40.} Is a free pro-$p\hs$-group representable as an
abstract group by matrices over a~commutative-associative ring with
$1$? \hfill \raisebox{-1ex}{\sl A.\,N.\,Zubkov}
 \emp

\bmp \textbf{14.41.} A group $G$ is said to be {\it para-free\/} if
all factors $\gamma _i(G)/\gamma _{i+1}(G)$ of its lower central
series are isomorphic to the corresponding lower central factors
of some free group, and $\bigcap_{i=1}^{\infty}\gamma _i(G)\,=1$.
Is an arbitrary para-free group representable by matrices over a
commutative-associative ring with $1$? \hfill \raisebox{-1ex}{\sl
A.\,N.\,Zubkov}
 \emp

\bmp \textbf{14.42.} Is a free pro-$p$-group representable by
matrices over an associative-com\-mu\-ta\-tive profinite ring
with~1?

\makebox[15pt][r]{}A negative answer is equivalent to the fact
that every linear pro-$p$-group satisfies a non-trivial
pro-$p$-identity. This is known to be true in dimension $2$
for $p\ne 2$ (A.\,N.\,Zubkov, {\it Siberian Math.~J.}, \textbf{28},
no.\,5 (1987), 742--747). It is also proved that a 2-dimensional linear pro-$2$-group in characteristic 2 satisfies a non-trivial pro-$2$-identity (D.\,E.-C.\,Ben-Ezra, E.\,Zelmanov, \textit{Trans. Amer. Math. Soc.}, \textbf{374}, no.\,6 (2021), 4093--4128).
\hfill \raisebox{-1ex}{\sl A.\,N.\,Zubkov,
V.\,N.\,Remeslennikov}
 \emp

\bmp \textbf{14.43.} Suppose that a finite group $G$ has the form
$G=AB $, where $A$ and $B$ are nilpotent subgroups of classes
$\alpha$ and $\beta$ respectively; then $G$ is soluble
(O.\,H.\,Kegel, H.\,Wielandt, 1961). Although the derived length
${\rm dl}(G)$ of $G$ need not be bounded by $\alpha +\beta$ (see
Archive, 5.17), can one bound ${\rm dl}(G)$ by a (linear) function of
$\alpha$ and $\beta$?

\hfill \raisebox{-1ex}{\sl L.\,S.\,Kazarin}

 \emp

\bmp \textbf{14.44.} Let $k(X)$ denote the number of conjugacy
classes of a finite group $X$. Suppose that a finite group $G=AB$
is a product of two subgroups $A,B$ of coprime orders. Is it true
that $k(AB) \leq k(A)k(B)$?

\makebox[15pt][r]{}Note that one cannot drop the
coprimeness condition and the answer is positive
if one of the subgroups is normal, see Archive,~11.43.
\hfill
\raisebox{-1ex}{\sl L.\,S.\,Kazarin, J.\,Sangroniz}

\emp

\bmp \textbf{14.45.} Does there exist a (non-abelian simple) linearly
right-orderable group all of whose proper subgroups are cyclic?
\hfill \raisebox{-1ex}{\sl U.\,E.\,Kaljulaid}
 \emp

\bmp

\textbf{14.46.} A finite group $G$ is said to be {\it almost
simple\/} if
 $T\leq G\leq {\rm Aut}(T)$ for some nonabelian simple group
$T$. By definition, a {\it finite linear space\/} consists of a
set $V$ of points, together with a collection of $k$-element
subsets of $V$, called {\it lines\/} ($k\geq 3$), such that every
pair of points is contained in exactly one line. Classify the
finite linear spaces which admit a line-transitive almost simple
subgroup $G$ of automorphisms which acts transitively on points.

\makebox[15pt][r]{}A.\,Camina, P.\,Neumann and C.\,E.\,Praeger have
solved this problem in the case where $T$ is an alternating group. In
(A.\,Camina, C.\,E.\,Praeger, {\it Aequat. Math.}, \textbf{61} (2001),
221--232) it is shown that a line-transitive group of automorphisms
of a finite linear space which is point-{\it quasiprimitive\/}
(i.\,e. all of whose non-trivial normal subgroups are
point-transitive) is almost simple or affine. \

 \makebox[15pt][r]{}{\it Remarks of
2005:} \
The cases of the following groups, but not almost simple groups with this
socle, were dealt with: $PSU(3,q)$ (W.\,Liu,
{\it Linear Algebra Appl.}, \textbf{374} (2003), 291--305);
$PSL(2,q)$
(W.\,Liu, {\it J.~Combin. Theory (A)}, \textbf{103} (2003),
209--222); $Sz(q)$ (W.\,Liu, {\it Discrete Math.}, \textbf{269} (2003),
181--190); $Ree(q)$ (W.\,Liu, {\it Europ. J.~Combin.}, \textbf{25} (2004),
311--325); sporadic simple groups were done in
(A.\,R.\,Ca\-mi\-na, F.\,Spiezia, {\it J.~Combin. Des.}, \textbf{8}
(2000), 353--362).
\ {\it Remark of 2009:} It was shown (N.\,Gill,
{\it Trans. Amer. Math. Soc.}, {\bf 368} (2016), 3017--3057)
that for non-desarguesian
projective planes a point-quasiprimitive group would be affine.
\ {\it Remark of 2013:}
The problem is solved for large-dimensional classical groups
(A.\,Camina, N.\,Gill, A.\,E.\,Zalesski, {\it Bull. Belg. Math. Soc. Simon Stevin}, {\bf 15}, no.\,4 (2008), 705--731).
\hfill \raisebox{-1ex}{\sl A.\,Camina, C.\,E.\,Praeger}

\emp

\bmp
\textbf{14.47.}
 Is the lattice of all soluble Fitting classes of finite groups modular?

\hfill \raisebox{-1ex}{\sl S.\,F.\,Kamornikov, A.\,N.\,Skiba}
 \emp

\bmp \textbf{14.51.} (Well-known problems). Is there a finite basis
for the identities of any

\makebox[25pt][r]{a)} abelian-by-nilpotent group?

\makebox[25pt][r]{b)} abelian-by-finite group?

\makebox[25pt][r]{c)} abelian-by-(finite nilpotent) group? \hfill
\raisebox{0ex}{\sl A.\,N.\,Krasil'nikov}

 \emp

 \bmp \textbf{14.53.} {\it Conjecture:\/} Let $G$ be a profinite group
such that the set of solutions of the equation $x^n=1$ has
positive Haar measure. Then $G$ has an open subgroup $H$ and an
element $t$ such that all elements of the coset $tH$ have order
dividing $n$.

\makebox[15pt][r]{}This is true in the case $n=2$. It would be
interesting to see whether similar results hold for profinite
groups in which the set of solutions of some equation has positive
measure.

\makebox[15pt][r]{}{\it Comment of 2021}:
This is also proved for $n=3$ (A.\,Abdollahi, M.\,S.\,Malekan, \textit{Adv. Group Theory Appl.}, \textbf{13} (2022), 71--81). {\it Comment of 2025}:
Significant progress was obtained in (S.\,Kionke, N.\,Otmen, T.\,Toti, M.\,Vannacci,
T.\,Weigel, {\it Preprint}, 2025,
 \url{https://arxiv.org/pdf/2507.19086}).
\hfill \raisebox{-1ex}{\sl L.\,Levai, L.\,Pyber}
 \emp

\bmp \textbf{\zv 14.54.}
Let $k(G)$ denote the number of conjugacy
classes of a finite group~$G$. Is it true that $k(G)\leq |N|$ for
some nilpotent subgroup
 $N$ of~$G$? This is true if $G$ is simple; besides, always $k(G)\leq
|S|$ for
some soluble subgroup $S\leq G$.

\hfill \raisebox{-1ex}{\sl M.\,W.\,Liebeck, L.\,Pyber}

\ul

\otv \! No, not always (S.\,Tertooy, {\it Preprint}, 2026,

\url{https://arxiv.org/abs/2606.21404}).
 \emp

\bmp
\textbf{14.55.} b)
 Prove that the Nottingham group $J=N({\Bbb Z}/p{\Bbb Z})$ (as
defined in Archive, 12.24) is finitely presented
 for $p=2$. \hfill \raisebox{-1ex}{\sl
C.\,R.\,Leedham-Green}
 \emp

\bmp
\textbf{14.56.}
 Prove that if $G$ is an infinite pro-$p\hs$-group with
$G/\gamma_{2p+1}(G)$ isomorphic to
$J/\gamma_{2p+1}(J)$ then $G$ is isomorphic to $J$, where $J$ is the
Nottingham group.

\hfill \raisebox{-1ex}{\sl C.\,R.\,Leedham-Green}

 \emp

\bmp \textbf{14.57.} Describe the hereditarily just infinite
pro-$p\hs$-groups of finite width.

\makebox[15pt][r]{}Definitions: A pro-$p\hs$-group $G$ has finite
{\it width\/} $p^d$ if $|\gamma _i(G)/\gamma _{i+1}(G)|\leq p^d$
for all~$i$, and $G$ is {\it hereditarily just infinite\/} if $G$
is infinite and each of its open subgroups has no closed normal
subgroups of infinite index. \hfill \raisebox{-1ex}{\sl
C.\,R.\,Leedham-Green}

 \emp

\bmp \textbf{14.58.} a) Suppose that $A$ is a periodic group of regular
automorphisms of an abelian group. Is $A$ cyclic if $A$ has prime
exponent? \hfill \raisebox{-1ex}{\sl V.\,D.\,Mazurov}

 \emp

\bmp
\textbf{14.59.}
Suppose that $G$ is a triply transitive
group in which a stabilizer of two points contains no involutions,
and a stabilizer of three points is trivial.
Is it true that $G$ is similar to $PGL_2(P)$ in its natural action
on the projective line
 $P\cup \{ \infty \}$, for some field
 $P$ of characteristic~$2$? This is true under the condition that the
stabilizer of two points is periodic. \hfill \raisebox{-1ex}{\sl
V.\,D.\,Mazurov}
 \emp

\bmp \textbf{14.61.} Determine all pairs $({\cal S},G)$, where $\cal
S$ is a semipartial geometry and $G$ is an almost simple
flag-transitive group of automorphisms of $\cal S$. A~system of
points and lines $(P,B)$ is a {\it semipartial geometry\/} with
parameters $(\alpha,s,t,\mu)$ if every point belongs to exactly
 $t+1$
lines (two different points belong to at most one line); every
line contains exactly $s+1$ points; for any anti-flag $(a,l)\in
(P,B)$ the number of lines containing
 $a$ and intersecting
 $l$ is either $0$ or $\alpha$; and for any
non-collinear points
 $a,b$ there are exactly
$\mu$ points collinear with $a$ and with $b$. \hfill
\raisebox{-1ex}{\sl A.\,A.\,Makhn\"ev}

 \emp

\bmp \textbf{14.64.} (M.\,F.\,Newman). Classify the finite $5$-groups of
maximal class; computer calculations suggest some conjectures
(M.\,F.\,Newman, in: {\it Groups--Canberra, 1989 $($Lecture Notes
Math.}, \textbf{1456}), Springer, Berlin, 1990, 49--62).

\makebox[15pt][r]{}{\it M.\,F.\,Newman's comment of 2025}:
Recent work on related conjectures can be found in (A.\,Cant,  H.\,Dietrich, B.\,Eick, T.\,Moede, {\it J.~Algebra}, {\bf 604} (2022), 429--450).

\hfill
\raisebox{-1ex} {\sl A.\,Moret\'o}

 \emp

\bmp
\textbf{14.65.}
(Well-known problem).
For a finite group $G$ let $\rho(G)$ denote the
set of prime numbers dividing the order of some conjugacy class, and
$\sigma(G)$ the maximum number of primes dividing the order of some
conjugacy class. Is it true that $|\rho(G)|\leq3\sigma(G)$?

\makebox[15pt][r]{}A possible linear bound for $|\rho(G)|$ in
terms of $\sigma(G)$ cannot be better than $3\sigma(G)$, since
there is a family of groups $\{G_n\}$ such that
$\lim_{n\to\infty}|\rho(G_n)|/\sigma(G_n)=3$ (C.\,Casolo,
S.\,Dolfi, {\it Rend. Sem. Mat. Univ. Padova}, \textbf{96} (1996),
121--130).

\makebox[15pt][r]{}{\it Comment of 2009:\/}
A quadratic bound was obtained in (A.\,Moret\'o, {\it Int. Math.
Res. Not.}, \textbf{2005}, no.\,54, 3375--3383), and a linear bound
in (C.\,Casolo, S.\,Dolfi, {\it J.~Group Theory}, \textbf{10}, no.\,5
(2007), 571--583). \hfill \raisebox{-1ex} {\sl A.\,Moret\'o}

 \emp

\bmp
\textbf{\zv 14.67.}
Suppose that $a$ is a non-trivial element
of a finite
 group $G$ such that
 $|C_{G}(a)| \geq |C_{G}(x)|$ for every non-trivial element
$x\in G$, and let
$H$ be a nilpotent
subgroup of~$G$ which is normalized by $C_{G}(a)$.
Is it true that $H \leq C_{G}(a)$?
This is true if $H$ is abelian, or if $H$ is a
$p'$-group for some prime number
$p\in \pi(Z(C_G(a)))$.

\hfill \raisebox{-1ex}{\sl I.\,T.\,Mukhamet'yanov, A.\,N.\,Fomin}

\ul

\otv No, not always (AI solution confirmed by K.\,Muliarchyk: \  Achyuth Jayadevan,  {\it Preprint}, 2026, \url{https://zenodo.org/records/22861469}).
 \emp

\bmp \textbf{14.68.} (Well-known problem). Suppose that $F$ is an
automorphism of order $2$ of the polynomial ring $R_n= {\Bbb C}
[x_1, \ldots ,x_n]$, $n> 2$. Does there exist an automorphism $G$
of~$R_n$ such that $G^{-1}FG$ is a linear automorphism?

\makebox[15pt][r]{}This is true for $n=2$.
\hfill \raisebox{0ex}{\sl
M.\,V.\,Neshchadim}
 \emp

\bmp \textbf{14.69.} For every finite simple group find the minimum of
the number of generating involutions satisfying an additional
condition, in each of the following cases.

\makebox[25pt][r]{\zva a)}
The product of the generating involutions
equals~$1$.

\makebox[25pt][r]{b)} (Malle--Saxl--Weigel). All generating
involutions are conjugate.

\makebox[25pt][r]{c)} (Malle--Saxl--Weigel). The conditions a) and
b) are simultaneously satisfied.

{\it Editors' comment:}
partial progress was obtained in
(J.\,Ward, PhD Thesis, QMW, 2009, \url{http://www.maths.qmul.ac.uk/~raw/JWardPhD.pdf}).

\makebox[25pt][r]{d)}
All generating
involutions are conjugate and two of them commute.

\hfill \raisebox{-1ex}{\sl Ya.\,N.\,Nuzhin}

\ul

\otv a) These numbers are found (R.\,I.\,Gvozdev, Ya.\,N.\,Nuzhin, {\it Siberian MAth.~J.}, {\bf 64}, no.\,6 (2023), 1160--1171; \
M.\,A.\,Vsemirnov, R.\,I.\,Gvozdev, Ya.\,N.\,Nuzhin, {\it
Abstracts of Int. Conf. Mal'cev Meeting 2023}, Novosibirsk, 2023, p.\,148).
\emp

\bmp \textbf{14.70.} A group $G$ is called {\it $n$-Engel\/} if it
satisfies the identity $[x,y,\ldots ,y]=1$, where $y$ is taken $n$
times. Are there non-nilpotent finitely generated $n$-Engel
groups?

{\it Comment of 2017}:
For $n=2,3,4$, finitely generated
 $n$-Engel
groups are nilpotent. A course of lectures devoted to constructing non-nilpotent finitely generated $n$-Engel
groups for sufficiently large $n$ is given by E.\,Rips and is available on the Internet.

\hfill \raisebox{-1ex}{\sl B.\,I.\,Plotkin}

 \emp

\bmp \textbf{\zv 14.72.}
Let $X$ be a regular algebraic variety over a
field of arbitrary characteristic, and $G$ a finite cyclic group
of automorphisms of $X$. Suppose that the fixed point variety
$X^G$ of $G$ is a regular hypersurface of $X$ (of
codimension~$1$). Is the quotient variety $X/G$ regular? \hfill
\raisebox{-1ex}{\sl K.\,N.\,Ponomar\"ev}

\ul

\otv No, not always (AI solution confirmed by K.\,N.\,Ponomar\"ev: \  M.\,Kinyon, J.\,Urban, {\it Preprint}, 2026,
\url{https://kourovkanotebookorg.wordpress.com/wp-content/uploads/2026/09/kourovka-mathematics.pdf}).
 \emp 

\bmp
\textbf{14.73.} {\it Conjecture:\/} There is a function $f$ on
the natural numbers such that, if $\Gamma$ is a finite,
vertex-transitive, locally-quasiprimitive graph of valency $v$,
then the number of automorphisms fixing a given vertex is at most
$f(v)$. (By definition, a vertex-transitive graph $\Gamma$ is {\it
locally-quasiprimitive\/} if the stabilizer in $Aut(\Gamma)$ of a
vertex $\alpha$ is quasiprimitive (see 14.46) in its action on the
set of vertices adjacent to $\alpha$.)

\makebox[15pt][r]{}To prove the conjecture above one need only
consider the case where $Aut(\Gamma)$ has the property that every
non-trivial normal subgroup has at most two orbits on vertices
(C.\,E.\,Praeger, {\it Ars Combin.}, \textbf{19\,A} (1985),
149--163). The analogous conjecture for finite, vertex-transitive,
locally-primitive graphs was made by R.\,Weiss in 1978 and is
still open. For non-bipartite graphs, there is a ``reduction'' of
Weiss' conjecture to the case where the automorphism group is
almost simple (see 14.46 for definition) (M.\,Conder, C.\,H.\,Li,
C.\,E.\,Praeger, {\it Proc. Edinburgh Math. Soc. (2)}, \textbf{43},
no.\,1 (2000), 129--138).

\makebox[15pt][r]{}{\it Comment of 2013:}
These conjectures have been proved in the case where there is an upper bound on the degree of any alternating group occurring as a quotient of a subgroup (C.\,E.\,Praeger, L.\,Pyber, P.\,Spiga, E.\,Szab\'{o},
{\it Proc. Amer. Math. Soc.}, {\bf 140}, no.\,7 (2012), 2307--2318). {\it Comment of 2021:}
This has been proved in the case where the group induced on the neighbourhood of a vertex has an abelian regular normal subgroup (P.\,Spiga, {\it Bull. London Math. Soc.}, {\bf 48}, no.\,1 (2016), 12--18.)
\hfill \raisebox{-1ex}{\sl
C.\,E.\,Praeger}
 \emp

\bmp \textbf{14.74.} Let $k(G)$ denote the number of conjugacy classes
of a finite group~$G$. Is it true that $k(G)\leq k(P_1)\cdots
k(P_s)$, where $P_1,\ldots ,P_s$ are Sylow subgroups of $G$ such that
$|G|=|P_1|\cdots |P_s|$? \hfill \raisebox{-1ex}{\sl L.\,Pyber}

 \emp

\bmp \textbf{14.75.} Suppose that ${\frak G}=\{ G_1,\,G_2,\ldots \}$ is
a family of finite 2-generated groups which generates the variety of
all groups. Is it true that a free group of rank 2 is residually in
${\frak G}$? \hfill \raisebox{-1ex}{\sl L.\,Pyber}
 \emp

\bmp
\textbf{14.76.} Does there exist an absolute constant $c$ such
that any finite $p\hs$-group $P$ has an abelian section $A$
satisfying $|A|^c>|P|$?

\makebox[15pt][r]{}By a result of A.\,Yu.\,Olshanskii ({\it
Math. Notes}, \textbf{23} (1978), 183--185) we cannot require $A$ to
be a subgroup. By a result of J.\,G.\,Thompson ({\it J.~Algebra},
\textbf{13} (1969), 149--151) the existence of such a section $A$
would imply the existence of a class $2$ subgroup $H$ of $P$ with
$|H|^c
>|P|$.
 \hfill \raisebox{-1ex}{\sl L.\,Pyber}

\emp

\bmp \textbf{14.78.} Suppose that ${\frak H}\subseteq {\frak
F}_{1}\subseteq {\frak M}$ \ and \ ${\frak \ H} \subseteq {\frak
F}_{2}\subseteq {\frak M}$ where ${\frak H}$ and ${\frak M}$ are
local formations of finite groups and ${\frak F}_{1}$ is a
complement for ${\frak F}_{2}$ in the lattice of all formations
between ${\frak H}$ and $ {\frak M}$. Is it true that ${\frak
F}_{1}$ and ${\frak F}_{2}$ are local formations? This is true in
the case ${\frak H}=(1)$. \hfill \raisebox{-1ex}{\sl
A.\,N.\,Skiba}

 \emp

\bmp \textbf{14.79.} Suppose that $ {\frak F}={\frak MH}={\rm lFit}( G)$
is a soluble one-generated local Fitting class of finite groups where
${\frak H}$ and ${\frak M}$ are Fitting classes and ${\frak M}\neq
{\frak F}$. Is ${\frak H}$ a local Fitting class? \hfill
\raisebox{-1ex}{\sl A.\,N.\,Skiba}
 \emp

\bmp
\textbf{\zv 14.81.}
Prove that the formation generated by a
finite group has only finitely many
 $S_{n}$-closed subformations.
\hfill \raisebox{-1ex}{\sl A.\,N.\,Skiba, L.\,A.\,Shemetkov}

\ul

\otv A counterexample is constructed (C.\,Lin, {\it Preprint}, 2026, \url{https://arxiv.org/abs/2610.37462}).
 \emp

\bmp \textbf{14.83.} We say that an infinite simple group $G$ is a {\it
monster of the third kind\/} if for every non-trivial elements
 $a,b$, of which at least one is not an involution, there are
infinitely many elements
$g\in G$ such that
$\left< a,b^g\right> =G$. (Compare with V.\,P.\,Shunkov's
definitions in Archive,~6.63, 6.64.) Is it true that every simple
quasi-Chernikov group
is a monster of the third kind? This is true for
quasi-finite groups
 (A.\,I.\,Sozutov, {\it Algebra and Logic}, \textbf{36}, no.\,5 (1997),
336--348). (We say that a non-$\s$ group is {\it quasi-$\s$\/} if
all of its proper subgroups have the property $\s$.)
 \hfill \raisebox{-1ex}{\sl A.\,I.\,Sozutov}

 \emp

\bmp \textbf{14.84.} An element $g$ of a (relatively) free group
$F_r({\frak M})$ of rank $r$ of a variety ${\frak M}$ is said to
be {\it primitive\/} if it can be included in a basis of
$F_r({\frak M})$.

\makebox[25pt][r]{a)} Do there exist a group $F_r({\frak M})$ and
a non-primitive element $h\in F_r({\frak M})$ such that for some
monomorphism $\alpha$ of $F_r({\frak M})$ the element $\alpha (h)$
is primitive?

\makebox[25pt][r]{b)} Do there exist a group $F_r({\frak M})$ and a
non-primitive element $h\in F_r({\frak M})$ such that for some $n>r$
the element $h$ is primitive in $F_n({\frak M})$? \hfill
\raisebox{-1ex}{\sl E.\,I.\,Timoshenko}

 \emp

\bmp
\textbf{\zv 14.85.}
Suppose that an endomorphism
 $\f$
 of a free metabelian group of rank $r$ takes every primitive
element to a primitive one.
 Is $\f$ necessarily an automorphism? This is true for
$r\leq 2$. \hfill \raisebox{-1ex}{\sl E.\,I.\,Timoshenko,
V.\,Shpilrain}

\ul

\otv Yes, it is (V.\,Ionin, A.\,Semidetnov,  {\it Preprint}, 2026, \url{https://arxiv.org/pdf/2608.29219}).
 \emp

\bmp \textbf{14.89.} (E.\,A.\,O'Brien, A.\,Shalev). Let $P$ be a
finite $p\hs$-group of order $p^m$ and let $m=2n+e$ with $e=0$
or~$1$. By a theorem of P.\,Hall the number of conjugacy classes
of $P$ has the form $n(p^2-1) + p^e + a(p^2-1)(p-1)$ for some
integer $a\geq 0$, which is called the {\it abundance\/} of~$P$.

\makebox[25pt][r]{a)} Is there a bound for the coclass of $P$
which depends only on $a$? Note that $a=0$ implies coclass $1$ and
that all known examples with $a=1$ have coclass $\le 3$. (The
group $P$ has {\it coclass\/} $r$ if $|P|=p^{c+r}$ where $c$ is
the nilpotency class of $P$.)

\makebox[25pt][r]{b)} Is there an element $s\in P$ such that
$|C_P(s)|\le p^{f(a)}$ for some $f(a)$ depending only on $a$?
We already know that we can take $f(0)=2$ and it seems that $f(1)=3$.

\makebox[15pt][r]{}Note that A.\,Jaikin-Zapirain ({\it J.~Group
Theory}, \textbf{3}, no.\,3 (2000), 225--231) has proved that $|P|\le
p^{f(p,a)}$ for some function $f$ of $p$ and $a$ only.

\hfill \raisebox{-1ex}{\sl G.\,Fern\'andez--Alcober}

 \emp

\bmp \textbf{14.90.} Let $P$ be a finite $p\hs$-group of
abundance~$a$ and nilpotency class~$c$. Does there exist an
integer $t=t(a)$ such that $\gamma_i(G)=\zeta _{c-i+1}(G)$ for
$i\geq t$? This holds for $a=0$ with $t=1$, since $P$ has maximal
class; it can be proved that for $a=1$ one can take $t=3$. \hfill
\raisebox{-1ex}{\sl G.\,Fern\'andez--Alcober}

 \emp

\bmp \textbf{14.91.} Let $p$ be a fixed prime. Do there exist finite
$p\hs$-groups of abundance $a$ for any $a\geq 0$? \hfill
\raisebox{-1ex}{\sl G.\,Fern\'andez--Alcober}

 \emp

\bmp \textbf{14.93.}
Let $N({\Bbb Z}/p{\Bbb Z})$ be the group defined
in Archive,~12.24 (the so-called ``Nottingham group'', or the ``Wild
group'').
 Find relations of $N({\Bbb Z}/p{\Bbb Z})$ as a pro-$p\hs$-group
(it has two generators, e.\,g., $x+x^2$ and $x/(1-x)$). {\it Comment
of 2009:} for $p>2$ see progress in (M.\,V.\,Ershov, {\it J. London Math. Soc.}, \textbf{71} (2005),
362--378).\hfill
\raisebox{-1ex}{\sl I.\,B.\,Fesenko}
 \emp

\bmp \textbf{14.94.} For each positive integer $r$ find the
$p\hs$-cohomological dimension ${\rm cd}_p(H_r)$ where $H_r$ is
the closed subgroup of $N({\Bbb Z}/p{\Bbb Z})$ consisting of
the series \ $x\left( 1 + \sum\limits_{i=1}^{\infty} a_ix^{p^ri}\right)
$, $a_i\in {\Bbb Z}/p{\Bbb Z}$. \hfill \raisebox{-1ex}{\sl
I.\,B.\,Fesenko}

 \emp

\bmp
\textbf{14.95.} (C.\,R.\,Leedham-Green, P.\,M.\,Neumann,
J.\,Wiegold). For a finite $p\hs$-group~$P$, denote by $c=c(P)$
its nilpotency class and by $b=b(P)$ its {\it breadth}, that is,
$p^b$ is the maximum size of a conjugacy class in $P$. {\it
Class-Breadth Problem:\/} Is it true that $c\le b+1$ if $p\ne 2$?

\makebox[15pt][r]{}So far, the best known bound is $c <
\frac{p}{p-1}b+1$ (C.\,R.\,Leedham-Green, P.\,M.\,Neu\-mann,
J.\,Wiegold, {\it J.~London Math.~Soc. (2)}, \textbf{1} (1969),
409--420). For $p=2$ for every $n\in{\Bbb N}$ there exists a
$2$-group $T_n$ such that $c(T_n)\geq b(T_n)+n$ (W.\,Felsch,
J.\,Neu\-b\"user, W.\,Plesken, {\it J.~London Math.~Soc. (2)},
\textbf{24} (1981), 113--122). \hfill \raisebox{-1ex}{\sl
A.\,Jaikin-Zapirain}
 \emp

\bmp \textbf{14.97.} Is it true that for any two different prime
numbers $p$ and $q$ there exists a non-primary periodic locally
soluble $\{ p,q\}$-group that can be represented as the product of
two of its $p\hs$-sub\-groups? \hfill \raisebox{-1ex}{\sl
N.\,S.\,Chernikov}

 \emp

\bmp \textbf{\zv 14.98.} We say that a metric space is a {\it $2$-end\/} one
(a {\it narrow\/} one), if it is quasiisometric to the real line
${\Bbb R}$
 (respectively, to a subset of ${\Bbb R}$). All other spaces are said
to be {\it wide}. Suppose that the Caley graph $\Gamma = \Gamma
(G, A)$ of a group $G$ with a finite set of generators~$A$ in the
natural metric contains a 2-end subset, and suppose that there is
 $\varepsilon > 0$ such that the complement in $\Gamma$ to the
$\varepsilon$-neighbourhood of any connected 2-end subset contains
exactly two wide connected components. Is it true that the group
$G$ in the word metric is quasiisometric to the Euclidean or
hyperbolic plane?

\hfill \raisebox{-1ex}{\sl V.\,A.\,Churkin}

\ul

\otv Yes, it is true (J.\,MacManus, {\it Preprint}, 2025, \url{https://arxiv.org/abs/2511.10759}).
 \emp

\bmp \textbf{14.99.}
a) A formation ${\frak F}$ of finite groups is
called {\it superradical\/} if it is $S_{n}$-closed and contains
every finite group of the form $G=AB$ where $A$ and $B$ are $
{\frak F}$-subnormal ${\frak \ F}$-sub\-groups. Find all superradical local formations.
 \hfill
\raisebox{-1ex}{\sl L.\,A.\,Shemetkov}
 \emp

\bmp \textbf{14.100.} Is it true that in a Shunkov group
(i.~e.~conjugately biprimitively finite group, see 6.59) having
infinitely many elements of finite order every element of prime
order is contained in some infinite locally finite subgroup? This
is true under the additional condition that any two conjugates of
this element generate a soluble subgroup (V.\,P.\,Shunkov, {\it
$M_p$-groups}, Moscow, Nauka, 1990 (Russian)). \hfill
\raisebox{-1ex}{\sl A.\,K.\,Shl\"epkin}

 \emp

\bmp \textbf{14.101.} A group $G$ is {\it saturated\/} with groups from
a class
 $X$ if every finite subgroup $K\leq G$
is contained in a subgroup
 $L\leq G$
 isomorphic to some group from~$X$. Is it true that
a periodic group saturated with finite simple groups of Lie type of
uniformly bounded ranks is itself a simple group of Lie type of
finite rank? \hfill \raisebox{-1ex}{\sl A.\,K.\,Shl\"epkin}

 \emp

\bmp
\textbf{\zv 14.102.}
(V.\,Lin). Let $B_n$ be the braid group
on
 $n$
strings, and let $n>4$.

\makebox[25pt][r]{a)} Does $B_n$ have any non-trivial
non-injective endomorphisms with non-cyclic images?

\makebox[25pt][r]{b)}
Is it true that every non-trivial endomorphism of the derived subgroup
 $[B_n,\,B_n]$ is an automorphism?
 \hfill \raisebox{-1ex}{\sl
V.\,Shpilrain}

\ul

\otv a) No, for $n\geqslant 5$ any non-injective homomorphism $B_n\to B_n$ has  cyclic image (F.\,Castel, {\it Geometric representations of the braid groups}, ({\it Ast\'erisque}, {\bf 378}), Paris, 2016, for $n\geq 6$, and
L.\,Chen, K.\,Kordek, D.\,Margalit, {\it Preprint}, 2019, \url{https://arxiv.org/abs/1910.00712} for $n=5$).

\otv b) Yes, it is true for $n\geq 5$ (S.\,Orevkov,
{\it Ann. Fac. Sci. Toulouse. Math.}, {\bf 33}, no.\,1 (2024), 105--121, extending the result for $n\geq 7$ in
K.\,Kordek, D.\,Margalit, {\it Bull. London Math. Soc.}, {\bf 54}, no.\,1 (2022), 95--111).
\emp

\raggedbottom
\newpage

\newpage

\pagestyle{myheadings} \markboth{15th Issue (2002)}{15th Issue
(2002)} \thispagestyle{headings} ~ \vspace{2ex}

\centerline {\Large \textbf{Problems from the 15th Issue (2002)}}
\phantomsection\label{15izd}
\vspace{4ex}

 \bmp \textbf{15.1.} (P.\,Longobardi,
M.\,Maj, A.\,H.\,Rhemtulla). Let $w = w(x_1, \ldots , x_n)$ be a
group word in $n$ variables $x_1, \dots , x_n$, and
 $V(w)$ the variety of groups defined by the law $w= 1$.
Let $V(w^*)$
 (respectively, $V(w^{\# })$) be the class of all groups $G$ in
 which for every $n$ infinite subsets $S_1,\ldots , S_n$
 there exist $s_i \in S_i$ such that
 $w(s_1,\ldots ,s_n)=1$ (respectively, $\left< s_1,\ldots ,s_n
\right> \in
 V(w)$).

\makebox[15pt][r]{}a) Is there some word $w$ and an infinite group
$G$ such that $G \in V(w^\#)$ but $G \notin V(w)$?

\makebox[15pt][r]{}b) Is there some word $w$ and an infinite group
$G$ such that $G \in V(w^*)$ but $G\notin V(w^{\# })$?

\makebox[15pt][r]{}The answer to both of these questions is likely to
be ``yes''. It is known that in a) \ $w$ cannot be any of several
words such as $x_1^n$, \ $[x_1,\ldots ,x_n]$, \ $[x_1,x_2]^2$, \
 $(x_1x_2)^3x_2^{-3}x_1^{-3}$, and $x_1^{a_1}\cdots
 x_n^{a_n}$ for any non-zero integers $a_1,\ldots ,a_n$.
\hfill \raisebox{-1ex}{\sl A.\,Abdollahi}

\emp

 \bmp \textbf{15.2.}
By a theorem of W.~Burnside, if $\chi \in {\rm Irr}(G)$ and $\chi(1)
> 1$, then there exists $x \in G$ such that $\chi(x) = 0$, that is,
only the linear characters are ``nonvanishing''. It is interesting to
consider the dual notion of {\it nonvanishing
 elements\/} of a finite group~$G$, that is, the elements $x \in G$ such
that $\chi(x) \ne 0$ for all $\chi \in
{\rm Irr}(G)$.

\makebox[15pt][r]{}a) It is proved (I.\,M.\,Isaacs, G.\,Navarro,
T.\,R.\,Wolf,
 {\it J.~Algebra}, \textbf{222}, no.\,2 (1999), 413--423) that if $G$ is
solvable and $x \in G$ is a nonvanishing element of odd order,
then $x$ belongs to the Fitting subgroup $\textbf{F}( G)$. Is this
true for elements of even order too? (M.\,Miyamoto showed in 2008
that every nontrivial abelian normal subgroup of a finite group
contains a nonvanishing element.)

\makebox[15pt][r]{}b) Which nonabelian simple groups have nonidentity
nonvanishing elements? (For example, ${\Bbb A}_7$ has.) \hfill
\raisebox{-1ex}{\sl I.\,M.\,Isaacs}

\emp

 \bmp \textbf{15.3.}
Let $\alpha$ and $\beta$ be faithful non-linear irreducible
characters of a finite group $G$. There are non-solvable groups
$G$ giving examples when the product $\alpha\beta$ is again an
irreducible character (for some of such $\alpha,\beta$). One
example is $G=SL_2(5)$ with two irreducible characters of degree
2. In (I.\,Zisser, {\it Israel J.~Math.}, \textbf{84}, no.\,1--2
 (1993), 147--151) it is proved that such an example exists in an
alternating group ${\Bbb A}_n$
if and only if $n$ is a square exceeding~$4$. But do solvable examples
exist? Evidence (but no proof) that they do not is given in
 (I.\,M.\,Isaacs, {\it J.~Algebra}, \textbf{223}, no.\,2 (2000), 630--646).
\hfill \raisebox{-1ex}{\sl I.\,M.\,Isaacs}

\emp

\bmp \textbf{15.8.}
b) Let $G$ be any Lie group (indeed, any separable continuous
group) made discrete: can $G$ act faithfully on a countable set?

\makebox[15pt][r]{}{\it Comment of 2021}:
an affirmative answer was obtained for the nilpotent case (N.\,Monod, {\it J.~Group Theory}, {\bf 25}, no.\,5 (2022), 851--865).
{\it Comment of 2025}:
this problem is reduced to the case of simple Lie groups; in particular, the
solvable case is settled (A.\,Conversano, N.\,Monod, {\it J.~Algebra}, {\bf 640} (2024), 106--116).
\hfill \raisebox{-1ex}{\sl P.~de~la~Harpe} \emp

 \bmp
 \textbf{15.9.}
 An automorphism $\f$ of the free group
$F_n$ on the free generators $x_1, x_2,\ldots , x_n$ is called
 {\it conjugating\/} if $x_i^{\f}=
t_i^ {-1} x_ {\pi (i)} t_i$, $i=1,\, 2,\ldots ,n$, for some
permutation $\pi \in {\Bbb S}_n$ and some elements $t_i\in F_n$.
The set of conjugating automorphisms fixing the product
 $x_1x_2\cdots x_n$ forms the braid group $B_n$. The group
$B_n$ is linear for any $n\geq 2$, while the group of all
automorphisms ${\rm Aut}\, F_n$ is not linear for $n\geq 3$. Is the
group of all conjugating automorphisms linear for
 $n\geq 3$?
\hfill \raisebox{-1ex}{\sl V.\,G.\,Bardakov}

\emp

 \bmp
 \textbf{15.11.}
 (M.\,Morigi). An automorphism of a group is called a {\it
power automorphism\/} if it leaves every subgroup invariant. Is every
finite abelian $p\hs$-group the group of all power automorphisms of
some group?
\hfill \raisebox{-1ex}{\sl V.\,G.\,Bardakov}

\emp

\bmp

\textbf{15.12.} Let $G$ be a group acting faithfully and
level-transitively by automorphisms on a rooted tree $\cal T$. For a
vertex $v$ of $\cal T$, the {\it rigid vertex stabilizer\/} at $v$
consists of those elements of $G$ whose support in $\cal T$ lies
entirely in the subtree ${\cal T}_v$ rooted at $v$. For a
non-negative integer $n$, the {\it $n$-th rigid level stabilizer\/}
is the subgroup of $G$ generated by all rigid vertex stabilizers
corresponding to the vertices at the level $n$ of the tree~$\cal T$.
The group $G$ is a {\it branch group\/} if all rigid level
stabilizers have finite index in $G$. For motivation, examples and
known results see (R. I. Grigorchuk, {\it in: New horizons in pro-$p$
groups}, Birkh\"auser, Boston, 2000, 121--179).

\makebox[15pt][r]{}Do there exist branch groups with Kazhdan's
$T$-property? (See A14.34.)

\hfill \raisebox{-1ex}{\sl L.\,Bartholdi, R.\,I.\,Grigorchuk,
Z.\,{\v S}uni\'k}

\emp

 \bmp

\textbf{15.13.} Do there exist finitely presented branch groups?

\hfill \raisebox{-1ex}{\sl L.\,Bartholdi, R.\,I.\,Grigorchuk,
Z.\,{\v S}uni\'k}
\emp

 \bmp
\textbf{\zv 15.14.}
b) Do there exist finitely generated branch groups
that are non-amenable and do not contain the free group $F_2$ on
two generators?

\hfill \raisebox{-1ex}{\sl L.\,Bartholdi, R.\,I.\,Grigorchuk,
Z.\,{\v S}uni\'k}

\ul

\otv b) Yes, there exist finitely generated non-amenable torsion branch groups (S.\,Kionke, E.\,Schesler, {\it Bull. London Math. Soc.}, {\bf 56}, no.\,2 (2024), 536--550).
\emp

 \bmp

\textbf{15.16.} Do there exist groups whose rate of growth is
$e^{\sqrt{n}}$

\makebox[15pt][r]{}a) in the class of finitely generated branch
groups?

\makebox[15pt][r]{}b) in the whole class of finitely generated
groups? This question is related to~9.9.

\hfill \raisebox{-1ex}{\sl L.\,Bartholdi, R.\,I.\,Grigorchuk, Z.\,{\v
S}uni\'k}

\emp

 \bmp \textbf{15.19.}
Let $p$ be a prime, and ${\cal F}_p$ the class of finitely
generated groups acting faithfully on a $p$-regu\-lar rooted
tree by finite automata. Any group in ${\cal F}_p$ is
residually-$p$ (residually in the class of finite $p$-groups) and
has word problem that is solvable in (at worst) exponential time.
There exist therefore groups that are residually-$p$, have a
solvable word problem, and do not belong to ${\cal F}_p$; though
no concrete example is known. For instance:\sm

\makebox[25pt][r]{a)} Is it true that some (or even all) the
groups given in (R.\,I.\,Grigorchuk, {\it Math. USSR--Sb.}, {\bf
54} (1986), 185--205) do not belong to ${\cal F}_p$ when the
sequence $\omega$ is computable, but not periodic?\sm

\makebox[25pt][r]{\zva b)}
Does ${\Bbb Z}\wr({\Bbb Z}\wr {\Bbb Z})$
belong to ${\cal F}_2$? (Here the wreath products are restricted.)\sm

\makebox[15pt][r]{}See (A.\,M.\,Brunner, S.\,Sidki, {\it J.~Algebra},
\textbf{257} (2002), 51--64) and (R.\,I.\,Grigorchuk,
V.\,V.\,Nekrashevich, V.\,I.\,Sushchanski\u{\i}, {\it in: Dynamical
systems, automata, and infinite groups. Proc. Steklov Inst. Math.},
\textbf{231} (2000),
 128--203).
 \hfill \raisebox{-1ex}{\sl L.\,Bartholdi, S.\,Sidki}

\ul

\otv b) Yes, it does (A.\,C.\,Dantas, J.\,R.\,Oliveira, T.\,M.\,G.\,Santos, {\it Preprint}, 2024, \url{https://arxiv.org/pdf/2405.16678}).
\emp

 \bmp
\textbf{15.20.} (B.\,Hartley). An infinite transitive permutation group
is said to be {\it barely transitive\/} if each of its proper
subgroups has only finite orbits. Can a locally finite barely
transitive group coincide with its derived subgroup?

\makebox[15pt][r]{}Note that there are no simple locally finite
barely transitive groups (B.\,Hartley, M.\,Kuzucuo\v{g}lu, {\it Proc.
Edinburgh Math. Soc.}, \textbf{40} (1997), 483--490), any locally finite
barely transitive group is a $p\hs$-group for some prime $p$, and if
the stabilizer of a point in a locally finite barely transitive group
$G$ is soluble of derived length $d$, then $G$ is soluble of derived
length bounded by a function of $d$ (V.\,V.\,Belyaev,
M.\,Kuzucuo\v{g}lu, {\it Algebra and Logic}, \textbf{42} (2003),
147--152).
 \hfill \raisebox{-1ex}{\sl V.\,V.\,Belyaev,
M.\,Kuzucuo\v{g}lu}

\emp

 \bmp \textbf{15.21.} (B.\,Hartley). Do there exist
torsion-free barely transitive groups?

\hfill \raisebox{-1ex}{\sl V.\,V.\,Belyaev}

\emp

 \bmp \textbf{15.22.}
A permutation group is said to be {\it finitary\/} if each of its
elements moves only finitely many points. Do there exist finitary
barely transitive groups? \hfill \raisebox{-1ex}{\sl
V.\,V.\,Belyaev}

\emp

 \bmp \textbf{15.23.} A transitive permutation group is
said to be {\it totally imprimitive\/} if every finite set of points
is contained in some finite block of the group. Do there exist
totally imprimitive barely transitive groups that are not locally
finite? \hfill \raisebox{-1ex}{\sl V.\,V.\,Belyaev}

\emp

 \bmp
\textbf{\zv 15.26.}
A {\it partition\/} of a group is a representation of it as
 a set-theoretic union of some of its proper subgroups ({\it
components\/}) that intersect pairwise trivially. Is it true that
every nontrivial partition
 of a finite $p$-group has an abelian component?

 \hfill
\raisebox{-1ex}{\sl Ya.\,G.\,Berkovich}

\ul

\otv Yes, it is true (AI solution confirmed by E.\,I.\,Khukhro: \ G.\,Cat, {\it Preprint}, 2026,

\url{https://kourovkanotebookorg.wordpress.com/wp-content/uploads/2026/08/kourovka_15-26-ai.pdf}).
\emp

 \bmp
 \textbf{15.28.}
Suppose that a finite $p\hs$-group $G$ is the product of two
subgroups: $G=AB$.

\makebox[15pt][r]{}a) Is the exponent of $G$ bounded in terms of the
exponents of $A$ and $B$?

\makebox[15pt][r]{}b) Is the exponent of $G$ bounded if $A$ and $B$
are groups of exponent $p$?

\hfill \raisebox{-1ex}{\sl Ya.\,G.\,Berkovich}

\emp

 \bmp
 \textbf{15.29.}
(A.\,Mann). The dihedral group of order $8$ is isomorphic to its own
automorphism group. Are there other non-trivial finite $p\hs$-groups
with this property?

\hfill \raisebox{-1ex}{\sl Ya.\,G.\,Berkovich}

\emp

\bmp
 \textbf{\ai 15.30.} 
Is it true that every finite abelian $p$-group is isomorphic to the
 Schur multiplier of some nonabelian finite $p$-group?
 \hfill \raisebox{-1ex}{\sl Ya.\,G.\,Berkovich}

\ul

\oai Yes, it is true (unverified AI solution: S.\,Cao, {\it Preprint}, 2026, \url{https://kourovkanotebookorg.wordpress.com/wp-content/uploads/2026/11/15_30.pdf}).
\emp 

\bmp
 \textbf{15.31.}
M.\,R.\,Vaughan-Lee and J.\,Wiegold ({\it Proc. R.~Soc. Edinburgh
Sect. A}, \textbf{95} (1983), 215--221) proved that if a finite
$p\hs$-group $G$ is generated by elements of breadth $\leq n$
(that is, having at most $p^n$ conjugates), then $G$ is nilpotent
of class $\leq n^2+1$; the bound for the class was later improved
by A.\,Mann ({\it J.~Group Theory}, \textbf{4}, no.\,3 (2001),
241--246) to $\leq n^2-n+1$. Is there a linear bound for the class
of $G$ in terms of $n$? \hfill \raisebox{-1ex}{\sl
Ya.\,G.\,Berkovich}

\emp

\bmp \textbf{\zv 15.32.}
Does there exist a function $f(k)$ (possibly depending also on $p$) such that if a finite $p\hs$-group $G$ of order $p^m$ with $m\geq f(k)$ has an automorphism of order
$p^{m-k}$, then $G$ possesses a cyclic subgroup of index $p^k$? \hfill
\raisebox{-1ex}{\sl Ya.\,G.\,Berkovich}

\ul

\otv Yes, one can take $f(k)=64k+2$ (AI solution confirmed by A.\,Caranti: \   Achyuth Jayadevan,  {\it Preprint}, 2026,
\url{https://zenodo.org/records/22910914}).
 \emp 

\bmp \textbf{15.36.} For a class $M$ of groups let
$L(M)$ be the class of groups $G$ such that the normal subgroup
$\left< a^G\right>$ generated by an element
 $a$ belongs to $M$ for any $a\in G$. Is it true that the class
$L(M)$ is finitely axiomatizable if $M$ is the quasivariety
generated by a finite group? \hfill \raisebox{-1ex}{\sl
A.\,I.\,Budkin}

\emp

\bmp \textbf{15.37.} Let $G$ be a group satisfying
the minimum condition for centralizers. Suppose that $X$ is a normal
subset of $G$ such that $[x,y,\ldots ,y]=1$ for any $x,y\in X$ with
$y$ repeated $f(x,y)$ times. Does $X$ generate a locally nilpotent
(whence hypercentral) subgroup?

\makebox[15pt][r]{}If the numbers $f(x,y)$ can be bounded, the
answer is affirmative (F.\,O.\,Wagner, {\it J.~Algebra}, {\bf
217}, no.\,2 (1999), 448--460).
\hfill \raisebox{-1ex}{\sl F.\,O.\,Wagner} \emp

\bmp \textbf{\zv 15.40.}
Let $N$ be a
nilpotent subgroup of a
 finite simple group~$G$. Is it true that there exists a subgroup
$N_1$ conjugate to $N$ such that $N\cap N_1=1$?

\makebox[15pt][r]{}The answer is known to be affirmative if
$N$ is a $p$-group. {\it Comment of 2013:}
an affirmative answer for alternating groups is obtained in (R.\,K.\,Kurmazov, {\it Siberian Math.~J.}, {\bf 54}, no.\,1 (2013), 73--77). 
\hfill \raisebox{-0ex}{\sl E.\,P.\,Vdovin}

\ul

\otv Yes it is true; moreover, for any two nilpotent subgroups $H,K$ of a (non-abelian) finite simple group $G$ there is $g\in G$ such that $H \cap K^g = 1$ (T.\,C.\,Burness, H.\,Y.\,Huang, {\it Preprint}, 2025, \url{https://arxiv.org/abs/2508.03479}). \emp

\bmp
 \textbf{15.41.} Let $R(m,p)$ denote the largest finite
$m$-gene\-ra\-tor group of prime exponent $p$.

\makebox[15pt][r]{}a) Can the nilpotency class of $R(m,p)$ be bounded
by a polynomial in $m$?
(This is true for $p=2,\,3,\,5,\,7$.)

\makebox[15pt][r]{}b) Can the nilpotency class of $R(m,p)$ be bounded
by a linear function in $m$? (This is true for $p=2,\,3,\,5$.)

\makebox[15pt][r]{}c) In particular, can the nilpotency class of
$R(m,7)$ be bounded by a linear function in $m$?

\makebox[15pt][r]{}My guess is ``no'' to the first two questions for
general $p$, but ``yes'' to the third. By contrast, a beautiful and
simple argument of Mike Newman
 shows that if $m\geq 2$ and $k\geq 2$ ($k\geq 3$ for $p=2$), then the
order of $R(m,p^k)$ is at least $p^{p^{.^{.^{.^{p^m}}}}}$, with
$p$ appearing $k$ times in the tower; see (M.\,Vaughan-Lee,
E.\,I.\,Zelmanov, {\it J.~Austral. Math. Soc. (A)}, \textbf{67},
no.\,2 (1999), 261--271). \hfill \raisebox{-1ex}{\sl
M.\,R.\,Vaughan-Lee}

\emp \vskip3ex

   \bmp \textbf{\zv 15.42.}
 Is it true that the group algebra $k[F]$ of R.\,Thompson's group
 $F$ (see 12.20) over a field $k$ satisfies the Ore condition, that
is, for
any $a,b \in
 k[F]$ there exist $u,v \in k[F]$ such that $au=bv$ and either $u$ or
 $v$ is nonzero? If the answer is negative, then $F$ is not amenable.

 \makebox[15pt][r]{}{\it Comment of 2026}: If the answer is positive, then $F$ is amenable (so
 problem 15.42 is equivalent to 12.20), because the ``if and only if'' statement holds for groups satisfying Kaplansky's zero-divisor conjecture (D.\,Kielak, \emph{Prepint}, 2026, \url{https://arxiv.org/pdf/1605.09133v2}), and $F$ does satisfy it.
\hfill \raisebox{-1ex}{\sl V.\,S.\,Guba}

\ul

\otv No, it is not true because R.\,Thompson's group
 $F$ is not amenable (AI solution confirmed by S.-H.\,Kim and Y.\,Lodha and by M.\,de~la~Salle: \  OpenAI, {\it Preprint}, 2026,
\url{https://kourovkanotebookorg.wordpress.com/wp-content/uploads/2026/10/12-20.pdf}; \ S.-H.\,Kim and Y.\,Lodha, {\it Preprint}, 2026, \url{https://arxiv.org/abs/2610.18360}; \  M.~de~la~Salle,  {\it Preprint}, 2026, \url{http://math.univ-lyon1.fr/homes-www/delasalle/thompson.pdf}).
\emp

\bmp
 \textbf{15.44.} a) Let $G$ be a reductive group over an algebraically
closed field $K$ of arbitrary characteristic. Let
 $X$ be an affine $G$-vari\-ety such that, for a fixed Borel
subgroup $B\leq G$, the coordinate algebra $K[X]$ as a $G$-module
is the union of an ascending chain of submodules each of whose
factors is an induced module ${\rm Ind}^G_BV$ of some
one-dimensional $B$-module $V$. (See S.\,Donkin, {\it Rational
representations of algebraic groups. Tensor products and
filtration $($Lect. Notes Math.}, \textbf{1140}), Springer, Berlin,
1985). Suppose in addition that $K[X]$ is a Cohen--Macaulay ring,
that is, a free module over the subalgebra generated by any
homogeneous system of parameters. Is then the ring of invariants
$K[X]^G$ Cohen--Macaulay? \ \ {\it Comment of 2005:} This is
proved in the case where $X$ is a rational $G$-module
(M.\,Hashimoto, {\it Math. Z.}, \textbf{236} (2001), 605--623).
\hfill \raisebox{-1ex}{\sl A.\,N.\,Zubkov} \emp

 \bmp
 \textbf{15.45.}
 We define the class of {\it hierarchically
decomposable\/} groups in the following way. First, if ${\frak X}$ is
any class of groups, then let $\textbf{H}_1{\frak X}$ denote the class
of groups which admit an admissible action on a finite-dimensional
contractible complex in such a way that every cell stabilizer belongs
to ${\frak X}$. Then the ``big'' class $\textbf{H}{\frak X}$ is defined
to be the smallest $\textbf{H}_1$-closed class containing ${\frak X}$.

\makebox[25pt][r]{\zv a)} Let ${\frak F}$ be the class of all finite
groups. Find an
 example of an $\textbf{H}{\frak F}$ group which is not in ${\bf
H}_3{\frak F}$ $(=\textbf{H}_1\textbf{H}_1\textbf{H}_1F)$.

\makebox[15pt][r]{}b) Prove or disprove that there is an ordinal
$\alpha$ such that $\textbf{H}_{\alpha}{\frak F}=\textbf{H}{\frak F}$,
where $\textbf{H}_{\alpha}$ is the operator on classes of groups defined
by transfinite induction in the obvious way starting from ${\bf
H}_1$. 
\ {\it Editor's comment}: It is proved that such an ordinal cannot be countable (T.\,Januszkiewicz, P.\,H.\,Kropholler, I.\,J.\,Leary, {\it Bull. London Math. Soc.}, {\bf 42}, No. 5 (2010), 896--904).
\hfill \raisebox{-1ex}{\sl P.\,Kropholler}

\ul

\otv a) Such an example is found; moreover, for any countable ordinal $\alpha$ there are groups in $\textbf{H}_{\alpha+1}{\frak F}$ that are not in ${\bf
H}_{\alpha}{\frak F}$ (T.\,Januszkiewicz, P.\,H.\,Kropholler, I.\,J.\,Leary, {\it Bull. London Math. Soc.}, {\bf 42}, No. 5 (2010), 896--904). Furthermore,
torsion-free
examples with similar properies have been constructed (F.\,Fournier-Facio, B.\,Sun, {\it Preprint}, 2025, \url{https://arxiv.org/abs/2503.01987v3}).
\emp

 \bmp
 \textbf{15.46.} Can the question 7.28 on conditions for
admissibility of an elementary carpet ${\goth A} = \{ {\goth A}_r\mid
r \in \Phi \}$
 be reduced to Lie rank 1 if $K$ is a field?
A carpet ${\goth A}$ of type $\Phi$
 of additive subgroups of $K$ is called {\it addmissible\/} if
in the
Chevalley group over $K$ associated with the root system
 $\Phi$ the subgroup
 $\left< x_r({\goth A}_r)\mid r \in \Phi \right>$ intersects
 $x_r(K)$ in $x_r({\goth A}_r)$. More precisely, is it true that
the carpet ${\goth A}$ is admissible if and only if the subcarpets
 $\{ {\goth A}_r, \,{\goth A}_{-r}
\}$, $r \in \Phi$, of rank 1 are admissible? The answer is known to
be affirmative if the field $K$ is locally finite
 (V.\,M.\,Levchuk, {\it Algebra and Logic}, \textbf{22}, no.\,5 (1983),
 362--371).
\hfill \raisebox{-1ex}{\sl V.\,M.\,Levchuk}

\emp

 \bmp
\textbf{\ai 15.47.} 
Let $M<G\leq {\rm Sym}(\Omega)$, where $\Omega$ is
finite, be such that $M$ is transitive on $\Omega$ and there is a
$G$-in\-va\-ri\-ant partition $\cal{P}$ of
$\Omega\times\Omega\setminus\{(\alpha,\alpha)\mid \alpha\in\Omega\}$
such that $G$ is transitive on the set of parts of $\cal{P}$ and $M$
fixes each part of $\cal{P}$ setwise. (Here $\cal{P}$ can be
identified with a decomposition of the complete directed graph with
vertex set $\Omega$ into edge-disjoint isomorphic directed graphs.)
If $G$ induces a cyclic permutation group on $\cal{P}$, then we
showed (\!{\it Trans. Amer. Math. Soc.}, \textbf{355}, no.\,2
(2003),
637--653) that the numbers $n=|\Omega|$ and
$k=|\cal{P}|$ are such that the $r$-part $n_r$ of $n$ satisfies
$n_r\equiv 1\pmod{k}$ for each prime $r$. Are there examples with $G$
inducing a non-cyclic permutation group on $\cal{P}$ for any $n, k$
not satisfying this congruence condition? \hfill \raisebox{-1ex}{\sl
C.\,H.\,Li, C.\,E.\,Praeger}

\ul

\oai No, there are no such examples (unverified AI solution: B.\,G.\,Di Siena, {\it Preprint}, 2026, \url{https://www.proofforum.org/abs/PF:2026.000036}).
\emp

 \bmp \textbf{15.48.} Let $G$ be any non-trivial finite group,
and let $X$ be any generating set for $G$. Is it true
 that every element of $G$ can be obtained from $X$ using fewer than
$2\log_2|G|$ multiplications? (When counting the number of
mutiplications on a path from the generators to a given element, at
each step one can use the elements obtained at previous steps.)
\hfill \raisebox{-1ex}{\sl C.\,R.\,Leedham-Green}

\emp

 \bmp \textbf{\zv 15.50.}
 Let $G$ be a group of automorphisms of an abelian
group of prime exponent. Suppose that there exists $x\in G$ such that
$x$ is regular of order 3 and the order of $[x,g]$ is finite for
every $g\in G$. Is it true that $\langle x^G\rangle$ is locally
finite? \hfill \raisebox{-1ex}{\sl V.\,D.\,Mazurov}

\ul

\otv No, not always (AI solution confirmed by V.\,D.\,Mazurov: \  B.\,G.\,Di~Siena, {\it Preprint}, 2026, \url{https://kourovkanotebookorg.wordpress.com/wp-content/uploads/2026/10/15-50.pdf}).
\emp

 \bmp \textbf{15.51.} Suppose that $G$ is a periodic group
satisfying the identity $[x,y]^5=1$. Is then the derived subgroup
$[G,G]$ a $5$-group? \hfill \raisebox{-1ex}{\sl V.\,D.\,Mazurov}

\emp

 \bmp \textbf{15.52.} (Well-known problem). By a famous
theorem of
 Wielandt the sequence $G_0=G,\,G_1,\ldots $, where $
G_{i+1}={\rm Aut\,}G_i$, stabilizes for any finite group $G$ with
trivial centre. Does there exist a function $f$ of natural argument such
 that $|G_i|\leq
 f(|G|)$ for all $i=0,\,1,\ldots $
 for an arbitrary finite group $G$ and the same kind of sequence?
\hfill \raisebox{-1ex}{\sl V.\,D.\,Mazurov}

\emp

\bmp \textbf{\zv 15.53.}
 Let $S$ be the set of all
 prime numbers $p$ for which there exists a finite
simple group $G$ and an absolutely irreducible $G$-module $V$ over a
field of characteristic $p$ such that the order of any element in the
natural semidirect product $VG$ coincides with the order of some
element in $G$. Is $S$ finite or infinite? \hfill \raisebox{-1ex}{\sl
V.\,D.\,Mazurov}

\ul

\otv The set $S$ is finite and consists of 2 and 3. This follows from the
complete list of finite simple groups $G$ with a module $V$ in
characteristic $p$ satisfying those properties, which was determined
in a number of papers, the last of which is (M.\,A.\,Grechkoseeva,
S.\,V.\,Skresanov, {\it Sibirsk. Elektron. Matem. Izv.}, {\bf 17}
(2020), 585--589). Then the group $G$ must be either $3D_4(2)$ with
$p=2$ (V.\,D.\,Mazurov, {\it Algebra Logic}, {\bf 52}, no.\,5 (2013),
400--403), or $U_5(2)$ with $p=3$ (A.\,V.\,Zavarnitsine, {\it Siberian
Math. J.}, {\bf 49}, no.\,2 (2008), 246--256, and
M.\,A.\,Grechkoseeva, {\it J. Algebra}, {\bf 339}, no.\,1 (2011),
304--319), or $U_3(q)$, where $q$ is a special Mersenne prime, with
$p=2$ (A.\,V.\,Zavarnitsine, {\it Algebra Logic}, {\bf 45}, no.\,2
(2006), 106--116).
\emp

 \bmp \textbf{15.54.} Suppose that $G$ is a periodic group
containing an involution $i$ such that the centralizer $C_G(i)$ is a
locally cyclic $2$-group. Does the set of all elements of odd order
in $G$ that are inverted by $i$ form a subgroup? \hfill
\raisebox{-1ex}{\sl V.\,D.\,Mazurov}

\emp

 \bmp \textbf{15.55.} a) The Monster, $M$, is a
$6$-trans\-po\-si\-tion
group. Pairs of
 Fischer transpositions generate 9 \ $M$-classes of dihedral
groups. The
 order of the product of a pair is the coefficient of the highest
 root of affine type $E_8$. Similar properties hold for Baby $B$, and
$F_{24}'$ with respect to $E_7$ and $E_6$ when the product is read
modulo centres ($2.B$, $3.F_{24}'$). Explain this. \ {\it Editors'
Comment of 2005:\/} Some progress was made in
 (C.\,H.\,Lam, H.\,Yamada, H.\,Yamauchi, {\it Trans. Amer. Math. Soc.}, {\bf 359}, no.\,9 (2007), 4107--4123).

\makebox[15pt][r]{}b) Note that the
Schur multiplier of the sporadic simple groups $M$, $B$, $F_{24}'$ is
the fundamental group of type $E_8$, $E_7$, $E_6$, respectively. Why?

\makebox[15pt][r]{} See (J.\,McKay, in: {\it Finite groups, Santa
Cruz Conf. 1979 $($Proc. Symp. Pure Math., \textbf{37}}), Amer. Math.
Soc., Providence, RI, 1980, 183--186) and (R.\,E.\,Borcherds, {\it
Doc. Math., J. DMV Extra Vol.
 ICM Berlin 1998}, Vol. I (1998), 607--616).
\hfill\raisebox{-1ex}{\sl J.\,McKay}

\emp

 \bmp \textbf{15.56.} Is there a spin manifold such that the
Monster acts on its loop space? Perhaps 24-dimensional?
hyper-K\"{a}hler, non-compact? See (F.\,Hirzebruch, T.\,Berger,
R.\,Jung, {\it Manifolds and modular forms}, Vieweg, Braunschweig,
1992). \hfill\raisebox{-1ex}{\sl J.\,McKay} \emp

 \bmp \textbf{15.57.} Suppose that $H$ is a subgroup of
$SL_2({\Bbb Q} )$ that is dense in the Zariski topology and has no nontrivial finite quotients. Is then $H=SL_2({\Bbb Q} )$?
\hfill\raisebox{-1ex}{\sl J.\,McKay, J.-P.\,Serre}
\emp

 \bmp \textbf{15.58.} Suppose that a free
profinite product $G*H$ is a free profinite group of finite rank.
Must $G$ and $H$ be free profinite groups?

\makebox[15pt][r]{}By (J.\,Neukirch, {\it Arch. Math.}, \textbf{22},
no.\,4 (1971), 337--357) this may not be true if the rank of $G*H$
is infinite. \hfill \raisebox{-1ex}{\sl O.\,V.\,Mel'nikov}

\emp

 \bmp \textbf{15.59.} Does there exist a
profinite group $G$ that is not free but can be represented as
a~projective limit $G=\lim\limits_{\longleftarrow}(G/N_{\alpha})$, where
all the $G/N_{\alpha}$ are free profinite groups of finite ranks?

\makebox[15pt][r]{}The finiteness condition on the ranks of the
$G/N_{\alpha}$ is essential. Such a group $G$ cannot satisfy the first
axiom of
 countability (O.\,V.\,Mel'nikov, {\it Dokl. AN BSSR}, {\bf
24}, no.\,11 (1980), 968--970 (Russian)). \hfill
\raisebox{-1ex}{\sl O.\,V.\,Mel'nikov}

\emp

 \bmp
 \textbf{15.61.}
Is it true that \ $l_ {\pi} ^ {n} (G) \le n (G_ {\pi} ) -1+\max_
{p\in \pi} l_p
 (G) $ \ for any $\pi$-solu\-ble group~$G$? Here
$n (G_ {\pi} ) $ is the nilpotent length of a Hall
$\pi$-sub\-group
 $G_ {\pi} $ of the
group $G$ and $l_ {\pi} ^ {n} (G) $ is the nilpotent $\pi$-length of
$G$, that is, the minimum number of
 $\pi$-fac\-tors in those normal series of~$G$ whose factors are either
 $\pi '$-groups, or
 nilpotent $\pi$-groups. The answer is known to be affirmative in the
case when all proper subgroups of
 $G_ {\pi} $ are supersoluble.
\hfill \raisebox{-1ex}{\sl V.\,S.\,Monakhov}

\emp

 \bmp
 \textbf{15.63.}
a) Let $F_n$ be the free group of finite rank $n$ on the free
generators $x_1,\ldots ,x_n$. An element $u \in F_n$ is called
{\it positive\/} if $u$ belongs to the semigroup generated by the
$x_i$. An element $u \in F_n$ is called {\it potentially
positive\/} if $\alpha (u)$ is positive for some automorphism
$\alpha$ of~$F_n$.
Is the property of an
element to be potentially positive algorithmically recognizable?

\makebox[15pt][r]{}{\it Comment of 2013:}
In (R.\,Goldstein, {\it Contemp. Math.}, Amer. Math. Soc., {\bf 421}
(2006), 157--168) the problem was solved in the affirmative in the special case $n=2$.

\hfill \raisebox{-1ex}{\sl A.\,G.\,Myasnikov,~V.\,E.\,Shpilrain}

\emp

 \bmp
 \textbf{15.64.}
For finite groups $G$, $X$ define $r(G;X)$ to be
the number of inequivalent actions of $G$ on $X$,
that is, the number of equivalence classes of
homomorphisms $G \to \hbox{\rm Aut}\, X$,
where equivalence is defined by conjugation by an element
of $\hbox{\rm Aut}\,X$.
Now define $
r_G(n) := \max\{r(G;X) : |X| = n\}.
$

\makebox[15pt][r]{}Is it true that $r_G(n)$ may be bounded as a
function of $\lambda(n)$, the total number (counting
multiplicities) of prime factors of $n$? \hfill
\raisebox{-1ex}{\sl P.\,M.\,Neumann}

\emp

 \bmp
 \textbf{\ai 15.65.} 
A square matrix is said to be {\it separable\/} if its minimal polynomial
has no repeated roots,
 and {\it cyclic\/} if its minimal and characteristic polynomials are equal.
For a matrix group $G$ over a finite field define $s(G)$ and
$c(G)$ to be the proportion of separable and of cyclic elements
respectively in $G$. For a classical group ${\rm X}(d,q)$ of
dimension $d$ defined over the field with $q$ elements let \
$S(\hbox{\rm X}; q) := \lim_{d \to \infty}s(\hbox{\rm X}(d,q))$ \
and \ $ C(\hbox{\rm X}; q) := \lim_{d \to \infty}c(\hbox{\rm
X}(d,q)).$ Independently G.\,E.\,Wall ({\it Bull. Austral. Math.
Soc.}, \textbf{60}, no.\,2 (1999), 253--284)
 and J.\,Fulman ({\it J. Group Theory}, \textbf{2}, no.\,3 (1999),
251--289) have evaluated $S(\hbox{\rm GL}; q)$ and $C(\hbox{\rm
GL}; q)$, and have found them to be rational functions of $q$. Are
$S(\hbox{\rm X}; q)$ and $C(\hbox{\rm X}; q)$ rational functions
of $q$ also for the unitary, symplectic, and orthogonal groups?
\hfill \raisebox{-1ex}{\sl P.\,M.\,Neumann}

\ul

\oai No, they are not (unverified AI solution:  M.\,Kinyon, J.\,Urban, {\it Preprint}, 2026,

\url{https://kourovkanotebookorg.wordpress.com/wp-content/uploads/2026/09/kourovka-mathematics.pdf}).

\emp 

 \bmp
 \textbf{15.66.}
For a class $\mathfrak X$ of groups let $g_{\mathfrak X}(n)$ be
the number of groups of order $n$ in the class $\mathfrak X$ (up
to isomorphism). Many years ago I formulated the following
problem: find good upper bounds for the quotient $g_{\mathfrak
V}(n)/g_{\frak U}(n)$, where $\frak V$ is a variety that is
defined by its finite groups and $\frak U$ is a subvariety of
$\mathfrak V$. (This quotient is not defined for all $n$ but only
for those for which there are groups of order $n$ in $\mathfrak
U$.) Some progress has been made by G.\,Venkataraman ({\it Quart.
J. Math. Oxford (2)}, \textbf{48}, no.\,189 (1997), 107--125) when
$\mathfrak V$ is a variety generated by finite groups all of whose
Sylow subgroups are abelian. 
{\it Conjecture:\/} if $\mathfrak V$ is a locally finite
variety
of $p\hs$-groups and $\mathfrak U$ is a non-abelian subvariety of ${\frak
V}$, then
$g_{\mathfrak V}(p^m)/g_{\mathfrak U}(p^m) < p^{O(m^2)}$.

\makebox[15pt][r]{}Moreover, this seems a possible way to attack
the Sims Conjecture that when we write the number of groups of
order $p^m$ as $p^{{2\over27}m^3 + \varepsilon(m)}$ the error term
$\varepsilon(m)$ is $O(m^2)$.

\hfill \raisebox{-1ex}{\sl P.\,M.\,Neumann}

\emp

 \bmp
 \textbf{15.67.} Which adjoint Chevalley groups (of
normal type) over the integers are generated by three involutions
two of which commute?

 \makebox[15pt][r]{}The groups $SL_n({\Bbb Z})$, $n>13$, satisfy
this condition (M.\,C.\,Tamburini, P.\,Zucca, {\it J.~Algebra},
\textbf{195}, no.\,2 (1997), 650--661). The groups
 $PSL_n({\Bbb Z})$ satisfy it if and only if
$n>4$ (N.\,Ya.\,Nuzhin, {\it Vladikavkaz. Mat. Zh.}, \textbf{10},
no.\,1 (2008), 68--74 (Russian)). The group
 $PSp_4({\Bbb Z})$ does not satisfy it, which follows from the
corresponding fact for $PSp_4(3)$, see Archive, 7.30.
The group $G_2$
satisfies it (I.\,A.\,Timofeenko, {\it J.~Siberian Fed. Univ. Ser. Math. Phys.}, {\bf 8}, no.\,1 (2015), 104--108. The group $E_l$ satisfies it (I.\,A.\,Timofeenko, {\it Sibrsk. \`Electron. Mat. Izv.}, {\bf 14} (2017), 807--820 (Russian)).

\hfill
\raisebox{-1ex}{\sl Ya.\,N.\,Nuzhin}

\emp

 \bmp
 \textbf{15.68.} Does there exist an infinite
finitely generated 2-group (of finite exponent) all of whose
proper subgroups are locally finite? \hfill \raisebox{-1ex}{\sl
A.\,Yu.\,Olshanskii}

\emp

 \bmp
 \textbf{15.69.}
Is it true that every hyperbolic group has a free normal subgroup
with the factor-group of finite exponent? \hfill
\raisebox{-1ex}{\sl A.\,Yu.\,Olshanskii}

\emp

 \bmp \textbf{\zv 15.70.}
Do there exist groups of arbitrarily large cardinality that satisfy
the minimum condition for subgroups? \hfill \raisebox{-1ex}{\sl
A.\,Yu.\,Olshanskii}

\ul

\otv A positive answer is consistent with the ZFC axioms of set theory, and the consistency of a negative answer follows from the consistency of certain large cardinal assumptions (S.\,M.\,Corson, S.\,Shelah, {\it Preprint}, 2024, \url{https://arxiv.org/pdf/2408.03201}).
\emp

\bmp
 \textbf{15.71.}
 (B.\,Huppert). Let $G$ be a finite solvable group, and let $\rho(G)$
denote the set of prime divisors of the degrees of irreducible
characters of $G$. Is it true that there always exists an
irreducible character of $G$ whose degree is divisible by at least
$|\rho(G)|/2$ different primes? \hfill \raisebox{-1ex}{\sl
P.\,P.\,P\'alfy}\emp

\bmp
 \textbf{15.72.}
For a fixed prime $p$ does there exist a sequence of groups $P_n$ of
order $p^n$ such that the number of conjugacy classes $k(P_n)$
satisfies $\lim_{n\to\infty} \log k(P_n) / \sqrt{n} = 0$?

\makebox[15pt][r]{}Note that J.\,M.\,Riedl ({\it J.~Algebra}, {\bf
218} (1999), 190--215) constructed $p\hs$-groups for which the above
limit is $2 \log p$.
 \hfill \raisebox{-1ex}{\sl P.\,P.\,P\'alfy}\emp

\bmp
 \textbf{15.73.} Is it true that for every finite lattice $L$ there
exists a finite group $G$ and a subgroup $H \le G$ such that the
interval $Int(H;G)$ in the subgroup lattice of $G$ is isomorphic to
$L$? (Probably not.) \hfill \raisebox{-1ex}{\sl P.\,P.\,P\'alfy} \emp

\bmp
 \textbf{15.74.}
For every prime $p$ find a finite $p\hs$-group of nilpotence class
$p$ such that its lattice of normal subgroups cannot be embedded
into the subgroup lattice of any abelian group. (Solved for $p=2$,
$3$.) \hfill \raisebox{-1ex}{\sl P.\,P.\,P\'alfy}\emp

 \bmp
\textbf{15.75.}
b) Consider the sequence $u_1 = [x,y],\ldots$, $ u_{n+1} =
[[u_n, x],[u_n,y]]$. Is it true that an arbitrary finite group is
soluble if and only if it satisfies one of these identities $u_n =
1$? \hfill \raisebox{-1ex}{\sl B.\,I.\,Plotkin} \emp

 \bmp \textbf{\ai 15.76.} 
 b)
If $\Theta$ is a variety of groups, then let $\Theta ^0$ denote the
category of all free groups of finite rank in $\Theta$. It is proved
(G.\,Mashevitzky, B.\,Plotkin, E.\,Plotkin {\it J.~Algebra}, {\bf
282} (2004), 490--512) that if $\Theta$ is the variety of all groups,
then every automorphism of the category $\Theta ^0$ is an inner one.
The same is true if $\Theta$ is the variety of all
nilpotent groups of class $\leq d$, for $d\geq 2$ (A.\,Tsurkov, {\it Int. J. Algebra Comput.}, {\bf 17} (2007), 1273--1281). Is
this true for the varieties of solvable groups?
of metabelian groups?

\makebox[15pt][r]{}An automorphism $\varphi$ of a category is called
{\it inner\/} if it is isomorphic to the identity automorphism. Let
$s:1\to\varphi$ be a function defining this isomorphism. Then for
every object $A$ we have an isomorphism $s_A: A\to \varphi(A)$ and
for any morphism of objects $\mu: A\to B$ we have
$\varphi(\mu)=s_B\mu s_A^{-1}$.
 \hfill \raisebox{-1ex}{\sl B.\,I.\,Plotkin}

\ul

\oai Yes, it is true (unverified AI solution:  M.\,Kinyon, J.\,Urban, {\it Preprint}, 2026,

\url{https://kourovkanotebookorg.wordpress.com/wp-content/uploads/2026/09/kourovka-mathematics.pdf}).
\emp 

\bmp
 \textbf{\zv 15.77.}
Elements $a,b$ of a group $G$ are said to be {\it symmetric with
respect to an element\/} $ g\in G$ if $a=gb^{-1}g$. Let $G$ be an
infinite abelian group, $\alpha$ a cardinal, $\alpha <|G|$. Is it
true that for any $n$-colour\-ing $\chi : G\rightarrow\{ 0,\,1,\ldots
,n-1\}$ there exists a monochrome subset of cardinality $\alpha$
that is symmetric with respect to some element of $G$? This is
known to be true for $n\leq 3$. \hfill \raisebox{-1ex}{\sl
I.\,V.\,Protasov}

\ul

\otv Yes, it is true (AI solution confirmed by I.\,V.\,Protasov: \  Achyuth Jayadevan,  {\it Preprint}, 2026,
\url{https://zenodo.org/records/22909420}).
 \emp 

 \bmp \textbf{15.78.} (R.\,I.\,Grigorchuk). Is it true that for any
$n$-colouring of a free group of any rank there exists a
monochrome subset that is symmetric with respect to some element
of the group? This is true for $n=2$. \hfill \raisebox{-1ex}{\sl
I.\,V.\,Protasov}
\emp

 \bmp \textbf{\zv 15.80.}
 A sequence $\{ F_n\}$
of pairwise disjoint finite subsets of a topological group is
called {\it expansive\/} if for every open subset $U$ there is a
number $m$ such that $ F_n\cap U\ne \varnothing$ for all $n>m$.
Suppose that a group $G$ can be partitioned into countably many
dense subsets. Is it true that in $ G$ there exists an expansive
sequence? \hfill \raisebox{-1ex}{\sl I.\,V.\,Protasov}

\ul

\otv Not necessarily. One example is $G = \prod_{\mathbb{N}}\mathbb{R}$ under the box topology. It is well-known that $G$ does not have a countable dense subset (hence has no expansive sequence). However, the topology on $G$ has a basis having $2^{\aleph_0}$ elements, each element of the basis having $2^{\aleph_0}$ points, and one can partition $G$ into countably many dense subsets (via a transfinite induction of length $2^{\aleph_0}$). ({\it Letter of S.\,Corson of 19 May 2025}.) Another example: any infinite abelian group $G$ with finite set of elements of order $2$ can be partitioned into countably many subsets dense in every non-discrete group topology on $G$ according to Corollary~12.21 in (Y.\,Zelenyuk, {\it Ultrafilters and topologies on groups}, De Gruyter, 2011) ({\it Letter of I.\,V.\,Protasov of 19 May 2025}).
\emp

 \bmp
 \textbf{15.82.}
Suppose that a periodic group $G$ contains a strongly isolated
 $2$-sub\-group $U$. Is it true that either
 $G$ is locally finite, or $U$ is a normal subgroup of~$G$?

\hfill \raisebox{-1ex}{\sl A.\,I.\,Sozutov, N.\,M.\,Suchkov}

\emp

 \bmp
 \textbf{15.83.} (Yu.\,I.\,Merzlyakov). Does there
exist a rational number $\alpha $ such that $|\alpha |<2$ and the
matrices
 {\small $\left(
\begin{array}{cc}1&\alpha \\0&1\end{array}\right)$} and {\small
$\left(
\begin{array}{cc}1&0\\\alpha &1\end{array}\right)$} generate a free
group? \hfill \raisebox{-1ex}{\sl Yu.\,V.\,Sosnovski\u{\i}}

\emp

 \bmp
 \textbf{15.84.} Yu.\,I.\,Merzlyakov ({\it Sov. Math. Dokl.}, \textbf{19}
(1978), 64--68) proved that
 if the complex numbers $\alpha ,\beta ,\gamma $ are each at
least 3 in absolute value,
then the matrices
{\small $\left(
\begin{array}{cc}1&\alpha \\0&1\end{array}\right)$}, \
{\small $\left(
\begin{array}{cc}1&0 \\\beta&1\end{array}\right)$}, \ and \
{\small $\left(
\begin{array}{cc}1-\gamma &-\gamma \\\gamma &1+\gamma \end{array}\right)$}
generate a free group of rank three. Are there rational numbers
$\alpha ,\beta , \gamma $, each less than 3 in absolute value,
with the same property? \hfill \raisebox{-1ex}{\sl
Yu.\,V.\,Sosnovski\u{\i}}

\emp

 \bmp
 \textbf{15.85.} A torsion-free group all of whose subgroups
are subnormal is nilpotent (H.\,Smith, {\it Arch. Math.\/}, {\bf
76}, no.\,1 (2001), 1--6). Is a torsion-free group with the
normalizer condition

\makebox[15pt][r]{}a) hyperabelian?

\makebox[15pt][r]{}b) hypercentral? \hfill {\sl
Yu.\,V.\,Sosnovski\u{\i}}

\emp

 \bmp
 \textbf{15.87.}
Suppose that a $2$-group $G$ admits a regular periodic locally
cyclic group of automorphisms that transitively permutes the set
of involutions of $G$.
 Is $G$ locally finite?
\hfill \raisebox{-1ex}{\sl N.\,M.\,Suchkov}

\emp

 \bmp
 \textbf{15.88.} Let ${\goth A}$ and ${\goth
N}_c$ denote the varieties of abelian groups and nilpotent groups of
class $\leq c$, respectively. Let $F_r= F_r({\goth A}{\goth N}_c)$ be
a free group of rank $r$ in $ {\goth A}{\goth N}_c$. An element of
 $ F_r$ is called {\it primitive\/} if it can
be included in a basis of~$ F_r$. Does there exist an algorithm
recognizing primitive elements in~$ F_r$? \hfill\raisebox{-1ex}{\sl
E.\,I.\,Timoshenko} \emp

\bmp
 \textbf{\zv 15.89.}
Let $\Gamma$ be an infinite undirected connected vertex-symmetric
graph of finite valency without loops or multiple edges. Is it true
that every complex number is an eigenvalue of the adjacency matrix of
$\Gamma$ under its natural action as a linear operator on the complex
vector space of all complex-valued functions on the vertex
set of $\Gamma$?

 \hfill \raisebox{-1ex}{\sl V.\,I.\,Trofimov}

 \ul

 \otv No, not always (AI solution confirmed by V.\,I.\,Trofimov: \ Achyuth Jayadevan, {\it Preprint},  2026, \url{https://zenodo.org/records/22722077}).  \emp

 \bmp
 \textbf{15.90.}
 Let $\Gamma$ be an infinite directed graph, and $\overline{ \Gamma}$ the
underlying
undirected graph. Suppose that the graph $\Gamma$ admits a vertex-transitive
group of automorphisms, and the graph $ \overline{ \Gamma}$ is connected and
of finite valency. Does there exist a positive integer $k$ (possibly
depending
 on $\Gamma$) such that for any positive integer $n$
there is a directed path of length at most $k\cdot n$ in the graph
$\Gamma$ whose initial and terminal vertices are at distance at least
$n$
 in the graph $\overline{
\Gamma}$?

 \makebox[15pt][r]{}{\it Comment of 2005:\/} It is proved
that there exists a positive integer $k'$ (depending only on the
valency of $\overline{ \Gamma}$) such that for any positive
integer $n$ there is a directed path of length at most $k'\cdot
n^2$ in the graph $\Gamma$ whose initial and terminal vertices are
at distance at least $n$ in the graph $\overline{ \Gamma}$
(V.\,I.\,Trofimov, {\it Europ. J. Combinatorics}, \textbf{27}, no.\,5
(2006), 690--700). \hfill \raisebox{-1ex}{\sl V.\,I.\,Trofimov}

\emp

 \bmp \textbf{15.91.} Is it true that any irreducible faithful
representation of a linear group $G$ of finite rank over a finitely
generated field of characteristic zero is induced from an irreducible
representation of a finitely generated dense subgroup of the group
$G$?

\hfill \raisebox{-1ex}{\sl A.\,V.\,Tushev}

\emp

 \bmp \textbf{\zv 15.92.}
A group $\left< x,y\mid x^l=y^m=(xy)^n\right>$ is called the {\it
triangle group\/} with parameters $(l,\,m,\,n)$. It is proved in
(A.\,M.\,Brunner, R.\,G.\,Burns, J.\,Wiegold, {\it Math. Scientist},
\textbf{4} (1979), 93--98) that the triangle group $ (2,\,3,\,30)$ has
uncountably many non-isomorphic homomorphic
 images that are residually finite alternating groups. Is the same true
 for the triangle group $(2,\,3,\,n)$ for all $ n > 6$?
\hfill \raisebox{-1ex}{\sl J.\,Wiegold}

\ul

\otv Yes, it is true (AI solution confirmed by D.\,Kielak: \  M.\,Kinyon, J.\,Urban, {\it Preprint}, 2026,
\url{https://kourovkanotebookorg.wordpress.com/wp-content/uploads/2026/09/kourovka-mathematics.pdf}).
\emp 

\bmp \textbf{15.93.} Let $G$ be a pro-$p$ group, $p>2$, and $\varphi$ an
automorphism of $G$ of order~2. Suppose that the centralizer
$C_G(\varphi )$ is abelian. Is it true that $G$ satisfies a pro-$p$
identity?

\makebox[15pt][r]{}An affirmative answer would generalize a result
of A.\,N.\,Zubkov ({\it Siberian Math.~J.}, \textbf{28}, no.\,5
(1987), 742--747) saying that non-abelian free pro-$p$ groups
cannot be represented by $2\times 2$ matrices. \hfill
\raisebox{-1ex}{\sl A.\,Jaikin-Zapirain}

\emp

 \bmp \textbf{15.94.}
Define the {\em weight} of a group $G$ to be the minimum number of
generators of $G$ as a normal subgroup of itself.
Let $G=G_1*\cdots *G_n$ be a free product of $n$ nontrivial groups.

\makebox[15pt][r]{}a) Is it true that the
weight of $G$ is at
least $n/2$?

\makebox[15pt][r]{}b) Is it true if the $G_i$ are cyclic?

\makebox[15pt][r]{}c) Is it true for $n=3$ (without assuming the $G_i$
cyclic)?

\makebox[15pt][r]{}Problem 5.53 in Archive (now answered in the affirmative) is
the special case with $n=3$ and all the $G_i$ cyclic.\hfill
\raisebox{-1ex}{\sl J.\,Howie}

\emp

 \bmp \textbf{15.95.} (A.\,Mann, Ch.\,Praeger). Suppose
that all fixed-point-free elements of a~transitive permutation
 group $G$ have prime order $p$. If $G$ is a finite $p\hs$-group,
must $G$ have exponent $p$?

\makebox[15pt][r]{}This is true for $p=2,\,3$; it is also known
 that the
exponent of $G$ is bounded in terms of $p$ (A.\,Hassani,
M.\,Khayaty, E.\,I.\,Khukhro, C.\,E.\,Praeger, {\it J.~Algebra},
\textbf{214}, no.\,1 (1999), 317--337). \hfill \raisebox{-1ex}{\sl
E.\,I.\,Khukhro}

\emp

 \bmp \textbf{15.96.} 
  An automorphism $\f $ of a group $G$ is
called a {\it splitting automorphism of order
 $n$\/} if $\f ^n=1$ and $xx^{\f}x^{\f
^{2}}\cdots x^{\f ^{n-1}}=1$ for any $x\in G$.

\makebox[25pt][r]{\zva a)} Is it true that the derived length of a
$d$-gene\-ra\-ted nilpotent $p$-group admitting a splitting
automorphism of order $p^n$ is bounded by a function of $d$,
 $p$, and $n$? This is true for $n=1$, see Archive,~7.53.

\makebox[25pt][r]{b)} The same question for $p^n=4$. \hfill
\raisebox{0ex}{\sl E.\,I.\,Khukhro}

\ul

\otv (a) No, not always (J.\,Yang,  {\it Preprint}, 2026, \url{https://doi.org/10.13140/RG.2.2.34838.10568}).
\emp

 \bmp
 \textbf{15.97.} Let $p$ be a prime. A group $G$ satisfies the
{\it $p\hs$-mini\-mal condition\/} if there are no infinite
descending chains $G_1>G_2>\ldots $ of subgroups of $G$ such that
each difference $G_i\setminus G_{i+1}$ contains a $p\hs$-ele\-ment
(S.\,N.\,Chernikov). Suppose that a locally finite group $G$
satisfies the $p\hs$-mini\-mal condition and has a subnormal series
each of whose factors is finite or a $p'$-group. Is it true that all
$p\hs$-ele\-ments of $G$ generate a Chernikov subgroup?

\hfill \raisebox{-1ex}{\sl N.\,S.\,Chernikov} \emp

 \bmp
 \textbf{15.98.}
 Let $\frak{F}$ be a saturated formation, and $G$
a finite soluble minimal non-$\frak{F}$-group such that
$G^{\frak{F}}$ is a Sylow $p\hs$-sub\-group of $G$. Is it true
that $G^{\frak{F}}$ is isomorphic to a factor group of the Sylow
$p\hs$-sub\-group of $PSU(3, p^{2n})$? This is true for the
formation of finite nilpotent groups (V.\,D.\,Mazurov,
S.\,A.\,Syskin, {\it Math. Notes}, \textbf{14}, no.\,2 (1973),
683--686; \ A.\,Kh.\,Zhurtov, S.\,A.\,Syskin, {\it Siberian
Math.~J.}, \textbf{26}, no.\,2 (1987), 235--239).

\hfill \raisebox{-1ex}{\sl L.\,A.\,Shemetkov}

\emp

 \bmp
 \textbf{15.99.}
Let $f(n)$ be the number of isomorphism classes of finite groups
of order $n$. Is it true that the equation $f(n)=k$ has a solution
for any positive integer $k$? The answer is affirmative for all $k
\leq 1000$ (G.\,M.\,Wei, {\it Southeast Asian Bull. Math.}, {\bf
22}, no.\,1 (1998), 93--102). \hfill \raisebox{-1ex}{\sl
W.\,J.\,Shi}

\emp

 \bmp \textbf{15.100.} Is a periodic group locally finite if it has
a non-cyclic subgroup of order 4 that coincides with its centralizer?
\hfill \raisebox{-1ex}{\sl A.\,K.\,Shl\"epkin}

\emp
 \bmp \textbf{15.101.}
Is a periodic group locally finite if it has an involution whose
centralizer is a locally finite group with Sylow 2-subgroup of order
2? \hfill \raisebox{-1ex}{\sl A.\,K.\,Shl\"epkin} \emp

 \bmp
 \textbf{15.102.} (W.\,Magnus). An
element $r$ of a free group $F_n$ is called a {\it normal root\/} of
an element $u\in F_n$ if $u$ belongs to the normal closure of $r$
in~$F_n$.
 Can an element $u$ that lies outside the commutator
subgroup $[F_n, \,F_n]$ have infinitely many non-conjugate normal
roots? \hfill \raisebox{-1ex}{\sl V.\,E.\,Shpilrain}

\emp

 \bmp
 \textbf{15.103.} (Well-known problem). Is the group ${\rm
Out}(F_3)$ of outer automorphisms of a free group of rank $3$ linear?
\hfill \raisebox{-1ex}{\sl V.\,E.\,Shpilrain} \emp

 \bmp \textbf{15.104.} Let $n$ be a positive integer and let
$w$ be a
group word in the variables $x_1,\,x_2,\ldots $.
Suppose that a residually finite group $G$ satisfies the identity $
w^n=1$.
 Does it follow that the verbal subgroup $w(G)$ is
locally finite?

\makebox[15pt][r]{}This is the Restricted Burnside Problem if
$w=x_1$. A positive answer was obtained also in a number of other
particular cases (P.\,V.\,Shumyatsky, {\it Quart. J.~Math.}, {\bf
51}, no.\,4 (2000), 523--528).
 \hfill \raisebox{-1ex}{\sl P.\,V.\,Shumyatsky}

\emp

 \raggedbottom

\newpage

\pagestyle{myheadings} \markboth{16th Issue (2006)}{16th Issue
(2006)} \thispagestyle{headings} ~ \vspace{2ex}

\centerline {\Large \textbf{Problems from the 16th Issue (2006)}}
\phantomsection\label{16izd}
\vspace{4ex}

\bmp \textbf{16.1.} Let $G$ be a finite non-abelian group, and
$Z(G)$ its centre. One can associate a graph $\Gamma_G$ with $G$
as follows: take $G\backslash Z(G)$ as vertices of $\Gamma_G$ and
join two vertices $x$ and $y$ if $xy\not=yx$. Let $H$ be a finite
non-abelian group such that $\Gamma_G \cong \Gamma_H$.\sm

\makebox[25pt][r]{b)} If $H$ is nilpotent,
 is it true that $G$ is nilpotent? \ \

\makebox[15pt][r]{}{\it Comment of 2009:\/}
This is true if $|H|=|G|$ (A.\,Abdollahi, S.\,Akbari, H.\,R.\,Maimani, {\it J.~Algebra}, \textbf{298} (2006),
468--492).

\makebox[15pt][r]{}{\it Comments of 2023:\/}
This is true if $|Z(H)|\geq |Z(G)|$
(H.\,Shahverdi, {\it J.~Algebra}, \textbf{642} (2024), 60--64), using the earlier result under the additional assumption that the centralizer of every non-central
element in $H$ is abelian (V.\,Grazian, C.\,Monetta, {\it J.~Algebra}, \textbf{633} (2023), 389--402).
\sm

\makebox[25pt][r]{c)} If $H$ is solvable,
 is it true that $G$ is solvable?

 \hf
\raisebox{-1ex}{\sl A.\,Abdollahi, S.\,Akbari, H.\,R.\,Maimani}

\emp

\bmp \textbf{16.3.} Is it true that if $G$ is a finite group with
all conjugacy classes of distinct sizes, then $G\cong {\Bbb
S}_3$?

\makebox[15pt][r]{}This is true if $G$ is solvable (J.\,Zhang,
{\it J.\,Algebra}, \textbf{170} (1994), 608--624); it
is also known that $F(G)$ is nontrivial (Z.\,Arad,
M.\,Muzychuk, A.\,Oliver, {\it J.~Algebra}, \textbf{280} (2004), 537--576). \hf
\raisebox{-1ex}{\sl Z.\,Arad}
 \emp

 \bmp \textbf{16.4.} Let $G$ be a finite group with $C,D$ two nontrivial
conjugacy classes such that $CD$ is also a conjugacy class. Can
$G$ be a non-abelian simple group? \hf \raisebox{-1ex}{\sl
Z.\,Arad}

\emp

\bmp \textbf{16.5.} A group is said to be {\it perfect\/} if it
coincides with its derived subgroup. Does there exist a perfect
locally finite $p$-group

\makebox[25pt][r]{a)} all of whose proper subgroups are
hypercentral?

\makebox[25pt][r]{b)} all of whose proper subgroups are solvable?
 \hf \raisebox{0ex}{\sl A.\,O.\,Asar}

\emp

 \bmp \textbf{16.6.} Can a perfect locally finite $p\hs$-group be
generated by a subset of bounded exponent

\makebox[25pt][r]{a)} if all of its proper subgroups are
hypercentral?

\makebox[25pt][r]{b)} if all of its proper subgroups are
solvable? \hf \raisebox{0ex}{\sl A.\,O.\,Asar}

\emp

\bmp \textbf{16.7.} Is it true that the membership problem is
undecidable for any semidirect product $F_n \leftthreetimes F_n$
of non-abelian free groups $F_n$?

\makebox[15pt][r]{}An affirmative answer would imply an answer
to~6.24. Recall that the membership problem is decidable for
$F_n$ (M.\,Hall, 1949) and undecidable for $F_n \times F_n$
(K.\,A.\,Mikhailova, 1958). \hf \raisebox{-1ex}{\sl
V.\,G.\,Bardakov}

\emp

\bmp \textbf{\zv 16.9.}
An element $g$ of a free group $F_n$ on the
free generators $x_1,\ldots,x_n$ is called a {\it palindrome\/}
with respect to these generators if the reduced word
representing~$g$ is the same when read from left to right or from
right to left. The {\it palindromic length\/} of an element $w \in
F_n$ is the smallest number of palindromes in $F_n$ whose product
is~$w$. Is there an algorithm for finding the palindromic length
of a given element of~$F_n$?

\hf \raisebox{-1ex}{\sl
V.\,G.\,Bardakov, V.\,A.\,Tolstykh, V.\,E.\,Shpilrain}

\ul

\otv Yes, there is (V.\,Diekert, A.\,E.\,Frid, S.\,Natterer, A.\,Thumm, {\it Preprint of 10 October  2026}, \url{https://arxiv.org/abs/2610.14536}; \ AI solution confirmed by V.\,Bardakov and V.\,A.\,Tolstykh: \  Achyuth Jayadevan,  {\it Preprint of 11 September 2026}, \url{https://zenodo.org/records/22722079}; \ unverified AI solutions: J.\,Liao,   {\it Preprint of 15 September 2026}, \url{https://arxiv.org/abs/2609.17027}; \  M.\,Kinyon, J.\,Urban, {\it Preprint of 17 September 2026},
 \url{https://kourovkanotebookorg.wordpress.com/wp-content/uploads/2026/09/kourovka-mathematics.pdf}).

\emp 

\bmp \textbf{16.10.} Is there an algorithm for finding the primitive
length of a given element of~$F_n$? The definition of a primitive
element is given in 14.84; the primitive length is defined
similarly to the palindromic length in 16.9.

\hf
\raisebox{-1ex}{\sl V.\,G.\,Bardakov, V.\,A.\,Tolstykh,
V.\,E.\,Shpilrain}

\emp

\bmp \textbf{\zv 16.11.}
Let $G$ be a finite $p$-group. Does there
always exist a finite $p$-group $H$ such that $\Phi(H)\cong
[G,G]$? \hf \raisebox{-1ex}{\sl Ya.\,G.\,Berkovich}

\ul

\otv No, not always (V.\,Ionin, A.\,Semidetnov,  {\it Preprint}, 2026, \url{https://arxiv.org/pdf/2608.29219}).
\emp

\bmp \textbf{16.14.} Let $G$ be a finite $2$-group such that
$\Omega_1(G)\leq Z(G)$. Is it true that the rank of $G/G^2$ is at
most double the rank of $Z(G)$?\hf \raisebox{-1ex}{\sl
Ya.\,G.\,Berkovich}

\emp

\bmp \textbf{16.15.} An element $g$ of a group $G$ is an {\it Engel
element\/} if for every $h\in G$ there exists $k$ such that
$[h,g,\ldots ,g]=1$, where $g$ occurs $k$ times; if there is such
$k$ independent of $h$, then $g$ is said to be {\it boundedly
Engel}.

\makebox[25pt][r]{b)} Does the set of boundedly
Engel elements  form a subgroup in a torsion-free group?

\makebox[25pt][r]{c)} The same question for right-ordered groups.

\makebox[25pt][r]{d)} The same question for linearly ordered
groups. \hf \raisebox{0ex}{\sl V.\,V.\,Bludov}
\emp

\bmp \textbf{16.16.} {a)} Does the set of (not necessarily
boundedly) Engel elements of a group without elements of order~2
form a subgroup?

\makebox[25pt][r]{b)} The same question for torsion-free groups.

\makebox[25pt][r]{c)} The same question for right-ordered groups.

\makebox[25pt][r]{d)} The same question for linearly ordered
groups.

There are examples of 2-groups where a product of two
(unboundedly) Engel elements is not an Engel element. \hf
\raisebox{-1ex}{\sl V.\,V.\,Bludov}

\emp

\bmp \textbf{16.18.} Does there exist a linearly orderable soluble
group of derived length exactly~$n$ that has a single proper
normal relatively convex subgroup

 \makebox[15pt][r]{}a) for $n=3$?

 \makebox[15pt][r]{}b) for $n > 4$?

Such groups do exist for $ n=2$ and for
$n=4$.\hf \raisebox{0ex}{\sl V.\,V.\,Bludov}

\emp

\bmp \textbf{16.19.} Is the variety of lattice-ordered groups
generated by nilpotent groups finitely based? (Here the variety is
considered in the signature of group and lattice operations.)
\hfill \raisebox{-0.8ex}{\sl V.\,V.\,Bludov}

\emp

\bmp \textbf{\zv 16.20.}
 Let $M$ be a quasivariety of groups. The {\it
 dominion\/}
 ${\rm dom}^M_A(H)$ of a subgroup $H$ of a group $A$ (in
$M$) is the set of all elements $a\in A$ such that for any two
homomorphisms $f,g:A\rightarrow B\in M$, if $f,g$ coincide on $H$,
then $f(a)=g(a)$. Suppose that the set $\{ {\rm dom}^N_A(H)\mid N
\mbox{ is a quasivariety}, N\subseteq M\} $ forms a lattice with
respect to set-theoretic inclusion. Can this lattice be modular
and non-distributive? \hf \raisebox{-0.8ex}{\sl A.\,I.\,Budkin}

\ul

\otv Yes, it can (AI solution confirmed by A.\,I.\,Budkin: \  M.\,Kinyon, J.\,Urban, {\it Preprint}, 2026,
\url{https://kourovkanotebookorg.wordpress.com/wp-content/uploads/2026/09/kourovka-mathematics.pdf}).
\emp 

\bmp \textbf{16.21.} Given a non-central matrix $\alpha \in
SL_n(F)$ over a field $F$ for $n>2$, is it true that every
non-central matrix in
 $SL_n(F)$ is a product of $n$ matrices, each similar to~$\alpha
$?

\hf \raisebox{-0.8ex}{\sl L.\,Vaserstein}

\emp

 \bmp \textbf{16.22.} (Well-known problem.) Let
$E_n(A)$ be the subgroup of $GL_n(A)$ generated by elementary
matrices. Is $SL_2(A) = E_2(A)$ when $A= {\Bbb Z}[x, 1/x]$? \hf
\raisebox{-0.8ex}{\sl L.\,Vaserstein}

\emp

\bmp \textbf{16.23.} Is there, for some $n>2$ and a ring $A$ with
$1$, a matrix in $E_n(A)$ that is nonscalar modulo any proper
ideal and is not a commutator?
 \hf \raisebox{-0.8ex}{\sl L.\,Vaserstein}
\emp

\bmp \textbf{16.26.} We say that the prime graphs of finite groups
$G$ and $H$ coincide if the sets of primes dividing their orders
are the same, $\pi (G)=\pi (H)$, and for any distinct $p,q\in \pi
(G)$ there is an element of order $pq$ in $G$ if and only if there
is such an element in $H$. Does there exist a positive integer $k$
such that there are no $k$ pairwise non-isomorphic finite
non-abelian simple groups with the same graphs of primes? {\it
Conjecture:\/} $k=5$. \hf \raisebox{-0.8ex}{\sl A.\,V.\,Vasil'ev}

\emp

\bmp \textbf{16.28.}
Let $G$ be a connected linear reductive
algebraic group over a field of positive characteristic, $X$ a
closed subset of~$G$, and let $X^k=\{x_1\cdots x_k\mid x_i\in
X\}$.

\makebox[25pt][r]{\zv a)} Is it true that there always exists a
positive integer $c=c(X)>1$ such that $X^c$ is closed?

\makebox[25pt][r]{b)} If $X$ is a conjugacy class of $G$ such that
$X^2$ contains an open subset of $G$, then is $X^2=G$? \hf
\raisebox{-1ex}{\sl E.\,P.\,Vdovin}

\ul

\otv No, not always (AI solution confirmed by E.\,P.\,Vdovin: \ Achyuth Jayadevan, {\it Preprint of 12 September 2026}, \url{https://zenodo.org/records/22729229}; \ unverified AI solution: M.\,Kinyon, J.\,Urban, {\it Preprint of 17 September 2026},    \url{https://kourovkanotebookorg.wordpress.com/wp-content/uploads/2026/09/kourovka-mathematics.pdf}).
\emp

 \bmp \textbf{16.29.}
Which finite simple groups of Lie type $G$ have the following property:
for every semisimple abelian subgroup $A$ and proper subgroup $H$ of $G$ there exists $x\in G$ such that $A^x\cap H=1$?

\makebox[15pt][r]{}This property holds when $A$ is a cyclic subgroup
(J.\,Siemons, A.\,Zalesskii, {\it J.\,Alge\-bra}, \textbf{256} (2002), 611--625), as well as when $A$ is contained in some maximal torus and $G=PSL_n(q)$ (J.\,Siemons, A.\,Zalesskii, {\it J.\,Alge\-bra}, \textbf{226} (2000), 451--478). Note that if $G=L_2(5)$, $A=2\times 2$, and $H=5:2$ (in Atlas notation), then $A^x\cap H>1$ for every $x\in G$ (this example was communicated to the author by V.\,I.\,Zenkov).
\hf \raisebox{-1ex}{\sl E.\,P.\,Vdovin}
\emp

\bmp \textbf{\ai 16.30.} 
Suppose that $A$ and $B$ are subgroups of a
group $G$ and $G = AB$. Will $G$ have composition (principal)
series if $A$ and $B$ have composition (respectively, principal)
series?\hf \raisebox{-1ex}{\sl V.\,A.\,Vedernikov}

\ul

\oai No, not always, both for composition and principal series (unverified AI solution:  Achyuth Jayadevan,  {\it Preprint}, 2026,
\url{https://zenodo.org/records/23421559}).
\emp 

\bmp \textbf{16.32.} Suppose that a group $G$ has a composition
series and let ${\rm Fit}(G)$ be the Fitting class generated by
$G$. Is the set of all Fitting subclasses of ${\rm Fit}(G)$
finite? Cf.\,14.31.

\hf \raisebox{-1ex}{\sl V.\,A.\,Vedernikov}
\emp

\bmp \textbf{16.33.} Suppose that a finite $p\hs$-group $G$ has
 an abelian subgroup $A$ of order $p^n$. Does $G$ contain an
abelian subgroup $B$ of order $p^n$ that is normal in $\langle
B^G\rangle$

\makebox[25pt][r]{a)} if $p=3$?

\makebox[25pt][r]{b)} if $p=2$?

\makebox[25pt][r]{c)} If $p=3$ and $A$ is elementary abelian, does
$G$ contain an elementary abelian subgroup $B$ of order $3^n$ that
is normal in $\langle B^G\rangle$?

\makebox[15pt][r]{}The corresponding results have been proved for
greater primes $p$, based on extensions of Thompson's Replacement
Theorem (and the dihedral group of order $32$ shows that the third
question has negative answer for $p=2$). See (G.\,Glauberman, {\it
J.~Algebra}, \textbf{196} (1997), 301--338; \ J.~Alperin,
G.\,Glauberman, {\it J.~Algebra}, \textbf{203} (1998), 533--566; \
G.\,Glauberman, {\it J.~Algebra}, \textbf{272} (2004), 128--153). \hf
\raisebox{-1ex}{\sl G.\,Glauberman}
\emp

\bmp \textbf{\zv 16.34.} Suppose that $G$ is a finitely generated
group acting faithfully on a regular rooted tree by finite-state
automorphisms. Is the conjugacy problem decidable for $G$?

\makebox[15pt][r]{}See the definitions in (R.\,I.\,Grigorchuk,
V.\,V.\,Nekrashevich,
V.\,I.\,Sushchanski\u{\i}, 
{\it Proc. Steklov Inst. Math.}, \textbf{2000}, no.\,4 (231),
128--203).

\hf \raisebox{-1ex}{\sl R.\,I.\,Grigorchuk,
V.\,V.\,Nekrashevich, V.\,I.\,Sushchanski\u{\i}}

\ul

\otv No, not always (Z.\,\v{S}uni\'{c}, E.\,Ventura,  {\it J.~Algebra}, {\bf  364} (2012), 148--154).
\emp

 \bmp \textbf{\zv 16.38.}
 Let $G$ be a soluble group, and let $A$ and $B$ be periodic
subgroups of $G$. Is it true that any subgroup of $G$ contained in
the set $AB=\{ ab\;|\; a\in A,\; b\in B\}$ is periodic? This is
known to be true if $AB$ is a subgroup of $G$. \hf
\raisebox{-1ex}{\sl F.\,de\,Giovanni}

\ul

\otv Yes, it is true (AI solution confirmed by M.\,De~Falco: \  M.\,Kinyon, J.\,Urban, {\it Preprint}, 2026,
\url{https://kourovkanotebookorg.wordpress.com/wp-content/uploads/2026/09/kourovka-mathematics.pdf}).
\emp 

 \bmp \textbf{16.39.} (J.\,E.\,Humphreys, D.\,N.\,Verma).
Let $G$ be a semisimple algebraic group over an algebraically
closed field $k$ of characteristic $p>0$. Let ${\frak g}$ be the
Lie algebra of $G$ and let $u=u({\frak g})$ be the restricted
enveloping algebra of ${\frak g}$. By a theorem of Curtis every
irreducible restricted $u$-module (i.\,e.~every irreducible
restricted ${\frak g}$-module)
 is the restriction to ${\frak g}$ of a (rational) $G$-module. Is
it also true that every projective indecomposable $u$-module is
the restriction of a rational $G$-module? This is true if $p\geq
2h-2$ (where $h$ is the Coxeter number of $G$) by results of
Jantzen. \hf \raisebox{-1ex}{\sl S.\,Donkin}

\emp

\bmp \textbf{16.40.} Let $\Delta$ be a subgroup of the automorphism
group of a free pro-$p$ group of finite rank $F$ such that
$\Delta$ is isomorphic (as a profinite group) to the group $Z_p$
of $ p\hs$-adic integers. Is the subgroup of fixed points of
$\Delta$ in $F$ finitely generated (as a profinite group)?\hf
\raisebox{-1ex}{\sl P.\,A.\,Zalesski\u{\i}}

\emp

\bmp \textbf{16.41.} Let $F$ be a free pro-$p$ group of finite rank
$n>1$. Does ${\rm Aut}\,F$ possess an open subgroup of finite
cohomological dimension?\hf \raisebox{-1ex}{\sl
P.\,A.\,Zalesski\u{\i}}

\emp

\bmp \textbf{16.44.} Is there in ZFC a countable non-discrete
topological group not containing discrete subsets with a single
accumulation point? (Such a group is known to exist under Martin's
Axiom.)\hf \raisebox{-1ex}{\sl E.\,G.\,Zelenyuk}

\emp

 \bmp
 \textbf{\zv 16.45.} Let $G$ be a permutation group on a set $\Omega$. A sequence of
points of $\Omega$ is a {\it base\/} for $G$ if its pointwise
stabilizer in $G$ is the identity; it is {\it minimal\/} if no
point may be removed. Let $b(G)$ be the maximum, over all
permutation representations of the finite group $G$, of the
maximum size of a minimal base for $G$. Let $\mu'(G)$ be the
maximum size of an {\it independent set\/} in $G$, a set of
elements with the property that no element belongs to the subgroup
generated by the others. Is it true that $b(G)=\mu'(G)$? (It is
known that $b(G)\le\mu'(G)$, and that equality holds for the
symmetric groups.)

\makebox[15pt][r]{}{\it Remark.} An equivalent question is the
following. Suppose that the Boolean lattice $B(n)$ of subsets of
an $n$-element set is embeddable as a meet-semilattice of the
subgroup lattice of $G$, and suppose that $n$ is maximal with this
property. Is it true that then there is such an embedding of
$B(n)$ with the property that the least element of $B(n)$ is a
normal subgroup of $G$? \hf \raisebox{-1ex}{\sl P.\,J.\,Cameron}

\ul

\otv No, not always (AI solution confirmed by P.\,J.\,Cameron:  \ A.\,Kriebel,  {\it Preprint}, 2026, \url{https://zenodo.org/records/22929486}).
\emp

 \bmp \textbf{16.46.} Among the finitely-presented groups that act arc-transitively on
the (infinite) 3-valent tree with finite vertex-stabilizer are the
two groups $G_3 = \left< h,\,a,\,P,\,Q \mid h^3\right.$, $ a^2$, $
P^2$, $ Q^2$, $ [P,Q]$, $ [h,P]$, $ (hQ)^2$, $\left. a^{-1}PaQ
\right>$ \ and \ $G_4 = \left< h,\,a,\,p,\,q,\,r \mid
h^3\right.$, $ a^2$, $ p^2$, $ q^2$, $ r^2$, $ [p,q]$, $ [p,r]$, $
p(qr)^2$, $ h^{-1}phq$, $ h^{-1}qhpq$, $ (hr)^2$, $ [a,p]$,
$\left. a^{-1}qar \right>$, each of which contains the modular
group $G_1 = \left< h,a \mid h^3,\, a^2 \right> \cong PSL_2({\Bbb
Z})$ as a subgroup of finite index. The free product of $G_3$ and
$G_4$ with subgroup $G_1 = \left< h,a \right>$ amalgamated has a
normal subgroup $K$ of index 8 generated by $A = h$, \ $B = aha$,
\ $C = p$, \ $D = PpP$, \ $E = QpQ$, and $F = PQpPQ$, with
dihedral complement $\left< a,\,P,\,Q \right>$. The group $K$ has
presentation $\left< A,\,B,\,C,\,D,\,E,\,F \mid A^3\right.$, $
B^3$, $ C^2$, $ D^2$, $ E^2$, $ F^2$, $ (AC)^3$, $ (AD)^3$, $
(AE)^3$, $ (AF)^3$, $ (BC)^3$, $ (BD)^3$, $ (BE)^3$, $ (BF)^3$, $
(ABA^{-1}C)^2$, $ (ABA^{-1}D)^2$, $ (A^{-1}BAE)^2$, $
(A^{-1}BAF)^2$, $ (BAB^{-1}C)^2$, $ (B^{-1}ABD)^2$, $
(BAB^{-1}E)^2$, $\left. (B^{-1}ABF)^2 \right>$. Does this group
have a non-trivial finite quotient? \hf \raisebox{-1ex}{\sl
M.\,Conder}
\emp

\bmp \textbf{16.47.} (P.\,Conrad).
Is it true that every torsion-free
abelian group admits an Archi\-me\-dean lattice ordering?\hfill
\raisebox{-1ex}{\sl V.\,M.\,Kopytov, N.\,Ya.\,Medvedev}
\emp

\bmp \textbf{\zv 16.48.}
A group $H$ is said to have {\it generalized
torsion\/} if there exists an element $h\ne 1$ such that
$h^{x_1}h^{x_2}\cdots h^{x_n}= 1$ for some $n$ and some $x_i\in
H$. Is every group without generalized torsion right-orderable?
\hf \raisebox{-1ex}{\sl V.\,M.\,Kopytov, N.\,Ya.\,Medvedev}

\ul

\otv No, not every (T.\,W.\,Cai, A.\,Clay, {\it Preprint}, 2025, \url{https://arxiv.org/abs/2504.08084})
\emp

\bmp \textbf{\zv 16.51.}
Do there exist groups that can be
right-ordered in infinitely countably many ways?

\makebox[15pt][r]{}There exist groups that can be fully ordered in
infinitely countably many ways (R.\,Buttsworth).\hf
\raisebox{-1ex}{\sl V.\,M.\,Kopytov, N.\,Ya.\,Medvedev}

\ul

\otv No, there are no such groups (A.\,Clay, \emph{Preprint}, 2009, \url{http://arxiv.org/abs/0909.0273}; \ A.\,Navas, C.\,Rivas, \emph{Algebr. Geom. Topol.}, \textbf{9} (2009), 2079--2100, with an appendix by A.\,Clay; \ P.\,A.\,Linnell, \emph{Bull. London Math. Soc.}, \textbf{43} (2011), 200--202, based on (A.\,V.\,Zenkov, \emph{Sibirsk. Mat. Zh.}, \textbf{38} (1997), 90--92 (Russian), English transl., \emph{Siberian Math.~J.}, \textbf{38}, no.\,1 (1997), 76--77), where this
was proved for locally indicable groups.
\emp

\bmp \textbf{\zv 16.53.}
Let $d(G)$ denote the smallest cardinality of a
generating set of the group~$G$. Suppose that
$G=\left<A,B\right>$, where $A$ and $B$ are two $d$-generated
finite groups of coprime orders. Is it true that $d(G)\leq d+1$?

\makebox[15pt][r]{}See Archive~12.71 and (A.\,Lucchini, {\it
J.~Algebra}, \textbf{245} (2001), 552--561).

\makebox[15pt][r]{}{\it Comment of 2021}:
the answer is yes if $A$ and $B$ have prime power order (E.\,Detomi, A.\,Lucchini, \emph{J.~London Math. Soc. (2)}, \textbf{87}, no.\,3 (2013), 689--706).
\hf \raisebox{-1ex}{\sl A.\,Lucchini}

\ul

\otv No, not always (K.\,Muliarchyk, {\it Preprint}, 2026, \url{https://kourovkanotebookorg.wordpress.com/wp-content/uploads/2026/09/kourovka_16_53_solution.pdf}).
\emp

\bmp \textbf{16.56.} The {\it spectrum\/} $\omega(G)$ of a group
$G$ is the set of orders of elements of $G$. Suppose that
$\omega(G)=\{1,2,3,4,5,6\}$. Is $G$ locally finite?\hf
\raisebox{-1ex}{\sl V.\,D.\,Mazurov}
\emp

\bmp \textbf{\zv 16.60.}
If $G$ is a finite group, let $T(G)$ be
the sum of the degrees of the irreducible complex representations
of $G$, $T(G)=\sum\limits_{i=1}^{k(G)}d_i$, where $G$ has $k(G)$
conjugacy classes. If $\alpha \in {\rm Aut\,} G$, let $S_{\alpha
}=\{g\in G\mid \alpha (g)=g^{-1}\}$. Is it true that $T(G)\geq
|S_{\alpha }|$ for all $\alpha \in {\rm Aut\,} G$?\hf
\raisebox{-1ex}{\sl D.\,MacHale}

\ul

\oai Yes, it is true (AI solution confirmed by T.\,Wilde: \ A.\,Ferudun, {\it Preprint}, 2026, \url{https://doi.org/10.5281/zenodo.23445966}).
\emp 

\bmp \textbf{\ai 16.63.} 
Is there a non-trivial finite $p$-group $G$
of odd order such that $|{\rm Aut\,} G| =|G|$?

\makebox[15pt][r]{}See also Archive 12.77.\hf \raisebox{-1ex}{\sl
D.\,MacHale}

\ul

\oai Yes, there is (unverified AI solution:  A.\,Kriebel,  {\it Preprint}, 2026,

\url{https://kourovkanotebookorg.wordpress.com/wp-content/uploads/2026/09/kourovka_16_63_v1.1.0.pdf}).
\emp

\bmp \textbf{16.64.}
A non-abelian variety in which all finite
groups are abelian is called {\it pseudo-abelian}. A group variety
is called a {\it $t$-variety\/} if for all groups in this variety
the relation of being a normal subgroup is transitive. By
(O.\,Macedo\'nska, A.\,Storozhev, {\it Commun. Algebra}, \textbf{25},
no.\,5 (1997), 1589--1593) each non-abelian $t$-variety is
pseudo-abelian, and the pseudo-abelian varieties constructed in
(A.\,Yu.\,Olshanskii, {\it Math. USSR--Sb.},
\textbf{54} (1986), 57--80) 
are $t$-varieties. Is every pseudo-abelian variety a $t$-variety?

\hf \raisebox{-1ex}{\sl O.\,Macedo\'nska}

\emp

 \bmp
\textbf{16.66.} For a group $G$, let $D_n(G)$ denote the $n$-th
dimension subgroup of $G$, and $\zeta_n(G)$ the $n$-th term of its
upper central series. For a given integer $n\geq 1$, let
$f(n)=\max\{ m\mid \exists\;\text{a nilpotent group}\;G\;\text{of
class}\;n\;\text{with}\;D_m(G)\neq 1\}$ and $g(n)=\max\{ m\ |\
\exists\ \text{a nilpotent group}\ G\ \text{of\ class}\ n\
\text{such that}\ D_n(G)\not\subseteq \zeta_m(G)\}.$

\makebox[15pt][r]{}a) What is $f(3)$?

\makebox[15pt][r]{}b) (B.\,I.\,Plotkin). Is it true that $f(n)$ is
finite for all $n$?

\makebox[15pt][r]{}c) Is the growth of $f(n)$ and $g(n)$
polynomial, exponential, or intermediary?

\makebox[15pt][r]{}It is known that both $f(n)-n$ and $g(n)$ tend
to infinity as $n\rightarrow \infty$ (N.\,Gupta, Yu.\,V.\,Kuz'min,
{\it J.~Pure Appl. Algebra}, \textbf{104} (1995), 191--197).

\hf
\raisebox{-1ex}{\sl R.\,Mikhailov, I.\,B.\,S.\,Passi} \emp

\bmp \textbf{\zv 16.68.}
Let $W(x, y)$ be a non-trivial reduced group
word, and $G$ one of the groups $PSL(2,{\Bbb R})$, \ $PSL(2,
{\Bbb C})$, or $SO(3, {\Bbb R})$. Are all the maps $W:G \times
G\rightarrow G$ surjective?

\makebox[15pt][r]{}For $SL(2, {\Bbb R})$, \ $SL(2, {\Bbb C})$, \
$GL(2, {\Bbb C})$, and for the group $S^3$ of quaternions of norm
1 there exist non-surjective words; see (J.\,Mycielski, {\it Amer.
Math. Monthly}, \textbf{84} (1977), 723--726; \ \textbf{85} (1978),
263--264).\hf \raisebox{-1ex}{\sl J.\,Mycielski}

\ul

\otv These maps are not surjective for $PSL(2,{\Bbb R})$ and $SO(3, {\Bbb R})$, and surjective for $PSL(2, {\Bbb C})$ (AI solution confirmed by A.\,Thom: \   Nilradical,  {\it Preprint}, 2026, \url{https://kourovkanotebookorg.wordpress.com/wp-content/uploads/2026/09/nilradical-kourovka-16-68-69.pdf}). For $SO(3, {\Bbb R})$ the negative answer follows from  (A.\,Thom, {\it Canad. Math. Bull.}, {\bf 56} (2013), 424--433).
\emp 

\bmp \textbf{\zv 16.69.}
Let $W$ be as above.

\makebox[25pt][r]{a)} Must the range of the function ${\rm
Tr}(W(x, y))$ for $x, y\in GL(2, {\Bbb R})$ include the interval
$[-2, +\infty)$?

\makebox[25pt][r]{b)} For $x, y$ being non-zero quaternions,
 must the range of the function ${\rm Re}(W(x, y))$ include the interval
$[-5/27, 1]$?

\makebox[15pt][r]{}(For some comments see {\it ibid}.)\hf
\raisebox{0ex}{\sl J.\,Mycielski}

\ul

\otv a) No, it need not (AI solution confirmed by A.\,Thom: \  Nilradical,  {\it Preprint}, 2026,

\url{https://kourovkanotebookorg.wordpress.com/wp-content/uploads/2026/09/nilradical-kourovka-16-68-69.pdf}).

\otv b) No, it need not  (A.\,Thom, {\it Canad. Math. Bull.}, {\bf 56} (2013), 424--433).
\emp 

\bmp \textbf{16.70.} Suppose that a finitely generated group $G$ acts
freely on an $A$-tree, where $A$ is an ordered abelian group. Is
it true that $G$ acts freely on a ${\Bbb Z}^n$-tree for some
$n$?

\makebox[15pt][r]{}{\it Comment of 2013:}
It was proved that every finitely presented group acting freely on an $A$-tree acts freely on some $\mathbb{R}^n$-tree for a suitable $n$, where $\mathbb{R}$ has the lexicographical order (O.\,Kharlampovich, A.\,Myasnikov, D.\,Serbin, {\it Int. J. Algebra Comput.}, {\bf 23} (2013), 325--345).
\hf \raisebox{-1ex}{\sl A.\,G.\,Myasnikov,
V.\,N.\,Remeslennikov, O.\,G.\,Kharlampovich}

\emp

\bmp \textbf{16.72.} Does there exist an exponential-time
algorithm for obtaining a JSJ-decom\-po\-si\-tion of a finitely
generated fully residually free group?

\makebox[15pt][r]{}Some algorithm was found in (O.\,Kharlampovich,
A.\,Myasnikov, in: {\it Contemp. Math. AMS $($Algorithms,
Languages, Logic\/}), \textbf{378} (2005), 87--212).

\hf\raisebox{-1ex}{\sl A.\,G.\,~Myasnikov,
O.\,G.\,Kharlampovich}

\emp

\bmp \textbf{16.73.} b) Let $G$ be a group generated by a finite set
$S$, and let $l(g)$ denote the word length function of $g\in G$
with respect to $S$. The group $G$ is said to be {\it
contracting\/} if there exist a faithful action of $G$ on the set
$X^*$ of finite words over a finite alphabet $X$ and constants
$0<\lambda <1$ and $C>0$ such that for every $g\in G$ and $x\in X$
there exist $h\in G$ and $y\in X$ such that $l(h)<\lambda l(g)+C$
and $g(xw)=yh(w)$ for all $w\in X^*$.

\makebox[15pt][r]{}Do there exist non-amenable contracting
groups? \hf \raisebox{-1ex}{\sl V.\,V.\,Nekrashevich}
\emp

\bmp \textbf{16.74.} b) Is it true that all groups
generated by automata of polynomial growth in the sense of
S.\,Sidki ({\it Geom. Dedicata}, \textbf{108} (2004), 193--204) are
amenable?

\hf \raisebox{-1ex}{\sl V.\,V.\,Nekrashevich}
\emp

\bmp \textbf{16.76.} We call a group $G$ {\it strictly real\/} if
each of its non-trivial elements is conjugate to its inverse by
some involution in~$G$. In which groups of Lie type over a field
of characteristic 2 the maximal unipotent subgroups are
strictly real?\hf \raisebox{-1ex}{\sl Ya.\,N.\,Nuzhin}
\emp

\bmp \textbf{16.77.} It is known that in every Noetherian group the
nilpotent radical coincides with the collection of all Engel
elements (R.\,Baer, {\it Math. Ann.}, \textbf{133} (1957), 256--270;
B.\,I.\,Plotkin, {\it Izv. Vys\v{s}. U\v{c}ebn. Zaved. Mat.}, {\bf
1958}, no.\,1(2), 130--135 (Russian)). It would be nice to find a
similar characterization of the solvable radical of a finite
group. More precisely, let $u=u(x,y)$ be a sequence of words
satisfying 15.75. We say that an element $g\in G$ is $u$-Engel if
there exists $n=n(g)$ such that $u_n(x,g)=1$ for every element
$x\in G$. Does there exist a sequence $u=u(x,y)$ such that the
solvable radical of a finite group coincides with the set of all
$u$-Engel elements? \hf \raisebox{-1ex}{\sl B.\,I.\,Plotkin}

\emp

\bmp \textbf{\ai 16.78.} 
Do there exist linear non-abelian simple
groups without involutions?

\hf \raisebox{-1ex}{\sl B.\,Poizat}

\ul

\oai Yes, such groups do exist (unverified AI solution:  I.\,S.\,Nnuga Mongo, {\it Preprint}, 2026, \url{https://arxiv.org/abs/2610.37703}).
\emp 

\bmp \textbf{16.80.} Suppose that a group $G$ is obtained from the
free product of torsion-free groups $A_1, \ldots , A_n$ by
imposing $m$ additional relations, where $m<n$. Is it true that
the free product of some $n-m$ of the $A_i$ embeds into~$G$? \hf
\raisebox{-1ex}{\sl N.\,S.\,Romanovski\u{\i}}
\emp

\bmp \textbf{16.83.} b) Let $E_n$ be a free locally nilpotent $n$-Engel
group on countably many generators, and let $\pi(E_n)$ be the set
of prime divisors of the orders of elements of the periodic part
of $E_n$. It is known that $2,\,3,\,5\in \pi(E_4)$.

\makebox[15pt][r]{}Is it true that $\pi(E_n)=\pi(E_{n+1}) $
for all sufficiently large~$ n$? \hf \raisebox{-1ex}{\sl
Yu.\,V.\,Sosnovski\u{\i}}
\emp

\bmp \textbf{16.84.} Can the braid group ${\frak B}_n$, $n\geq 4$,
act faithfully on a regular rooted tree by finite-state
automorphisms? Such action is known for ${\frak B}_3$.

\makebox[15pt][r]{}See the definitions in (R.\,I.\,Grigorchuk,
V.\,V.\,Nekrashevich,
V.\,I.\,Sushchanski\u{\i}, 
{\it Proc. Steklov Inst. Math.}, \textbf{2000}, no.\,4 (231),
128--203).\hf \raisebox{-1ex}{\sl V.\,I.\,Sushchanski\u{\i}}

\emp

\bmp \textbf{\zv 16.87.}
Let ${\frak M}$ be a variety of groups and
let $G_r$ be a free $r$-generator group in $\frak M$. A subset
$S\subseteq G_r$ is called a {\it test set\/} if every
endomorphism of $G_r$ identical on
 $S$ is an automorphism. The minimum of the cardinalities of test
sets is called the test rank of $G_r$. Suppose that the test rank
of $G_r$ is $r$ for every $r \geq 1$.

\makebox[15pt][r]{}a) Is it true that ${\frak M}$ is an abelian
variety?

\makebox[15pt][r]{}b) Suppose that ${\frak M}$ is not a periodic
variety. Is it true that ${\frak M}$ is the variety of all abelian
groups?\hf \raisebox{-1ex}{\sl E.\,I.\,Timoshenko}

\ul

\otv a), b) No, it is not true in both cases (AI solution confirmed by E.\,I.\,Timoshenko: \  M.\,Kinyon, J.\,Urban, {\it Preprint}, 2026,
\url{https://kourovkanotebookorg.wordpress.com/wp-content/uploads/2026/09/kourovka-mathematics.pdf}).
\emp

\bmp \textbf{16.88.} (G.\,M.\,Bergman). The {\it width} of a group $G$
with respect to a generating set $X$ means the supremum over all
$g\in G$ of the least length of a group word in elements of $X$
expressing $g$. A group $G$ has {\it finite width} if the width
of $G$ with respect to every
generating set is
finite. Does there exist a countably infinite group of finite
width?

\makebox[15pt][r]{}All known infinite groups of finite width (infinite
permutation groups, infinite-dimen\-si\-o\-nal general linear
groups, and some other groups) are uncountable
(G.\,M.\,Bergman,
{\it Bull. London Math. Soc.}, {\bf 38} (2006), 429--440, and
references therein).
 \hfill \raisebox{-1ex}{\sl V.\,Tolstykh}

\emp

\bmp \textbf{16.89.} (G.\,M.\,Bergman). Is it true that the automorphism
group of an infinitely generated free group $F$ has finite width?
The answer is affirmative if $F$ is countably generated. \hf
\raisebox{-1ex}{\sl V.\,Tolstykh}

\emp

\bmp \textbf{16.90.} Is it true that the automorphism group
$\text{Aut}\,F$ of an infinitely generated free group $F$ is

\makebox[15pt][r]{}a) the normal closure of a single element?

 \makebox[15pt][r]{}b) the normal closure of some involution in $\text{Aut}\,F$?

 \makebox[15pt][r]{}For a free group $F_n$ of
finite rank $n$ \ M.\,Bridson and K.\,Vogtmann
have recently shown
that $\text{Aut}\,F_n$ is the normal closure of some involution,
which permutes some basis of $F_n$. It is also known that the
automorphism groups of infinitely generated free nilpotent groups
have such involutions.
 \hf \raisebox{-1ex}{\sl V.\,Tolstykh}
\emp

\bmp \textbf{16.91.} Let $F$ be an infinitely generated free group.
Is there an $IA$-automorphism of $F$ whose normal closure in
$\text{Aut}\,F$ is the group of all $IA$-automorphisms of~$F?$

 \hf
\raisebox{-1ex}{\sl V.\,Tolstykh} \emp

\bmp \textbf{16.92.} Let $F$ be an infinitely generated free group.
Is $\text{Aut}\,F$ equal to its derived subgroup? This is true if
 $F$ is countably generated (R.\,Bryant,
V.\,A.\,Roman'kov, {\it J.~Algebra}, \textbf{209} (1998), 713--723).
 \hf \raisebox{-1ex}{\sl V.\,Tolstykh}

\emp

\bmp \textbf{16.93.} Let $F_n$ be a free group of finite rank $n \geq
2$. Is the group $\text{Inn}\,F_n$ of inner automorphisms of $F_n$
a first-order definable subgroup of $\text{Aut}\,F_n$?
 It is known that the set of inner automorphisms induced by
powers of primitive elements is definable in $\text{Aut}\,F_n$.
 \hf \raisebox{-1ex}{\sl V.\,Tolstykh}

\emp

\bmp \textbf{16.94.} If $G=[G,G]$, then the {\it commutator width\/}
of the group $G$ is its width relative to the set of commutators.
Let $V$ be an infinite-dimensional vector space over a division
ring. It is known that the commutator width of $GL(V)$ is finite.
Is it true that the commutator width of $GL(V)$ is one?\hf
\raisebox{-1ex}{\sl V.\,Tolstykh}

\emp

\bmp \textbf{16.95.}
{\it Conjecture:\/} If $F$ is a field and
$A$ is in $GL(n,F)$, then there is a permutation matrix $P$ such
that $AP$ is cyclic, that is, the minimal polynomial of $AP$ is
also its characteristic polynomial. \hf \raisebox{-1ex}{\sl
J.\,G.\,Thompson}
\emp

\bmp \textbf{16.96.} Let $G$ be a locally finite $n$-Engel
$p\hs$-group where $p$ a prime greater than~$n$. Is $G$ then a
Fitting group?
 (Examples of N.\,Gupta and F.\,Levin show that the
condition $p>n$ is necessary in general.) \hf \raisebox{-1ex}{\sl
G.\,Traustason}

\emp

\bmp \textbf{16.97.} Let $G$ be a torsion-free group with all
subgroups subnormal of defect at most~$n$. Must $G$ then be
nilpotent of class at most~$n$? (This is known to be true for $n <
5$). \hf \raisebox{-1ex}{\sl G.\,Traustason}

\emp

\bmp \textbf{16.98.} Suppose that $G$ is a solvable finite group
and $A$ is a group of automorphisms of $G$ of relatively prime
order. Is there a bound for the Fitting height $h(G)$ of $G$ in
terms of $A$ and $h(C_G(A))$, or even in terms of the length
$l(A)$ of the longest chain of nested subgroups of $A$ and
$h(C_G(A))$?

\makebox[15pt][r]{}When $A$ is solvable, it is proved in
(A.\,Turull, {\it J.~Algebra}, \textbf{86} (1984), 555--566) that
$h(G)\leq h(C_G(A)) + 2 l(A)$ and this bound is best possible for
$h(C_G(A))> 0$. For $A$ non-solvable some results are in
(H.\,Kurzweil, {\it Manuscripta Math.}, \textbf{41} (1983),
233--305).\hf \raisebox{-1ex}{\sl A.\,Turull}

\emp

\bmp \textbf{16.99.} Suppose that $G$ is a finite solvable group, $A
\leq {\rm Aut}\, G$, \ $C_G(A)=1$, the orders of $G$ and $A$ are
coprime, and let $l(A)$ be the length of the longest chain of
nested subgroups of $A$. Is the Fitting height of $G$ bounded
above by~$l(A)$?

\makebox[15pt][r]{}For $A$ solvable the question coincides
with~5.30. It is proved that for any finite group $A$, first,
there exist $G$ with~$h(G)=l(A)$ and, second, there is a finite
set of primes $\pi$ (depending on $A$) such that if $|G|$ is
coprime to each prime in $\pi$, then $h(G)\leq l(A)$. See
(A.\,Turull, {\it Math.~Z.}, \textbf{187} (1984), 491--503). \hf
\raisebox{-1ex}{\sl A.\,Turull}

\emp \bmp \textbf{16.100.} Is there an (infinite) 2-generator
simple group $G$ such that ${\rm Aut}\,F_2$ is transitive on the
set of normal subgroups $N$ of the free group $F_2$ such that
$F_2/N\cong G$? Cf.~6.45.

\hf \raisebox{-1ex}{\sl J.\,Wiegold}
\emp

 \bmp \textbf{16.102.} We say that a group $G$ is {\it rational\/}
 if any two elements $x,y\in
G$ satisfying $\left< x \right>=\left< y \right>$ are conjugate.
Is it true that for any $d$ there exist only finitely many finite
rational $d$-generated $2$-groups?
\hf \raisebox{-1ex}{\sl A.\,Jaikin-Zapirain}
\emp

\bmp \textbf{16.105.} Is it true that a locally graded group which
is a product of two almost polycyclic subgroups (equivalently, of
two almost soluble subgroups with the maximal condition) is almost
polycyclic?

\makebox[15pt][r]{}By (N.\,S.\,Chernikov, {\it 
Ukrain. Math. J.}, \textbf{32} (1980) 476--479) a locally graded
group which is a product of two Chernikov subgroups (equivalently,
of two almost soluble subgroups with the minimal condition) is
Chernikov.\hf \raisebox{-1ex}{\sl N.\,S.\,Chernikov}
\emp

\bmp \textbf{16.108.} Do braid groups $B_n,\; n>4$, have
non-elementary hyperbolic factor groups?

\hf \raisebox{-1ex}{\sl V.\,E.\,Shpilrain}
\emp

\bmp \textbf{16.110.} (I.\,Kapovich, P.\,Schupp). Is there an
algorithm which, when given two elements $u, v$ of a free group
$F_n$, decides whether or not the cyclic length of $\phi(u)$
equals the cyclic length of $\phi(v)$ for every automorphism
$\phi$ of~$F_n$?

\makebox[15pt][r]{}{\it Comment of 2009:\/} The answer was shown
to be positive for $n=2$ in (D.\,Lee, {\it J.~Group Theory}, {\bf
10} (2007), 561--569).
 \hf \raisebox{-1ex}{\sl V.\,E.\,Shpilrain}
\emp

\bmp \textbf{16.111.} Must an infinite simple periodic group with a
dihedral Sylow 2-sub\-group be isomorphic to $L_2(P)$ for a
locally finite field $P$ of odd characteristic? \hfill
\raisebox{-1ex}{\sl V.\,P.\,Shunkov}

\emp

 \raggedbottom

\newpage

\pagestyle{myheadings} \markboth{17th Issue (2010)}{17th Issue
(2010)} \thispagestyle{headings} ~ \vspace{2ex}

\centerline {\Large \textbf{Problems from the 17th Issue (2010)}}
\phantomsection\label{17izd}
\vspace{4ex}

\bmp \textbf{17.1.} (I.\,M.\,Isaacs). Does there exist a finite group
partitioned by subgroups of equal order not all of which are
abelian? (Cf.~15.26.) It is known that such a group must be of
prime exponent (I.\,M.\,Isaacs, {\it Pacific J.~Math.}, \textbf{49}
(1973), 109--116).\hf \raisebox{-1ex}{\sl A.\,Abdollahi}
 \emp

\bmp \textbf{17.4.} Let $x$ be a right $4$-Engel element of a group
$G$.

\makebox[15pt][r]{}a) Is it true that the normal closure $\langle
x \rangle ^G$ of $x$ in $G$ is nilpotent if $G$ is locally
nilpotent?

\makebox[15pt][r]{}b) If the answer to a) is affirmative, is there
a bound on the nilpotency class of~$\langle x \rangle ^G$?

\makebox[15pt][r]{}c) Is it true that $\langle x \rangle ^G$ is
always nilpotent? \hfill \raisebox{0ex}{\sl A.\,Abdollahi}
 \emp

\bmp \textbf{17.5.} Is the nilpotency class of a nilpotent group
generated by $d$ left $3$-Engel elements bounded in terms of~$d$?
A group generated by two left $3$-Engel elements is nilpotent of
class at most $4$ ({\it J.~Pure Appl. Algebra}, \textbf{188} (2004),
1--6.)\hf \raisebox{-1ex}{\sl A.\,Abdollahi}
 \emp

 \bmp \textbf{17.6.}
a) Is there a function $f:\mathbb{N}\times\mathbb{N}\rightarrow
\mathbb{N}$ such that every nilpotent group generated by $d$ left
$n$-Engel elements is nilpotent of class at most $f(n,d)$?

\makebox[15pt][r]{}b) Is there a function
$g:\mathbb{N}\times\mathbb{N}\rightarrow \mathbb{N}$ such that
every nilpotent group generated by $d$ right $n$-Engel elements is
nilpotent of class at most $g(n,d)$? \hf \raisebox{-1ex}{\sl
A.\,Abdollahi}
 \emp

 \bmp \textbf{17.7.} a) Is there a group in which the set of right Engel elements does
not form a subgroup?

\makebox[15pt][r]{}b) Is there a group in which the set of
bounded right Engel elements does not form a subgroup?\hf
\raisebox{-1ex}{\sl A.\,Abdollahi}
 \emp

 \bmp \textbf{17.8.} a) Is there a group containing a right Engel element which is not
a left Engel element?

\makebox[15pt][r]{}b) Is there a group containing a bounded right
Engel element which is not a left Engel element?\hf
\raisebox{-1ex}{\sl A.\,Abdollahi}
 \emp

 \bmp \textbf{17.9.} Is there a group containing a left
Engel element whose inverse is not a left Engel element?\hf
\raisebox{-1ex}{\sl A.\,Abdollahi}
 \emp

 \bmp \textbf{17.10.} Is there a group containing a right
$4$-Engel element which does not belong to the Hirsch--Plotkin
radical?\hf \raisebox{-1ex}{\sl A.\,Abdollahi}
 \emp

 \bmp \textbf{17.11.} Is there a group containing a left
$3$-Engel element which does not belong to the Hirsch--Plotkin
radical?\hf \raisebox{-1ex}{\sl A.\,Abdollahi}
 \emp

\bmp \textbf{17.13.} Let $G$ be a totally imprimitive $p$-group of
finitary permutations on an infinite set. Suppose that the support
of any cycle in the cyclic decomposition of every element of $G$
is a block for $G$. Does $G$ necessarily contain a minimal
non-$FC$-subgroup?

\hf \raisebox{-1ex}{\sl A.\,O.\,Asar}
 \emp

\bmp \textbf{17.15.}
Construct an algorithm which, for a given polynomial $f \in
{\mathbb Z}[x_1, x_2, \ldots, x_n]$, finds (explicitly, in terms of
generators) the maximal subgroup $G_f$ of the group $Aut ({\mathbb
Z}[x_1, x_2, \ldots, x_n])$ that leaves $f$ fixed. This question
is related to description of the solution set of the Diophantine
equation $f = 0$.\hf \raisebox{-1ex}{\sl V.\,G.\,Bardakov}
 \emp

\bmp \textbf{17.16.}
Let $A$ be an Artin group of finite type, that is, the
corresponding Coxeter group is finite. It is known that the group
$A$ is linear (S.\,Bigelow, {\it J.~Amer. Math. Soc.}, \textbf{14} 
(2001), 471--486; \ D.\,Krammer, {\it Ann. Math.}, \textbf{155} 
(2002), 131--156; \ F.\,Digne, {\it J.~Algebra}, \textbf{268} 
(2003), 39--57; \ A.\,M.\,Cohen, D.\,B.\,Wales, {\it Israel
J.~Math.}, \textbf{ 131} (2002), 101--123). Is it true that the
automorphism group $Aut (A)$ is also linear?

\makebox[15pt][r]{}An affirmative answer for $A$ being a braid
group $B_n$, $n \geq 2$, was obtained in (V.\,G.\,Bardakov,
O.\,V.\,Bryukhanov, {\it Vestnik Novosibirsk. Univ. Ser. Mat. Mekh.
Inform.}, \textbf{7}, no.\,3 (2007), 45--58 (Russian)).
 \hfill
\raisebox{-1ex}{\sl V.\,G.\,Bardakov}
 \emp

\bmp \textbf{17.18.} Let $\mathbf{A}$ be the class of compact groups
$A$ with the property that whenever two compact groups $B$ and $C$
contain $A$, they can be embedded in a common compact group $D$ by
embeddings agreeing on $A$. I showed ({\it Manuscr. Math.}, {\bf
58} (1987) 253--281) that all members of $\mathbf{A}$ are (not
necessarily connected) finite-dimensional compact Lie groups
satisfying a strong ``local simplicity'' property, and that all
finite groups do belong to $\mathbf{A}$.

\makebox[15pt][r]{}a)~Is it true that
$\mathbb{R}/\mathbb{Z}\in\mathbf{A}$?

\makebox[15pt][r]{}b)~Do any {\it nonabelian connected} compact
Lie groups belong to $\mathbf{A}$?

\makebox[15pt][r]{}c)~If $A$ belongs to $\mathbf{A}$, must the
connected component of the identity in $A$ belong to $\mathbf{A}$?
\hf \raisebox{-1ex}{\sl G.\,M.\,Bergman}

 \emp

\bmp \textbf{\ai 17.22.} 
Suppose $A$ is a group which belongs to a
variety $\mathbf V$ of groups, and which is embeddable in the full
symmetric group $S$ on an infinite set. Must the coproduct in
$\mathbf V$ of two copies of $A$ also be embeddable in $S$?
(N.\,G.\,de~Bruijn proved that this is true if $\mathbf V$ is the
variety of all groups.)\hf \raisebox{-1ex}{\sl G.\,M.\,Bergman}

 \ul

 \oai  No, it need not (unverified AI solution:  Achyuth Jayadevan,  {\it Preprint},  2026,

\url{https://zenodo.org/records/22909835}).
 \emp 

\bmp \textbf{17.23.} Suppose the full symmetric group $S$ on a
countably infinite set is generated by the union of two subgroups
$G$ and $H$. Must $S$ be finitely generated over one of these
subgroups?\hf \raisebox{-1ex}{\sl G.\,M.\,Bergman}
 \emp

\bmp \textbf{17.25.} (S.\,P.\,Farbman). Let $X$ be the set of
complex numbers $\alpha$ such that the group generated by the two
$2\times 2$ matrices $I+\alpha e_{12}$ and $I+ e_{21}$ is not free
on those generators.

\makebox[15pt][r]{}a) Does $X$ contain all rational numbers in the
interval $(-4,4)$?

\makebox[15pt][r]{}b) Does $X$ contain any rational number in the
interval $[27/7,4)$?

\makebox[15pt][r]{}Cf. 4.41 in Archive, 15.83, and (S.\,P.\,Farbman, {\it Publ. Mat.}, \textbf{39}
(1995) 379--391).

\hf \raisebox{-1ex}{\sl G.\,M.\,Bergman}

 \emp

 \bmp \textbf{17.26.} Are the classes of right-orderable and right-ordered groups closed
under taking solutions of equations
$w(a_1,\dots,a_k,x_1,\dots,x_n)=1$? (Here the closures are under
group embeddings and order-preserving embeddings, respectively.)
This is true for equations with a single constant, when $k=1$
({\it J. Group Theory}, \textbf{11} (2008), 623--633).

\hf
\raisebox{-1ex}{\sl V.\,V.\,Bludov, A.\,M.\,W.\,Glass}
 \emp

 \bmp \textbf{\ai 17.27.} 
 Can the free product of two ordered groups with order-isomorphic
amalgamated subgroups be lattice orderable but not orderable?

\hf
\raisebox{-1ex}{\sl V.\,V.\,Bludov, A.\,M.\,W.\,Glass}

\ul

\oai Yes, it can (unverified AI solution:  Achyuth Jayadevan,  {\it Preprint}, 2026,

\url{https://zenodo.org/records/22937468}).
 \emp

 \bmp \textbf{17.30.} Is there a constant $l$ such that every finitely presented soluble group
has a subgroup of finite index of nilpotent length at most $l$?
Cf.~16.35 in Archive.

 \hf \raisebox{-1ex}{\sl V.\,V.\,Bludov, J.\,R.\,J.\,Groves}
 \emp

\bmp \textbf{17.31.} Can a soluble right-orderable group have finite
quotient by the derived subgroup? \hfill \raisebox{-1ex}{\sl
V.\,V.\,Bludov, A.\,H.\,Rhemtulla}
 \emp

 \bmp \textbf{\zv 17.32.}
 Is the following analogue of the Cayley--Hamilton theorem
 true for the free group $F_n$ of rank $n$:
 If $w\in F_n$ and $\f\in {\rm Aut\,}F_n$ are such
 that $\langle w,\f (w),\dots,\f^{n}(w)\rangle= F_n$, then $\langle w,\f (w),\dots,\f^{n-1}(w)\rangle=F_n$?
 \hf \raisebox{-1ex}{\sl O.\,V.\,Bogopolski}

 \ul

 \otv No, it is not true (V.\,Ionin, A.\,Semidetnov,  {\it Preprint}, 2026, \url{https://arxiv.org/pdf/2608.29219}; \  O.\,Kharlampovich, {\it Preprint}, 2026,
 \url{https://kourovkanotebookorg.wordpress.com/wp-content/uploads/2026/08/kourovka_14_23_17_32_note_with_problem_authors.pdf}).
 \emp

  \bmp \textbf{\zv 17.33.}
Can the quasivariety generated by the group $\langle a,b\mid
a^{-1}ba=b^{-1}\rangle$ be defined by a set of quasi-identities in
a finite set of variables? The answer is known for all other
 groups with one defining relation. \hf \raisebox{-1ex}{\sl A.\,I.\,Budkin}

 \ul

\otv No, it cannot (AI solution confirmed by A.\,I.\,Budkin: \   M.\,Kinyon, J.\,Urban, {\it Preprint}, 2026,
  \url{https://kourovkanotebookorg.wordpress.com/wp-content/uploads/2026/09/kourovka-mathematics.pdf}).
 \emp 

\bmp \textbf{\zv 17.34.}
Let $N_c$ be the quasivariety of nilpotent torsion-free groups of
class at most~$c$. Is it true that the dominion in $N_c$ (see the
definition in~16.20) of a divisible subgroup in every group
 in $N_c$ is equal to this subgroup?
\hf \raisebox{-1ex}{\sl A.\,I.\,Budkin}

\ul

\otv Yes, it is true (AI solution confirmed by A.\,I.\,Budkin: \     Achyuth Jayadevan,  {\it Preprint}, 2026,
\url{https://zenodo.org/records/22752789}).
 \emp 

\bmp \textbf{\zv 17.37.}
 Is there an integer $n$ such that for all $m
> n$ the alternating group $A_{m}$ has no
non-trivial $A_{m}$-permutable subgroups?
(See the definition in 17.112.)

 \hf \raisebox{-1ex}{\sl
A.\,F.\,Vasil'ev, A.\,N.\,Skiba}

\ul

\otv Yes, there is (N.\,Yang, A.\,Galt, {\it Bull. Malays. Math. Sci. Soc. (2)}, {\bf 46}, no.\,5 (2023), Paper no.\,177, 16 p.).
 \emp

\bmp \textbf{17.38.}
A formation $\mathfrak F$ is called {\it radical in a class
$\mathfrak H$\/} if $\mathfrak F \subseteq \mathfrak H$ and in
every $\mathfrak H$-group the product of any two normal $\mathfrak
F$-subgroups belongs to $\mathfrak F$. Let $\mathfrak M$ be the
class of all saturated hereditary formations of finite groups such
that the formation of all finite supersoluble groups is radical in
every element of $\mathfrak M$. Is it true that $\mathfrak M$ has
the largest (by inclusion) element? \hf \raisebox{-1ex}{\sl
A.\,F.\,Vasil'ev, L.\,A.\,Shemetkov}
 \emp

 \bmp \textbf{\ai 17.39.} 
Is there a positive integer $n$ such that the hypercenter of any
finite soluble group coincides with the intersection of $n$ system
normalizers of that group? What is the least number with this
property?

\makebox[15pt][r]{}{\it Comment of 2017}:
This is proved for groups with trivial Frattini subgroup (A.\,Ballester-Bolinches, J.\,Cossey, S.\,F.\,Kamornikov, H.\,Meng, {\it J.~Algebra}, \textbf{499} (2018), 438--449).
 \hf \raisebox{-1ex}{\sl A.\,F.\,Vasil'ev, L.\,A.\,Shemetkov}

 \ul

\oai No, there is no such number (unverified AI solution:  M.\,Kinyon, J.\,Urban, {\it Preprint}, 2026,
\url{https://kourovkanotebookorg.wordpress.com/wp-content/uploads/2026/09/kourovka-mathematics.pdf}).

 \emp 

\bmp \textbf{17.41.}
Let $S$ be a solvable subgroup of a finite group $G$ that has no
nontrivial solvable normal subgroups.

\makebox[25pt][r]{(a)} (L.\,Babai, A.\,J.\,Goodman, L.\,Pyber). Do
there always exist seven conjugates of $S$ whose intersection is
trivial?

\makebox[25pt][r]{(b)}
Do there always exist five conjugates of $S$
whose intersection is trivial?

\makebox[15pt][r]{}{\it Editors' comments}: Reduction of part (b) to almost simple group $G$ is obtained in (E.\,P.\,Vdovin, {\it J.~Algebra Appl.}, {\bf 11}, no.\,1 (2012), 1250015 (14 pages)). An affirmative answer to (b) has been obtained for almost simple groups with socle isomorphic to an alternating group (A.\,A.\,Baikalov, {\it Algebra Logic}, {\bf 56} (2017), 87--97), to a sporadic simple group (T.\,C.\,Burness, {\it Israel J.~Math.}, {\bf 254} (2023), 313--340), to $PSL(n,q)$ (A.\,A.\,Baykalov, {\it J.~Group Theory}, {\bf 28} (2025), 1003--1077), \ to $PSU(n,q)$ or $PSp(n,q)'$ (A.\,A.\,Baykalov, {\it Int. J.~Algebra Comput.}, {\bf 35}, no.\,06 (2025), 823--908). Part (b)
also has an affirmative answer if $S$ is a maximal subgroup of $G$ (T.\,C.\,Burness, {\it Algebra Number Theory}, {\bf 15}, no.\,7 (2021), 1755--1807).\hf \raisebox{-1ex}{\sl
E.\,P.\,Vdovin}
 \emp

\bmp \textbf{\zv 17.42.}
Let $\overline{G}$ be a simple algebraic group of adjoint type
over the algebraic closure $\overline{F}_p$ of a finite field
$F_p$ of prime order $p$, and $\sigma$ a Frobenius map (that is, a
surjective homomorphism such that
$\overline{G}_\sigma=C_{\overline{G}}(\sigma)$ is finite). Then
$G=O^{p'}(\overline{G}_\sigma)$ is a finite group of Lie type. For
a maximal $\sigma$-stable torus $\overline{T}$ of $\overline{G}$,
let $N=N_{\overline{G}}(\overline{T})\cap G$. Assume also that $G$
is simple and $G\not\cong \mathrm{SL}_3(2)$. Does there always
exist $x\in G$ such that $N\cap N^x$ is a $p$-group? \hf
\raisebox{-1ex}{\sl E.\,P.\,Vdovin}

\ul

\otv A definitive answer for a stronger question on the size of a base has been obtained in (T.\,C.\,Burness, A.\,R.\,Thomas, {\it J.~Algebra}, {\bf 619} (2023), 459--504), which implies the following answer (in the notation therein): there is $x\in G$ such that $N\cap N^x$ is a $p$-group if and only if $(G,N)$ is not one of the following: $(\operatorname L_3(2),7{:}3)$, $(\operatorname U_4(2),3^3{:}S_4)$, $(\operatorname U_5(2),3^4{:}S_5)$.
 \emp

\bmp \textbf{17.43.}
Let $\pi$ be a set of primes. Find all finite simple
$D_\pi$-groups (see Archive, 3.62) in which

 \makebox[25pt][r]{(a)} every subgroup is a $D_\pi$-group (H.\,Wielandt);

 \makebox[25pt][r]{(b)} every subgroup possessing a Hall $\pi$-subgroup is a $D_\pi$-group.

All finite simple $D_\pi$-groups are known (D.\,O.\,Revin, {\it
Algebra Logic}, \textbf{47}, no.\,3 (2008), 210--227). {\it Comment of 2013:}
Alternating and sporadic groups for both parts are found in (N.\,Ch.\,Manzaeva,
{\it Siberian Electron. Math. Rep.}, {\bf 9} (2012), 294--305 (Russian)). {\it Comment of 2025}: Simple groups of Lie type of rank 1 satisfying condition (a) are known (D.\,O.\,Revin, V.\,D.\, Shepelev, {\it Siberian Math. J.}, \textbf{65}, no.\,5 (2024), 1187--1194, V.\,D.\, Shepelev, {\it Siberian Math. J.}, \textbf{66}, no.\,4 (2025), 1049--1062).
\hf \raisebox{-1ex}{\sl E.\,P.\,Vdovin, D.\,O.\,Revin}
 \emp

\bmp \textbf{17.45.}
c) A subgroup $H$ of a group $G$ is called {\em pronormal\/} if $H$
and $H^g$ are conjugate in $\langle H,H^g\rangle$ for every $g\in
G$. We say that $H$ is {\em strongly pronormal\/} if $L^g$ is
conjugate to a subgroup of $H$ in $\langle H,L^g\rangle$ for every
$L\leq H$ and $g\in G$.

 \makebox[15pt][r]{}In a finite group, is a Hall subgroup with a Sylow tower always strongly pro\-normal?

\makebox[15pt][r]{}Notice that there
exist finite (non-simple) groups with a non-pronormal Hall
subgroup. Hall subgroups with a Sylow tower are known to be
pronormal.

 \hf\raisebox{-1ex}{\sl E.\,P.\,Vdovin, D.\,O.\,Revin}

 \emp

  \bmp \textbf{\zv 17.47.}
 Let $G$ be a nilpotent group in which every
coset $x[G,G]$ for $x\not\in [G,G]$ coincides with the conjugacy
class $x^G$. Is there a bound for the nilpotency class of~$G$?

\makebox[15pt][r]{}If
$G$ is finite, then its class is at most~3 (R.\,Dark,
C.\,Scoppola, \emph{J. Algebra}, {\bf 181}, no.\,3 (1996), 787--802).
\hfill \raisebox{-1ex}{\sl S.\,H.\,Ghate, A.\,S.\,Muktibodh}

 \ul

 \otv No, there is no such bound (V.\,Ionin, A.\,Semidetnov,  {\it Preprint}, 2026, \url{https://arxiv.org/pdf/2608.29219}).
  \emp

\bmp \textbf{17.48.} Does the free product of two groups with stable
first-order theory also have stable first-order theory? \hf
\raisebox{-1ex}{\sl E.\,Jaligot}
 \emp

\bmp \textbf{17.49.} Is it consistent with ZFC that every
non-discrete topological group contains a nonempty nowhere dense
subset without isolated points? Under Martin's Axiom, there are
countable non-discrete topological groups in which each nonempty
nowhere dense subset has an isolated point. \hf
\raisebox{-1ex}{\sl E.\,G.\,Zelenyuk}
 \emp

\bmp \textbf{17.51.} Is it true that every non-discrete topological
group containing no countable open subgroup of exponent 2 can be
partitioned into two (infinitely many) dense subsets? If the group
is Abelian or countable, the answer is yes.

 \hf \raisebox{-1ex}{\sl E.\,G.\,Zelenyuk, I.\,V.\,Protasov}
 \emp

\bmp \textbf{17.52.} (N.\,Eggert). Let $R$ be a commutative
associative finite-dimensional nilpotent algebra over a finite
field $F$ of characteristic $p$. Let $R^{(p)}$ be the subalgebra
of all elements of the form $r^p$ ($r\in R$). Is it true that
$\dim R\geq p\, \dim R^{(p)}$?

\makebox[15pt][r]{}An affirmative answer gives a
solution to 11.44 for the case where the factors $A$ and $B$ are
abelian. \hf \raisebox{-1ex}{\sl L.\,S.\,Kazarin}
 \emp

\bmp \textbf{17.53.} A finite group $G$ is called {\it simply
reducible\/} (SR-group) if every element of $G$ is conjugate to
its inverse, and the tensor product of any two irreducible
representations of $G$ 
decomposes
into a sum of irreducible representations of $G$ with coefficients
$0$ or~$1$.

\makebox[15pt][r]{}a) Is it true that the nilpotent length of a
soluble SR-group is at most $5$?

\makebox[15pt][r]{}b) Is it true that the derived length of a
soluble SR-group is bounded by some constant $c$?
 \hfill \raisebox{-1ex}{\sl L.\,S.\,Kazarin}

 \emp

\bmp \textbf{17.54.}
Does there exist a non-hereditary local formation ${\frak F}$ of
finite groups such that in every finite group the set of all
${\frak F}$-subnormal subgroups is a sublattice of the subgroup
lattice? %

\makebox[15pt][r]{}A negative answer would solve 9.75, since all the hereditary local
formations with the same property are already known.\hf
\raisebox{-1ex}{\sl S.\,F.\,Kamornikov}
 \emp

\bmp \textbf{\zv 17.56.}
Suppose that a subgroup $H$ of a finite group
$G$ is such that $HM=MH$ for every minimal non-nilpotent subgroup
$M$ of $G$. Must $H/H_G$ be nilpotent, where $H_G$ is the largest
normal subgroup of $G$ contained in $H$?
\hf \raisebox{-1ex}{\sl V.\,N.\,Knyagina, V.\,S.\,Monakhov}

\ul

\otv No, not always (A.\,Ballester-Bolinches, R.\,Esteban-Romero, S.\,F.\,Kamornikov, V.\,P\'erez-Calabuig, {\it J.\,Algebra}, {\bf 608} (2022), 382--387).
 \emp

\bmp \textbf{17.57.} Let $r(m)=\{r+km\mid k\in \Z\}$ for integers
$0\leq r<m$. For \(r_1(m_1)\cap r_2(m_2) =\nobreak\varnothing\) define the
\emph{class transposition} \(\tau_{r_1(m_1),r_2(m_2)}\) as the
involution which interchanges\break \(r_1+km_1\) and \(r_2+km_2\) for
each integer \(k\) and fixes everything else. The group \({\rm
CT}(\mathbb{Z})\) generated by all class transpositions is simple
({\it Math. Z.}, {\bf 264}, no.\,4 (2010), 927--938). Is \({\rm Out}({\rm CT}(\mathbb{Z})) = \langle \sigma \mapsto
\sigma^{n \mapsto -n-1} \rangle \cong {\rm C}_2\)?
\hfill \raisebox{-1ex}{\sl S.\,Kohl}
 \emp

\bmp \textbf{17.58.} Does \({\rm CT}(\mathbb{Z})\) have subgroups of
intermediate (word-) growth? \hf \raisebox{0ex}{\sl S.\,Kohl}
 \emp

\bmp \textbf{17.59.} \ A permutation of \(\mathbb{Z}\) is called
\emph{residue-class-wise affine} if there is a positive integer
\(m\) such that its restrictions to the residue classes (mod
\(m\)) are all affine. Is \({\rm CT}(\mathbb{Z})\) the group of
all residue-class-wise affine permutations of \(\mathbb{Z}\) which
fix the nonnegative integers setwise? \hf \raisebox{-1ex}{\sl
S.\,Kohl}
 \emp

\bmp \textbf{17.60.} Given a set \(\mathcal{P}\) of odd primes, let
\({\rm CT}_{\mathcal{P}}(\mathbb{Z})\) denote the subgroup of
\({\rm CT}(\mathbb{Z})\) which is generated by all class
transpositions which interchange residue classes whose moduli have
only prime factors in \(\mathcal{P} \cup \{2\}\). The groups
\({\rm CT}_{\mathcal{P}}(\mathbb{Z})\) are simple ({\it Math. Z.}, {\bf 264}, no.\,4 (2010), 927--938).
Are they pairwise nonisomorphic?\hf \raisebox{-1ex}{\sl
S.\,Kohl}
 \emp

\bmp \textbf{17.61.} The group \({\rm CT}_{\mathcal{P}}(\mathbb{Z})\)
is finitely generated if and only if \(\mathcal{P}\) is finite. If
\(\mathcal{P} =\nobreak \varnothing\), then it is isomorphic to the
finitely presented (first) Higman--Thompson group
(J.\,P.\,McDermott, see Remark 1.4 in S.\,Kohl, {\it J.~Group Theory}, {\bf 20}, no.\,5 (2017), 1025--1030).
Is it always finitely presented if \(\mathcal{P}\) is finite?
\hfill \raisebox{-1ex}{\sl S.\,Kohl}
 \emp

\bmp \textbf{\zv 17.62.}
Given a free group $F_n$ and a proper
characteristic subgroup $C$, is it ever possible to generate the
quotient $F_n/C$ by fewer than $n$ elements? \hf
\raisebox{-1ex}{\sl J.\,Conrad}

\ul

\otv Yes, it is possible; moreover, for all $n\geq 2$, the free group $F_n$ admits continuum many pairwise non-isomorphic infinite simple characteristic 2-generated quotients (R.\,Coulon, F.\,Fournier-Facio, \url{https://arxiv.org/abs/2312.11684}).
\emp

\bmp \textbf{17.63.} Prove that if $p$ is an odd prime and $s$ is a
positive integer, then there are only finitely many $p$-adic space
groups of finite coclass with point group of coclass~$s$.

\hf
\raisebox{-1ex}{\sl C.\,R.\,Leedham-Green}
 \emp

 \bmp \textbf{17.64.} Say that a group $G$ is an {\it
 $n$-approximation\/}
to the Nottingham group $J=N(p)$ (as defined in Archive, 12.24) if
$G$ is an infinite pro-$p$ group, and $G/\gamma_n(G)$ is
isomorphic to $J/\gamma_n(J)$. Does there exist a function $f(p)$
such that, if $G$ is an $f(p)$-approximation to the Nottingham
group, then $\gamma_i(G)/\gamma_{i+1}(G)$ is isomorphic to
$\gamma_i(J)/\gamma_{i+1}(J)$ for all $i$?

\makebox[15pt][r]{}Cf. 14.56. Note that an affirmative solution to
this problem trivially implies the now known fact that $J$ is
finitely presented as a pro-$p$ group (if $p>2$), see
14.55.

\hf\raisebox{-1ex}{\sl C.\,R.\,Leedham-Green} \emp

\bmp \textbf{17.65.} If Problem 17.64 has an affirmative answer,
does it follow that, for some function $g(p)$, the number of
isomorphism classes of $g(p)$-approximations to $J$ is

\makebox[15pt][r]{}a) countable?

\makebox[15pt][r]{}b) one?

Note that if, for some $g(p)$, there are only finitely many
isomorphism classes, then for some $h(p)$ there is only one.
\hf\raisebox{-1ex}{\sl C.\,R.\,Leedham-Green} \emp

\bmp \textbf{17.66.} (R.\,Guralnick). Does there exist a
positive integer $d$ such that $\dim H^1(G,V)\leq d$ for any
faithful absolutely irreducible module $V$ for any finite group
$G$? Cf. 16.55 in Archive. \hf \raisebox{-1ex}{\sl V.\,D.\,Mazurov}

 \emp

\bmp \textbf{17.68.} Let $C$ be a cyclic subgroup of a group $G$ and
$G=CFC$ where $F$ is a finite cyclic subgroup. Is it true that
$|G:C|$ is finite? Cf.\,12.62. \hf \raisebox{-1ex}{\sl
V.\,D.\,Mazurov} \emp

\bmp \textbf{17.69.} Let $G$ be a group of prime exponent acting
freely on a non-trivial abelian group. Is $G$ cyclic? \hf
\raisebox{-1ex}{\sl V.\,D.\,Mazurov} \emp

\bmp \textbf{\zv 17.70.}
Let $\alpha$ be an automorphism of prime order
$q$ of an infinite free Burnside group $G=B(q,p)$ of prime
exponent $p$ such that $\alpha$ cyclically permutes the free
generators of $G$. Is it true that $\alpha$ fixes some non-trivial
element of $G$?

 \makebox[15pt][r]{}{\it Editors' comment}:
this is proved for $q=2$ in (V.\,S.\,Atabekyan, H.\,T.\,Aslanyan, {\it Proc. Yerevan State Univ. Physic. Math. Sci.}, \textbf{53}, no.~3 (2019), 147--149).

 \hf \raisebox{-1ex} {\sl V.\,D.\,Mazurov}

 \ul

 \otv No, it does not  for all $q\geq 3$ and all sufficiently large $p$ (V.S. Atabekyan,  {\it Preprint}, 2026, \url{https://arxiv.org/abs/2609.33941}).
 \emp

\bmp \textbf{\zv 17.71.} 
b) Let $\alpha$ be a fixed-point-free automorphism
of prime order $p$ of a periodic group~$G$. Suppose that $G$ does not contain a
non-trivial $p$-element. Is $\alpha$ a splitting automorphism?
\hfill {\raisebox{-1ex}{\sl V.\,D.\,Mazurov}}

\ul

\otv No, not always (V.\,S.\,Atabekyan, {\it Preprint}, 2026, \url{https://arxiv.org/abs/2610.15360}).
\emp

\bmp \textbf{17.72.} b) Let $AB$ be a Frobenius group with kernel $A$
and complement $B$. Suppose that $AB$ acts on a finite group $G$
so that $GA$ is also a Frobenius group with kernel $G$ and
complement $A$.
Is the exponent of $G$ bounded in terms of
$|B|$ and the exponent of $C_G(B)$?
\hf \raisebox{-1ex}{\sl V.\,D.\,Mazurov}
\emp

 \bmp \textbf{\zv 17.75.}
Can the Monster $M$ act on a nontrivial finite $3$-group $V$ so
that all elements of $M$ of orders $41$, $47$, $59$, $71$ have no
nontrivial fixed points in~$V$?
\hf \raisebox{-1ex}{\sl V.\,D.\,Mazurov}

\ul

\otv No, it cannot; moreover, $M$ cannot act on a $3$-group in such a way that an element of order 41 has no nontrivial fixed points, as shown in \S\,3 of (M.\,Lee, T.\,Popiel, {\it  J.~Group Theory}, {\bf  26}  (2023), 193--205).
 \emp

\bmp \textbf{\zv 17.76.}
 Does there exist a finite group $G$, with $|G|
> 2$, such that there is  exactly one element in $G$ which is not
a commutator?
 \hf \raisebox{-1ex} {\sl D.\,MacHale}

 \ul

 \otv Yes, such group do exist, and there are infinitely many examples (S.\,V.\,Skresanov, {\it Siberian Electron. Math. Rep.}, {\bf 23}, no.\,1 (2026), 603--607); two other examples are also found in
 (O.\,Hatem, D.\,Siniora, {\it Int. J. Group Theory},
{\bf  15},  no.\,3 (2026), 161--167).
 \emp

\bmp \textbf{17.78.} Does there exist a finitely generated group
without free subsemigroups generating a proper variety {containing
${\frak A}_p{\frak A}$?}
\hfill { \raisebox{-1ex}{\sl O.\,Macedo\'nska}}
 \emp

\bmp \textbf{17.80.} Let
$[u,\,_nv]:=[u,\overbrace{v,\dots ,v}^{n}\,]$. Is the group
$M_n=\langle x, y\mid x=[x,\,_ny],\;\, y=[y,\,_nx]\,\rangle$
infinite for every $n>2$? (See also 11.18.)
\hf \raisebox{-1ex}{\sl O.\,Macedo\'nska}
 \emp

\bmp \textbf{17.81.} Given normal subgroups $R_1,\dots, R_n$ of a
group $G$, let $\displaystyle [[R_1,\dots, R_n]]:=\prod
\Big[\bigcap_{i\in I}R_i,\bigcap_{j\in J}R_j\Big] $, where the
product is over all $I\cup J=\{1,\dots,n\}$, $ I\cap J=\varnothing$.
Let $G$ be a free group, and let $R_i=\langle r_i\rangle^G$ be the
normal closures of elements $r_i\in G$. It
is known (B.\,Hartley, Yu.\,Kuzmin,
{\it J. Pure Appl. Algebra}, \textbf{74} (1991), 247--256) that the
quotient $(R_1\cap R_2)/{[R_1,\,R_2]}$ is a free abelian group. On
the other hand, the quotient $(R_1\cap \dots \cap
R_n)/[[R_1,\dots, R_n]] $ has, in general, non-trivial torsion for
$n\geq 4$. Is this quotient always torsion-free for $n=3$? This is
known to be true if $r_1,r_2,r_3$ are not proper powers in~$G$.
 \hf\raisebox{-1ex}{\sl R.\,Mikhailov}
 \emp

\bmp \textbf{17.83.} Does there exist a group such that every
central extension of it is residually nilpotent, but there exists
a central extension of a central extension of it which is not
residually nilpotent? \hf \raisebox{-1ex}{\sl R.\,Mikhailov}
 \emp

\bmp \textbf{17.84.}
An associative algebra $A$ is said to be {\it Calabi--Yau of
dimension $d$} (for short, CY$d$) if there is a natural
isomorphism of $A$-bimodules ${\rm
Ext}_{A\text{-bimod}}^d(A,\,A\otimes A)\cong \nobreak A$ and ${\rm
Ext}_{A\text{-bimod}}^n(A,\,A\otimes A)=0$ for $n\ne d$. By
Kontsevich's theorem, the complex group algebra ${\Bbb C}G$ of the
fundamental group $G$ of a 3-dimensional aspherical manifold is
CY$3$.
Is every group with CY$3$ complex group algebra residually
finite?
 \hf
\raisebox{-1ex}{\sl R.\,Mikhailov}

 \emp

\bmp \textbf{17.85.} Let $\mathcal V$ be a variety of groups such
that any free group in $\mathcal V$ has torsion-free integral
homology groups in all dimensions. Is it true that $\mathcal V$ is
abelian?

\makebox[15pt][r]{}It is known that free metabelian groups may
have 2-torsion in the second homology groups, and free nilpotent
groups of class 2 may have 3-torsion in the fourth homology
groups. \hf \raisebox{-1ex}{\sl R.\,Mikhailov}
 \emp

\bmp \textbf{17.86.} (Simplest questions related to the
 Whitehead asphericity conjecture). Let $\langle x_1,x_2,x_3\mid
r_1,r_2,r_3\rangle$ be a presentation of the trivial group.

\makebox[25pt][r]{a)} Prove that the group $\langle
x_1,x_2,x_3\mid r_1,r_2\rangle$ is torsion-free.
 \hf
\raisebox{0ex}{\sl R.\,Mikhailov}
 \emp

\bmp \textbf{17.87.} Construct a group of intermediate growth with
finitely generated Schur multiplier. \hfill \raisebox{-1ex}{\sl
R.\,Mikhailov}

 \emp

\bmp \textbf{17.88.} Compute $K_0({\Bbb F}_2G)$ and $K_1({\Bbb
F}_2G)$, where $G$ is the first Grigorchuk group, ${\Bbb F}_2G$
its group algebra over the field of two elements, and $K_0, K_1$
the zeroth and first $K$-functors. \hfill \raisebox{-1ex}{\sl
R.\,Mikhailov}

 \emp

\bmp \textbf{17.89.} By Bousfield's theorem, the free pronilpotent
completion of a non-cyclic free group has uncountable Schur
multiplier. Is it true that the free prosolvable completion (that
is, the inverse limit of quotients by the derived subgroups) of a
non-cyclic free group has uncountable Schur multiplier?
\hf \raisebox{-1ex}{\sl R.\,Mikhailov}

 \emp

\bmp \textbf{17.90.} (G.\,Baumslag). A group is {\it
para\-free\/} if it is residually nilpotent and has the same lower
central quotients as a free group. Is it true that $H_2(G)=0$ for
any finitely generated para\-free group $G$? \hf \raisebox{-1ex}{\sl
R.\,Mikhailov}\emp

\bmp \textbf{\ai 17.91.} 
Let $d(X)$ denote the derived length of a group
$X$.

\makebox[15pt][r]{}a) Does there exist an absolute constant $k$
such that $d(G)-d(M)\leq k$ for every finite soluble group $G$
and any maximal subgroup $M$?

\makebox[15pt][r]{}b) Find the minimum $k$ with this property.
\hf \raisebox{-0ex}{\sl V.\,S.\,Monakhov}

\ul

\oai No, there is no such constant (unverified AI solution: L.\,Sabatini, {\it Preprint}, 2026, \url{https://kourovkanotebookorg.wordpress.com/wp-content/uploads/2026/11/17-91-1.pdf}).
\emp 

\bmp \textbf{17.93.} (Well-known problem). Let $G$ be a compact
topological group which has elements of arbitrarily high orders.
Must $G$ contain an element of infinite order?

\hf \raisebox{-1ex}{\sl J.\,Mycielski}
\emp

\bmp \textbf{\zv 17.94.}
 (Well-known problem). Can the free product $Z*G$
of the infinite cyclic group $Z$ and a nontrivial group $G$ be the
normal closure of a single element?

 \makebox[15pt][r]{}{\it Editors' comments}:
the negation is also known as Kervaire's conjecture. It was proved for torsion-free groups (A.\,Klyachko, {\it Commun. Algebra}, {\bf 21}, no.\,7 (1993), 2555--2575) and for residually finite groups (M.\,Gerstenhaber, O.\,S.\,Rothaus, {\it Proc. Nat. Acad.
Sci. U.S.A.}, {\bf 48} (1962) 1531--1533).
\hfill \raisebox{-1ex}{\sl J.\,Mycielski}

\ul

\otv No, it cannot (AI solution confirmed by N.\,V.\,Bogachev: \ OpenAI, {\it Preprint}, 2026, \url{https://kourovkanotebookorg.wordpress.com/wp-content/uploads/2026/10/17-94.pdf}).
\emp 

\bmp \textbf{17.95.} Let $G$ be a permutation group on the finite
 set $\Omega$. A partition $\rho$ of $\Omega$
 is said to be 
 $G$-{\it regular} if there exists
 a subset $S$ of $\Omega$ such that $S^g$ is a 
 transversal 
 of $\rho$ for all $g \in G$. The group $G$ is said to be {\it synchronizing} if $|\Omega| > 2$
 and there are no non-trivial proper $G$-regular partitions on $\Omega$.

\makebox[15pt][r]{}a) Are the following affine-type primitive
groups synchronizing:\vs

$2^{p}.{\rm PSL}(2,2p+1)$ where
 both $p$ and $2p + 1$ are prime,
 $p \equiv 3\;(\mathrm{mod}\;4)$ and $p > 23$?\vs

 $2^{101}.He$?\vs

\makebox[15pt][r]{}b) For which finite simple groups $S$
 are the groups $S \times S$ acting
 on $S$ by $(g,h): x \mapsto g^{-1}xh$
 non-synchronizing?
 \hf
\raisebox{-1ex}{\sl P.\,M.\,Neumann}

 \emp

 \bmp \textbf{17.96.} Does there exist a variety of groups
 which contains only countably many subvarieties
 but in which there is an infinite properly
 descending chain of subvarieties?

 \hf
\raisebox{-1ex}{\sl P.\,M.\,Neumann}

 \emp

 \bmp \textbf{17.97.} Is every variety of groups of exponent $4$
 finitely based?
 \hf
\raisebox{-1ex}{\sl P.\,M.\,Neumann}

 \emp

 \bmp \textbf{17.98.} A variety of groups is said to be {\it small\/}
 if it contains only countably many non-isomorphic finitely
 generated groups.

\makebox[15pt][r]{}a) Is it true that if $G$ is a finitely
generated group and the variety
 ${\rm Var}(G)$ it generates is
 small then $G$ satisfies the
 maximal condition on normal subgroups?

\makebox[15pt][r]{}b) Is it true that a variety is
 small if and only if all its finitely
 generated groups have the Hopf property?
 \hf
\raisebox{-1ex}{\sl P.\,M.\,Neumann}
 \emp

\bmp \textbf{\zv 17.99.}
Consider the group $B = \langle a,b \mid
(bab^{-1}) a (bab^{-1})^{-1} = a^2\rangle $ introduced by Baumslag
in 1969. The same relation is satisfied by the functions $ f(x) =
2x$ and $g(x) = 2^x$ under the operation of composition in the
group of germs of monotonically increasing to $\infty $ continuous
functions on $(0, \infty )$, where two functions are identified if
they coincide for all sufficiently large arguments. Is the
representation $a\rightarrow f$, $b\rightarrow g$ of the group $B$
faithful? \hf \raisebox{-1ex}{\sl A.\,Yu.\,Olshanskii}

\ul

\otv No, it is not (C.-F.\,Nyberg-Brodda, {\it Preprint}, 2026, \url{https://arxiv.org/abs/2606.27408}).
\emp

\bmp \textbf{17.100.} {\it Conjecture}\/: A finite group is not
simple if it has an irreducible complex character of odd degree
vanishing on a class of odd length.

\makebox[16pt][r]{}If true, this implies the solvability of groups
of odd order, so a proof independent of CFSG is of special
interest. \hf \raisebox{-1ex}{\sl V.\,Pannone}
 \emp

\bmp \textbf{\ai 17.101.} 
 According to P.\,Hall, a group $G$ is said to be
{\it homogeneous\/} if every isomorphism of its finitely generated
subgroups is induced by an automorphism of $G$. It is known
(Higman--Neumann--Neumann) that every group is embeddable into a
homogeneous one. Is the same true for the category of
representations of groups? A~representation $(V,G)$ is finitely
generated if $G$ is a finitely generated group, and $V$ a
finitely generated module over the group algebra of~$G$.
 \hf
\raisebox{-1ex}{\sl B.\,I.\,Plotkin}

\ul

\oai Yes, it is true (unverified AI solutions:  Achyuth Jayadevan,  {\it Preprint of 11 September 2026},
\url{https://zenodo.org/records/22722081}; \ M.\,Kinyon, J.\,Urban, {\it Preprint of 17 September 2026},    \url{https://kourovkanotebookorg.wordpress.com/wp-content/uploads/2026/09/kourovka-mathematics.pdf}).

 \emp 

\bmp \textbf{\zv 17.102.}
We say that two subsets $A, B$ of an infinite
group $G$ are {\it separated\/} if there exists an infinite subset
$X$ of $G$ such that $1\in X$, $X=X^{-1}$, and $XAX\cap
B=\varnothing$. Is it true that any two disjoint subsets $A, B$ of
an infinite group $G$ satisfying $|A|<|G|$, $|B|<|G|$ are
separated? This is so if $A,B$ are finite, or $G$ is Abelian. \hf
\raisebox{-1ex}{\sl I.\,V.\,Protasov}

\ul

\otv No, not always (S.\,Wan, {\it Preprint}, 2026, \url{https://arxiv.org/abs/2608.00504}).
 \emp

\bmp \textbf{17.103.} Does there exist a continuum of sets $\pi$ of primes for which
every finite group possessing a Hall $\pi$-subgroup is a
$D_\pi$-group?

 \makebox[15pt][r]{}{\it Comment of 2025}:
It has been proved that, for any integer $x$, if $\pi$ is the set of all prime numbers $p>x$, then any finite group containing a Hall $\pi$-subgroup is a $D_\pi$-group (K.\,A.\,Ilenko, N.\,V.\,Maslova, {\it Siberian Math. J.}, \textbf{62}, no.\,1 (2021), 44--51).

\hf \raisebox{-1ex}{\sl D.\,O.\,Revin}
 \emp

\bmp \textbf{17.104.} Let $\Gamma$ be a finite non-oriented graph on
the set of vertices $\{x_1,\dots,x_n\}$ and let
$S_\Gamma=\langle x_1,\dots,x_n\mid x_ix_j=x_jx_i \Leftrightarrow
(x_i,x_j) \in \Gamma;\;\, {\frak A}^2\rangle$ be a presentation of
a partially commutative metabelian group $S_\Gamma$ in the variety
of all metabelian groups. Is the universal theory of the group
$S_\Gamma$ decidable?
 \hf \raisebox{-1ex}{\sl
V.\,N.\,Remeslennikov, E.\,I.\,Timoshenko} \emp

\bmp \textbf{17.105.} An equation over a pro-$p$-group $G$ is an
expression $v(x)=1$, where $v(x)$ is an element of the free
pro-$p$-product of $G$ and a free pro-$p$-group with basis $\{
x_1, \ldots , x_n \}$; solutions are sought in the affine space
$G^n$. Is it true that a free pro-$p$-group is equationally Noetherian,
 that is, for any $n$ every system of
equations in $x_1, \ldots , x_n$ over this group is equivalent to
some finite subsystem of it?
 \hf
\raisebox{-1ex}{\sl N.\,S.\,Romanovski\u{\i}}

 \emp

\bmp \textbf{17.107.} Does $G=\mathrm{SL}_2(\mathbb{C})$ contain a
2-generated free subgroup that is conjugate in~$G$ to a proper
subgroup of itself?

\makebox[15pt][r]{}A 6-generated free subgroup with this property
has been found (D.\,Calegari, N.\,M.\,Dunfield,
{\it Proc. Amer. Math. Soc.}, \textbf{134}, no.\,11 (2006), 3131--3136).
 \hf
\raisebox{-1ex}{\sl M.\,Sapir}

 \emp

\bmp \textbf{17.109.} A non-trivial group word $w$ is \emph{uniformly
elliptic} in a class $\mathcal{C}$ if there is a function
$f:\mathbb{N}\rightarrow\mathbb{N}$ such that the width of $w$ in
every $d$-generator $\mathcal{C}$-group $G$ is bounded by $f(d)$
(i.\,e. every element of the verbal subgroup $w(G)$ is equal to a
product of $f(d)$ values of $w$ or their inverses). If
$\mathcal{C}$ is a class of finite groups, this is equivalent to
saying that in every finitely generated pro-$\mathcal{C}$ group
$G$ the verbal subgroup $w(G)$ is closed. A.\,Jaikin-Zapirain
(\emph{Revista Mat. Iberoamericana}, \textbf{24} (2008), 617--630)
proved that $w$ is uniformly elliptic in finite $p$-groups if and
only if $w\notin F^{\prime\prime}(F^{\prime})^{p}$, where $F$ is
the free group on the variables of $w$. Is it true that $w$ is
uniformly elliptic in $\mathcal{C}$ if and only if $w\notin
F^{\prime\prime}(F^{\prime})^{p}$ for every prime $p$ in the case
where $\mathcal{C}$ is the class of all

\makebox[15pt][r]{}a) finite soluble groups?

\makebox[15pt][r]{}b) finite groups?

See also Ch.\,4 of (D.\,Segal, \emph{Words: notes on verbal width
in groups}, LMS Lect. Note Series, \textbf{361}, Cambridge Univ.
Press, 2009. \hfill \raisebox{-1ex}{\sl D.\,Segal}
 \emp

\bmp \textbf{17.110.} a) Is it true that for each word $w$, there is
a function $h:\mathbb{N}\times \mathbb{N}\rightarrow\mathbb{N}$
such that the width of $w$ in every finite $p$-group of Pr\"{u}fer
rank $r$ is bounded by $h(p,r)$?

\makebox[15pt][r]{}b) If so, can $h(p,r)$ be made independent of
$p$?

\makebox[15pt][r]{}An affirmative answer to a) would imply
A.\,Jaikin-Zapirain's result ({\it ibid.}) that every word has
finite width in each $p$-adic analytic pro-$p$ group, since a
pro-$p$ group is $p$-adic analytic if and only if the Pr\"{u}fer
ranks of all its finite quotients are uniformly bounded.
 \hf \raisebox{-1ex}{\sl D.\,Segal}
 \emp

\bmp \textbf{17.112.} A subgroup $A$ of a group $G$ is said to be
{\it $G$-permutable in $G$} if for every subgroup $B$ of $G$
there exists an element $x\in G $ such that $AB^x=B^xA$. A
subgroup $A$ is said to be {\it hereditarily $G$-permutable in
$G$} if $A$ is $E$-permutable in every subgroup $E$ of $G$
containing $A$. Which finite non-abelian simple groups $G$
possess

\makebox[15pt][r]{}a) a non-trivial $G$-permutable subgroup?

\makebox[15pt][r]{}b) a non-trivial hereditarily
$G$-per\-mu\-table subgroup?

\makebox[15pt][r]{}{\it Comments of 2025}:
Among sporadic groups, only  $J_1$ has a proper $G$-permutable subgroup (A.\,A.\,Galt, V.\,N.\,Tyutyanov, {\it Siberian Math.~J.}, {\bf 63}, no.\, 4 (2022), 691--698). Sporadic, alternating, and exceptional groups of Lie type have no proper hereditarily $G$-per\-mu\-table subgroups (A.\,A.\,Galt, V.\,N.\,Tyutyanov, {\it Siberian Math.~J.}, {\bf 63}, no.\, 4 (2022), 691--698; \ A.\,F.\,Vasilyev, V.\,N.\,Tyutyanov, {\it
Izv. Gomel. Gos. Univ.  F. Skoriny}, {\bf 2012}, no.\,5(74) (2012), 148--150 (Russian)); \ A.\,A.\,Galt, V.\,N.\,Tyutyanov, {\it Siberian Math.~J.}, {\bf 64}, no.\,5 (2023), 1110--1116).
 \hf \raisebox{-1ex}{\sl A.\,N.\,Skiba, V.\,N.\,Tyutyanov}
 \emp

\bmp \textbf{\zv 17.113.}
Do there exist for $p>3$ \ 2-generator finite
$p$-groups
 with deficiency zero (see 8.12) of arbitrarily high
nilpotency class? \hf\raisebox{-1ex}{\sl J.\,Wiegold}

\ul

\otv Yes, they do exist (AI solution confirmed by \'{O}.\,Fern\'{a}ndez-Ayala and M.\,Trombetti,  {\it Preprint}, 2026,
  \url{https://kourovkanotebookorg.wordpress.com/wp-content/uploads/2026/10/17.113.pdf}).
\emp

 \bmp \textbf{17.114.}
Does every generalized free product (with amalgamation) of two
non-trivial groups have maximal subgroups?
 \hf \raisebox{-1ex}{\sl J.\,Wiegold}
\emp

\bmp \textbf{17.115.} Can a locally free non-abelian group have
non-trivial Frattini subgroup?

\hf \raisebox{-1ex}{\sl J.\,Wiegold}
\emp

\bmp \textbf{17.117.} (Well-known problem). If groups $A$ and $B$
have decidable elementary theories $Th(A)$ and $Th(B)$, must
$Th(A*B)$ be decidable? \hfill \raisebox{-1ex}{\sl
O.\,Kharlampovich}
 \emp

\bmp \textbf{17.118.} Suppose that a finite $p$-group $G$ has a
subgroup of exponent $p$ and of index~$p$. Must $G$ also have a
characteristic subgroup of exponent $p$ and of index bounded in
terms of $p$? \hf \raisebox{-1ex}{\sl E.\,I.\,Khukhro}
\emp

 \bmp \textbf{17.119.} Suppose that a finite soluble group $G$ admits a soluble group of automorphisms
 $A$ of coprime order such that $C_G(A)$ has rank $r$. Let $|A|$ be the product of $l$ not necessarily distinct primes.
Is there a linear function $f$ such that $G/F_{f(l)}(G)$ has
$(|A|,r)$-bounded rank, where $F_{f(l)}(G)$ is the $f(l)$-th
Fitting subgroup?

\makebox[15pt][r]{}An exponential function $f$ with this property
was found in (E.\,Khukhro, V.\,Mazurov, {\it Groups St.\,Andrews
2005}, vol.\,II, Cambridge Univ. Press, 2007, 564--585). It is
also known that $|G/F_{2l+1}(G)|$ is bounded in terms of
$|C_G(A)|$ and $|A|$ (Hartley--Isaacs). \hfill \raisebox{-1ex}{\sl
E.\,I.\,Khukhro}

 \emp

\bmp \textbf{17.120.}
 Is a residually finite group all of whose subgroups of infinite
index are finite necessarily cyclic-by-finite? \hfill
\raisebox{-1ex}{\sl N.\,S.\,Chernikov}
 \emp

\bmp \textbf{17.121.} Let $ G $ be a group whose set of elements is
the real numbers and which is ``nicely definable'' (see below).
Does
 $ G $ being $ { \aleph_\omega} $-free
imply it being free if ``nicely definable'' means

\makebox[15pt][r]{}a) being $ F_ \sigma $?

\makebox[15pt][r]{}b) being Borel?

\makebox[15pt][r]{}c) being analytic?

\makebox[15pt][r]{}d) being projective $ {\Bbb L} [ {\Bbb R} ]
$?

See \url{https://arxiv.org/pdf/math/0212250.pdf} for justification of the
restrictions.

\hf \raisebox{-0ex}{\sl S.\,Shelah}
 \emp

\bmp \textbf{17.122.} The same questions as 17.121 for an abelian
group $G$.
 \hf \raisebox{-0ex}{\sl
S.\,Shelah}
 \emp

\bmp \textbf{17.123.} Do there exist finite groups $G_1,\,G_2$ such
that $\pi _e(G_1)=\pi _e(G_2)$, \ $h(\pi _e(G_1))<\infty$, and
each non-abelian composition factors of each of the groups
$G_1,\,G_2$ is not isomorphic to a section of the other?

\makebox[15pt][r]{}For definitions see Archive, 13.63. \hfill
\raisebox{0ex}{\sl W.\,J.\,Shi}\emp

\bmp \textbf{\ai 17.124.} 
Is the set of finitely presented metabelian
groups recursively enumerable?

\hf \raisebox{-1ex}{\sl
V.\,Shpilrain}

\ul

\oai Yes, it is (unverified AI solution:  Achyuth Jayadevan,  {\it Preprint}, 2026, \url{https://arxiv.org/abs/2609.10281}).
 \emp 

\bmp \textbf{17.125.} Does every finite group $G$ contain a pair of
conjugate elements $a,b$ such that $\pi(G)=\pi(\langle
a,b\rangle)$? This is true for soluble groups.

\makebox[15pt][r]{}{\it Comment of 2013:}
it was proved in (A.\,Lucchini, M.\,Morigi, P.\,Shumyatsky, {\it Forum Math.}, {\bf 24} (2012), 875--887) that every finite group $G$ contains a 2-generator subgroup $H$ such that $\pi(G)=\pi(H)$.\hf
\raisebox{-1ex}{\sl P.\,Shumyatsky}
 \emp

\bmp \textbf{17.126.} Suppose that $G$ is a residually finite group
satisfying the identity $[x,y]^n=1$. Must $[G,G]$ be locally
finite?

\makebox[15pt][r]{}An equivalent question: let $G$ be a finite
soluble group satisfying the identity $[x,y]^n=1$; is the Fitting
height of $G$ bounded in terms of $n$? Cf. Archive, 13.34

\hf
\raisebox{-1ex}{\sl P.\,Shumyatsky}
 \emp

\bmp \textbf{17.127.} Suppose that a finite soluble group $G$ of
derived length $d$ admits an elementary abelian $p$-group of
automorphisms $A$ of order $p^n$ such that $C_G(A)=1$. Must $G$
have a normal series of $n$-bounded length with nilpotent factors
of $(p,n,d)$-bounded nilpotency class?

\makebox[15pt][r]{}This is true for $p=2$. An affirmative answer would follow from
an affirmative answer to~11.125.\hf
\raisebox{-1ex}{\sl P.\,Shumyatsky}
 \emp

\bmp \textbf{17.128.} Let $T$ be a finite $p$-group admitting an
elementary abelian group of automorphisms $A$ of order $p^{2}$
such that in the semidirect product $P=TA$ every element of
$P\setminus T$ has order $p$. Does it follow that $T$ is of
exponent~$p$? \hf \raisebox{-1ex}{\sl E.\,Jabara}

\emp

\newpage

\pagestyle{myheadings} \markboth{18th Issue (2014)}{18th Issue
(2014)} \thispagestyle{headings} ~ \vspace{2ex}

\centerline {\Large \textbf{Problems from the 18th Issue (2014)}}
\phantomsection\label{18izd}
\vspace{4ex}


\bmp \textbf{\zv 18.1.}
Given a group $G$ of finite order $n$, does there necessarily
exist a bijection $f$ from $G$ onto a cyclic group of order $n$
such that for each element $x\in G$, the order of $x$ divides
the order of $f(x)$?
An affirmative answer in the case where $G$ is solvable was given by F.\,Ladisch.
\hf \raisebox{-1ex}{\sl I.\,M.\,Isaacs}
\ul

\otv Yes, it does (M.\,Amiri, {\it J.~Pure Applied Algebra}, \textbf{228}, no.~7 (2024), 107632).
\emp

\bmp \textbf{18.2.} Let the group $G=AB$ be the product of two Chernikov subgroups $A$ and $B$ each of which
has an abelian subgroup of index at most 2. Is $G$ soluble?
\hf \raisebox{-1ex}{\sl B.\,Amberg}
 \emp

\bmp \textbf{18.3.} Let the group $G=AB$ be the product of two central-by-finite subgroups $A$ and $B$. By a theorem of N.\,Chernikov, $G$ is soluble-by-finite. Is $G$ metabelian-by-\linebreak\text{finite?} \text{} \hf \raisebox{-1ex}{\sl B.\,Amberg}
 \emp

\bmp \textbf{18.4.} (A.\,Rhemtulla, S.\,Sidki). a) Is a group of the form $G=ABA$ with cyclic subgroups $A$ and $B$ always soluble? (This is known to be true if $G$ is finite.)

\makebox[15pt][r]{}b) The same question when in addition $A$ and $B$ are conjugate in $G$.
\hf \raisebox{-1ex}{\sl B.\,Amberg}
 \emp

\bmp \textbf{18.5.}
Let the (soluble) group $G=AB$ with finite
torsion-free rank $r_0(G)$ be the product of two subgroups $A$ and $B$.
Is the equation $r_0(G) = r_0(A) + r_0(B) - r_0(A \cap B)$ valid in this case?
It is known that the inequality $\leq$ always holds.
\hf \raisebox{-1ex}{\sl B.\,Amberg}
 \emp

\bmp
 \textbf{18.6.} Is every finitely generated one-relator group residually amenable?

\hfill \raisebox{-1ex}{\sl G.\,N.\,Arzhantseva}\emp

\bmp \textbf{18.7.} Let $n\geq 665$ be an odd integer. Is it true that the group of outer automorphisms
${\rm Out}(B(m,n))$ of the free Burnside group $B(m,n)$ for $m>2$ is complete, that is, has trivial center and all its automorphisms are inner?
\hf \raisebox{-1ex}{\sl V.\,S.\,Atabekyan}
 \emp

\bmp
 \textbf{18.8.} Let $\mathcal X$ be a class of finite simple groups such that $\pi(\mathcal X) = char(\mathcal X)$. A~formation of finite groups $\mathfrak F$ is said to be $\mathcal X$-saturated if a finite group $G$ belongs to $\mathfrak F$ whenever the factor group $G/\Phi(O_{\mathcal X}(G))$ is in $\mathfrak F$, where $O_{\mathcal X}(G)$ is the largest normal subgroup of $G$ whose composition factors are in $\mathcal X$. Is every $\mathcal X$-saturated formation $\mathcal X$-local in the sense of F\"orster? (See the definition in
(P.\,F\"orster, 
{\it Publ. Sec. Mat. Univ. Aut\`onoma Barcelona}, {\bf 29}, no.\,2--3 (1985), 39--76).)
\hfill \raisebox{-1ex}{\sl A.\,Ballester-Bolinches}\emp

\bmp
 \textbf{18.10.} A formation $\mathfrak F$ of finite groups satisfies the {\it Wielandt property for residuals\/} if whenever $U$ and $V$ are subnormal subgroups of $\langle U, V\rangle $ in a finite group $G$, then the $\mathfrak F$-residual $\langle U, V\rangle ^{\mathfrak F}$ of $\langle U, V\rangle $ coincides with $\langle U^{\mathfrak F}, V^{\mathfrak F}\rangle $. Does every Fitting formation $\mathfrak F$ satisfy the Wielandt property for residuals? \hfill \raisebox{-1ex}{\sl A.\,Ballester-Bolinches}\emp

\bmp
 \textbf{\zv 18.11.}
 Let $G$ and $H$ be subgroups of the automorphism group $Aut (F_n)$ of a free group $F_n$ of rank
 $n \geq 2$. Is it true that the free product $G * H$ embeds in the automorphism group
$Aut (F_m)$ for some $m$?
 \hfill \raisebox{-1ex}{\sl V.\,G.\,Bardakov}

 \ul

\otv Yes, it is true (AI solution confirmed by V.\,G.\,Bardakov:  \ Achyuth Jayadevan,  {\it Preprint}, 2026,
\url{https://zenodo.org/records/22722085}).
 \emp 

\bmp
 \textbf{18.12.}
Let $G$ be a finitely generated group of intermediate growth. Is it true that there is a positive integer $m$ such that every element of the derived subgroup $G'$ is a product of at most $m$ commutators?

\makebox[15pt][r]{}This assertion is valid for groups of polynomial growth, since they are almost nilpotent by Gromov's theorem. On the other hand, there are groups of exponential growth (for example, free groups) for which this assertion is not true. \hfill \raisebox{-1ex}{\sl V.\,G.\,Bardakov}
 \emp

\bmp
 \textbf{18.13.} (D.\,B.\,A.\,Epstein).
Is it true that the group $ H = (\mathbb{Z}_3 \times \mathbb{Z}) * (\mathbb{Z}_2 \times \mathbb{Z}) =
 \langle x, y, z, t \mid x^3 = z^2 = [x, y] = [z, t] = 1 \rangle$ cannot be defined by three relators in the generators $x, y, z, t$?

 \makebox[15pt][r]{}It is known that the relation module of the group $H$ has rank~3 \ (K.\,W.\,Gruenberg, P.\,A.\,Linnell, {\it J. Group Theory}, {\bf 11}, no.\,5 (2008), 587--608). An affirmative answer would give a solution of the relation gap problem.
 \hfill \raisebox{-1ex}{\sl V.\,G.\,Bardakov, M.\,V.\,Neshchadim}
 \emp

\bmp
 \textbf{\zv 18.14.}
For an automorphism $\varphi \in Aut(G)$ of a group $G$, let $[e]_{\varphi} =\{g^{-1}g^{\varphi}\mid g\in G\}$.

 \makebox[15pt][r]{}\emph{Conjecture}: If $[e]_{\varphi}$ is a subgroup for every $\varphi \in Aut(G)$, then
the group $G$ is nilpotent. If in addition $G$ is finitely generated, then $G$ is abelian.

 \hfill \raisebox{-1ex}{\sl V.\,G.\,Bardakov, M.\,V.\,Neshchadim, T.\,R.\,Nasybullov}

 \ul

 \otv The conjecture is disproved (C.\,Nicotera, {\it Arch. Math.}, {\bf 123}, no.\,3 (2024), 225--232).
 \emp

\bmp
 \textbf{\zv 18.16.}
Is it true that any definable endomorphism of any ordered abelian group
is of the form $x\mapsto rx$, for some rational number $r$?\hfill \raisebox{-1ex}{\sl O.\,V.\,Belegradek}

\ul

\otv Yes, it is true (H.\,Ji, F.\,Xu, D.\,Qiu,  {\it Preprint}, 2026, \url{https://arxiv.org/abs/2608.17861}).
\emp

\bmp
 \textbf{\zv 18.18.}
 Mal'cev proved that the set of sentences that hold in all finite groups is not
computably enumerable, although its complement is. Is it true that both
the set of sentences that hold in almost all finite groups and its complement
are not computably enumerable?\hfill \raisebox{-1ex}{\sl O.\,V.\,Belegradek}

\ul

\otv Yes, it is true (AI solution confirmed by H.\,D.\,Macpherson and B.\,I.\,Zilber:  \ M.\,Kinyon, J.\,Urban, {\it Preprint}, 2026,
\url{https://kourovkanotebookorg.wordpress.com/wp-content/uploads/2026/09/kourovka-mathematics.pdf}).
\emp 

\bmp
 \textbf{18.19.} Is any torsion-free, relatively free group of infinite rank not $\aleph_1$-homo\-ge\-neous?
This is true for group varieties in which all free groups are residually finite
(O.\,Belegradek, \emph{Arch. Math. Logic}, {\bf 51}
(2012), 781--787).\hfill \raisebox{-1ex}{\sl O.\,V.\,Belegradek}
\emp

 \bmp
 \textbf{18.20.} Characters $\varphi$ and $\psi$ of a finite group $G$
are said to be {\it semiproportional\/} if they are not proportional and there is a
normal subset $M$ of $G$ such that
$\varphi |_M$ is proportional to $\psi |_M$ and
$\varphi |_{G\setminus M}$ is proportional to $\psi |_{G\setminus M}$.

 \makebox[15pt][r]{}{\it Conjecture\/}: If $\varphi$ and $\psi$ are semiproportional
irreducible characters of a finite group, then $\varphi(1) = \psi(1)$.
\hfill \raisebox{-1ex}{\sl V.\,A.\,Belonogov}
\emp

\bmp
 \textbf{18.21.} (B.\,H.\,Neumann and H.\,Neumann).
Fix an integer $d\geq 2$. If $\mathfrak{V}$ is a variety of groups such that
all $d$-generated groups in $\mathfrak{V}$ are finite,
must $\mathfrak{V}$ be locally finite?

\hfill \raisebox{-1ex}{\sl G.\,M.\,Bergman}
\emp

\bmp
 \textbf{18.22.} (H.\,Neumann).
Is the Kostrikin variety of all locally finite groups
of given prime exponent $p$
determined by finitely many identities?
\hfill \raisebox{-1ex}{\sl G.\,M.\,Bergman}
\emp

\bmp
 \textbf{18.24.} For a group $G$, a function $\phi : G\to {\Bbb R}$ is a \emph{quasimorphism} if there is a least non-negative number $D(\phi)$ (called the \emph{defect}) such that $|\phi(gh) -\phi(g) - \phi(h)| \leq D(\phi)$ for all $g, h\in G$. A quasimorphism is \emph{homogeneous} if in addition $\phi(g^n) = n\phi(g)$ for all $g\in G$.
 Let $\phi$ be a homogeneous quasimorphism of a free group $F$. For any quasimorphism $\psi$ of $F$ (not required to be homogeneous) with $|\phi - \psi| < {\rm const}$, we must have $D(\psi) \ge D(\phi)/2$.
 Is it true that there is some $\psi$ with $D(\psi) = D(\phi)/2$?

\hfill \raisebox{-1ex}{\sl M.\,Burger, D.\,Calegari}
 \emp

\bmp
 \textbf{\ai 18.25.} 
Let ${\rm form}(G)$ be the formation generated by a finite group $G$. Suppose that $G$ has a unique composition series
$ 1\vartriangleleft G_1\vartriangleleft G_2\vartriangleleft G$ and the consecutive factors of this series are $Z_p$, $X$, $Z_q$,
where $p$ and $q$ are primes and $X$ is a non-abelian simple group. Is it true that ${\rm form}(G)$ has infinitely many subformations if and only if $p=q$?

\hfill \raisebox{-1ex}{\sl V.\,P.\,Burichenko}

\ul

\oai No, it is not true (unverified AI solution: B.\,G.\,Di Siena, {\it Preprint}, 2026, \url{https://www.proofforum.org/abs/PF:2026.000038}).
\emp

\bmp
 \textbf{18.26.} Suppose that a finite group $G$ has a normal series $ 1< G_1<G_2<G $
such that the groups $G_1$ and $G_2/G_1$
are elementary abelian $p$-groups, $G/G_2\cong A_5\times A_5$
(where $A_5$ is the alternating group of degree 5), $G_1$ and
$G_2/G_1$ are minimal normal subgroups of $G$ and $G/G_1$, respectively.
Is it true that ${\rm form}(G)$ has finitely many subformations?
\hfill \raisebox{-1ex}{\sl V.\,P.\,Burichenko}
\emp

\bmp
 \textbf{18.27.} Does there exist an algorithm that determines whether there are finitely many subformations in ${\rm form}(G)$ for a given finite group $G$? \hfill \raisebox{-1ex}{\sl V.\,P.\,Burichenko}
\emp

\bmp
 \textbf{18.28.} Is it true that any subformation of every one-generator formation ${\rm form}(G)$ is also one-generator?
\hfill \raisebox{-1ex}{\sl V.\,P.\,Burichenko}
\emp

\bmp
 \textbf{18.29.} Let $\mathfrak {F}$ be a Fitting class of finite soluble groups which contains every soluble group $G=AB$, where $A$ and $B$ are abnormal $\mathfrak {F}$-subgroups of $G$. Is $\mathfrak {F}$ a formation?

\hfill \raisebox{-1ex}{\sl A.\,F.\,Vasil'ev}
\emp

\bmp
 \textbf{18.34.} Let $G$ be a topological group and $S={G_0\leq \dots \leq G_n=G}$ a subnormal series of closed subgroups. The \emph{infinite-length} of $S$ is the cardinality of the set of infinite factors $G_i/G_{i-1}$. The \emph{virtual length} of $G$ is the supremum of the infinite lengths taken over all such series of $G$. It is known that if $G$ is a pronilpotent group with finite virtual length then every closed subnormal subgroup is topologically finitely generated (N.\,Gavioli, V.\,Monti, C.\,M.\,Scoppola, \emph{J.~Austral. Math. Soc.}, {\bf 95}, no.\,3 (2013), 343--355). Let $G$ be a pro-$p$ group (or more generally a pronilpotent group). If every closed subnormal subgroup of $G$ is tolopologically finitely generated, is it true that $G$ has finite virtual length?
\hfill \raisebox{-1ex}{\sl N.\,Gavioli, V.\,Monti, C.\,M.\,Scoppola}\emp

\bmp
 \textbf{18.35.} A family ${\mathscr F}$ of group homomorphisms $A\to B$ is {\it separating\/} if for every nontrivial $a\in A$ there is $f\in {\mathscr F}$ such that $f(a)\ne 1$, and {\it discriminating\/} if for any finitely many nontrivial elements $a_1,\dots ,a_n\in A$ there is $f\in {\mathscr F}$ such that $f(a_i)\ne 1$ for all $i=1,\dots , n$.
Let $G$ be a relatively free group of rank 2 in the variety of metabelian groups. Let $T$ be a metabelian group
in which the centralizer of every nontrivial element is abelian. If $T$ admits a separating family of surjective
homomorphisms $T\to G$ must it also admit a discriminating family of surjective homomorphisms $T\to G$?

 \makebox[15pt][r]{}A positive answer would give a metabelian analogue of a classical theorem of B.\,Baumslag.
\hfill \raisebox{-1ex}{\sl A.\,Gaglione, D.\,Spellman}
 \emp

\bmp
 \textbf{\zv 18.36.}
There are groups of cardinality at most $2^{\aleph_{0}}$, even nilpotent of
class $2$, that cannot
be embedded in $S:=Sym(\mathbb{N})$ (V.\,A.\,Churkin,
\emph{Algebra and model theory} (Novosibirsk State Tech. Univ.), {\bf  5} (2005), 39--43 (Russian)).
One condition which might characterize subgroups of $S$ is the
following. Consider the metric $d$ on $S$ where for all $x\neq y$
we define $d(x,y):=2^{-k}$ when $k$ is the least element of the set $\mathbb{N}$ such that $k^{x}\neq
k^{y}$. Then $(S,d)$ is a separable topological group, and so
each subgroup of $S$ (not necessarily a closed subgroup of $(S,d)$) is a
separable topological group under the induced metric. Is it true that every
group $G$ for which there is a metric $d^{\prime}$ such that $(G,d^{\prime})$
is a separable topological group is (abstractly) embeddable in $S$?

\makebox[15pt][r]{}Note that any such group $G$ can be embedded as a \emph{section} in $S$.
 \hfill \raisebox{-1ex}{\sl  J.\,D.\,Dixon}

 \ul

 \otv No, there exists an uncountable connected Polish group whose every abstract homomorphism into $S$ is trivial
 (C.\,Rosendal, S.\,Solecki, {\it Israel J. Math.},
 {\bf  162} (2007), 349--371); see also (A.\,Ka\"{i}chouh, F.\, Le Ma\^{i}tre, \emph{Bull. London Math. Soc.}, \textbf{47} (2015), 996--1009) (S.\,Corson, {\it Letter  of 5 June 2024}).
\emp

\bmp \textbf{18.37.} (Well-known problem). Is every locally graded group of finite rank almost locally soluble? \hf \raisebox{-1ex}{\sl M.\,Dixon}
 \emp

\bmp \textbf{\zv 18.38.}
 Let $E_\pi$ denote the class of finite groups that contain a Hall $\pi$-subgroup. Does the inclusion $E_{\pi_1} \cap E_{\pi_2} \subseteq E_{\pi_1 \cap \pi_2}
$ hold for arbitrary sets of primes $\pi_1$ and $\pi_2$?

 \hf \raisebox{-1ex}{\sl A.\,V.\,Zavarnitsine}

 \ul

 \otv Yes, it does (A.\,Buturlakin, N.\,Yang, {\it Preprint}, 2025,
 \url{https://arxiv.org/abs/2501.05865}).
 \emp

\bmp
 \textbf{18.39.} \textit{Conjecture}: Let $G$ be a hyperbolic group. Then every 2-dimensional rational homology class is virtually represented by a sum of closed surface subgroups, that is, for any $\alpha \in H_2(G;\Q)$ there are finitely many closed oriented surfaces $S_i$ and injective homomorphisms $\rho_i:\pi_1(S_i) \to G$ such that $\sum_i [S_i] = n \alpha$, where $[S_i]$ denotes the image of the fundamental class of $S_i$ in $H_2(G)$.\hfill \raisebox{-1ex}{\sl D.\,Calegari}
 \emp

\bmp
 \textbf{18.40.} The commutator length $cl(g)$ of $g \in [G,G]$ is the least number of commutators in $G$ whose product is $g$, and the stable commutator length is $scl(g):=
\lim_{n \to \infty} cl(g^n)/n$.

 \makebox[15pt][r]{} \textit{Conjecture}: Let $G$ be a hyperbolic group. Then the stable commutator length takes on rational values on $[G,G]$.
\hfill \raisebox{-1ex}{\sl D.\,Calegari}
 \emp

\bmp
 \textbf{18.41.} Let $F$ be a free group of rank 2.

 \makebox[15pt][r]{}a) It is known that $scl(g)=0$ only for $g=1$, and $scl(g)\ge 1/2$ for all $1\ne g\in [F,F]$. Is $1/2$ an isolated value?

 \makebox[15pt][r]{}b) Are there any intervals $J$ in ${\Bbb R}$ such that the set of values of $scl$ on $[F,F]$ is dense in $J$?

 \makebox[15pt][r]{}c) Is there some $T\in {\Bbb R}$ such that every rational number ${}\ge T$ is a value of $scl(g)$ for $g\in [F,F]$?
\hfill \raisebox{-1ex}{\sl D.\,Calegari}
 \emp

\bmp
 \textbf{18.42.}
 For an orientation-preserving homeomorphism $f:S^1\to S^1$ of a unit circle $S^1$, let $\tilde f$ be its lifting to a homeomorphism of ${\Bbb R}$; then the \emph{rotation number} of $f$ is defined to be the limit $\lim _{n\to \infty} (\tilde f^n(x)-x)/n$ (which is independent of the point $x\in S^1$). Let $F$ be a free group of rank 2 with generators $a,b$. For $ w\in F$ and $r,s\in {\Bbb R}$, let $R(w,r,s)$ denote the maximum value of the (real-valued) rotation number of $w$, under all representations from $F$ to the universal central extension of $Homeo^+(S^1)$ for which the rotation number of $a$ is $r$, and the rotation number of $b$ is $s$.

\makebox[25pt][r]{a)} If $r,s$ are rational, must $R(w,r,s)$ be rational?

\makebox[25pt][r]{b)} Let $R(w,r-,s-)$ denote the supremum of the rotation number of $w$ (as above) under all representations for which $a$ and $b$ are
conjugate to rotations through $r$ and $s$, respectively (in the universal central extension of $Homeo^+(S^1)$). Is $R(w,r-,s-)$ always rational if $r$ and $s$ are rational?

\makebox[25pt][r]{c)} \textit{Weak Slippery Conjecture}: Let $w$ be a word containing only positive powers of $a$ and $b$. A pair of values $(r,s)$ is \textit{slippery} if there is a strict inequality $R(w,r',s') < R(w,r-,s-)$ for all $r'<r$, $s'<s$. Is it true that then $R(w,r-,s-) = h_a(w)r + h_b(w)s$, where $h_a(w)$ counts the number of $a$'s in $w$, and $h_b(w)$ counts the number of $b$'s in $w$?

\makebox[25pt][r]{d)} \textit{Slippery Conjecture}: More precisely, is it always (without assuming $(r,s)$ being slippery) true that if $w$ has the form $w=a^{\alpha_1}b^{\beta_1}\cdots a^{\alpha_m}b^{\beta_m}$ (all positive) and $R(w,r,s) = p/q$ where $p/q$ is reduced, then $|p/q - h_a(w)r - h_b(w)s| \le m/q$?

\hfill \raisebox{-1ex}{\sl D.\,Calegari, A.\,Walker}
 \emp

\bmp
 \textbf{18.43.}
a) Do there exist elements $u,v$ of a 2-generator free group $F(a,b)$ such that $u,v$ are not conjugate in $F(a,b)$ but for any matrices $A,B\in GL(3, \mathbb{C})$ we have $trace(u(A,B))= trace(v(A,B))$?

 \makebox[15pt][r]{}b) The same question if we replace $GL(3,\mathbb{C})$ by $SL(3,\mathbb{C})$. \hfill \raisebox{0ex}{\sl I.\,Kapovich}
\emp

\bmp
 \textbf{\zv 18.44.} Let $G$ be a finite non-abelian group and $V$ a finite faithful irreducible $G$-module. Suppose that $M=|G/G'|$ is the largest orbit size of $G$ on $V$, and among orbits of $G$ on $V$ there are exactly two orbits of size $M$.
Does this imply that $G$ is dihedral of order 8, and $|V|=9$?\hfill \raisebox{-1ex}{\sl T.\,M.\,Keller}

\ul

\otv Yes, it does (N.\,A.\,Jones,\, T.\,M.\,Keller, {\it  Ann. Mat. Pura Appl.}, {\bf 199} (2020), 2323--2340).
\emp

\bmp
 \textbf{18.45.} Let $G$ be a finite $p$-group of maximal class such that all its irreducible characters are induced from linear characters of normal subgroups. Let $S$ be the set of the derived
subgroups of the members of the central series of $G$. Is there a logarithmic bound for the derived length of $G$
in terms of $|S|$?\hfill \raisebox{-1ex}{\sl T.\,M.\,Keller}\emp

\bmp
 \textbf{18.46.}
Every finite group $G$ can be embedded in a group $H$ in such a way
that every element of $G$ is a square of an element of $H$.
The overgroup $H$ can be chosen such that
$|H| \le 2|G|^2$.
Is this estimate sharp?

\makebox[15pt][r]{}It is known that the best possible estimate cannot be better than $|H| \le |G|^2$ (D.\,V.\,Baranov, Ant.\,A.\,Klyachko, \emph{Siber. Math.~J.}, {\bf 53}, no.\,2 (2012), 201--206).

\hfill \raisebox{-1ex}{\sl Ant.\,A.\,Klyachko}\emp

\bmp \textbf{18.47.}
Is it algorithmically decidable whether a group generated by three given
class transpositions (for the definition, see 17.57)

\makebox[25pt][r]{a)} has only finite orbits on \(\mathbb{Z}\)?

\makebox[25pt][r]{b)} acts transitively on the set of nonnegative integers in its support?

A difficult case is the group
\(\langle \tau_{1(2),4(6)}, \tau_{1(3),2(6)}, \tau_{2(3),4(6)} \rangle\),
which acts transitively on \(\mathbb{N} \setminus 0(6)\) if and only if
Collatz' \(3n+1\) conjecture is true. \hf \raisebox{-1ex}{\sl
S.\,Kohl}\emp

\bmp \textbf{\zv 18.48.}
Is it true that there are only finitely many integers which occur as orders of products
of two class transpositions? (For the definition, see 17.57.)
\hf \raisebox{-1ex}{\sl
S.\,Kohl}

\ul

\otv Yes, it is true (A.\,Iskra,  {\it Preprint}, 2026, \url{https://arxiv.org/abs/2610.01793}).
\emp

\bmp \textbf{\zv 18.50.}
Let \(n \in \mathbb{N}\). Is it true that for every \(k \in \{1, \dots, n!\}\) there is
some group \(G\) and pairwise distinct elements \(g_1, \dots, g_n \in G\) such that the set
\(\{g_{\sigma(1)} \cdots g_{\sigma(n)} \ | \ \sigma \in {\rm S}_n\}\)
of all products of the \(g_i\) obtained by permuting the factors has cardinality \(k\)?
\hf \raisebox{-1ex}{\sl
S.\,Kohl}

\ul

\otv Yes, it is true (AI solution confirmed by S.\,Kohl: \  P.\,Monticone, {\it Preprint}, 2026,  \url{https://kourovkanotebookorg.wordpress.com/wp-content/uploads/2026/05/18_50.pdf}); see also (W.\,van\,Doorn, E.\,Judin, P.\,Monticone, D.\,Morrison,  {\it Preprint}, 2026, \url{https://arxiv.org/abs/2607.17477});
see also S.\,Kohl's post at
\url{https://mathoverflow.net/questions/123890/numbers-of-distinct-products-obtained-by-permuting-the-factors/511354#511354}).
\emp

\bmp \textbf{18.51.}
Given a prime \(p\) and \(n \in \mathbb{N}\), let \(f_p(n)\) be the smallest
number such that there is a group of order \(p^{f_p(n)}\) into which every group of
order \(p^n\) embeds. Is it true that \(f_p(n)\) grows faster than polynomially
but slower than exponentially when \(n\) tends to infinity?
\hf \raisebox{-1ex}{\sl
S.\,Kohl}
\emp

\bmp \textbf{18.53.}
Is there a non-linear simple locally finite group in which every
centralizer of a non-trivial element is almost soluble, that is, has a
soluble subgroup of finite index?

\hf \raisebox{-1ex}{\sl M.\,Kuzucuo\u{g}lu}
 \emp

\bmp \textbf{18.54.} Can a group be equal to the union of conjugates of a proper finite nonabelian simple subgroup?
\hf \raisebox{-1ex}{\sl G.\,Cutolo} \emp

\bmp \textbf{18.55.} a)
Can a locally finite $p$-group $G$ of finite exponent be the union of conjugates of an abelian proper subgroup?

\makebox[25pt][r]{b)} Can this happen when $G$ is of exponent $p$?
\hf \raisebox{0ex}{\sl G.\,Cutolo} \emp

\bmp \textbf{18.56.} Let $G$ be a finite 2-group, of order greater than 2, such that $|H/H_G| \le 2$ for all $H\le G$, where $H_G$ denotes the largest normal subgroup of $G$ contained in $H$. Must $G$ have an abelian subgroup of index 4?
\hf \raisebox{-1ex}{\sl G.\,Cutolo} \emp

 \bmp \textbf{18.58.}
Let $G$ be a group generated by finite number $n$ of involutions
in which $(uv)^4=1$ for all involutions $u,v\in G$. Is it true
that $G$ is finite? is a 2-group?
This is true for $n\leq 3$. 

\makebox[15pt][r]{}{\it Editors' comment of 2021}: the previously published claim
of an affirmative solution of this problem proved to be erroneous. \hfill \raisebox{-1ex}{\sl D.\,V.\,Lytkina}
 \emp

  \bmp \textbf{\zv 18.59.}
Does there exist a periodic group $G$ such that $G$ contains an
involution, all involutions in $G$ are conjugate, and the
centralizer of every involution $i$ is isomorphic to $\langle i\rangle \times L_2(P)$,
where $P$ is some infinite locally finite field of characteristic~$2$?

\hf \raisebox{-1ex}{\sl D.\,V.\,Lytkina}

\ul

\otv No, there are no such groups (K.\,Muliarchyk,  {\it Preprint}, 2026, \url{https://kourovkanotebookorg.wordpress.com/wp-content/uploads/2026/09/kourovka_18_59-18_60_solution.pdf}).
 \emp

\bmp \textbf{\zv 18.60.}
Let $V$ be an infinite countable elementary abelian additive 2-group.
Does ${\rm Aut}\, V$ contain a subgroup $G$ such that

\makebox[25pt][r]{(a)} $G$ is transitive on the set of non-zero elements of $V$, and

\makebox[25pt][r]{(b)} if $H$ is the stabilizer in $G$ of a non-zero element $v\in V$,
then $V=\langle v\rangle \oplus V_v$, where $V_v$ is
$H$-invariant, $H$ is isomorphic to the multiplicative group $P^*$
of a locally finite field $P$ of characteristic $2$, and the
action $H$ on $V_v$ is similar to the action of $P^*$ on $P$ by
multiplication?

\makebox[15pt][r]{}{\it Conjecture:} such a group $G$ does not exist. If so, then the
group $G$ in the previous problem does not exist too.
\hf \raisebox{-1ex}{\sl D.\,V.\,Lytkina}

\ul

\otv Yes, there are no such groups (K.\,Muliarchyk,  {\it Preprint}, 2026, \url{https://kourovkanotebookorg.wordpress.com/wp-content/uploads/2026/09/kourovka_18_59-18_60_solution.pdf}).
 \emp

\bmp \textbf{18.61.} Is a 2-group nilpotent if all its finite subgroups are nilpotent
of class at most~3? This is true if 3 is replaced by 2.
\hf \raisebox{-1ex}{\sl D.\,V.\,Lytkina.}
 \emp

\bmp \textbf{18.62.} 
The {\it spectrum\/} of a finite group is the
set of orders of its elements. Let $\omega$ be a finite set of positive integers. A group $G$ is
said to be
{\it $\omega$-critical\/}
if the spectrum of $G$ coincides with $\omega$, but the spectrum
of every proper
section of $G$
is not equal to $\omega$.

\makebox[25pt][r]{\zva a)} Does there exist a number $n$ such that
for every finite simple group $G$ the number of
$\omega(G)$-critical groups is less than $n$?

\makebox[25pt][r]{b)} For every finite simple
 group $G$, find all $\omega(G)$-cri\-ti\-cal groups.
\hf \raisebox{-1ex}{\sl V.\,D.\,Mazurov}

\ul

\otv a) No, there is no such number  (AI solution confirmed by V.\,D.\,Mazurov: \  I.\,S.\,Nnuga Mongo, {\it Preprint}, 2026, \url{https://arxiv.org/abs/2611.18465}).
 \emp

 \bmp \textbf{\zv 18.63.} 
Let $G$ be a periodic group generated by two fixed-point-free
automorphisms of order 5 of an abelian group. Is $G$ finite?
\hf \raisebox{-1ex}{\sl V.\,D.\,Mazurov}

\ul

\otv No, not necessarily (AI solution confirmed by V.\,D.\,Mazurov: \   I.\,S.\,Nnuga Mongo, {\it Preprint}, 2026, \url{https://arxiv.org/abs/2611.18465}).
 \emp

\bmp \textbf{18.64.} (K.\,Harada). {\it Conjecture\/}: Let $G$ be a finite group, $p$ a prime, and
$B$ a $p$-block of $G$. If $J$ is a non-empty subset of $Irr(B)$
such that $\sum\limits_{\chi \in J}\chi (1)\chi (g)=0$
for
every \smash[t]{$p$-singular element $g\in G$, then $J=Irr(B)$.}
\hf \smash[t]{\raisebox{-1ex}{\sl V.\,D.\,Mazurov}}
 \emp

  \bmp \textbf{\zv 18.65.} 
 (R.\,Guralnick, G.\,Malle). {\it Conjecture\/}: Let $p$ be a prime different from~5, and $C$ a class of conjugate $p$-elements in a finite group $G$. If $[c,d]$ is a $p$-element for any $c,d\in C$, then $C\subseteq O_p(G)$.
 \hfill \raisebox{-1ex}{\sl V.\,D.\,Mazurov}

\ul

\otv Conjecture is refuted (AI solution confirmed by G.\,Malle: \   I.\,S.\,Nnuga Mongo, {\it Preprint}, 2026, \url{https://arxiv.org/abs/2611.18465}).
\emp

\bmp \textbf{18.66.} Suppose that a finite group $G$ admits a Frobenius group of automorphisms $FH$ with kernel $F$ and complement $H$ such that $GF$ is also a Frobenius group with kernel $G$ and complement $F$.
Is the derived length of $G$ bounded in terms of $|H|$ and the derived length of $C_G(H)$?
\hf \raisebox{-1ex}{\sl N.\,Yu.\,Makarenko, E.\,I.\,Khukhro, P.\,Shumyatsky}
 \emp

\bmp \textbf{18.67.} Suppose that a finite group $G$ admits a Frobenius group of automorphisms $FH$ with kernel $F$ and complement $H$ such that $C_G(F)=1$. Is the exponent of $G$ bounded in terms of $|F|$ and the exponent of $C_G(H)$?

\makebox[15pt][r]{}This was proved when $F$ is cyclic (E.\,I.\,Khukhro, N.\,Y.\,Makarenko, P.\,Shumyatsky, \textit{Forum Math.},
{\bf 26} (2014),
73--112) and for $FH\cong {\Bbb A}_4$ when $(|G|,3)=1$ (P.\,Shumyatsky, {\it J.~Algebra}, {\bf 331} (2011), 482--489).
\hf \raisebox{-1ex}{\sl N.\,Yu.\,Makarenko, E.\,I.\,Khukhro, P.\,Shumyatsky}
 \emp

\bmp
 \textbf{\ai 18.68.} 
 What are the nonabelian composition factors of a finite nonsoluble group all of whose maximal subgroups have complements? \hfill

\makebox[15pt][r]{}Note that finite simple groups with all maximal subgroups having complements are, up to isomorphism, $L_2(7)$, $L_2(11)$, and $L_5(2)$ (V.\,M.\,Levchuk, A.\,G.\,Likharev, \emph{Siberian Math. J.}, {\bf 47}, no.\,4 (2006), 659--668).
The same groups exhaust finite simple groups with Hall maximal subgroups and composition factors of groups with Hall maximal subgroups (N.\,V.\,Maslova, \emph{Siberian Math. J.}, {\bf 53}, no.\,5 (2012), 853--861).
It is also proved that in a finite group with Hall maximal subgroups all maximal subgroups have complements
(N.\,V.\,Maslova, D.\,O.\,Revin, \emph{Siberian Adv. Math.}, {\bf 23}, no.\,3 (2013), 196--209).\hfill
 \raisebox{-1ex}{\sl N.\,V.\,Maslova, D.\,O.\,Revin}

\ul

\oai These factors are described (unverified AI solution: \ R.\,Sater,  {\it Preprint}, 2026, \url{https://kourovkanotebookorg.wordpress.com/wp-content/uploads/2026/09/kourovka-18-68-v1.1.1.pdf
}).
 \emp 

\bmp
 \textbf{18.69.} Does there exist a relatively free group $G $
containing a free subsemigroup and having $[G,G]$ finitely generated?

 \makebox[15pt][r]{}Note that $G$ cannot be locally graded ({\it Publ. Math. Debrecen}, {\bf 81}, no.\,3-4 (2012), 415--420.) \hfill \raisebox{-1ex}{\sl O.\,Macedo\'nska}
\emp

\bmp \textbf{18.70.}
Is every finitely generated Coxeter group
conjugacy separable?

\makebox[15pt][r]{}Some special cases were considered in
(P.-E.\,Caprace, A.\,Minasyan, {\it Illinois J. Math.}, \textbf{57}, no.\,2 (2013), 499--523).
\hf \raisebox{-1ex}{\sl A.\,Minasyan}
 \emp

\bmp
 \textbf{\zv 18.71.}
 (N.\,Aronszajn). Suppose that $W(x, y) = 1$ in an open subset of $G\times G$, where $G$ is a connected topological group. Must $W(x, y) = 1$ for all $(x, y)\in G\times G$?

\makebox[15pt][r]{}For locally compact groups the answer is yes.
\hfill \raisebox{-1ex}{\sl J.\,Mycielski}

\ul

\otv No, not necessarily: for any odd integer $n>10^{10}$ there is a connected topological group such that the identity $x^n=1$ holds in some neighborhood of unity, but not in the entire group (E.\,Reznichenko, I.\,Zyabrev, \textit{Preprint}, 2024, \url{https://arxiv.org/abs/2406.05203}).
\emp

\bmp
 \textbf{18.72.}
 Is it true that an existentially closed subgroup of a nonabelian free group of finite rank is a nonabelian free factor of this group?
\hfill \raisebox{-1ex}{\sl A.\,G.\,Myasnikov, V.\,A.\,Roman'kov}
\emp

\bmp \textbf{18.73.}
a) Does every finitely generated solvable group of derived length $l\ge 2$ embed
into a $2$-generated solvable group of length $l+1$?

\makebox[15pt][r]{}{\it Comment of 2021}: It is proved that any countable solvable group of derived length $l$ with torsion-free abelianization embeds in a 2-generated solvable group of derived length $l+1$ (V.\,A.\,Roman'kov, {\it Proc. Amer. Math. Soc.}, \textbf{149} (2021), 4133--4143).

\hf \raisebox{-1ex}{\sl A.\,Yu.\,Olshanskii}
 \emp

\bmp \textbf{18.74.} Let $G$ be a finitely generated elementary amenable group which
is not virtually nilpotent. Is there a finitely generated metabelian non-virtually-nilpotent section
in $G$?
\hf \raisebox{-1ex}{\sl A.\,Yu.\,Olshanskii}
 \emp

\bmp \textbf{\zv 18.76.}
Let $A$ be a division ring, $G$ a subgroup of the
multiplicative group of $A$, and $E$ an extension of the additive
group of $A$ by $G$ such that $G$ acts by multiplication in~$A$.
Is it true that $E$ splits? This is true if $G$ is finite.
\hf \raisebox{-1ex}{\sl E.\,A.\,Palyutin.}

\ul

\otv No, not always (AI solutions confirmed by B.A.F.Wehrfritz: \  Achyuth Jayadevan,  {\it Preprint of 14 September 2026}, \url{https://zenodo.org/records/22751423}; \ M.\,Kinyon, J.\,Urban, {\it Preprint of 17 September 2026},    \url{https://kourovkanotebookorg.wordpress.com/wp-content/uploads/2026/09/kourovka-mathematics.pdf}).
).
 \emp 

\bmp
 \textbf{18.77.}
Let $G$ be a finite $p$-group and let $p^e$ be the largest
degree of an irreducible complex representation of $G$. If $p>e$, is it
necessarily true that $\bigcap \ker \Theta= 1$,
where the intersection runs over all irreducible complex representations $\Theta$ of
$G$ of degree $p^e$?

\hfill \raisebox{-1ex}{\sl D.\,S.\,Passman}
 \emp

\bmp
 \textbf{18.78.} Let $K^tG$ be a twisted group algebra of the finite group $G$ over the field $K$.
If $K^tG$ is a simple $K$-algebra, is $G$ necessarily solvable? This is known to be
true if $K^tG$ is central simple.
\hfill \raisebox{-1ex}{\sl D.\,S.\,Passman}
 \emp

\bmp
 \textbf{18.79.} Let $K[G]$ be the group algebra of the finitely generated group $G$
over the field $K$. Is the Jacobson radical
${\cal J}K[G]$ equal to the join of all nilpotent ideals of the ring? This is
known to be true if $G$ is solvable or linear.
\hfill \raisebox{-1ex}{\sl D.\,S.\,Passman}
 \emp

\bmp
 \textbf{18.80.} (I.\,Kaplansky). For $G\neq 1$, show that
the augmentation ideal of the group algebra $K[G]$ is equal to the Jacobson radical of the ring
if and only if ${\rm char}\, K=p>0$ and
$G$ is a locally finite $p$-group. \hfill \raisebox{-1ex}{\sl D.\,S.\,Passman}
 \emp

\bmp
 \textbf{\ai 18.81.} 
 Let $G$ be a finitely generated $p$-group that is residually finite.
Are all maximal subgroups of $G$ necessarily normal?
\hfill \raisebox{-1ex}{\sl D.\,S.\,Passman}

\ul

\oai No, not always (unverified AI solution: B.\,G.\,Di Siena, {\it Preprint}, 2026, \url{https://www.proofforum.org/abs/PF:2026.000030}).
\emp

\bmp
\textbf{18.83.}
A generating system $X$ of a group $G$ is \emph{fast} if there is an integer $n$
such that every element of $G$ can be expressed as a product of at most
$n$ elements of $X$ or their inverses. If not, we say that it is \textit{slow}.
For instance, in $(\Z ,+)$, the squares are fast, but the powers of 2 are slow.

\makebox[15pt][r]{}Do there exist countable infinite groups without an infinite slow generating set? Uncountable ones do exist.
\hfill \raisebox{-1ex}{\sl B.\,Poizat}

\emp

\bmp \textbf{18.84.} Let $\pi$ be a set of primes. We say that a finite group is a $BS_\pi$-group if every conjugacy class in this group any two elements of which generate a $\pi$-subgroup itself generates a
$\pi$-subgroup. Is every normal subgroup of a $BS_\pi$-group a $BS_\pi$-group?

\makebox[15pt][r]{}In the case
$2\not\in\pi$, an affirmative answer follows from
(D.\,O.\,Revin, {\it Siberian Math. J.}, \textbf{52}, no.\,2 (2011), 340--347).
\hf \raisebox{-1ex}{\sl D.\,O.\,Revin}
 \emp

\bmp
 \textbf{18.85.}
 A subset of a group is said to be \emph{rational} if it can be obtained from finite subsets by finitely many rational operations, that is, taking union, product, and the submonoid generated by a set.

 \makebox[15pt][r]{}{\it Conjecture}\/: every finitely generated solvable group in which all rational subsets form a Boolean algebra is virtually abelian.

 \makebox[15pt][r]{}This is known to be true if the group is metabelian, or polycyclic, or of finite rank (G.\,A.\,Bazhenova).
\hfill \raisebox{-1ex}{\sl V.\,A.\,Roman'kov} 
\emp

\bmp
 \textbf{18.87.}
 A system of equations with coefficients in a group $G$ is said to be \emph{independent} if the matrix composed of the sums of exponents of the unknowns has rank equal to the number of equations.

 \makebox[25pt][r]{\zv a)} The {\it Kervaire--Laudenbach Conjecture (KLC)\/}: every independent system of equations with coefficients in an arbitrary group $G$ has a solution in some overgroup~$\bar{G}$.
This is true for every locally residually finite group $G$ (M.\,Gerstenhaber and O.\,S.\,Rothaus).

 \makebox[25pt][r]{b)} KLC --- nilpotent version: every independent system of equations with coefficients in an arbitrary nilpotent group $G$ has a solution in some nilpotent overgroup~$\bar{G}$.

 \makebox[25pt][r]{c)} KLC --- solvable version: every independent system of equations with coefficients in an arbitrary solvable group $G$ has a solution in some solvable overgroup~$\bar{G}$.

\hfill \raisebox{-1ex}{\sl V.\,A.\,Roman'kov}

\ul

\otv a) The conjecture is proved (AI solution confirmed by N.\,V.\,Bogachev:

  OpenAI, {\it Preprint}, 2026, \url{https://kourovkanotebookorg.wordpress.com/wp-content/uploads/2026/10/18-87.pdf}).
\emp 

\bmp \textbf{18.88.} 
Can a finitely generated infinite group of finite exponent be the quotient of a residually finite group by a locally finite normal subgroup?

\makebox[15pt][r]{}If not, then there exists a hyperbolic group that is not residually finite.

\makebox[15pt][r]{}{\it Comment of 2026}: there exist hyperbolic groups that are not residually finite (unverified AI solution: OpenAI, {\it Preprint}, 2026, \url{https://kourovkanotebookorg.wordpress.com/wp-content/uploads/2026/10/partial-12-64_18-88.pdf}).

\hf \raisebox{-1ex}{\sl M.\,Sapir}
\emp 

\bmp \textbf{\zv 18.89.}
Consider the set of balanced presentations $\langle x_1,\dots ,x_n\mid r_1,\dots, r_n\rangle $ with fixed $n$ generators $x_1,\dots, x_n$.  By definition, the group ${AC}_n$ of Andrews--Curtis moves on this set of balanced presentations  is generated by the Nielsen transformations together with conjugations of relators. Is  ${AC}_n$  finitely presented?
  \hf \raisebox{-1ex}{\sl J.\,Swan, A.\,Lisitsa}

  \ul

  \otv The group ${AC}_2$ is not finitely presented (S.\,Krsti\'c,  J.\,McCool, \emph{J.~London Math. Soc. (2)}, {\bf 56}, no.\,2 (1997),  264--274; also V.\,A.\,Roman'kov, {\it Preprint}, 2023, \url{https://arxiv.org/pdf/2305.11838.pdf}), while the groups ${AC}_n$, $n\geq 3$, are finitely presented, which follows from  (G.\,Kiralis, S.\,Krsti\v{c},  J.\,McCool, {\it Proc. London Math. Soc. (3)}, {\bf 73}, no.\,3 (1996), 481--720) (M.\,Ershov, {\it Letter of 16 October 2025}, \url{https://kourovkanotebookorg.wordpress.com/wp-content/uploads/2025/11/kourovka1889.pdf}).
  \emp

\bmp
 \textbf{18.90.} Prove that the group
$Y(m,n) =\langle a_{1},a_{2},\dots,a_{m} \mid
a_{k}^{n} =1 \text{ }(1\leq k\leq m)\text{, } ( a_{k}^{i}a_{l}^{i}) ^{2} =e\text{ }(1\leq k<l\leq m\text{, }
1\leq i\leq \frac{n}{2})\rangle $
is finite for every pair $( m,n) $.

\makebox[15pt][r]{}Some known cases are (where $C_n$ denotes a cyclic group of order $n$): $Y(m,2)=C_{2}^{m}$; $Y(m;3)={A}_{m+2}$ (presentation by Carmichael); $Y(m;4)$
has order $2^{\frac{m(m+3)}{2}}$, nilpotency class $3$ and exponent $4$;
\ $Y(2,n)$ is the natural extension of the augmentation ideal $\omega ( C_{n})$ of $GF(2)[C_{n}]$ by $C_{n}$ (presentation by Coxeter);
$Y(3,n)\cong SL(2,\omega ( C_{n}) )$.

\makebox[15pt][r]{}These groups have connections with classical groups in characteristic $2$. When $n$ is odd, the group $Y(m,n)$ has the following presentation, where $\tau _{ij}$ are transpositions:
$y(m,n) =\langle a,\text{ }{\Bbb S} _{m}\text{ }\mid a^{n}=e,
[ \tau _{12}^{a^{i}},\tau _{12}] =e\;\, ( 1\leq i\leq
\frac{n}{2}) ,
\tau _{12}^{1+a+\dots+a^{n-1}} =e,\text{ }
\tau _{i,i+1}a\tau _{i,i+1} =a^{-1}\text{ }( 2\leq i\leq m-1) \rangle$, and it is known that, for example,
$y(3,5) \cong SL(2,16)\cong \Omega ^{-}(4,4)$, $y(4,5)\cong Sp(2,16)\cong
\Omega (5,4)$, $y(5,5) \cong SU(4,16)\cong \Omega ^{-}(6,4)$, $y(6,5)\cong 4^{6}\Omega
^{-}(6,4)$, $ y(3,7) \cong SL(2,8)^{2}\cong \Omega ^{+}(4,8)$, $y(4,7)\cong Sp(4,8)\cong
\Omega (5,8)$. Extensive computations by Felsch, Neub\"user and O'Brien confirm this trend.
Different types reveal Bott periodicity and a connection with Clifford
algebras.\hfill \raisebox{-1ex}{\sl S.\,Sidki}\emp

\bmp
 \textbf{\zv 18.92.}
A non-empty set $\theta$ of formations is called a \emph{complete lattice} of
formations if the intersection of any set of formations in $\theta$
belongs to $\theta$ and $\theta$ has the largest element (with respect to inclusion). If $L$ is a complete lattice, then an element $a\in L$ is
said to be compact if $ a\leq \bigvee X$ for any $X\subseteq L$ implies that $a\leq\bigvee X_1$ for some finite $X_1\subset X$. A complete lattice is called \emph{algebraic} if every element is the join of a (possibly infinite) set of compact elements.

\makebox[15pt][r]{}a) Is there a non-algebraic complete lattice of formations of finite groups?

\makebox[15pt][r]{}b) Is there a non-modular complete lattice of formations of finite groups?

\hfill \raisebox{-1ex}{\sl A.\,N.\,Skiba}

\ul

\otv a), b) Yes, there are, in both cases (AI solution confirmed by A.\,N.\,Skiba: \    M.\,Kinyon, J.\,Urban, {\it Preprint}, 2026,   \url{https://kourovkanotebookorg.wordpress.com/wp-content/uploads/2026/09/kourovka-mathematics.pdf}).
\emp 

\bmp
 \textbf{18.93.}
Let $\cal M$ be a one-generated saturated formation, that is, the
intersection of all saturated formations containing some fixed finite
group. Let $\cal F$ be a subformation of $\cal M$ such that $\cal F
\ne {\cal F}{\cal F}$.

\makebox[15pt][r]{}a) Is it true that then $\cal F$ can be written in the form
 ${\cal F}={\cal F}_{1} \cdots {\cal F}_{t}$, where ${\cal F}_{i}$ is a
non-decomposable formation for every $i=1, \ldots , t$?

\makebox[15pt][r]{}b) Suppose that ${\cal F}={\cal F}_{1} \cdots {\cal F}_{t}$, where
${\cal F}_{i}$ is a non-decomposable formation for every $i=1, \ldots , t$. Is it
true that then all factors ${\cal F}_{i}$ are uniquely determined?

\hfill \raisebox{-1ex}{\sl A.\,N.\,Skiba}
\emp

\bmp
 \textbf{18.96.} Suppose that a periodic group $G$ contains involutions and the centralizer of each involution is locally finite. Is it true that $G$ has a nontrivial locally finite normal subgroup?
 \hfill \raisebox{-1ex}{\sl N.\,M.\,Suchkov}
 \emp

\bmp
 \textbf{18.97.} Let $G$ be a periodic Zassenhaus group, that is, a two-transitive permutation group with trivial stabilizer of every three points. Suppose that the stabilizer of a point is a Frobenius group with locally finite kernel $U$ containing an involution. Is it true that $U$ is a $2$-group? This was proved for finite groups by Feit.
 \hfill \raisebox{-1ex}{\sl N.\,M.\,Suchkov}
 \emp

\bmp
 \textbf{18.98.} The work of many authors shows that most finite
simple groups are generated by two elements of orders 2 and~3; for example, see the survey \url{http://math.nsc.ru/conference/groups2013/slides/MaximVsemirnov_slides.pdf}.
 Which finite simple groups cannot be generated by two elements of orders 2 and~3?

\makebox[15pt][r]{}In particular, is it true that, among classical simple groups of Lie type, such exceptions, apart from $PSU(3,5^2)$, arise only when the characteristic is $2$ or $3$?

\makebox[15pt][r]{}The question of which finite simple groups are (2,3)-generated remains open only
for the orthogonal groups of even dimension $2m>8$ (M.\,A.\,Pellegrini, M.\,C.\,Tamburini Bellani, {\it J.~Austral. Math. Soc.},
{\bf 117} (2024), 130--148).
\hfill \raisebox{-1ex}{\sl M.\,C.\,Tamburini}
\emp

\bmp
 \textbf{18.99.} Let ${\rm cd}(G)$ be the set of all complex irreducible character degrees of a finite group $G$. \textit{Huppert's Conjecture:}
 if $H$ is a nonabelian simple group and $G$ is any finite group such that ${\rm cd}(G)={\rm cd(H)}$, then $G\cong H\times A$, where $A$ is an abelian group.

\hfill \raisebox{-1ex}{\sl H.\,P.\,Tong-Viet}
\emp

\bmp \textbf{18.100.}
For a given class of groups $\mathfrak{X}$, let $(\mathfrak{X},\infty)^{\ast }$ denote the class of groups in which every infinite subset contains two distinct elements $x$ and $y$ that satisfy $\langle x,x^{y}\rangle\in\mathfrak{X}$.
Let $G$ be a finitely generated soluble-by-finite group in the class $(\mathfrak{X},\infty)^{\ast }$, let $m$ be a positive integer, and let $\mathfrak{F}$, $\mathfrak{E}_{m}$, $\mathfrak{E}$, $\mathfrak{N}$, and $\mathfrak{P}$ denote the classes of finite groups, groups of exponent dividing $m$, groups of finite exponent, nilpotent groups, and polycyclic groups, respectively.

\makebox[25pt][r]{a)} If $\mathfrak{X}=\mathfrak{E}\mathfrak{N}$, then is $G$ in $\mathfrak{E}\mathfrak{N}$?

\makebox[25pt][r]{b)} If $\mathfrak{X}=\mathfrak{E}_{m}\mathfrak{N}$, then is $G$ in $\mathfrak{E}_{m}(\mathfrak{F}\mathfrak{N})$?

\makebox[25pt][r]{c)} If $\mathfrak{X}=\mathfrak{N}(\mathfrak{P}\mathfrak{F})$, then is $G$ in $\mathfrak{N}(\mathfrak{P}\mathfrak{F})$?
\hf \raisebox{0ex}{\sl N.\,Trabelsi} 
 \emp

\bmp
 \textbf{18.102.} (E.\,C.\,Dade). Let $C$ be a Carter subgroup of a finite solvable group $G$, and let $\ell(C)$ be the number of primes dividing $|C|$ counting multiplicities. It was proved in (E.\,C.\,Dade, {\it Illinois J.~Math.}, {\bf 13} (1969), 449--514) that there is an exponential function $f$ such that the nilpotent length of $G$ is at most $f(\ell (C))$. Is there a linear (or at least a polynomial) function $f$ with this property?

 \makebox[15pt][r]{}{\it Editors' comment}:
 A quadratic function with that property was proved to exist in the special case where $G=H\rtimes C$ and $C$ is a cyclic subgroup such that $C_H(C)=1$, (E.\,Jabara, {\it J.~Algebra}, {\bf 487} (2017), 161--172).
 \hfill \raisebox{-1ex}{\sl A.\,Turull}
\emp

\bmp
 \textbf{18.103.} A group is said to be \emph{minimax} if it has a finite subnormal series each of whose factors satisfies either the minimum or the maximum condition on subgroups.
Is it true that in the class of nilpotent minimax groups only finitely generated groups may have faithful irreducible primitive representations over a finitely generated field of characteristic zero?
\hfill \raisebox{-1ex}{\sl A.\,V.\,Tushev}
\emp

\bmp
 \textbf{18.104.} Construct an example of a pro-$p$ group $G$ and a proper abstract normal subgroup $K$ such that $G/K$ is perfect.\hfill \raisebox{-1ex}{\sl John S.\,Wilson}
\emp

\bmp
 \textbf{18.105.} Let $R$ be the formal power series algebra in commuting indeterminates $x_1,\dots x_n, \dots$ over the field with $p$ elements, with $p$ a prime. (Thus its elements are (in general infinite) linear combinations of monomials $x_1^{a_1}\dots x_n^{a_n}$ with $n, a_1,\dots,a_n$ non-negative integers.)
Let $I$ be the maximal ideal of $R$ and $J$ the (abstract) ideal generated by all products of two elements of $I$. Is there an ideal $U$ of $R$ such that $U<I$ and $U+J=I$?

 \makebox[15pt][r]{}If so, then ${\rm SL}_3(R)$
and the kernel of the map to ${\rm SL}_3(R/U)$ provide an answer to the previous question.\hfill \raisebox{-1ex}{\sl John S.\,Wilson}
\emp

\bmp
 \textbf{18.106.} Let $R$ be a group that acts coprimely on the finite group $G$.
Let $p$ be a prime, let $P$ be the unique maximal $RC_{G}(R)$-invariant $p$-subgroup
of $G$ and assume that $C_{G}(O_{p}(G)) \leq O_{p}(G)$.
If $p>3$ then $P$ contains a nontrivial characteristic subgroup that is normal in $G$ (P.\,Flavell, {\it J.~Algebra}, {\bf 257} (2002), 249--264). Does the same result hold for the primes 2 and 3?
 \hfill \raisebox{-1ex}{\sl P.\,Flavell}
\emp

\bmp \textbf{18.107.} Is there a finitely generated infinite
residually finite $p$-group such that every subgroup of infinite
index is cyclic?

\makebox[15pt][r]{}The answer is known to be negative
for $p=2$. Note that in (M.\,Ershov, A.\,Jaikin-Zapirain,
{\it J.~Reine Angew. Math.}, {\bf 677} (2013), 71--134) it is shown that for every prime $p$ there is a finitely generated infinite residually finite $p$-group such that every finitely generated subgroup of infinite index is finite.
\hf \raisebox{-1ex}{\sl A.\,Jaikin-Zapirain}
 \emp
\bmp

\textbf{18.108.}
A group $G$ is said to have \emph{property $(\tau )$} if
its trivial representation is an isolated point in the subspace of irreducible representations with finite images (in the topological unitary dual space of $G$).
Is it true that there exists a finitely generated group satisfying property $(\tau )$ that homomorphically maps onto every finite simple group?

\makebox[15pt][r]{}The motivation is the theorem in (E.\,Breuillard, B.\,Green, M.\,Kassabov, A.\,Lubotzky, N.\,Nikolov, T.\,Tao) saying that there are $\epsilon>0$ and $k$ such that every finite simple group $G$ contains a generating set $S$ with at most $k$ elements such that the Cayley graph of $G$ with respect to $S$ is an $\epsilon$-expander. In (M.\,Ershov, A.\,Jaikin-Zapirain, M.\,Kassabov, {\it Mem. Amer. Math. Soc.}, {\bf 1186}, 2017) it is proved that there exists a group with Kazhdan's property (T) that maps onto every simple group of Lie type of rank ${}\geq 2$. (A group $G$ is said to have \emph{Kazhdan's property} (T) if
its trivial representation is an isolated point in the unitary dual space of $G$.)
\hfill \raisebox{-1ex}{\sl A.\,Jaikin-Zapirain}
\emp

\bmp
\textbf{18.109.} Is it true that a group satisfying Kazhdan's property (T) cannot homomorphically map onto infinitely many simple groups of Lie type of rank 1?

\makebox[15pt][r]{}It is known that a group which maps onto $PSL_2(q)$ for infinitely many $q$ does not have Kazhdan's property (T).
\hf \raisebox{-1ex}{\sl A.\,Jaikin-Zapirain}\emp

\bmp
 \textbf{18.110.} The \emph{non-$p$-soluble length} of a finite group $G$ is the number of non-$p$-soluble factors in a shortest normal series each of whose factors either is $p$-soluble or is a direct product of non-abelian simple groups of order divisible by $p$. For a given prime $p$ and a given proper group variety ${\frak V}$, is there a bound for the non-$p$-soluble length of finite groups whose Sylow $p$-subgroups belong to ${\frak V}$?

 \makebox[15pt][r]{}{\it Comment of 2021}:
 the existence of such a bound is proved for~$p=2$ (F.\,Fumagalli, F.\,Leinen, O.\,Puglisi, {\it Proc. London Math. Soc. (3)}, {\bf 125}, no.\,5 (2022), 1066--1082).

\hfill \raisebox{-1ex}{\sl E.\,I.\,Khukhro, P.\,Shumyatsky}
\emp

\bmp
 \textbf{\zv 18.111.}
 Let $G$ be a discrete countable group, given as a central extension $0\to {\Z}\to G\to Q\to 0$. Assume that either $G$ is quasi-isometric to $\Z\times Q$, or that $\Z\to G$ is a quasi-isometric embedding. Does that imply that $G$ comes from a bounded cocycle on $Q$?
\hfill \raisebox{-1ex}{\sl I.\,Chatterji, G.\,Mislin}

\ul

\otv No, it does not, in both cases (R.\,Frigerio, A.\,Sisto, {\it Ann. Henri Lebesgue}, {\bf  6} (2023), 225--258).
\emp

\bmp
 \textbf{18.113.}
Let $\mathfrak{M}$ be a finite set of finite simple nonabelian
groups. Is it true that a periodic group saturated with groups from $\mathfrak{M}$ (see 14.101) is isomorphic to one of groups in~$\mathfrak{M}$? The case where $\mathfrak{M}$ is one-element is of special interest.
\hfill \raisebox{-1ex}{\sl A.\,K.\,Shl\"epkin}\emp

\bmp
 \textbf{18.114.}
 Does there exist an irreducible $5'$-subgroup $G$ of $GL(V)$ for some finite ${\Bbb F} _5$-space $V$ such that the number of conjugacy classes of the semidirect product $VG$ is equal to $|V|$ but $G$ is not cyclic?

\makebox[15pt][r]{}Such a subgroup exists if 5 is replaced by $p=2, 3$ but does not exist for primes $p>5$ (\emph{J.~Group Theory}, {\bf 14} (2011), 175--199).
\hfill \raisebox{-1ex}{\sl P.\,Schmid}
\emp

\bmp
 \textbf{18.116.}
 Given a finite $2$-group $G$ of order $2^n$, does there exist a finite set $S$ of rational primes such that $|S|\le n$ and $G$ is a quotient group of the absolute Galois group
of the maximal $2$-extension of $\Q$ unramified outside $S\cup \{\infty \}$?

\makebox[15pt][r]{}This is true for $p$-groups, $p$ odd, as noted by Serre.\hfill \raisebox{0ex}{\sl P.\,Schmid}
\emp
\bmp \textbf{18.117.}
Is every element of a nonabelian finite simple group a commutator of two elements of coprime orders?

\makebox[15pt][r]{}The answer is known to be ``yes'' for the alternating groups (P.\,Shumyatsky, {\it Forum Math.}, {\bf 27}, no.\,1 (2015), 575--583) and for $PSL(2,q)$ (M.\,A.\,Pellegrini, P.\,Shumyatsky, {\it Arch. Math.}, {\bf 99} (2012), 501--507). Without the coprimeness condition, this is Ore's conjecture proved in (M.\,W.\,Liebeck, E.\,A.\,O'Brien, A.\,Shalev, Pham Huu Tiep, {\it J.~Eur. Math. Soc.}, {\bf 12}, no.\,4 (2010), 939--1008).
\hf \raisebox{-1ex}{\sl P.\,Shumyatsky}
 \emp

\bmp \textbf{18.118.} Let $G$ be a profinite group such that $G/Z(G)$ is periodic. Is $[G,G]$ necessarily periodic?
\hf \raisebox{-1ex}{\sl P.\,Shumyatsky}
 \emp

\bmp \textbf{18.119.} A multilinear 
commutator word is any commutator of weight $n$ in $n$ distinct variables. Let $w$ be a multilinear commutator word and let $G$ be a finite group. Is it true that every Sylow $p$-subgroup of the verbal subgroup $w(G)$ is generated by $w$-values?
\hf \raisebox{-1ex}{\sl P.\,Shumyatsky}
 \emp

\bmp
 \textbf{18.120.} Let $P=AB$ be a finite $p$-group
factorized by an abelian subgroup $A$ and a class-two subgroup $B$. Suppose, if
 necessary, $A \cap B=1$. Then is it true that $ \langle A, [B,B]
 \rangle =AB_{0}$, where $B_{0}$ is an abelian subgroup of $ B$?
\hfill \raisebox{-1ex}{\sl E.\,Jabara}\emp

\newpage

\pagestyle{myheadings} \markboth{19th Issue (2018)}{19th Issue
(2018)} \thispagestyle{headings} ~ \vspace{2ex}

\centerline {\Large \textbf{Problems from the 19th Issue (2018)}}
\phantomsection\label{19izd}
\vspace{4ex}

\bmp \textbf{\ai 19.1.} 
An element $g$ of a group $G$ is a
\textit{non-near generator} of $G$ if for every subset $S\subseteq G$ such that
$|G:\langle g,S\rangle|<\infty$ it follows that $|G:\langle S\rangle |<\infty$. The set of all non-near generators of $G$ forms a characteristic subgroup of $G$ called the
\textit{lower near Frattini subgroup} of $G$, denoted by $\lambda(G)$. A~subgroup $M\leq G$ is \textit{nearly maximal} in $G$ if it is maximal
with respect to being of infinite index in $G$. The intersection of all
nearly maximal subgroups forms a characteristic subgroup called the
\textit{upper near Frattini subgroup} of $G$, denoted by $\mu(G)$. In
general, $\lambda(G)\leq\mu(G)$. If $\lambda(G)=\mu(G)$, then this subgroup is called the \textit{near Frattini subgroup} of $G$, denoted by
$\psi(G)$.

\makebox[25pt][r]{a)} Is it true that $\psi(G)=1$ if $G$ is the knot group of any product of knots?

\makebox[25pt][r]{b)} Is it true that $\psi(G)=1$ if $G$ is a cable knot group?
\hfill \raisebox{-1ex}{\sl M.\,K.\,Azarian}

\ul

\oai a), b) Yes, it is true, in both cases (unverified AI solution: M.\,Trombetti,
 {\it Preprint}, 2026, \url{https://kourovkanotebookorg.wordpress.com/wp-content/uploads/2026/10/19-1.pdf}).
\emp

\bmp \textbf{19.2.} For a subgroup $L$ of a group $G$, let $L_G$ denote the largest normal subgroup of $G$ contained in $L$. If $M\vartriangleleft G$, then we say that $G$ \textit{nearly splits}
over $M$ if there is a subgroup $N\leq G$ such that $|G:N|=\infty$,
$|G:MN|<\infty$, and $(M\cap N)_G=1$.

\makebox[15pt][r]{}Let $G$ be any group, and $H$ a normal subgroup of prime order. Is it true that $\psi(G)\cap H=1$
 if and only if $G$ nearly splits over $H$?
\hfill \raisebox{-1ex}{\sl M.\,K.\,Azarian}
\emp

\bmp \textbf{19.3.} Let $G=A\ast_{H}B$ be the generalized free product of groups $A$ and $B$ with amalgamated subgroup $H$. Is $\psi(G)=1$ in the following cases?\vspace{-1.5ex}
\begin{itemize}
 \item[{a)}] $H$ is finite cyclic and $H_G=1$.\vspace{-1.5ex}

 \item[{b)}] $H$ is finite cyclic and either $\lambda(A)\cap H_G=1$ or $\lambda(B)\cap H_G=1$.\vspace{-1.5ex}

 \item[{c)}] $H_G=1$ and $H$ satisfies the minimum condition on
subgroups.\vspace{-1.5ex}

 \item[d)] $\lambda(G)\cap H=1$.\vspace{-1.5ex}

 \item[{e)}] $A$ and $B$ are free
groups, $H$ is finitely generated and at least one
of $|A:H|$ or $|B:H|$ is infinite.\vspace{-1.5ex}

 \item[{f)}] $H$ is infinite cyclic and is a retract
of $A$ and $B$.\vspace{-1.5ex}

 \item[{g)}] $G$ nearly splits over $H$, and $H$ is a normal subgroup of $G$ of prime order.\vspace{-1.5ex}
\end{itemize}

\hfill \raisebox{-1ex}{\sl M.\,K.\,Azarian}
\emp

\bmp \textbf{19.4.} a) If $G$ is the free product of infinitely many finitely
generated free groups with cyclic amalgamation, then
is $\psi(G)=1$?

\makebox[25pt][r]{b)} If $G$ is the free product of infinitely many finitely
generated free abelian groups with cyclic amalgamation, then
is $\psi(G)=1$?
\hfill \raisebox{-1ex}{\sl M.\,K.\,Azarian}
\emp

\bmp \textbf{19.5.}
If $G$ is the free product of infinitely many finitely generated
abelian groups with amalgamated subgroup $H$, then is $\psi(G)$
equal to the torsion subgroup of $H$?

\hfill \raisebox{-1ex}{\sl M.\,K.\,Azarian}
\emp

\bmp \textbf{19.6.} Let $G=A\ast_{H}B$ be the generalized free product of groups $A$ and $B$ with amalgamated subgroup $H$.

\makebox[77pt][r]{a) \textit{Conjecture}:} If both $A$ and $B$ are
nilpotent, then $\mu(G)\leq H$.

\makebox[77pt][r]{b) \textit{Conjecture}:} If $G$ is residually finite, and $H$
satisfies a nontrivial identical relation, then $\lambda(G)\leq H$.
\hfill \raisebox{-1ex}{\sl M.\,K.\,Azarian}
\emp

\bmp \textbf{19.7.}
 The group of virtual pure braids $VP_n$, $n \geq 2$, is generated by elements $\lambda_{ij}$, $1 \leq i \not= j \leq n$, and is defined by the relations
$
\lambda_{ij}\lambda_{kl}=\lambda_{kl}\lambda_{ij}$, \ $\lambda_{ki}\lambda_{kj}\lambda_{ij}=\lambda_{ij}\lambda_{kj}\lambda_{ki}
$,
where different letters denote different indices.

 \makebox[25pt][r]{a)} Construct a normal form of words in the group $VP_n$ for $n \geq 4$.

 \makebox[25pt][r]{b)} Is the group $VP_n$ linear for $n \geq 4$? It is known that the group $VP_3$ is linear.

\hfill \raisebox{-1ex}{\sl V.\,G.\,Bardakov}
\emp

\bmp \textbf{\zv 19.8.}
 A word in an alphabet $A = \{ a_1, a_2, \ldots, a_n \}$ is called a palindrome if it reads the same from left to right and from right to left. Let $k$ a non-negative integer. A~word in the alphabet $A$ is called an {\it almost $k$-palindrome} if it can be transformed into a palindrome by changing $\leq k$ letters in it. (So an almost $0$-palindrome is a palindrome.) Let elements of a free group $F_2 = \langle x, y \rangle$ be represented as words in the alphabet $\{ x^{\pm 1}, y^{\pm 1} \}$. Do there exist positive integers $m$ and $c$ such that every element in $F_2$ is a product of $\leq c$ almost $m$-palindromes?

 \makebox[15pt][r]{}It is known that for $m=0$ there is no such a number $c$. \hfill \raisebox{0ex}{\sl V.\,G.\,Bardakov}

 \ul

 \otv No, there are no such integers (M.\,Staiger, {\it J.~Algebra}, {\bf 659} (2024), 475--481).
 \emp

\bmp \textbf{\ai 19.9.} 
 Let $G$ be a finitely generated linear group.

 \makebox[15pt][r]{}a) Is it true that the non-abelian tensor square $G \otimes G$ is a linear group?

 \makebox[15pt][r]{}b) In particular, is this true for the braid group $B_n$ for $n > 3$?

 \makebox[15pt][r]{}It is known that $B_3 \otimes B_3$ is linear, and there is a countable group $G$ such that $G \otimes G$ is not linear. See the definition of the non-abelian tensor square in (R.\,Brown, J.-L.\,Loday, {\it Topology}, {\bf 26}, no.\,3 (1987), 311--335).
\hfill \raisebox{-1ex}{\sl V.\,G.\,Bardakov, M.\,V.\,Neshchadim}

\ul

\oai a) No, not always (unverified AI solution:  M.\,Trombetti, {\it Preprint}, 2026, \url{https://kourovkanotebookorg.wordpress.com/wp-content/uploads/2026/10/19-9.pdf}).

\oai b) Yes, this is true (unverified AI solution:  M.\,Trombetti, {\it Preprint}, 2026, \url{https://kourovkanotebookorg.wordpress.com/wp-content/uploads/2026/10/19-9.pdf}).
\emp 

\bmp \textbf{19.10.} \zva a)
Is it possible to embed a finitely generated non-abelian free pro-$p$ group as an open subgroup of a simple totally disconnected locally compact group?

 \makebox[15pt][r]{}b) If the answer is yes, can we require that the simple envelope is also compactly generated?

 \makebox[15pt][r]{}{\it Comment of 2025}: Partial progress on part b) is made in (Y.\,Barnea, M.\,Ershov, A.\,Le\,Boudec, C.\,D.\,Reid,
M.\,Vannacci, T.\,Weigel, {\it Preprint}, 2025, \url{https://arxiv.org/pdf/2507.04120}).
\hfill \raisebox{-1ex}{\sl Y.\,Barnea}

\ul

\otv a) Yes, it is possible (Y.\,Barnea, M.\,Ershov, A.\,Le\,Boudec, C.\,D.\,Reid,
M.\,Vannacci, T.\,Weigel, {\it Preprint}, 2025, \url{https://arxiv.org/pdf/2507.04120}).
\emp

\bmp \textbf{19.11.}
 Does there exist a constant $c$ such that the
number of conjugacy classes in a finite group $G$ is always at least $c \log_2 |G|$?

 \makebox[15pt][r]{}{\it Editors' comment of 2021}: It is proved that every group $G$ contains at least $\varepsilon \log |G|/(\log \log |G|)^8$ conjugacy classes for some fixed $\varepsilon > 0$ (L.\,Pyber, \emph{J.~London Math. Soc. (2)}, \textbf{46}, no.\,2 (1992), 239--249). It is also proved that for every $\varepsilon > 0$ there exists $\delta > 0$ such that every finite group $G$ of order at least 3 has at least $\delta \log _2|G|/(\log_2 \log _2|G|)^{3+\varepsilon}$ conjugacy classes (B.\,Baumeister, A.\,Mar\'oti, H.\,P.\,Tong-Viet, \emph{Forum Math.}, \textbf{29}, no.\,2 (2017), 259--275).
\hfill \raisebox{-1ex}{\sl E.\,Bertram} 
\emp

\bmp \textbf{19.12.} A finite group $G$ is called conjugacy-expansive if for every normal
subset $N$ and conjugacy class $C$ of $G$ the normal set $NC$ contains at least
as many conjugacy classes of $G$ as $N$ does. Is it true that every finite
simple group is conjugacy-expansive?

\hfill \raisebox{-1ex}{\sl M.\,Bezerra, Z.\,Halasi, A.\,Mar\'oti, S.\,Sidki}
\emp

\bmp \textbf{19.13.}
 The well-known Baer--Suzuki theorem
 states that if every two conjugates of an element $a$ of
a finite group $G$ generate a finite $p\hskip0.2ex$-sub\-group,
then $a$ is contained in a normal $p\hskip0.2ex$-sub\-group. Does
such a theorem hold in the class of periodic groups for $p>2$?

\makebox[15pt][r]{}A counterexample is known for $p=2$; see Archive 11.11a).
\hfill \raisebox{-0ex}{\sl A.\,V.\,Borovik}
\emp

\bmp \textbf{19.14.}
Let $G$ be a finite group such that the quasivariety generated by all its Sylow subgroups contains only finitely many subquasivarieties. Is it true that the quasivariety generated by $G$ also contains only finitely many subquasivarieties?

\hfill \raisebox{-1ex}{\sl A.\,I.\,Budkin}
\emp

\bmp \textbf{19.15.}
For a class of groups $M$, let $L(M)$ denote the class of all groups $G$ in which
the normal closure $\langle a^G\rangle $ of any element $a\in G$ is contained in $M$. Let $qA$ be the quasivariety generated by a finite nilpotent group $A$. Can the quasivariety $L(qA)$ contain a non-nilpotent group?
\hfill \raisebox{-1ex}{\sl A.\,I.\,Budkin}
\emp

\bmp \textbf{19.16.} (C.\,Drutu and M.\,Sapir). Is every free-by-cyclic group of the form $F_n \rtimes \Z$ linear?
\hfill \raisebox{-1ex}{\sl J.\,O.\,Button}
\emp

\bmp \textbf{19.17.} (Well-known problem). Suppose that $G$ is a finitely presented group such that the set of first Betti numbers over all the
finite index subgroups of $G$ is unbounded above. (Here, the first Betti number is the maximum $n$ for which there is a surjective homomorphism onto $\Z^n$.)

 \makebox[25pt][r]{a)} Must $G$ have a finite index subgroup with a surjective
homomorphism to a non-abelian free group?

 \makebox[25pt][r]{b)} Can $G$ even be soluble?
\hfill \raisebox{-0ex}{\sl J.\,O.\,Button}
\emp

\bmp \textbf{19.18.} Given an infinite set $\Omega$, define an
algebra $A$ (the {\it reduced incidence algebra of finite
subsets\/}) as follows. Let $V_n$ be the set of functions from the
set of $n $-element subsets of $\Omega$ to the rationals $\Bbb Q$.
Now let $A=\bigoplus V_n$, with multiplication as follows: for
$f\in V_n$, $g\in V_m$, and $|X|=m+n$, let $(fg)(X)=\sum
f(Y)g(X\setminus Y)$, where the sum is over the $n $-element
subsets $Y$ of $X$. If $G$ is a permutation group on $\Omega$, let
$A^G$ be the algebra of $G $-inva\-ri\-ants in $A$.

 \makebox[15pt][r]{}Suppose that $G$ has no finite orbits on $\Omega$. It is known that then $A^G$ is an integral domain (M.\,Pouzet, \textit{Theor. Inf. App.}
\textbf{42}, no.\,1 (2008), 83--103). Is it true that the quotient of $A^G$ by the ideal generated by constant functions is also an integral domain?
\hfill \raisebox{-1ex}{\sl P.\,J.\,Cameron}
\emp

\bmp \textbf{19.19.} A finite transitive
permutation group $G$ is said to have the \emph{road closure property} if,
given any orbit $O$ of $G$ on $2$-sets, and any proper block of imprimitivity
for $G$ acting on $O$, the graph with edge set $O\setminus B$ is connected.
Such a group must be primitive, and basic (not contained in a wreath product
with the product action); it cannot have an imprimitive subgroup of index~$2$.
In addition, such a group cannot be one of the permutation groups arising from triality (whose socle is $D_4(q)$ and intersects the point stabiliser in the parabolic subgroup corresponding to the three leaves in the Coxeter--Dynkin diagram for $D_4$).

 \makebox[15pt][r]{}Classify the basic primitive groups $G$ which do not have the road closure property. In particular, is it true that such a group either has a subgroup of index at most 3 or is almost simple?
\hfill \raisebox{-1ex}{\sl P.\,J.\,Cameron}
\emp

\bmp \textbf{\zv 19.20.}
For a finite group $G$, let $\mathrm{End}(G)$ denote
the semigroup of endomorphisms of~$G$, and $\mathrm{PIso}(G)$ the semigroup
of partial isomorphisms of $G$ (isomorphisms between subgroups of $G$). If
$G$ is abelian, then $|\mathrm{End}(G)|=|\mathrm{PIso}(G)|$. Is the converse
true?

\hfill \raisebox{-1ex}{\sl P.\,J.\,Cameron}

\ul

\otv No, it is not (K.\,Muliarchyk,  {\it Preprint}, 2026, \url{https://kourovkanotebookorg.wordpress.com/wp-content/uploads/2026/09/kourovka_19.20_solution.pdf}).
\emp

\bmp \textbf{19.21.}
Can a non-discrete, compactly generated, topologically simple, locally compact group be amenable and non-compact?

 \makebox[15pt][r]{}(A positive answer implies a negative answer to 19.107.)
\hfill \raisebox{0ex}{\sl P.-E.\,Caprace}
\emp

\bmp \textbf{\zv 19.22.}
Let $\mathscr{S}$ denote the class of nondiscrete compactly generated, topologically simple totally disconnected locally compact groups.
Two topological groups are called {\em locally isomorphic} if they contain isomorphic open subgroups. Is the number of local isomorphism classes of groups in $\mathscr{S}$ uncountable?
\hfill \raisebox{-1ex}{\sl P.-E.\,Caprace}

\ul

\otv Yes, it is (I.\,Castellano, J.\,Fari\~{n}a-Asategui, M.\,E.\,Garciarena,
B.\,Marchionna, M.\,Quick, C.\,D.\,Reid, S.\,M.\,Smith, S.\,Tornier,
M.\,Vannacci, J.\,S.\,Wilson, {\it Preprint}, 2026, \url{https://arxiv.org/abs/2609.17410}).
\emp

\bmp \textbf{\zv 19.25.}
Let $G$ and $H$ be finite groups of the same order with $\sum_{g\in G}\varphi (|g|) = \sum_{h\in H} \varphi (|h|)$, where $\varphi$ is the Euler totient function. Suppose that $G$ is simple. Is $H$ necessarily simple?
\hf \raisebox{-1ex}{\sl B.\,Curtin, G.\,R.\,Pourgholi} 

\ul

\otv No, it need not be (AI solution confirmed by GAP:  \ P.\,Monticone, {\it Preprint}, 2026,

\url{https://kourovkanotebookorg.wordpress.com/wp-content/uploads/2026/05/19_25.pdf}); see also (W.\,van\,Doorn, E.\,Judin, P.\,Monticone, D.\,Morrison,  {\it Preprint}, 2026, \url{https://arxiv.org/abs/2607.17477}).
\emp

\bmp \textbf{19.26.} (Y.\,O.\,Hamidoune).
Suppose that $A$ and $B$ are finite subsets of a group $G$
such that $\vert A \vert \ge 2$ and $\vert B \vert \ge 2$, and let
$A\cdot_2B$ denote the set of elements of $G$ which can be expressed in the form
 $ab$ for at least two different $(a,b) \in A \times B$. Is it true that
$\vert A \vert + \vert B \vert - (1/2) \vert AB \vert - (1/2) \vert A\cdot_2B\vert
 \le \max\{2, \vert gH \vert \,:\, H \le G,\; g \in G,\; gH \subseteq A\cdot_2B\}$?

 \makebox[15pt][r]{}This is proved if $G$ is abelian.
\hfill \raisebox{-1ex}{\sl W.\,Dicks}
\emp

\bmp \textbf{19.27.} (Well-known question). A finitely generated group $G$ that acts on a tree in such a way that all vertex and edge stabilizers are infinite cyclic groups is called a generalized Baumslag--Solitar group. The Bass--Serre theory gives finite presentations of such groups. Is the isomorphism problem soluble for
 generalized Baumslag--Solitar groups?
\hfill \raisebox{-1ex}{\sl F.\,A.\,Dudkin}

\emp

\bmp \textbf{19.28.}
 Let $G$ be a group, and $\varphi$ an automorphism of $G$. Elements $x,y\in G$ are said to be $\varphi$-conjugate if $x=z^{-1}y\varphi (z)$ for some $z\in G$. The $\varphi$-conjugacy is an equivalence relation; the number of $\varphi$-conjugacy classes is denoted by $R(\varphi )$.

 \makebox[15pt][r]{}{\it Conjecture}: If a finitely generated residually finite group $G$ has an automorphism $\varphi$ such that $R(\varphi )$ is finite, then $G$ has a soluble subgroup of finite index.

 \makebox[15pt][r]{}The conjecture is proved if $\varphi$ has prime order (E.\,Jabara, {\it J.~Algebra}, {\bf 320}, no.\,10 (2008), 3671--3679). If $G$ is infinitely generated, then $G$ does not have to be almost solvable (K.\,Dekimpe, D.\,Gon\c{c}alves, {\it Bull. London Math. Soc.}, {\bf 46}, no.\,4 (2014), 737--746), but if it is of finite upper rank, then it has to be almost solvable
(E.\,Troitsky, {\it J. Group Theory}, {\bf 28}, no.\,1 (2025), 151--164).
 \hf \raisebox{-1ex}{\sl A.\,L.\,Fel'shtyn, E.\,V.\,Troitsky}
\emp

\bmp \textbf{19.29.}
 Let $\Phi \in {\rm Out\,} G = {\rm Aut}\, G/{\rm Inn}\, G$. Two automorphisms $\varphi,\psi\in\Phi$ are said to be isogradient if $\varphi=\hat g^{-1}\psi\hat g$ for some inner automorphism $\hat g$. The number of isogradiency classes in $\Phi$ is denoted by $S(\Phi )$.

 \makebox[15pt][r]{}{\it Conjecture}: If a finitely generated residually finite group $G$ has an outer automorphism $\Phi$ such that $S(\Phi )$ is finite, then $G$ has a soluble subgroup of finite index.

\hf \raisebox{-1ex}{\sl A.\,L.\,Fel'shtyn, E.\,V.\,Troitsky}
\emp

\bmp \textbf{19.30.} An element $g$ of a finite group $G$ is said to be vanishing if $\chi(g)=0$ for some irreducible complex character $\chi\in {\rm
Irr}(G)$. Must a finite group and a finite simple group be isomorphic if they have equal orders and the same set of orders of vanishing elements?
\hfill \raisebox{-1ex}{\sl M.\,Foroudi Ghasemabadi, A.\,Iranmanesh, }
\emp

\bmp \textbf{\zv 19.31.}
Let $\omega(G,S)$ be the exponential growth rate of
a group $G$ with a finite generating set $S$ (see 14.7).
Let $MCG(\Sigma_g)$ be the mapping class group
of the orientable surface $\Sigma_g$ of genus $g$.
Is there a constant $C>1$ such that
$\omega(MCG(\Sigma_g), S) \ge C$
for every $g >0$ and every finite generating set $S$ of $MCG(\Sigma_g)$?
\hfill \raisebox{-1ex}{\sl K.\,Fujiwara}

\ul

\otv  No, there is no such constant that works for all genera; but if the genus $g$ is fixed, then such a constant $C=C(g)>1$ does exist (J.\,Mangahas,
{\it Geom. Funct. Anal.} {\bf 19}, no.\,5 (2010), 1468--1480).
\emp

\bmp \textbf{19.32.}
Let $p\geq 673$ be a prime and $r \geq 2$ be an integer. Let $B$ be the free group in the Burnside variety of exponent $p$ on $a_1, \dots , a_r$. Must every non-identical one-variable equation $w(a_1, \dots , a_r, x) = 1$ over $B$ have at most finitely many solutions in~$B$?\hf \raisebox{-1ex}{\sl A.\,Gaglione}
\emp

\bmp \textbf{19.33.}
{\it Conjecture}: Let $G$ be a finite group, $p$ a prime number, and $P$ a Sylow $p$-subgroup of $G$. Suppose that an irreducible ordinary character $\chi$ of $G$ has degree divisible by $p$. If the restriction $\chi_P$ of $\chi$ to $P$ has a linear constituent, then $\chi_P$ has at least $p$ different linear constituents.

 \makebox[15pt][r]{}This conjecture has been verified for symmetric, alternating, $p$-solvable, and sporadic simple groups.
\hfill \raisebox{-1ex}{\sl E.\,Giannelli}
\emp

\bmp \textbf{19.34.}
Let $G$ be a finitely generated group such that for every element $g\in G$ the set of commutators $\{[x,g]\mid x\in G\}$ is a subgroup of $G$. Is it true that $G$ is residually nilpotent?

 \makebox[15pt][r]{}The answer is affirmative if $G$ is finite.
\hfill \raisebox{0ex}{\sl D.\,Gon\c{c}alves, T.\,Nasybullov}
\emp

\bmp \textbf{\zv 19.38.}
 Suppose that $H$ is a
 subgroup of a finite soluble group $G$ that covers all Frattini chief factors of $G$ and avoids all complemented chief factors of $G$. Is it true that there are elements $x,y\in G$ such that $H\cap H^x\cap H^y=H_G$, where $H_G$ is the largest normal subgroup of $G$ contained in $H$?

 \makebox[15pt][r]{}This is true if $H$ is a prefrattini subgroup of $G$.
\hfill \raisebox{0ex}{\sl S.\,F.\,Kamornikov}

\ul

\otv Yes, it is true (S.\,F.\,Kamornikov, O.\,L.\,Shemetkova, {\it
J.~Algebra}, {\bf 641} (2024), 1--8).
\emp

\bmp \textbf{19.39.} A group is said to be $\frak X$-critical if it does not belong to $\frak X$ but all its proper subgroups belong to~$\frak X$.
Let $\frak F$ be a soluble hereditary formation of finite groups for which all $\frak F$-critical groups are either groups of prime order or $\frak U$-critical groups, where $\frak U$ is the formation of all supersoluble finite groups. Must $\frak F$ be a saturated formation?

\hfill \raisebox{-1ex}{\sl S.\,F.\,Kamornikov}
\emp

\bmp \textbf{\zv 19.41.}
Let $\phi$ be an automorphism of a free group $F_n$ of finite rank, and let $F_n\rtimes_\phi \mathbb Z$ be the split extension of $F_n$ by $\mathbb Z$ with $\mathbb Z$ acting as $\langle \phi\rangle$. Is the group $F_n\rtimes_\phi \mathbb Z$ conjugacy separable?
\hfill \raisebox{-1ex}{\sl I.\,Kapovich}

\ul

\otv Yes, it is (F.\,Dahmani, S.\,Hughes, M.\,Kudlinska, N.\, Touikan, {\it Preprint}, 2026, \url{https://arxiv.org/abs/2509.22346v2}).
\emp

\bmp \textbf{19.42.}
Suppose that $H$ is a word-hyperbolic subgroup of a word-hyperbolic group $G$ such that the inclusion of $H$ to $G$ extends to a continuous $H$-equivariant map $j:\partial H\to \partial G$ between their hyperbolic boundaries. If such an extension exists, it is unique and $j$ is called the \emph{Cannon--Thurston map}.

 \makebox[25pt][r]{a)} Is it true that for every point $p\in \partial G$ its full preimage $j^{-1}(p)$ is finite?

 \makebox[25pt][r]{b)} Moreover, is it true that there is a number $N=N(G,H)<\infty$ such that for every $p\in \partial G$ the full preimage $j^{-1}(p)$ consists of at most $N$ points?

 \makebox[15pt][r]{}The answer to both questions is ``yes'' in all the cases where the Cannon--Thurston map is known to exist and where it has been possible to analyze the multiplicity of this map. This includes the original set-up considered by Cannon and Thurston where $H$ is the surface group, and $G$ is the fundamental group of a closed hyperbolic 3-manifold fibering over the circle with that surface as a fiber. In this case, $j:\mathbb S^1\to\mathbb S^2$ is a uniformly finite-to-one continuous surjective `Peano curve'.
\hfill \raisebox{-1ex}{\sl I.\,Kapovich}
\emp

\bmp \textbf{\zv 19.43.}
Suppose that $\varphi $ is an automorphism of a finite soluble group $G$. Is the Fitting height of $G$ bounded in terms of $|\varphi |$ and $|C_G(\varphi )|$?
\hfill \raisebox{-1ex}{\sl E.\,I.\,Khukhro}

\ul

\otv Yes, it is (E.\,Khukhro,
\emph{Bull. London Math. Soc.}, \textbf{58}, no.~1 (2026), Article ID e70206).
\emp

\bmp \textbf{19.44.} By definition a profinite group has finite rank at most $r$ if every subgroup of it can be (topologically) generated by $r$ elements.
Suppose that for every element $g$ of a profinite group $G$ there is a closed subgroup $E_g$ of finite rank such that for every $x\in G$ all sufficiently long Engel commutators $[x,g,\dots, g]$ belong to $E_g$, that is, for every $x\in G$ there is a positive integer $n(x,g)$ such that $[x,g,\dots, g]\in E_g$ whenever $g$ is repeated $\geq n(x,g)$ times. Is it true that $G$ has a normal subgroup $N$ of finite rank with locally nilpotent quotient $G/N$?
\hfill \raisebox{-1ex}{\sl E.\,I.\,Khukhro}
\emp

\bmp \textbf{19.45.}
Let \(G\) be a group generated by 3 class transpositions (see the definition in 17.57), and let \(m\) be the least common multiple of the moduli of the residue classes interchanged by the generators of~\(G\). Assume that \(G\) does not setwisely stabilize any union of residue classes modulo~\(m\) except for~\(\varnothing\) and~\(\mathbb{Z}\), and assume that the integers
\(0, 1, \dots, 42\) all lie in the same orbit under the action of \(G\) on \(\mathbb{Z}\). Is the action of \(G\) on~\(\mathbb{N}\cup \{0\}\) necessarily transitive?

 \makebox[15pt][r]{}The bound 42 cannot be replaced by a smaller number, since the finite group \(\langle \tau_{0(2),1(2)}, \tau_{0(3),2(3)}, \tau_{0(7),6(7)} \rangle\) acts transitively on the set \(\{0, \dots, 41\}\), as well as on the set of residue classes modulo~42.
\hfill \raisebox{-1ex}{\sl S.\,Kohl}
\emp

\bmp \textbf{19.46.}
Does the group \({\rm CT}(\mathbb{Z})\) have finitely generated infinite periodic subgroups? (See the definition of \({\rm CT}(\mathbb{Z})\) in 17.57).
\hfill \raisebox{-1ex}{\sl S.\,Kohl}
\emp

\bmp \textbf{\ai 19.47.} 
Let $K$ be a finite extension of degree $n$ of a field $k$ of odd characteristic. The multiplicative group $K^{*}$ embeds into the group of all $k$-linear automorphisms $\mathrm{Aut}_{k}(K)$ by the rule $t(x) = tx$ for all $x\in K$. For a fixed basis of $K$ over $k$, the group $\mathrm{Aut}_k(K)$ is isomorphic to
$GL(n,k)$. The image of $K^{*}$ is called a \textit{nonsplit
 maximal torus} corresponding to the extension $K/k$ and is denoted by $T(K/k)$. A subgroup of $GL(n,k)$ is said to be \textit{rich in transvections} if it contains all elementary
transvections.

 \makebox[15pt][r]{}Let $H$ be a subgroup of $GL(n, k)$ which contains $T(K/k)$ and a
one-dimensional transformation. Is $H$ rich in transvections?
\hfill \raisebox{-1ex}{\sl V.\,A.\,Ko\u{\i}baev}

\ul

\oai No, not always (unverified AI solution: B.\,G.\,Di Siena, {\it Preprint}, 2026, \url{https://www.proofforum.org/abs/PF:2026.000035}).
\emp

\bmp \textbf{\zv 19.48.}
A system of additive subgroups $\sigma_{ij}$, $1\leq i,j\leq
n$, of a field $K$ is called a \textit{net} (or a \textit{carpet}) of order
$n$ if $\sigma_{ir}\sigma_{rj}\subset\sigma_{ij}$ for all $i, r, j$.
A net that does not contain the diagonal is called an \textit{elementary net}. A net $\sigma=(\sigma_{ij})$ is said to be \textit{irreducible} if all $\sigma_{ij}$ are nontrivial. An elementary net
$\sigma$ is said to be \textit{closed} if the elementary net subgroup
$E(\sigma)$ does not contain additional elementary transvections.

 \makebox[15pt][r]{}Let $R$ be a principal ideals domain with $1\in R$, let $k$ be the field of fractions of $R$, and $K$ an algebraic extension of the field $k$. Let $\sigma=(\sigma_{ij})$ be an irreducible elementary net
$\sigma=(\sigma_{ij})$ over $K$ such that all $\sigma_{ij}$ are
$R$-modules. Is the net $\sigma$ closed?

\hfill \raisebox{-1ex}{\sl V.\,A.\,Ko\u{\i}baev}

\ul

\otv  No, not always  (AI solution confirmed by V.\,A.\,Ko\u{\i}baev:  \ M.\,Brescia, {\it Preprint}, 2026,
\url{https://kourovkanotebookorg.wordpress.com/wp-content/uploads/2026/10/19-48.pdf}). A more general result was obtained earlier (V.\,A.\,Ko\u{\i}baev, {\it unpublished} (in Russian), \url{https://kourovkanotebookorg.wordpress.com/wp-content/uploads/2026/10/19.48-koib.pdf}).
\emp 

\bmp \textbf{19.51.} A finite group is \textit{monomial} if each irreducible character of it is induced from a linear character of some subgroup. Monomial groups are soluble (Taketa). The group is \textit{normally} (\textit{subnormally}) \textit{monomial} if each irreducible character is induced from a linear character of some normal (respectively, subnormal) subgroup. Metabelian groups are normally monomial, and abelian-by-nilpotent groups are subnormally monomial. In the other direction, it is known that there exist normally monomial groups of arbitrarily large derived length.

\makebox[15pt][r]{}a) For a given prime $p$, do there exist normally monomial finite $p$-groups of arbitrarily large derived length?

\makebox[15pt][r]{}b) Do there exist subnormally monomial groups of arbitrarily large nilpotence length?
\hfill \raisebox{-1ex}{\sl A.\,Mann}
\emp

\bmp \textbf{19.52.}
 The Gruenberg--Kegel graph (or the prime graph) $GK(G)$ of a finite group $G$ has vertex set consisting of all prime divisors of the order of $G$, and different vertices $p$ and $q$ are adjacent in $GK(G)$ if and only if the number $pq$ is an element order in~$G$. Is there a finite non-solvable group $G$ such that $GK(G)$
 does not contain 3-cocliques
 and is not isomorphic to the Gruenberg--Kegel graph of any finite solvable group?

 \makebox[15pt][r]{}It is known that there are no examples of such groups $G$ among almost simple groups.
\hfill \raisebox{-1ex}{\sl N.\,V.\,Maslova}
\emp

\bmp \textbf{19.53.}
Let $G$ be a group generated by elements $x,y,z$ such that
$$x^3=y^2=z^2=(xy)^3=(yz)^3=1\quad \text{and}\ \quad g^{12}=1 \text{ for all }g\in G.$$
Is it true that $|G|\leq 12$?

\makebox[15pt]{}This is true if $G$ is finite (D.\,V.\,Lytkina, V.\,D.\,Mazurov, {\it Siberian Math. J.}, {\bf 56}, no.\,3 (2015), 471--475).
\hfill \raisebox{-1ex}{\sl V.\,D.\,Mazurov}
\emp

\bmp \textbf{19.54.}
What are the chief factors of a finite group in which every
2-maximal subgroup is not $m$-maximal for any~$m\ge 3$?

\makebox[15pt]{}A subgroup $H$ of a group $G$ is said to be $m$-maximal if there is a chain of subgroups $ H=H_0<H_1< \ldots < H_{m-1}< H_m=G$ in which
$H_i$ is maximal in $H_{i+1}$ for every~$i$. Note that for every $m\ge 3$ there is a finite group in which some 2-maximal subgroup is $m$-maximal.
\hfill \raisebox{-1ex}{\sl V.\,S.\,Monakhov}
\emp

\bmp \textbf{\ai 19.56.} 
A $\star$-commutator is a commutator $[x,y]$ of two elements $x,y$ of coprime prime-power orders. Is a finite group~$G$ soluble if $|ab|\ge |a||b|$ for any $\star$\nobreakdash-\hspace{0pt}commutators~$a$ and $b$ of coprime orders?
\hfill \raisebox{-1ex}{\sl V.\,S.\,Monakhov}

\ul

\oai Yes, it is (unverified AI solution:   M.\,Kinyon, J.\,Urban, {\it Preprint}, 2026,
\url{https://kourovkanotebookorg.wordpress.com/wp-content/uploads/2026/09/kourovka-mathematics.pdf}).
\emp 

\bmp \textbf{\zv 19.57.}
What are the non-abelian composition factors of a finite group $G$
in which every maximal subgroup is simple or $p$\nobreakdash-\hspace{0pt}nilpotent for some fixed odd prime $p\in \pi (G)$?

 \makebox[15pt][r]{}{\it Comment of 2026}: These factors are determined in the cases where $G$ is  non-simple, or alternating, or sporadic simple  (A.\,Beltr\'an, F.\,Viudez,  {\it Preprint}, 2026, \url{https://arxiv.org/abs/2607.06434}).
\hfill \raisebox{-1ex}{\sl V.\,S.\,Monakhov, V.\,N.\,Tyutyanov}

\ul

\otv These factors are described (AI solution confirmed by  V.\,N.\,Tyutyanov: \ R.\,Sater,  {\it Preprint}, 2026, \url{https://kourovkanotebookorg.wordpress.com/wp-content/uploads/2026/09/kourovka-19-57-19-58-v1.2.1.pdf }).
\emp 

\bmp \textbf{\zv 19.58.}
What are the non-abelian composition factors of a finite group
in which every maximal subgroup is simple or $p$-decomposable
for some fixed odd prime $p\in \pi (G)$?

 \makebox[15pt][r]{}{\it Comment of 2026}: These factors are determined in the cases where $G$ is  non-simple, or alternating, or sporadic simple  (A.\,Beltr\'an, F.\,Viudez,  {\it Preprint}, 2026, \url{https://arxiv.org/abs/2607.06434}).
\hfill \raisebox{-1ex}{\sl V.\,S.\,Monakhov, V.\,N.\,Tyutyanov}

\ul

\otv These factors are described (AI solution confirmed by  V.\,N.\,Tyutyanov: \ R.\,Sater,  {\it Preprint}, 2026, \url{https://kourovkanotebookorg.wordpress.com/wp-content/uploads/2026/09/kourovka-19-57-19-58-v1.2.1.pdf }).
\emp 

\bmp \textbf{19.59.}
 Does the group of isometries of $\Q^3$ have a free non-abelian subgroup such that every non-trivial element acts without nontrivial fixed points in $\Q^3$?

 \makebox[15pt][r]{}The answer is yes if the field of rational numbers $\Q$ is replaced by the field $\Bbb R$ of real numbers (see G.\,Tomkowicz, S.\,Wagon, {\it The Banach--Tarski Paradox}, Cambridge Univ. Press, 2016.)
 \hfill \raisebox{-1ex}{\sl J.\,Mycielski}
\emp

\bmp \textbf{19.60.}
Let $R$ be a principal ideal domain, let $d$ be a positive integer, and let $X_1,\, \ldots,\, X_m$ be members of $\mathrm{SL}(d,R)$ (or of $\mathrm{GL}(d,R)$). Is it decidable whether or not these $m$ matrices freely generate a free group? If it is, design an effective algorithm for computing the answer.
\hfill \raisebox{-1ex}{\sl P.\,M.\,Neumann}
\emp

\bmp \textbf{\zv 19.61.}
Let $\mathfrak{A}=\{\mathfrak{A}_r\ |\ r\in \Phi\}$ be an elementary carpet of type $\Phi$ over a commutative ring $K$
(see 7.28), and let $\Phi(\mathfrak{A})=\langle x_r(\mathfrak{A}_r)\ |\ r\in\Phi\rangle$ be its carpet subgroup. Define the {\it closure} of the carpet~$\mathfrak{A}$ to be the set of additive subgroups
$\overline{\mathfrak{A}}=\{\overline{\mathfrak{A}}_r\ |\ r\in \Phi\},
$ where $ \overline{\mathfrak{A}}_r=\{t\in K\ |\ x_r(t)\in \Phi(\mathfrak{A})\}.
$ Is the closure $\overline{\mathfrak{A}}$ of a carpet $\mathfrak{A}$ always a carpet?

\makebox[15pt][r]{}An affirmative answer is known if $\Phi=A_l, D_l, E_l$.
\hfill \raisebox{-1ex}{\sl Ya.\,N.\,Nuzhin}

\ul

\otv Yes, it is (AI solution confirmed by Ya.\,N.\,Nuzhin: \   M.\,Kinyon, J.\,Urban, {\it Preprint}, 2026,
\url{https://kourovkanotebookorg.wordpress.com/wp-content/uploads/2026/09/kourovka-mathematics.pdf}).
\emp 

\bmp \textbf{\zv 19.62.} 
Let $\mathfrak{A}=\{\mathfrak{A}_r\ |\ r\in \Phi\}$ be an elementary carpet of type $\Phi$ of rank $l\geq 2$ (see 7.28). For $p\in \Phi$, define a set of additive subgroups
 $ \mathfrak{B}_p=\sum C_{ij,rs}\mathfrak{A}_r^i\mathfrak{A}_s^j$, where the sum is taken over all natural numbers $i$, $j$ and roots $r,s\in \Phi$ such that $ir+js=p$. It is known that the set $\mathfrak{B}=\{\mathfrak{B}_p\ |\ p\in \Phi\}$ is a carpet called the {\it derived carpet} of $\mathfrak{A}$.
 It is also known that for $\Phi=A_l$ the set $\mathfrak{B}$ is a closed (admissible) carpet, which means that its carpet subgroup does not contain new root elements.
Is every derived carpet of type $\Phi$ over a commutative ring
closed (admissible)?

\makebox[15pt][r]{}{\it Comment of 2025}: An affirmative answer
was obtained for $\Phi$ of type $B_l$, $C_l$, or $F_4$ when $\mbox{GCD}(p,2)=1$, and for $\Phi$ of type $G_2$ when $\mbox{GCD}(p,6)=1$ (Ya.\,N.\,Nuzhin, {\it J.~Siberian Fed. Univ. Ser. Math. Phys.}, {\bf 16}, no.\,6 (2023), 732--737). \hfill \raisebox{-1ex}{\sl Ya.\,N.\,Nuzhin}

\ul

\otv Yes, it is (AI solution confirmed by Ya.\,N.\,Nuzhin: \  M.\,Kinyon, J.\,Urban, {\it Preprint}, 2026,
\url{https://kourovkanotebookorg.wordpress.com/wp-content/uploads/2026/09/kourovka-mathematics.pdf}).
\emp 

\bmp \textbf{\zv 19.63.}
Let $\mathfrak{A}=\{\mathfrak{A}_r\ |\ r\in \Phi\}$ be an elementary carpet of type $\Phi$ over a commutative ring $K$
(see 7.28) and let $\mathfrak{A}_r^2=\{t^2\mid t\in\mathfrak{A}_r\}$. Are the inclusions
$\mathfrak{A}_r^2\mathfrak{A}_{-r}\subseteq\mathfrak{A}_r$, $r\in\Phi$, sufficient for the carpet $\mathfrak{A}$ to be closed (admissible)?
 \hfill \raisebox{-1ex}{\sl Ya.\,N.\,Nuzhin}

 \ul

 \otv Yes, they are sufficient (Ya.\,N.\,Nuzhin, {\it J.~Siberian Fed. Univ. Ser. Math. Phys.}, {\bf 16}, no.\,6 (2023), 732--737; \ Ya.\,N.\,Nuzhin, {\it to appear in Siberian Math.~J.}).

\emp

\bmp \textbf{19.64.}
Let $G$ be a group, and $(g_1, \ldots , g_n)$ a tuple of its elements. The {\em type} of this tuple in $G$, denoted $Tp^G(g_1, \ldots , g_n)$, is the set of all first order formulas in free variables $x_1, \ldots ,x_n$ in the standard group theory language which are true on $(g_1, \ldots , g_n)$ in $G$.
Two groups $G$ and $H$ are called {\it isotypic} if for every tuple of elements $\bar h=(h_1,\ldots,h_n)$ in $H$ there is a tuple $\bar g =(g_1,\ldots,g_n)$ in $G$, such that $Tp(\bar h)=Tp(\bar g)$ and vice versa, for every tuple $\bar g$ in $G$ there is a tuple $\bar h$ in $H$ such that $Tp(\bar h)=Tp(\bar g)$.

 \makebox[15pt][r]{}Is it true that every two isotypic finitely generated groups are isomorphic?

 \makebox[15pt][r]{}The answer is positive if one of the finitely generated groups is free (R.\,Sklinos), abelian (G.\,Zhitomirski), virtually polycyclic, metabelian, free solvable (A.\,Myasnikov, N.\,Romanovskii), co-Hopfian (in particular homogeneous), finitely presented Hopfian or geometrically Noetherian Hopfian (R.\,Sklinos). In particular, torsion-free hyperbolic groups, braid groups, linear groups, mapping class groups of compact surfaces, limit groups, quasi-cyclic groups, $Out(F_n)$, $n > 2$, irreducible arithmetic lattices (which are not virtually free groups) in semi-simple Lie groups, are of such kind (R.\,Sklinos).
\hfill \raisebox{-1ex}{\sl B.\,I.\,Plotkin}
\emp

\bmp \textbf{19.65.}
It is well known that varieties of groups form a free semigroup $N$ (A.\,L.\,Shmel'kin, {\it DAN SSSR}, {\bf 149} (1963), 543--545 (Russian)); \ B.\,H.\,Neuman, H.\,Neumann, P.\,M.\,Neumann, {\it Math. Z.}, {\bf 80} (1962), 44--62). Varieties of linear representations over an infinite field of zero characteristic also form a free semi\-group~$M$, and $N$ acts freely on $M$ (B.\,I.\,Plotkin, {\it Siberian Math. J.}, {\bf 13}, no.\,5 (1972), 713--729) 
A similar theorem holds for Lie algebras (V.\,A.\,Parfenov, Algebra Logika, {\bf 6}, no.\,4 (1967), 61--73 (Russian); \ L.\,A.\,Simonyan, {\it Siberian Math.~J.}, {\bf 29}, no.\,2 (1988), 276--283). 

 \makebox[15pt][r]{}Are there other varieties of algebras $\Theta$ where the same situation takes place?

\hfill \raisebox{-1ex}{\sl B.\,I.\,Plotkin}
\emp

\bmp \textbf{19.66.}
We say that a variety $\Theta$ is of \textit{Tarski type} if any two non-abelian $\Theta$-free groups of finite rank are elementarily equivalent.

\makebox[25pt][r]{a)} Find examples of Tarski type varieties distinct from the variety of all groups.

\makebox[25pt][r]{b)} Is it true that the {Burnside variety $B_n$} of all groups of exponent $n$, where $n$ is big enough, is of Tarski type?

\makebox[25pt][r]{c)} Is it true that the {\it $n$-Engel variety $E_n$} of all groups satisfying the identity $[[[x,\underbrace{y],y],\ldots,y}_{n}]\equiv 1$, where $n$ is big enough, is of Tarski type?

\hfill \raisebox{2ex}{\sl \smash{B.\,I.\,Plotkin, E.\,B.\,Plotkin}}
\emp

\bmp \textbf{19.68.} For a finite group $G$ and a permutation group $K$, let $b_G(K)$ denote the number of conjugacy classes of regular subgroups of $K$ isomorphic to~$G$. Does there exist a function~$f$ such that
$b_G(K)\le n^{f(r)}$ for every abelian group $G$ of order $n$ and rank~$r$, and every group~$K$ such that $K^{(2)}=K$? (See Archive 19.67 for the definition of~$K^{(2)}$.)
\hfill \raisebox{-1ex}{\sl I.\,Ponomarenko}
\emp

\bmp \textbf{19.69.} (P.\,Wesolek).
The acronym tdlc stands for totally disconnected and locally compact, and tdlcsc for tdlc and second countable. Let $\rm{Res}(G)$ denote the intersection of all open normal subgroups of a topological group $G$. The class of elementary tdlcsc groups is defined as the smallest class $\mathscr{E}$ of tdlcsc groups such that\vspace{-1.5ex}
\begin{itemize}
 \item[(1)] $\mathscr{E}$ contains all second countable profinite groups and countable discrete groups;\vspace{-1.5ex}
 \item[(2)] $\mathscr{E}$ is closed under taking closed subgroups, Hausdorff quotients, directed unions of open subgroups, and group extensions.\vspace{-1.5ex}
\end{itemize}
This class admits a well-behaved, ordinal-valued decomposition rank $\xi$ defined recursively as follows: $\xi(\{1\}) = 1$, and if $G \in \mathscr{E}$ is a union of an increasing sequence $(O_i)$ of compactly generated open subgroups, then $\xi(G) = \sup_i \{\xi(\rm{Res}(O_i))\} + 1$.

 \makebox[15pt][r]{}What is the supremum of the decomposition ranks of elementary tdlcsc groups? In particular, is it countable?
\hfill \raisebox{-1ex}{\sl C.\,Reid}
\emp

\bmp \textbf{19.70.} (P.\,Wesolek and P.-E.\,Caprace).
Let $\mathscr{S}$ denote the class of nondiscrete compactly generated, topologically simple tdlc groups. Let $G$ be a non-elementary tdlcsc group. Is there a compactly generated closed subgroup $H \leq G$ such that $H$ has a continuous quotient in $\mathscr{S}$?
\hfill \raisebox{-1ex}{\sl C.\,Reid}
\emp

\bmp \textbf{19.71.} (G.\,Willis).
Can a group $G$ in $\mathscr{S}$ be such that every element of $G$ normalizes a~compact open subgroup of $G$?
\hfill \raisebox{-1ex}{\sl C.\,Reid}
\emp

\bmp \textbf{19.72.}
Let $G$ be a group in $\mathscr{S}$. Can every element of $G$ have trivial contraction group?
\hfill \raisebox{-1ex}{\sl C.\,Reid}
\emp

\bmp \textbf{19.73.} (P.-E.\,Caprace and N.\,Monod).
Is there a compactly generated, locally compact group that is topologically simple, but not abstractly simple?

\hfill \raisebox{-1ex}{\sl C.\,Reid, S.\,M.\,Smith}
\emp

\bmp \textbf{\zv 19.74.}
A subgroup $H$ of a group $G$ is called {\em pronormal\/} if $H$
and $H^g$ are conjugate in $\langle H,H^g\rangle$ for every $g\in
G$. A subgroup $H$ of a group $G$ is called {\em abnormal\/} if $g\in \langle H,H^g\rangle$ for every $g\in G$.

 \makebox[15pt][r]{}Does there exist an infinite group that does not contain nontrivial proper pronormal subgroups?

 \makebox[15pt][r]{}The question is equivalent to the following: does there exist an infinite simple group that does not contain proper abnormal subgroups?\hfill \raisebox{-1ex}{\sl D.\,O.\,Revin}

\ul

\otv Yes, it does (S.\,Corson, {\it Monatsh. Math.} (2025) \url{https://doi.org/10.1007/s00605-025-02116-8}).
\emp

\bmp \textbf{19.76.}
A semigroup presentation is called {\em tree-like} if all
relations have the form $a=bc$ where $a, b, c$ are letters and no two
relations share the left-hand side or the right-hand side.
Is it decidable whether the semigroup given by a finite tree-like
presentation contains an idempotent?

 \makebox[15pt][r]{}This is equivalent to the question whether the
closure of a finitely generated subgroup of R.\,Thompson's group $F$
contains an isomorphic copy of $F$.
\hfill \raisebox{-1ex}{\sl M.\,Sapir}
\emp

\bmp \textbf{19.77.} {\it The Stable Small Cancellation Conjecture}: If $ F$ is a free group of rank $r\geq 2$, then, for any fixed $n > 0$, there exists a generic subset $S$ of $F^n$ such that,
for any automorphism $\phi$ of $F$, the set $\phi(S)$ satisfies the small cancellation condition $C'(1/6)$.

\hfill \raisebox{-1ex}{\sl P.\,E.\,Schupp}
\emp

\bmp \textbf{19.78.}
 Does there exist an infinite simple subgroup of $SL(2,{\Bbb Q})$?

 \makebox[15pt][r]{}This is a special case of a question which Serge Cantat asked me about ``non-algebraic'' simple subgroups of $GL(n,{\Bbb C})$. It is also a special case of 15.57.
\hfill \raisebox{-1ex}{\sl J.-P.\,Serre}
\emp

\bmp \textbf{19.79.} (T.\,Springer). Let $G$ be a group. Suppose that for every integer $n > 0$ the group $G$ has a unique (up to isomorphism) irreducible complex linear representation of dimension $n$. (Note that $G = SL(2,{\Bbb Q})$ has these properties, by a theorem of Borel--Tits ({\it Ann. Math.}, {\bf 97} (1973), 499--571).) Is it true that $G$ has a normal subgroup $N$ such that
 $G/N$ is isomorphic to $SL(2,{\Bbb Q})$?
\hfill \raisebox{-1ex}{\sl J.-P.\,Serre}
\emp

\bmp \textbf{19.82.} Let $h^*(G)$ denote the generalized Fitting height of a finite group $G$ defined as the minimum number $k$ such that $F^*_k(G)=G$, where $F^*_1(G)=F^*(G)$ is the generalized Fitting subgroup of $G$, and by induction $F^*_{i+1}(G)$ is the inverse image of $F^*(G/F^*_i(G))$. If $G$ is soluble, then $h^*(G)=h(G)$ is the Fitting height of $G$.

 \makebox[15pt][r]{}Does every finite group $G$ contain a soluble subgroup $K$ such that $h^*(G)=h(K)$?

\hf \raisebox{-1ex}{\sl P.\,Shumyatsky}
\emp

\bmp \textbf{19.83.}
An element $g$ of a group $G$ is almost Engel if there is a finite set ${\mathscr E}(g)$ such that for every $x\in G$ all sufficiently long commutators $[x,\,{}_n g]$ belong to ${\mathscr E}(g)$, that is, for every $x\in G$ there is a positive integer $n(x,g)$ such that $[x,\,{}_n g]\in {\mathscr E}(g)$ whenever $n\geq n(x,g)$. By a linear group we understand a subgroup of $GL(m,F)$ for some field $F$ and a positive integer $m$.

 \makebox[15pt][r]{}Is the set of almost Engel elements in a linear group always a subgroup?

\hfill \raisebox{-1ex}{\sl P.\,Shumyatsky}
\emp

\bmp \textbf{19.86.}
Suppose that a finite group $G$ has a
$\sigma_{i}$-Hall subgroup for every $i \in I$. Suppose that a subgroup $A\leq G$ is such that $A\cap H$ is a
 $\sigma_{i}$-Hall subgroup of $A$ for every $i\in I$ and every $\sigma_{i}$-Hall subgroup $H$ of $G$ (see 19.84 in Archive). Is it true that then $A$ is ${\sigma}$-subnormal in $G$?

\makebox[15pt][r]{}An affirmative answer is known if $\sigma =\{\{2\}, \{3 \}, \ldots \}$
and if $G$ is $\sigma$-soluble.

\hfill \raisebox{-1ex}{\sl A.\,N.\,Skiba}
\emp

\bmp \textbf{\zv 19.89.}
A permutation group is {\em subdegree-finite} if every orbit of every point stabiliser is finite. Let $\Omega$ be countably infinite, and let $G$ be a closed and subdegree-finite subgroup of ${\rm Sym}(\Omega)$ regarded as a topological group under the topology of pointwise convergence.

 \makebox[15pt][r]{}{\it Conjecture}: If the minimal degree of $G$ is infinite, then there is some subset of $\Omega$ whose setwise stabiliser in $G$ is trivial.

 \makebox[15pt][r]{}Note that this conjecture implies Tom Tucker's well-known Infinite Motion Conjecture for graphs.
\hfill \raisebox{-1ex}{\sl S.\,M.\,Smith}

\ul

\otv The conjecture is proved (F.\,Lehner, {\it Bull. London Math. Soc.}, {\bf  55} (2023), 1874--1879).
\emp

\bmp \textbf{19.90.}
A {\it skew brace} is a set $B$ equipped with two operations $+$ and $\cdot$ such that $(B,+)$ is an additively written (but not necessarily abelian) group, $(B,\cdot )$ is a multiplicatively written group, and $a\cdot (b+c)=ab-a+ac$ for any $a,b,c\in B$.

 \makebox[25pt][r]{b)} Is there a skew brace with non-soluble additive group but nilpotent multiplicative group?

 \makebox[25pt][r]{c)} Is there a finite skew brace with soluble additive group but non-soluble multiplicative group?
\hfill \raisebox{-1ex}{\sl A.\,Smoktunowicz, L.\,Vendramin}
\emp

\bmp \textbf{19.91.}
Let $G$ be a finite group with an abelian Sylow $p$-subgroup $A$. Suppose that
$B$ is a strongly closed elementary abelian subgroup of $A$. Without invoking the
Classi\-fication Theorem for Finite Simple Groups (CFSG), prove that $G$ has a normal
subgroup $N$ such that $B =\Omega _1(A\cap N)$.

 \makebox[15pt][r]{}For $p = 2$, this is a corollary of a theorem of Goldschmidt. For $p$ odd, this has been proved by Flores and Foote ({\it Adv. Math.}, {\bf 222} (2009), 453--484), but their proof relies on CFSG. A~CFSG-free proof in the special case when $p = 3$ and $A$ has 3-rank 3 would already be quite interesting. This case arises in Aschbacher's treatment of the $e(G) = 3$ problem, and a proof would provide an alternative to part of his argument. (For this application, one could assume that all proper simple sections of $G$ are known.)
\hfill \raisebox{-1ex}{\sl R.\,Solomon}
\emp

 \bmp \textbf{\zv 19.92.}
Let $A$ be the algebra of $3\times 3$ skew-Hermitian matrices over the real octonions, where multiplication is given by bracket product. Determine ${\rm Aut}(A)$.

\makebox[15pt][r]{}The question might be of interest for octonion algebras over a different base field or ring. The question was inspired by an observation of John Faulkner concerning a possible connection between $A$ or a related algebra and the Dwyer--Wilkerson 2-compact group $BDI(4)$.
\hfill \raisebox{-1ex}{\sl R.\,Solomon}

\ul

\otv It is determined: $Aut(A)\cong G_2 \times SO(3)$, which follows from the result about derivations in (H.\,Petyt, {\it Commun. Algebra}, {\bf 47}, no.\,10 (2019), 4216--4223).
\emp

\bmp \textbf{19.93.}
{\it Conjecture}:
There exists a function $f:\mathbb{N}\times\mathbb{N}\to \mathbb{N}$ such that, if $G$ is a finite $p$-group with $d$ generators and $G$ has no epimorphic images isomorphic to the wreath product $C_p\wr C_p$, then each factor of the $p$-lower central series of $G$ has order bounded above by $f(p,d)$

 \makebox[25pt][r]{\zv a)} for $p\geq 7$;

 \makebox[25pt][r]{b)} for $p =2,3,5$.

 \makebox[15pt][r]{}It is interesting to compare this conjecture with the celebrated characterization of Aner Shalev of finitely generated $p$-adic analytic pro-$p$-groups.
\hfill \raisebox{-1ex}{\sl P.\,Spiga}

\ul

\otv a)  The conjecture is refuted for $p\geq 7$  (AI solution confirmed by P.\,Spiga: \   M.\,Kinyon, J.\,Urban, {\it Preprint}, 2026,
\url{https://kourovkanotebookorg.wordpress.com/wp-content/uploads/2026/09/kourovka-mathematics.pdf}).
\emp 

\bmp \textbf{\zv 19.94.}
By definition two elements $g_1$ and $g_2$ of a free metabelian group $F$ with basis $\{x_1,\ldots, x_n\}$ have the same distribution if the equations $g_1(x_1,\ldots,
 x_n) = g$ and $g_2(x_1,\ldots, x_n) = g$ have the same number of solutions in any finite metabelian group $G$ for every $g \in G$. Is it true that elements $g_1$ and $g_2$ have the same distribution if and only if they are conjugate by some automorphism of $F?$
 \hfill \raisebox{-1ex}{\sl E.\,I.\,Timoshenko}

 \ul

\otv No, it is not true (V.\,Ionin, A.\,Semidetnov,  {\it Preprint}, 2026, \url{https://arxiv.org/pdf/2608.29219}).
 \emp

\bmp \textbf{19.95.}
Let $G_\Gamma$ be a partially commutative soluble group of derived length $n \geq 3$ with defining graph $\Gamma$ (the definition is similar to the case of $n=2$, see 17.104). Is it true that the centralizer of any vertex of the graph is generated by its adjacent vertices?

 \hfill \raisebox{-1ex}{\sl E.\,I.\,Timoshenko}
 \emp

\bmp \textbf{19.96.}
Let $M$ be a free metabelian group of finite rank $r \geq 2$, and let $P$ be the set of its primitive elements. (An element of a relatively free group is said to be primitive if it can be included in a basis of the group.)

 \makebox[25pt][r]{a)} Is $P$ a first-order definable set?

 \makebox[25pt][r]{b)} Is the set of bases
 of the group
 $M$ a first-order definable set?

 \hfill \raisebox{-1ex}{\sl E.\,I.\,Timoshenko}
\emp

\bmp \textbf{19.97.} Let $G$ be an $m$-generated group the elementary theory of which $Th(G)$ coincides with the elementary theory
$Th(F)$ of a free soluble group $F$ of derived length $n
\geq 3$ of finite rank $r \geq 2$. Must the groups $G$ and
$F$ be isomorphic?

 \makebox[15pt][r]{}The answer is affirmative if $m \leq r$ or $n=2$
 \hfill \raisebox{-1ex}{\sl E.\,I.\,Timoshenko}
\emp

\bmp \textbf{19.99.}
 Is it true that for any positive integer $t$ there is a positive integer $n(t)$ such that any finite group with at least $n(t)$ conjugate classes of soluble maximal subgroups of Fitting height at most $t$ is itself a soluble group of Fitting height at most $t$?

 \hfill \raisebox{-1ex}{\sl A.\,F.\,Vasil'ev, T.\,I.\,Vasil'eva}
\emp

\bmp \textbf{\zv 19.100.}
 Suppose that a finite group $G$ admits a factorization $G=AB=BC=CA$, where $A, B, C$ are abnormal supersoluble subgroups. Is $G$ supersoluble?

 \hfill \raisebox{-1ex}{\sl A.\,F.\,Vasil'ev, T.\,I.\,Vasil'eva}

 \ul

 \otv Yes, it is (L.\,S.\,Kazarin, V.\,N.\,Tyutyanov, in: {\it Collection of surveys to the 50th anniversary of the Mathematical Faculty of Yaroslavl' University}, 2026, 123--125 (Russian);
 \url{https://kourovkanotebookorg.wordpress.com/wp-content/uploads/2026/10/19.100.pdf}).
\emp

\bmp \textbf{\zv 19.102.}
A subgroup $H$ of a free group $F$ is called \textit{inert} if $r(H\cap K)\leq r(K)$ for every $K\leq F$; and \textit{compressed} if $r(H)\leq r(K)$ for every $H\leq K\leq F$.

 \makebox[15pt][r]{}Is it true that compressed subgroups are inert?
\hfill \raisebox{0ex}{\sl E.\,Ventura}

\ul

\otv Yes, it is true (A.\,Jaikin-Zapirain, \emph{Preprint}, 2024, \url{https://arxiv.org/pdf/2403.09515 }).
\emp

\bmp \textbf{\zv 19.103.}
Is it true that an intersection of compressed subgroups is compressed?

\makebox[15pt][r]{}It is known that arbitrary intersections of inert subgroups are inert.
\hfill \raisebox{-1ex}{\sl E.\,Ventura}

\ul

\otv Yes, it is true (A.\,Jaikin-Zapirain, \emph{Preprint}, 2024, \url{https://arxiv.org/pdf/2403.09515 }).
\emp

\bmp \textbf{\zv 19.104.}
Is there an algorithm which decides whether a given subgroup of $F$ is inert?

\makebox[15pt][r]{}An algorithm to decide whether a given $H$ is compressed is known.
\hfill \raisebox{-1ex}{\sl E.\,Ventura}

\ul

\otv Yes, there is (A.\,Jaikin-Zapirain, \emph{Preprint}, 2024, \url{https://arxiv.org/pdf/2403.09515 }).
\emp

\bmp \textbf{\zv 19.105.}
Is it true that the fixed subgroups of endomorphisms of $F$ are inert?

\hfill \raisebox{-1ex}{\sl E.\,Ventura}

\ul

\otv Yes, it is true (A.\,Jaikin-Zapirain, \emph{Preprint}, 2024, \url{https://arxiv.org/pdf/2403.09515 }).
\emp

\bmp \textbf{\ai 19.106.} 
Let $G$ be a uniformly locally finite group (which means that there is a function $f$ on the natural numbers such that any subgroup generated by $n$ elements has size at most $f(n)$). Suppose that for any two definable subgroups $H$ and $K$, the intersection $H\cap K$ has finite index either in $H$ or in $K$. Is $G$ necessarily nilpotent-by-finite? Or even finite-by-abelian-by-finite?
\hfill \raisebox{-1ex}{\sl F.\,Wagner}

\ul

\oai No, it is not always nilpotent-by-finite or  finite-by-abelian-by-finite (unverified AI solution: B.\,G.\,Di Siena, {\it Preprint}, 2026, \url{https://www.proofforum.org/abs/PF:2026.000031}).
\emp

\bmp \textbf{19.107.}
Let $\mathscr{AE}$ be the smallest class of locally compact groups such that\vspace{-1.5ex}
\begin{itemize}
 \item[(1)] $\mathscr{AE}$ contains the compact groups and the discrete amenable groups;\vspace{-1.5ex}
 \item[(2)] $\mathscr{AE}$ is closed under taking closed subgroups, Hausdorff quotients, directed unions of open subgroups, and group extensions.\vspace{-1.5ex}
\end{itemize}
Is every amenable locally compact group an element of $\mathscr{AE}$?

\makebox[15pt][r]{}(A negative answer would follow from a positive answer to 19.21.)
\hfill \raisebox{-0ex}{\sl P.\,Wesolek}
\emp

\bmp \textbf{19.108.}
Let $P$ be a finite $p$-group, where $p$ is an odd prime. Let $\chi$ be a
complex irreducible character of $P$. If $\chi(x)\ne 0$ for some
$x\in P$, is it true that the order of $x$ must divide $|P|/\chi(1)^2$?
\hfill \raisebox{-1ex}{\sl T.\,Wilde}
\emp

\bmp \textbf{19.110.}
 (G.\,M.\,Bergman and A.\,Magidin). Do there exist varieties of groups in which the relatively free
group of rank 2 is finite, and the relatively free group of rank 3 is
infinite?
\hf \raisebox{-1ex}{\sl P.\,Zusmanovich} 
\emp

\bmp \textbf{19.111.} Do there exist infinite groups
all of whose proper subgroups are cyclic of order~$p$ that do not satisfy any
nontrivial group identity except $x^p = 1$ and its consequences?

\makebox[15pt][r]{}Note
that there exist infinite
groups all of whose proper subgroups are infinite cyclic that do not satisfy
any nontrivial group identity (due to A.\,Olshanskii; see P.\,Zusmanovich, {\it J.~Algebra}, {\bf 388} (2013), 268--286, Remark after Theorem~6.1).

\hf \raisebox{-1ex}{\sl P.\,Zusmanovich} 
\emp

\newpage
\pagestyle{myheadings} \markboth{20th Issue (2022)}{20th Issue
(2022)} \thispagestyle{headings} ~ \vspace{2ex}

\centerline {\Large \textbf{Problems from the 20th Issue (2022)}}
\phantomsection\label{20izd}
\vspace{4ex}

\bmp \textbf{\zv 20.1.}
For $k\geq1$, a group $G$ is said to be \emph{totally $k$-closed} if in each of its faithful permutation representations, say on a set $\Omega$, $G$ is the largest subgroup of ${\rm Sym}(\Omega)$ which leaves invariant each of the $G$-orbits in the induced action on the set of ordered $k$-tuples $\Omega^k$.

\makebox[15pt][r]{}Are there any finite insoluble totally $2$-closed groups with nontrivial Fitting subgroup?

\makebox[15pt][r]{}The finite totally $2$-closed groups which are either soluble or have trivial Fitting subgroup are known.
 \hf\raisebox{-1ex}{\sl M.\,Arezoomand, M.\,A.\,Iranmanesh, C.\,E.\,Praeger, G.\,Tracey}

 \ul

\otv Yes, there are (AI solution confirmed by G.\,Tracey: \  M.\,Kinyon, J.\,Urban, {\it Preprint}, 2026,
\url{https://kourovkanotebookorg.wordpress.com/wp-content/uploads/2026/09/kourovka-mathematics.pdf}). Another solution follows from earlier work: by Theorem 1.1 in (M.\,Arezoomand, M.\,A.\,Iranmanesh, C.\,E.\,Praeger, G.\,Tracey, {\it Bull. Math. Sci.}, {\bf 14}, no.\,1 (2024), 2350004), the sporadic groups $J_1, J_3, J_4, Ly, Th$ and $M$ are totally 2-closed. Cyclic groups are totally 2-closed by Theorem 2 in (A.\,Abdollahi, M.\,Arezoomand, {\it J.~Algebra Appl.}, {\bf  17}, no.\, 4 (2018), 1850065). For $K$ one of the listed sporadic groups, and  $H$ any cyclic group of order coprime to $|K|$, let $G=H\times K$, and let $\Omega$ be a set on which $G$ acts faithfully. Total 2-closure forces $H^{(2),\Omega}=H$ and $K^{(2),\Sigma}=K$, where $\Sigma=\{\alpha^H : \alpha\in\Omega\}$. Then by Proposition~2 in (A.\,Abdollahi, M.\,Arezoomand, {\it J.~Algebra Appl.}, {\bf  17}, no.\, 4 (2018), 1850065), $G^{(2),\Omega}=G$. Hence, $G$ is totally 2-closed. For instance, $J_1 \times  C_p$ is totally 2-closed for all primes $p$ in $\{13,17,p : p>19\}$.
\emp

\bmp \textbf{\zv 20.2.}
Are there any nonabelian simple groups of Lie type which are totally $3$-closed?

\makebox[15pt][r]{}There are exactly six nonabelian simple totally $2$-closed groups --- all are sporadic groups, the largest being the Monster.

 \hf\raisebox{-1ex}{\sl M.\,Arezoomand, M.\,A.\,Iranmanesh, C.\,E.\,Praeger, G.\,Tracey}

 \ul

 \otv  Yes, there are (T.\,Gong, M.\,R.\,Zeng, Y.\,Yang,  {\it Preprint}, 2026, \url{https://arxiv.org/abs/2608.17878}).
\emp

\bmp \textbf{20.4.} A finite group $G$ is said to be \emph{cut} (or \emph{inverse semi-rational}) if $\langle x \rangle = \langle y \rangle$ implies that $x$ is conjugate to $y$ or to $y^{-1}$ for all $x, y \in G$.

\makebox[25pt][r]{a)} Let $\mathbb{Q}(G)$ denote the field extension of the rationals obtained by adjoining all entries of the ordinary character table of $G$. Is there $c > 0$ such that $|\mathbb{Q}(G):\mathbb{Q}| \leqslant c$ for all cut groups? This is true if one assumes in addition that $G$ is solvable (J.\,F.\,Tent, \emph{J.~Algebra}, \textbf{363} (2012), 73--82).

\makebox[25pt][r]{b)} Is a Sylow $3$-subgroup of a cut group also a cut group?

\makebox[25pt][r]{c)} Let $\operatorname{O}_p(G)$ denote the largest normal $p$-subgroup of $G$. Let $G$ be a solvable cut group. Is it true that for $p \in \{5, 7\}$ the exponent of $\operatorname{O}_p(G)$ divides~$p$?

\makebox[15pt][r]{}For additional information and some positive results see (\emph{Adv. Group Theory Appl.}, \textbf{8}B, 2020, 157--160 or \texttt{https://arxiv.org/abs/2001.02637}).
\hf \raisebox{-1ex}{\sl A.\,B{\"a}chle}
\emp

\bmp \textbf{20.5.}
(Well-known problem). Let $G$ be a finite group, and $\mathrm{V}(\mathbb{Z}G)$ the group of normalized units of the integral group ring of~$G$. Do the spectra of $G$ and $\mathrm{V}(\mathbb{Z} G)$ coincide? That is, is it true that, for any integer $n$, there is an element of order $n$ in $\mathrm{V}(\mathbb{Z} G)$ if and only if there is an element of order $n$ in $G$?

\makebox[15pt][r]{}The answer is ``yes'' if $G$ is solvable (M.\,Hertweck, \emph{Comm. Algebra} \textbf{36} (2008), 3585--3588).
\hf \raisebox{-1ex}{\sl A.\,B{\"a}chle, L.\,Margolis}
\emp

\bmp \textbf{20.6.}
 (W.~Kimmerle).
 The prime graph (or Gruenberg--Kegel graph) $\Gamma(X)$ of a group $X$ has vertices labeled by primes appearing as orders of elements in $X$; two distinct primes $p$ and $q$ are adjacent in $\Gamma(X)$ if and only if $X$ contains an element of order $pq$.
Denote by $\mathrm{V}(\mathbb{Z}G)$ the group of normalized units of the integral group ring of a group $G$. Is it true that for each finite group $G$ the prime graphs of $G$ and $\mathrm{V}(\mathbb{Z}G)$ coincide?

\makebox[15pt][r]{}This question has been reduced to almost simple groups (W.~Kimmerle, A.~Konovalov, \emph{Internat. J. Algebra Comput.}, \textbf{27} (2017), 619--631).
\hf \raisebox{-1ex}{\sl A.\,B{\"a}chle, L.\,Margolis}
\emp

\bmp \textbf{\ai 20.7.} 
(W.\,Boone and G.\,Higman) Does every finitely generated group with solvable word problem embed into a finitely presented simple group? It is known that every such group embeds into a simple subgroup of a finitely presented group (W.\,Boone, G.\,Higman, \emph{J.~Austral. Math. Soc.}, \textbf{18}, no.\,1 (1974), 41--53).
\hf
\raisebox{-1ex}{\sl J.\,Belk}

\ul

\oai Yes, it does (unverified AI solution: OpenAI, {\it Preprint}, 2026,

\url{https://kourovkanotebookorg.wordpress.com/wp-content/uploads/2026/10/20-7.pdf}).
\emp 

\bmp \textbf{20.8.}
Suppose $K < H < F$ are free groups of finite rank such that
$\operatorname{rank}(H) < \operatorname{rank}(K)$, but all proper subgroups of $H$ which contain $K$ have ranks ${}\geq \operatorname{rank}(K)$. Then is the inclusion of $H$ in $F$ the only homomorphism $H\rightarrow F$ fixing all elements of $K$?
 \hf\raisebox{-1ex}{\sl G.\,M.\,Bergman}
\emp

\bmp \textbf{20.9.}
In the group algebra of a free group over a field, does every element whose
support in the group has cardinality more than 1 generate a proper 2-sided ideal?

\makebox[15pt][r]{}This question is one of several related
questions in (G.\,M.\,Bergman, \emph{Commun. Algebra}, \textbf{49}, no.~9 (2021), 3760--3776).
 \hf\raisebox{-1ex}{\sl G.\,M.\,Bergman}
\emp

\bmp \textbf{20.10.}
(a) If ${\mathscr U}$ and ${\mathscr U'}$ are nonprincipal ultrafilters on $\N$, can every group
which can be written as a homomorphic image of an
ultraproduct of groups with respect to ${\mathscr U}$ also be written as a homomorphic image of an
ultraproduct of groups with respect to ${\mathscr U'}$?

\makebox[25pt][r]{\zv (b)}
If the answer to (a) is negative, is it at least true that for any
two nonprincipal ultrafilters ${\mathscr U}$ and ${\mathscr U'}$ on $\N$, there exists a nonprincipal ultrafilter ${\mathscr U''}$ on $\N$
such that every group which can be written as a homomorphic image of an ultraproduct of
groups with respect to ${\mathscr U}$ or with respect to ${\mathscr U'}$ can be written as a homomorphic image of
an ultraproduct with respect to ${\mathscr U''}$?

\makebox[15pt][r]{}A positive answer to (b) would imply that the class of groups which can be written
as homomorphic images of nonprincipal countable ultraproducts of groups is closed under
finite direct products. See (G.\,M.\,Bergman, \emph{Pacific J.~Math.}, \textbf{274} (2015), 451--495).
 \hf\raisebox{-1ex}{\sl G.\,M.\,Bergman}

 \ul

 \otv (b) The affirmative answer is consistent with the ZFC axioms of set theory (S.\,M.\,Corson, {\it Preprint}, 2025, \url{https://arxiv.org/pdf/2503.09228}).
 \emp

\bmp \textbf{20.11.}
Let $F \leq H$ be free groups such that there exists a free group
$G$ of finite rank with $F\leq G\leq H$, and let $r$
be the least of the ranks of such groups $G$. Which,
if any, of the following statements must hold?

\makebox[25pt][r]{(i)\,} There is a largest $G$ of rank $r$ between $F$ and $H$.

\makebox[25pt][r]{(i$'$)\,} For any two $G_1$, $G_2$ of rank $r$ between $F$ and $H$,
 the subgroup $\langle G_1, G_2\rangle$ has rank $r$.

\makebox[25pt][r]{(ii)\,} There is a smallest $G$ of rank $r$ between $F$ and $H$.

\makebox[25pt][r]{(ii$'$)\,} For any two $G_1$, $G_2$ of rank $r$ between $F$ and $H$,
 the subgroup $G_1\cap G_2$ has rank $r$.

\makebox[15pt][r]{}If (i) holds for all such $F$ and $H$, then so does (i$'$). If (ii) holds for all $F$, $H$, then so
does (ii$'$). The converse of the former implication holds because
subgroups of a free group of any fixed finite rank satisfy ACC,
but I don't see a way to get the converse of the other implication.
 \hf\raisebox{-1ex}{\sl G.\,M.\,Bergman}
\emp

\bmp \textbf{20.12.}
Let us say that a group $G$ has the \emph{unique $n$-fold product} property (u.-$n$-p.)
if for every $n$-tuple of finite nonempty subsets $A_1, \dots , A_n\subseteq G$ there
exists $g\in G$ which can be written in one and only one way as $g = a_1\cdots a_n$ with $a_i\in A_i$. It is easy to see that u.-$n$-p. implies u.-$m$-p. for $n \geq m$ (since some of the $A_i$ can be $\{1\}$).

\makebox[15pt][r]{}Are the conditions u.-$n$-p. ($n\geq 2$) all equivalent
 to u.-2-p., the usual unique product condition?
 \hf\raisebox{-1ex}{\sl G.\,M.\,Bergman}
\emp

\bmp \textbf{\zv 20.13.}
{(a)} If an abelian group can be written as a homomorphic
image of a nonprincipal countable ultraproduct of \emph{not necessarily} abelian groups $G_i$, must
it be a homomorphic image of a nonprincipal countable ultraproduct of abelian groups? See (G.\,M.\,Bergman, \emph{Pacific J.~Math.}, \textbf{274} (2015), 451--495).

\makebox[25pt][r]{(b)} If an abelian group can be written as a homomorphic image of a \emph{direct product} of an infinite family of not necessarily abelian \emph{finite} groups, can it be written as a homomorphic image of a direct product of finite abelian groups?

\makebox[15pt][r]{}To see that neither question is trivial, choose for each $n > 0$ a finite group $G_n$ which is
perfect but has elements which cannot be written as products of fewer than $n$ commutators.
Then both the direct product of the $G_n$ and any nonprincipal ultraproduct of those groups
will have elements which are not products of commutators; hence its abelianization $A$ will
be nontrivial. There is no evident family of abelian groups from which to obtain $A$ as an
image of a nonprincipal ultraproduct, nor a family of finite abelian groups from which to
obtain $A$ as an image of a direct product.

 \hf\raisebox{-1ex}{\sl G.\,M.\,Bergman}

 \ul

 \otv (a) Yes, it must (S.\,M.\,Corson, {\it Preprint}, 2025, \url{https://arxiv.org/pdf/2503.09228}).

 \otv (b) Yes it can (S.\,M.\,Corson, {\it Preprint}, 2025, \url{https://arxiv.org/pdf/2503.09228}).
\emp

\bmp \textbf{20.14.}
(a) Do there exist a variety $\mathfrak{V}$ of groups and a
group $G\in \mathfrak{V}$ such that the coproduct in $\mathfrak{V}$ of two copies of $G$ is embeddable in $G$, but the coproduct of three such copies is not? See (G.\,M.\,Bergman, \emph{Indag. Math.}, \textbf{18} (2007), 349--403).

\makebox[15pt][r]{}Given an embedding $G *_\mathfrak{V} G \to G$, one might expect the induced map \vspace{-1ex}
\begin{center}
 $G*_\mathfrak{V} (G *_\mathfrak{V} G) \to G *_\mathfrak{V} G \to G $
\end{center}
\vspace{-1ex} to be an embedding. But this is not automatic, because in a general group
variety $\mathfrak{V}$, a map $G *_\mathfrak{V} A \to G *_\mathfrak{V} B$ induced by an embedding $A\to B$ is not necessarily again an embedding.

 \makebox[25pt][r]{(b)} If there exist $\mathfrak{V}$ and $G$ as in (a), does there
in fact exist an example with $\mathfrak{V}$ the variety of groups generated by $G$? For this
and related questions, see (G.\,M.\,Bergman, \emph{Algebra Number Theory}, \textbf{3} (2009), 847--879).
 \hf\raisebox{-1ex}{\sl G.\,M.\,Bergman}
\emp

\bmp \textbf{20.15.}
Let $\kappa$ be an infinite cardinal. If a residually finite
group $G$ is embeddable in the full permutation group of a set of cardinality $\kappa$, must it be
embeddable in the direct product of $\kappa$ finite groups? See (G.\,M.\,Bergman, \emph{Indag. Math.}, \textbf{18} (2007), 349--403).

\makebox[15pt][r]{}The converse is true: any such direct product is residually finite and embeddable
in the indicated permutation group. Also, the statement asked for becomes true if one
replaces ``direct product of $\kappa$ finite groups'' by ``direct product of $2^\kappa$ finite groups'', since $G$ has cardinality at most $2^\kappa$, hence homomorphisms to
 that many finite groups can be chosen which, together,
 separate each element of $G$ from $e$.
 \hf\raisebox{-1ex}{\sl G.\,M.\,Bergman}
\emp

\bmp \textbf{20.16.} By Hartley's theorem (based on CFSG), if a simple locally finite group does not
contain a particular finite group as a section, then it is isomorphic to a group of Lie type over a locally finite field.

\makebox[25pt][r]{a)} Is the same true for an infinite
definably simple locally finite group? A group is said to be \textit{definably simple} if it does not contain proper definable normal subgroups.

\makebox[25pt][r]{b)} The same question for an infinite
definably simple locally finite group with bounded derived lengths of its solvable subgroups.
\hf \raisebox{-1ex}{\sl A.\,V.\,Borovik}
\emp

\bmp \textbf{20.17.}
 The \textit{normal covering number} $\gamma(G)$ of a finite non-cyclic group $G$ is the minimum number of proper subgroups of $G$ such that $G$ is the union of their conjugates.\vs

 \makebox[25pt][r]{\zv (ae)} What is the exact value of \ $\displaystyle\liminf_{n\ \mathrm{even}}\gamma(S_n)/n$?\vs

 \makebox[25pt][r]{(ao)} What is the exact value of \ $\displaystyle\liminf_{n\ \mathrm{odd}}\gamma(S_n)/n$?
 \vs
 \makebox[15pt][r]{}{\it Comment of 2025}: It belongs to $[1/12,\,1/6]$ for $n$  odd (S.\,Eberhard, C.\,Mellon, {\it Bull. London Math. Soc.} (2025), \url{https://doi.org/10.1112/blms.70154}).
 \vs\vs

 \makebox[25pt][r]{(b)} What are the exact values of \ $\displaystyle\liminf_{n\ \mathrm{even}}\gamma(A_n)/n$ \ and \ $\displaystyle\liminf_{n\ \mathrm{odd}}\gamma(A_n)/n$?
 \vs

 \makebox[15pt][r]{}{\it Comment of 2025}: These belong to $[1/18,\,1/6]$ for $n$ even, and to $[1/6,\,4/15]$ for $n$ odd (S.\,Eberhard, C.\,Mellon, {\it Bull. London Math. Soc.} (2025), \url{https://doi.org/10.1112/blms.70154}).

\vs\vs

 \makebox[25pt][r]{\zv (c)} What is the exact value of $\displaystyle\limsup_{n\ \mathrm{even}}\gamma(S_n)/n$?
 \vs

\makebox[15pt][r]{}It is known that $\displaystyle\limsup_{n\ \mathrm{odd}}\gamma(S_n)/n=1/2$.
\vs\vs

 \makebox[25pt][r]{\zv (d)} What are the exact values of \ $\displaystyle\limsup_{n\ \mathrm{odd}}\gamma(A_n)/n$ \ and \ $\displaystyle\limsup_{n\ \mathrm{even}}\gamma(A_n)/n$?
 \vs\vs
 \hf\raisebox{-1ex}{\sl D.\,Bubboloni, C.\,E.\,Praeger, P.\,Spiga}

\ul

\otv (ae) It is $1/6$ (S.\,Eberhard, C.\,Mellon, {\it Bull. London Math. Soc.} (2025), \url{https://doi.org/10.1112/blms.70154}).

\otv (c) It is 1/4 (S.\,Eberhard, C.\,Mellon, {\it Bull. London Math. Soc.} (2025), \url{https://doi.org/10.1112/blms.70154}).

\otv (d) These are 1/3 for the odd case, and 1/4 for the even case (S.\,Eberhard, C.\,Mellon, {\it Bull. London Math. Soc.} (2025), \url{https://doi.org/10.1112/blms.70154}).
\emp

\bmp \textbf{20.18.}
 Let $R_{p^k}$ be the variety of class $2$ nilpotent groups of exponent
$p^k$, where $p$ is a prime number and $k\geq 2$. It is true that for every $p$ and $k$ there are infinitely many subquasivarieties of $R_{p^k}$ each of which is generated by a finite group with derived subgroup of exponent $p^k$ and does not have an independent basis of quasi-identities?

 \hf
\raisebox{-1ex}{\sl A.\,I.\,Budkin}
\emp

\bmp \textbf{20.19.}
A subgroup $H$ of a group $G$ is called \textit{commensurated} if for all $g\in G$, the index $|H : H \cap gHg^{-1}|$ is finite. Can a non-abelian free group (or a non-elementary hyperbolic group) contain two infinite commensurated subgroups $A$, $B$ with a trivial intersection? The answer is negative if $A$ or $B$ is normal.
\hf
\raisebox{-1ex}{\sl P.-E.\,Caprace}
\emp

 \bmp \textbf{20.20.} Let $E(l,m,n)$ be a group on generators $a,b$ with the presentation
 $\langle a,b\mid a^l=\nobreak 1,\;(ab)^m=b^n\rangle$, where $l,m,n$ satisfy $1/l + 1/m + 1/n < 1$.
Are there any values for $l,m,n$ such that $E(l,m,n)$ is shortlex automatic?

\makebox[15pt][r]{}It is known that $E(6,2,6)$ is not shortlex automatic.
\hf \raisebox{0ex}{\sl C.\,Chalk}
\emp

 \bmp \textbf{20.21.} (G.\,Verret).
 Does there exist a finite group $G$ with two normal subgroups $K$ and $L$, each with index 12 in $G$,
 such that $K$ is isomorphic to $L$, but $G/K$ is isomorphic to $C_{12}$, while $G/L$ is isomorphic to $A_4$?
 \hf\raisebox{-1ex}{\sl M.\,Conder}
\emp

\bmp \textbf{\zv 20.22.}
Let $G$ be a non-free and non-cyclic one-relator group in which every subgroup of infinite index is free. Is $G$ a surface group?
 \hf\raisebox{-1ex}{\sl B.\,Fine, G.\,Rosenberger, L.\,Wienke}

 \ul

 \otv Yes, it is (H.\,Wilton, {\it Preprint}, 2024, \url{https://arxiv.org/abs/2406.02121}).
\emp

\bmp \textbf{\zv 20.23.}
Is there a constant $\delta$ with $0 < \delta < 1$ and an integer $N$
such that whenever $A$, $B$, and $C$ are conjugacy classes in the
alternating group $Alt(n)$ each of size at least $|Alt(n)|^{\delta}$
with $n \geq N$, then $ABC = Alt(n)$?
\hf \raisebox{-1ex}{\sl M.\,Garonzi, A.\,Mar\'oti}

\ul

\otv Yes, there are (D.\,Dona, {\it Preprint}, 2025, \url{https://arxiv.org/abs/2505.06012}).
\emp

\bmp \textbf{20.24.}
Does the group $SL_2(\Z [\sqrt 2])$ admit a faithful transitive
amenable action?

\makebox[15pt][r]{}This is the same as asking if it admits any co-amenable subgroup except
the finite index subgroups. \hf\raisebox{-1ex}{\sl Y.\,Glasner, N.\,Monod}
\emp

\bmp \textbf{20.25.}
Let $n=k^2$. Let $\sigma\in S_n$ be the product of $k$ disjoint cycles, of lengths
$1,3,5,\dots , 2k-1$. Find the limiting proportion of odd entries in the
$\sigma$-column of the character table of $S_n$.

\makebox[15pt][r]{}For some background, see pages 1005--1006 in (D.\,Gluck, \emph{Proc. Amer. Math. Soc.}, \textbf{147} (2019), 1005--1011).
\hf
\raisebox{-1ex}{\sl D.\,Gluck}
\emp

\bmp \textbf{20.26.}
 Let $\phi$ be the Euler totient function. Does there exist a constant $a > 0$ such that $| Aut(G)| \geq \phi (|G|)^a$

 \makebox[25pt][r]{a)} for every finite group~$G$?

\makebox[25pt][r]{b)} for every finite nilpotent group $G$?

\makebox[15pt][r]{}If such a constant $a$ exists, then $a\le \frac{40}{41}$ (J.\,Gonz\'alez-S\'anchez, A.\,Jaikin-Zapirain, {\it Forum Math. Sigma}, {\bf 3}, Article ID e7, 11 p., electronic only (2015)). Also cf. Archive 12.77.
\hf \raisebox{-1ex}{\sl J.\,Gonz\'alez-S\'anchez, A.\,Jaikin-Zapirain}
\emp

\bmp \textbf{\zv 20.27.}
Let $G$ be a finite group, $p$ a prime number, and let $|x^G|_p$ denote the maximum power of $p$ that divides the class size of an element $x\in G$. Suppose that there exists a $p$-element $g\in G$ such that $|g^G|_p=\max_{x\in G}|x^G|_p$. Is it true that $G$ has a normal $p$-complement?

\makebox[15pt][r]{}A partial answer is in (\url{https://arxiv.org/abs/1812.03641}).
\hf
\raisebox{0ex}{\sl I.\,B.\,Gorshkov}

\ul

\otv No, not necessarily, a counterexample is given by \verb#SmallGroup(192,945)#
 (B.\,Sam\-bale, {\it Letter of 16 February 2022}).
\emp

\bmp \textbf{20.28.}
Let $L$ be a non-abelian finite simple group, and let $H(L)=M(L).L$ be the universal perfect central extension, where $M(L)$ is the Schur multiplier of $L$. Suppose that $G$ is a finite group such that the set of class sizes of $G$ is the same as the set of class sizes of $H(L)$. Is it true that $G\simeq H(L)\times A$, where $A$ is an abelian group?

\makebox[15pt][r]{}This is proved for $L=Alt_5$ in ({\it J.~Algebra Appl.}, {\bf 21}, no.\,11 (2022), Article ID 2250226, 8 p.).
\hf
\raisebox{-1ex}{\sl I.\,B.\,Gorshkov}
\emp

\bmp \textbf{20.29.}
Let $S$ be a non-abelian finite simple group. Is it true that for any $n\in \mathbb{N}$, if the set of class sizes of a centreless finite group $G$ is the same as the set of class sizes of the direct power $S^n$, then $G\simeq S^n$?

\makebox[15pt][r]{}This is proved for $S=Alt_5^2$ in ({\it J.~Algebra Appl.}, {\bf 21}, no.\,11 (2022), Article ID 2250226, 8 p.).
\hf
\raisebox{0ex}{\sl I.\,B.\,Gorshkov}
\emp

\bmp \textbf{\zv 20.30.}
(P.\,M.\,Neumann and M.\,R.\,Vaughan-Lee).
Let $G$ be a perfect and centreless finite group, and let $n$ be the maximum size of a conjugate class in $G$. Is it true that $|G|\leq n^2$?
 \hf\raisebox{-1ex}{\sl I.\,B.\,Gorshkov}

\ul

\otv No, not always (K.\,Muliarchyk, {\it Preprint}, 2026, \url{https://kourovkanotebookorg.wordpress.com/wp-content/uploads/2026/09/kourovka_20_30_solution.pdf}).
\emp

\bmp \textbf{20.31.} 
Let $\omega(G)$ denote the set of element orders of a finite group $G$,
and $h(G)$ the number of pairwise nonisomorphic finite groups $H$ with $\omega(H)=\omega(G)$. Find $h(L)$, where

\makebox[25pt][r]{\zva (a)\;} $L$ is the symmetric group of degree $10$;

\makebox[25pt][r]{(b)\;} $L$ is the automorphism group of the simple sporadic
 Janko group $J_2$;

\makebox[25pt][r]{(c)\;} $PSL(2,q)<L\leq \operatorname{Aut}(PSL(2,q))$.

See the current status in (M.\,A.\,Grechkoseeva, V.\,D.\,Mazurov, W.\,J.\,Shi, A.\,V.\,Vasil'ev, N.\,Yang, {\it Commun. Math. Stat.}, {\bf 11}, no.\,2 (2023), 169--194; Subsection 4.2).

 \hf\raisebox{-1ex}{\sl M.\,A.\,Grechkoseeva, A.\,V.\,Vasil'ev}

 \ul

 \otv (a) $h(S_{10})=\infty$ (D.\,V.\,Lytkina, V.\,D.\,Mazurov, {\it Preprint}, 2026, \url{https://arxiv.org/abs/2611.25766}; \ \ \ F.\,Scott, {\it PhD Thesis}, Univ. of  New South Wales,  2026,
 
 \url{https://kourovkanotebookorg.wordpress.com/wp-content/uploads/2026/11/20-31a.pdf}).
\emp

\bmp \textbf{20.32.} (E.\,Rapaport Strasser). The Lov\'asz conjecture states that a vertex-transitive connected graph is Hamiltonian. In the special case of Cayley graphs, one can ask the following. Consider a set $A$ which generates a group $G$ and is symmetric ($x\in A$ implies $x^{-1} \in A$). Is there a list $a_1,a_2,\dots, a_n$ of elements of $A$ such that $a_1, a_1a_2, a_1a_2a_3,\dots, a_1a_2\cdots a_n$ is a complete list of the elements of $G$? \hf
\raisebox{-1ex}{\sl B.\,Green}
\emp

 \bmp \textbf{\zv 20.33.} (A.\,Bauer).
Does Higman's Embedding Theorem relativize in the following way? Is it the case that for every subset $X \subseteq \mathbb{N}$, there is a finitely generated group $G_X$ that has an $X$-computable presentation (that is, there is a finite generating set relative to which the set of relations is computably enumerable with an $X$ oracle), and such that any finitely generated group has an $X$-computable presentation if and only if it can be embedded as a finitely generated subgroup of a quotient of a free product of finitely many copies of $G_X$ by the normal closure of a finite subset?

\makebox[15pt][r]{}Higman's Embedding Theorem says that for computable $X$, one may take $G_X=\mathbb{Z}$.

\hf \raisebox{-1ex}{\sl J.\,Grochow}

\ul

\otv Yes, it does (AI solution confirmed by E.\,Ingrosso:  {\it Preprint}, 2026,

\url{https://kourovkanotebookorg.wordpress.com/wp-content/uploads/2026/10/20.33.pdf}).
\emp

\bmp \textbf{20.34.}
Let $f$ be an inner screen of a saturated formation $\mathfrak{F}$ and suppose that a finite group $A$ acts faithfully and $f$-hypercentrally on a finite group $G$. Is it true that $G\rtimes A\in \mathfrak{F}$ for all $G\in \mathfrak{F}$ if and only if $\mathfrak{F}$ is a Fitting class?
 \hf\raisebox{-1ex}{\sl W.\,Guo}
\emp

\bmp \textbf{20.35.} (Well-known problem).
Let $S$ be a connected orientable hyperbolic surface of finite type and complexity at least 2. Does its mapping class group ${MCG}(S)$ have a non-elementary hyperbolic quotient?

\makebox[15pt][r]{}An affirmative answer is known when $S$ is a closed genus-2 surface (\url{https://arxiv.org/pdf/2005.00567.pdf}), and in the same paper it is shown that the answer will be affirmative provided certain hyperbolic groups are residually finite. See also a closely related question 16.108 about braid groups.
\hf
\raisebox{-1ex}{\sl M.\,Hagen}
\emp

\bmp \textbf{20.36.}
The free Burnside group of exponent four on three generators,
$B(3,4)$, has order $2^{69}$ as shown by Bayes, Kautsky, and Wamsley (1974).
Their proof is based on a theorem of Sanov (1940) which shows
that $B(n,4)$ is finite. Sanov's proof for $B(3,4)$ uses more than $2^{32}$ fourth powers,
because the subgroup of $B(3,4)$ generated by two of its generators and the
square of the third has order $2^{32}$. It is also known that $B(3,4)$ needs
at least 105 relations to define it, as shown by Havas and Newman (1980).

\makebox[25pt][r]{a)} Can $B(3,4)$ be defined with fewer than a million fourth powers?

\makebox[25pt][r]{b)} Can $B(3,4)$ be defined with fewer than a thousand fourth powers?

\makebox[25pt][r]{c)} What is the smallest number of fourth powers which define $B(3,4)$?

\hf\raisebox{-1ex}{\sl G.\,Havas, M.\,F.\,Newman}
\emp

\bmp \textbf{20.37.}
Is it true that for every finite group $G$ and every factorization $|G|=ab$ there exist subsets $A,B\subseteq G$ with $|A|=a$ and $|B|=b$ such that $G=AB$? Cf. 19.35 in Archive.

\makebox[15pt][r]{}{\it Comment of 2025}:
any minimal counterexample $G$ must be a simple group without prime-index subgroups (R.\,R.\,Bildanov, V.\,A.\,Goryachenko, A.\,V.\,Vasil'ev, {\it Siberian Electron. Math. Rep.}, \textbf{17} (2020), 683--689;
\ M.\,H.\,Hooshmand, \emph{Commun. Algebra}, \textbf{49}, no.\,7 (2021), 2927--2933).
\hf
\raisebox{-1ex}{\sl M.\,H.\,Hooshmand}
\emp

\bmp \textbf{20.38.}
Suppose that $H$ is an almost simple group of Lie type and $G$ is a finite group such that
$G$ and $H$ have the same sets of degrees of irreducible complex characters. Must there exist an abelian normal subgroup $A$ of $G$ such that $G/A$ is isomorphic to~$H$?

\makebox[15pt][r]{}{\it Comment of 2025}: The conjecture was confirmed for projective general linear and unitary groups of dimension 3 (F.\,Shirjian, A.\,Iranmanesh, {\it Illinois J. Math.}, {\bf  64}, no.\,1 (2020), 49--69).
\hf\raisebox{-1ex}{\sl A.\,Iranmanesh}
\emp

\bmp \textbf{20.39.}
Let $F$ be a non-abelian finitely generated free group, $1\ne w\in F$, and $n\ge 1$. Is the group $\langle F,t \mid t^n=w \rangle $ linear of degree 2 over a field of characteristic 0 if $w$ is not a proper power in~$F$?

\makebox[15pt][r]{}This question is motivated by the following well-known question: is it true that the free $\Q$-group $F^\Q$ is linear over a field? (See 13.39(b).)
\hf \raisebox{-1ex}{\sl A.\,Jaikin-Zapirain}
\emp

\bmp \textbf{\zv 20.40.}
Let $G$ be a $\kappa$-existentially closed group of cardinality $\lambda > \kappa$, where $\kappa$ is a regular cardinal. Is it true that $|Aut(G)|=2^{\lambda}$?

 \makebox[15pt][r]{}The answer is known to be affirmative if $\lambda=\kappa$ (B.\,Kaya, M.\,Kuzucuo\u{g}lu, P.\,Longobardi, M.\,Maj,
 {\it J.~Algebra}, {\bf 666} (2025), 840--849).
\hf\raisebox{-1ex}{\sl B.\,Kaya, M.\,Kuzucuo\u{g}lu}

\ul

\otv No, not necessarily (S.\,Corson, S.\,Shelah, {\it Preprint}, 2026, \url{https://arxiv.org/abs/2610.39700}).
\emp

\bmp \textbf{\zv 20.41.}
(L.\,Ciobanu, B.\,Fine, G.\,Rosenberger). Suppose that $G$ is a non-cyclic residually finite group in which every subgroup of finite index (including the group itself) is defined by a single defining relation, while all infinite index subgroups are free. Is it true that $G$ is either free or isomorphic to the fundamental group of a compact surface?

\makebox[15pt][r]{}See Archive 7.36 for a negative solution of a similar question without the assumption on infinite index subgroups. \hfill\raisebox{-1ex}{\sl D.\,Kielak}

 \ul

 \otv Yes, it is true (H.\,Wilton, {\it Preprint}, 2024, \url{https://arxiv.org/abs/2406.02121}).
\emp

\bmp \textbf{20.42.} (Y.\,Cornulier, A.\,Mann, A.\,Thom).
Does there exist a finitely generated residually finite group that is not (elementary) amenable and satisfies a nontrivial group law?
\hf\raisebox{-1ex}{\sl S.\,Kionke}
\emp

\bmp \textbf{20.43.}
Is every family of finite groups satisfying a common group law uniformly amenable?
\hf\raisebox{-1ex}{\sl S.\,Kionke, E.\,Schesler}
\emp

\bmp \textbf{20.44.}
The definition of \({\rm CT}(\mathbb{Z})\) is given in 17.57. Is it true that a finitely generated subgroup of \({\rm CT}(\mathbb{Z})\) either has
only finitely many orbits on \(\mathbb{Z}\) or there is a set of representatives for its orbits
on \(\mathbb{Z}\) which has positive density?
 \hf\raisebox{-1ex}{\sl S.\,Kohl}
\emp

\bmp \textbf{20.45.}
Let \(n\) be a positive integer, and let \(G \leq {\rm GL}(n,\mathbb{Z})\)
be finitely generated. Given a bound \(b \in \mathbb{N}\), let \(e_b\) be
the number of elements of \(G\) all of whose matrix entries have absolute value
\(\leq b\). Does the limit \ \(\displaystyle\lim_{b \rightarrow \infty} \ln{e_b}/\ln{b}\) \ always exist?
 \hf\raisebox{-1ex}{\sl S.\,Kohl}
\emp

\bmp \textbf{20.46.}
A group action on a compact space is said to be \emph{topologically free} if the set of points with trivial stabilizer is dense. Let $G$ be a locally compact group, and $\partial_{sp} G$ its Furstenberg boundary (the largest minimal and strongly proximal compact $G$-space). Let $\Gamma_1$ and $\Gamma_2$ be two lattices in $G$, both acting faithfully on $\partial_{sp} G$.

\makebox[25pt][r]{a)} Is it possible that the $\Gamma_1$-action on $\partial_{sp} G$ is topologically free, but the $\Gamma_2$-action on $\partial_{sp} G$ is not topologically free?

\makebox[25pt][r]{b)} If yes, can this also happen if $\partial_{sp} G = G/H$ is a homogeneous $G$-space?

 \hf\raisebox{-1ex}{\sl A.\,Le Boudec}
\emp

\bmp \textbf{20.47.}
Let $G = F_n$ be a finitely generated free group. Let $X$ be a compact $G$-space on which the $G$-action is faithful, minimal, and strongly proximal. Does it follow that the action is topologically free?
 \hf\raisebox{-1ex}{\sl A.\,Le Boudec, N.\,Matte Bon}
\emp

\bmp \textbf{20.48.}
Suppose that $P$ is a finite $p$-group with a non-trivial partition (which is equivalent to having proper Hughes subgroup $H_p(P):=\langle g\in P\mid g^p\ne 1\rangle\ne P$). If $P$ admits a fixed-point-free automorphism of prime order, must $P$ be of exponent~$p$?

\hf \raisebox{-1ex}{\sl M.\,Lewis}
\emp

\bmp \textbf{20.49.}
Is it true that any finite group contains a 2-generated subgroup with the
same exponent?

\makebox[15pt][r]{}An affirmative answer is known for soluble groups. It is also known that any finite group contains a 3-generated subgroup with the same exponent (E.\,Detomi, A.\,Lucchini, \emph{J.~London Math. Soc. (2)}, \textbf{87}, no.\,3 (2013), 689--706).
\hf
\raisebox{-1ex}{\sl A.\,Lucchini}
\emp

\bmp \textbf{20.50.}
For a positive integer $m$, let $G_m$ be the largest group generated by $m$ involutions such that $(xy)^4=1$ for any two involutions $x,y\in G_m$. What is the order of $G_4$?

 \makebox[15pt][r]{}It is known that the order of $G_3$ is equal to $2^{11}$. Also cf. 18.58.
\hf\raisebox{0ex}{\sl D.\,V.\,Lytkina}
\emp

\bmp \textbf{20.51.}
Is it true that for every odd integer $t > 3$ there exists a finite non-abelian group $G$ of odd order with exactly $t$ conjugacy classes?
\hf
\raisebox{-1ex}{\sl D.\,MacHale}
\emp

\bmp \textbf{20.52.} 
A famous result of Burnside states that if $k(G)$ is the number of conjugacy classes of a finite group $G$ of odd order, then $|G|- k(G)$ is divisible by 16. Is it true that for every integer $m > 0$ there exists a finite non-abelian group $G(m)$ of odd order such that $|G(m)| - k(G(m)) = 16m$?

\makebox[15pt][r]{}{\it Comment of 2026}: The existence of such a group $G(m)$ is shown for every odd $m$ and for every $m$ in which $2$ occurs to an odd power  (M.\, Mart\'{i} Salvador, {\it Preprint}, 2026, \url{https://palomar-registry.org/entry.html?id=PALOMAR-2026-11-27-000085&version=1}; \url{https://github.com/marcosmartinsalvador/conjecture-proofs}).
\hf
\raisebox{-1ex}{\sl D.\,MacHale}
\emp

\bmp \textbf{20.53.}
A nonisotropic unitary graph $\Gamma$ is distance-regular with intersection array $\{q(q-1),\,(q+1)(q-2),\,q+1;\,1,\,1,\,q(q-2)\}$ for some prime power $q$. The group $G=Aut (\Gamma)$ acts transitively on the vertex set and on the edge set of $\Gamma$.
It is known that $\Gamma$ is distance-transitive if $q=3$. Does there exist a distance-regular graph with such an intersection array if $q$ is not a prime power?
 \hf\raisebox{-1ex}{\sl A.\,A.\,Makhnev}
\emp

\bmp \textbf{20.54.}
A distance-regular graph of diameter 3 with the second eigenvalue $\theta_1=a_3$ is called a \emph{Shilla graph}.
For a Shilla graph $\Gamma$ the number $a=a_3$ divides $k$ and we set $b=b(\Gamma)=k/a$. Koolen and Park proved that there are 12 feasible intersection arrays of Shilla graphs with $b=3$. At present it is proved that a Shilla graph with $b=3$ has intersection array $\{12,10,3;1,3,8\}$ (Doro graph), $\{12,10,5;1,1,8\}$ (nonisotropic unitary graph for $q=4$), or $\{15,12,6;1,2,10\}$. The automorphisms of the last graph were found by A.\,Makhnev and N.\,Zyulyarkina (\emph{Doklady Maths.}, {\bf 84}, no.\,1 (2011), 510--514). Does the graph with intersection array $\{15,12,6;1,2,10\}$ exist?

 \hf\raisebox{-1ex}{\sl A.\,A.\,Makhnev, D.\,O.\,Revin}
\emp

\bmp \textbf{20.55.}
Are there soluble finite groups $G$ and $H$, of derived lengths 2 and 4 and having identical character tables?

\makebox[15pt][r]{}Pairs of nonisomorphic soluble finite groups with identical character tables and with derived
lengths $n$ and $n+1$ for any $n\ge 2$ were constructed in (S.~Mattarei, {\it J.~Algebra}, {\bf 175} (1995) 157--178).
\hf \raisebox{-1ex}{\sl S.\,Mattarei}
\emp

\bmp \textbf{\zv 20.56.}
Is a periodic group $G$ locally finite if it is generated by involutions and the centralizer of every involution in $G$ is locally finite?
 \hf\raisebox{-1ex}{\sl V.\,D.\,Mazurov}

 \ul

\otv No, not always  (AI solution confirmed by V.\,D.\,Mazurov: \  B.\,G.\,Di~Siena, {\it Preprint}, 2026, \url{https://kourovkanotebookorg.wordpress.com/wp-content/uploads/2026/10/20-56.pdf}).
\emp

\bmp \textbf{\zv 20.57.}
Is a periodic group $G$ locally finite if every finite subgroup of $G$ is contained in a subgroup of $G$ isomorphic to a finite alternating group?
 \hf\raisebox{-1ex}{\sl V.\,D.\,Mazurov}

 \ul

\otv No, not always  (AI solution confirmed by V.\,D.\,Mazurov: \  B.\,G.\,Di~Siena, {\it Preprint}, 2026, \url{https://kourovkanotebookorg.wordpress.com/wp-content/uploads/2026/10/20.57.pdf}).
\emp

\bmp \textbf{20.58.}
Let $\omega(G)$ denote the set of element orders of a finite group $G$.
A finite group~$G$ is said to be \textit{recognizable} (\textit{by spectrum})
if every finite group $H$ with $\omega(H)=\omega(G)$ is isomorphic to $G$.

\makebox[25pt][r]{\zva a)}
Is it true that for every $n$ there is a recognizable group that is the $n$-th direct power of a~nonabelian simple group?

\makebox[25pt][r]{b)} Is it true that there is a nonabelian simple group $L$ such that for every $n$ there is a recognizable group whose socle is the $k$-th direct power of $L$ for some $k\geqslant n$?

 \hf\raisebox{-1ex}{\sl V.\,D.\,Mazurov, A.\,V.\,Vasil'ev}

 \ul

 \otv a) Yes, it is true (N.\,Yang, I.\,Gorshkov, A.\,Staroletov, A.\,V.\,Vasil'ev, {\it Annali Matem. Pura Appl.}, {\bf 202} (2023), 2699--2714). \emp

\bmp \textbf{20.59.}
 A subgroup $H$ is called a \emph{virtual retract} of a group $G$ if $H$ is a retract of a finite-index subgroup of $G$. Is it true that every finitely generated subgroup of a finitely generated virtually free group is a virtual retract?

 \makebox[15pt][r]{}For motivation and partial results see (A. Minasyan, {\it Int. Math. Res. Notes}, \textbf{2021}, no.\,17 (2021), 13434--13477).
 \hf \raisebox{-1ex}{\sl A.\,Minasyan}
\emp

\bmp \textbf{20.60.} We say that $G$ is a \textit{virtually compact special group} if $G$ has a finite-index subgroup which is isomorphic to the fundamental group of a compact special complex (in the sense of F.\,Haglund, D.\,T.\,Wise, {\it Geom. Funct. Anal.}, \textbf{17}, no.\,5 (2008), 1551--1620).
Let $G$ be a virtual retract of a finitely generated right-angled Artin group. Must $G$ be a virtually compact special group?

 \makebox[15pt][r]{}An affirmative answer would provide an algebraic characterization of the class of virtually compact special groups as the class groups admitting finite index subgroups that are virtual retracts of right-angled Artin groups.
 \hf \raisebox{-1ex}{\sl A.\,Minasyan}
\emp

\bmp \textbf{20.61.} Suppose that $G$ is a virtually compact special group. Is it true that the centralizer of any element in $G$ is itself virtually compact special?
\hf\raisebox{-1ex}{\sl A.\,Minasyan}
\emp

\bmp \textbf{\ai 20.62.} 
(J.\,Dixmier)
A group is said to be \emph{unitarisable} if its every uniformly bounded representation on a Hilbert space is
unitarisable. Is every unitarisable
group amenable?
 \hf\raisebox{-1ex}{\sl N.\,Monod}

\ul

\oai Yes, it is (unverified AI solution: OpenAI, {\it Preprint}, 2026,

\url{https://kourovkanotebookorg.wordpress.com/wp-content/uploads/2026/10/20-62.pdf}).
\emp 

\bmp \textbf{20.63.}
a) Prove that every unitarisable group has trivial cost.

\makebox[15pt][r]{}This is known to be true for residually finite groups (I.\,Epstein, N.\,Monod, \emph{Int. Math. Res. Notes}, \textbf{2009}, no.\,22 (2009), 4336--4353).

\makebox[15pt][r]{}b) At least prove that every unitarisable group has vanishing
first $L^2$ Betti number.

\hf\raisebox{-1ex}{\sl N.\,Monod}
\emp

\bmp \textbf{20.64.}
Let $G$ be a non-amenable group.

\makebox[25pt][r]{a)} Prove that $G^n$ is non-unitarisable for some $n$.

\makebox[25pt][r]{b)} At least prove that $G^{\infty}$, the
direct sum (restricted product), is non-unitarisable.

\hf\raisebox{-1ex}{\sl N.\,Monod}
\emp

 \bmp \textbf{20.65.}
Given a complex irreducible character $\chi\in \rm{Irr}(G)$ of a finite group $G$, let $\mathbb{Q}(\chi)$ denote the field extension of $\mathbb{Q}$ obtained by adjoining to $\mathbb{Q}$ all the values of $\chi$. We say that a finite group $G$ is $k$-\textit{rational} if $|\mathbb{Q}(\chi):\mathbb{Q}|$ divides $k$ for every $\chi\in \rm{Irr}(G)$. Does there exist a real-valued function $f$ such that if $p$ is the order of a cyclic composition factor of a $k$-rational group $G$, then $p\leq f(k)$?

\makebox[15pt][r]{}If $k=1$, then we know that $p\leq 11$ by a theorem of J.\,G.\,Thompson ({\it J.~Algebra}, {\bf 319} (2008), 558--594).
\hf \raisebox{-1ex}{\sl A.\,Moret\'o}
\emp

\bmp \textbf{\ai 20.66.} 
A Schmidt $(p, q)$-group is a finite non-nilpotent group all of whose proper subgroups are nilpotent and whose Sylow $p$-subgroup is normal.
The $N$-critical graph $\Gamma_{Nc}(G)$ of a finite group $G$ is a directed graph on the vertex set
of all prime divisors of $|G|$ in which $(p, q)$ is an edge of $\Gamma_{Nc}(G)$ if and only if $G$ has a Schmidt $(p, q)$-subgroup. Suppose that a finite group $ G$ is such that $G=AB=AC=BC$, where $ A, B, C$ are subgroups of $G$.
Is $\Gamma_{Nc}(G)=\Gamma_{Nc}(A)\cup\Gamma_{Nc}(B)\cup\Gamma_{Nc}(C)$?
This is true if $A, B, C$ are soluble.
 \hf\raisebox{-1ex}{\sl V.\,I.\,Murashka, A.\,F.\,Vasil'ev}

\ul

\oai No, not always (unverified AI solution: B.\,G.\,Di Siena, {\it Preprint}, 2026, \url{https://www.proofforum.org/abs/PF:2026.000034}).
\emp

\bmp \textbf{20.67.} Suppose that $G$ is a finite group, $p$ is
a prime, and $B$ is a Brauer $p$-block of $G$
with defect group $D$. Let ${\rm cd}(B)$ be the
set of degrees of the irreducible complex characters in $B$.

\makebox[15pt][r]{}a) Is it true that the derived length of $D$ is bounded by $|{\rm cd}(B)|$?

\makebox[15pt][r]{}b) (A.\,Jaikin-Zapirain). Is is even true that $|{\rm cd}(D)|\le |{\rm cd}(B)|$?
\hf\raisebox{0ex}{\sl G.\,Navarro}
\emp

\bmp \textbf{20.68.} The \emph{submonoid membership problem} for a group $G$ generated by an alphabet~$A$ asks, for given words $x_1, x_2, \dots , x_n$ and a word $w$ over $A$, whether $w$ belongs to the submonoid of $G$ generated by the~$x_i$. Does there exist a hyperbolic one-relator group with undecidable submonoid membership problem?
\hf\raisebox{-1ex}{\sl C.-F.\,Nyberg-Brodda}
\emp

\bmp \textbf{20.69.} Is the submonoid membership problem decidable for every one-relator group with torsion?
\hf\raisebox{-1ex}{\sl C.-F.\,Nyberg-Brodda}
\emp

\bmp \textbf{20.70.}
As a strengthening of the Burnside restriction, for every pair $(k,n)$ of positive integers, let a group $G$ satisfy condition $C_{k,n}$ if every $k$-generated subgroup of $G$ is finite of order at most~$n$.

\makebox[25pt][r]{a)} Does there exist $k\ge 2$ such that for any $n$ all groups with condition $C_{k,n}$ are locally finite?

\makebox[25pt][r]{b)} In particular, is it true that for any $n$ the condition $C_{2,n}$ implies local finiteness?

\makebox[25pt][r]{c)} Find possibly more pairs $(k,n)$ for which groups with condition $C_{k,n}$ are locally finite. (For example, all groups with condition $C_{2,20}$ are metabelian, and therefore locally finite.)
\hf
\raisebox{-1ex}{\sl A.\,Yu.\,Olshanskii}
\emp

\bmp \textbf{20.71.} (J.\,Lauri).
A \emph{card} of a finite simple undirected graph $G$ of order $n=|V(G)|$ is an induced subgraph of order $n-1$. For a connected graph $G$ with $k$ isomorphism types of cards, can $Aut(G)$ have more than $k$ orbits on the vertex set $V(G)$

\makebox[25pt][r]{a)}   for $k=2$?

\makebox[25pt][r]{b)}  for $k=3$?

\makebox[25pt][r]{\zv c)}
  for $k=4$?

\makebox[15pt][r]{}If a graph $G$ has $k$ isomorphism types of cards, then the group $Aut(G)$ of automorphisms of $G$ has obviously at least $k$ orbits on $V(G)$. It is known that if all cards are mutually isomorphic, then $Aut(G)$ is transitive on $V(G)$. Examples are known of graphs $G$ with 5 isomorphism types of cards for which $Aut(G)$ has $6$ orbits on $V(G)$, and of graphs with 6 isomorphism types for which $Aut(G)$ has 7 orbits on $V(G)$.

\hf\raisebox{-1ex}{\sl V.\,Pannone}

\ul

\otv c) Yes, it can for $k=4$  (J.\,Prajapati, {\it Preprint}, 2026,

\url{https://kourovkanotebookorg.wordpress.com/wp-content/uploads/2026/10/kourovka-20-71-prajapati.pdf}).
\emp

\bmp \textbf{\zv 20.72.}
Is there a variety of groups $\Theta$ such that the group of automorphisms of the \linebreak category of free finitely generated groups $\Theta^0$ contains outer automorphisms?

\hf\raisebox{-1ex}{\sl E.\,Plotkin}

\ul

\otv Yes, there is (R.\,Barbosa-Fernandes, A.\,Tsurkov, {\it Preprint}, 2019,  \url{https://arxiv.org/abs/1909.05955}).
\emp

 \bmp \textbf{\zv 20.73.}
Can every countable group $G$ be factorized $G=AB$ into infinite subsets $A, B$ such that every
element $g\in G$ has a unique representation $g = ab$ for $a\in A$, $b \in B$?

\makebox[15pt][r]{}This is true if $G$ is topologizable. \hf
\raisebox{-0ex}{\sl I.\,V.\,Protasov}

\ul

\otv Yes, it can (AI solution confirmed by I.\,V.\,Protasov: \  B.\,G.\,Di~Siena, {\it Preprint}, 2026, \url{https://kourovkanotebookorg.wordpress.com/wp-content/uploads/2026/10/20.73.pdf}).
\emp 

\bmp \textbf{\zv 20.74.}
a) {\it Conjecture}: there are only finitely many nonabelian finite simple groups $G$ that have a normal subset $S$ closed under inversion such that $|S | > |G|/\log_2|G|$ and $S^2\ne G$. {\it Comment of 2024}: This is true for the class of alternating groups (M.\,Larsen, P.\,H.\,Tiep, {\it Preprint}, 2023, \url{arXiv:2305.04806}) and for any class of groups of Lie type of bounded Lie rank (S.\,V.\,Skresanov, {\it Bull. London Math. Soc.}, {\bf 57}, no.\,11 (2025), 3429--3440).

\makebox[25pt][r]{b)} The same conjecture for the groups $PSL(2,p)$.
\hf\raisebox{0ex}{\sl L.\,Pyber}

\ul

\otv a) No, there are infinitely many such groups (O.\,Sheinfeld,  {\it Preprint}, 2026,
\url{https://arxiv.org/abs/2609.22004}).

\otv b) The conjecture is proved for the groups $PSL(2,p)$ (S.\,V.\,Skresanov, {\it Bull. London Math. Soc.}, {\bf 57}, no.\,11 (2025), 3429--3440).
\emp

\bmp \textbf{\zv 20.75.}
 Let $G$ be a finite $p$-group and assume that all abelian normal
subgroups of $G$ can be generated by $k$ elements.
Is it true that every abelian subgroup of $G$ can be generated by $2k$ elements?

\makebox[15pt][r]{}The $k$-th direct power of $D_{16}$ shows that this bound would be best
possible.

\hf\raisebox{-1ex}{\sl L.\,Pyber}

\ul

\otv No, it is not (K.\,Muliarchyk,  {\it Preprint}, 2026, \url{https://kourovkanotebookorg.wordpress.com/wp-content/uploads/2026/09/kourovka_20.75_solution.pdf}).
 \emp

\bmp \textbf{\zv 20.76.}
Let $G$ be a finite $p$-group and assume that all abelian normal
subgroups of $G$ have order at most $p^k$.
Is it true that every abelian subgroup of $G$ has order at most $p^{2k}$?
 \hf\raisebox{-1ex}{\sl L.\,Pyber}

 \ul

\otv No, it is not (K.\,Muliarchyk,  {\it Preprint}, 2026, \url{https://kourovkanotebookorg.wordpress.com/wp-content/uploads/2026/09/kourovka_20_75-20_76_solution.pdf}).
 \emp

\bmp \textbf{20.77.} (O.\,I.\,Tavgen').
A group is said to be \emph{boundedly generated} if it is a
product of finitely many cyclic groups.
Is it true that every boundedly generated residually finite group is linear over a field
of characteristic zero?

\makebox[15pt][r]{}This is true for residually finite-soluble groups (L.\,Pyber, D.\,Segal, \emph{J.~Reine Angew. Math.}, \textbf{612} (2007), 173--211).
\hf\raisebox{-1ex}{\sl L.\,Pyber}
\emp

\bmp \textbf{20.78.}
 For an irreducible complex character $\chi$ of a finite group
$G$, the \emph{codegree} of $\chi$ is defined by ${\rm
cod}(\chi)=|G:\ker\chi|/\chi(1)$. Let ${\rm Cod}(G)$ be the set of
irreducible character codegrees of $G$.
\emph{Conjecture}: If $G$ has an
element of order $m$, then $m$ divides some member of ${\rm
Cod}(G)$.

 \makebox[15pt][r]{}The conjecture is proved when $m$ is a prime power (G.\,Qian, \emph{Arch. Math.}, \textbf{97} (2011), 99--103); when $m$ is square-free (I.\,M.\,Isaacs, \emph{Arch. Math.}, \textbf{97} (2011), 499--501); or when $G$ is solvable (G. Qian, \emph{Bull. London Math. Soc.}, \textbf{53} (2021) 820--824); or when $G$ is a symmetric or alternating group
 (E. Giannelli ({\it J.~Algebra Appl.}, {\bf 23}, no.\,9 (2024),
 article ID 2450144);
 or when $G$ is an almost simple group (S.\,Y.\,Madanha, {\it Commun. Algebra}, {\bf 51}, no.\,7 (2023), 3143--3151); or when $\mathbf{F}(G)=1$ (Z.\,Akhlaghi, E.\,Pacifici, L.\,Sanus, {\it J.~Algebra}, {\bf 644} (2024), 428--441).
\hf\raisebox{-1ex}{\sl G.\,Qian}
\emp

\bmp \textbf{20.79.}
\emph{Conjecture}: Suppose that $G$ is a non-abelian simple group and $H$ is a
finite group such that ${\rm Cod}(G)={\rm Cod}(H)$. Then $G\cong
H$.

 \makebox[15pt][r]{}{\it Comment of 2025}: The conjecture is proved for
$G$ of type $\mathrm{PSL}(2,q)$ (A.\,Bahri, Z.\,Akhlaghi, B.\,Khosravi, {\it
Bull. Austral. Math. Soc.}, {\bf 104}, no.\,2 (2021), 278--286); \
 $\mathrm{PSL}(3, q)$ or $\mathrm{PSU}(3, q)$ (Y.\,Liu, Y.\,Yang, {\it Results Math.}, {\bf 78}, no.\,1 (2023), article No.\,7); \ a sporadic simple group (M.\,Dolorfino, L.\,Martin, Z.\,Slonim, Y.\,Sun, Y.\,Yang, {\it Bull. Austral. Math. Soc.}, {\bf 109}, no.\,1 (2024), 57--66); \ an alternating group (M.\,Dolorfino, L.\,Martin, Z.\,Slonim, Y.\,Sun, Y.\,Yang, {\it Bull. Austral. Math. Soc.}, {\bf 110}, no.\,1 (2024), 115--120); \ $^2F_4(q^2)$ (Y.\,Yang, {\it J.~Group
Theory}, {\bf 27}, no.\,1 (2024), 141--155); \ $\mathrm{PSL}(n,q)$ or a simple exceptional group of Lie type (H.\,P.\,Tong Viet, {\it Math. Nachr.}, {\bf 298}, no.\,4 (2025), 1356--1369).
\hf\raisebox{-1ex}{\sl G.\,Qian}
\emp

\bmp \textbf{20.80.}
Let $G$ be a finite group, and $\chi$ an irreducible complex
character of $G$. We call $\chi$ a $\mathcal{P}$-character if
$\chi$ is a constituent of $(1_H)^G$ for some maximal subgroup $H$
of~$G$; and $\chi$ is said to be monomial if $\chi$ is induced by a
linear character of a subgroup of $G$.

\makebox[15pt][r]{}\emph{Conjecture}: $G$ is solvable
if and only if all its $\mathcal{P}$-characters are monomial.

 \makebox[15pt][r]{}The necessity part of the conjecture is true (G.\,Qian, Y.\,Yang,
\emph{Commun. Algebra}, \textbf{46}, no.\,1 (2018), 167--175).
\hf\raisebox{-1ex}{\sl G.\,Qian}
\emp

\bmp \textbf{20.81.}
Does there exist a class $\mathfrak{X}$ of finite groups satisfying the following conditions:

\makebox[25pt][r]{(1)} $\mathfrak{X}$ is closed with respect to taking subgroups and homomorphic images;

\makebox[25pt][r]{(2)} the product of two normal $\mathfrak{X}$-subgroups of an arbitrary group is always an $\mathfrak{X}$-group;

\makebox[25pt][r]{(3)} for every positive integer $n$, there exist a finite group and its conjugacy class $D$ such that any $n$ elements of $D$ generate an $\mathfrak{X}$-subgroup, whereas $\langle D\rangle\notin\mathfrak{X}$.

\makebox[15pt][r]{}It
is known that there are no classes $\mathfrak{X}$ satisfying (1)--(3) that are closed with respect to extensions (D.\,O.\,Revin, {\it Algebra i Analiz}, \textbf{37}, no.\,1 (2025), 141--176 (Russian)).
 \hf\raisebox{-1ex}{\sl D.\,O.\,Revin}
\emp

\bmp \textbf{20.82.}
Two groups are said to be {\it isospectral} if they have the same set of element
orders. Suppose that $G$ is a finite group such that every finite group isospectral to $G$ is isomorphic to $G$. Is it true that the quotient of $G$ by its socle is solvable?

 \hf\raisebox{-1ex}{\sl D.\,O.\,Revin}
\emp

\bmp \textbf{20.83.}
Let $\mathfrak{X}$ be a class of finite groups that is closed with respect to taking subgroups, homomorphic images, and extensions. A subgroup $H$ of a finite group $G$ is said to be \emph{$\mathfrak{X}$-submaximal} if there exists an embedding of $G$ into a group $G^*$ such that $G$ is subnormal in $G^*$ and $H$ coincides with the intersection of $G$ and an $\mathfrak{X}$-maximal subgroup of $G^*$. Suppose that all $\mathfrak{X}$-submaximal subgroups of a characteristic subgroup $N$ of a finite group $G$ are conjugate in $N$. Does it follow that $HN/N$ is an $\mathfrak{X}$-submaximal subgroup of $G/N$ for every $\mathfrak{X}$-submaximal subgroup $H$ of $G$?

\makebox[15pt][r]{}For normal subgroups $N$, this is not true even if $N$ is an $\mathfrak{X}$-group or if $N$ does not contain nontrivial $\mathfrak{X}$-subgroups. A positive answer is known in the case where $N$ coincides with the $\mathfrak{F}$-radical of $G$ for a Fitting class $\mathfrak{F}$.

 \hf\raisebox{-1ex}{\sl D.\,O.\,Revin, A.\,V.\,Zavarnitsine}
\emp

\bmp \textbf{20.84.}
Does there exist a nilpotent group of class 3 with a non-trivial 5-th dimension subgroup?
\hf\raisebox{-1ex}{\sl E.\,Rips}
\emp

\bmp \textbf{20.85.} Let $F$ be a free group. An element $\omega\in F$ is said to be \emph{primitive} if there is a minimal generating system of $F$ that contains $\omega$, \emph{almost primitive} if it is primitive in each finitely generated proper subgroup of $F$ containing $\omega$, \emph{tame almost primitive} if, whenever $\omega^\alpha$ is contained in a subgroup $H$ of $F$ with $\alpha\ge 1$ minimal, either $\omega^\alpha$ is primitive in $H$ or the index of $H$ in $F$ is just $\alpha$. In (B.\,Fine, A.\,Moldenhauer, G.\,Rosenberger, L.\,Wienke, \emph{Topics in Infinite Group Theory: Nielsen Methods, Covering Spaces, and Hyperbolic Groups}, De Gruyter, Berlin, 2021) it is shown that
 $u=[a_1,b_1][a_2,b_2]\cdots [a_g,b_g]$ is tame almost primitive in the free group on $a_1, b_1,\dots ,a_g,b_g$ with $g\ge 1$, and
 $v=c_1 ^2\cdots c_p ^2$ is tame almost primitive in the free group on $c_1,\dots,c_p$ with $p\ge 2$.

 \makebox[15pt][r]{}Are there tame almost primitive elements in free groups other than $u$, $v$, and
 their product $uv$ in the free group on $a_1, b_1,\dots ,a_g,b_g,c_1,\dots,c_p$?

 \hf\raisebox{-1ex}{\sl G.\,Rosenberger, L.\,Wienke}
\emp

\bmp \textbf{20.86.}
Suppose $G$ is a finite group and $p$ a prime such that the number $s_p(G)=1+kp$ of Sylow $p$-subgroups $P$ of $G$ is greater than~1. By a theorem of Frobenius the set $G_p$ of all $p$-elements in $G$ has cardinality $|G_p|=|P|\cdot f_p(G)$ for some positive integer $f_p(G)$, and one easily gets that
 $f_p(G)=1+\ell(p-1)\le k-{k-1\over p-1}$.
Usually, $\ell <k$ (but $\ell =k$ if $k<p$). Is always $\ell^p\ge k^{p-1}$?

 \makebox[15pt][r]{}Recently P.\,Gheri showed that $f_p(G)^p\ge s_p(G)^{p-1}$ if $G$ is $p$-solvable (\emph{Ann. Mat. Pura Appl. (4)}, {\bf 200} (2021), 1231--1243).
\hf\raisebox{-1ex}{\sl P.\,Schmid}
\emp

\bmp \textbf{\ai 20.87.}
What are the non-abelian composition factors of finite groups in which the order of every element is divisible by at most two primes?

 \makebox[15pt][r]{}For the case of one prime, see (\url{https://arxiv.org/pdf/2003.09445.pdf}).

 \hf\raisebox{-1ex}{\sl W.\,J.\,Shi}

 \ul

\oai These factors are described (unverified AI solution: AutoPasha, {\it Preprint}, 2026, \url{https://kourovkanotebookorg.wordpress.com/wp-content/uploads/2026/10/20-87.pdf}); necessary conditions are also given in (J.\,Li, X.\,Xu, Y.\,Xu, Y.\,Yang, {\it
Commun. Algebra}, {\bf  53}, no.\,6 (2025), 2204--2212).
\emp

\bmp \textbf{\zv 20.88.}
Let $u$ be an element of a free group $F_r$. Is it true that there is $v \in F_r$ (that depends on $u$) that cannot be a subword of any cyclically reduced word $\varphi(u)$, where $\varphi$ is an automorphism of $F_r$?
\hf\raisebox{-1ex}{\sl V.\,Shpilrain}

\ul

\otv Yes, it is true (L.\,Koch-Hyde, S.\,O'Connor, E.\,Olive, V.\,Shpilrain, {\it Preprint}, 2025, \url{https://arxiv.org/abs/2505.00477}).
\emp

\bmp \textbf{20.89.}
An element $g$ of a group $G$ is said to be \emph{almost Engel} if there is a finite subset $E(g)$ of $G$ such that for every $x\in G$ all sufficiently long commutators $[...[[x, g], g], \dots , g]$ belong to $E(g)$, that is,
there is a positive integer $n(x, g)$ such that $[...[[x, g], g], \dots ,g]\in E(g)$ if $g$ is repeated ${}\geq n(x, g)$ times. An element $g$ is Engel if we can take $E(g)=\{1\}$. The set of Engel elements of a linear group is a subgroup by a well-known result of Gruenberg (\emph{J.~Algebra}, \textbf{3} (1966), 291--303). Is the set of almost Engel elements of a linear group a subgroup?

 \makebox[15pt][r]{}A linear group in which all elements are almost Engel is finite-by-hypercentral (P.\,Shumyatsky, \emph{Monatsh. Math.}, \textbf{186} (2018), 711--719).
 \hf\raisebox{-1ex}{\sl P.\,Shumyatsky}
\emp

 \bmp \textbf{\zv 20.90.}
 A profinite group in which all centralizers of non-trivial elements are pronilpotent is called a CN-group. Find an example of a finitely generated infinite profinite
CN-group which is not prosoluble.

\makebox[15pt][r]{}The structure of profinite CN-groups is described in (P.\,Shumyatsky,
\textit{Israel J.~Math.}, \textbf{235}, no.\,1 (2020), 325--347).
\hf \raisebox{-1ex}{\sl P.\,Shumyatsky}

\ul

\otv An example is found (AI solution confirmed by P.\,Shumyatsky:  M.\,Kinyon, J.\,Urban, {\it Preprint of 17 September 2026},
\url{https://kourovkanotebookorg.wordpress.com/wp-content/uploads/2026/09/kourovka-mathematics.pdf}; \ L.\,Mendonc,  T.\,S.\,Weigel,  T.\,Zapata, {\it Preprint of 30 October 2026},
\url{https://kourovkanotebookorg.wordpress.com/wp-content/uploads/2026/11/20-90.pdf}).
\emp

\bmp \textbf{20.91.}
A group $G$ is said to be \textit{stable} if for any first-order formula $\phi(\bar{x}, \bar{y})$ (in the first-order language of groups $\{\,\cdot\,, {}^{-1}, 1\}$) there exists a natural number $n$ such that whenever there exist sequences of tuples of $G$, $(a_i)_{i<m}$, $(b_i)_{i<m}$ with $G\models \phi(a_i, b_j)$ if and only if $i<j$, then $m\leq n$. Sela proved that all torsion-free hyperbolic groups are stable.

\makebox[25pt][r]{a)} Is there a non-stable hyperbolic group?

\makebox[25pt][r]{b)} Is there a non-stable virtually free group?

\makebox[25pt][r]{c)} Is $SL_2(\mathbb{Z})$ non-stable?
 \hf\raisebox{0ex}{\sl R.\,Sklinos}
\emp

\bmp \textbf{\ai 20.92.} 
This question is about existence of an analogue of the Lazard correspondence for pre-Lie algebras and braces. A {\em pre-Lie algebra} $A$ is a vector space with a bilinear operation $(x, y) \rightarrow xy$ satisfying $(xy)z -x(yz) = (yx)z - y(xz)$ for every $x,y,z\in A$.
A~pre-Lie algebra $A$ is said to be {\em left nilpotent} if, for some $n\in \mathbb N$, $A\cdot (A\cdot( A\cdots A)))=0$ where $A$ appears $n$ times in the product.
Recall that a set $A$ with binary operations $+$ and $\circ $ is a {\em left brace} if $(A, +)$ is an abelian group, $(A, \circ )$ is a group, and $a\circ (b+c)+a=a\circ b+a\circ c$ for every $a,b,c\in A$. Let $p$ be a prime, and $\mathbb F_{p}$ the field of $p$ elements. A~left brace $A$ is called an $\mathbb F_p$-brace if its additive group is an $\mathbb F_p$-vector space such that
$a*({\alpha }b)={\alpha }(a* b)$ for all $a,b\in A$,
$ {\alpha }\in \mathbb F_p$, where $ a*b=a \circ b -a -b$. The idea of a connection between braces and pre-Lie algebras comes from a paper by W.\,Rump (2014).

\makebox[25pt][r]{\ai a)} Let $A$ be an $\mathbb F_{p}$-brace of cardinality $p^{k}$ for some $k$. Is it true that when $p$ is sufficiently large relative to $k$, the set $A$ with the same additive operation $+$ and with the operation $\cdot $ defined as $a\cdot b=-\sum_{i=0}^{p-2}{\frac 1{2^{i}}}((2^{i}a)* b)$ is a pre-Lie algebra?

\makebox[15pt][r]{}{\it Comment of 2025}:
If $2\equiv \xi ^{p^{n-1}}$ when $\xi $ is a primitive root modulo $p$, and if $A$ is strongly nilpotent of nilpotency index less than $p$, then the result follows by (A.\,Smoktunowicz, {\em Adv. Math.}, {\bf 409}, part~B (2022), 108683). It is not known if the result follows when $A$ is not strongly nilpotent.\vskip1ex

\makebox[25pt][r]{\zva b)}
 Let $k$ be a natural number, and let $p$ be a prime number such that $p>2^{k}$. Is there a bijective correspondence between $\mathbb F_{p}$-braces of cardinality $p^{k}$ and left nilpotent pre-Lie algebras over ${\mathbb F}_{p}$ of cardinality $p^{k}$?

 \makebox[15pt][r]{}An affirmative answer to any of the above questions would have consequences for the theory of set-theoretic solutions of the Yang--Baxter equation and for the theory of Hopf--Galois extensions.
\hf\raisebox{-1ex}{\sl A.\,Smoktunowicz}

\ul

\oai a) Yes, it is (unverified AI solution:  M.\,Kinyon, J.\,Urban, {\it Preprint}, 2026,

\url{https://kourovkanotebookorg.wordpress.com/wp-content/uploads/2026/09/kourovka-mathematics.pdf}).

\otv b) Yes, there is (S.\,Trappeniers, {\em Preprint}, 2024, \url{https://arxiv.org/abs/2406.02475}).
\emp 

\bmp \textbf{20.93.}
 A group $G$ is called a \emph{Shunkov group} if for any finite subgroup $H\leq G$ any two conjugate elements of prime order in $N_G(H)/H$ generate a finite subgroup. Are the following well-known results of the theory of finite groups true in the class of (periodic) Shunkov groups?

 \makebox[25pt][r]{a)} The Baer--Suzuki theorem (see 11.11 in Archive).

 \makebox[25pt][r]{b)} The Burnside--Brauer--Suzuki theorem on the existence of a normal section of order~2 in a group with a non-trivial Sylow 2-subgroup containing only one involution (see 4.75).

 \makebox[25pt][r]{c)} Glauberman's $Z^*$-theorem (see 10.62 and Archive, 11.13).
\hf\raisebox{0ex}{\sl A.\,I.\,Sozutov}
\emp

\bmp \textbf{20.94.}
Does there exist an infinite periodic simple group saturated (see the definition in~14.101) with finite Frobenius groups?
\hf\raisebox{-1ex}{\sl A.\,I.\,Sozutov}
\emp

\bmp \textbf{20.95.}
Is a periodic group a Frobenius group (see 6.53) if it is saturated with finite Frobenius groups, contains an involution, and does not contain non-cyclic subgroups of order~4?
\hf\raisebox{-1ex}{\sl A.\,I.\,Sozutov}
\emp

\bmp \textbf{\zv 20.96.}
Is a periodic group a Frobenius group if it has a proper non-trivial normal abelian subgroup
that contains the centralizer of each of its non-identity elements?

\hf\raisebox{-1ex}{\sl A.\,I.\,Sozutov}

\ul

\otv Yes, it is (D.\,V.\,Lytkina, V.\,D.\,Mazurov, {\it Siberian Math. J.}, {\bf 64}, no.\,6 (2023), 1350--1353).
\emp

\bmp \textbf{20.97.}
 Is a 2-group locally finite if the centralizer of every involution is locally finite?

\hf\raisebox{-1ex}{\sl N.\,M.\,Suchkov}
\emp

\bmp \textbf{\zv 20.98.}
Let $G$ be a group of permutations of the set of positive integers $\N$ isomorphic to the additive group of rational numbers. Must there be an element $g\in G$ such that the set $\{a-ag\mid a\in \N\}$ is infinite?
\hf\raisebox{-1ex}{\sl N.\,M.\,Suchkov}

\ul

\otv Yes, it must (N.\,M.\,Suchkov, A.\,A.\,Shlepkin,\, D.\,A.\,Taysnyov, {\it Siberian Math.~J.}, {\bf 65} (2024), 1390--1394).
\emp

\bmp \textbf{20.99.}
\emph{Conjecture}: Let $a_1G_1,\ldots,a_kG_k$, $k>1$, be finitely many pairwise disjoint left cosets in
a group $G$ with $[G:G_i]<\infty$ for all $i=1,\ldots,k$. Then
\vs \vs \centerline{$\gcd([G:G_i],[G:G_j])\ge k\qquad \text{for some }1\leq i<j\leq k$.}
\vs \vs
\makebox[15pt][r]{}The conjecture is known to hold for $k=2,3,4$.
 \hf\raisebox{0ex}{\sl Z.-W.\,Sun}
\emp

\bmp \textbf{20.100.}
\emph{Conjecture}:
Let $n$ be a positive integer, and let $G$ be a group
containing no elements of order among $2,\ldots,n+1$.
Then, for any $A\subseteq G$ with $|A|=n$, we may write $A=\{a_1,\ldots,a_n\}$
with $a_1,a_2^2,\ldots,a_n^n$ pairwise distinct.

\makebox[15pt][r]{}The conjecture is known to hold for $n\leq 3$.
 \hf\raisebox{0ex}{\sl Z.-W.\,Sun}
\emp

\bmp \textbf{\zv 20.101.}
Let $G$ be a finitely generated branch group. Are all finite-index maximal subgroups of $G$ necessarily normal?

\makebox[15pt][r]{}Cf. 18.81.
\hf\raisebox{0ex}{\sl A.\,Thillaisundaram}

\ul

\otv No, not necessarily (P.\,Neumann, \emph{Illinois J. Math.}, \textbf{30} (1986), 301--316). Moreover, there are finitely generated branch groups all of whose maximal subgroups are of finite index but not all of them are normal (M.\,Garciarena, M.\,Petschick, \emph{Preprint}, 2024, \url{https://arxiv.org/abs/2410.06783}).
\emp

\bmp \textbf{20.102.} (Y.\,Barnea, A.\,Shalev). Let $G$ be a finitely generated pro-$p$ group. Must $G$ be $p$-adic analytic if it has finite Hausdorff spectrum with respect to

\makebox[25pt][r]{a)} the $p$-power series?

\makebox[25pt][r]{b)} the iterated $p$-power series?

\makebox[25pt][r]{c)} the lower $p$-series?

\makebox[25pt][r]{d)} the Frattini series?

\makebox[25pt][r]{e)} the dimension subgroup series?
\hf\raisebox{0ex}{\sl A.\,Thillaisundaram}
\emp

\bmp \textbf{20.103.}
Let $G_\Gamma$ be a partially commutative soluble group
of derived length $n \geq 3$ with defining graph $\Gamma$ (the
definition is similar to the case of $n=2$, see 17.104). Is
$G_\Gamma$ a torsion-free group?
\hf\raisebox{-1ex}{\sl E.\,I.\,Timoshenko}
\emp

\bmp \textbf{20.104.} (Well-known questions).
Suppose that a group $G$ is finitely generated and decomposable into a direct product $G = G_1 \times G_2$.

\makebox[25pt][r]{a)} Is it true that the elementary theory of $G$ is decidable if and only if the elementary theories of the groups $G_1$ and $G_2$ are decidable?

\makebox[25pt][r]{b)} Is it true that the universal theory of $G$ is decidable if and only if the universal theories of the groups $G_1$ and $G_2$ are decidable?

 \makebox[15pt][r]{}An affirmative answer to question a) for almost soluble groups follows from (G.\,A.\,Noskov, {\it Izv. Akad. Nauk SSSR, Ser. Mat.}, {\bf 47}, no. 3 (1983), 498--517).

\hf\raisebox{-1ex}{\sl E.\,I.\,Timoshenko}
\emp

\bmp \textbf{20.105.}
Is there a perfect locally nilpotent $p$-group, for some prime~$p$, whose proper subgroups are hypercentral?
 \hf\raisebox{-1ex}{\sl N.\,Trabelsi}
\emp

\bmp \textbf{20.106.}
Let $G$ be a residually finite $2$-group and let $x\in G$ be a left 3-Engel element of order~2. Is $\langle x^{G}\rangle$ locally nilpotent?

\makebox[15pt][r]{}It is known that in any group a 3-Engel element of odd order belongs to the locally nilpotent radical (E.\,Jabara, G.\,Traustason, \emph{Proc. Amer. Math. Soc.}, \textbf{147}, no.\,5 (2019), 1921--1927).

\makebox[15pt][r]{}{\it Comment of 2025}:
An element $a\in G$ is called a strong left $3$-Engel element if $\langle a,a^{g}\rangle$ is nilpotent of class at most $2$ and $\langle a,a^{g},a^{h}\rangle$ is nilpotent of class at most $3$ for all $g,h\in G$. (This is equivalent to $a$ being left $3$-Engel when $a$ is of odd order.) It is proved that if $a$ is a strong left $3$-Engel element in an arbitrary group $G$, then $\langle a\rangle^{G}$ is locally nilpotent (A.\,Hadjievangelou, G.\,Traustason, {\it Proc. Amer. Math. Soc.},
 {\bf 152}, no.\,4 (2024), 1467--1477).\hfill \raisebox{-1ex}{\sl G.\,Traustason}
\emp

\bmp \textbf{20.107.}
Let $G$ be a group of exponent 8 and $x\in G$ a left 3-Engel element of order~2. Is $\langle x^{G}\rangle$ locally finite?
\hf \raisebox{-1ex}{\sl G.\,Traustason}
\emp

\bmp \textbf{20.108.}
The \emph{holomorph} $\mathrm{Hol}(G)$ of a group $G$ can be defined as the normalizer of the subgroup of left translations in the group of all permutations of the set $G$. The \emph{multiple holomorph} $\mathrm{NHol}(G)$ of $G$ is the normalizer of the holomorph. Set $T(G) = \mathrm{NHol}(G)/\mathrm{Hol}(G)$.

\makebox[25pt][r]{\zv a)} Is there a centerless group $G$ for which $T(G)$ is not an elementary abelian 2-group?

\makebox[25pt][r]{\zv b)} Is there a finite centerless group $G$ for which $T(G)$ is not an elementary abelian 2-group?

\makebox[25pt][r]{c)} (A.\,Caranti). Is there a finite $p$-group $G$ for which the order of $T(G)$ has a prime divisor not dividing $(p-1)p$?
\hf\raisebox{-1ex}{\sl C.\,Tsang}

\ul

\otv a), b) Yes, there is, in both cases  (AI solution confirmed by C.\,Tsang:  \ M.\,Kinyon, J.\,Urban, {\it Preprint}, 2026,
\url{https://kourovkanotebookorg.wordpress.com/wp-content/uploads/2026/09/kourovka-mathematics.pdf}).
\emp 

\bmp \textbf{\zv 20.109.}
A {\it skew brace} is a set $B$ equipped with two operations $+$ and $\cdot$ such that $(B,+)$ is an additively written (but not necessarily abelian) group, $(B,\cdot )$ is a multiplicatively written group, and $a\cdot (b+c)=ab-a+ac$ for any $a,b,c\in B$.

\makebox[15pt][r]{}Is there a finite skew brace with perfect additive group and non-perfect almost simple multiplicative group?
\hf\raisebox{-1ex}{\sl C.\,Tsang}

\ul

\otv Yes, there is (M.\,Di\,Matteo, M.\,Ferrara, {\it Preprint}, 2026,  \url{https://arxiv.org/abs/2607.14059}
\emp

\bmp \textbf{20.110.}
Are there residually finite hereditarily just infinite groups that are

\makebox[25pt][r]{a)} amenable but not solvable?

\makebox[25pt][r]{b)} amenable but not elementary amenable?

\makebox[25pt][r]{c)} of intermediate word growth?

\makebox[25pt][r]{d)} of intermediate subgroup growth?

\makebox[15pt][r]{}All examples that we know are either linear (hence Tits Alternative applies), or have a quotient with property (T) (hence cannot be amenable).
\hf\raisebox{-1ex}{\sl M.\,Vannacci}
\emp

\bmp \textbf{20.111.}
Are there residually finite hereditarily just infinite groups admitting a self-similar action on a rooted tree that are not linear?
\hf\raisebox{-1ex}{\sl M.\,Vannacci}
\emp

\bmp \textbf{20.112.}
Let $\mathfrak{F}$ be a hereditary saturated formation and let ${\rm w}^{*}\mathfrak{F}$ denote the class of all
finite groups $G$ for which
$\pi(G)\subseteq\pi(\mathfrak{F})$ and the normalizers of all Sylow subgroups of $G$ are $\mathfrak{F}$-subnormal in $G$.
It is known that ${\rm w}^{*}\mathfrak{F}$ is a formation.
Must ${\rm w}^{*}\mathfrak{F}$ be a~saturated formation?
 \hf\raisebox{-1ex}{\sl A.\,F.\,Vasil'ev, T.\,I.\,Vasil'eva}
\emp

\bmp \textbf{20.113.}
Let $H(q,c) = \langle a,b,c,d \mid [b,a]=[d,c], \text{ of exponent } q, \text{ nilpotent of class } c\rangle$.

 \makebox[25pt][r]{a)} Is it true that the Schur multiplier $M(H(8,12))$ has exponent 32?

 \makebox[25pt][r]{b)} Is it true that the Schur multiplier $M(H(7,13))$ has exponent 49?

 \makebox[15pt][r]{}The difficulty is that these groups are too big to compute using current versions of the $p$-Quotient Algorithm, which use 32 bit arithmetic. So to tackle these groups
it would help to have a version of the $p$-Quotient Algorithm using 64 bit arithmetic.

\hf\raisebox{-1ex}{\sl M.\,R.\,Vaughan-Lee}
\emp

\bmp \textbf{20.114.}
Is there a positive integer $c$ such that for any $q$ the exponent of the Schur multiplier of a finite group of exponent $q$ divides $q^c$?
\hf\raisebox{-1ex}{\sl M.\,R.\,Vaughan-Lee}
\emp

\bmp \textbf{20.115.}
Let $\chi$ be a complex irreducible character of a finite group $G$. If $\chi(x)\ne 0$ for some $x\in G$, must the order $o(x)$ of $x$ divide $|G|/\chi(1)$?

\makebox[15pt][r]{}This is known to be true if $G$ is solvable, and it is known that $(o(x)\chi(1))^4$ divides $|G|^5$ for arbitrary $G$.
 \hf\raisebox{-1ex}{\sl T.\,Wilde}
\emp

\bmp \textbf{20.116.} (Well-known problem).
Does a profinite torsion group have finite exponent?

 \makebox[15pt][r]{}This problem reduces to the case of pro-$p$ groups (cf. W.\,Herfort, {\it Arch. Math.}, \textbf{33} (1980), 404--410). Also cf. Archive 3.41.
\hf\raisebox{-1ex}{\sl John S.\,Wilson}
\emp

\bmp \textbf{20.117.}
Let $G$ be a pro-$p$ group that is a 3-dimensional Poincar\'e duality group. Is $G$ coherent?

 \makebox[15pt][r]{}A group is said to be \textit{coherent} if each of its finitely generated subgroups is finitely presented, and in the question the coherency is used in the pro-$p$ sense. The question is a famous problem for abstract groups, but has a positive answer for 3-manifold groups.
\hf\raisebox{-1ex}{\sl P.\,Zalesskii}
\emp

\bmp \textbf{20.118.}
Let $G$ be a pro-$p$ group that is a 3-dimensional Poincar\'e duality group. Can it contain a direct product $F_2\times F_2$ of free pro-$p$ groups of rank 2?

 \makebox[15pt][r]{}If it can, then the answer to 20.117 is negative, but in the abstract case it is known that it can not (P.\,H.\,Kropholler, M.\,A.\,Roller, {\it J.~London
Math. Soc. (2)}, {\bf 39} (1989), 271--284).
\hf\raisebox{-1ex}{\sl P.\,Zalesskii}
\emp

\bmp \textbf{20.119.} (G.\,Wilkes).
A pro-$p$ group is said to be \textit{accessible} if there is a
number $n = n(G)$ such that any finite proper reduced graph of pro-$p$ groups with finite edge groups having fundamental group isomorphic to $G$ has at most $n$ edges (\textit{J.~Algebra}, {\bf 525} (2019), 1--18).

 \makebox[15pt][r]{}Is every finitely presented pro-$p$ group accessible?
\hf\raisebox{0ex}{\sl P.\,Zalesskii}
\emp

\bmp \textbf{\zv 20.120.}
Let $G$ be a locally nilpotent group with an automorphism (of infinite order). Can $G$ be simple as a group with automorphism?
\hf \raisebox{-1ex}{\sl E.\,I.\,Zelmanov}

\ul

\otv  Yes, it can; as noted by  I.\,S.\,Nnuga Mongo in ({\it Preprint}, 2026, \url{https://arxiv.org/abs/2611.18465}), an example is provided by Theorem~D in (P.\,Hall, {\it Proc. Cambridge Philos.
Soc.}, {\bf  58} (1962), 170--184): for the cyclic group $C_p$ of prime order $p$, the wreath power $W=\operatorname{Wr}C_p^{\Z}$ admits an infinite cyclic group of automorphisms  $\langle\f\rangle$ induced by the order-preserving permutations of $\Z$,  the derived subgroup $W'$ has no normal $\f$-invariant subgroups, and $W'$ is a locally finite $p$-group.
\emp

\bmp \textbf{\zv 20.121.}
Let $G$ be a finite almost simple group such that $Soc(G)$ is a non-soluble group of Lie type over a field of characteristic~$p$. Let $R$ be a Sylow $q$-subgroup of $G$ (for some $q$) and let $\operatorname{Min}_G(R)$ be the subgroup of $R$ generated by all minimal by inclusion intersections of the form $R\cap R^g$, where $g\in G$.
 Is it true that for $p>3$ the subgroup $\operatorname{Min}_G(R)$ is non-trivial if and only if the following hold: $G=\operatorname{Aut}(L_2(p))$, where $p$ is a Mersenne prime, $\operatorname{Min}_G(R)=R$, and $q=2$.

 \makebox[15pt][r]{}It is known that this is not always the case for $p=2,3$.
\hf\raisebox{0ex}{\sl V.\,I.\,Zenkov}

\ul

\otv Yes, it is true (T.\,C.\,Burness, H.\,Y.\,Huang, {\it Preprint}, 2025, \url{https://arxiv.org/pdf/2506.19745}).
\emp

\bmp \textbf{\zv 20.122.}
For nilpotent subgroups $A,B, C$ of a finite group $G$, let
 $\operatorname{Min}_G(A,B,C)$ be the subgroup of $A$ generated by all minimal by inclusion
 intersections of the form $A\cap B^x\cap C^y$, where $ x, y\in G$, and let $\operatorname{min}_G(A,B,C)$ be the subgroup of $\operatorname{Min}_G(A,B,C)$ generated by all intersections of this kind of minimal order.

 \makebox[25pt][r]{a)} Is it true that $\operatorname{min}_G(A,B,C)\leq F(G)$?

\makebox[25pt][r]{b)} Is it true that $\operatorname{Min}_G(A,B,C)\leq F(G)$?

\makebox[25pt][r]{c)} The same questions for soluble groups.
\hf\raisebox{0ex}{\sl V.\,I.\,Zenkov}

\ul

\otv  a), b), c) No, it is not always true in all three cases (AI solution confirmed by V.\,I.\,Zenkov:   \ M.\,Kinyon, J.\,Urban, {\it Preprint}, 2026,
\url{https://kourovkanotebookorg.wordpress.com/wp-content/uploads/2026/09/kourovka-mathematics.pdf}).
\emp 

\bmp \textbf{20.123.} A finite group is called a {\it $D_{\pi} $-group\/} if any two of its maximal $\pi
 $-sub\-groups
are conjugate.

 \makebox[25pt][r]{a)} Is it true that for any finite $D_\pi$-group $G$ and a $\pi$-Hall subgroup $H$ of $G$, there are elements $x,y,z\in G$ such that $O_{\pi }(G)=H\cap H^x\cap H^y\cap H^z$?

\makebox[25pt][r]{\zva b)}
 Suppose that $G$ is a finite $D_{\pi} $-group in which all simple non-abelian composition factors are sporadic or alternating groups, and let $H$ be a Hall $\pi $-subgroup of $G$. Is it true that $H\cap H^x\cap H^y= O_\pi(G)$ for some $x,y\in G$?

\makebox[25pt][r]{\zva c)}
Suppose that $G$ is a finite $D_{\pi} $-group with trivial soluble radical in which all simple non-abelian composition factors are sporadic groups, and let $H$ be a Hall $\pi $-sub\-group of $G$. Is it true that $H\cap H^g= O_\pi(G)$ for some $g\in G$? \hfill \raisebox{-1ex}{\sl V.\,I.\,Zenkov}

\ul

\otv b) Yes, it is true (I.\,N.\,Belousov, V.\,I.\,Zenkov, Preprint, 2025 (Russian),

\url{https://kourovkanotebookorg.wordpress.com/wp-content/uploads/2025/09/20.123b.pdf})

\otv c) Yes, it is true (I.\,N.\,Belousov, V.\,I.\,Zenkov, {\it Trudy Inst. Mat. Mekh. Ural. Otdel. Ross. Akad. Nauk}, {\bf 31}, no.\,1 (2025), 19--35 (Russian)).
\emp

\bmp \textbf{20.124.}
A {\it Rota--Baxter operator} on a group $G$ is a mapping $B : G \to G$ such that
$B(g) B(h) = B\left(gB(g)hB(g)^{-1} \right)$ for all $g, h \in G$.
Let $F$ be a non-abelian free group. Is there a Rota--Baxter operator on $F$ such that its image is equal to the derived subgroup $[F, F]$?
\hf\raisebox{-1ex}{\sl V.\,G.\,Bardakov}
\emp

\bmp \textbf{\zv 20.125.}
Does there exist a non-abelian group $G$ and a Rota--Baxter operator $B : G \to G$ such that $B$ is surjective but not injective?
\hf\raisebox{-1ex}{\sl V.\,G.\,Bardakov}

\ul

\otv Yes, it does exist (AI solution confirmed by V.\,G.\,Bardakov: \  P.\,Monticone, {\it Preprint}, 2026,  \url{https://kourovkanotebookorg.wordpress.com/wp-content/uploads/2026/03/20_125-1.pdf});

 see also (W.\,van\,Doorn, E.\,Judin, P.\,Monticone, D.\,Morrison,  {\it Preprint}, 2026, \url{https://arxiv.org/abs/2607.17477}).
\emp

\bmp \textbf{\zv 20.126.} 
A \textit{brace} $(G; +, \circ)$ is non-empty set $G$ with two binary operations $+$, $\circ$ such that
 $(G, +)$ is an additively written abelian group, $(G, \circ)$ is a multiplicatively written group, and $a \circ (b + c) + a = (a \circ b) + (a \circ c)$ for all $a, b, c \in G$. Does there exist a brace with finitely generated group $(G, +)$ such that

 \makebox[25pt][r]{(a)} the group $(G, \circ)$ is non-solvable?

 \makebox[25pt][r]{(b)} the group $(G, \circ)$ contains a non-abelian free group?

\hf\raisebox{-1ex}{\sl V.\,G.\,Bardakov, M.\,V.\,Neshchadim, M.\,K.\,Yadav}

\ul

\otv (a) and (b) Yes, such a brace does exist (L.\,Gong, {\it Preprint}, 2026, \url{https://arxiv.org/abs/2611.14972}).
\emp
\newpage

\pagestyle{myheadings}
\markboth{\protect\vphantom{(y}New Problems (21st issue, 2026)}{\protect\vphantom{(y}New Problems (21st issue, 2026)}
\thispagestyle{headings}

~

\vspace{2ex}

\centerline {\Large \textbf{ New Problems}}

\phantomsection\label{21izd}
\parindent=0pt
\vspace{4ex}

\bmp \textbf{21.1.}
Let $n$ be a positive integer. For a finite group $K$  and an automorphism $\phi$ of $K$ of order dividing $n$, let $X_{n,\phi}(K):=\{x \in K \mid  xx^\phi \cdots x^{\phi^{n-1}} = 1\}$.
  Let  $c_n$ be the supremum of the ratios $|X_{n,\phi}(H)|/|H|$ over all finite groups $H$ and their automorphisms $ \phi\in{\rm Aut}(H)$  such that   $\phi^n={\rm id}$ and $X_{n,\phi}(H)\ne H$.

 \makebox[25pt][r]{(a)} Let $n>1$ be a positive integer such that $c_d < 1$ for all prime power divisors $d$ of $n$. Is it true that $c_n < 1$?

 \makebox[25pt][r]{(b)} For a finite group $G$ and a positive integer $n$,  the generalized Hughes--Thompson subgroup is defined as $H_n(G)=\langle x\in G \mid  x^n \neq 1\rangle$.  Suppose that $n$ is a positive integer for which there is  a positive integer $k_n$ depending only on $n$ such that $|G : H_n(G)|\leq k_n$ for all finite groups $G$ with $ H_n(G)\ne 1$. Is it true that then $c_n < 1$?
 This question is open even when $n \geq 5$ is prime.
 \hf\raisebox{-1ex}{\sl A.\,Abdollahi, M.\,S.\,Malekan}
\emp

\bmp \textbf{21.2.}
Let $S$ be a finite simple group, and let $G$ be a finite group for which there exists a bijection $f:G\to S$ such that $|x|$ divides $|f(x)|$ for all $x\in G$. Must $G$ necessarily be simple?
 \hfill \raisebox{-1ex}{\sl M.\,Amiri}
\emp

\bmp \textbf{21.3.}
Let $G = A_n$ or $S_n$ and let $H,K$ be soluble subgroups of $G$.

\makebox[25pt][r]{\zva a)} For all sufficiently large $n$, can we always find an element $x \in G$ such that $H \cap K^x = 1$?

\makebox[25pt][r]{b)} Does this hold for all $n \geqslant 21$?

\makebox[15pt][r]{}Note that the conclusion is false when $G = S_{20}$ and $H = K = (S_4 \wr S_4) \times S_4$.

\hf \raisebox{-1ex}{\sl M.\,Anagnostopoulou-Merkouri, T.\,C.\,Burness}

\ul

\otv a) Yes, we can (A.\,Rizzoli, {\it Preprint}, 2026, \url{https://kourovkanotebookorg.wordpress.com/wp-content/uploads/2026/09/kourovka_21_3a.pdf}).
\emp

\bmp \textbf{21.4.} Let $G$ be a finite group with trivial solvable radical and let $H_1, \ldots, H_5$ be solvable subgroups of $G$. Then do there always exist elements $x_i \in G$ such that $\bigcap_i H_i^{x_i} = 1$? Cf.~17.41(b).
\hf \raisebox{-1ex}{\sl M.\,Anagnostopoulou-Merkouri, T.\,C.\,Burness}
\emp

\bmp \textbf{21.5.}
Let $p$ be a prime. Let $G$ be a transitive subgroup of the group of finitary permutations $FSym(\Omega)$ of a set $\Omega$, let $N$ be a normal subgroup of $G$, and let $S$ be a transitive Sylow $p$-subgroup of $G$.

\makebox[25pt][r]{\zva a)} Is it true that $S\cap  N$ is a Sylow $p$-subgroup of $N$?

\makebox[25pt][r]{\zva b)} Is it true that $SN/N$ is a Sylow $p$-subgroup of $G/N$?

\makebox[25pt][r]{c)} Are any two transitive  Sylow $p$-subgroups of $G$ locally conjugate in $G$?

\makebox[15pt][r]{}Two subgroups $X,Y$ of a group $G$ are said to be \emph{locally conjugate} if there is a~locally inner  automorphism $\varphi$ of $G$ such that $X^\varphi=Y$. An automorphism $\varphi$ of $G$ is said to be \emph{locally inner} if for every finite subset $A\subseteq  G$ there is an element $g=g(A)\in G$ such that $a^\varphi=g^{-1}a{g}$ for all $a\in A$.
 \hfill \raisebox{-1ex}{\sl A.\,O.\,Asar}

 \ul

 \otv a) Yes, it is true (AI solution confirmed by F.\,Leinen: \ P.\,Monticone, {\it Preprint}, 2026,
\url{https://kourovkanotebookorg.wordpress.com/wp-content/uploads/2026/08/21_5.pdf}).

 \otv b) Yes, it is true (AI solution confirmed by F.\,Leinen: \  P.\,Monticone, {\it Preprint}, 2026,
\url{https://kourovkanotebookorg.wordpress.com/wp-content/uploads/2026/08/21_5.pdf}).
\emp

\bmp \textbf{21.6.}
Let $p$ be a prime. A totally imprimitive $p$-group $H$ of finitary permutations is said to have the \emph{cyclic-block property} if in the cycle decomposition of every element the support of every cycle is a block for $H$. Let $G$ be a transitive subgroup of the group of finitary permutations $FSym(\Omega)$ of a set $\Omega$.
Does every transitive Sylow $p$-subgroup of $G$ contain a transitive subgroup which has the cyclic-block property?

 \hfill \raisebox{-1ex}{\sl A.\,O.\,Asar}
\emp

\bmp \textbf{21.7.} 
 (Well-known problem).
 A finite group $G$ is called an {\it IYB-group} if it is isomorphic to the permutation group of a finite involutive non-degenerate set-theoretic solution of the Yang--Baxter equation, or equivalently, $G$ is isomorphic to the multiplicative group of a finite left brace.  Assume that the Sylow subgroups of a finite soluble group $G$ are IYB-groups. Is $G$ an IYB-group?
 \hfill \raisebox{-1ex}{\sl A.\,Ballester-Bolinches}


 \emp

\bmp \textbf{\zv 21.8.}
 As in 17.57, let $r(m)=\{r+km\mid k\in \Z\}$ for integers
$0\leq r<m$; for \(r_1(m_1)\cap r_2(m_2) =\nobreak\varnothing\) let the
\emph{class transposition} \(\tau_{r_1(m_1),r_2(m_2)}\) be the
involution which interchanges  \(r_1+tm_1\) and \(r_2+tm_2\) for
each integer \(t\) and fixes everything else, and let  \({\rm
CT}(\mathbb{Z})\) be the group generated by all class transpositions.

\makebox[15pt][r]{}Let ${\rm CT}_k$ be the subgroup of  \({\rm
{\rm CT}}(\mathbb{Z})\) generated by the class transpositions
$\tau_{r_1(k),r_2(k)}$ for $0 \leq r_1 \not= r_2 < k$. Since    $\tau_{r_1(k),r_2(k)}$ permutes the residue classes modulo $k$, the group  ${\rm CT}_k$ is isomorphic to the symmetric group $S_k$. Let ${\rm CT}_{(k)}=\langle {\rm CT}_2, {\rm CT}_3, \ldots, {\rm CT}_k\rangle $. Is it true that for $k > 3$ the group ${\rm CT}_{(k)}$ is isomorphic to the symmetric group $S_N$, where $N$ is the least common multiple of the numbers $2, 3, \ldots, k$?

 \hfill \raisebox{-1ex}{\sl V.\,G.\,Bardakov, A.\,L.\,Iskra}

\ul

\otv Yes, it is true (Junyao Pan, {\it Preprint of 14 April 2026}, \url{https://arxiv.org/abs/2604.12553}; \ AI solution confirmed by A.\,L.\,Iskra: \ P.\,Monticone, {\it Preprint of 30 April 2026}, \url{https://kourovkanotebookorg.wordpress.com/wp-content/uploads/2026/05/21_8-1.pdf}); see also (W.\,van\,Doorn, E.\,Judin, P.\,Monticone, D.\,Morrison,  {\it Preprint}, 2026, \url{https://arxiv.org/abs/2607.17477}).
\emp

\bmp \textbf{21.9.} Let $F$ be a non-abelian free pro-$p$ group of finite rank. Can one find a finite collection $U_1, \dots  , U_n$ of open subgroups of $F$ including $F$ itself such that the only subgroup of $F$ which is contained in $U_i$ and is characteristic in $U_i$ for every $i$ is the trivial subgroup?
\hf \raisebox{-1ex}{\sl Y.\,Barnea}
\emp

\bmp \textbf{\zv 21.10.}
We call a group presentation \emph{finite} if it represents a finite group. We say that a presentation is \emph{just finite} if it is finite and is no longer finite on removal of any relation from it.
Is it true that every finite group has a just finite presentation?

\makebox[15pt][r]{}Note that if a group has a balanced presentation, then it is just finite. A similar argument can be applied for some $p$-groups using the  Golod--Shafarevich inequality.

\hf \raisebox{-1ex}{\sl Y.\,Barnea}

\ul

\otv Yes, it is true (M.\,Lackenby, {\it Preprint}, 2026, \url{https://arxiv.org/abs/2605.10402}).
\emp

\bmp \textbf{21.11.}
Can some or all groups of the following
sorts be written as homomorphic
images of nonprincipal ultraproducts of countable families of groups? This is Question~18 in (G.\,M.\,Bergman, {\it Pacific J.~Math.}, {\bf 274} (2015) 451--495).

\makebox[25pt][r]{(a)} Infinite finitely generated  groups of finite exponent.

\makebox[25pt][r]{(b)} For an infinite set $X$, the group of those permutations of $X$ that move only finitely many elements.

\makebox[15pt][r]{}It is known that no group of permutations containing an element with exactly one infinite orbit can be written as an image of such an ultraproduct (\emph{ibid.}).

\hfill \raisebox{-1ex}{\sl G.\,M.\,Bergman}
\emp

\bmp \textbf{\zv 21.12.}
 Suppose that $\mathscr{U}$ is a nonprincipal ultrafilter on $\omega$, and $B$ is a group such that every element $b\in B$ belongs to a subgroup of $B$ that is a  homomorphic image  of $\Z^{\omega}/\mathscr{U}$. Must $B$ then be a homomorphic image of an ultraproduct group $\big(\prod_{i\in \omega}G_i\big)/\mathscr{U}$ for some groups $G_i$?

\makebox[15pt][r]{}This is Question~19 in (G.\,M.\,Bergman, {\it Pacific J.~Math.}, {\bf 274} (2015) 451--495). An affirmative answer would imply that every torsion group was such a homomorphic
image for every $\mathscr{U}$, and so would give positive answers to both parts of 21.11.

\hfill \raisebox{-1ex}{\sl G.\,M.\,Bergman}

\ul

\otv No, it need not (S.\,M.\,Corson, \emph{J.~Algebra}, \textbf{681} (2025), 306--317).
\emp

\bmp \textbf{21.13.}
It is known
that the group $\Z^\omega$ has a subgroup whose dual is free abelian of rank $2^{\aleph_0}$ (see 17.24 in Archive). Does $\Z^\omega$ have a subgroup whose dual is free abelian
of still larger rank (the largest possible being $2^{2^{\aleph_0}}$)? This is Question~11 in  (G.\,M.\,Bergman, {\it Portugaliae Math.}, {\bf  69} (2012) 69--84).
\hfill \raisebox{-1ex}{\sl G.\,M.\,Bergman}
\emp

\bmp \textbf{\zv 21.14.}
Suppose $\alpha$ is an endomorphism of a group $G$ such that for
every group $H$ and every  homomorphism $f : G \to H$, there exists an endomorphism $\beta _f$ of $H$
such that $\beta _f f = f\alpha$.
Must  $\alpha$ then be either an inner
automorphism of $G$ or the trivial endomorphism? This is Question~5 in  (G.\,M.\,Bergman, {\it Publ.  Matem.}, {\bf 56} (2012), 91--126).
\hfill \raisebox{-1ex}{\sl G.\,M.\,Bergman}

\ul

\otv Yes, it must (F.\,Fournier-Facio, {\it Preprint}, 2026,
\url{https://arxiv.org/abs/2604.05728})
\emp

\bmp \textbf{\zv 21.15.}
Suppose $B$ is a subgroup of the symmetric group $S_{\Omega}$ on an infinite set $\Omega$. Will the amalgamated free product $S_{\Omega}*_B S_{\Omega}$ of two copies of $S_{\Omega}$ with amalgamation of
$B$ be embeddable in $S_{\Omega}$? This is a weakened form of
the group case of Question~4.4 in (G.\,M.\,Bergman, {\it Indag. Math.}, {\bf 18} (2007), 349--403).

\makebox[15pt][r]{}It is known that $S_{\Omega}*_B S_{\Omega}$ need not be so embeddable by a map
respecting $B$ (\emph{Algebra Number Theory}, {\bf 3}
(2009), 847--879, \S\,10).
\hfill \raisebox{-1ex}{\sl G.\,M.\,Bergman}

\ul

\otv No, not necessarily: take $\Omega$ countably infinite, $B$ a subgroup of $S_{\Omega}$ which is not Borel, and $\phi: S_{\Omega} *_B S_{\Omega} \rightarrow S_{\Omega}$ an injective group homomorphism.  Let $\phi_1$ be the restriction of $\phi$ to the first copy of $S_{\Omega}$ and $\phi_2$ that to the second.  Each $\phi_i$ is continuous (A.\,S.\,Kechris, C.\,Rosendal, \emph{Proc. Lond. Math. Soc.}, {\bf 94}, no.\,2 (2007), 302--350),  so $\operatorname{im}(\phi_i)$ is Borel in $S_{\Omega}$ (as a continuous injective image of a Polish space).  Now $B = \phi_1^{-1}(\phi(B)) = \phi_1^{-1}(\operatorname{im}(\phi_1) \cap \operatorname{im}(\phi_2))$ is Borel, a contradiction. (S.\,M.\,Corson, {\it Letter of 26 January 2026}; see also \url{https://kourovkanotebookorg.wordpress.com/wp-content/uploads/2026/03/solution-of-21.15.pdf}.)
\emp

\bmp \textbf{21.16.}
Let the \textit{width} of a group (respectively, a monoid) $H$ with respect
to a generating set $X$ mean the supremum over $h\in  H$ of the least length of a group word (respectively, a monoid word) in elements of $X$ expressing $h$. A group (or monoid) is said to have finite width if its width with respect to every generating set is finite.  (A~common finite
bound for these widths is not required.) Do there exist groups $G$ having finite width as groups, but not as monoids? This is Question~9  in (G.\,M.\,Bergman, {\it Bull. London Math. Soc.}, {\bf 38} (2006), 429--440).
\hfill \raisebox{-1ex}{\sl G.\,M.\,Bergman}
\emp

\bmp \textbf{21.17.}
If $\frak{X}$ is a class of groups, let $\textbf{H}(\frak{X})$ denote the class of homomorphic images of groups in $\frak{X}$, let $\textbf{S}(\frak{X})$ denote the class of groups isomorphic to subgroups of groups in~$\frak{X}$, let $\textbf{P}(\frak{X})$
denote the class of groups isomorphic to (unrestricted) direct 
products of families of groups in $\frak{X}$, and let $\textbf{P}_f (\frak{X})$
denote the class of groups isomorphic to direct products of finite families of groups in $\frak{X}$. By Birkhoff's theorem, $\textbf{H} (\textbf{S} (\textbf{P}(\frak{X})))$ is the variety of groups generated by $\frak{X}$.

 \makebox[15pt][r]{}If $\frak{M}$ is a class of metabelian groups, must $\textbf{H} ( \textbf{S} (\textbf{P}_f (\frak{M})))\subseteq \textbf{S} (\textbf{H} (\textbf{P} (\textbf{S}(\frak{M}))))$? This is Question~27  in (G.\,M.\,Bergman, {\it Algebra Universalis}, {\bf 26} (1989), 267--283).

 \hfill \raisebox{-1ex}{\sl G.\,M.\,Bergman}
\emp

\bmp \textbf{\zv 21.18.}
Suppose that $G$ is a finite group, and $A_1, A_2, A_3$ are subsets of $G$ such that
the multiplication map $A_1\times A_2\times A_3 \to G$ is bijective. Must the subgroup $\langle A_2\rangle$ generated
by $A_2$ have order divisible by the cardinality $|A_2|$? This is Question~8 in (G.\,M.\,Bergman, {\it J.~Iranian Math. Soc.}, {\bf 1} (2020), 157--161).

\makebox[15pt][r]{}It is known (\emph{ibid.)} that the corresponding statement is true for the subgroups
$\langle A_1\rangle$ and $\langle A_3\rangle$. Moreover, $|A_2|$ will at least divide the order of the
least subgroup containing $A_2$ and closed under conjugation by members of $A_1$, and similarly of the
least subgroup containing $A_2$ and closed under conjugation by members of $A_3$.

 \hfill \raisebox{-1ex}{\sl G.\,M.\,Bergman}

 \ul

 \otv No, it need not (M.\,I.\,Kabenyuk, {\it Preprint}, 2021, \url{https://arxiv.org/abs/2102.08605}).
 \emp

\bmp \textbf{21.19.}
 Suppose that $S$ and $M$ are groups of finite Morley rank,
$S$ is an infinite group,  and $M$ is a non-trivial connected group definably and faithfully acting on~$S$. This action is  said to be \emph{irreducible} if
$M$ does not leave invariant any definable non-trivial proper subgroup of $S$.

\makebox[15pt][r]{}Prove that if $S$ is a simple group such that every proper definable subgroup of $S$ is nilpotent,  and the action of $M$ on $S$ is irreducible, then this action
is equivalent to the action of $S$ on itself by conjugation.
\hfill \raisebox{-1ex}{\sl A.\,V.\,Borovik}
\emp

\bmp \textbf{21.20.}
 Prove that  a simple group of finite Morley rank without involutions cannot act definably, faithfully, and irreducibly on a connected group  other than on itself acting by conjugation.
\hfill \raisebox{-1ex}{\sl A.\,V.\,Borovik}
\emp

\bmp \textbf{21.21.}
Prove that a simple (that is, without proper non-trivial connected normal subgroups) algebraic group $M$ over an algebraically closed
field cannot act definably, faithfully, and irreducibly on a simple group of finite Morley rank other than $M/Z(M)$.

\hfill \raisebox{-1ex}{\sl A.\,V.\,Borovik}
\emp

\bmp \textbf{21.22.}
Is the (standard, restricted) wreath product $G \wr H$ of two finitely generated Hopfian groups Hopfian?

 \makebox[15pt][r]{}The same question where $G$ is assumed to be abelian or nilpotent is equivalent to Kaplansky's direct finiteness conjecture; see (H.\,Bradford, F.\,Fournier-Facio, {\it Math.~Z.}, {\bf 308}, no.\,4 (2024), Paper no.\,58).
\hf\raisebox{-1ex}{\sl H.\,Bradford, F.\,Fournier-Facio}
\emp

\bmp \textbf{21.23.}
A graph is called a \emph{cograph} if it has no induced subgraph isomorphic to a path with $4$ vertices. A graph is said to be \emph{chordal} if it has
no induced cycles with $n$ vertices for every $n\geq 4$. For a finite group $G$, the \emph{enhanced power graph} $\mathcal{E}(G)$
is the graph with vertex set $G$ and edges $\{x,y\}$ for all $x\neq y\in G$ such that $\langle x,y\rangle$ is cyclic.

\makebox[25pt][r]{(a)} For a given integer  $n\geq 4$, determine the set of all finite nonabelian simple groups $G$ such that $\mathcal{E}(G)$ has
no induced cycles with $n$ vertices.

\makebox[25pt][r]{(b)} Determine the set of all  finite nonabelian simple groups $G$ such that $\mathcal{E}(G)$ is chordal.

\makebox[15pt][r]{}In ({\it Preprint}, 2025, \url{https://arxiv.org/abs/2510.18073}) we proved that if the enhanced power graph of a given finite group is a cograph, then it is also chordal. Also the finite nonabelian simple groups whose enhanced power graph is a cograph are described, and additional information is obtained on finite nonabelian simple groups whose  enhanced power graph has no induced cycles with $4$ vertices.

\hf \raisebox{-1ex}{\sl D.\,Bubboloni, F.\,Fumagalli, C.\,E.\,Praeger}
\emp

\bmp \textbf{\zv 21.24.}
For a finite group $G$, the \emph{power graph} $\mathcal{P}(G)$
is the graph with vertex set $G$ and edges $\{x,y\}$ for all $x\neq y\in G$ such that either
$x\in \langle y\rangle$ or $y\in \langle x\rangle$.
Is it true that,  for every finite group $G$, if $\mathcal{P}(G)$ is a cograph, then $\mathcal{P}(G)$ is  chordal? Cf. 21.23.

\makebox[15pt][r]{}This holds if every element of $G$ has prime power order (D.\,Bubboloni, F.\,Fumagalli, C.\,E.\,Praeger, {\it Preprint}, 2025, \url{https://arxiv.org/abs/2510.18073})  and  if $G$ is a nonabelian simple group
(J.\,Cameron, P.\,Manna, R.\,Mehatari, \emph{J.~Algebra},
{\bf 591} (2022), 59--74; \ J.\,Brachter, E.\,Kaja, \emph{J.~Algebr. Comb.},  {\bf 58} (2023), 1095--1124).

\hf \raisebox{-1ex}{\sl D.\,Bubboloni, F.\,Fumagalli, C.\,E.\,Praeger}

\ul

\otv Yes, it is true (M.\,Rundstr\"om, {\it Preprint of 30 January 2026}, \url{https://kourovkanotebookorg.wordpress.com/wp-content/uploads/2026/04/21.24-runds.pdf}; \ AI solution confirmed by C.\,E.\,Praeger: \ P.\,Monticone, {\it Preprint of 29 March 2026},  \url{https://kourovkanotebookorg.wordpress.com/wp-content/uploads/2026/04/21_24-1.pdf}); see also (W.\,van\,Doorn, E.\,Judin, P.\,Monticone, D.\,Morrison,  {\it Preprint}, 2026, \url{https://arxiv.org/abs/2607.17477}).
\emp

\bmp \textbf{21.25.} (T.\,Breuer, R.\,M.\,Guralnick). Let $G$ be a finite simple group and let $p_1,p_2$ be any (not necessarily distinct) prime divisors of $|G|$. Then can we always find  Sylow $p_i$-subgroups $H_i$ such that $G = \langle H_1, H_2\rangle$?
\hf \raisebox{-1ex}{\sl T.\,C.\,Burness}
\emp

\bmp \textbf{21.26.} (F.\,Lisi, L.\,Sabatini). Let $G$ be a non-trivial finite group and let $p_1, \ldots, p_k$ be the distinct prime divisors of $|G|$. For each $i$, let $H_i$ be a Sylow $p_i$-subgroup of $G$. Is it true that there exists an element $x \in G$ such that for all $i$ the subgroup $H_i \cap H_i^x$ is inclusion-minimal in $\{H_i \cap H_i^g \,: \, g \in G\}$?
\hf \raisebox{-1ex}{\sl T.\,C.\,Burness}
\emp

\bmp \textbf{21.27.} (M.\,Larsen, A.\,Shalev, P.\,H.\,Tiep).  A permutation on a set $\Omega$ is called a \emph{derangement} if it has no fixed points in $\Omega$.
Let $G $ be a finite simple transitive permutation group.
Is it true that every element in $G$ is the product of two derangements?

\hf \raisebox{-1ex}{\sl T.\,C.\,Burness}
\emp

\bmp \textbf{21.28.} Let $G $ be a finite simple transitive permutation group, and let $\delta(G)$ be the proportion of derangements in $G$. Is it true that $\delta(G) \geqslant 89/325$?

\makebox[15pt][r]{}Note that $\delta(G) = 89/325$ for the action of the Tits group $G = {}^2F_4(2)'$ on the cosets of a maximal parabolic subgroup of the form $2^2.[2^8].S_3$ ({\it Forum Math. Sigma}, \textbf{13} (2025), paper no.\,e98, 62 pp.).
\hf \raisebox{-1ex}{\sl T.\,C.\,Burness, M.\,Fusari}
\emp

\bmp \textbf{\zv 21.29.}
  Let $G \leqslant {\rm Sym}(\Omega)$ be a finite primitive permutation group with a regular suborbit (that is, $G$ has a trivial $2$-point stabiliser). Then is it true that for all $\alpha,\beta \in \Omega$, there exists $\gamma \in \Omega$ such that the $2$-point stabilisers $G_{\alpha,\gamma}$ and $G_{\beta,\gamma}$ are both trivial?
\hf \raisebox{-1ex}{\sl T.\,C.\,Burness, M.\,Giudici}

\ul

\otv No, not always (A.\,Rizzoli, A.\,R.\,Thomas,  {\it Preprint}, 2026, \url{https://arxiv.org/abs/2609.01367}).
\emp

\bmp \textbf{21.30.} (Well-known question). A discrete group $G$ is said to have the  \emph{Haagerup property} (also known as  \textit{Gromov's a-T-menability property}) if there exists a metrically proper isometric action of $G$ on a (possibly infinite-dimensional) Hilbert space. Are all 1-relator groups Haagerup groups?
\hf \raisebox{-1ex}{{\sl J.\,O.\,Button}}
\emp

\bmp \textbf{21.31.}
{\it Conjecture}: If $N$ is a finite soluble group, then any regular subgroup
 in the holomorph $Hol(N)$ of $N$ is also soluble.
 \hf\raisebox{-1ex}{\sl N.\,Byott}
\emp

\bmp \textbf{\zv 21.32.}
Is the following problem decidable, and if so, what is its complexity?
Given a finite group $G$, is there a finite group $H$ such that the derived
subgroup of $H$ is isomorphic to $G$?
 \hfill \raisebox{-1ex}{\sl P.\,J.\,Cameron}

 \ul

 \otv Yes, it is decidable, in nondeterministic quasi-polynomial time (S.\,Kalithasan, V.\,Z.\,Thomas,  {\it Preprint}, 2026, \url{https://arxiv.org/abs/2602.18829}).
\emp

\bmp \textbf{21.33.}
Does an analogue of Dunwoody's theorem hold for totally disconnected locally compact groups, that is, must a tdlc group of rational discrete cohomological dimension  at most 1  be topologically isomorphic to the fundamental group of a  graph of profinite groups?
 \hf\raisebox{-1ex}{\sl I.\,Castellano}
\emp

\bmp \textbf{21.34.} (Well-known problem). A group $G$ is a \emph{unique product group} if, for any nonempty finite subsets $A, B$ of $G$, there exists an element of $G$ which can be written uniquely as $ab$ with $a \in A$ and $b \in B$.  A group $G$ is \emph{locally invariant orderable} if $G$ admits a partial order $<$ such that for all $g, h \in G$ with $h \neq 1$, we have either $gh > g$ or $gh^{-1} >g$.
Does there exist a unique product group which is not locally invariant orderable?
 \hfill \raisebox{-1ex}{\sl A.\,Clay}
\emp

\bmp \textbf{\ai 21.35.} 
Let $G$ be a finite group,  $w$  a multilinear commutator group-word, and $p$ a prime. Suppose that $p$ divides the order  $|xy|$ whenever $x$ is a $w$-value  of $p'$-order in $G$ and $y$ is a  $w$-value  in $G$ of order divisible by~$p$. Is it true that then the verbal subgroup $w(G)$ must be $p$-nilpotent?

\makebox[15pt][r]{}Without the assumption that $w$ be multilinear, the answer is negative. An affirmative answer has been obtained in several special cases (\emph{J.~Algebra}, {\bf 609} (2022),  926--936).
\hfill \raisebox{-1ex}{\sl Y.\,Contreras Rojas, V.\,Grazian,  C.\,Monetta}

\ul

\oai Yes, it is true (unverified AI solution:  Nilradical,  {\it Preprint}, 2026,

\url{https://kourovkanotebookorg.wordpress.com/wp-content/uploads/2026/09/nilradical-kourovka-21-35.pdf}).
\emp 

\bmp \textbf{21.36.}
Kropholler's hierarchy (see 15.45)  is closed under finite extensions, that is,
$(\textbf{H}_\alpha\mathfrak{F})\mathfrak{F} \subseteq \textbf{H}_\alpha\mathfrak{F}$ for every $\alpha$ (P.\,Kropholler, {\it  J.~Pure Appl. Algebra}, {\bf 90} (1993), 55--67). Let a hierarchy of tdlc groups $\textbf{H} \mathfrak{K}$ be defined analogously to Kropholler's hierarchy in 15.45, with $\mathfrak{K}$ being the class of profinite groups and with the cell stabilisers of the admissible action required to be open.

\makebox[15pt][r]{}Is it true that $\textbf{H}\mathfrak{K}$ is closed under profinite extensions, that is, $(\textbf{H}_\alpha\mathfrak{K})\mathfrak{K} \subseteq \textbf{H}_\alpha\mathfrak{K}$ for every $\alpha$?
 \hf\raisebox{-1ex}{\sl G.\,C.\,Cook}
\emp

\bmp \textbf{21.37.}
By definition, a \textit{constructible} totally disconnected, locally compact (tdlc) group is the result of a sequence of profinite extensions and ascending HNN-extensions starting from the trivial group. As in the discrete case, soluble constructible tdlc groups have type $FP_\infty$ (G.\,C.\,Cook,  I.\,Castellano, {\it J.~Algebra}, {\bf 543} (2020), 54--97).

 \makebox[15pt][r]{}Are soluble tdlc groups of type $FP_\infty$ constructible?
 \hf\raisebox{-1ex}{\sl G.\,C.\,Cook}
\emp

\bmp \textbf{\zv  21.38.}
(S.\,Harper, C.\,Donoven).
The \emph{spread} of a group $G$ is the greatest nonnegative integer $k$ such that for all nontrivial elements $x_1, \ldots, x_k \in G$ there exists $y \in G$ such that $\langle x_1, y \rangle = \cdots = \langle x_k, y \rangle = G$, or is $\infty$ in case there is no such maximum.  Does there exist a group with spread equal to $1$?

  \makebox[15pt][r]{}Such a group must be infinite if it exists (T.\,C.\,Burness, R.\,M.\,Guralnick, S.\,Harper, \emph{Ann. Math.}, {\bf  193} (2021), 619--687).
\hf \raisebox{-1ex}{\sl  S.\,Corson}

\ul

\otv Yes, such a group exists (AI solution confirmed by S.\,Corson: \  Nilradical,  {\it Preprint}, 2026, \url{https://kourovkanotebookorg.wordpress.com/wp-content/uploads/2026/09/nilradical-kourovka-21-38.pdf}).
\emp 

\bmp \textbf{21.39.}
It is known that there exist residually finite, locally finite, characteristically simple groups with finitely many orbits under automorphisms
(A.\,B.\,Apps, \emph{J.~Algebra}, {\bf 81} (1983), 320--339).
Are there any locally finite, characteristically simple groups with finitely many orbits  under automorphisms that are not residually finite?

\hfill \raisebox{-1ex}{\sl A.\,Dantas, E.\,de\,Melo}
\emp

\bmp \textbf{\ai 21.40.} 
 Let $G$ be a subgroup of $GL(n, \mathbb{Q})$ with finitely many orbits under automorphisms. Is $G$ a virtually soluble group?
 \hfill \raisebox{-1ex}{\sl A.\,Dantas, E.\,de\,Melo}

 \ul

\oai Yes, it is (unverified AI solutions:  Achyuth Jayadevan,  {\it Preprint of 14 September 2026},
\url{https://kourovkanotebookorg.wordpress.com/wp-content/uploads/2026/09/kourovka-21-40-paper.pdf}; \ M.\,Kinyon, J.\,Urban, {\it Preprint of 17 September 2026},    \url{https://kourovkanotebookorg.wordpress.com/wp-content/uploads/2026/09/kourovka-mathematics.pdf}).
 \emp 

 \bmp \textbf{21.41.} A group is said to be \emph{self-similar} if it admits a faithful state-closed representation by automorphisms  of a regular one-rooted $m$-tree for some $m$.
  Can a torsion-free finitely presented metabelian group which is self-similar
  contain  a subgroup isomorphic to the restricted wreath product $H = \mathbb{Z} \wr\mathbb{Z}$?

\makebox[15pt][r]{}It is known that $\mathbb{Z} \wr \mathbb{Z}$  itself is
self-similar
(A.\,C.\,Dantas, T.\,M.\,G.\,Santos, S.\,N.\,Sidki, \emph{J.~Algebra}, {\bf 567} (2021), 564--581).
\hfill \raisebox{-1ex}{\sl A. \,Dantas, S. \,Sidki}\emp

 \bmp \textbf{21.42.} Let $T_{d, c}$ denote the class of $d$-generated, torsion-free nilpotent groups having class $c$.  It is known that $T_{d, 2}$-groups are self-similar for all $d$,  that  $T_{2, 3}$-groups are self-similar, and that there are $T_{4, 3}$-groups that are not self-similar (A.\,Berlatto, T.\,Santos, {\it Preprint}, 2025, \url{https://arxiv.org/abs/2509.16947}). Are there $T_{3,3}$-groups that are not self-similar?
 \hfill \raisebox{-1ex}{\sl A.\,Dantas, S.\,Sidki}
 \emp

\bmp \textbf{\zv 21.43.}
{\it Conjecture}: Suppose that for a fixed positive integer $k$ at least half of the elements of a finite group $G$ have order $k$. Then $G$ is solvable.
 \hf\raisebox{-1ex}{\sl M.\,Deaconescu}

 \ul

 \otv This is not always true. In a direct product of $N$ copies of $A_5$ there are $60^N$ elements, while the number of elements of order exactly 30 is $60^N - 45^N - 40^N - 36^N + 25^N + 21^N + 16^N - 1$. As $N\to \infty$, the ratio of elements of order 30 converges to~1. (L. Tae Young, \emph{Letter of 9 January 2026};  see also
 (R.\,McCulloch, L.\,Tae~Young, {\it Preprint}, 2026, \url{https://arxiv.org/abs/2602.19340}).)
\emp

\bmp \textbf{\zv 21.44.}
 Let $W_n = A_5 \wr \cdots \wr A_5$ be the $n$-times iterated permutational wreath product of $A_5$ in its natural action (so $W_n$ acts on $5^n$ points), and let $W = \displaystyle \lim_{\longleftarrow} W_n$ be the inverse limit (infinite iterated wreath product of $A_5$). Does $W$ contain a finitely generated dense subgroup of subexponential growth?
\hf \raisebox{-1ex}{{\sl S.\,Eberhard}}

\ul

\otv Yes, it does (AI solution confirmed by S.\,Eberhard: \  Nilradical,  {\it Preprint}, 2026,

\url{https://kourovkanotebookorg.wordpress.com/wp-content/uploads/2026/09/nilradical-kourovka-21-44.pdf}).
\emp

\bmp \textbf{21.45.} (Well-known problem).
Does there exist a finitely presented (infinite) simple group requiring more than two generators? Cf. 6.44 in Archive.
\hf\raisebox{-1ex}{\sl F.\,Fournier-Facio}
\emp

\bmp \textbf{21.46.} (Well-known problem).
Does there exist a finitely presented (infinite) simple group of finite cohomological dimension greater than $2$?
\hf\raisebox{-1ex}{\sl F.\,Fournier-Facio}
\emp

\bmp \textbf{21.47.} (Well-known problem).
Does there exist a finitely presented group $G$ such that $G \cong G \times H$ for some non-trivial group $H$?

 \makebox[15pt][r]{}The first finitely generated example was constructed in (J.\,M.\,Tyrer Jones, {\it J.~Austral. Math. Soc.}, {\bf 17} (1974), 174--196).  A finitely presented group that surjects onto its own direct square was  constructed in (G.\,Baumslag, C.\,F.\,Miller, III, {\it Bull. London Math. Soc.},
 {\bf 20}, no.\,3 (1988), 239--244).
 \hf\raisebox{-1ex}{\sl F.\,Fournier-Facio}
 \emp

\bmp \textbf{21.48.} A \emph{quasimorphism} on a group $G$ is a function $f : G \to  {\Bbb R}$ such that the quantity $\ds\sup_{g,h}
|f(g) + f(h) - f(gh)|$ is finite. A quasimorphism is \emph{homogeneous} if it
restricts to a homomorphism on every cyclic subgroup of $G$.
Let $G$ be a group admitting an unbounded homogeneous quasimorphism $G \to\nobreak \mathbb{R}$ that is not a homomorphism. Must $G$ contain a non-abelian free subgroup?
\hf\raisebox{-1ex}{\sl F.\,Fournier-Facio}
\emp

\bmp \textbf{21.49.}
An isometric action of a group $G$ on a metric space $S$ is called \emph{acylindrical} if for every $\e  > 0$ there exist
$R, N > 0$ such that for every two points $x, y$ with $d(x, y) \geq  R$, there are at most $N$ elements $g\in  G$ satisfying $d(x, gx) \leq  \e$ and $d(y, gy) \leq \e$. A  group
 is said to be \emph{acylindrically hyperbolic} if it is not virtually cyclic and  admits an  acylindrical action on a hyperbolic space with
unbounded orbits. Is the automorphism group of a finitely generated acylindrically hyperbolic group  also acylindrically hyperbolic?
\hf\raisebox{-1ex}{\sl A.\,Genevois}
\emp

\bmp \textbf{21.50.}  Does every finite $3$-group $T$ have a nontrivial characteristic subgroup $C$ such that if $T$ is a Sylow $3$-subgroup of a finite group $G$, then   $T\cap G'=T\cap H'$, where $H=N_G(C)$?

 \makebox[15pt][r]{}Such a characteristic subgroup is known to exist in $p$-groups for $p\geq 5$ (G.\,Glauberman, {\it Math. Z.}, {\bf 117} (1970), 46--56), and for $p=3$ there are two characteristic subgroups $K_1,K_2$ such that $T\cap G'=(T\cap H_1')(S\cap H_2')$, where $H_i=N_G(K_i)$  (G.\,Glauberman, {\it  J.~Algebra}, {\bf 648} (2024),  62--86). The group $S_4$ shows that no such characteristic subgroups can be found in some Sylow $2$-subgroups.
 \hf\raisebox{-1ex}{\sl G.\,Glauberman}
\emp

\bmp \textbf{21.51.} 
Let $p$ be a prime, and $P$ a finite $p$-group.

\makebox[25pt][r]{(a)} Suppose that $P$ has an abelian subgroup of  order $p^n$. For which $n$ does $P$ necessarily have a normal abelian subgroup of order~$p^n$?

\makebox[25pt][r]{(b)} Suppose that $P$ has an elementary abelian subgroup of  order~$p^n$. For which $n$ does $P$ necessarily have a normal elementary abelian subgroup of order $p^n$?

\makebox[15pt][r]{}It is easy to see that for $p=2$, the answer to (b) is ``yes'' only for $n=1$. The answer to both questions is ``yes'' for  $n <  (p+2)/2$ (G.\,G.\,Glau\-berman, {\it J.~Algebra}, {\bf 319}, no.\,2 (2008),  800--805), as well as for $n\leq 5$ when  $p\ne 2$ (M.\,Konvisser,  D.\,Jonah, {\it J.~Algebra}, {\bf 34} (1975), 309--330). The answer to both questions is ``no'' for  $n\geq  (p+9)/2$ when $p\geq 5$ (for $p\geq 7$ due to G.\,Glauberman, {\it  Contemp. Math.}, {\bf 524} (2010), 61--65;  for $p=3, 5$  due to  Ya.\,G.\,Berkovich, {\it J.~Algebra}, {\bf 248}, no.\,2 (2002),  472--553).  Thus, the only open cases for $p\geq 5$ are  $n= 6$ for $p =5$ , and $n=(p+3)/2$, \ $(p+5)/2$,  \ $(p+7)/2$  for $p > 5$.

\makebox[15pt][r]{}{\it Comment of 2026}: The answer is ``yes'' for $p=5$ and $n=6$ for both parts (a) and~(b) (L.\,Gong, J.\,Chen, B.\,Li, {\it Preprint}, 2026, \url{https://arxiv.org/abs/2611.34205}).
 \hfill \raisebox{-1ex}{\sl G.\,Glauberman}
\emp

\bmp \textbf{21.52.}
Let $L$ be a finite non-abelian simple group, and let $D$ be a conjugacy class of involutions in $L$. Consider the complete graph $\Gamma$ with vertex set $D$. Define an equivalence relation $\sim$ (graph coloring) on the set of edges as follows: $(a,b) \sim (c,d)$ if and only if  $|ab| = |cd|$. An automorphism of the coloured graph $\Gamma$ is a permutation $\tau \in \nobreak S_D$ such that $(a,b)\sim (a^\tau, b^\tau)$ for every edge $(a,b)$. Is it true that the automorphism group of $\Gamma$ is a subgroup of $\operatorname{Aut}(L)$?
\hf \raisebox{-1ex}{\sl I.\,B.\,Gorshkov}
\emp

\bmp \textbf{\zv 21.53.}
In the notation of 21.52, let $\operatorname{Aut}_t(\Gamma)$ be the set of  permutations $\tau\in \nobreak S_D$ such that  $(a,b)\sim(a^{\tau},b^{\tau})$ whenever $|ab|=t$  for  $a,b\in D$. Clearly, $\operatorname{Aut}(\Gamma)=\bigcap _t\operatorname{Aut}_t(\Gamma)$.
Is it true that for every finite simple group $G$ we have $\operatorname{Aut}(\Gamma)= \operatorname{Aut}_2(\Gamma)\cap \operatorname{Aut}_p(\Gamma)$, where $\{2,p\}$ are the two minimal prime divisors of~$|G|$?
\hf \raisebox{-1ex}{\sl I.\,B.\,Gorshkov}

\ul

\otv No, not for every (AI solution confirmed by I.\,B.\,Gorshkov: \

Nilradical,  {\it Preprint}, 2026,
\url{https://kourovkanotebookorg.wordpress.com/wp-content/uploads/2026/09/nilradical-kourovka-21-53.pdf}).
\emp 

\bmp \textbf{21.54.}
Let $G$ be a finite soluble \textit{group with triality}, which means that $G$ admits a~group of automorphisms $S$ isomorphic to the symmetric group of degree~$3$ given by the presentation $S=\langle \sigma, \rho \mid \sigma ^2=\rho ^3=1;\;\, \sigma\rho\sigma=\rho^2\rangle $ such that
 $m\cdot m^{\rho}\cdot m^{\rho ^2}=1$ for all $m$ in the set of commutators $ M(G):=\{[g,\sigma]\mid g\in G\}$.

  \makebox[15pt][r]{}Suppose in addition that $G=[G,S]$, the group $G$ is generated by $d$ elements of $M(G)$ and their images under $S$, and $x^n=1$ for all $x\in M(G)$. Is it true that the Fitting height of $G$ is bounded in terms of $d$ and $n$?

  \makebox[15pt][r]{}An affirmative answer would provide a reduction of the analogue of the Restricted Burnside Problem for Moufang loops to the nilpotent case.

\hfill \raisebox{-1ex}{\sl A.\,N.\,Grishkov, A.\,V.\,Zavarnitsine}
\emp

\bmp \textbf{\ai 21.55.}   
Let $q$ be  a power of a prime $p$, and let $m_n(q)$ be the maximum  $p$-length of  $p$-solvable subgroups of~$GL(n,q)$. Is it true that $\displaystyle\lim_{n\to\infty} m_n(q)/\log_2n=1$?
\hfill \raisebox{-1ex}{\sl \.{I}.\,G\"{u}lo\u{g}lu}

\ul

\oai No, it is not true (unverified AI solution: B.\,G.\,Di Siena, {\it Preprint}, 2026, \url{https://www.proofforum.org/abs/PF:2026.000033}).
\emp

\bmp \textbf{21.56.} Let $\ell (X)$ denote the composition length of a finite group $X$. Let $A$ be a finite nilpotent group acting by automorphisms on a finite soluble group~$G$. Let $c(G,A)$ be the number of trivial $A$-modules in a given $A$-composition series of $G$. (Note that $c(G,A)=\ell (C_ {G} (A)) $ if $(|A|,|G|)=1$.)

\makebox[15pt][r]{}{\it Conjecture}: there are absolute constants $C_1$ and $C_2$ such that the Fitting height of $G$ is at most $C_1\ell(A) +C_2c(G,A)$.
  \hfill \raisebox{-1ex}{\sl \.I.\,G\"ulo\v{g}lu}
\emp

\bmp \textbf{\ai 21.57.} 
Let $\mathfrak{X}$ be a non-empty class of finite groups of odd order closed under taking subgroups, homomorphic images, and extensions. Let  $H$ be an $\mathfrak{X}$-maximal subgroup of a finite group $G$, and $N$ a normal subgroup of $G$. Must $H\cap N$ be an $\mathfrak{X}$-maximal subgroup of $N$?
\hf \raisebox{-1ex}{\sl  W.\,Guo, D.\,O.\,Revin}

\ul

\oai No, it need not (unverified AI solution: B.\,G.\,Di~Siena, {\it Preprint}, 2026, \url{https://kourovkanotebookorg.wordpress.com/wp-content/uploads/2026/10/21-57.pdf}).
\emp

\bmp \textbf{\zv 21.58.}
We say that a product $XY=\{xy\mid x\in X,\;y\in Y\}$ of two subsets $X,Y$ of a group $G$ is \emph{direct} if for every $z\in XY$ there are unique $x\in X$, $y\in Y$ such that $z=xy$. Is there an infinite group $G$ such that every subset $A \subseteq G$ satisfies the following property: all the maximal subsets $B$ for which the product $AB$ is direct have the same cardinality?

\makebox[15pt][r]{}Note that for checking the property for a given infinite group $G$, it suffices to consider only those subsets $A \subseteq G$ for which $|A| = |G \setminus A|$. Indeed, the property is equivalent to $A^{-1}A \cap BB^{-1} = \{1\}$ and $A^{-1}AB = G$, and these imply $|G| = |A||B|$, since $G$ is infinite. Now, if $|A| < |G \setminus A|$, then   $|A| < |G | $, and so   $|B| = |G|$; and if
 $|A| > |G \setminus A|$, then $A^{-1}A = G$, and so $|B| = 1$, for all $B$ satisfying the property.

\hf \raisebox{-1ex}{\sl  M.\,H.\,Hooshmand}

\ul

\otv No, there are no such groups (M.\,I.\,Kabenyuk, {\it Preprint}, February 2026, \url{https://arxiv.org/abs/2602.22876}; \ M.\,H.\,Hooshmand,
{\it Preprint}, April 2026, \url{https://arxiv.org/abs/2604.08724}, Remark 1.16).
\emp

 \bmp \textbf{21.59.}
 For a finite group $G$,
let $\chi _1(G)$ denote the totality of the degrees of all
irre\-du\-cible complex characters of $G$ with allowance for their
multiplicities.  Suppose that $H$ is a finite group with $\chi_1(H)=\chi_1(G)$.

\makebox[25pt][r]{a)}  If  $G$ is an almost simple group, must  $H$ be  isomorphic to $G$?

\makebox[25pt][r]{\zva b)} If  $G$ is a quasisimple group, must  $H$ be  isomorphic to $G$?

\makebox[15pt][r]{}This is true if $G$ is a simple group (see~11.8(a) in Archive).

 \hfill \raisebox{-1ex}{\sl A.\,Iranmanesh, F. Shirjian}

 \ul

 \otv b) Yes, it must, by Theorem~B in (C.\,Bessenrodt, H.\,N.\,Nguyen, J.\,B.\,Olsson, H.\,P.\,Tong-Viet, {\it Algebra Number Theory}, {\bf 9}, no.\,3 (2015), 601--628).  \emp

\bmp \textbf{\zv 21.60.} Let $G$ be a finite group, ${\Bbb Z}_{(p)}$ the
localization at $p$, and $\F_p$ the field of $p$ elements. Let $X$ be the class of $\F_pG$-modules obtained by reduction of simple $\Q G$-modules. Is it true that $\Z_{(p)}G$ is semiperfect if and only if each projective indecomposable $\F_pG$-module can be written as an $\N$-linear combination of modules in $X$ inside the Grothendieck group of $\F_pG$?

 \makebox[15pt][r]{}The ``only if'' direction is proved, and the affirmative answer is obtained if $p$ does not divide $|G|$ (D.\,Johnston, D.\,Rumynin, {\it J.~Algebra}, {\bf 687}, no.\,1 (2026),  776--791).

\hfill \raisebox{-1ex}{\sl D.\,Johnston, D.\,Rumynin}

\ul

\otv No, not always (AI solution confirmed by D.\,Rumynin: \ M.\,Kinyon, J.\,Urban, {\it Preprint}, 2026,
 \url{https://kourovkanotebookorg.wordpress.com/wp-content/uploads/2026/09/kourovka-mathematics.pdf}).
\emp

\bmp \textbf{\ai 21.61.} 
For a fixed (finitely generated free)-by-cyclic group  $G=F_n\rtimes \mathbb Z$, is there an algorithm that, given a finite subset $S$ of $G$, finds a finite presentation for the subgroup $H=\langle S\rangle$? Cf. 4.8 in Archive.
\hfill \raisebox{-1ex}{\sl I.\,Kapovich}

\ul

\oai No, there is no such algorithm  (unverified AI solution: T.\,Stengel, {\it Preprint}, 2026, \url{https://doi.org/10.5281/zenodo.24557836}).
\emp 

\bmp \textbf{21.62.}
 Is the uniform subgroup membership problem decidable for (finitely generated free)-by-cyclic groups? That is, for a fixed  group $G=F_n\rtimes \mathbb Z$, is there an algorithm that, given elements $w, h_1, \dots , h_k\in  G$, decides whether or not $w$ belongs to the subgroup $H=\langle h_1,\dots ,h_k\rangle$?
\hfill \raisebox{-1ex}{\sl I.\,Kapovich}
\emp

\bmp \textbf{21.63.} (E.\,Zelmanov). Let $F$ be a field of characteristic $p>0$, and let $\Gamma$
 be the principal congruence subgroup of $\mathop{\rm Aut}(F[x_1,\dots , x_n])$ consisting of all automorphisms that send each variable $x_i$ to $x_i$ modulo terms of higher degree. Then $\Gamma$ is a residually-$p$ group. Does $\Gamma$ satisfy a pro-$p$ identity?
 \hf\raisebox{-1ex}{\sl E.\,I.\,Khukhro}
\emp

\bmp \textbf{21.64.} Is it true that if a normal subgroup $A$ of a Sylow $p$-subgroup of a $p$-soluble finite group $G$ has exponent $p^e$, then the normal closure of  $A$ in $G$  has $(p,e)$-bounded (or even $e$-bounded) $p$-length?
 \hf\raisebox{-1ex}{\sl E.\,I.\,Khukhro}
\emp

\bmp \textbf{21.65.}
Suppose that  $\varphi $ is an automorphism of a finite soluble group $G$. Must $G$ contain a subgroup of index bounded in terms of $|\varphi |$ and $|C_G(\varphi )|$ whose Fitting height is bounded

\makebox[25pt][r]{(a)} in terms of $|\varphi |$?

\makebox[25pt][r]{(b)} or even in terms of the composition length of $\langle \varphi \rangle$?

\makebox[15pt][r]{}An affirmative answer to part~(b) is known when $|\varphi |$ is a prime power (B.\,Hartley--V.\,Turau),
or when $(|G|,|\varphi |)=1$ (A.\,Turull, B.\,Hartley--I.\,M.\,Isaacs).
Also cf. 19.43.

\hfill \raisebox{-1ex}{\sl E.\,I.\,Khukhro}
\emp

\bmp \textbf{21.66.}
Suppose that  $A$ is a nilpotent group of automorphisms of a finite soluble group $G$. Is the Fitting height of $G$ bounded in terms of $|A|$ and $|C_G(A)|$?

\makebox[15pt][r]{}An affirmative answer is known when $(|G|,|A|)=1$ (J.\,G.\,Thompson, even for soluble $A$, with improved bounds in subsequent papers of H.\,Kurzweil, A.\,Turull, B.\,Hartley--I.\,M.\,Isaacs), when  $A$ is cyclic (see 19.43), or when $C_G(A)=1$ (E.\,C.\,Dade). Note that for any non-nilpotent finite group $A$ there are finite soluble groups $G$ of unbounded Fitting height with $C_G(A)=1$ (S.\,D.\,Bell--B.\,Hartley).
\hfill \raisebox{-1ex}{\sl E.\,I.\,Khukhro}
\emp

\bmp \textbf{21.67.}
Suppose that  $\varphi $ is an automorphism of a finite soluble group $G$, and let $r$ be the (Pr\"ufer) rank of the fixed-point subgroup $C_G(\varphi )$. Is the Fitting height of $G$ bounded terms of $|\varphi |$ and $r$?

\makebox[15pt][r]{}This is known to be true in the case where $|\varphi |$ is a product of at most two prime powers (B.\,Hartley, {\it Preprint}, 1994,
\url{https://kourovkanotebookorg.wordpress.com/wp-content/uploads/2025/08/hartley94-prepr-mims.pdf}). An affirmative answer to this question would imply an affirmative answer to 13.8(b).
\hfill \raisebox{-1ex}{\sl E.\,I.\,Khukhro}
\emp

\bmp \textbf{\zv 21.68.}
A finite group $G$  is said to be \emph{semi-abelian} if it has a sequence of subgroups $1=G_0\leq G_1\leq \cdots \leq G_n=G$ such that for every $i$ the subgroup $G_{i+1}$ is isomorphic to a quotient of a semidirect product $A_i\rtimes G_i$ for some abelian group $A_i$.

\makebox[15pt][r]{}{\it Conjecture}:  Semi-abelian finite groups are monomial.
 \hf\raisebox{0ex}{\sl M.\,Kida}

\ul

\otv Conjecture is refuted in (A.\,Rizzoli, {\it Preprint}, 2026,

\url{https://kourovkanotebookorg.wordpress.com/wp-content/uploads/2026/09/kourovka_21_68.pdf}).
\emp

\bmp \textbf{21.69.}
 Is there an algorithm deciding if a given one-relator group is hyperbolic?

 \hf\raisebox{-1ex}{\sl D.\,Kielak}
\emp

\bmp \textbf{21.70.}
A group $G$ is called an \emph{orientable Poincar\'e duality group}  of dimension $n$ over a ring $R$ if it is of type $\mathtt{FP}$ over $R$ and   $\mathrm{H}^i(G;RG) =  0$ for $i\ne n$, while
$\mathrm{H}^n(G;RG) =  R$ as an $RG$-module, where the action on $R$ is trivial. (Note that $G$ is not required to be finitely presented.) 
If $G$ is an orientable Poincar\'e duality group of dimension $n$ over all fields, is it an orientable Poincar\'e duality group over the integers?
\hfill \raisebox{-1ex}{\sl D.\,Kielak}
\emp

\bmp \textbf{21.71.}
For a ring $R$, we say that a group $G$ is of type $\mathtt{FL}(R)$ if the trivial $RG$-module $R$ admits a finite resolution by finitely generated free modules. If $R$ is a field, we define the Euler characteristic of $G$ over $R$ to be the alternating sum of $R$-ranks of homology groups $\mathrm{H}_i(G;R)$.

 \makebox[15pt][r]{}Does there exist a group of type $\mathtt{FL}(\mathbb F_i)$ for two fields $\mathbb F_1$ and $\mathbb F_2$ such that the Euler characteristics of the group over the fields $\mathbb F_1$ and $\mathbb F_2$ differ?
 \hf\raisebox{-1ex}{\sl D.\,Kielak}
\emp

\bmp \textbf{21.72.}
We say that a group is a \emph{Tarski monster} if it is finitely
generated, not cyclic, and all of its proper non-trivial subgroups are
isomorphic to each other.

 \makebox[25pt][r]{a)} Do there exist amenable torsion-free Tarski monsters?

 \makebox[25pt][r]{b)} Do there exist amenable Tarski monsters of prime exponent?
\hf \raisebox{-1ex}{{\sl D.\,Kielak}}
\emp

\bmp \textbf{21.73.}
Is the conjugacy problem in \({\rm CT}(\mathbb{Z})\) algorithmically decidable?

\makebox[15pt][r]{}See the definition of \({\rm CT}(\mathbb{Z})\) in 17.57.
 \hfill \raisebox{0ex}{\sl S.\,Kohl}
\emp

\bmp \textbf{21.74.}
Is it algorithmically decidable whether a given element \(g \in {\rm CT}(\mathbb{Z})\)

\makebox[25pt][r]{(a)} permutes a nontrivial partition of \(\mathbb{Z}\) into residue classes?

\makebox[25pt][r]{(b)} has only finite cycles?

  \makebox[25pt][r]{(c)} has no finite cycles?
 \hfill \raisebox{0ex}{\sl S.\,Kohl}
\emp

\bmp \textbf{\zv 21.75.}
Given two distinct sets \(\mathcal{P}_1\) and \(\mathcal{P}_2\) of odd primes none
of which is a subset of the other, is it true that
\(\langle {\rm CT}_{\mathcal{P}_1}(\mathbb{Z}), {\rm CT}_{\mathcal{P}_2}(\mathbb{Z}) \rangle \lneq
{\rm CT}_{\mathcal{P}_1 \cup \mathcal{P}_2}(\mathbb{Z})?\)

\makebox[15pt][r]{}See the definition of \({\rm CT}_{\mathcal{P}}(\mathbb{Z})\) in 17.60.
 \hfill \raisebox{0ex}{\sl S.\,Kohl}

 \ul

 \otv No, in fact, equality holds (A.\,L.\,Iskra, {\it Preprint}, 2026,
 \url{http://arxiv.org/abs/2610.04274}). A simplified solution is in (S.\,Kohl, {\it Preprint}, 2026,
 \url{ https://kourovkanotebookorg.wordpress.com/wp-content/uploads/2026/09/21.75-kohl.pdf}).
\emp

\bmp \textbf{\zv 21.76.}
Let $\sigma=(\sigma_{ij})$, $1\leq i\neq j\leq n$, be an irreducible elementary net  (carpet) of order $n\geq 3$ over a field $K$ (see 19.48). The
net $\sigma$ is said to be \textit{closed} if the elementary net
subgroup $E(\sigma)$ does not contain new elementary transvections. The net $\sigma$ is said to be \textit{completable} if its diagonal can be supplemented with subgroups  to a complete net. Completable  elementary nets are closed. It is known that over fields of characteristic $0$ and $2$ there exist irreducible closed elementary nets that are not completable (V.\,A.\,Koibaev, {\it Trudy Inst. Mat. Mekh. Ural Div. Ross. Akad. Nauk}, \textbf{17}, no.\,4 (2011), 134--141 (Russian); \ V.\,A.\,Koibaev, {\it
Siberian Math.~J.}, \textbf{62}, no.\,2 (2021),   262--266).

 \makebox[15pt][r]{}Do there exist irreducible closed elementary nets of order $n\geq 3$ over a field of odd characteristic that are not completable?
 \hfill \raisebox{-1ex}{\sl V.\,A.\,Koibaev}

 \ul

\otv Yes, such nets do exist (Ya.\,N.\,Nuzhin, {\it  Siberian Math. J.}, {\bf 67},  no.\,4 (2026), 840--845).
\emp

\bmp \textbf{21.77.}
Let $d$  be an integer that is not divisible by $n$-th powers of primes, let  $x^n-d$ be an irreducible polynomial over $\mathbb{Q}$,  let $\theta=\sqrt[n]{d}$,  and let   $K=\mathbb{Q}(\theta)$ be the radical extension of degree $n$  of the field $\mathbb{Q}$. The multiplicative group  $K^\ast$ of the field $K$ is canonically embedded into the group   $Aut_{\mathbb{Q}}(K)$ of all invertible $\mathbb{Q}$-linear mappings of the $\mathbb{Q}$-space~$K$; let $T$ be the image of $K^\ast$ under this embedding. In the natural basis
$1, \theta, \theta^2,  \ldots , \theta^{n-1}$  of the $\mathbb{Q}$-space $K$ the   group $Aut_{\mathbb{Q}}(K)$ corresponds to $G=GL(n, \mathbb{Q})$, and the subgroup $T$  to a subgroup $T(d)$ (unsplit maximal torus).
Every subgroup $H$ of $G$ containing $T(d)$ and some one-dimensional transformation is rich in elementary transvections (V.\,A.\,Koibaev,  {\it St. Petersbg. Math.~J.}, {\bf  21}, no.\,5 (2010), 731--742) and thus defines a net $\sigma=\sigma(H)$ (V.\,A.\,Koibaev, A.\,V.\,Shilov, {\it J.~Math. Sci.  New York}\, {\bf 171}, no.\, 3 (2010), 380--385). Let $E(\sigma)$ denote the subgroup generated by all transvections in the net group $G(\sigma)$. Is it true that
$H\leq N_{G}(E(\sigma))$?

 \makebox[15pt][r]{}This inclusion was proved in the case $n=2$ (V.\,A.\,Koibaev, {\it Dokl. Math.}, {\bf 41}, no.\,3 (1990), 414--416).
\hfill \raisebox{-1ex}{\sl V.\,A.\,Koibaev}
  \emp

\bmp \textbf{21.78.}
Let $p$ be a prime and let $G$ be a pro-$p$ group. Suppose that all of the (continuous Galois) cohomology groups $H^n(G, \mathbb F_p)$ of $G$ with coefficients in the field of $p$ elements are finite. Does it necessarily follow that the cohomology ring $H^*(G,\mathbb F_p)$ is finitely generated?

\makebox[15pt][r]{}The answer is known to be `yes' if $G$ is abelian-by-($p$-adic analytic), as follows from (J.\,King, {\it  Commun. Algebra}, {\bf 27}, no.\,10 (1999), 4969--4991).
\hfill \raisebox{-1ex}{\sl P.\,Kropholler}
\emp

\bmp \textbf{21.79.}  Let $G$ be a finitely generated group with a fixed finite generating set $S$ and the corresponding word metric $L_S(*)$. An element $g$ is said to be \emph{distorted} in $G$ if
$ {L_S(g^n)}/{n}\to 0$ as $n\to \infty$; this notion is independent of the choice of the generating set~$S$. For any, not necessarily finitely generated, group $H$,  an element $g\in H$ is said to be \emph{distorted} if there is a finitely generated subgroup $G$ of $H$ containing $g$ in which $g$ is distorted. Do there exist finitely generated left-orderable groups in which every nontrivial element is distorted?

\makebox[15pt][r]{}Note that it is straightforward to construct countable (not finitely generated) left-orderable groups with this property using $HNN$-extensions and applying results of V.\,V.\,Bludov and A.\,M.\,W.\,Glass.
\hfill \raisebox{-1ex}{\sl Y.\,Lodha,  A.\,Navas}

\emp

\bmp \textbf{21.80.}
Do there exist finitely generated left-orderable groups with only one nontrivial conjugacy class?

\makebox[15pt][r]{}A positive answer to this question implies a positive answer to 21.78. Note that D.\,Osin constructed torsion-free finitely generated groups with only one nontrivial conjugacy class; see 9.10 in Archive.
\hfill \raisebox{-1ex}{\sl Y.\,Lodha,  A.\,Navas}
\emp

\bmp \textbf{21.81.} (J.\,Wiegold).
Let $\Gamma$ be a finite simple group and let $\mathcal{N}_n(\Gamma)$ denote the set of normal subgroups of the free group $F_n$ of rank $n$ whose quotient is isomorphic to $\Gamma$.

 \makebox[15pt][r]{}{\it Conjecture}: $\mathrm{Aut}(F_n)$ acts transitively on $\mathcal{N}_n(\Gamma)$ for $n \geq 3$.

 \makebox[15pt][r]{}This is not true for $n = 2$ (B.\,H.\,Neumann, H.\,Neumann, {\it Math. Nachr.}, {\bf 4} (1951),  106--125).
\hf\raisebox{-1ex}{\sl A.\,Lubotzky}
\emp

\bmp \textbf{21.82.}  {\it Conjecture:\/} For $n \geq 3$, there are no finite simple characteristic quotients of the free group~$F_n$.

 \makebox[15pt][r]{}This is not true for $n = 2$ (W.\,Y.\,Chen, A.\,Lubotzky, P.\,H.\,Tiep, {\it to appear in Comment. Math. Helvetici}, 2025).
\hf\raisebox{-1ex}{\sl A.\,Lubotzky}
\emp

\bmp \textbf{21.83.}
The function $d_n (\sigma , \tau) =(1/n) \cdot  |\{x\in \{1,\dots,n\}\mid \sigma  (x) \ne \tau  (x)\}|$ is a distance  on  the symmetric group~$S_n$. For a finitely generated group $G$, an \emph{almost-homomorphism} is a sequence of set-theoretic maps $f_n : G \to  S_n$ satisfying
$d_n (f_n(g)f_n(h), f_n(gh)) \to 0$ as  $n\to\infty$ for all $g, h\in G$.
An almost-homomorphism $\{f_n\}$ is said to be \emph{close to a homomorphism} if there is a sequence
of group homomorphisms $\rho _n : G \to S_n$ such that $d_n (\rho _n(g), f_n(g)) \to 0$ as  $n\to\infty$ for all $g\in G$.
The group $G$ is said to be \emph{permutation stable} if every
almost-homomorphism of $G$ is close to a homomorphism.

 \makebox[15pt][r]{}{\it Conjecture}: Metabelian groups are permutation-stable.
 \hf\raisebox{-0ex}{\sl A.\,Lubotzky}
\emp

\bmp \textbf{21.84.}
For $\sigma\in S_n$ and $\tau\in S_m$, where $n\leq  m$, let
\vspace{-2ex}
$$d^{\rm flex}_n (\sigma , \tau) =(1/n) \cdot  (|\{x\in \{1,\dots,n\}\mid \sigma  (x) \ne \tau  (x)\}| + (m - n)).\vspace{-2ex}$$
An almost-homomorphism $\{f_n\}$ is said to be \emph{flexibly close to a homomorphism} if there is a sequence
of group homomorphisms $\rho _n : G \to S_{m_n}$ with $n\leq m_n$ such that $d_n^{\rm flex} (\rho _n(g), f_n(g)) \to 0$ as  $n\to\infty$ for all $g\in G$.
The group $G$ is said to be \emph{ flexibly permutation stable} if every
almost-homomorphism of $G$ is  flexibly close to a~homomorphism.

 \makebox[15pt][r]{}Is $SL_n(\Z )$ flexibly permutation-stable?
 \hf\raisebox{0ex}{\sl A.\,Lubotzky}
\emp

\bmp \textbf{21.85.}
Is a flexibly permutation-stable group always permutation-stable?

 \hf\raisebox{-1ex}{\sl A.\,Lubotzky}
\emp

\bmp \textbf{\zv 21.86.} (M.\,Gromov, B.\,Weiss).
A group $G$ is said to be \emph{sofic} if for every finite set
$F\subseteq G$ containing $1$ and every $\e > 0$ there exist $n\in \N$ and a map
$\varphi  : F \to S_n$ such that $\varphi (1) = 1$, \ $d(\varphi (gh), \varphi(g)\varphi(h)) < \e $ for all $g, h$
such that $gh \in F$, \  $\varphi (g)$  does not have fixed points for every
$g\in  F\setminus \{1\}$.

\makebox[15pt][r]{}Is every group sofic?
 \hf\raisebox{0ex}{\sl A.\,Lubotzky}

\ul

\otv No, not every (AI solution by OpenAI:  \url{https://cdn.openai.com/pdf/ten-proofs-oai.pdf}).  More examples were found in (G.\,Kun, A.\,Thom, {\it Preprint}, 2026, \url{https://arxiv.org/abs/2608.06222}).
\emp

\bmp \textbf{\zv 21.87.}
Assume that a finite group $G$ has a family of $d$-generator subgroups
whose indices have no common divisor. Is it true that $G$ can be
generated by $d + 1$ elements?

\makebox[15pt][r]{}The answer is positive if $G$ is solvable (L.\,G.\,Kov\'acs,
H.-S.\,Sim, {\it  Indag. Math.}, {\bf 2} (1991), 229--232). For an arbitrary finite group $G$ it is proved  that $G$ can be generated by $d+2$ elements  (A.\,Lucchini, {\it  Commun. Algebra}, {\bf 28}, no.\,4
(2000), 1875--1880).

 \hfill \raisebox{-1ex}{\sl A.\,Lucchini}

 \ul

 \otv Yes, it is true (AI solutions confirmed by A.\,Lucchini: \ J.\,DeCaro,
 {\it Preprint}, July 2026,
 \url{https://kourovkanotebookorg.wordpress.com/wp-content/uploads/2026/08/21-97-decaro-op-0717-preprint.pdf});

 R.\,Sater, {\it Preprint}, August 2026, \url{https://arxiv.org/abs/2608.12432}).
\emp

\bmp \textbf{\zv 21.88.}
Is there a finite non-abelian group $G$ of odd order, with $k(G)$ conjugacy classes, such that  $k(G)/|G|  = 1/17$?
\hf \raisebox{-1ex}{\sl D.\,MacHale}

\ul

\otv No, there are no such groups (B.\,Beyer\,de\,Ryke,  {\it Preprint}, 2026, \url{https://arxiv.org/abs/2608.03003}).
\emp

\bmp
   \textbf{21.89.}
 For $n>39$, is it true that the number of conjugacy classes in the  symmetric group $S_n$ of degree $n$ is never a divisor of the order of $S_n$? In other words, is it true that, for $n>39$,  the number $ p(n)$  of integer partitions of $n$ is never a divisor of $n!$?

\hfill \raisebox{-1ex}{\sl D.\,MacHale}
\emp

\bmp \textbf{21.90.}
Let $\Gamma$ be a graph of diameter $d$. For  $i\in \{1,2,\dots ,d\}$, let $\Gamma_i$ be the graph on the same
vertex set as $\Gamma$ with vertices $u,w$ adjacent in $\Gamma_i$ if and only if  $d_{\Gamma}(u,w)=i$. Does there exist a $Q$-polynomial distance-regular graph $\Gamma$ of diameter 3 such that $\Gamma_2$ and $\Gamma_3$ are strongly regular?
\hfill \raisebox{-1ex}{\sl A.\,A.\,Makhn\"ev}
\emp

\bmp \textbf{21.91.} (W.\,Willems). {\it Conjecture}:
 The sum of squares of the degrees of the irreducible $p$-Brauer characters of a finite
group $G$ is at least the $p'$-part of $|G|$.

\makebox[15pt][r]{}This is known to be true for $p=2$ (G.\,Malle, {\it Adv. Math.}, {\bf 380} (2021), Paper no.\,107609, 15 pp.).
\hf \raisebox{-1ex}{\sl G.\,Malle}
\emp

 \bmp \textbf{\zv 21.92.}
{\it Conjecture}: The number of  irreducible $p$-Brauer characters of a finite group $G$
is bounded above by the maximum of the number of conjugacy classes $k(H)$ in
$p'$-subgroups $H$ of $G$.
\hf \raisebox{-1ex}{\sl G.\,Malle, G.\,Navarro, G.\,Robinson}

\ul

\otv The conjecture is refuted: for $G=PSL_3(4).6$, an almost simple group with simple socle $PSL_3(4)$, and for $p=5$ or $7$, the group $G$ has 23 irreducible $p$-Brauer characters, but $k(H)$ is at most 22 for all $p'$-subgroups $H$  
(H.\,P.\,Tong-Viet, {\it Letter of 28 October 2026}).
\emp

\bmp \textbf{21.93.} (M.\,Herzog, J.\,Sch\"onheim). Let $G$ be a group and let $k \geq 2$. Let $H_1,\dots,H_k$ be subgroups of $G$, and $g_1,\dots,g_k$ elements of $G$ such that the cosets $g_1H_1,\dots,g_kH_k$ form a partition of $G$. Is
it true that $|G:H_i| = |G:H_j|$ for some $i \neq j$?

\makebox[15pt][r]{}This is known to be true for groups with a Sylow tower (M.\,A.\,Berger, A.\,Felzen\-baum, A.\,Fraenkel, {\it  Fund. Math.}, {\bf  128}, no.\,3 (1987), 139--144). Also cf. 20.99.

\hfill \raisebox{-1ex}{\sl L.\,Margolis}
\emp

 \bmp \textbf{21.94.} The Gruenberg--Kegel graph (or the prime graph) $GK(G)$ of a finite group~$G$ is a labelled graph with vertex set consisting of all prime divisors of the order of $G$ in which different vertices $p$ and $q$ are adjacent if and only if $G$ contains an element of order~$pq$. Let $\overline{GK}(G)$ denote the abstract graph obtained from $GK(G)$ by removing all labels.  A finite group $G$ is said to be \emph{recognizable by the isomorphism type of its Gruenberg--Kegel graph} if there are no finite groups $H\not\cong G$ with $\overline{GK}(H)$  isomorphic to $\overline{GK}(G)$.

\makebox[15pt][r]{}Are there infinitely many (pairwise non-isomorphic) finite  groups which are recog\-ni\-zable by the isomorphism type of the Gruenberg--Kegel graph?
 \hfill \raisebox{-1ex}{\sl N.\,V.\,Maslova}
 \emp

 \bmp \textbf{21.95.} Is there an almost simple but not simple group which is recognizable by the isomorphism type of its Gruenberg--Kegel graph?

\makebox[15pt][r]{}Note that if a group $G$ is recognizable by the isomorphism type of its Gruenberg--Kegel graph, then $G$ is recognizable by its Gruenberg--Kegel graph, and therefore $G$ is almost simple (P.\,J.\,Cameron, N.\,V.\,Maslova, {\it  J.~Algebra}, {\bf 607} (2022), 186--213).

 \hfill \raisebox{-1ex}{\sl N.\,V.\,Maslova}
 \emp

  \bmp \textbf{\zv 21.96.}
  Is it true that a periodic group containing an involution is locally finite if the centralizer of every element of even order is locally finite?
  \hfill \raisebox{-1ex}{\sl V.\,D.\,Mazurov}

  \ul

\otv No, not always (AI solution confirmed by V.\,D.\,Mazurov: \  B.\,G.\,Di~Siena, {\it Preprint}, 2026, \url{https://kourovkanotebookorg.wordpress.com/wp-content/uploads/2026/10/21-96.pdf}).
   \emp

\bmp \textbf{\zv 21.97.} 
(M.\,T\u{a}rn\u{a}uceanu). Is it true that for every positive rational number $r$ there exists  a finite group $G$ such that ${|\operatorname{Aut}(G)|}/{|G|}=r$?

\makebox[15pt][r]{}A similar question is answered in the positive for graphs, monoids, partial groups, and posets (R.\,Molinier, {\it Preprint}, 2025, \url{https://arxiv.org/abs/2504.21059}). It is also known that the set $\{|\operatorname{Aut}(G)| / |G|\mid G\text{ is a finite abelian group}\}$ is dense in $[0, +\infty)$ (M.\,T\u{a}rn\u{a}uceanu, {\it Elemente Math.} (2025), \url{https://ems.press/journals/em/articles/14298544}).
\hf \raisebox{-1ex}{\sl R.\,Molinier}

\ul

\otv Yes, it is true (S.\,Sureaux, {\it Preprint}, 2026, \url{https://kourovkanotebookorg.wordpress.com/wp-content/uploads/2026/08/21.97-solution.pdf}); \ later also with abelian groups  $G$ (R.\,McCulloch, {\it J.~Algebra Appl.}, 2026, \url{https://doi.org/10.1142/S0219498828500077}).
\emp

\bmp \textbf{21.98.}
Let $w$ be a multilinear commutator word, and assume that $G$ is a group where the set of $w$-values is covered by finitely many cyclic subgroups. Is it true that the verbal subgroup $w(G)$ is finite-by-cyclic?

 \makebox[15pt][r]{}This is true for lower central words (G.\,Cutolo, C.\,Nicotera,  {\it J. Algebra}, {\bf  324}, no.\,7 (2010), 1616--1624).
\hfill \raisebox{-1ex}{\sl M.\,Morigi}
\emp

\bmp \textbf{\zv 21.99.}
  {\it Conjecture}: If $G$ is a transitive permutation group  on a finite set $\Omega$, then for any distinct $\alpha,\beta$ in $\Omega$ there is an element $g\in G$ with
$\alpha^g=\beta$ whose number of fixed points is different from $1$.
\hfill \raisebox{-1ex}{\sl P.\,M\"uller}

\ul

\otv The conjecture is refuted in  (K.\,Muliarchyk,  {\it Preprint}, 2026, \url{https://kourovkanotebookorg.wordpress.com/wp-content/uploads/2026/09/kourovka_21.99_solution.pdf}).
 \emp

\bmp \textbf{\oai 21.100.} 
Suppose that $A$ and $G$ are finite groups such that $A$ acts coprimely on $G$ by automorphisms.
Let $C = {\bf C}_G(A)$ be the fixed-point subgroup, and let $C^\prime$ denote its derived subgroup.

\makebox[15pt][r]{}Is it true that the number of $A$-invariant irreducible characters $\chi$ of $G$ whose restriction
$\chi_C$ is never zero is exactly $|C/C'|$?
This would follow if one could show that $\chi_C$ is never zero if and only if the
Glauberman--Isaacs correspondent $\chi^{*}$ of $\chi$ is linear.
\hfill \raisebox{-1ex}{\sl G.\,Navarro}

\ul

\oai No, not always (unverified AI solution: E.\,Hou, {\it Preprint}, 2026, \url{https://arxiv.org/abs/2609.16227}).
\emp

\bmp \textbf{21.101.}
Which finite almost simple groups are the automorphism groups of regular polytopes of rank~3? In other words, which  finite almost simple groups are
generated by three involutions two of which commute?

\makebox[15pt][r]{}This question has been answered for finite simple groups; see 7.30 in Archive.

 \hfill \raisebox{-1ex}{\sl Ya.\,N.\,Nuzhin}
\emp

\bmp \textbf{21.102.}
Let $\mathscr V$ be a variety generated by a finite group, and let $f(n)$ be the order of the free group in $\mathscr V$  on $n$ generators. Is it true that the sequence $\sqrt[n]{\log f(n)}$ has a limit as $n \to\infty$, and this limit is an integer?
\hf \raisebox{-1ex}{\sl A.\,Yu.\,Olshanskii}
\emp

\bmp \textbf{21.103.} (V.\,V.\,Uspenskii).
A Hausdorff topological group $G$ is called \emph{minimal} if it does not admit a strictly coarser Hausdorff group topology. A topological group is called \emph{Raikov complete} if its two-sided uniform structure is complete. It is known that a finite direct product of Raikov complete minimal topological groups is again minimal. Is it true that an arbitrary Cartesian product of Raikov complete minimal topological groups remains minimal?

\makebox[15pt][r]{}It is known that an arbitrary Cartesian product of centre-free minimal topological groups is minimal (M.\,Megrelishvili, {\it Topology Appl.}, {\bf 62}, no.\,1  (1995), 1--19).
\hfill \raisebox{-1ex}{\sl D.\,Peng}
\emp

\bmp \textbf{21.104.}
For a group word $w(x_1, \dots, x_n)$  on $n$ letters, define $e_{0}(x_1, \dots, x_n) = x_1$ and $e_{k+1}(x_1, \dots, x_n) = w(e_k(x_1, \dots, x_n), \dots, x_n)$ for all $k \in \mathbb{N}$. A group $G$ is said to satisfy the \emph{Engel type iterated identity $w$} if for all $x_1, \dots, x_n \in G$ there exists $m \in \mathbb{N}$ such that $e_m(x_1, \dots, x_n)=1$.

\makebox[15pt][r]{}{\it Conjecture}: For every non-trivial word $w$, if a finitely generated branch group $G$ (see 15.12) satisfies the iterated identity $w$, then $G$ is a torsion group.

\makebox[15pt][r]{}This is true in the case of  the commutator word $w=[x_1,x_2]$  (G.\,Fern\'andez-Alcober, M.\,Noce, G.\,Tracey, \emph{J.~Algebra}, \textbf{554} (2020), 54--77).
\hf \raisebox{-1ex}{\sl  M.\,Petschick}
\emp

\bmp \textbf{21.105.}  (D.\,Segal). A group word $w$ is said to be \emph{concise} in a class $\mathfrak{C}$ of groups if for every group $G$ in $\mathfrak{C}$ such that the set $G_w$ of word values of $w$ in $G$ is finite, the verbal subgroup $w(G)=\langle G_w\rangle$ is also finite (see also Archive 2.45). Is every word concise in the class of residually finite groups?

\makebox[15pt][r]{}See (D.\,Segal, \emph{Words. Notes on verbal width in groups}, Cambridge Univ. Press, 2009) for a partial solution.
	\hf \raisebox{-1ex}{\sl  M.\,Petschick}
\emp

\bmp \textbf{\zv 21.106.}
A first order formula $\varphi(x)$ in the group language with one free variable is said to be  \emph{concise} in a class $\mathfrak{C}$ of groups if for every group $G$ in $\mathfrak{C}$ such that the set $G_\varphi$ of elements in $G$ satisfying $\varphi$ is finite, the subgroup $\varphi(G)$ generated by $G_\varphi$ is also finite  (cf. 21.105 and  Archive 2.45).  Is every formula with one free variable concise in the class of residually finite groups?
\hf \raisebox{-1ex}{\sl  M.\,Petschick}

\ul

\otv No, not every  (AI solution confirmed by M.\,Petschick: \ Achyuth Jayadevan, {\it Preprint of 12 September 2026},
 \url{https://zenodo.org/records/22723520}; \ unverified AI solution: M.\,Kinyon, J.\,Urban, {\it Preprint of 17 September 2026},
 \url{https://kourovkanotebookorg.wordpress.com/wp-content/uploads/2026/09/kourovka-mathematics.pdf}).
\emp

\bmp \textbf{21.107.}
 A sequence $\{ F_n\}$
of pairwise disjoint finite subsets of a topological group is
called {\it expansive\/} if for every open subset $U$ there is a
number $m$ such that $ F_n\cap U\ne \varnothing$ for all $n>m$.
Suppose that a countable group $G$ can be partitioned into countably many
dense subsets. Is it true that in $ G$ there exists an expansive
sequence? Cf.~15.80.

\hfill \raisebox{-1ex}{\sl I.\,V.\,Protasov}
\emp

\bmp \textbf{21.108.}
For a finite group $G$ let $\mathrm{Cod}(G)$ denote the set of irreducible character codegrees of $G$ (see 20.78). Define $\sigma(G)=\max\{ |\pi(m)|: \, m\in \mathrm{Cod}(G)\}$, where $\pi(m)$ denotes the set of prime divisors of an integer $m$. It is proved that there exists a~constant $k$ such that $|\pi(G)| \leq k \cdot \sigma(G)$  for every  finite group $G$
(Y.\,Yang, G.\,Qian, {\it J. Algebra}, {\bf 478} (2017), 215--219), but the estimate provided for $k$ is very crude. Can the constant $k$ be taken as $4$?
  \hfill \raisebox{-1ex}{\sl G.\,Qian}
\emp

\bmp \textbf{21.109.}
{\it Conjecture}: The derived length of  a finite solvable group $G$  does not exceed $|\mathrm{Cod}(G)|-1$.

\makebox[15pt][r]{}The Fitting height of  $G$ is known to be at most
$\min\{|\mathrm{Cod}(G)|-\nobreak 1, \, |\mathrm{Cod}(G)|/2+1\}$ (G.\,Qian,  Y.\,Zeng, {\it J.~Group Theory}, to appear). \hfill \raisebox{-1ex}{\sl G.\,Qian}
\emp

\bmp \textbf{21.110.}
Let $S$ be a nonabelian finite simple group, and $x$ a nonidentity auto\-mor\-phism of~$S$. Let $\alpha(x)$ be the smallest number of conjugates of $x$ in $G=\langle x, \operatorname{Inn}S\rangle$ that generate $G$. The values of $\alpha(x)$ had been studied in (R.\,Guralnick, J.\,Saxl,  {\it J. Algebra}, \textbf{268}, no.\,2 (2003), 519--571).

\makebox[25pt][r]{(a)} (R.\,Guralnick, J.\,Saxl). {\it Conjecture}:  If $S$ is an exceptional group of Lie type, then $\alpha(x)\leqslant 5$ for every nonidentity automorphism $x$ of $S$.

\makebox[25pt][r]{(b)} For each exceptional group $S$ of Lie type, find the largest value of $\alpha(x)$.

\hf \raisebox{-1ex}{{\sl D.\,O.\,Revin}}
\emp

\bmp \textbf{21.111.}  Let $S$ be a finite simple nonabelian group that is not isomorphic to any group ${}^2B_2(q)$. A nonidentity automorphism $x$ of $S$ is called a \emph{$\tau$-automorphism} if every two conjugates of $x$ in $\langle x, \operatorname{Inn}(S)\rangle$  generate a subgroup of order not divisible by~$3$. If $S$ admits a $\tau$-automorphism, we call $S$ a \emph{$\tau$-group}.

\makebox[25pt][r]{(a)} List all $\tau$-groups up to isomorphism.

\makebox[25pt][r]{(b)} Do $\tau$-automorphisms of odd order exist?
\hf \raisebox{0ex}{{\sl D.\,O.\,Revin, N.\,Y.\,Yang}}
\emp

\bmp \textbf{21.112.} A nonempty class $\mathfrak{X}$ of finite groups is said to be \emph{complete} if $\mathfrak{X}$ is closed under taking subgroups, homomorphic images, and extensions. The \emph{symmetric boundary} of a complete class $\mathfrak{X}$ other than the class of all finite groups is defined as the largest integer $n$ such that $\mathbb{S}_n \in \mathfrak{X}$. Every positive integer $n \ne 3$ coincides with the symmetric boundary of some complete class. It is proved (mod CFSG, D.\,O.\,Revin, \emph{Algebra i Analiz}, \textbf{37}, no.\,1 (2025), 141--176 (Russian)) that, for every complete class $\mathfrak{X}$, there exists a nonnegative integer $m$ with the following property: for every finite group $G$ and each conjugacy class $D$ of $G$, if every $m$ elements of $D$ generate a subgroup belonging to $\mathfrak{X}$, then $\langle D\rangle \in \mathfrak{X}$. The smallest such $m$ is called the \emph{Baer--Suzuki width} of~$\mathfrak{X}$
denoted by $\operatorname{BS}(\mathfrak{X})$.
  It is also proved (mod CFSG, \emph{ibid.})  that, for a complete class $\mathfrak{X}$ of symmetric boundary $n$, the value of $\operatorname{BS}(\mathfrak{X})$ is at least $n$ and is bounded above in terms of~$n$. For every positive integer $n \ne 3$, let $f_+(n)$ and $f_-(n)$ be respectively the maximum and the minimum of $\operatorname{BS}(\mathfrak{X})$, where $\mathfrak{X}$ runs over all complete classes of symmetric boundary~$n$.

\makebox[25pt][r]{(a)} Find $f_+(n)$ for $n=4,5,6$. It is known that $f_+(1)=2$, $f_+(2)=3$, and $f_+(n)=2(n-1)$ for $n\geqslant 7$.

\makebox[25pt][r]{(b)}
 Is it true that $f_-(n)=n$ for all $n \ne 3$? This is known to be true for $n=1,2,4$.

\hf \raisebox{-1ex}{{\sl D.\,O.\,Revin, N.\,Y.\,Yang}}
\emp

\bmp \textbf{21.113.}
 Let $G$ be a finite group and $p$ be a prime. Let $\Psi_{p,G}$ be the class function of~$G$ which vanishes on all $p$-singular elements of $G$ and whose value at each $p$-regular element $x$ of $G$ is the number of $p$-elements of $C_{G}(x)$.

 \makebox[25pt][r]{(a)} Is it true that $\Psi_{p,G}$ is a character of $G$?

\makebox[25pt][r]{(b)} If yes, can $\Psi_{p,G}$ be afforded by a projective $RG$-module, where $R$ is a complete discrete valuation ring of characteristic zero such that the field of fractions of $R$ is a splitting field for $G$ and its subgroups, and the residue field $R/J(R)$ is a splitting field of characteristic $p$ for $G$ and its subgroups?

\makebox[15pt][r]{}It is known that $\Psi_{p,G}$ is a character when $G \cong S_{n}$ for any positive integer $n$ and any prime $p$ (T.\,Scharf, {\it J.~Algebra}, {\bf 139}, no.\,2 (1991), 446--457).
\hfill \raisebox{-1ex}{\sl G.\,Robinson}
\emp

\bmp \textbf{21.114.}
A finite group $G$ is called {\itshape weakly ab-maximal} if $|H:[H,H]| \leq |G:[G,G]|$ for all $H \leqslant G$. Do weakly ab-maximal groups have bounded derived length?

 \makebox[15pt][r]{}It is known that weakly ab-maximal groups are direct products of weakly ab-maximal $p$-groups (F.\,Lisi,  L.\,Sabatini, {\itshape J.~Group Theory}, \textbf{27} (2024), 1203--1217).

\hfill \raisebox{-1ex}{\sl L.\,Sabatini}
\emp

\bmp \textbf{\zv 21.115.}
Let $C_1,\ldots,C_n$ be (left or right) cosets of a finite group $G$ such that $U:=C_1\cup\cdots\cup C_n$ is not $G$. Is it always true that $|G\setminus U|\ge |G|/2^n$?

\makebox[15pt][r]{}Affirmative answers are known  in some special cases (B.\,Sambale,  M.\,T\u{a}rn\u{a}uceanu, {\it J.~Algebraic Combin.}, {\bf 55} (2022), 979--987).
\hfill \raisebox{-1ex}{\sl B.\,Sambale}

\ul

\otv Yes, it is true (Anonymous AI-proof confirmed by B.\,Sambale: {\it Preprint}, August 2026, \url{https://www.unsolvedmath.com/problems/KOU-21.115}; \ B.\,Sambale, {\it Preprint}, September 2026, \url{https://arxiv.org/abs/2609.09052}).
\emp

\bmp \textbf{21.116.} A group is \textit{boundedly acyclic} if its bounded cohomology with trivial real coefficients vanishes in all positive degrees.
Is every branch group boundedly acyclic?

\hfill \raisebox{-1ex}{\sl E.\,Schesler}
\emp

 \bmp \textbf{\zv 21.117.}
(a)  Does there exist a finitely generated simple group that is of exponential growth but not of uniformly exponential growth?

\makebox[15pt][r]{}(b) Does there exist a finitely generated hereditarily just-infinite group that is of exponential growth but not of uniformly exponential growth?
\hfill \raisebox{-1ex}{\sl E.\,Schesler}

\ul

\otv Yes, Thompson's group $V$ provides an affirmative answer to both questions (R.\,Sauer, E.\,Schesler, {\it Preprint}, 2026, \url{https://arxiv.org/abs/2605.30163}).
\emp

\bmp \textbf{21.118.} (A.\,Thom).  Is there any group which is not isomorphic to the quotient of a residually finite group by an amenable normal subgroup?
\hfill \raisebox{-1ex}{\sl E.\,Schesler}
\emp

\bmp \textbf{21.119.} Does there exist a group $G$ that contains a family $(G_n)_{n \in \mathbb{N}}$ of finite-index subgroups  such that for every $n$ there is a homomorphism $f_n : G_n \rightarrow \mathbb{Z}$ whose kernel is of type $F_n$, but not of type $F_{n+1}$?
\hfill \raisebox{-1ex}{\sl E.\,Schesler}
\emp

\bmp \textbf{21.120.}  A pro-$p$ group is (\emph{relatively}) \emph{strictly finitely presented} if it is the pro-$p$ completion of a group that is finitely presented (respectively, finitely presented in some finitely-based variety of groups).
A pro-$p$ group is \emph{finitely axiomatizable} if it is determined up to isomorphism by a single sentence in the first-order
language of group theory.

 \makebox[25pt][r]{(a)} Does there exist a  (\emph{relatively}) \emph{strictly finitely presented} pro-$p$ group that is not finitely axiomatizable in the class of all pro-$p$ groups?

 \makebox[25pt][r]{(b)} In particular, is every finitely generated free pro-$p$ group finitely axiomatizable?

\makebox[15pt][r]{}See (A.\,Nies, K.\,Tent, D.\,Segal, \emph{Proc. London Math. Soc. (3)},
\textbf{123} (2021), 597--635; D.\,Segal, \emph{Preprint}, 2025, \url{https://arxiv.org/abs/2505.04816}).
 \hf\raisebox{-1ex}{\sl D.\,Segal}
\emp

\bmp \textbf{21.121.}
 Let $p$ be a prime number. A group $\Gamma$ is called \emph{$p$-Jordan}
if there exist constants $J$ and $e$ such that any finite subgroup
$G\subset\Gamma$ contains a normal abelian subgroup of order coprime to $p$
and of index at most~\mbox{$J\cdot |G_{(p)}|^e$}.
(For example by the results of Brauer--Feit and Larsen--Pink, for any field $\mathbb{K}$ of characteristic $p$
the group $\mathrm{GL}_n(\mathbb{K})$ is $p$-Jordan with $e=3$.)
Let the \emph{$p$-Jordan exponent}~$e(\Gamma)$ of the group $\Gamma$
be the infimum of all constants $e$ for which the above bound holds
for some constant $J=J(e)$.

\makebox[25pt][r]{\zva  a)}
Is it true that this infimum is always attained?

\makebox[25pt][r]{b)} Is it true that $e(\Gamma)\leqslant 3$
for any $p$-Jordan group~$\Gamma$?
 \hfill \raisebox{0ex}{\sl C.\,Shramov}

\ul

\otv No, not always (AI solution confirmed by C.\,Shramov:  \ Achyuth Jayadevan,  {\it Preprint of 14 September 2026},  \url{https://zenodo.org/records/22751568}; \ M.\,Kinyon, J.\,Urban, {\it Preprint of 17 September 2026},  \url{https://kourovkanotebookorg.wordpress.com/wp-content/uploads/2026/09/kourovka-mathematics.pdf}).
  \emp  

\bmp \textbf{21.122.}
Let $w$ be a group word, and $G$ a profinite group. Is it true that the cardinality of the set of $w$-values in $G$ is either finite or at least continuum?

\makebox[15pt][r]{}An affirmative answer is known for several important words; see (E.\,Detomi, B.\,Klopsch, P.\,Shumyatsky, {\it  J. London Math. Soc. (2)}, {\bf 102} (2020),  977--993).

\hfill \raisebox{-1ex}{\sl P.\,Shumyatsky}
\emp

\bmp \textbf{21.123.}
Is it true that the extension of the  A.\,Agrachev--R.\,Gamkrelidze construction of groups from pre-Lie rings in ({\it J.~Soviet Math.}, {\bf 17} (1981), 1650--1675) suggested in Definition~66 of (A.\,Smoktunowicz, {\it  J.~Pure Appl. Algebra}, {\bf  229}, no.\,12 (2025), 108128)  produces groups from pre-Lie rings?

\makebox[15pt][r]{}If this extended construction does give a group, then it also gives a brace, and so an  affirmative answer to this question would have consequences for the theory of set-theoretic solutions of the Yang--Baxter equation and for the theory of braces.

\hfill \raisebox{-1ex}{\sl A.\,Smoktunowicz}
\emp

   \bmp \textbf{21.124.}
A group $G$ is said to be \emph{virtually special} if $G$ has a finite-index
subgroup isomorphic to the fundamental group of a special complex
(in the sense of F.\,Haglund, D.\,T.\,Wise, \emph{Geom.\,Funct.\,Anal.}, \textbf{17}, no.\,5 (2008), 1551--1620; cf 20.60.)
A group $G$ is called a \emph{CAT(0) group} if it acts properly discontinuously and cocompactly by isometries on a CAT(0) metric space.

\makebox[25pt][r]{a)} Is every CAT(0) free-by-cyclic group virtually special?

\makebox[25pt][r]{b)} A weaker question: does every CAT(0) free-by-cyclic group virtually embed into a right-angled Artin group?
\hfill \raisebox{-1ex}{\sl I.\,Soroko}
\emp

   \bmp \textbf{21.125.} (M.\,Bridson). Let $F_m$ be a free group of rank $m$ and let $\phi\in\operatorname{Aut}(F_m)$ be a polynomially growing automorphism of maximal degree $m-1$, which means that for some (equivalently, any) free basis $\{x_1,\dots,x_m\}$ of $F_m$, the sequence $\max_i|\phi^n(x_i)|$ grows at the rate of $n^{m-1}$, where $|g|$ denotes the minimal length of $g$ in the $x_i$ and their inverses.

\makebox[25pt][r]{a)} Is the free-by-cyclic group $F_m\rtimes_\phi{\mathbb{Z}}$ virtually special?

\makebox[25pt][r]{b)}~In particular, are the Hydra groups
$G_m=F_m\rtimes\mathbb{Z}=\langle a_1,\dots,a_m,t\mid t^{-1}a_1t=\nobreak
a_1,\; t^{-1}a_it=a_ia_{i-1} \text{ for all }i>1\rangle$ virtually special?
\hfill \raisebox{-1ex}{\sl I.\,Soroko}
\emp

   \bmp \textbf{21.126.} (N.\,Brady). Do there exist finitely presented subgroups of right-angled Artin groups whose Dehn functions are super-exponential, or sub-exponential but not polynomial?

   \makebox[15pt][r]{}Such subgroups are known to exist in general CAT(0) groups, whereas the only Dehn functions currently realized for subgroups of right-angled Artin groups are exponential and polynomial of arbitrary degree.
\hfill \raisebox{-1ex}{\sl I.\,Soroko}
\emp

   \bmp \textbf{21.127.}
 Let $G$ be a right-angled Artin group. Is the stable commutator length $\operatorname{scl}(g)$ a rational number for every $g\in[G,G]$? For free groups this is true by Calegari's Rationality Theorem. (See {18.40} for the definition of $\operatorname{scl}(g)$.)
\hfill \raisebox{-1ex}{\sl I.\,Soroko}
\emp

   \bmp \textbf{21.128.}
Two groups $G_1$ and $G_2$ are said to be \emph{commensurable} if there exist finite-index subgroups $H_1\leqslant G_1$ and $H_2\leqslant G_2$ (not necessarily of the same index) such that $H_1\cong H_2$. Let $A[F_4]$ and $A[H_4]$ denote the Artin groups of spherical types $F_4$ and $H_4$, respectively. Are these two groups commensurable? This is the most difficult case in the classification of Artin groups of spherical type up to commensurability.
\hfill \raisebox{-1ex}{\sl I.\,Soroko}
\emp

   \bmp \textbf{21.129.}
 If two Artin groups of spherical type are quasi-isometric, must they be commensurable? (This is not true for right-angled Artin groups.)
 \hfill \raisebox{-1ex}{\sl I.\,Soroko}
\emp

\bmp \textbf{21.130.}
{\it Conjecture}:  Let $G$ be a finite additive abelian group with $|G|$ odd. Then any subset $A$ of $G$ with $|A|=n>2$ can be written as $\{a_1,\ldots,a_n\}$ in such a way that all the sums
$a_1+a_2,\ a_2+a_3,\ \ldots,\ a_{n-1}+a_n,\ a_n+a_1$
are distinct.
 \hf\raisebox{-1ex}{\sl Z.\,W.\,Sun}
\emp

\bmp \textbf{21.131.}
Construct a homomorphism of a subgroup of a Golod group (see 9.76) onto an infinite {\it $AT$-group} as defined in  (A.\,V.\,Rozhkov, {\it Math. Notes}, \textbf{40}, no.\,5 (1986), 827--836). Cf. 13.55
\hfill \raisebox{-1ex}{\sl A.\,V.\,Timofeenko}
\emp

\bmp \textbf{21.132.}
Based on the development of E.\,S.\,Golod's construction (see, for example, {\it Discrete Math. Appl.}, {\bf 23}, no.\,5--6 (2013), 491--501), for each prime number $p$, construct a finitely generated residually finite $p$-group with a non-trivial finite centre. Such groups with infinite and trivial centres are known (see 9.76 and Archive 11.101).

\hfill \raisebox{-1ex}{\sl A.\,V.\,Timofeenko}
\emp

\bmp \textbf{21.133.}
Does a group need to have a subnormal abelian series if every countable subgroup of it has such a series?
 \hf\raisebox{-1ex}{\sl M.\,Trombetti}
 \emp

  \bmp \textbf{\zv 21.134.}
 For a finite group $G$, let the \emph{type} of $G$ be the function on positive integers whose value at $n$ is the number of solutions of the equation $x^n=1$ in $G$.

\makebox[25pt][r]{a)} Is it true that a group having the same type as a group with trivial solvable radical must also have trivial solvable radical? Note that there are solvable and nonsolvable groups with the same type (see 12.37).

\makebox[25pt][r]{b)} Is it true that a group having the same type as an almost simple group must be isomorphic to it?
This is true for a group having the same type as a simple group,  as follows from the affirmative answer to~12.39 (in Archive).
 \hfill \raisebox{-1ex}{\sl A.\,V.\,Vasil'ev}

 \ul

 \otv It turned out that a negative answer to both questions had already been given by J.\,G.\,Thompson; see, for example, (Y.\,Li, W.\,Shi, {\it Ric. Mat.}, {\bf 74}, no.~1 (2025), 559--563)  ({\it Letter of A.\,V.\,Vasil'ev of 18 May 2026}).
  \emp

  \bmp \textbf{21.135.} For a finite group $G$,
let $\chi _1(G)$ denote the totality of the degrees of all
irre\-du\-cible complex characters of $G$ with allowance for their
multiplicities.  Suppose that $H$ is a finite group with $\chi_1(H)=\chi_1(G)$.
If
$G$ has trivial solvable radical, must $H$ also have trivial solvable radical?
 \hfill \raisebox{-1ex}{\sl A.\,V.\,Vasil'ev}
  \emp

\bmp \textbf{21.136.}
Let $G$ be a profinite group with fewer than $2^{\aleph_0}$ conjugacy classes of elements of infinite order.  Must $G$ be a torsion group?

\makebox[15pt][r]{}This holds in the case when $G$ is finitely generated (J.\,S.\,Wilson, {\em Arch. Math.}, {\bf 120} (2023), 557--563).
\hfill \raisebox{-1ex}{\sl John S.\,Wilson}
\emp

\bmp \textbf{\zv 21.137.}
 If the $p$-th powers in a finite $p$-group form a subgroup, must that subgroup be powerful? That is, for $p\ne 2$, if the $p$-th powers in a $p$-group of exponent $p^2$ form a subgroup, must that subgroup be abelian? For a $2$-group of exponent $8$, if the squares form a subgroup, must that subgroup be abelian?
 \hfill \raisebox{-1ex}{\sl L.\,Wilson}

 \ul

 \otv No, it need not be. For $p=3$ a counterexample is constructed in (K.\,Muliarchyk,  {\it Preprint}, 2026,
 \url{https://kourovkanotebookorg.wordpress.com/wp-content/uploads/2026/08/kourovka_21_137_counterexample.pdf}). For $p=\nobreak 2$, the Sylow 2-subgroup of $S_8$ has exponent 8, its squares form a subgroup, and that subgroup $C_2\times D_8$ is nonabelian and not powerful  ({\it Letter of A.\,Chang of 1 August 2026}). Further conterexamples for all odd primes $p$ are in (K.\,Muliarchyk,  {\it Preprint}, 2026,
 \url{https://kourovkanotebookorg.wordpress.com/wp-content/uploads/2026/09/kourovka_21.137_solution.pdf}).
\emp

\bmp \textbf{21.138.}
Let $G$ be an infinite finitely presented group such that every subgroup of infinite index is free. Must $G$ be isomorphic to either a free group or a surface group?

\hfill \raisebox{-1ex}{\sl H.\,Wilton}
\emp

\bmp \textbf{21.139.}
Let $G$ be a hyperbolic group which is virtually compact special in the sense of Haglund--Wise. Suppose that the set of second Betti numbers of the finite-index subgroups of $G$ is bounded. Must $G$ be virtually either a free group or a surface group?
\hfill \raisebox{-1ex}{\sl H.\,Wilton}
\emp

\bmp \textbf{21.140.}
Let $G$ be a torsion-free group of type $F_\infty$ of infinite cohomological dimension. Must $G$ contain a copy of Thompson's group $F$?
\hf\raisebox{-1ex}{\sl S.\,Witzel, M.\,C.\,B.\,Zaremsky}
\emp

\bmp \textbf{21.141.}
Let $G=G_1\amalg_H G_2$ be a free pro-$p$ product of coherent pro-$p$
groups with polycyclic amalgamation. Is $G$ coherent?

\makebox[15pt][r]{}For abstract groups this is known to be true.  A group is said to be \textit{coherent} if each of its finitely generated subgroups is finitely presented, and in the question the coherency is used in the pro-$p$ sense.
\hfill \raisebox{-1ex}{\sl P.\,A.\,Zalesskii}
\emp

\bmp \textbf{\zv 21.142.}
A group $G$ is said to be \emph{invariably generated} by $a$ and $b$ if $G$ is generated by the conjugates $a^g$, $b^h$ for every $g,h$.
Let $p\neq q$ be fixed primes.  Does every finite group embed into a
finite  group invariably generated by an element of order $p$ and an
element of order~$q$?
\hfill \raisebox{-1ex}{\sl P.\,A.\,Zalesskii}

\ul

\otv No, for any two primes $p\neq q$ there is a finite group that does not embed into a
finite  group invariably generated by an element of order $p$ and an
element of order~$q$ (T.\,Gong, M.\,R.\,Zeng, Y.\,Yang, {\it Preprint of  1 August 2026}, \url{https://arxiv.org/abs/2608.00703}. Examples for $p=2$, $q=3$ are also constructed in  (I.\,Capdeboscq, C.\,Parker, {\it Preprint of  4 August 2026}, \url{https://arxiv.org/abs/2608.03935}).
\emp

\bmp \textbf{21.143.} (Well-known problem).
 Is Thompson's group $F$ automatic?

\hfill \raisebox{-1ex}{\sl M.\,C.\,B.\,Zaremsky}

\emp
\bmp \textbf{21.144.}  (M.\,Brin,  M.\,Sapir). \emph{Conjecture}:
Every subgroup of Thompson's group~$F$ is either elementary amenable or else contains a subgroup isomorphic to $F$.

\hfill \raisebox{-1ex}{\sl M.\,C.\,B.\,Zaremsky}

\emp
\bmp \textbf{21.145.} (M.\,Bridson).  Is Thompson's group $F$ quasi-isometric

\makebox[25pt][r]{(a)} to $F \times \mathbb{Z}$?

\makebox[25pt][r]{(b)} to $F \times F$?
\hfill \raisebox{0ex}{\sl M.\,C.\,B.\,Zaremsky}

\emp
\bmp \textbf{21.146.} (Well-known problem).
A classifying space for a group $G$ is a connected $CW$-complex with fundamental group $G$ and all higher homotopy groups trivial. A group is of type $\mathrm{F}_n$ if it has a classifying space with finite $n$-skeleton. For example, type $\mathrm{F}_1$ is equivalent to finite generation, and type $\mathrm{F}_2$ is equivalent to finite presentability. Type $\mathrm{F}_\infty$ means type $\mathrm{F}_n$ for all $n$.

\makebox[15pt][r]{}For $n\ge 3$, does every group of type $\mathrm{F}_{n-1}$ embed as a subgroup of a group of type~$\mathrm{F}_n$? Or even in a group of type $\mathrm{F}_\infty$?
\hfill \raisebox{-1ex}{\sl M.\,C.\,B.\,Zaremsky}
\emp

\bmp \textbf{\zv 21.147.}
(V.\,M.\,Kopytov, N.\,Ya.\,Medvedev).
A subgroup $H$ of a right-orderable group~$G$ is said to be {\it right-relatively convex} if it is convex under some right ordering on $G$.
Is the lattice of right-relatively convex subgroups of a right-orderable group
always a sublattice of the lattice of its subgroups?
\hfill \raisebox{-1ex}{\sl A.\,V.\,Zenkov}

\ul

\otv No, not always (AI solution confirmed by A.\,V.\,Zenkov: \ P.\,Monticone, {\it Preprint}, 2026,
\url{https://kourovkanotebookorg.wordpress.com/wp-content/uploads/2026/06/21_147.pdf});

see also (W.\,van\,Doorn, E.\,Judin, P.\,Monticone, D.\,Morrison,  {\it Preprint}, 2026, \url{https://arxiv.org/abs/2607.17477}).

\emp

\bmp \textbf{\zv 21.148.}
(V.\,M.\,Kopytov, N.\,Ya.\,Medvedev).
Is it true that the lattice of right-relatively convex subgroups of a right-orderable group is distributive if and only if it is a chain?
\hfill \raisebox{-1ex}{\sl A.\,V.\,Zenkov}

\ul

\otv No, it is not  (K.\,Muliarchyk,  {\it Preprint}, 2026, \url{https://kourovkanotebookorg.wordpress.com/wp-content/uploads/2026/09/kourovka_21.148_solution.pdf}).
\emp

\bmp \textbf{21.149.}
(V.\,M.\,Kopytov, N.\,Ya.\,Medvedev).
Are there order automorphisms of Dlab groups that are not induced by conjugation by elements of a (possibly bigger) Dlab group?
\hfill \raisebox{-1ex}{\sl A.\,V.\,Zenkov}
\emp

\bmp \textbf{\zv 21.150.}
Let $G$ be an extension of a normal elementary abelian subgroup $A$ by an elementary abelian group $B\cong G/A$ such that $A$ contains an element $a$ with $C_B(a)=\nobreak 1$. Is it true that the rank of the subgroup
$Z(\langle a, B\rangle)\cap(\langle a, B\rangle) ' $ is at most the rank of $B$?
\hfill \raisebox{-1ex}{\sl V.\,I.\,Zenkov}

\ul

\otv No, not always. Let $A=\F_3[x, y]/I^3$ be the additive group of the quotient of the polynomial algebra $\F_3[x, y]$ by the ideal $I^3$, where $I$ is the ideal generated by $x,y$. Let $B$ be the group of automorphisms of $A$ generated by  multiplication by $1+x$ and $1+y$. Let $G=A\rtimes B$. For $a=1\in A$ we have $C_B(a)=1_B$. For $H=\langle a, B\rangle$ we have $Z(H)\cap H'={}_+\langle x^2, xy,y^2\rangle$, which has rank 3, while the rank of $B$ is 2. (P.\,Monticone, {\it Letter of 21 March 2026}); see also (P.\,Monticone,  {\it Preprint}, 2026, \url{https://kourovkanotebookorg.wordpress.com/wp-content/uploads/2026/03/solution_21_150-3.pdf})  and (W.\,van\,Doorn, E.\,Judin, P.\,Monticone, D.\,Morrison,  {\it Preprint}, 2026, \url{https://arxiv.org/abs/2607.17477}).
\emp

\newpage

\pagestyle{myheadings}
\markboth{\protect\vphantom{(y}{Archive of solved
problems (1st ed., 1965)}}{\protect\vphantom{(y}{Archive of solved
problems (1st ed., 1965)}} \thispagestyle{headings}

\parindent=0cm

~
\vspace{2ex}

\centerline {\Large {\bf Archive of solved problems}}
\phantomsection\label{archiv}
\vspace{4ex}

\mb{}This section contains all the solved problems that had already
been commented on, with a reference
to a detailed publication containing a complete answer, by the day of the  first appearance of the previous edition in 2022. The solutions that appeared in one of the updates after that day remain in the main part of the ``Kourovka
Notebook'', among the unsolved problems of the corresponding
section. \vskip4ex

\bmp \textbf{1.1.} Do there exist non-trivial finitely-generated
divisible groups? Equivalently, do there exist non-trivial
finitely-generated divisible simple groups? \hf {{\sl
Yu.\,A.\,Bogan}}\vs

 Yes, such groups do exist (V.\,S.\,Guba, {\it Math.
USSR--Izv.}, {\bf 29} (1987), 233--277). \emp

\bmp \textbf{1.2.} Let $G$ be a group, $F$ a free group with free
generators $x_1,\ldots , x_n$, and $R$ the free product of $G$ and
$F$. An {\it equation $($in the unknowns $x_1,\ldots ,x_n)$
over\/} $G$ is an expression of the form $v(x_1,\ldots \,x_n)=1$,
where on the left is an element of $R$ not conjugate in $R$ to any
element of $G$. We call $G$ {\it algebraically closed\/} if every
equation over $G$ has a solution in $G$. Do there exist
algebraically closed groups? \hfill {\sl L.\,A.\,Bokut'}\vs

Yes, such groups do exist (S.\,D.\,Brodski\u{\i}, Dep. no.\,2214-80,
VINITI, Moscow, 1980 (Russian)). \emp

\bmp \textbf{1.4.}
 (A.\,I.\,Mal'cev). Does there exist a ring without zero divisors
which is not em\-bed\-dable in a skew field, while the
multiplicative semigroup of its non-zero elements is embeddable in
a group? \hf {{\sl L.\,A.\,Bokut'}} \vs

Yes, such groups do exist (L.\,A.\,Bokut', {\it Soviet Math. Dokl.}, {\bf 8}
(1967), 901--904; \ A.\,Bowtell, {\it J.~Algebra}, {\bf 7} (1967),
126--139; \ A.\,Klein, {\it J.~Algebra}, {\bf 7} (1967),
100--125). \vskip3ex

\emp \bmp \textbf{1.9.}
 Can the factor-group of a locally normal group by the second
term of its upper central series be embedded (isomorphically) in
a direct product of finite groups?

\hf {\sl Yu.\,M.\,Gorchakov}\vs

Yes, it can (Yu.\,M.\,Gorchakov, {\it Algebra and Logic}, {\bf 15}
(1976), 386--390).

\emp \bmp \textbf{1.10.}
 An automorphism $\f$ of a group $G$ is called {\it splitting\/}
if $gg^{\f}\cdots g^{\f^{n-1}}=1$ for any element $g\in G$,
where $n$ is the order of $\f$. Is a soluble group admitting a
regular splitting automorphism of prime order necessarily
nilpotent?\hf {\sl Yu.\,M.\,Gorchakov}\vs

Yes, it is (E.\,I.\,Khukhro, {\it Algebra and Logic}, {\bf 17}
(1978), 402--406) and even without the regularity condition on the
automorphism (E.\,I.\,Khukhro, {\it Algebra and Logic},
 {\bf 19} (1980), 77--84).

\emp \bmp \textbf{1.11.}
 (E.\,Artin). The conjugacy problem for the braid groups ${\cal
B}_n$, $n>4$.

\hf {\sl M.\,D.\,Greendlinger}\vs

The solution is affirmative (G.\,S.\,Makanin, {\it Soviet Math.
Dokl.}, {\bf 9} (1968), 1156--1157; \ F.\,A.\,Garside, {\it Quart.
J. Math.}, {\bf 20} (1969), 235--254).

\emp

\bmp \textbf{1.13.}
(J.\,Stallings). If a finitely presented group is trivial, is it
always possible to replace one of the defining words by a
primitive element without changing the group?

\hfill {\sl M.\,D.\,Greendlinger}

\vs

No, not always
(S.\,V.\,Ivanov, {\it Invent. Math.}, \textbf{165}, no.\,3 (2006),
525--549).

\emp

\bmp \textbf{1.14.}
 (B.\,H.\,Neumann). Does there exist an infinite simple finitely
presented group?

\hf {\sl M.\,D.\,Greendlinger}\vs

Yes, such groups do exist (G.\,Higman, {\it Finitely presented infinite
simple groups $($Notes on Pure Mathematics, no.\,$8)$} Dept. of
Math., Inst. Adv. Studies Austral. National Univ., Canberra,
1974), where there is also a reference to R.\,Thompson.

\emp \bmp \textbf{1.17.}
 Write an explicit set of generators and defining relations for
a universal finitely presented group.\hf {\sl
M.\,D.\,Greendlinger}\vs

This was done (M.\,K.\,Valiev, {\it Soviet Math. Dokl.}, {\bf 14}
(1973), 987--991).

\emp \bmp \textbf{1.18.} (A.\,Tarski). a) Does there exist an
algorithm for determining
the solubility of equations in a free group?

\mb{b)} Describe the
structure of all solutions of an equation when it has at least
one solution.\hf {\sl Yu.\,L.\,Ershov}\vs

{a)} Yes, there does (G.\,S.\,Makanin, {\it Math. USSR--Izv.}, {\bf
21} (1982), 546--582; {\bf 25} (1985), 75--88).

{b)} This was described (A.\,A.\,Razborov, {\it Math.
USSR--Izv.}, {\bf 25} (1984), 115--162).

\emp

\bmp \textbf{1.19.}
(A.\,I.\,Mal'cev). Which subgroups (subsets) are
first order definable in a free group? Which subgroups are
relatively elementarily definable in a free group? In particular,
is the derived subgroup first order definable (relatively
elementarily definable) in a free group? \hfill
{\sl Yu.\,L.\,Ershov}
\vs

All definable (and relatively elementarily definable) sets are described; in particular, a subgroup is relatively elementarily definable if and only if it is cyclic (O.\,Kharlampovich, A.\,Myasnikov, {\it Int. J. Algebra Comput.}, {\bf 23}, no.\,1 (2013), 91--110).

\emp

\bmp \textbf{1.21.}
 Are there only finitely many finite simple groups of a given
exponent $n$?

\hf {\sl M.\,I.\,Kargapolov}\vs

Yes, modulo CFSG, see, for example, (G.\,A.\,Jones, {\it J. Austral.
Math. Soc.}, {\bf 17} (1974), 162--173).

\emp \bmp \textbf{1.23.}
 Does there exist an infinite simple locally finite group of
finite rank?

\hf {\sl M.\,I.\,Kargapolov}\vsh

No (V.\,P.\,Shunkov, {\it Algebra and Logic}, {\bf 10} (1971),
127--142).

\emp \bmp \textbf{1.24.}
 Does every infinite group possess an infinite abelian
subgroup?

\hf {\sl M.\,I.\,Kargapolov}\vs

No (P.\,S.\,Novikov, S.\,I.\,Adian, {\it Math. USSR--Izv.}, {\bf
2} (1968), 1131--1144).

\emp \bmp \textbf{1.25.}
 a) Is the universal theory of the class of finite groups
decidable?

\mb{b)} Is the universal theory of the class of finite nilpotent
groups decidable?

\hf {\sl M.\,I.\,Kargapolov}\vsh

a) No (A.\,M.\,Slobodskoi, {\it Algebra and Logic}, {\bf 20} (1981),
139--156). \

b) No (O.\,G.\,Kharlampovich, {\it Math. Notes}, {\bf 33} (1983),
254--263).

\emp \bmp \textbf{1.26.}
 Does elementary equivalence of two finitely generated
nilpotent groups imply that they are isomorphic?\hf {\sl
M.\,I.\,Kargapolov}\vs

No (B.\,I.\,Zil'ber, {\it
Algebra and Logic}, {\bf 10} (1971), 192--197).

\emp

\bmp \textbf{1.29.} (A.\,Tarski). Is
 the elementary theory of a free group decidable?

\hfill {\sl M.\,I.\,Kargapolov}

\vs

Yes, it is
(O.\,Kharlampovich, A.\,Myasnikov, {\it J.~Algebra}, \textbf{302},
no.\,2 (2006), 451--552).

\emp

\bmp \textbf{1.30.}
 Is the universal theory of the class of soluble groups
decidable?

\hf {\sl M.\,I.\,Kargapolov}\vsh

No (O.\,G.\,Kharlampovich, {\it Math. USSR--Izv.}, {\bf 19}
(1982), 151--169).

\emp \bmp \textbf{1.32.}
 Is the Frattini subgroup of a finitely generated matrix group
over a field nilpotent?\hfill {\sl M.\,I.\,Kargapolov}\vs

Yes, it is (V.\,P.\,Platonov, {\it Soviet Math. Dokl.}, {\bf 7}
(1966), 1557--1560).

\emp \bmp \textbf{1.34.}
 Has every orderable polycyclic group a faithful
representation by matrices over the integers?\hf {\sl
M.\,I.\,Kargapolov}\vs

Yes, it has (L.\,Auslander, {\it Ann. Math. (2)}, {\bf 86} (1967),
112--117; \ R.\,G.\,Swan, {\it Proc. Amer. Math. Soc.}, {\bf 18}
(1967), 573--574).

\emp \bmp \textbf{1.35.}
 A group is called {\it pro-orderable\/} if every partial
ordering of the group extends to a linear ordering.

\mb{a)} Is the wreath product of two arbitrary pro-orderable
groups again pro-orderable?

\mb{b)} (A.\,I.\,Mal'cev). Is every subgroup of a pro-orderable group
again pro-orderable?

\hf {\sl M.\,I.\,Kargapolov}\vs

Not always, in both cases (V.\,M.\,Kopytov, {\it Algebra i
Logika}, {\bf 5}, no.\,6 (1966), 27--31 (Russian)).

\emp \bmp \textbf{1.36.}
 If a group $G$ is factorizable by $p\hs$-sub\-groups, that
is, $ G = AB$, where $A$ and $B$ are $p\hs$-sub\-groups, does it
follows that $G$ is itself a $p\hs$-group?\hf {\sl
Sh.\,S.\,Kemkhadze}\vs

Not always (Ya.\,P.\,Sysak, {\it Products of infinite groups},
 Math. Inst. Akad. Nauk Ukrain. SSR, Kiev, 1982 (Russian)).\emp

 \bmp \textbf{1.37.} Is it true that every subgroup of a locally nilpotent group is quasi-invariant?

\hf {\sl Sh.\,S.\,Kemkhadze}\vs

No, not always; for example, in the inductive limit of Sylow $p$-subgroups of symmetric groups of degrees $p^i$, $i=1,2,\dots$. (V.\,I.\,Sushchanski\u{\i}, \emph{Abstracts of 15th All-USSR Algebraic Conf.}, Krasnoyarsk, 1979, 154 (Russian)).

\emp \bmp \textbf{1.38.}
An {\it $N^0$-group\/} is a group in which every cyclic
subgroup is a term of some normal system of the group. Is every
$N^0$-group an $\widetilde N$-group?\hf {\sl
Sh.\,S.\,Kemkhadze}\vs

No. Every group $G$ with a central system $\left\{
H_i\right\}$ is an $N^0$-group, since for any $g\in G$ the system
$\left\{ \left< g,H_i\right> \right\}$ refined by the trivial
subgroup and the intersections of
all the subsystems is a normal system of $G$ containing $\left<
g\right>$. Free groups have central systems but
are not $\widetilde N$-groups. (Yu.\,I.\,Merzlyakov, 1973.)

\emp

\bmp \textbf{1.39.} Is a group binary
nilpotent if it is the
product of two
normal binary nilpotent subgroups? \hf
{\sl Sh.\,S.\,Kemkhadze}\vs

No, not always (A.\,I.\,Sozutov, {\it Algebra and Logic}, {\bf
30}, no.\,1 (1991), 70--72).

\emp

\bmp \textbf{1.41.}
 A subgroup of a linearly orderable group is called {\it
relatively convex\/} if it is convex with respect to some linear
ordering of the group. Under what conditions is a subgroup of an
orderable group relatively convex?\hf {\sl
A.\,I.\,Kokorin}\vs

Necessary and sufficient conditions are given in (A.I. Kokorin,
V.\,M.\,Kopytov, {\it Fully ordered groups}, John Wiley, New
York, 1974).

\emp \bmp \textbf{1.42.}
 Is the centre of a relatively convex subgroup relatively
convex?\hf {\sl A.\,I.\,Kokorin}\vs

Not always (A.\,I.\,Kokorin, V.\,M.\,Kopytov, {\it Fully
ordered groups}, John Wiley, New York, 1974).

\emp \bmp \textbf{1.43.}
 Is the centralizer of a relatively convex subgroup relatively
convex?

\hf {\sl A.\,I.\,Kokorin}\vs

Not always (A.\,I.\,Kokorin, V.\,M.\,Kopytov, {\it Siberian
Math.~J.}, {\bf 9} (1968), 622--626).

\emp \bmp \textbf{1.44.}
 Is a maximal abelian normal subgroup relatively convex?\hf {\sl
A.\,I.\,Kokorin}\vs

Not always (D.\,M.\,Smirnov, {\it Algebra i Logika}, {\bf 5},
no.\,6 (1966), 41--60 (Russian)).

\emp \bmp \textbf{1.45.}
 Is the largest locally nilpotent normal subgroup relatively
convex?

\hf {\sl A.\,I.\,Kokorin}\vs

Not always (D.\,M.\,Smirnov, {\it Algebra i Logika}, {\bf 5},
no.\,6 (1966), 41--60 (Russian)).

\emp \bmp \textbf{1.47.}
 A subgroup $H$ of a group $G$ is said to be {\it strictly
isolated\/} if, whenever $xg_1^{-1}xg_1\cdots g_n^{-1}xg_n$
belongs to $H$, so do $x$ and each $g_i^{-1}xg_i$. A group in
which the identity subgroup is strictly isolated is called an
{\it $S$-group}. Do there exist $S$-groups that are not
orderable groups?\hf {\sl A.\,I.\,Kokorin}\vs

Yes, such groups do exist (V.\,V.\,Bludov, {\it Algebra and Logic}, {\bf 13}
(1974), 343--360).

\emp \bmp \textbf{1.48.}
 Is a free product of two orderable groups with an
amalgamated subgroup that is relatively convex in each of the
factors again an orderable group?\hf {\sl
A.\,I.\,Kokorin}\vs

Not always (M.\,I.\,Kargapolov, A.\,I.\,Kokorin, V.\,M.\,Kopytov,
{\it Algebra i Logika}, {\bf 4}, no.\,6 (1965), 21--27 (Russian)).

\emp \bmp \textbf{1.49.}
 An automorphism $\f$ of a linearly ordered group is called {\it
order preserving\/} if $ x < y$ implies $x^{\f} <y^{\f}$. Is it
possible to order an abelian strictly isolated normal subgroup of
an $S$-group (see Archive 1.47) in such a way that the order is
preserved under the action of the inner automorphisms of the whole
group?\hf {\sl A.\,I.\,Kokorin}\vs

Not always (V.\,V.\,Bludov, {\it
Algebra and Logic}, {\bf 11} (1972), 341--349).

\emp \bmp \textbf{1.50.}
 Do the order-preserving automorphisms of a linearly ordered
group form an orderable group?\hf {\sl
A.\,I.\,Kokorin}\vs

Not always (D.\,M.\,Smirnov, {\it Algebra i Logika}, {\bf 5},
no.\,6 (1966), 41--60 (Russian)).

\emp \bmp \textbf{1.52.}
 Describe the groups that can be ordered linearly in a unique
way (reversed orderings are not regarded as different).\hf {\sl
A.\,I.\,Kokorin}\vs

This was done (V.\,V.\,Bludov, {\it Algebra and Logic}, {\bf 13}
(1974), 343--360).

\emp \bmp \textbf{1.53.}
 Describe all possible linear orderings of a free nilpotent
group with finitely many generators.\hf {\sl
A.\,I.\,Kokorin}\vs

This was done (V.\,F.\,Kleimenov, in: {\it Algebra, Logic and
Applications}, Irkutsk, 1994, 22--27 (Russian)). \emp

\bmp \textbf{1.56.}
 Is a torsion-free group pro-orderable (see Archive 1.35) if the
factor-group by its centre is pro-orderable?\hf {\sl
A.\,I.\,Kokorin}\vs

Not always (M.\,I.\,Kargapolov, A.\,I.\,Kokorin, V.\,M.\,Kopytov,
{\it Algebra i Logika}, {\bf 4}, no.\,6 (1965), 21--27 (Russian)).

\emp \bmp \textbf{1.60.}
 Can an orderable metabelian group be embedded in a radicable orderable
 group?\hf {\sl A.\,I.\,Kokorin}\vs

Yes, it can (V.\,V.\,Bludov, N.\,Ya.\,Medvedev, {\it Algebra and
Logic}, {\bf 13} (1974), 207--209).

\emp

\bmp \textbf{1.61.} Can any orderable group be embedded into an
orderable group

\mb{a)} with a radicable maximal locally nilpotent
normal subgroup?

\mb{b)} with a radicable maximal abelian normal
subgroup? \hf{\sl A.\,I.\,Kokorin}\vs

Yes, it can, in both cases (S.\,A.\,Gurchenkov, {\it Math. Notes},
{\bf 51}, no.\,2 (1992), 129--132). \emp

\bmp \textbf{1.63.}
 A group $G$ is called {\it dense\/} if it has no proper
isolated subgroups other than its trivial subgroup.

\mb{a)} Do there exist dense torsion-free groups that are not
locally cyclic?

\mb{b)} Suppose that any two non-trivial elements $x$ and $y$ of a
torsion-free group $G$ satisfy the relation $ x^k = y^l$, where
$k$ and $l$ are non-zero integers depending on $x$ and~$y$. Does
it follow that $G$ is abelian?\hf {\sl
P.\,G.\,Kontorovich}\vs

a) Yes, such groups do exist; b) Not necessarily (S.\,I.\,Adian, {\it Math.
USSR--Izv.}, {\bf 5} (1971), 475--484).

\emp \bmp \textbf{1.64.}
 A torsion-free group is said to be {\it separable\/} if it can
be represented as the set-theoretic union of two of its proper
subsemigroups. Is every $R$-group separable?

\hf {\sl
P.\,G.\,Kontorovich}\vs

No (S.\,J.\,Pride, J.\,Wiegold, {\it Bull. London Math. Soc.}, {\bf 9}
(1977),
36--37).

\emp

\bmp \textbf{1.66.}
Suppose that $T$ is a periodic abelian group, and
${\goth m}$ an uncountable cardinal number. Does there always
exist an abelian torsion-free group $U(T,{\goth m})$ of
cardinality ${\goth m}$ with the following property: for any
abelian torsion-free group $A$ of cardinality $\leq {\goth m}$,
the equality ${\rm Ext}\, (A,T)=0$ holds if and only if $A$ is
embeddable in $U(T,{\goth m})$?

\hfill {\sl L.\,Ya.\,Kulikov}\vs

No, not always. There is a model of ZFC in which for a certain
class of cardinals ${\goth m}$ the answer is negative (S.\,Shelah,
L.\,Str\"ungmann, {\it J.~London Math. Soc.}, {\bf 67}, no.\,3
(2003), 626--642). On the other hand, under the assumption of
G\"odel's constructivity hypothesis ($V=L$) the answer is
affirmative for any cardinal if $T$ has only finitely many
non-trivial bounded $p$-components (L.\,Str\"ungmann, {\it Ill.
J.~Math.}, {\bf 46}, no.\,2 (2002), 477--490). \emp

\bmp \textbf{1.68.}
 (A.\,Tarski). Let ${\frak K}$ be a class of groups and $Q{\frak
K}$ the class of all homomorphic images of groups from ${\frak
K}$. If ${\frak K}$ is axiomatizable, does it follow that $Q{\frak
K}$ is?

\hf {\sl Yu.\,I.\,Merzlyakov}\vsh

Not always (F.\,Clare, {\it Algebra Universalis}, {\bf 5}
(1975), 120--124).

\emp \bmp \textbf{1.69.}
 (B.\,I.\,Plotkin). Do there exist locally nilpotent torsion-free
groups without the property $RN^*$?\hf {\sl
Yu.\,I.\,Merzlyakov}\vs

Yes, such groups do exist (E.\,M.\,Levich, A.\,I.\,Tokarenko, {\it Siberian
Math.~J.}, {\bf 11} (1970), 1033--1034; \ A.\,I.\,Tokarenko, {\it
Trudy Rizhsk. Algebr. Sem.}, Riga Univ., Riga, {\bf 1969},
280--281 (Russian)).

\emp \bmp \textbf{1.70.}
 Let $p$ be a prime number and let $G$ be the group of all matrices of the
form {\small $\left( \begin{array}{cc} 1+p\alpha & p\beta\\
p\gamma & 1+p\delta \end{array}\right)$}, where $\alpha$, $\beta$,
$\gamma$, $\delta$ are rational numbers with denominators coprime
to $p$. Does $G$ have the property $\,\overline{\! RN}$?\hf {\sl
Yu.\,I.\,Merzlyakov}\vs

Yes, it does (G.\,A.\,Noskov, {\it Siberian Math.~J.}, {\bf 14}
(1973), 475--477).

\emp \bmp \textbf{1.71.}
 Let $G$ be a connected algebraic group over an algebraically
closed field. Is the number of conjugacy classes of maximal
soluble subgroups of $G$ finite?

\hf {\sl V.\,P.\,Platonov}\vsh

Yes, it is (V.\,P.\,Platonov, {\it Siberian Math.~J.}, {\bf 10}
(1969), 800--804).

\emp \bmp \textbf{1.72.}
 D.\,Hertzig has shown that a connected algebraic group over an
algebraically closed field is soluble if it has a rational regular
automorphism. Is this result true for an arbitrary field?\hf {\sl
V.\,P.\,Platonov}\vs

No (V.\,P.\,Platonov, {\it Soviet Math. Dokl.}, {\bf 7} (1966),
825--829).

\emp \bmp \textbf{1.73.}
 Are there only finitely many conjugacy classes of maximal
periodic subgroups in a finitely generated linear group over the
integers?\hf {\sl V.\,P.\,Platonov}\vs

Not always. The extension $H$ of the free group on free generators
$a,b$ by the automorphism $ \f :a\rightarrow
a^{-1},\;\;b\rightarrow b$ is a linear group over the integers.
For every $n\in \Z$ the element $c_n=\f b^{-n}ab^n$ has order 2
and $C_H(c_n)=\left< c_n\right>$. Let the dash denote an
isomorphism of $H$ onto its copy $ H'$. As a counterexample one
can take the subgroup $G\leq H\times H'$ generated by the elements
$ a,\,\f ,\,a'\f ',\,bb'$. Indeed, $G$ contains the subgroups
$T_n=\left< c_0,c_n'\right>$, \ $n\in \Z$, which are maximal
periodic (even in $H\times H'$). If $T_n$ and $T_m$ are conjugate
by an element $xy'\in G$ (where $x\in H$ and $y'\in H'$), then
$c_0^x=c_0$ and $c_n^y=c_m$, whence $ x\in \left< c_0\right>$, $
y\in b^{m-n}\left< c_m\right>$. Since $xy'\in G$, the sums of the
exponents at the occurrences of $b$ in $x$ and $y$ (in any
expression) must coincide; hence $m=n$. (Yu.\,I.\,Merzlyakov,
1973.)

\emp \bmp \textbf{1.75.} Classify the infinite simple periodic
linear groups over a
field of characteristic $p > 0$.\hfill {\sl
V.\,P.\,Platonov}\vs

They are classified for every $p$ mod CFSG (V.\,V.\,Belyaev, in: {\it
 Investigations in Group Theory}, Sverdlovsk, UNC AN SSSR, 1984,
39--50 (Russian); \ A.\,V.\,Borovik, {\it Siberian Math.~J.}, {\bf
24}, no.\,6 (1983), 843--851; \
 B.\,Hartley, G.\,Shute, {\it Quart. J.~Math. Oxford (2)}, {\bf 35}
(1984), 49--71; \ S.\,Thomas, {\it Arch. Math.}, {\bf 41} (1983),
103--116). For $p>2$ there is a proof without using the CFSG
(A.\,V.\,Borovik, {\it Siberian Math.~J.}, {\bf 25}, no.\,2
(1984), 221--235). For all $p$, without using CFSG, the classification also follows from the paper (M.\,J.\,Larsen,
R.\,Pink, {\it J.~Amer. Math. Soc.}, {\bf 24} (2011), 1105--1158) also known as a preprint of
1998.

\emp

\bmp \textbf{1.76.}
 Does there exist a simple, locally nilpotent, locally compact,
topological group?

\hf {\sl V.\,P.\,Platonov}\vsh

No (I.\,V.\,Protasov, {\it Soviet Math. Dokl.}, {\bf 19} (1978),
487--489).

\emp \bmp \textbf{1.78.} Let a group $G$ be the product of two
divisible abelian $p\hs$-groups of finite rank. Is then $G$ itself a
divisible abelian $p\hs$-group of finite rank?\hf {\sl
N.\,F.\,Sesekin}\vs

Yet, it is (N.\,F.\,Sesekin, {\it Siberian Math.~J.}, {\bf 9}
(1968), 1070--1072).

\emp \bmp \textbf{1.80.}
 Does there exist a finite simple group whose Sylow 2-subgroup
 is a direct product of quaternion groups?\hf {\sl
A.\,I.\,Starostin}\vs

No (G.\,Glauberman, {\it J.~Algebra}, {\bf 4} (1966), 403--420).

\emp \bmp \textbf{1.81.}
 The {\it width\/} of a group $G$ is, by definition, the
smallest cardinal $ m = m(G)$ with the property that any
subgroup of $G$ generated by a finite set $S\subseteq G$ is
generated by a subset of $S$ of cardinality at most $m$.

\mb{a)} Does a group of finite width satisfy the minimum condition
for subgroups?

\mb{b)} Does a group with the minimum condition for subgroups have
finite width?

\mb{c)} The same questions under the additional condition of local
finiteness. In particular, is a locally finite group of finite
width a Chernikov group?\hfill {\sl L.\,N.\,Shevrin}\vs

a) Not always (S.\,V.\,Ivanov, {\it Geometric methods in the study
of groups with given subgroup properties}, Cand. Diss.,
Moscow Univ., Moscow, 1988 (Russian)).

b) Not always
(G.\,S.\,Deryabina, {\it Math. USSR--Sb.}, {\bf 52} (1985),
481--490).

c) Yes, it is (V.\,P.\,Shunkov, {\it Algebra and Logic}, {\bf 9}
(1970), 137--151; {\bf 10} (1971), 127--142).

\emp \bmp \textbf{1.82.}
 Two sets of identities are said to be {\it equivalent\/} if they
determine the same variety of groups. Construct an infinite set
of identities which is not equivalent to any finite one.\hf {\sl
A.\,L.\,Shmel'kin}\vs

This was done (S.\,I.\,Adian, {\it Math. USSR--Izv.}, {\bf 4}
(1970), 721--739).

\emp \bmp \textbf{1.83.}
 Does there exist a simple group, the orders of whose elements
are unbounded, in which a non-trivial identity relation
holds?\hf {\sl A.\,L.\,Shmel'kin}\vs

Yes, such groups do exist (V.\,S.\,Atabekyan, Dep. no.\,5381-V86, VINITI,
Moscow, 1986 (Russian)).

\emp \bmp \textbf{1.84.} Is it true that
 a polycyclic group $G$ is
residually a finite $p\hs$-group if and only if $G$ has a
nilpotent normal torsion-free subgroup of $p\hs$-power index?\hf
{\sl A.\,L.\,Shmel'kin}\vs

No, it is not (K.\,Seksenbaev, {\it Algebra i Logika}, {\bf 4},
no.\,3 (1965), 79--83 (Russian)).

\emp \bmp \textbf{1.85.}
 Is it true that the identity relations of a metabelian group
have a finite basis?

\hf {\sl A.\,L.\,Shmel'kin}\vsh

Yes, it is (D.\,E.\,Cohen, {\it J.~Algebra}, {\bf 5} (1967),
267--273).

\emp \bmp \textbf{1.88.}
 Is it true that if a matrix group over a field of
characteristic 0 does not satisfy any non-trivial identity
relation, then it contains a non-abelian free subgroup?

\hf {\sl
A.\,L.\,Shmel'kin}\vsh

Yes, it is (J.\,Tits, {\it J.~Algebra}, {\bf 20} (1972),
250--270).

\emp \bmp \textbf{1.89.}
 Is the following assertion true? Let $G$ be a free soluble
group, and $a$ and $b$ elements of $G$ whose normal closures
coincide. Then there is an element $ x\in G$ such that $ b^{\pm
1} = x^{-1}ax$.\hf {\sl A.\,L.\,Shmel'kin}\vs

No (A.\,L.\,Shmel'kin, {\it Algebra i Logika}, {\bf 6}, no.\,2
(1967), 95--109 (Russian)).

\emp

\bmp \textbf{1.90.}
 A subgroup $H$ of a group $G$ is called {\it $2$-infi\-ni\-te\-ly
isolated\/} in $G$ if, whenever the centralizer $C_G(h)$ in $G$
of some element $h\ne 1$ of $H$ contains at least one involution
 and intersects $H$ in an infinite subgroup, it follows
that $C_G(h)\leq H$. Let $G$ be an infinite simple locally
finite group whose Sylow 2-subgroups are Chernikov groups, and
suppose $G$ has a proper 2-infinitely isolated subgroup $H$
containing some Sylow 2-subgroup of $G$. Does it follow that $G$
is isomorphic to a group of the type $PSL_2(k)$, where $k$ is a
field of odd characteristic?\hf {\sl V.\,P.\,Shunkov}\vs

Yes, it does (V.\,P.\,Shunkov, {\it Algebra and Logic}, {\bf 11} (1972),
260--272).

\emp

\bmp \textbf{2.1.}
 Classify the finite groups having a Sylow $p\hs$-sub\-group as a
maximal subgroup.

\hf {\sl V.\,A.\,Belonogov,
A.\,I.\,Starostin}\vs

A description can be extracted from (B.\,Baumann, {\it J.
Algebra}, {\bf 38} (1976), 119--135) and (L.\,A.\,Shemetkov, {\it
Math. USSR--Izv.}, {\bf 2} (1968), 487--513).

\emp

\markboth{\protect\vphantom{(y}{Archive of solved
problems (2nd ed., 1966)}}{\protect\vphantom{(y}{Archive of solved
problems (2nd ed., 1966)}}

\bmp \textbf{2.2.}
 A {\it quasigroup\/} is a groupoid $Q(\cdot )$ in which the
equations $ ax = b$ and $ ya = b$ have a unique solution for any
$ a, b\in Q$. Two quasigroups $Q(\cdot )$ and $Q(\circ )$ are
{\it isotopic\/} if there are bijections $\alpha$, $\beta$,
$\gamma$ of the set $Q$ onto itself such that $x \circ y=\gamma
(\alpha x\cdot \beta y)$ for all $ x, y \in Q$. It is well-known
that all quasigroups that are isotopic to groups form a variety
${\frak G}$. Let ${\frak V}$ be a variety of quasigroups.
 Characterize the class of groups isotopic to quasigroups in $
{\frak G}\cap {\frak V}$.
 For which identities characterizing ${\frak V}$ is every
group isotopic to a quasigroup in $ {\frak G}\cap {\frak V}$?
 Under what conditions on $ {\frak V}$ does any group
isotopic to a quasigroup in $ {\frak G}\cap {\frak V}$ consist of a
single element? \hf {\sl V.\,D.\,Belousov,~A.\,A.\,Gvaramiya}\vs

Every part is answered (A.\,A.\,Gvaramiya, Dep. no.\,6704-V84,
VINITI, Moscow, 1984 (Russian)).

\emp \bmp \textbf{2.3.}
 A finite group is called {\it quasi-nilpotent\/} ($\Gamma$-{\it
quasi-nil\-po\-tent\/}) if any two of its (maximal) subgroups $A$ and
$B$ satisfy one of the conditions 1) $ A \leq B$, \ \ 2) $ B \leq
A$, \ \ 3) $N_A(A\cap B) \ne A\cap B \ne N_B(A\cap B)$. Do the
classes of quasi-nilpotent and $\Gamma$-quasi-nil\-po\-tent groups
coincide? \hfill {\sl Ya.\,G.\,Berkovich,~M.\,I.\,Kravchuk}\vs

No. The group $G=\left< x,y,z,t\mid
x^4=y^4=z^2=t^3=1,\;\;[x,y]=z,\;\; [x,z]=[y,z]=1,\right.$ $\left.
x^t=y,\;\; y^t=x^{-1}y^{-1}\right>$ is $\Gamma$-quasi-nil\-po\-tent,
but not quasi-nilpotent. Since $ G/\Phi (G)\cong {\Bbb A}_4$,
 the intersection of any two maximal subgroups $A$ and $B$ of
$G$ equals $ \Phi (G)$, whence $N_A(A\cap B)\ne A\cap B\ne
N_B(A\cap B)$; thus $G$ is $\Gamma$-quasi-nil\-po\-tent. On the
other hand, if $ A_1=\left< zx^2,\,zy^2,\,t\right>$ and $B_1=\left<
z,t\right>$, then $ N_{A_1}(A_1\cap B_1)=$ $A_1\cap B_1=\left<
t\right>$; hence $G$ is not quasi-nilpotent.
 (V.\,D.\,Mazurov, 1973.)

\emp \bmp \textbf{2.4.}
 (S.\,Chase). Suppose that an abelian group $A$ can be written as
the union of pure subgroups $A_{\alpha}$, $\alpha\in\Omega$,
where $\Omega$ is the first non-denumerable ordinal,
$A_{\alpha}$ is a free abelian group of denumerable rank, and,
 if $ \beta < \alpha$, then $A_{\beta}$ is a direct summand of~$
A_{\alpha}$. Does it follow that $A$ is a free abelian
group?\hf {\sl Yu.\,A.\,Bogan}\vs

Not always (P.\,A.\,Griffith, {\it Pacific J.~Math.}, {\bf 29} (1969),
279--284).

\emp \bmp \textbf{2.7.}
 Find the cardinality of the set of all polyverbal operations
acting on the class of all groups.\hf {\sl O.\,N.\,Golovin}\vs

It has the cardinality of the continuum (A.\,Yu.\,Olshanskii,
{\it Math. USSR--Izv.}, {\bf 4} (1970), 381--389).

\emp

\bmp \textbf{2.8.} (A.\,I.\,Mal'cev). Do there exist regular
associative operations having the heredity property with respect
to transition from the factors to their subgroups? \hfill {\sl
O.\,N.\,Golovin}\vs

Yes, there do (S.\,V.\,Ivanov, {\it Trans. Moscow Math. Soc.}, {\bf
1993}, 217--249).
 \emp

\bmp \textbf{2.10.}
Prove an analogue of the Remak--Shmidt theorem for decompositions of a group into nilpotent products.
 \hfill {\sl O.\,N.\,Golovin}
 \vs

Such an analogue is proved (V.\,V.\,Limanski\u{\i}, \emph{Trudy Moskov. Mat. Obshch.}, {\bf 39} (1979), 135--155 (Russian)).
\emp

\bmp \textbf{2.13.}
 (Well-known problem). Let $G$ be a periodic group in which
every $\pi$-ele\-ment commutes with every $\pi '$-ele\-ment. Does
$G$ decompose into the direct product of a maximal
$\pi$-sub\-group and a maximal $\pi '$-sub\-group?\hf {\sl
S.\,G.\,Ivanov}\vs

Not always. S.\,I.\,Adian ({\it Math.
USSR--Izv.}, {\bf 5}
(1971), 475--484) has
 constructed a group $A=A(m,n)$
which is torsion-free and has a central element $d$ such that
$A(m,n)/\left< d\right>\cong B(m,n)$, the free $m$-gene\-ra\-tor
Burnside group of odd exponent $ n\geq 4381$. Given a
 prime $p$ coprime to $n$, a
 counterexample with $ \pi =\left\{ p\right\}$ can be found in the form
$G=A/\langle d^{p^k}\rangle$ for some positive integer $k$. Indeed,
$\left< d\right> /\langle d^{p^k}\rangle$ is a maximal
$p\hs$-sub\-group of $G$. Suppose that $A/\langle d^{p^k}\rangle =
\left< d\right> /\langle d^{p^k}\rangle\times H_k/\langle
d^{p^k}\rangle$ for every $k$. Then $ \left< d\right> \cap
H_k=\langle d^{p^k}\rangle$ for all $k$ and therefore $ \left<
d\right>\cap H=1$, where $H=\bigcap_kH_k$. Since $H$ is
torsion-free and is isomorphic to a subgroup of $ A/\left<
d\right>\cong B(m,n)$, we obtain $H=1$. This implies
that $A$ is isomorphic to a subgroup of the Cartesian product of
the abelian groups $A/H_k$, a contradiction. (Yu.\,I.\,Merzlyakov,
1973.)

\emp \bmp \textbf{2.15.}
 Does there exist a torsion-free group such that the factor
group by some term of its upper central series is nontrivial
periodic, with a bound on the orders of the elements?\hf {\sl
G.\,A.\,Karas\"ev}\vs

Yes, such groups do exist (S.\,I.\,Adian, {\it Proc. Steklov Inst. Math.}, {\bf 112}
(1971),
61--69).

\emp \bmp \textbf{2.16.}
 A group $G$ is called {\it conjugacy separable\/} if any two of
its elements are conjugate in $G$ if
and only if their images are conjugate in every finite
homomorphic image of~$G$. Is $G$ conjugacy separable in the
following
cases:

\mb{a)} $G$ is a polycyclic group,

\mb{b)} $G$ is a free soluble group,

\mb{c)} $G$ is a group of (all) integral matrices,

\mb{d)} $G$ is a finitely generated group of matrices,

\mb{e)} $G$ is a finitely generated metabelian group?\hf {\sl
M.\,I.\,Kargapolov}\vs

a) Yes (V.\,N.\,Remeslennikov, {\it Algebra and Logic}, {\bf 8}
(1969), 404--411); \ E.\,Formanek, {\it J.~Algebra}, {\bf 42}
(1976), 1--10).

b) Yes (V.\,N.\,Remeslennikov, V.\,G.\,Sokolov, {\it Algebra and
Logic}, {\bf 9} (1970), 342--349).

c), d) Not always (V.\,P.\,Platonov, G.\,V.\,Matveev, {\it Dokl.
Akad. Nauk BSSR}, {\bf 14} (1970), 777--779 (Russian); \
 V.\,N.\,Remeslennikov, V.\,G.\,Sokolov, {\it Algebra and Logic},
{\bf 9} (1970), 342--349; \ V.\,N.\,Remeslennikov, {\it Siberian
Math.~J.}, {\bf 12} (1971), 783--792).

e) Not always. Let $p$ be a prime and let $A_1,A_2$ be two
copies of the additive group $\{ m/p^k\mid m,k\in \Z \}$. Let
$b_1,b_2$ be the automorphisms of the direct sum $A=A_1\oplus A_2$
defined by $a^{b_2}=pa$ for any $a\in A$ and $a_2^{b_1}=a_1+a_2$
and $a_1^{b_1}=a_1$ for some fixed elements $a_1\in A_1$, $a_2\in
A_2$. Let $G$ be the semidirect product of $A$ and the direct
product $\langle b_1\rangle \times \langle b_2\rangle$ (of two
infinite cyclics). It can be shown that the elements $a_2$ and
$a_2+a_1/p$ are not
conjugate in $G$, but their images are conjugate in any finite
quotient of $G$. (M.\,I.\,Kargapolov, E.\,I.\,Timoshenko, {\it
Abstracts of the 4th All-Union Sympos. Group Theory,
Akademgorodok, 1973}, Novosibirsk, 1973,
86--88 (Russian).)

\emp \bmp \textbf{2.17.}
 Is it true that the wreath product $A\wr B$ of two groups
that are conjugacy separable is itself conjugacy separable if and
only if either $A$ is abelian or $B$ is finite?

\hf {\sl
M.\,I.\,Kargapolov}\vsh

No (V.\,N.\,Remeslennikov, {\it Siberian Math.~J.}, {\bf 12}
(1971), 783--792).

\emp

\bmp \textbf{2.18.}
 Compute the ranks of the factors of the lower central series of
a free soluble group.\hf {\sl M.\,I.\,Kargapolov}\vs

They were computed (V.\,G.\,Sokolov, {\it Algebra and Logic},
{\bf 8} (1969),
212--215; \ Yu.\,M.\,Gor\-chakov, G.\,P.\,Egorychev, {\it Soviet
Math. Dokl.}, {\bf 13} (1972), 565--568).

\emp \bmp \textbf{2.19.}
 Are finitely generated subgroups of a free soluble group
finitely separable?

\hf {\sl M.\,I.\,Kargapolov}\vs

Not always (S.\,A.\,Agalakov, {\it Algebra and Logic}, {\bf 22} (1983),
261--268).

\emp \bmp \textbf{2.20.}
 Is it true that the wreath products $A\wr B$ and $A_1\wr B_1$
are elementarily equivalent if and only if $A$,~$B$ are
elementarily equivalent to $A_1$,~$B_1$, respectively?

\hf {\sl
M.\,I.\,Kargapolov}\vs

No, but if the word ``elementarily" is replaced by the word
``universally," then it is true (E.\,I.\,Timoshenko, {\it Algebra
and Logic}, {\bf 7}, no.\,4 (1968), 273--276).

\emp \bmp \textbf{2.21.}
 Do the classes of Baer and Fitting groups coincide?\hf {\sl
Sh.\,S.\,Kemkhadze}\vs

No (R.\,S.\,Dark, {\it Math. Z.}, {\bf 105} (1968), 294--298).

\emp \bmp \textbf{2.22.}
 b) An abstract group-theoretical property $\Sigma$ is called
{\it radical $($in our sense\/$)$} if, in any group $G$, the join
$\Sigma(G)$ of all normal $\Sigma$-sub\-groups is a
$\Sigma$-sub\-group. Is the property $N^0$ (see Archive,~1.38)
radical?\hf {\sl Sh.\,S.\,Kemkhadze}\vs

No. The group $SL_n(\Z )$ for sufficiently large $n\geq 3$ contains a
non-abelian finite simple group and therefore is not an $N^0$-group.
On the other hand, as shown in (M.\,I.\,Kargapolov,
Yu.\,I.\,Merzlyakov, {\it in: Itogi Nauki. Algebra. Topologiya.
Geometriya, 1966}, VINITI, Moscow, 1968, 57--90 (Russian)), it is the
product of its congruence-subgroups mod\,2 and mod\,3, which have
central systems and hence are $N^0$-groups (see Archive, 1.38).
(Yu.\,I.\,Merzlyakov, 1973.)

\emp \bmp \textbf{2.23.}
 a) A subgroup $H$ of a group $G$ is called {\it
quasisubinvariant\/} if there is a normal system of $G$ passing
through $H$. Let ${\frak K}$ be a class of groups that is closed
with respect to taking homomorphic images. A group $G$ is called
an {\it $R^0({\frak K})$-group\/} if each of its non-trivial
homomorphic images has a non-trivial quasisubinvariant ${\frak
K}$-sub\-group. Do $R^0({\frak K})$ and $\,\overline{\! RN}$
coincide when ${\frak K}$ is the class of all abelian groups?

\hfill
{\sl Sh.\,S.\,Kemkhadze}\vs

No, they do not  (J.\,S.\,Wilson, {\it Arch. Math.}, {\bf 25} (1974), 574--577).

\emp

\bmp \textbf{2.25.}
 b) Do there exist groups that are linearly orderable in
countably many ways?

\hf {\sl A.\,I.\,Kokorin}\vs

Yes, there do (R.\,N.\,Buttsworth, {\it Bull. Austral. Math. Soc.},
{\bf 4}
(1971), 97--104).

\emp
\bmp
 {\bf 2.29.}
 Does the class of finite groups in which every proper abelian
subgroup is contained in a proper normal subgroup coincide with
the class of finite groups in which every proper abelian subgroup
is contained in a proper normal subgroup of prime index? \hf {\sl
P.\,G.\,Kontorovich,~V.\,T.\,Nagrebetski\u{\i}}\vs

No. Let $B$ be a finite group such that $B=[B,B]\ne 1$ and let $r$
be the rank of $B$. Put $ A=\bigoplus\limits_{p\mid |B|}(\Z /p\Z )
^{r+1}$ and $G=A\wr B$. We define a homomorphism $ \f
:G\rightarrow A$ by setting $(bf)^{\f}=\sum_{x\in B}f(x)$, where
$b\in B$ and $ f\in F={\rm Fun}(B,A)$. Suppose that $f^b=f$ for
$b\in B$ and $f\in F$. It is clear that $f^{\f}\in nA$, where
$n=|b|$. In particular, $ C_F(b)^{\f}\leq pA\leq O_{p'}(A)$ if
$p\mid n$. Every proper abelian subgroup $H$ of $G$ is contained
in a proper normal subgroup. Indeed, we may assume that $H\not\leq
F$. We fix an element $bf\in H$, where $ f\in F$, $b\in B$, $b\ne
1$. Let $p$ be a prime dividing $|b|$. For any $ h\in H\cap F$ we
have $bfh=hbf=bh^bf=bfh^b$, whence $h=h^b$. Hence $ (H\cap
F)^{\f}\leq C_F(b)^{\f}\leq O_{p'}(A)$. If $T$ is the full
preimage of $O_{p'}(A)$ in $G$, then $ G/T\cong (\Z
 /p\Z )^{r+1}$ and therefore $H/H\cap T$ is an elementary abelian
$p\hs$-group. The rank of it is $\leq r$, since
$H/H\cap T$ embeds into $B$ and $H\cap F\leq H\cap T$. Thus, $HT$ is a
proper normal subgroup containing $H$. On the other hand, $F$ is
a proper abelian subgroup that is not contained in any proper
normal subgroup of prime index, since $G/F\cong B=[B,B]$.
(G.\,M.\,Bergman, I.\,M.\,Isaacs, {\it Letter of June,
17, 1974.})

\emp

 \bmp \textbf{2.30.}
 Does there exist a finite group in which a Sylow $p\hs$-sub\-group
is covered by other Sylow $p\hs$-sub\-groups? \hf {\sl
P.\,G.\,Kontorovich,~A.\,L.\,Starostin}\vs

Yes, there does, for any prime $p$ \ (V.\,D.\,Mazurov, {\it Ural
Gos. Univ. Mat. Zap.}, {\bf 7}, no.\,3 (1969/70), 129--132
(Russian)).

\emp \bmp \textbf{2.31.}
 Can every group admitting an ordering with only finitely many
convex subgroups be represented by matrices over a field?\hfill
{\sl V.\,M.\,Kopytov}\vs

No, not every. For example, the group
 $G=\left< a_n,\,b_n\;(n\in \Z
),\;c,\,d\mid [a_n,b_n]=c ,\right.$ $\left. a_n^d=a_{n+1},\;\;
b_n^d=b_{n+1}, \;\; [a_n,a_m]=[b_n,b_m]=[a_n,c]=[b_n,c]=1\;
\;(n,m\in \Z ) \right>$ is soluble and orderable with finitely
many convex subgroups, but it is not residually finite and hence
is not linear over a field. (V.\,A.\,Churkin, 1969.)

\emp \bmp \textbf{2.33.}
 Is a direct summand of a direct sum of finitely generated
modules over a~Noetherian ring again a direct sum of finitely
generated modules?\hfill {\sl V.\,I.\,Kuz'minov}\vs

Not always (P.\,A.\,Linnell, {\it Bull. London Math. Soc.}, {\bf 14}
(1982),
124--126).
\emp

\bmp{\bf 2.35.} An inverse spectrum $\xi$ of abelian groups is
said to be {\it acyclic\/} if ${\lim\limits_{\longleftarrow}} ^{(p)}
\xi =0$ for $p > 0$. Here ${\lim\limits_{\longleftarrow}}^{(p)}$
denotes the right derived functor of the projective limit functor.
Let $\xi$ be an acyclic spectrum of finitely generated groups. Is
the spectrum $\bigoplus\xi _{\alpha}$ also acyclic, where each
spectrum $\xi _{\alpha}$ coincides with $\xi $?\hf{\sl
V.\,I.\,Kuz'minov}\vs

The answer depends on the axioms of Set Theory (A.\,A.\,Khusainov,
{\it Siberian Math.~J.}, {\bf 37}, no.\,2 (1996), 405--413).

\emp

\bmp \textbf{2.36.}
 (de Groot). Is the group of all continuous integer-valued functions
on a compact space free abelian?\hf {\sl
V.\,I.\,Kuz'minov}\vs

Yes, it is. By (G.\,N\"obeling, {\it Invent. Math.}, {\bf 6}
(1968) 41--55) the additive group of all bounded
integer-valued functions
on an arbitrary set is free. Hence the group of all
continuous integer-valued functions on the Chech
compactification of an arbitrary discrete space is also free.
For any compact space $X$ there exists a continuous mapping of
the Chech compactification $Y$ of a discrete space onto $X$.
This induces an
 embedding of the group of
continuous integer-valued functions on $X$ into the free group of
continuous integer-valued functions on $Y$.
(V.\,I.\,Kuz'minov, 1969.)

\emp \bmp \textbf{2.37.}
 Describe the finite simple groups whose Sylow $p\hs$-sub\-groups
are cyclic for all odd $p$.\hf {\sl V.\,D.\,Mazurov}\vs

This was done (M.\,Aschbacher, {\it J.~Algebra}, {\bf 54} (1978),
50--152).

\emp \bmp \textbf{2.38.}
 (Old problem). The class of rings embeddable in associative
division rings is universally axiomatizable. Is it finitely
axiomatizable?\hf {\sl A.\,I.\,Mal'cev}\vs

No (P.\,M.\,Cohn, {\it Bull. London Math. Soc.}, {\bf 6} (1974), 147--148).

\emp

\bmp \textbf{2.39.}
 Does there exist a non-finitely-axiomatizable variety of

 \mb{a)} (H.\,Neumann) groups?

\mb{b)} of associative rings (the
Specht problem)?

\mb{c)} Lie rings?\hf {\sl A.\,I.\,Mal'cev}\vs

a) Yes, there does (S.\,I.\,Adian, {\it Math. USSR--Izv.}, {\bf
4} (1970), 721--739; \ A.\,Yu.\,Ol'\-shan\-ski\u{\i}, {\it Math.
USSR--Izv.}, {\bf 4} (1970), 381--389).

b) Yes, there does
 (A.\,Ya.\,Belov,
 {\it Fundam. Prikl. Mat.},
{\bf 5},
 no.\,1 (1999), 47--66 (Russian), \ {\it Sb. Math.}, {\bf 191},
no.\,3--4 (2000), 329--340; \ A.\,V.\,Grishin, {\it Fundam. Prikl.
Mat.}, {\bf 5}, no.\,1 (1999), 101--118 (Russian); \ V.\,V.\,Shchigolev, {\it Fundam. Prikl.
Mat.}, {\bf 5}, no.\,1 (1999), 307--312 (Russian)).
 But every variety of associative algebras over a
 field of characteristic 0 has a finite basis of identities
 (A.\,R.\,Kemer, {\it Algebra and Logic}, {\bf 26}, no.\,5 (1987),
 362--397).

c) Yes, there does (M.\,R.\,Vaughan-Lee,
{\it
Quart. J.~Math.}, {\bf 21} (1970), 297--308).

\emp

\bmp \textbf{2.40.}
 The {\it $I$-theory\/} ({\it $Q$-theory\/}) of a class ${\frak
K}$ of universal algebras is the totality of all identities
(quasi-identities) that are true on all the algebras in ${\frak
K}$. Does there exist a finitely axiomatizable variety of

\mb{a)} groups,

\mb{b)} semigroups,

\mb{c)} (1) rings

\makebox[25pt][r]{(2)} of
 associative rings

\makebox[40pt][r]{(i)} whose $I $-theory is non-decidable?

\makebox[40pt][r]{(ii)} whose $Q\hskip0.1ex $-theory is non-decidable?

 \makebox[25pt][r]{(3)} of Lie rings

\makebox[40pt][r]{(ii)} whose $Q\hskip0.1ex $-theory is non-decidable?

whose $I$-theory ($Q\hs$-theory) is non-decidable?\hf {\sl
A.\,I.\,Mal'cev}\vs

a) Yes (Yu.\,G.\,Kleiman, {\it Trans. Moscow Math. Soc.}, {\bf
1983}, no.\,2, 63--110).

 b) Yes (V.\,L.\,Murski\u{\i}, {\it Math. Notes}, {\bf
3} (1968), 423--427).

c) (1) Yes (V.\,Yu.\,Popov,
{\it Math. Notes}, \textbf{67} 
(2000), 495--504).
(2i) No, it does not
(A.\,Ya.\,Belov, \textit{Izv. Math.}
\textbf{74}, no.\,1 (2010), 1--126).

(2ii) and (3ii): Yes, it
exists (A.\,I.\,Budkin, {\it Izv. Altai Univ.}, {\bf 65}, no.\,1 (2010), 15--17 (Russian)).
\emp

\bmp \textbf{2.41.} Is the variety generated
by \

\mb{a)}~a~finite associative ring;

\mb{b)}~a~finite Lie ring;

\mb{c)}~a~finite quasigroup

finitely axiomatizable?

\mb{d)} What is the cardinality $n$ of the smallest semigroup
generating a non-finitely-axiomatizable variety? \hfill {\sl
 A.\,I.\,Mal'cev}\vs

a) Yes (I.\,V.\,L'vov, {\it Algebra and Logic}, {\bf 12} (1973),
156--167; \ R.\,L.\,Kruse, {\it
 J.~Algebra}, {\bf 26} (1973), 298--318).

b) Yes (Yu.\,A.\,Bakhturin,
A.\,Yu.\,Olshanskii, {\it Math. USSR--Sb.}, {\bf 25} (1975), 507--523).

c) Not always (M.\,R.\,Vaughan-Lee, {\it Algebra
Universalis}, {\bf 9} (1979), 269--280).

d) It is proved that $n=6$ \ (P.\,Perkins, {\it J. Algebra}, {\bf 11} (1969),
298--314; \ A.\,N.\,Trakht\-man, {\it Semigroup Forum}, {\bf 27}
(1983), 387--389).

\emp

\bmp \textbf{2.43.} A group $G$ is called an {\it
$FN$-group\/} if the groups $\gamma_ iG/\gamma _{i+1}G$ are free
abelian and $\bigcap_{i=1}^{\infty} \gamma _iG = 1$, where $\gamma
_{i+1}G=[\gamma _{i}G, G]$. A variety of groups ${\frak M}$ is
called a {\it $\Sigma$-variety\/} (where $\Sigma$ is an abstract
property) if the ${\frak M}$-free groups have the property
$\Sigma$.

\mb{a)} Which properties $\Sigma$ are preserved
under multiplication and intersection of varieties? Is the
$FN$ property preserved under these operations?

\mb{b)} Are all varieties obtained by multiplication and
intersection from the nilpotent varieties
${\frak N}_1$, ${\frak N}_2, \ldots $ (where ${\frak N}_1$ is the
variety of abelian groups) $FN$-varie\-ties?

\hf {\sl
A.\,I.\,Mal'cev}\vs

a) The property $FN$ is preserved by multiplication of varieties
(A.\,L.\,Shmel'kin, {\it Trans. Moscow Math. Soc.}, {\bf 29}
(1973), 239--252). The property $FN$ is not preserved by the
intersection of varieties. The following example is due to
L.\,G.\,Kov\'acs. Let ${\goth U}$ and ${\goth V}$ be the varieties
of all nilpotent groups of class at most~4 satisfying the
identities $[x,y,y,x] \equiv 1$ and $[[x,y],[z,t]] \equiv 1$,
respectively. Then ${\goth U}$ and ${\goth V}$ are
$FN$-varie\-ties because 1)~both of them are nilpotent of class 4
and contain all nilpotent groups of class $\leq 3$, and
2)~relatively free groups in ${\goth U}$ and ${\goth V}$ are
torsion-free (well-known for ${\goth V}$ and follows for ${\goth
U}$ from (P.\,Fitzpatrick, L.\,G.\,Kov\'acs, {\it J.~Austral.
Math. Soc. Ser.~A}, {\bf 35}, no.\,1 (1983), 59--73)). On the
other hand, ${\goth U} \cap {\goth V}$ is not an $FN$-varie\-ty
because the relatively free group in ${\goth U} \cap {\goth V}$ of
rank~$3$ is not torsion-free. Indeed, every torsion-free group in
${\goth U} \cap {\goth V}$ is of class $\leq 3$ (follows from the
paper by Fitzpatrick and Kov\'acs cited above), and there is a
3-generated (exponent~2)-by-(exponent~2) group of class
precisely~4 in ${\goth U} \cap {\goth V}$. (A.\,N.\,Krasil'nikov,
{\it Letter of July, 17, 1998}.)

b) Yes (Yu.\,M.\,Gorchakov, {\it Algebra i Logika}, {\bf 6},
no.\,3 (1967), 25--30 (Russian)). \emp

\bmp \textbf{2.44.}
 Let ${\goth A}$ and ${\goth B}$ be subvarieties of a variety of
groups $
{\goth M}$; then $({\goth A}{\goth B})\cap {\goth M}$ is
called the {\it ${\goth M}$-pro\-duct of ${\goth A}$ by ${\goth
B}$\/}, where $ {\goth A}{\goth B}$ is the usual product. Does
there exist a non-abelian variety ${\goth M}$ with an infinite
lattice of subvarieties and commutative ${\goth
M}$-multi\-pli\-ca\-tion?\hfill {\sl A.\,I.\,Mal'cev}\vs

Yes, there does; for example, $ {\frak M}={\frak A}_p{\frak A}_p$,
where ${\frak A}_p$ is the variety of all abelian groups of prime
exponent $p$ (Yu.\,M.\,Gorchakov, {\it Talk of June 21, 1967,
Krasnoyarsk}).

\emp

\bmp \textbf{2.45.}
 (P.\,Hall). Prove or refute the following conjectures:

\mb{a)} If a word $v$ takes only finitely many values in a group $G$,
then the verbal subgroup $vG$ is finite.

\mb{c)} If $G$ satisfies the maximum condition and $vG$ is finite,
then $v^*G$ has finite index in $G$.\hf {\sl
Yu.\,I.\,Merzlyakov}\vs

These conjectures are refuted in:

a) (S.\,V.\,Ivanov, {\it
Soviet Math. (Izv. VUZ)}, {\bf 33}, no.\,6 (1989), 59--70;

c) I.\,S.\,Ashmanov, A.\,Yu.\,Olshanskii,
{\it Soviet Math. (Izv. VUZ)}, {\bf 29}, no.\,11 (1985), 65--82).

\emp

\bmp \textbf{2.46.}
 Find conditions under which a finitely-generated matrix group
is almost residually a finite $p\hs$-group for some prime $p$. \hfill
{\sl Yu.\,I.\,Merzlyakov}\vs

A finitely generated matrix group over a field of characteristic
zero has a subgroup of finite index which is residually a finite
$p\hs$-group for almost all primes $p$ (M.\,I.\,Kargapolov, {\it
Algebra i Logika}, {\bf 6}, no.\,5 (1967), 17--20 (Russian); \
Yu.\,I.\,Merzlyakov, {\it Soviet Math. Dokl.}, {\bf 8} (1967),
1538--1541; \ V.\,P.\,Platonov, {\it Dokl. Akad. Nauk Belarus.
SSR}, {\bf 12}, no.\,6 (1968), 492--494 (Russian)). \ A finitely
generated matrix group over a field of characteristic $ p>0$ has a
subgroup of finite index which is residually a finite $p\hs$-group
(V.\,P.\,Platonov, {\it Dokl. Akad. Nauk Belarus. SSR}, {\bf 12},
no.\,6 (1968), 492--494 (Russian); \ A.\,I.\,Tokarenko, {\it Proc.
Riga Algebr. Sem.}, Latv.
 Univ., Riga, 1969, 280--281 (Russian)). See also some
generalizations in (V.\,N.\,Remeslennikov, {\it Algebra i Logika},
{\bf 7}, no.\,4 (1968), 106--113 (Russian); \
B.\,A.\,F.\,Wehrfritz, {\it Proc. London Math. Soc.}, {\bf 20},
no.\,1 (1970), 101--122).

\emp \bmp \textbf{2.47.} In which abelian groups is the lattice of
all fully invariant subgroups a chain?

\hf {\sl
A.\,P.\,Mishina}\vs

Such groups were characterized (M.-P.\,Brameret, {\it S\'eminaire
P.\,Dubreil, M.-L.\,Dubreil-Jacotin, L.\,Lesieur, et C.\,Pisot},
{\bf 16} (1962/63), no.\,13, Paris, 1967).

\emp \bmp \textbf{2.49.}
 (A.\,Selberg). Let $G$ be a connected semisimple linear Lie
group whose corresponding symmetric space has rank greater than
1, and let $\Gamma$ be an irreducible discrete subgroup of $G$
such that $ G/\Gamma$ has finite volume. Does it follow that
$\Gamma$ is an arithmetic subgroup?\hf {\sl
V.\,P.\,Platonov}\vs

Yes, it does (G.\,A.\,Margulis, {\it Functional Anal. Appl.}, {\bf
8} (1974), 258--259).

\emp \bmp \textbf{2.50.}
 (A.\,Selberg). Let $\Gamma$ be an irreducible discrete subgroup
of a connected Lie group $G$ such that the factor space
$G/\Gamma$ is non-compact but has finite volume in the Haar
measure. Prove that $\Gamma$ contains a non-trivial unipotent
element.\hf {\sl V.\,P.\,Platonov}\vs

This is proved (D.\,A.\,Kazhdan, G.\,A.\,Margulis, {\it
Math. USSR--Sb.}, {\bf 4} (1968), 147--152).

\emp \bmp \textbf{2.51.}
 (A.\,Borel, R.\,Steinberg). Let $G$ be a semisimple algebraic
group and $R_G$ the set of classes of conjugate unipotent
elements of $ G$. Is $R_G$ finite?\hf {\sl
V.\,P.\,Platonov}\vs

Yes, it is (G.\,Lusztig, {\it Invent. Math.}, {\bf 34} (1976), 201--213).

\emp \bmp \textbf{2.52.}
 (F.\,Bruhat, N.\,Iwahori, M.\,Matsumoto). Let $G$ be a
semisimple algebraic group over a locally compact, totally
disconnected field. Do the maximal compact subgroups of $G$ fall
into finitely many conjugacy classes? If so, estimate this
number.

\hf {\sl V.\,P.\,Platonov}\vs

Yes, they do, and a formula for this number is found (F.\,Bruhat,
J.\,Tits, {\it Publ. Math. IHES}, {\bf 41} (1972), 5--251).

\emp

\bmp \textbf{2.54.} Can $SL(n,k)$ have maximal subgroups that are not closed in the Zariski topology?
\hfill {\sl V.\,P.\,Platonov}\vs

Yes, it can if $k$ is either $\Q$ or an algebraically closed field of characteristic zero (N.\,S.\,Ro\-ma\-novski\u{\i}, \emph{Algebra i Logika}, {\bf 6}, no.\,4 (1967), 75--82; {\bf 7}, no.\,3 (1968), 123 (Russian)).

\emp

\bmp \textbf{2.55.}
 b) Does $SL_n(\Z )$, $n\geq 2$, have maximal subgroups of infinite index?

\hf {\sl V.\,P.\,Platonov}\vs

Yes, it does (G.\,A.\,Margulis, G.\,A.\,Soifer, {\it Soviet Math.
Dokl.}, {\bf 18} (1977), 847--851).

\emp \bmp \textbf{2.58.}
 Let $G$ be a vector space and $\Gamma$ a group of automorphisms
of $G$. $\Gamma$ is called {\it locally finitely stable\/} if,
for any finitely generated subgroup $\Delta$ of $\Gamma$, \ $G$
has a finite series stable relative to $\Delta$. If the
characteristic of the field is zero and $\Gamma$ is locally
finitely stable, then $\Gamma$ is locally nilpotent and
torsion-free. Is it true that every locally nilpotent
torsion-free group can be realized in this way?\hf {\sl
B.\,I.\,Plotkin}\vs

No (L.\,A.\,Simonyan, {\it Siberian Math.~J.}, {\bf 12} (1971),
602--606).

\emp \bmp \textbf{2.59.}
 Let $\Gamma$ be any nilpotent group of class $ n-1$. Does $
\Gamma$ always admit a faithful representation as a group of
automorphisms of an abelian group with a series of length $n$
stable relative to $\Gamma$? \hf {\sl B.\,I.\,Plotkin}\vs

Not always (E.\,Rips, {\it Israel J. Math.}, {\bf 12} (1972),
342--346).

\emp \bmp \textbf{2.61.} Let $\G$ be a Noetherian group of
automorphisms of a vector space $G$ such that every element of
$\G$ is unipotent. Is $\G$ necessarily a stable group of
automorphisms?

\hf {\sl B.\,I.\,Plotkin}\vs

Not if the field has non-zero characteristic
(A.\,Yu.\,Olshanskii, {\it The geometry of defining relations
in groups}, Kluwer, Dordrecht, 1991). This is true for fields of
characteristic zero if the indices of unipotence are uniformly
bounded (B.\,I.\,Plotkin, S.\,M.\,Vovsi, {\it Varieties of group
representations}, Zinatne, Riga, 1983 (Russian)).

\emp \bmp \textbf{2.62.}
 If $G$ is a finite-dimensional vector space over a field and $
\G$ is a group of automorphisms of $G$ in which every element is
 stable, then the whole of $\G$ is stable (E.\,Kolchin). Is
Kolchin's theorem true for spaces over skew fields?\hf {\sl
B.\,I.\,Plotkin}\vs

Yes, it is, if the characteristic of the skew field
is zero or sufficiently large compared to the
dimension of the space (H.\,Y.\,Mochizuki, {\it Canad. Math. Bull.}, {\bf
21} (1978), 249--250).

\emp
\bmp \textbf{2.63.} Let $G$ be a group of automorphisms of a vector space over a field of characteristic zero, and suppose that all elements of $G$ are unipotent, with uniformly bounded unipotency indices. Must such a group be locally finitely stable?
\hf {\sl
B.\,I.\,Plotkin}\vs

Yes, it must, which follows from (H.\,Heineken, \emph{Arch. Math. (Basel)}, {\bf 13} (1962), 29--37). Moreover,
such a group is even finitely stable, as follows from (E.\,I.\,Zel'manov, \emph{Math. USSR-Sb.}, {\bf 66}, no.\,1 (1990), 159--168).
\emp

\bmp \textbf{2.64.}
 Does the set of nil-elements of a finite-dimensional linear
group coincide with its locally nilpotent radical?\hf {\sl
B.\,I.\,Plotkin}\vs

Not always. The non-nilpotent $3$-gene\-ra\-tor nilgroup of E.\,S.\,Golod
constructed by an algebra over $\Q$ is residually torsion-free
nilpotent; hence it is orderable and therefore it is embeddable into $
GL_n$ over some skew field by Mal'cev's theorem (V.\,A.\,Roman'kov, 1978).

\emp \bmp \textbf{2.65.}
 Does the adjoint group of a radical ring (in the sense of
Jacobson) have a central series?\hf {\sl
B.\,I.\,Plotkin}\vs

Not always (O.\,M.\,Neroslavski\u{\i}, {\it Vesci Akad. Navuk
BSSR,
Ser. Fiz.-Mat. Navuk}, {\bf 1973}, No. 2, 5--10 (Russian)).

\emp \bmp \textbf{2.66.}
 Is an $R$-group determined by its subgroup lattice?
 Is every lattice isomorphism of $R$-groups induced by a group
isomorphism?\hf {\sl L.\,E.\,Sadovski\u{\i}}\vs

Not always, in both cases (A.\,Yu.\,Olshanskii, {\it
C.~R.~Acad. Bulgare Sci.}, {\bf 32}, no.\,9 (1979), 1165--1166).

\emp \bmp \textbf{2.69.}
 Let a group $G$ be the product of two subgroups $A$ and $B$,
each of which is nilpotent and satisfies the minimum condition.
Prove or refute the following: \ a) $G$ is soluble; \ b) the
divisible parts of $A$ and $B$ commute elementwise.\hf {\sl
N.\,F.\,Sesekin}\vs

Both parts are proved (N.\,S.\,Chernikov, {\it Soviet Math.
Dokl.}, {\bf 21} (1980), 701--703).

\emp \bmp \textbf{2.70.}
 a) Let a group $G$ be the product of two subgroups $A$ and $B$,
each of which is locally cyclic and torsion-free. Prove that
either $A$ or $B$ has a non-trivial subgroup that is normal in
$G$.

\mb{b)} Characterize the groups that can be factorized in this way.
\hf {\sl N.\,F.\,Sesekin}\vs

a) This was proved (D.\,I.\,Zaitsev, {\it Algebra and Logic},
{\bf 19} (1980), 94--106).

b) They were characterized (Ya.\,P.\,Sysak, {\it
Algebra and Logic}, {\bf 25} (1986), 425--433).

\emp

\bmp \textbf{2.71.} Does there exist a finitely generated
right-orderable group which coincides with its derived subgroup
and, therefore, does not have the property $RN$? \hfill {\sl
D.\,M.\,Smirnov}\vs

Yes, there does (G.\,M.\,Bergman, {\it Pacific J.~Math.}, {\bf
147}, no.\,2 (1991), 243--248).

\emp

\bmp \textbf{2.72.} (G.\,Baumslag).
Suppose that $F$ is a finitely generated free group, $N$~its normal
subgroup and $V$ a fully invariant subgroup of $N$. Is $F/V$
necessarily Hopfian if $F/N$ is Hopfian? \hfill {\sl
D.\,M.\,Smirnov}

\vs

No, not necessarily
(S.\,V.\,Ivanov, A.\,M.\,Storozhev, {\it Geom. Dedicata}, \textbf{
114} (2005), 209--228). \emp

\bmp \textbf{2.73.} (Well-known problem). Does there exist an
infinite group all of whose proper subgroups have prime order?\hf {\sl
A.\,I.\,Starostin}\vs

Yes, there does (A.\,Yu.\,Olshanskii, {\it Algebra and Logic},
{\bf 21} (1982), 369--418). \emp

\bmp \textbf{2.75.}
 Let $G$ be a periodic group containing an infinite family of
finite subgroups whose intersection contains non-trivial
elements. Does $G$ contain a non-trivial element with infinite
centralizer?
\hf {\sl S.\,P.\,Strunkov}\vs

Not always (K.\,I.\,Lossov, Dep. no.\,5528-V88, VINITI, Moscow,
1988 (Russian)).

\emp \bmp \textbf{2.76.}
 Let $\G$ be the holomorph of an abelian group $A$. Find
conditions for $A$ to be maximal among the locally nilpotent
subgroups of $\G$.\hf {\sl D.\,A.\,Suprunenko}\vs

These are found (A.\,V.\,Yagzhev, {\it Mat. Zametki}, {\bf 46},
no.\,6 (1989), 118 (Russian)).

\emp \bmp \textbf{2.77.}
 Let $A$ and $B$ be abelian groups. Find conditions under which
every extension of $A$ by $B$ is nilpotent.\hf {\sl
D.\,A.\,Suprunenko}\vs

These are found (A.\,V.\,Yagzhev, {\it Math. Notes}, {\bf
43} (1988), 244--245).

\emp \bmp \textbf{2.79.}
 Do there exist divisible (simple) groups with maximal
subgroups?

\hf {\sl M.\,S.\,Tsalenko}\vs

Yes, there do  (V.\,G.\,Sokolov, {\it Algebra
Logic}, {\bf 7} (1968), 122--126).

\emp \bmp \textbf{2.83.} Suppose that a periodic group $G$ is the
product of two locally finite subgroups. Is then $G$ locally
finite?\hf {\sl V.\,P.\,Shunkov}\vs

Not always (V.\,I.\,Sushchanski\u{\i}, {\it Math. USSR--Sb.}, {\bf
67}, no.\,2 (1990), 535--553).

\emp \bmp \textbf{2.88.}
 Is every Hall $\pi$-sub\-group of an arbitrary group a maximal
$\pi$-sub\-group?

\hf {\sl M.\,I.\,\`Eidinov}\vs

Not always (R.\,Baer, {\it J. Reine Angew. Math.}, {\bf 239/240} (1969),
109--144).

\emp

\bmp \textbf{3.2.}
Classify the faithful irreducible
(infinite-dimensional) representations of the nilpotent group
defined by generators $a,b,c$ and relations $[a,b]=c$, $ac=ca$,
$bc=cb$. (A condition for a representation to be monomial is given
in (A.\,E.\,Zalesski\u{\i}, {\it Math. Notes}, \textbf{9} (1971),
117--123).) 
\hfill {\sl S.\,D.\,Berman, A.\,E.\,Zalesski\u{\i}}\vs

Classification of these representations up to equivalence does not seem to be feasible. D.\,Segal ({\it Math. Proc. Cambridge Phil. Soc.}, {\bf 81} (1977), 201--208) proved that there exist primitive irreducible representations of this group (that is, representations that are not induced from a representation of any proper subgroup). There are results concerning classification of primitive ideals of the groups algebra of this group (over various fields). For instance, Zalesskii (ibid.) showed that every primitive ideal is maximal.
\emp

\bmp \textbf{3.4.}
 (Well-known problems). a) Is there an algorithm that decides,
for any set of group words $f_1,\ldots ,f_m$ (in a fixed set of
variables $ x_1,x_2,\ldots $) and a separate word $f$, whether
$f=1$ is a consequence of $ f_1=1,\ldots , f_m=1$?

\mb{b)} Given words $f_1,\ldots ,f_m$, is there an algorithm that
decides, for any word $f$, whether $ f = 1$ is a consequence of
$ f_1=1,\ldots , f_m=1$?\hf {\sl L.\,A.\,Bokut'}\vs

No, in both cases (Yu.\,G.\,Kleiman, {\it Soviet Math. Dokl.},
{\bf 20} (1979), 115--119).

\emp

\markboth{\protect\vphantom{(y}{Archive of solved
problems (3rd ed., 1969)}}{\protect\vphantom{(y}{Archive of solved
problems (3rd ed., 1969)}}

\bmp \textbf{3.6.}
 Describe the insoluble finite groups in which every soluble
subgroup is either $2$-closed, $2'$-closed, or isomorphic to $\Sy
_4$.\hf {\sl V.\,A.\,Vedernikov}\vs

This follows from (D.\,Gorenstein, J.\,H.\,Walter, {\it Illinois
J.~Math.}, {\bf 6} (1962) 553--593; \ V.\,D.\,Mazurov, {\it
Soviet Math. Dokl.}, {\bf 7} (1966), 681--682; \
V.\,D.\,Mazurov, V.\,M.\,Sitnikov, S.\,A.\,Syskin, {\it Algebra
and Logic}, {\bf 9} (1970),
187--204).

\emp \bmp \textbf{3.7.}
 An automorphism $\s$ of a group $G$ is called {\it algebraic\/} if, for
any $g\in G$, the minimal $\s$-inva\-ri\-ant subgroup of $G$
containing $g$ is finitely generated. An automorphism $\s$ of
$G$ is called an {\it $e$-auto\-mor\-phism\/} if, for any
$\s$-inva\-ri\-ant subgroups $A$ and $B$, where $A$ is a proper
subgroup of $B$, $A\setminus B$ contains an element $ x$ such
that $[x, \s ] \in B$ . Is every $e$-auto\-mor\-phism algebraic?
\hf {\sl V.\,G.\,Vilyatser}\vs

Yes, it is (A.\,V.\,Yagzhev, {\it Algebra and Logic}, {\bf 28}
(1989), 83--85).

\emp \bmp \textbf{3.8.}
 Let $G$ be the free product of free groups $A$ and $B$, and $V$
the verbal subgroup of~$G$ corresponding to the equation $ x^4 =
1$. Is it true that, if $ a\in A\setminus V$ and $ b\in B\setminus
V$, then $ (ab)^2\notin V$? \hf {\sl V.\,G.\,Vilyatser}\vs

No, it is not. If $A$ and $B$ are free groups with bases $a_1,a_2$
and $b_1,b_2$, respectively, then, for example, the commutators
$a=[[a_1,a_2,a_2],[a_1,a_2]]$ and $b=[[b_1,b_2,b_2],[b_1,b_2]]$ do
not belong to $V$, but $(ab)^2\in V$. Indeed, it follows from
(C.\,R.\,B.\,Wright, {\it Pacif. J. Math.}, {\bf 11}, no.\,1
(1961), 387--394) that $a^2\in V$, $b^2\in V$, and $[a,b]\in V$.
(V.\,A.\,Roman'kov, {\it Talk at the seminar Algebra and Logic},
March, 17, 1970.)

\emp \bmp \textbf{3.9.} Does there exist an infinite periodic group
with a finite maximal subgroup?

\hf {\sl
Yu.\,M.\,Gorchakov}\vs

Yes, there does (S.\,I.\,Adian, {\it Proc. Steklov Inst. Math.},
{\bf 112} (1971), 61--69).

\emp \bmp \textbf{3.10.}
 Need the number of finite non-abelian simple groups contained
in a proper variety of groups be finite?\hf {\sl
Yu.\,M.\,Gorchakov}\vs

Yes, mod CFSG; see, for example, (G.\,A.\,Jones, {\it J. Austral.
Math. Soc.}, {\bf 17} (1974), 162--173).

\emp \bmp \textbf{3.11.}
 An element $g$ of a group $G$ is said to be {\it generalized
periodic\/} if there exist $ x_1,\ldots ,x_n\in G$ such that $
x_1^{-1}gx_1\cdots x_n^{-1}gx_n= 1$. Does there exist a finitely
generated torsion-free group all of whose elements are
generalized periodic?\hf {\sl Yu.\,M.\,Gorchakov}\vs

Yes, there does (A.\,P.\,Goryushkin, {\it Siberian Math.~J.}, {\bf
14} (1973), 146--148). Another example is $ G=\left< a,b\mid
(b^2)^{a}=b^{-2},\;\; (a^2)^{b}=a^{-2} \right>$. Then $N=\left<
a^2,\,b^2,\,(ab)^2\right>$ is an abelian normal subgroup of $G$
and $G/N$ is non-cyclic of order $4$. If $
a^{2l}b^{2m}(ab)^{2n}=1$, then after conjugating by $a$ we get $
a^{2l}b^{-2m}(ab)^{-2n}=1$, whence $a^{4l}=1$. In view of the
obvious homomorphism $ G\rightarrow \Z ({\rm Aut\,}\Z)$ that maps
$a$ to the number $1\in \Z$ we have $l=0$. Similarly, $m=n=0$.
Hence $N$ is free abelian of rank $3$. The squares of elements
outside of $N$ are non-trivial; for example, $\left(
aa^{2l}b^{2m}(ab)^{2n}\right)^2=
a^2a^{-1}(a^{2l}b^{2m}(ab)^{2n})aa^{2l}
b^{2m}(ab)^{2n}=a^{4l+2}\ne 1$. Hence $G$ is torsion-free. Since
the square of any element $x\in G$ belongs to $N$, we have $
x^2(x^2)^a(x^2)^b(x^2)^{ab}=1$. (V.\,A.\,Churkin,
 1973.)

\emp \bmp \textbf{3.13.}
 Is the elementary theory of subgroup lattices of finite abelian
groups decidable?

\hf {\sl Yu.\,L.\,Ershov,~A.\,I.\,Kokorin}\vs

No (G.\,T.\,Kozlov, {\it Algebra and Logic}, {\bf 9} (1970),
104--107).

\emp \bmp \textbf{3.14.}
 Let $G$ be a finite group of $ n\times n$ matrices over a
skew field $T$ of characteristic zero. Prove that $G$ has a
soluble normal subgroup $H$ whose index in $G$ is not greater
than some number $f(n)$ not depending on $G$ or $T$. This
problem is related to the representation theory of finite groups
(the Schur index).
\hf {\sl A.\,E.\,Zalesski\u{\i}}\vs

This is proved mod CFSG (B.\,Hartley, M.\,A.\,Shahabi Shojaei, {\it
Math. Proc. Cambridge Phil. Soc.}, {\bf 92} (1982), 55--64).

\emp

\bmp \textbf{3.15.}
 A group $G$ is said to be {\it $U$-embed\-dable\/} in a class $
{\frak K}$ of groups if, for any finite submodel $M\subset G$,
there is a group $ A\in {\goth K}$ such that $M$ is isomorphic
to some submodel of $A$. Are the following groups $U$-embed\-dable
in the class of finite groups:

\mb{a)} every group with one defining relation?

\mb{b)} every group defined by one relation in the variety of
soluble groups of a given derived length? \hf {\sl
M.\,I.\,Kargapolov}\vs

No, in both cases (A.\,I.\,Budkin, V.\,A.\,Gorbunov, {\it Algebra
and Logic}, {\bf 14} (1975), 73--84).
 Another example. The group $ G=\left< a,b\mid
(b^2)^a=b^3\right>$ is non-Hopfian as proved in (G.\,Baumslag,
D.\,Solitar, {\it Bull. Amer. Math. Soc.}, {\bf 68}, no.\,3
(1962), 199--201) and therefore is not metabelian. We choose
elements $a,b,x_1,x_2,x_3,x_4\in G$ such that
$w=[[x_1,x_2],\,[x_3,x_4]]\ne 1$, add to them 1, their inverses,
and all the initial segments of the word $w$ in the alphabet $\{
x_i\}$, and all the initial segments of the word
$a^{-1}b^2ab^{-3}$ and of each of the words $x_i$ in the alphabet
of $\{ a,b\}$. Let $M$ be the resulting model. If $M$ was
embeddable into a finite group $G_0$, then $G_0$ would have to be
metacyclic, and also would have to have elements satisfying
$[[x_1,x_2],\,[x_3,x_4]]\ne 1$, which is impossible. Hence $G$ is
not $U$-embed\-dable into finite groups. Since the group
$G/G^{(k)}$ is also non-Hopfian for a suitable $k$ ({\it
ibid.\/}), it is not $U$-embed\-dable into finite groups for
similar reasons. (Yu.\,I.\,Merzlyakov, 1969.)

\emp 

 \bmp \textbf{3.17.}
 a) Is a non-abelian group with a unique linear ordering
necessarily simple?

\mb{b)} Is a non-commutative orderable group simple if it
has no non-trivial normal relatively convex subgroups?
\hf {\sl A.\,I.\,Kokorin}\vs

Not always, in both cases (V.\,V.\,Bludov, {\it Algebra and
Logic}, {\bf 13} (1974), 343--360).

\emp \bmp \textbf{3.18.}
 (B.\,H.\,Neumann). A group $U$ is called {\it universal\/} for a
class ${\goth K}$ of groups if $U$ contains an isomorphic image
of every member of ${\goth K}$. Does there exist a countable
group that is universal for the class of countable orderable
groups?\hf {\sl A.\,I.\,Kokorin}\vs

No (M.\,I.\,Kargapolov, {\it Algebra and Logic}, {\bf 9} (1970),
257--261; \ D.\,B.\,Smith, {\it Pacific J. Math.}, {\bf 35}
(1970), 499--502).

\emp \bmp \textbf{3.19.}
 (A.\,I.\,Mal'cev). For an arbitrary linearly ordered group $G$,
does there exist a linearly ordered abelian group with the same
order-type as $G$?\hf {\sl A.\,I.\,Kokorin}\vs

No (W.\,C.\,Holland, A.\,H.\,Mekler, S.\,Shelah, {\it Order}, {\bf 1}
(1985), 383--397).

\emp

\bmp \textbf{3.21.}
 Let ${\goth K}$ be the class of one-based models of signature
$\s$, \ $\vartheta$ a~property that makes sense for models in
${\goth K}$, and ${\goth K}_0$ the class of two-based models whose
first base is a set $M$ taken from ${\goth K}$ and whose second
base consists of all submodels of $M$ with property $\vartheta$
and whose signature consists of the symbols in $\s$ together with
$\in$ and $\subseteq$ in the usual set-theoretic sense. The
elementary theory of ${\goth K}_0$ is called the {\it
element-$\vartheta$-sub\-model theory\/} of ${\goth K}$. Are
either of the following theories decidable:

\mb{a)} the element-pure-subgroup theory of abelian groups?

\mb{b)} the element-pure-subgroup
theory of abelian torsion-free groups?

\mb{c)} the element-$\vartheta$-sub\-group theory of abelian
groups, when the set of $\vartheta$-sub\-groups is linearly
ordered by inclusion?\hf {\sl A.\,I.\,Kokorin}\vs

No, in all cases (for a), c):~G.\,T.\,Kozlov, {\it Algebra and
Logic}, {\bf 9} (1970), 104--107, \ {\it Algebra}, no.\,1, Irkutsk
Univ., 1972, 21--23 (Russian); \ for b):~\'E.\,I.\,Fridman, {\it
Algebra}, no.\,1, Irkutsk Univ., 1972 , 97--100 (Russian)).

\emp

\bmp \textbf{\zv 3.22.}
 Let $\xi = \{ G_{\alpha}, \;\pi_{\beta}^{\alpha}\mid \alpha
,\beta \in I\}$ be a projective system (over a directed set $I$)
of finitely generated free abelian groups. If all the
projections $
\pi_{\beta}^{\alpha}$
 are epimorphisms and all the $G_{\alpha}$ are non-zero, does it
follow that $ \lim\limits _{\longleftarrow}\xi \ne 0$? Equivalently,
suppose every finite set of elements of an abelian group $A$ is
contained in a pure finitely-generated free subgroup of $A$.
Then does it follow that $A$ has a direct summand isomorphic to
the infinite cyclic group?\hf {\sl V.\,I.\,Kuz'minov}
\ul

\otv
Under the assumption of the continuum hypothesis, not always
(E.\,A.\,Palyutin, {\it Siberian Math.~J.}, {\bf 19} (1978),
1415--1417). Another example,
not using the continuum hypothesis, is given by the abelian group $A = \mathbb{Z}^{\aleph_1}/N$ where $N$ is the subgroup consisting of elements having countable support. Every countable subgroup of $A$ is free abelian, which means that each finite subset of $A$ is contained in a pure finitely generated free abelian subgroup. But $A$ has no direct summand isomorphic to $\mathbb{Z}$. Both properties can be verified using the standard ZFC axioms of set theory. (S.\,Corson, {\it Letter of 26 August 2024}.)
\emp

 \bmp \textbf{3.26.}
 (F.\,Gross). Is it true that finite groups of exponent
$p^{\alpha}q^{\beta}$ have nilpotent length $\leq\alpha
+\beta$?\hf {\sl V.\,D.\,Mazurov}\vs

No, not always (E.\,I.\,Khukhro, {\it Algebra and Logic}, {\bf 17} (1978), 473--482).

\emp \bmp \textbf{3.27.}
 (J.\,G.\,Thompson). Is every finite simple group with a nilpotent
maximal subgroup isomorphic to some $PSL_2(q)$? \hfill {\sl
V.\,D.\,Mazurov}\vs

Yes, it is (B.\,Baumann, {\it J.~Algebra}, {\bf 38} (1976),
119--135).

\emp \bmp \textbf{3.28.}
 If $G$ is a finite 2-group with cyclic centre and every abelian
normal subgroup 2-generated, then is every abelian subgroup
of $G$ \ $3$-gene\-ra\-ted?\hf {\sl V.\,D.\,Mazurov}\vs

Not always (Ya.\,G.\,Berkovich, {\it Algebra and Logic}, {\bf 9}
(1970), 75).

\emp \bmp \textbf{3.29.}
 Under what conditions can a wreath product of matrix groups
over a field be represented by matrices over a field?
\hf {\sl Yu.\,I.\,Merzlyakov}\vs

Conditions have been found (Yu.\,E.\,Vapne, {\it Soviet Math.
Dokl.}, {\bf 11} (1970), 1396--1399).

\emp \bmp \textbf{3.30.}
 A torsion-free abelian group is called {\it
factor-decomposable\/} if, in all its factor groups, the
periodic part is a direct summand. Characterize these groups.

\hf
{\sl A.\,P.\,Mishina}\vs

This is done (L.\,Bican, {\it Commentat. Math. Univ. Carolinae}, {\bf
19} (1978), 653--672).

\emp \bmp \textbf{3.31.}
 Find necessary and sufficient conditions under which every pure
subgroup of a completely decomposable torsion-free abelian group
is itself completely decomposable.

\hfill {\sl A.\,P.\,Mishina}\vs

These have been found (L.\,Bican, {\it Czech. Math.~J.}, {\bf 24}
(1974), 176--191; \ A.\,A.\,Krav\-chen\-ko, {\it Vestnik Moskov.
Univ. Ser. 1 Mat. Mekh.}, {\bf 1980}, no.\,3, p.\,104 (Russian)).

\emp

\bmp
\textbf{3.33.} Are two groups
necessarily isomorphic if each of them can be defined by a single
relation and is a homomorphic image of the other one? \hfill
 {\sl D.\,I.\,Moldavanski\u{\i}}

\vs

No, not necessarily
(A.\,V.\,Borshchev, D.\,I.\,Moldavanski\u{\i}, {\it Math. Notes},
\textbf{79}, no.\,1 (2006), 31--40). 
\emp

\bmp \textbf{3.35.}
 (K.\,Ross). Suppose that a group $G$ admits two topologies $\s$
and $\tau$ yielding locally compact topological groups $G_{\s}$ and $
G_{\tau}$. If the sets of closed subgroups in $G_{\s}$ and $G_{\tau}$
are the same, does it follow that $G_{\s}$ and $G_{\tau}$ are
topologically isomorphic?

\hf {\sl Yu.\,N.\,Mukhin}\vsh

Not always (A.\,I.\,Moskalenko, {\it Ukrain. Math.~J.}, {\bf
30} (1978), 199--201).

\emp \bmp \textbf{3.37.}
 Suppose that every finitely generated subgroup of a locally
compact group $G$ is pronilpotent. Then is it true that every
maximal closed subgroup of $G$ contains the derived subgroup
$G'$?
\hf {\sl Yu.\,N.\,Mukhin}\vs

Yes, it is (I.\,V.\,Protasov, {\it Soviet Math. Dokl.}, {\bf 19}
(1978), 1208--1210).
\emp

\bmp \textbf{3.39.} Describe the finite groups with a
self-centralizing subgroup of prime order.

\hf {\sl
V.\,T.\,Nagrebetski\u{\i}}\vs

They are described. Self-centralizing subgroups of prime order are
CC-subgroups. Finite groups with a CC-subgroup were fully classified
in (Z.\,Arad, W.\,Herfort, {\it Commun. in
 Algebra}, {\bf 32} (2004), 2087--2098).
\emp

\bmp \textbf{3.40.}
 (I.\,R.\,Shafarevich). Let $SL_2(\Z )\,\hat{}$ and $SL_2(\Z )^{-}$
denote the completions of $SL_2(\Z )$ determined by all subgroups
of finite index and all congruence subgroups, respectively, and
let $\psi : SL_2(\Z )\,\hat{}\rightarrow SL_2(\Z )^{-}$ be the
natural homomorphism. Is ${\rm Ker}\, \psi$ a free profinite
group? \hf {\sl V.\,P.\,Platonov}\vs

Yes, it is (O.\,V.\,Mel'nikov, {\it Soviet Math. Dokl.}, {\bf
17} (1976), 867--870).
\emp

\bmp \textbf{3.41.} Is every compact periodic group locally finite?
\hf {\sl V.\,P.\,Platonov} \vs

Yes, it is (E.\,I.\,Zel'manov, {\it Israel J.~Math.}, {\bf 77}
(1992), 83--95).
\emp

\bmp \textbf{3.42.} (Kneser--Tits conjecture). Let $G$ be a simply connected $k$-defined simple algebraic group, and $E_k(G)$ the subgroup generated by unipotent $k$-elements. If $E_k(G)\ne 1$, then $G_k=E_k(G)$. The proof is known for $k$-decomposable groups (C.\,Chevalley) and for local fields (V.\,P.\,Platonov).
\hf {\sl V.\,P.\,Platonov} \vs

The conjecture was refuted (V.\,P.\,Platonov, \emph{Math. USSR--Izv.}, {\bf 10}, no.\,2 (1976), 211--243).
\emp

\bmp \textbf{3.50.} Let $G$ be a group of order $p^{\alpha} \cdot
m$, where $p$ is a prime, $p$ and $m$ are coprime, and let $k$ be
an algebraically closed field of characteristic $p$. Is it true
that if the indecomposable projective module corresponding to the
1-representation of $G$ has $k $-dimen\-sion $p^{\alpha}$, then
$G$ has a Hall $p' $-sub\-group? The converse is trivially true.

\hfill {\sl A.\,I.\,Saksonov}
\vs

No, not always: see Example 4.5 in (W.\,Willems, {\it Math. Z.}, {\bf 171} (1980), 163--174).
\emp

\bmp \textbf{3.51.}
Is it true that every finite group with a group of
automorphisms $\Phi$ which acts regularly on the set of conjugacy
classes of $G$ (that is, leaves only the identity class fixed) is
soluble? The answer is known to be affirmative in the case where
$\Phi$ is a cyclic group generated by a regular automorphism. \hfill
{\sl A.\,I.\,Saksonov}
\vs

No, not always (Y.\,Fine,
\textit{J.~Group Theory}, {\bf 22}, no.\,6 (2019),
1077--1087).

\emp

\bmp \textbf{3.52.}
 Can the quasivariety generated by the free group of rank 2 be
defined by a system of quasi-identities in finitely many
variables?\hf {\sl D.\,M.\,Smirnov}\vs

No (A.\,I.\,Budkin, {\it Algebra and Logic}, {\bf 15} (1976),
25--33).

\emp \bmp \textbf{3.53.}
 Let $L({\goth N}_4)$ denote the lattice of subvarieties of the
variety ${\goth N}_4$ of nilpotent groups of class at most $4$. Is
$L({\goth N}_4)$ distributive?\hf {\sl D.\,M.\,Smirnov}\vs

No (Yu.\,A.\,Belov, {\it Algebra and Logic}, {\bf 9} (1970), 371--374).

\emp \bmp \textbf{3.56.}
 Is a 2-group with the minimum condition for abelian subgroups
locally finite?

\hf {\sl S.\,P.\,Strunkov}\vsh

Yes, it is (V.\,P.\,Shunkov, {\it Algebra and Logic}, {\bf 9}
(1970), 291--297).

\emp \bmp \textbf{3.58.}
 Let $G$ be a compact $0$-dimen\-sio\-nal topological group all of
whose Sylow $p\hs$-sub\-groups are direct products of cyclic groups
of order $p$. Then is every normal subgroup of $G$ complementable?
\hf {\sl V.\,S.\,Charin}\vs

Yes, it is (M.\,I.\,Kabenyuk, {\it Siberian Math.~J.}, {\bf 13}
(1972), 654--657).

\emp \bmp \textbf{3.59.}
 Let $H$ be an insoluble minimal normal subgroup of a finite
group $G$, and suppose that $H$ has cyclic Sylow $p\hs$-sub\-groups for
every prime $p$ dividing $ |G : H|$. Prove that $H$ has at least
one complement in $G$.\hf {\sl L.\,A.\,She\-met\-kov}\vs

This was proved (S.\,A.\,Syskin, {\it Siberian Math.~J.},
{\bf 12} (1971), 342--344; \ L.\,A.\,She\-met\-kov, {\it Soviet Math.
Dokl.}, {\bf 11} (1970), 1436--1438). \emp

\bmp \textbf{3.61.} Let $\sigma$ be an automorphism of prime order
$p$ of a finite group $G$, which has a Hall $\pi $-sub\-group with
cyclic Sylow subgroups. Suppose that $p \in \pi$. Does the
centralizer $C_G(\sigma )$ have at least one Hall $\pi
$-sub\-group?\hf {\sl L.\,A.\,Shemetkov}\vs

Yes, it does, mod CFSG (V.\,D.\,Mazurov, {\it Algebra and Logic},
{\bf 31}, no.\,6 (1992), 360--366). \emp

\bmp \textbf{3.62.}
(Well-known problem). A finite group is said to
be a {\it $D_{\pi} $-group\/} if any two of its maximal $\pi
 $-sub\-groups
are conjugate. Is an extension of a $D_{\pi} $-group by a $D_{\pi
} $-group always a $D_{\pi} $-group?
 \hfill
{\sl L.\,A.\,Shemetkov}\vs

Yes, it is mod CFSG
 (E.\,P.\,Vdovin, D.\,O.\,Revin, {\it Contemp.
Math.}, {\bf 402} (2006), 229--263). \emp

\bmp \textbf{3.64.} Describe the finite simple groups with a Sylow
2-subgroup of the following type: $\langle a,t\mid a^{2^n}=t^2=1,\;\;
tat=a^{2^{n-1}-1}\rangle$, \ $n>2$.\hf {\sl V.\,P.\,Shunkov}\vs

This was done (J.\,L.\,Alperin, R.\,Brauer, D.\,Gorenstein,
{\it Trans. Amer. Math. Soc.}, {\bf 151} (1970), 1--261).
\emp

\parindent=0cm

\bmp \textbf{4.1.} Find an infinite finitely generated group with
an identical relation of the form $x^{2^n}=1$. \hfill
{\sl S.\,I.\,Adian}\vs

It has been found (S.\,V.\,Ivanov, {\it Int. J.~Algebra Comput.},
{\bf 4}, no.\,1--2 (1994), 1--308; I.\,G.\,Ly\-s\"e\-nok, {\it
Izv. Math.}, {\bf 60}, no.\,3 (1996), 453--654). \emp

\markboth{\protect\vphantom{(y}{Archive of solved
problems (4th ed., 1973)}}{\protect\vphantom{(y}{Archive of solved
problems (4th ed., 1973)}}

 \bmp \textbf{4.3.} Construct a finitely presented group with
insoluble word problem and satisfying a non-trivial law.\hf {\sl
S.\,I.\,Adian}\vs

This has been done (O.\,G.\,Kharlampovich, {\it Math. USSR--Izv.},
{\bf 19} (1982), 151--169).

\emp

\bmp \textbf{4.4.} Construct a finitely presented group with
undecidable word problem all of whose non-trivial defining
relations have the form $A^2=1$. This problem is interesting for
topologists.\hf {\sl S.\,I.\,Adian}\vs

It is constructed (O.\,A.\,Sarkisyan, {\it The word problem for
some classes of groups and semigroups}, Candidate disser., Moscow
Univ., 1983 (Russian)). \emp

\bmp \textbf{4.5.}
 a) (J.\,Milnor). Is it true that an arbitrary finitely generated
group has either polynomial or exponential growth?

\mb{c)} Is it true that every finitely generated group with
undecidable
word problem has exponential growth?\hf {\sl
S.\,I.\,Adian}\vs

a) No, c) No (R.\,I.\,Grigorchuk, {\it Soviet Math. Dokl.}, {\bf
28} (1983), 23--26).

\emp

\bmp \textbf{4.8.}
Suppose $G$ is a finitely generated free-by-cyclic
group. Is $G$ finitely presented?

\hfill {\sl G.\,Baumslag}
\vs

Yes, it is;
moreover, every finitely generated subgroup is finitely presented. For infinite cyclic extensions this is proved in (M.\,Feighn, M.\,Handel, {\it Ann. Math.}, {\bf 149}, no.\,3 (1999), 1061--1077), while another easier argument applies for finite cyclic extensions.
\emp

\bmp \textbf{4.10.}
 A group $G$ is called {\it locally indicable\/} if every
non-trivial finitely generated subgroup of $G$ has an infinite
cyclic factor group. Is every torsion-free one-relator group
locally indicable? \hf {\sl G.\,Baumslag}\vs

Yes, it is (S.\,D.\,Brodski\u{\i}, Dep. no.\,2214-80, VINITI,
Moscow, 1980 (Russian)).

\emp \bmp \textbf{4.12.}
 Let $G$ be a finite group and $A$ a group of
automorphisms of $ G$ stabilizing a series of subgroups
beginning with $G$ and ending with its Frattini subgroup. Then
is $A$ nilpotent?\hf {\sl Ya.\,G.\,Berkovich}\vs

Yes, it is (P.\,Schmid, {\it Math. Ann.}, {\bf 202} (1973), 57--69).

\emp

\bmp \textbf{4.14.}
Let $p$ be a prime number. What are necessary and
sufficient conditions for a finite group $G$ in order that the
group algebra of $G$ over a field of characteristic $p$ be
indecomposable as a two-sided ideal? There exist some nontrivial
examples, for instance, the group algebra of the Mathieu group
$M_{24}$ is indecomposable when $p$ is~$2$.

\hfill
{\sl R.\,Brauer}
\vs

A necessary and sufficient condition can be extracted from
(G.\,R.\,Robinson, {\it J. Algebra}, {\bf 84} (1983),
493--502); another solution, which requires considering fewer subgroups, can be extracted from
(B.\,K\"ulshammer, {\it Arch. Math. (Basel)}, {\bf 56} (1991), 313--319).
\emp

\bmp \textbf{4.16.}
 Suppose that ${\goth K}$ is a class of groups meeting the
following requirements: 1)~sub\-groups and epimorphic images of
${\goth K}$-groups are ${\goth K}$-groups; \ 2) if the group $ G =
UV$ is the product of $ {\goth K}$-sub\-groups $U$ and $V$ (neither
of which need be normal), then $ G\in {\goth K}$. If $\pi$ is a
set of primes, then the class of all finite $\pi$-groups meets
these requirements. Are these the only classes ${\goth K}$ with
these properties? \hf {\sl R.\,Baer}\vs

No (S.\,A.\,Syskin, {\it Siberian Math.~J.}, {\bf 20} (1979),
475--476).

\emp

\bmp \textbf{4.20.}
a) Let $F$ be a free group, and $N$ a normal subgroup of it. Is it true that the Cartesian square of $N$ is $m$-reducible to $N$ (that is, there is an algorithm that from a pair of words $w_1,w_2\in F$ constructs a word $w\in F$ such that
$w_1\in N$ and $w_2\in N$ if and only if $w\in N$)?

\makebox[15pt][r]{}b) (Well-known problem). Do there exist finitely presented groups in which the word problem has an arbitrary pre-assigned recursively enumerable $m$-degree of unsolvability? \hfill{\sl M.\,K.\,Valiev}\vs

a) No, it is not true (O.\,V.\,Belegradek, \emph{Siberian Math. J.}, {\bf 19} (1978), 867--870).
\vs

b) No, there exist $m$-degrees which do not contain the word problem of any recursively presented cancellation semigroup (C.\,G.\,Jockusch, jr., \emph{Z.~Math. Logik Grundlag. Math.}, {\bf 26}, no.\,1 (1980), 93--95).
\emp

\bmp \textbf{4.21.}
 Let $G$ be a finite group, $p$ an odd prime number, and $P$ a
Sylow $p\hs$-sub\-group of $G$. Let the order of every non-identity
normal subgroup of $G$ be divisible by $p$. Suppose $P$ has an
element $x$ that is conjugate to no other from $P$. Does $x$
belong to the centre of $G$? For $ p = 2$, the answer is positive
(G.\,Glauberman, {\it J.~Algebra}, {\bf 4} (1966), 403--420). \hfill {\sl
G.\,Glauberman}\vs

Yes, it does, mod CFSG (O.\,D.\,Artemovich, {\it Ukrain. Math.
J.}, {\bf 40} (1988), 343--345).

\emp

\bmp \textbf{4.22.} (J.\,G.\,Thompson). Let $G$ be a
finite group, $A$ a group of automorphisms of $G$ such that
$|A|$ and $|G|$ are coprime. Does there exist an $A
$-inva\-riant soluble subgroup $H$ of $G$ such that
$C_A(H)=1$?\hf {\sl G.\,Glauberman}\vs

Yes, it does (S.\,A.\,Syskin, {\it Siberian Math.~J.}, {\bf 32},
no.\,6 (1991), 1034--1037). \emp

\bmp \textbf{4.23.}
 Let $G$ be a finite simple group, $\tau$ some element of prime
order, and $\alpha$ an automorphism of $G$ whose order is
coprime to $ |G|$. Suppose $\alpha$ centralizes $C_G(\tau )$.
Is $\alpha = 1$? \hf {\sl G.\,Glauberman}\vs

Not always; for example, $G=Sz(8)$, \ $|\tau |=5$, \ $|\alpha
|=3$ \
(N.\,D.\,Podufalov, {\it Letter of September, 3, 1975\/}).

\emp \bmp \textbf{4.24.}
 Suppose that $T$ is a non-abelian Sylow 2-subgroup of a finite
simple group $G$.

\mb{b)} Is it possible that $T$ is the direct product of two proper
subgroups?

\mb{c)} Is $T'=\Phi (T)$? \hf {\sl D.\,Goldschmidt}\vs

b) Yes, it is. For example, the Sylow 2-subgroups of the
alternating groups $ {\Bbb A}_{14}$ and $ {\Bbb A}_{15}$ are
isomorphic to the Sylow 2-subgroup of $ {\Bbb S}_{4} \times {\Bbb
S}_{8}$. The Sylow 2-subgroups of $D_4(q)$ for $q$ odd are also
decomposable into direct products. (A.\,S.\,Kondratiev, {\it Letter
of October, 13, 1977.\/})

 c) Not always; for example, for $G=PSL_3(q)$ with $ q\equiv 1\,({\rm
mod\,}4)$ (A.\,S.\,Kondratiev).

\emp

\bmp \textbf{4.27.} Describe all finite
simple groups $G$ which
can be
represented in the form $G=ABA$, where $A$ and $B$ are abelian
subgroups. \hf {\sl I.\,P.\,Doktorov}\vs

They are described mod CFSG (D.\,L.\,Zagorin, L.\,S.\,Kazarin,
{\it Dokl. Math.}, {\bf 53}, no.\,2 (1996), 237--239; \
D.\,L.\,Zagorin, {\it Some problems of $ABA$-factorizations of
permutation groups and classical groups $($Cand. Diss.}),
Yaroslavl', 1994 (Russian); \ D.\,L.\,Zagorin, L.\,S.\,Kazarin,
 {\it Problems in
algebra}, Gomel' Univ., Gomel', no.\,11 (1997), 27--41 (Russian);
\ E.\,P.\,Vdovin, {\it Algebra and Logic}, {\bf 38}, no.\,2
(1999), 67--83). \emp

\bmp \textbf{4.28.} For a given field $k$ of characteristic $p > 0$,
characterize the locally finite groups with semisimple group
algebras over $k$. \hf {\sl A.\,E.\,Zalesskii}\vs

They are characterized: simple groups in (D.\,S.\,Passman,
A.\,E.\,Zalesskii, {\it Proc. London Math. Soc. (3)}, {\bf 67}
(1993), 243--276); the general case in (D.\,S.\,Passman, {\it Proc.
London Math. Soc. (3)}, {\bf 73} (1996), 323--357). \emp

\bmp \textbf{4.29.}
 Classify the irreducible matrix groups over a finite field that
are generated by reflections, that is, by matrices with Jordan
form ${\rm diag}(-1, 1, \ldots , 1)$.\hf {\sl
A.\,E.\,Zalesski\u{\i}}\vs

These are classified (A.\,E.\,Zalesski\u{\i}, V.\,N.\,Ser\"ezhkin, {\it Math.
USSR--Izv.}, {\bf 17} (1981), 477--503; \ A.\,Wagner, {\it Geom.
Dedicata}, {\bf 9} (1980), 239--253; {\bf 10} (1981), 191--203, 475-523).

\emp \bmp \textbf{4.32.}
 The conjugacy problem for metabelian groups.\hf {\sl
M.\,I.\,Kargapolov}\vs

This is algorithmically soluble (G.\,A.\,Noskov, {\it Math.
Notes}, {\bf 31} (1982), 252--258).

\emp \bmp \textbf{4.35.}
 (Well-known problem). Is there an infinite locally finite
simple group satisfying the minimum condition for $p\hs$-sub\-groups
for every prime~$p$?\hf {\sl O.\,H.\,Kegel}\vs

No (V.\,V.\,Belyaev, {\it Algebra and Logic}, {\bf 20} (1981),
393--401).

\emp \bmp \textbf{4.36.}
 Is there an infinite locally finite simple group $G$ with an
involution $i$ such that the centralizer $C_G(i)$ is a Chernikov
group?
\hf {\sl O.\,H.\,Kegel}\vs

No (A.\,O.\,Asar, {\it Proc. London Math. Soc.}, {\bf 45} (1982), 337--364).

\emp \bmp \textbf{4.37.}
 Is there an infinite locally finite simple group $G$ that
cannot be represented by matrices over a field and is such that
for some prime $p$ the $p\hs$-sub\-groups of $G$ are either of
bounded derived length or of finite exponent?\hf {\sl
O.\,H.\,Kegel}\vs

No (mod CFSG) by Theorem 4.8 of (O.\,H.\,Kegel,
B.\,A.\,F.\,Wehrfritz, {\it Locally
finite groups}, North Holland, Amsterdam, 1973).

\emp \bmp \textbf{4.38.}
 Classify the composition factors of automorphism groups of
finite (nilpotent of class 2) groups of prime exponent.\hf {\sl
V.\,D.\,Mazurov}\vs

Any finite simple group can be such a composition factor
(P.\,M.\,Beletski\u{\i}, {\it Russian Math. Surveys}, {\bf 33},
no.\,6 (1980), 85--86).

\emp \bmp \textbf{4.39.}
 A countable group $U$ is said to be {\it ${\sc
SQ}$-uni\-ver\-sal\/} if
every countable group is isomorphic to a subgroup of a quotient
group of~$ U$. Let $G$ be a group that has a presentation with $ r
\geq 2$ generators and at most $r-2$ defining relations. Is $G$ \
${\sc SQ}$-uni\-ver\-sal? \hfill {\sl
A.\,M.\,Macbeath,~P.\,M.\,Neumann} \vs

Yes, it is (B.\,Baumslag, S.\,J.\,Pride, {\it J. London Math. Soc.}, {\bf
17} (1978), 425--426).

\emp \bmp \textbf{4.41.}
 A point $z$ in the complex plane is called {\it free\/} if the
matrices {\small $\left(
\begin{array}{cc}1&2\\0&1\end{array}\right)$} and
 {\small $\left( \begin{array}{cc}1&0\\z&1\end{array}\right)$} generate
a
free group. Are all points outside the rhombus with vertices
$\pm 2$, $\pm i$ free?\hf {\sl Yu.\,I.\,Merzlyakov}\vs

No (Yu.\,A.\,Ignatov, {\it Math. USSR--Sb.}, {\bf 35} (1979),
49--56; \ A.\,I.\,Shkuratski\u{\i}, {\it Math. Notes}, {\bf 24}
(1978), 720--721). \emp

\bmp \textbf{4.45.}
Let $G$ be a free product amalgamating proper
subgroups $H$ and $K$ of $A$ and~$B$, respectively.

\mb{a)} Suppose that $H,K$ are finite and $|A:H| > 2$, \ $|B:K| \geq
2$. Is $G$ \, ${\sc SQ}\hskip0.1ex $-uni\-ver\-sal?

\mb{b)} Suppose that $A,B,H,K$ are free groups of finite ranks. Can $G$
be simple?

 \hf {\sl P.\,M.\,Neumann}\vs

a) Yes, it is (K.\,I.\,Lossov, {\it Siberian Math.~J.}, {\bf 27},
no.\,6 (1986), 890--899).

b) Yes, it can (M.\,Burger, S.\,Mozes, {\it Inst. Hautes \'Etudes Sci. Publ. Math.}, {\bf 92} (2000), 151--194). But if one of the indices
$|A:H|$, $|B:K|$ is infinite, then no, it cannot (S.\,V.\,Ivanov,
P.\,E.\,Schupp, in: {\it Algorithmic problems in groups and
semigroups, Int. Conf., Lincoln, NE, 1998}, Boston MA,
Birkh\"auser, 2000, 139--142).

\emp

\bmp \textbf{4.46.}
b) We call a variety of groups a {\it limit
variety\/} if it cannot be defined by finitely many laws, while
each of its proper subvarieties has a finite basis of identities.
It follows from Zorn's lemma that every variety that has no finite
basis of identities contains a limit subvariety. Is the set of
limit varieties countable? \hfill {\sl
A.\,Yu.\,Olshanskii}\vs

No, there are continuum of such varieties (P.\,A.\,Kozhevnikov,
{\em On varieties of groups of large odd exponent},
Dep.\,1612-V00, VINITI, Moscow, 2000 (Russian); \ S.\,V.\,Ivanov,
A.\,M.\,Storozhev, {\it Contemp. Math.}, {\bf 360} (2004),
55--62). \emp

\bmp \textbf{4.47.}
 Does there exist a countable family of groups such that every variety
is generated by a subfamily of it?
\hfill {\sl A.\,Yu.\,Olshanskii}\vs

No (Yu.\,G.\,Kleiman, {\it Math. USSR--Izv.}, {\bf 22} (1984),
33--65).

\emp \bmp \textbf{4.49.}
 Let $G$ and $H$ be finitely generated torsion-free nilpotent
groups such that ${\rm Aut}\,G\cong {\rm Aut}\,H$. Does it
follow that $G\cong H$?\hf {\sl V.\,N.\,Remeslennikov}\vs

No (G.\,A.\,Noskov, V.\,A.\,Roman'kov, {\it Algebra and
Logic}, {\bf 13} (1974), 306--311).

\emp \bmp \textbf{4.51.}
 (Well-known problem). Are knot groups residually finite?\hf
{\sl V.\,N.\,Remeslennikov}\vs

Yes, they are (W.\,P.\,Thurston, {\it Bull. Amer. Math. Soc.}, {\bf 6}
(1982), 357--381).

\emp \bmp \textbf{4.52.}
 Let $G$ be a finitely generated torsion-free group that is an
extension of an abelian group by a nilpotent group. Then is $G$
almost a residually finite $p\hs$-group for almost all primes $p$?\hf
{\sl V.\,N.\,Remeslennikov}\vs

Yes, it is (G.\,A.\,Noskov, {\it Algebra and Logic}, {\bf 13}
(1974), 388--393; \ D.\,Segal, {\it J.~Algebra}, {\bf 32}
(1974), 389--399).

\emp \bmp \textbf{4.53.}
 P.\,F.\,Pickel has proved that there are only finitely many
non-isomorphic finitely-generated nilpotent groups having the
same family of finite homomorphic images. Can Pickel's theorem
be extended to polycyclic groups?\hf {\sl
V.\,N.\,Remeslennikov}\vs

Yes, it can (F.\,J.\,Grunewald, P.\,F.\,Pickel, D.\,Segal, {\it Ann. of
Math. (2)}, {\bf 111} (1980), 155--195).

\emp \bmp \textbf{4.54.}
 Are any two minimal relation-modules of a finite group
isomorphic?

\hf {\sl K.\,W.\,Roggenkamp}\vsh

Not always (P.\,A.\,Linnell, {\it Ph.\,D. Thesis}, Cambridge,
1979; \newline
A.\,J.\,Sieradski, M.\,N.\,Dyer, {\it J. Pure Appl. Algebra},
{\bf 15} (1979), 199--217).

\emp \bmp \textbf{4.57.}
 Let a group $G$ be the product of two of its abelian minimax
subgroups $A$ and~$B$. Prove or refute the following statements:

\mb{a)} $A_0B_0\ne 1$, where $A_0=\bigcap_{x\in G}A^x$ and similarly
for $ B_0$;

\mb{b)} the derived subgroup of $G$ is a minimax subgroup.\hf {\sl
N.\,F.\,Sesekin}\vs

Both parts were proved (D.\,I.\,Zaitsev, {\it Algebra and
Logic}, {\bf 19} (1980), 94--106).

\emp \bmp \textbf{4.58.}
 Let a finite group $G$ be the product of two subgroups $A$ and
$B$, where $A$ is abelian and $B$ is nilpotent. Find the
dependence of the derived length of $G$ on the nilpotency class
of $B$ and the order of its derived subgroup.\hf {\sl
N.\,F.\,Sesekin}\vs

This was found (D.\,I.\,Zaitsev, {\it Math. Notes}, {\bf
33} (1983), 414--419).

\emp \bmp \textbf{4.59.}
 (P.\,Hall). Find the smallest positive integer $n$ such that
every countable group can be embedded in a simple group with $n$
generators.\hf {\sl D.\,M.\,Smirnov}\vs

It is proved that $ n = 2$ (A.\,P.\,Goryushkin, {\it Math. Notes}, {\bf 16}
(1974), 725--727).

\emp \bmp \textbf{4.60.}
 (P.\,Hall). What is the cardinality of the set of simple groups
generated by two elements, one of order 2 and the other of order
3?\hf {\sl D.\,M.\,Smirnov}\vs

The cardinality of the continuum. Every 2-generated group $G$ is
embeddable into a simple group $H$ with two generators of orders
$2$ and $3$ (P.\,E.\,Schupp, {\it J.~London Math. Soc.}, {\bf 13},
no.\,1 (1976), 90--94). Such a group $H$ has at most countably
many 2-generator subgroups, while there are continuum of
 groups~$ G$. (Yu.\,I.\,Merzlyakov, 1976.)

\emp \bmp \textbf{4.61.}
 Does there exist a linear function $f$ with the following
property: if every abelian subgroup of a finite 2-group $G$ is
generated by $n$ elements, then $G$ is generated by $f(n)$
elements?\hf {\sl S.\,A.\,Syskin}\vs

No (A.\,Yu.\,Olshanskii, {\it Math. Notes}, {\bf 23} (1978),
183--185).

\emp \bmp \textbf{4.62.}
 Does there exist a finitely based variety of groups whose
universal theory is undecidable? \hf {\sl A.\,Tarski}\vs

Yes, there does; for example, ${\goth A}^5$. Indeed, it is shown
in (V.\,N.\,Remeslennikov, {\it Algebra and Logic}, {\bf 12},
no.\,5 (1975), 327--346) that there exists a finitely presented
group $G=\left< x_1,\ldots ,x_n\mid w_1,\ldots ,w_m,\;\, {\rm
mod}\,{\goth A}^5\right> $ in ${\goth A}^5$ with undecidable word
problem. We put $\Phi _w=(\forall x_1,\ldots
,x_n)((w_1=1)\wedge\ldots \wedge(w_m=1)\rightarrow (w=1))$, where
$w$ runs over all the words in $x_1,\ldots ,x_n$. Clearly, there
is no algorithm to decide whether a formula $\Phi _w$ is true in
${\goth A}^5$. (V.\,N.\,Remeslennikov, 1976.)

\emp \bmp \textbf{4.63.}
 Does there exist a non-abelian variety of groups (in
particular, one that contains the variety of all abelian groups)
whose elementary theory is decidable? \hf {\sl A.\,Tarski}\vs

No (A.\,P.\,Zamyatin, {\it Algebra and Logic}, {\bf 27} (1978), 13--17).

\emp \bmp \textbf{4.64.}
 Does there exist a variety of groups that does not admit an independent
system of defining identities? \hf {\sl A.\,Tarski}\vs

Yes, there does (Yu.\,G.\,Kleiman, {\it Math. USSR--Izv.}, {\bf 22} (1984),
33--65).

\emp \bmp \textbf{4.67.}
 Let $G$ be a finite $p\hs$-group. Show that the rank of the
multiplicator $M(G)$ of $G$ is bounded in terms of the rank of
$G$. \hf {\sl J.\,Wiegold}\vs

This was done (A.\,Lubotzky, A.\,Mann, {\it J.~Algebra}, {\bf
105} (1987), 484--505).

\emp \bmp \textbf{4.68.}
 Construct a finitely generated (infinite) characteristically
simple group that is not a direct power of a simple group. \hf
{\sl J.\,Wiegold}\vs

This was done (J.\,S.\,Wilson, {\it Math. Proc. Cambridge
Phil. Soc.}, {\bf 80} (1976), 19--35).

\emp

\bmp \textbf{4.69.}
Let $G$ be a finite $p\hskip0.2ex $-group, and
suppose that $|G'|
> p^{n(n-1)/2}$ for some non-negative integer $n$. Prove that
$G$ is generated by the elements of breadths $\geq n$. The {\it
breadth\/} of an element $x$ of $G$ is $b(x)$ where
$|G:C_G(x)|=p^{b(x)}$. \hfill {\sl J.\,Wiegold}
\vs

This has been proved (A.\,Skutin, {\it J.~Algebra}, {\bf 526} (2019), 1--5).
\emp

 \bmp \textbf{4.70.}
 Let $k$ be a field of characteristic different from $2$, and
$G_k$ the group of transformations $A = (a, \alpha) : x
\rightarrow ax +\alpha$ \ ($a, \alpha \in k$, \ $a\ne 0$). Extend
$G_k$ to the projective plane by adjoining the symbols $ (0,
\alpha )$ and a line at infinity. Then the lines are just the
centralizers $C_G(A)$ of elements $A\in G_k$ and their cosets. Do
there exist other groups $G$ complementable to the projective
plane such that the lines are just the cosets of the centralizers
of elements of $G$? \hf {\sl H.\,Schwerdtfeger}\vs

No (E.\,A.\,Kuznetsov, Dep. no.\,7028-V89, VINITI, Moscow, 1989
(Russian)).

\emp \bmp \textbf{4.71.}
 Let $A$ be a group of automorphisms of a finite group $G$ which
has a series of $A$-inva\-ri\-ant subgroups $G = G_0 > \cdots > G_k =
1$ such that every $|G_i:G_{i+1}|$ is prime. Prove that $A$ is
supersoluble.\hf {\sl L.\,A.\,Shemetkov}\vs

This was proved (L.\,A.\,Shemetkov, {\it Math. USSR--Sb.},
{\bf 23} (1974), 593--611).

\emp \bmp \textbf{4.73.} (Well-known problem). Does there exist a
non-abelian variety of groups

 \mb{a)} all of whose finite groups are abelian?

\mb{b)} all of whose periodic groups are abelian? \hf {\sl
A.\,L.\,Shmel'kin}\vs

a) Yes, there does (A.\,Yu.\,Olshanskii, {\it Math.
USSR--Sb.}, {\bf 54} (1986), 57--80).

b) Yes, there does (P.\,A.\,Kozhevnikov, {\em On varieties of
groups of large odd exponent}, Dep.\,1612-V00, VINITI, Moscow,
2000 (Russian); \ S.\,V.\,Ivanov, A.\,M.\,Storozhev, {\it Contemp.
Math.}, {\bf 360} (2004), 55--62). \emp

\bmp \textbf{4.74.} a) Is every 2-group of order greater than 2 non-simple?
\hf {\sl V.\,P.\,Shunkov}\vs

No, the existence of simple infinite 2-groups follows from the
existence of finitely generated infinite groups of exponent $2^n$
(S.\,V.\,Ivanov, {\it Int. J.~Algebra Comput.}, {\bf 4}, no.\,1--2
(1994), 1--308; \ I.\,G.\,Lys\"enok, {\it Izv. Math.}, {\bf 60},
no.\,3 (1996), 453--654) and from the positive solution to the
Restricted Burnside Problem for groups of exponent~$2^n$
(E.\,I.\,Zel'manov, {\it Math. USSR--Sb.}, {\bf 72} (1992),
543--565). \emp

\bmp \textbf{4.76.}
 Let $G$ be a locally finite group containing an element $a$ of
prime order such that the centralizer $C_G(a)$ is finite. Is $G$
almost soluble?\hf {\sl V.\,P.\,Shunkov}\vs

Yes, it is. P.\,Fong ({\it Osaka J. Math.}, {\bf 13} (1976),
483--489) has shown (mod CFSG) that $G$ is almost locally soluble.
Given this, B.\,Hartley and T.\,Meixner ({\it Arch. Math.}, {\bf
36} (1981), 211--213) have proved that $G$ is almost locally
nilpotent. Therefore $G$ is almost soluble by (J.\,L.\,Alperin,
{\it Proc. Amer. Math. Soc.}, {\bf 13} (1962), 175--180) and
(G.\,Higman, {\it J. London Math. Soc.}, {\bf 32} (1957),
321--334). Moreover, by (E.\,I.\,Khukhro, {\it Math. Notes}, {\bf
38}, no.\,5-6 (1985), 867--870; \ E.\,I.\,Khukhro, {\it Math.
USSR--Sb.}, {\bf 71}, no.\,1 (1992), 51--63)
 then $G$ is almost nilpotent.

\emp \bmp \textbf{4.77.}
 In 1972, A.\,Rudvalis discovered a new simple group $R$ of order
$2^{14}\cdot 3^3 \cdot 5^3
\cdot 7\cdot 13\cdot 29$. He has shown that $R$ possesses an
 involution $i$ such that $C_R(i)=V \times F$, where $V$ is a
4-group (an elementary abelian group of order $4$) and $ F\cong
Sz(8)$.

\mb{a)} Show that $R$ is the only finite simple group $G$ that
possesses an involution $i$ such that $C_G(i)=V\times F$, where
$V$ is a 4-group and $F\cong Sz(8)$.

\mb{b)} Let $G$ be a non-abelian finite simple group that possesses
an involution $i$ such that $C_G(i)=V\times F$, where $V$ is an
elementary abelian 2-group of order $2^n$, $ n \geq 1$, and
$F\cong Sz(2^m)$, $ m
\geq 3$. Show that $ n = 2$ and $ m = 3$. \hf {\sl
Z.\,Janko}\vs

a) This has been shown (V.\,D.\,Mazurov, {\it Math. Notes}, {\bf
31} (1982), 165--173).

b) This follows from the CFSG.

\emp
 \bmp \textbf{5.1.} a) Is every locally finite minimal non-$FC
$-group non-simple?

\hf {\sl V.\,V.\,Belyaev,~N.\,F.\,Sesekin}\vs

Yes, it is (M.\,Kuzucuo\=glu, R.\,E.\,Phillips, {\it Proc.
Cambridge Philos. Soc.}, {\bf 105}, no.\,3 (1989), 417--420). \emp
\markboth{\protect\vphantom{(y}{Archive of solved
problems (5th ed., 1976)}}{\protect\vphantom{(y}{Archive of solved
problems (5th ed., 1976)}}
 \bmp

{\bf 5.2.} Is it true that the growth function $f$ of any infinite
finitely generated group satisfies the inequality $f(n)\leq
(f(n-1)+f(n+1))/2$ for all sufficiently
large $n$ (for a fixed finite system of
generators)?\hf {\sl V.\,V.\,Belyaev,~N.\,F.\,Sesekin}\vs

No (P.\,\,de\,\,la\,\,Harpe, R.\,I.\,Grigorchuk, {\it Algebra and
Logic}, {\bf 37}, no.\,6 (1998), 353--356). \emp

\bmp \textbf{5.3.}
 (Well-known problem). Can every finite lattice $L$ be embedded
in the lattice of subgroups of a finite group?\hf {\sl
G.\,M.\,Bergman}\vs

Yet, it can (P.\,Pudl\'ak, J.\,T\raisebox{-0.1ex}{$\stackrel{\circ}{{\rm
u}}$}ma, {\it Algebra Universalis}, {\bf 10} (1980), 74--95).

\emp

\bmp \textbf{5.4.} Let
$g$ and $h$ be positive elements of a linearly ordered group $G$.
Can one always embed $G$ in a linearly ordered group
$\,\overline{\! G}$ in such a way that $g$ and $h$ are conjugate
in~$\,\overline{\! G}$?
\hfill {\sl V.\,V.\,Bludov}

\vs

No, not always
(V.\,V.\,Bludov, {\it Algebra and Logic}, \textbf{44}, no.\,6
(2005), 370--380). \emp 

\bmp \textbf{5.6.}
 If $R$ is a finitely generated integral domain of
characteristic $ p > 0$, does the profinite (ideal) topology of
$R$ induce the profinite topology on the group of units of~$R$?
It does for $ p = 0$.\hf {\sl B.\,A.\,F.\,Wehrfritz}\vs

Yes, it does (D.\,Segal, {\it Bull. London Math. Soc.}, {\bf 11}
(1979), 186--190). \emp

 \bmp \textbf{5.7.} An
algebraic variety $X$ over a field $k$ is called {\it rational\/}
if the field of functions $k(X)$ is purely transcendental over
$k$, and it is called {\it stably rational\/} if $k(X)$ becomes
purely transcendental after adjoining finitely many independent
variables. Let $T$ be a stably rational torus over a field $k$. Is
$ T$ rational? Other formulations of this question and some
related results see in (V.\,E.\,Voskresenski\u{\i}, {\it Russ.
Math. Surveys}, {\bf 28}, no.\,4 (1973), 79--105). \hfill {\sl
V.\,E.\,Voskresenski\u{\i}}\vs

Yes, it is (V.\,E.\,Voskresenskii, {\it J. Math. Sci. (N. Y.)}, {\bf 161}, no.\,1 (2009), 176--180).
 \emp

\bmp \textbf{5.8.}
 Let $G$ be a torsion-free soluble group of finite cohomological
dimension ${\rm cd}\,G$. If $hG$ denotes the Hirsch number of
$G$, then it is known that $ hG\leq {\rm cd}\,G\leq hG+1$. Find
a purely group-theoretic criterion for $hG={\rm cd}\,G$.\hf {\sl
K.\,W.\,Gruenberg}\vs

This has been found (P.\,H.\,Kropholler, {\it J. Pure Appl. Algebra},
{\bf 43} (1986), 281--287). \emp

\bmp \textbf{5.9.}
 Let $1\rightarrow R_i \stackrel{\pi _i}{\longrightarrow}F\rightarrow G\rightarrow 1$, \
$i=1,\,2$, be two exact sequences of groups with $G$ finite and
$F$ free of finite rank $d(F)$. If we assume that $d(F)=d(G)+1$
(where $d(G)$ is the minimum number of generators of $G$), are
the corresponding abelianized extensions isomorphic?\hf {\sl
K.\,W.\,Gruenberg}\vs

Yes, they are (P.\,A.\,Linnell, {\it J. Pure Appl. Algebra}, {\bf 22}
(1981), 143--166).
\emp

\bmp \textbf{5.10.} Let $E=H*K$ and denote the augmentation ideals of
the groups $E$, $H$, $K$ by $ {\goth e}$, ${\goth h}$, ${\goth
k}$, respectively. If $I$ is a right ideal in ${\Bbb Z}E$, let
$d_E(I)$ denote the minimum number of generators of $I$ as a right
ideal. Assuming $H$ and $K$ finitely generated, is it true that
$d_E({\goth e})=d_E({\goth h}E) + d_E({\goth k}E)?$ Here ${\goth
h}E$ is, of course, the right ideal generated by $ {\goth h}$;
similarly for ${\goth k}E$.\hf {\sl K.\,W.\,Gruenberg}\vs

Not always (P.\,A.\,Linnell, {\it unpublished\/}). A simple
example: $H$ is elementary abelian group of order $4$, $K$ is
elementary abelian of order $9$. Here $d_E({\goth e})=3$, although
$d_E({\goth h}E)=d_E({\goth k}E)=2$. (K.\,W.\,Gruenberg, {\it
Letter of July, 27, 1995}.) \emp 

\bmp \textbf{5.11.}
 Let $G$ be a finite group and suppose that there exists a
non-empty proper subset $\pi$ of the set of all primes dividing
$|G|$ such that the centralizer of every non-trivial
$\pi$-ele\-ment is a $\pi$-sub\-group. Does it follow that $G$
contains a subgroup $U$ such that $U^g\cap U=1$ or $U$ for
every $g\in G$, and the centralizer of every non-trivial
element of~$U$ is contained in $U$?\hf {\sl
K.\,W.\,Gruenberg}\vs

Yes, it does (mod CFSG) (J.\,S.\,Williams, {\it J.~Algebra}, {\bf
69} (1981), 487--513). \emp

\bmp \textbf{5.12.} Let $G$ be a finite group with trivial soluble
radical in which there are Sylow 2-subgroups having non-trivial
intersection. Suppose that, for any two Sylow 2-subgroups $P$ and
$Q$ of $G$ with $P \cap Q \ne 1$, the index $|P:P \cap Q|$ does
not exceed $2^n$. Is it then true that $|P| \leq 2^{2n}$?\hf {\sl
V.\,V.\,Kabanov}\vs

Yes, it is, mod CFSG (V.\,I.\,Zenkov, {\it Algebra and Logic},
{\bf 36}, no.\,2 (1997), 93--98). \emp

\bmp \textbf{5.13.}
 Suppose that $K$ and $L$ are distinct conjugacy classes of
involutions in a finite group $G$ and $\left< x,y\right>$ is a
$2$-group for all $x\in K$ and $y\in L$. Does it follow that
$G\ne [K,L]$?\hf {\sl V.\,V.\,Kabanov}\vs

Not always; for example, $G=Sp_4(2^n)$, $n\geq 2$, $K,L$ being
classes of involutions that have non-trivial intersections with
the centres of $N_G(M)'$ and $N_G(N)'$, respectively, where $M,N$
are distinct elementary abelian subgroups of order 8 in a Sylow
$2$-sub\-group of $G$ and the dash means taking the derived subgroup
(A.\,A.\,Makhn\"ev, {\it Letter of October, 10, 1981\/}). \emp

\bmp \textbf{5.17.} If the finite group $G$ has the form $G=AB$
where $A$ and $B$ are nilpotent of classes $\alpha $ and $\beta $,
respectively, then $G$ is soluble. Is $G^{(\alpha +\beta )}=1$?
One can show that $G^{(\alpha +\beta )}$ is nilpotent
(E.\,Pennington, 1973). One can ask the same question for
infinite groups (or Lie algebras), but there
is nothing known beyond Ito's theorem: $A'=B'=1$ implies
$G^{(2)}=1$.\hf {\sl O.\,H.\,Kegel}\vs

Not always (J.\,Cossey, S.\,Stonehewer, {\it Bull. London Math.
Soc.}, {\bf 30}, no.\,144 (1998), 247--250). See also new problem
14.43. \emp

 \bmp \textbf{5.18.} Let $G$ be an infinite locally finite simple group. Is
the centralizer of every element of $G$ infinite?\hfill {\sl
O.\,H.\,Kegel}\vs

Yes, it is (B.\,Hartley, M.\,Kuzucuo\=glu, {\it Proc. London Math.
Soc. (3)}, {\bf 62}, no.\,2 (1991), 301--324). \emp

\bmp \textbf{5.19.}
 a) Let $G$ be an infinite locally finite simple group
satisfying the minimum condition for 2-subgroups. Is $ G =
PSL_2(F)$, $F$ some locally finite field of odd characteristic, if
the centralizer of every involution of $G$ is almost locally
soluble?

\mb{b)} Can one characterize the simple locally finite groups with
the min-2 condition containing a maximal radical non-trivial
2-subgroup of rank $\leq 2$ as linear groups of small rank?\hf
{\sl O.\,H.\,Kegel}\vs

a) Yes, it is (N.\,S.\,Chernikov, {\it Ukrain. Math.~J.}, {\bf
35} (1983), 230--231).

b) Yes, one can (mod CFSG) by Theorem 4.8 of (O.\,H.\,Kegel,
B.\,A.\,F.\,Wehrfritz, {\it Locally finite groups}, North
Holland, Amsterdam, 1973).

\emp

\bmp \textbf{5.20.}
Is the elementary theory of lattices of $l $-ideals of
lattice-ordered abe\-li\-an groups decidable?\hfill {\sl
A.\,I.\,Kokorin}\vs

No, it is undecidable (N.\,Ya.\,Medvedev, {\it Algebra and Logic},
{\bf 44}, no.\,5 (2005), 302--312). \emp

 \bmp \textbf{5.21.}
Can every torsion-free group with solvable word problem be
embedded in a group with solvable conjugacy problem? An example
due to A.\,Macintyre shows that this question has a negative
answer when the condition of torsion-freeness is omitted.

\hfill {\sl D.\,J.\,Collins}

No, not every
(A.\,Darbinyan, \emph{Invent. Math.}, \textbf{224} (2021), 987--997).
 \emp

\bmp \textbf{5.22.}
Does there exist a version of the Higman embedding theorem in
which the degree of unsolvability of the conjugacy problem is
preserved? 
\hfill {\sl D.\,J.\,Collins}\vs

Yes, there does (A.\,Yu.\,Olshanskii, M.\,V.\,Sapir, {\it The
conjugacy problem and Higman embeddings\/} ({\it Memoirs AMS},
{\bf 804}), 2004, 133\,p.). \emp

\bmp \textbf{5.23.} Is
it true that a free lattice-ordered group of the variety of
lattice-ordered groups defined by the law $x^{-1}|y|x \ll |y|^2$
(or, equivalently, by the law $|[x,y]| \ll |x|$), is residually
linearly ordered nilpotent?\hfill
 {\sl V.\,M.\,Kopytov}

\vs

No, it is not
(V.\,V.\,Bludov, A.\,M.\,W.\,Glass, {\it Trans. Amer. Math. Soc.},
\textbf{358} (2006), 5179--5192).

\emp

\bmp \textbf{5.24.}
Is it true that a free lattice-ordered group of
the variety of the lattice-ordered groups which are residually
linearly ordered, is residually soluble linearly ordered?

 \hfill
{\sl V.\,M.\,Kopytov}\vs

No, it is not (N.\,Ya.\,Medvedev, {\it Algebra and Logic}, {\bf
44}, no.\,3 (2005), 197--204).\emp

\bmp \textbf{5.28.}
 Let $G$ be a group and $H$ a torsion-free subgroup of $G$ such
that the augmentation ideal $I_G$ of the integral group-ring $\Z
G$ can be decomposed as $ I_G=I_H \Z G\oplus M$
for some $\Z
G$-sub\-mo\-dule $ M$. Prove that $G$ is a free product of the form
$ G = H*K$.

\hf {\sl D.\,E.\,Cohen}\vs

This has been proved (W.\,Dicks, M.\,J.\,Dunwoody, {\it Groups acting on
graphs}, Cambridge Univ. Press, Cambridge--New York, 1989).

\emp \bmp \textbf{5.29.}
 Consider a group $G$ given by a presentation with $m$
generators and $n$ defining relations, where $m\geq n$. Do some
$m-n$ of the given generators generate a free subgroup of
$G$?\hf {\sl R.\,C.\,Lyndon}\vs

Yes, they do (N.\,S.\,Romanovski\u{\i}, {\it Algebra and Logic},
{\bf 16} (1977), 62--68).
\emp

\bmp \textbf{5.32.} Let $p$ be a prime, $C$ a conjugacy class of $
p\hskip0.1ex $-ele\-ments of a finite group $G$ and suppose that
for any two elements $x$ and $y$ of $C$ the product $xy^{-1}$ is a
$p\hskip0.1ex $-ele\-ment. Is the subgroup generated by the class $C$ a
$p\hskip0.1ex $-group? \hf {\sl V.\,D.\,Mazurov}\vs

Yes, it is (for $p>2$: W.\,Xiao, {\it Sci. in China}, {\bf 34},
no.\,9 (1991), 1025--1031; \ for $p=2$: V.\,I.\,Zenkov, in: {\it
Algebra. Proc. IIIrd Int. Conf. on Algebra, Krasnoyarsk, 1993},
Berlin, De~Gruyter, 1996, 297--302). \emp

\bmp \textbf{5.34.}
 Let ${\frak o}$ be a commutative ring with identity in which 2 is
invertible and which is not generated by zero divisors. Do there
exist non--standard automorphisms of $GL_n({\frak o})$ for $ n \geq
3$?\hf {\sl Yu.\,I.\,Merzlyakov}\vs

No (V.\,Ya.\,Bloshchitsyn, {\it Algebra and Logic}, {\bf 17}
(1978), 415--417).

\emp \bmp \textbf{5.40.}
 Let $G$ be a countable group acting on a set $\W$. Suppose
that $ G$ is $k$-fold transitive for every finite $k$, and $G$
contains no non-trivial permutations of finite support. Is it
true that $\W$ can be identified with the rational line $\Q$ in
such a way that $G$ becomes a group of autohomeomorphisms?\hf
{\sl P.\,M.\,Neumann}\vs

Not always (A.\,H.\,Mekler, {\it J. London Math. Soc.}, {\bf 33} (1986),
49--58).
\emp

\bmp \textbf{5.41.} Does every non-trivial finite group, which is
free in some variety, contain a non-trivial abelian normal
subgroup? \hf {\sl A.\,Yu.\,Olshanskii}\vs

Yes, it does (mod CFSG). Otherwise, if $G$ is a counterexample,
the centralizer of the product $E$ of all minimal normal subgroups
of $G$ is trivial and there is an element of $G/E$ whose order is
$m$, where $m$ is the maximum of the orders of 2-elements of~$G$.
By Theorem~1 in (M.\,Aschbacher, P.\,B.\,Kleidman,
M.\,W.\,Liebeck, {\it Math.~Z.}, {\bf 208}, no.\,3 (1991),
401--409), which is proved using CFSG, there is an element of
order $2m$ in $G$, which contradicts the choice of $m$.
(S.\,A.\,Syskin, {\it Letter of August, 20, 1992}.) \emp

\bmp \textbf{5.43.}
 Does there exist a soluble variety of groups that is not
generated by its finite groups?\hf {\sl
A.\,Yu.\,Olshanskii}\vs

Yes, there does (Yu.\,G.\,Kleiman, {\it Trans. Moscow Math. Soc.},
{\bf 1983}, no.\,2, 63--110).

\emp \bmp \textbf{5.46.}
 Is every recursively presented soluble group embeddable in a
group finitely presented in the variety of all soluble groups of
derived length $n$, for suitable $n$?

\hf {\sl
V.\,N.\,Remeslennikov}\vsh

Yes, it is (G.\,P.\,Kukin, {\it Soviet Math. Dokl.}, {\bf 21}
(1980), 378--381).

\emp \bmp \textbf{5.49.} Let $A_{m,n}$ be the group of all
automorphisms of
the free soluble group of derived length $n$ and rank $m$. Is it
true that

\mb{a)} $ A_{m,n}$ is finitely generated for any $m$ and $n$?

\mb{b)} every automorphism in $A_{m,n}$ is induced by some
automorphism in $A_{m,n+1}$?

\hf {\sl V.\,N.\,Remeslennikov}\vs

a) No: $A_{3,2}$ is not finitely generated, although $A_{m,2}$ is
 whenever $ m\ne 3$ (S.\,Bachmuth, H.\,Y.\,Mochizuki, {\it Trans. Amer.
Math. Soc.}, {\bf 270} (1982), 693--700, \ {\bf 292} (1985),
81--101). \ \ \ \

b) No (V.\,Shpilrain, {\it Int. J.~Algebra and Comput.}, {\bf 1},
no.\,2 (1991), 177--184). \emp

\bmp \textbf{5.50.}
 Is there a finite group whose set of quasi-laws does not have
an independent basis?\hf {\sl D.\,M.\,Smirnov}\vs

Yes, there is (A.\,N.\,Fedorov, {\it Siberian Math.~J.}, {\bf
21} (1980), 840--850).

\emp \bmp \textbf{5.51.}
 Does a non-abelian free group have an independent basis for its
quasi-laws?

\hf {\sl D.\,M.\,Smirnov}\vsh

Yes, it does (A.\,I.\,Budkin, {\it Math. Notes}, {\bf 31} (1982),
413--417).

\emp

 \bmp \textbf{5.53.}
(P.\,Scott). Let $p,q,r$ be distinct prime numbers. Prove that the
free product $G=C_p*C_q*C_r$ of cyclics of orders $p,q,r$ is not
the normal closure of a single element. By Lemma 3.1 of
(J.\,Wiegold, {\it J.~Austral. Math. Soc.}, {\bf 17}, no.\,2
(1974), 133--141), every soluble image, and every finite image, of
$G$ is the normal closure of a single element. \hfill {\sl
J.\,Wiegold}\vs

This is proved in (J.\,Howie, {\it J.~Pure Appl. Algebra}, {\bf
173}, no.\,2 (2002), 167--176). \emp

\bmp \textbf{5.56.} b)
 Does there exist a locally nilpotent group of prime exponent that
coincides with its derived subgroup (and hence has no maximal
subgroups)?\hf {\sl J.\,Wiegold}\vs

Yes, there does. Every non-soluble variety ${\goth V}$ contains
a non-trivial group coinciding with the derived subgroup, the
direct limit of the spectrum $ F\stackrel{\f}{\rightarrow}
F\stackrel{\f}{\rightarrow} \ldots $, where $F$ is the free group
in ${\goth V}$ on the free generators $x_1,x_2,\ldots $ and $\f$
is the homomorphism given by $x_i\rightarrow [x_{2i-1},x_{2i}]$.
As shown in (Yu.\,P.\,Razmyslov, {\it Algebra and Logic}, {\bf 10}
(1971), 21--29) the Kostrikin variety of locally nilpotent groups
of prime exponent $ p\geq 5$ is unsoluble. (E.\,I.\,Khukhro,
I.\,V.\,L'vov, {\it Letter of June, 19, 1976.}) The same was also
proved in (Yu.\,A.\,Kolmakov, {\it Math. Notes}, {\bf 35},
no.\,5--6 (1984), 389--391). An example answering the question was
also produced in (M.\,R.\,Vaughan-Lee, J.\,Wiegold, {\it Bull.
London Math. Soc.}, {\bf 13}, no.\,1 (1981), 45--46).

\emp \bmp \textbf{5.60.}
 Is an arbitrary soluble group that satisfies the minimum condition for
normal subgroups countable?\hf {\sl B.\,Hartley}\vs

Not always (B.\,Hartley, {\it Proc. London Math. Soc.}, {\bf 33} (1977),
55--75).

\emp \bmp \textbf{5.61.}
 (Well-known problem). Does an arbitrary uncountable locally
finite group have only one end? \hf {\sl B.\,Hartley}\vs

Yes, it does (D.\,Holt, {\it Bull. London Math. Soc.}, {\bf 13} (1981),
557--560).

\emp \bmp \textbf{5.62.}
 Is a group locally finite if it contains infinite abelian
subgroups and all of them are complementable?\hf {\sl
S.\,N.\,Chernikov}\vs

Not always (N.\,S.\,Chernikov, {\it Math. Notes}, {\bf 28}
(1980), 788--792).

\emp \bmp \textbf{5.63.}
 Prove that a finite group is not simple if it contains two
non-identity elements whose centralizers have coprime indices.\hf
{\sl S.\,A.\,Chunikhin}\vs

This has been proved mod CFSG (L.\,S.\,Kazarin, in: {\it Studies in
Group Theory}, Sver\-d\-lovsk, 1984, 81--99 (Russian)).

\emp \bmp \textbf{5.64.}
 Suppose that a finite group $G$ is the product of two subgroups
$A_1$ and $A_2$. Prove that if $A_i$ contains a nilpotent
subgroup of index $\leq 2$ for $ i = 1,\,2$, then $G$ is
soluble.\hf {\sl L.\,A.\,Shemetkov}\vs

This has been proved (L.\,S.\,Kazarin, {\it Math. USSR--Sb.}, {\bf
38} (1981), 47--59). \emp

\bmp \textbf{5.65.}
Is the class of all finite groups that have Hall
$\pi $-sub\-groups closed under taking finite subdirect products?
\hf {\sl L.\,A.\,Shemetkov}\vs

Yes, it is (D.\,O.\,Revin, E.\,P.\,Vdovin, {\it J. Group
Theory}, {\bf 14}, no.\,1 (2011), 93--101).
\emp

\bmp \textbf{5.68.}
 Let $G$ be a finitely-presented group, and assume that $G$ has
polynomial growth in the sense of Milnor. Show that $G$ has
soluble word problem.\hf {\sl P.\,E.\,Schupp}\vs

This has been shown (M.\,Gromov, {\it Publ. Math. IHES}, {\bf 53} (1981),
53--73).

\emp \bmp \textbf{5.69.}
 Is every lattice isomorphism between torsion-free groups having
no non-trivial cyclic normal subgroups induced by a group
isomorphism?\hf {\sl B.\,V.\,Yakovlev}\vs

No (A.\,Yu.\,Olshanskii, {\it C.~R.~Acad. Bulgare Sci.}, {\bf 32}
(1979), 1165--1166).

\emp \bmp \textbf{6.4.} A
group $G$ is called {\it of type $(FP)_{\infty}$} if the trivial
$G $-module ${\Bbb Z}$ has a resolution by finitely generated
projective $G $-modules.
 Is it true that every torsion-free group of type
$(FP)_{\infty}$ has finite cohomological dimension? \hf {\sl
R.\,Bieri}\vs

No (K.\,S.\,Brown, R.\,Geoghegan, {\it Invent. Math.}, {\bf 77}
(1984), 367--381). \emp
\markboth{\protect\vphantom{(y}{Archive of solved
problems (6th ed., 1978)}}{\protect\vphantom{(y}{Archive of solved
problems (6th ed., 1978)}}

\bmp \textbf{6.6.} Let $P,\, Q$ be
permutation representations of a finite group $G$ with the same
character. Suppose $P(G)$ is a primitive permutation group. Is
$Q(G)$ necessarily primitive? The answer is known to be
affirmative if $G$ is soluble. \hf {\sl H.\,Wielandt}\vs

No, not necessarily (mod CFSG) (R.\,M.\,Guralnick,
 J.\,Saxl, in: {\it
 Groups, Combinatorics, Geometry, Proc. Durham, 1990}, Cambridge
Univ. Press, 1992, 364--367). \emp

\bmp \textbf{6.7.} Suppose that $P$ is a finite 2-group. Does there
exist a characteristic subgroup $L(P)$ of $P$ such that $L(P)$ is
normal in $H$ for every finite group $H$ that satisfies the
following conditions: 1) $P$ is a Sylow 2-subgroup of $H$, \,
2)~$H$ is ${\Bbb S}_4 $-free, and 3)~$C_H(O_2(H)) \leq O_2(H)$?
\hf {\sl G.\,Glauberman}\vs

Yes, there does (B.\,Stellmacher, {\it Israel J.~Math.}, {\bf 94}
(1996), 367--379). \emp

\bmp \textbf{6.8.}
 Can a residually finite locally normal group be embedded in a
Cartesian product of finite groups in such a way that each
element of the group has at most finitely many central
projections?\hf {\sl Yu.\,M.\,Gorchakov}\vs

Yes, it can (M.\,J.\,Tomkinson, {\it Bull. London Math. Soc.}, {\bf 13}
(1981),
133--137).

\emp \bmp \textbf{6.12.}
 a) Is a metabelian group countable if it satisfies the weak
minimum condition for normal subgroups?

\mb{b)} Is such a group minimax if it is
torsion-free?\hf {\sl D.\,I.\,Zaitsev}\vs

Yes, it is, in both cases (D.\,I.\,Zaitsev, L.\,A.\,Kurdachenko,
A.\,V.\,Tushev, {\it Algebra and Logic}, {\bf 24} (1985), 412--436).

\emp \bmp \textbf{6.13.}
 Is it true that if a non-abelian Sylow 2-subgroup of a finite
group $G$ has a non-trivial abelian direct factor, then $G$ is
not simple?\hf {\sl A.\,S.\,Kondratiev}\vs

Yes, it is, mod CFSG (V.\,V.\,Kabanov, A.\,S.\,Kondratiev, {\it Sylow
$2$-sub\-groups of finite groups (a survey)}, Inst. Math. Mech. UNC
AN SSSR, Sverdlovsk, 1979 (Russian)).

\emp

 \bmp \textbf{6.14.} Are the following lattices locally finite: the lattice of
all locally finite varieties of groups? the lattice of all varieties
of groups? \hfill {\sl A.\,V.\,Kuznetsov}\vs

No, they are not (M.\,I.\,Anokhin, {\it Izv. Math.}, {\bf 63},
no.\,4 (1999),
 649--665).

\emp

\bmp \textbf{6.15.}
A variety is said to be {\it pro-locally-finite\/} if it is not
locally finite while all of its proper subvarieties are locally
finite. An example --- the variety of abelian groups. How many
pro-locally-finite varieties of groups are there? \hfill {\sl
A.\,V.\,Kuznetsov}\vs

There are continuum of
 such varieties (P.\,A.\,Kozhevnikov, {\em On varieties of groups of large odd
exponent}, Dep.\,1612-V00, VINITI, Moscow, 2000 (Russian)). \emp

\bmp \textbf{6.16.}
A variety is called {\it sparse\/} if it has at most countably
many subvarieties. How many sparse varieties of groups are there?
\hfill {\sl A.\,V.\,Kuznetsov}\vs

There are continuum of such varieties (P.\,A.\,Kozhevnikov, {\em
On varieties of groups of large odd exponent}, Dep.\,1612-V00,
VINITI, Moscow, 2000 (Russian); \ S.\,V.\,Ivanov,
A.\,M.\,Storozhev, {\it Contemp. Math.}, {\bf 360} (2004),
55--62). \emp

\bmp \textbf{6.17.}
 Is every variety of groups generated by its finitely generated
groups that have soluble word problem?\hf {\sl
A.\,V.\,Kuznetsov}\vs

No (Yu.\,G.\,Kleiman, {\it Math. USSR--Izv.}, {\bf 22} (1984), 33--65).

\emp

\bmp \textbf{6.18.}
 (Well-known problem). Suppose that a class ${\goth K}$ of
2-generator groups generates the variety of all groups. Is a
non-cyclic free group residually in ${\goth K}$? \hf {\sl
V.\,M.\,Levchuk}
\vs

Not always (S.\,J.\,Pride, \textit{Math.~Z.}, \textbf{131} (1973), 245--248).
\emp

\bmp \textbf{6.19.}
 Let $R$ be a nilpotent associative ring. Are the following two
statements for a subgroup $H$ of the adjoint group of $R$
always equivalent: \ 1) $H$ is a normal subgroup; \ 2) $H$ is an
ideal of the groupoid $R$ with respect to Lie multiplication?\hf
{\sl V.\,M.\,Levchuk}\vs

Not always. Let $R$ be the free nilpotent of index 3 associative
algebra over $\F_2$ on the free generators $x,y$. Let $S$ be the
subalgebra generated by the elements $[x,y^2]=x*y^2$, \ $
[x,xy]=x*(xy)=x(x*y)$, \ $[x,yx]=x*(yx)=(x*y)x$, where $[a,b]$
denotes the commutator in the adjoint group with multiplication
$a\circ b=a+b+ab$, and $a*b=ab-ba$ is Lie multiplication. Let
$M,N$ be the subalgebras generated by $S$ and the elements $x*y$
and $[x,y]$, respectively. The minimal subgroup of $(R,\circ )$
 that contains $x$ and is an ideal of the groupoid $(R,*)$ equals
$\left< x\right>\circ M=\left< x\right> +M$, where $ \left<
x\right> =\{ 0,\,x,\,x^2,\,x+x^2+x^3\}$, while the minimal normal
subgroup containing $x$ equals $\left< x\right>\circ N=\left<
x\right> +N$. Neither is contained in the other. (E.\,I.\,Khukhro,
{\it Letter of July, 23, 1979.})

\emp \bmp \textbf{6.20.}
 Does there exist a supersoluble group of odd order, all of
whose automorphisms are inner?\hf {\sl
V.\,D.\,Mazurov}\vs

Yes, there does (B.\,Hartley, D.\,J.\,S.\,Robinson, {\it Arch.
Math.}, {\bf
35} (1980), 67--74).

\emp \bmp \textbf{6.22.}
 Construct a braid that belongs to the derived subgroup of
the braid group but is not a commutator.\hf {\sl
G.\,S.\,Makanin}\vs

Let $x=\s _1,\,y=\s _2$ be the standard generators of the braid
group $B_3$; then
the braid $(xyxyx)^{12}(xy)^{-12}(xyx)^{-12}$ belongs to the
derived subgroup of $B_3$ but is not a commutator (Yu.\,S.\,Sem\"enov,
{\it Abstracts of
the 10th All-Union Sympos. on Group Theory}, Minsk, 1986, p.~207
(Russian)).

\emp

\bmp \textbf{6.23.}
 A braid $K$ of the braid group ${\frak
B}_{n+1}$ is said to be {\it smooth\/} if removing any of the threads
in $K$ transforms $K$ into a braid that is equal to $1$ in ${\frak
B}_{n}$. It is known that smooth braids form a free subgroup.
Describe generators of this subgroup.

\hf {\sl G.\,S.\,Makanin}\vs

They are described in (D.\,L.\,Johnson, {\it Math. Proc. Camb.
Philos. Soc.}, {\bf 92} (1982), 425--427). \emp

\bmp \textbf{6.25.}
 (Well-known problem). Find an algorithm for calculating the
rank of coefficient-free equations in a free group. The {\it rank\/} of an
equation is the maximal rank of the free subgroup generated by
a solution of this equation.\hf {\sl G.\,S.\,Makanin}\vs

This was found (A.\,A.\,Razborov, {\it Math. USSR--Izv.},
{\bf 25} (1984), 115--162).

\emp \bmp \textbf{6.31.}
 a) Suppose that $G$ is a finitely-generated residually-finite
group, $d(G)$ the minimal number of generators of $G$, and $\delta
(G)$ the minimal number of topological generators of the profinite
completion of $G$. Is $\delta (G)=d(G)$ always true?

\hf {\sl
O.\,V.\,Mel'nikov}\vsh

Not always (G.\,A.\,Noskov, {\it Math. Notes}, {\bf 33} (1983),
249--254). \emp

\bmp \textbf{6.34.} Let ${\goth o}$ be an associative ring with
identity. A system of its ideals ${\goth A}=\{ {\goth A}_{ij}\mid
i,j \in {\Bbb Z}\}$ is called a {\it carpet of ideals\/} if $
{\goth A}_{ik}{\goth A}_{kj} \subseteq {\goth A}_{ij}$ for all
$i,j,k \in {\Bbb Z}$. If ${\goth o}$ is commutative, then the set
$ \Gamma_n(A)=\{ x \in SL_n({\goth o})\mid x_{ij} \equiv \delta
_{ij}\, ({\rm mod}\, {\goth A}_{ij})\} $ is a group, the {\it
$($special\/$)$ congruenz-subgroup modulo the carpet} ${\goth A}$
(the {\it ``carpet subgroup"}). Under quite general conditions, it
was proved in (Yu.\,I.\,Merzlyakov, {\it Algebra i Logika}, {\bf
3}, no.\,4 (1964), 49--59 (Russian); see also M.\,I.\,Kargapolov,
Yu.\,I.\,Merzlyakov, {\it Fundamentals of the Theory of Groups},
3rd Ed., Moscow, Nauka, 1982, p.~145 (Russian))
 that in the groups $GL_n$ and $SL_n$ the mutual commutator
subgroup of the congruenz-subgroups modulo a carpet of ideals shifted
by $k$ and $l$ steps is again the congruenz-subgroup modulo the
same carpet shifted by $k + l$ steps. Prove analogous theorems \
 a) for orthogonal groups; \linebreak b) for unitary groups. \hf {\sl
Yu.\,I.\,Merzlyakov}\vs

These are proved (V.\,M.\,Levchuk, {\it Sov. Math. Dokl.}, {\bf
42}, no.\,1 (1991), 82--86; \ {\it Ukrain. Math.~J.}, {\bf 44},
no.\,6 (1992), 710--718). \emp

\bmp \textbf{6.35.} (R.\,Bieri, R.\,Strebel). Let ${\goth o}$ be an
associative ring with identity distinct from zero. A group $G$
is said to be {\it almost finitely presented
over} ${\goth o}$ if it has a presentation $G=F/R$ where $F$ is
a finitely generated free group and the ${\goth o}G $-module
$R/[R,R] \otimes _{{\Bbb Z}} {\goth o}$ is finitely generated. It
is easy to see that every finitely presented group $G$ is almost
finitely presented over ${\Bbb Z}$ and therefore also over an
arbitrary ring ${\goth o}$. Is the converse true? \hf {\sl
Yu.\,I.\,Merzlyakov}\vs

No, it is not true (M.\,Bestvina, N.\,Brady, {\it Invent. Math.},
{\bf 129}, no.\,3 (1997), 445--470). \emp

\bmp \textbf{6.36.}
 (J.\,W.\,Grossman). The {\it nilpotent-completion diagram\/} of a group
$G$ is as follows: $G/\g _1G\leftarrow G/\g _2G\leftarrow \cdots
$, where $\g _iG$ is the $i$th term of the lower central series
and the arrows are natural homomorphisms. It is easy to see that
every nilpotent-completion diagram $G_1\leftarrow G_2 \leftarrow
\cdots $ is a {\it $\g$-dia\-gram}, that is, every sequence
$1\rightarrow \g _sG_{s+1}\rightarrow G_{s+1}\rightarrow
G_s\rightarrow 1$, \ $s=1,\,2,\ldots$, is exact. Do $
\g$-dia\-grams exist that are not nilpotent-completion
diagrams?\hfill {\sl Yu.\,I.\,Merzlyakov}\vs

Yes, such  $
\g$-dia\-grams exist  (N.\,S.\,Romanovski\u{\i}, {\it Sibirsk. Mat.
Zh.}, {\bf 26}, no.\,4 (1985), 194--195 (Russian)).

\emp \bmp \textbf{6.37.}
 (H.\,Wielandt, O.\,H.\,Kegel). Is a finite group $G$ soluble if
it has soluble subgroups $A$, $B$, $C$ such that
$G=AB=AC=BC$?\hf {\sl V.\,S.\,Monakhov}\vs

Yes, it is (mod CFSG) (L.\,S.\,Kazarin, {\it Ukrain. Math.~J.}, \textbf{43} (1991), 883--886).
\emp

 \bmp \textbf{6.38.} a)
Let $k$ be a (commutative) field. Find all irreducible subgroups
$G$ of $GL_n(k)$ having the property that $G \cap C \ne
\varnothing$ for every conjugacy class $C$ of $GL_n(k)$.
I~con\-jec\-ture that $ G=GL_n(k)$ except in case $n={\rm char}\,
k=2$, the field $k$ is quadratically closed, and $G$ is conjugate
to the group of all matrices of the form
$\left(\!\!\begin{array}{cc}\alpha &0\cr
0&\beta\end{array}\!\!\right)$,
$\left(\!\!\begin{array}{cc}0&\alpha \cr \beta
&0\end{array}\!\!\right)$ where $\alpha \ne 0$ and $\beta \ne 0$.
\hfill {\sl P.\,M.\,Neumann}

\vs

a) The conjecture is refuted 
(S.\,A.\,Zyubin, {\it Algebra and Logic}, \textbf{45}, no.\,5 (2006),
296--305).
 \emp

\bmp \textbf{6.42.} Let $H$ be a strongly 3-embedded subgroup of a finite
group $G$. Suppose that $Z(H/O_{3'}(H))$ contains an element
of order 3. Does $Z(G/O_{3'}(G))$ necessarily contain an element
of order 3? \hf {\sl N.\,D.\,Podufalov}\vs

Yes, it does (mod CFSG) (W.\,Xiao, {\it Sci. China (A)}, {\bf 33}
(1990), 1172--1181). \emp

\bmp \textbf{6.43.}
 Does the set of quasi-identities holding in the class of all finite
groups possess a basis in finitely many variables? \hf
{\sl D.\,M.\,Smirnov}\vs

No (A.\,K.\,Rumyantsev, {\it Algebra and Logic}, {\bf 19}
(1980), 297--311).

\emp \bmp \textbf{6.44.}
 Construct a finitely generated infinite simple group requiring
more than two generators. \hf {\sl J.\,Wiegold}\vs

This has been done (V.\,S.\,Guba, {\it Siberian Math.~J.}, {\bf
27} (1986), 670--684).
\emp

\bmp \textbf{6.46.} If $G$ is $d $-gene\-ra\-tor group having no
non-trivial finite homomorphic images (in particular, if $G$ is an
infinite simple $d $-gene\-ra\-tor group) for some integer $d \geq
2$, must $G \times G$ be a $d $-gene\-ra\-tor group? \hf {\sl
John S.\,Wilson}\vs

No, it need not (V.\,N.\,Obraztsov, {\it Proc. Roy. Soc.
Edinburgh}, {\bf 123} (1993), 839--855).
\emp

\bmp \textbf{6.49.} Is the minimal condition for abelian normal subgroups inherited by subgroups of finite index?
 This is true for the minimal condition for (all) abelian subgroups (J.\,S.\,Wilson, {\it Math. Z.}, {\bf 114} (1970), 19--21).
 \hf {\sl
S.\,A.\,Chechin}\vs

No, not always. Let $G=[(A\times B\times C\times D)\rtimes (\langle g\rangle\times \langle t\rangle )]\rtimes (Y\times \langle x\rangle )$, where $A,B,C,D,Y$ are quasicyclic $p$-groups, $x^2=1$, while $g$ and $t$ are of infinite order. Let $A=\bigcup _{n=1}^{\infty}\langle a_n\rangle$, \ $B=\bigcup _{n=1}^{\infty}\langle b_n\rangle$, \ $C=\bigcup _{n=1}^{\infty}\langle c_n\rangle$, \ $D=\bigcup _{n=1}^{\infty}\langle d_n\rangle$, \ $Y=\bigcup _{n=1}^{\infty}\langle y_n\rangle$ with $a_{n+1}^p=a_n$, $b_{n+1}^p=b_n$, $c_{n+1}^p=c_n$, $d_{n+1}^p=d_n$, $y_{n+1}^p=y_n$. We impose the relations $[ABCD,Y]=[ABC, g]=[ABD,t]=1$; \ $[x,g]=gt^{-1}$; $[x,a_n]=a_nb_n^{-1}$; $[x,c_n]=c_nd_n^{-1}$; $[g,d_n]=b_n$; $[t,c_n]=a_n$; $[y_n,g]=c_n$; $[y_n,t]=d_n$. Then all abelian normal subgroups of $G$ satisfy the minimal condition for subgroups. The subgroup $H=[(A\times B\times C\times D)\rtimes (\langle g\rangle\times \langle t\rangle )]\rtimes Y$ does not satisfy the minimal condition for abelian normal subgroups, since the subgroups $E_n=A\times B\times C\times \langle g^{2^n}\rangle$ are normal in $H$ and form a strictly decreasing chain.
(S.\,A.\,Chechin, {\it Abstracts of 15th All-USSR Algebraic Conf.}, Krasnoyarsk, 1979).
\emp

\bmp \textbf{6.52.} Let $f$ be a local screen of a
formation which contains all finite nilpotent groups and let $A$ be a
group of automorphisms of a finite group $G$. Suppose that $A$ acts $
f $-stably on the socle of $G/\Phi (G)$. Is it true that $A$ acts $f
$-stably on $\Phi (G)$?

\hf {\sl L.\,A.\,Shemetkov}\vsh

No, not always (A.\,N.\,Skiba, {\it Siber. Math.~J.}, {\bf 34}
(1993), 953--958). \emp

\bmp \textbf{6.53.} A group $G$ of the form
$G=F\leftthreetimes H$ is said to be a {\it Frobenius group with
kernel $F$ and complement} $H$ if $H \cap H^g=1$ for any $g \in
G\setminus H$ and $F\setminus \{ 1\} =G\setminus \bigcup\limits_{g\in
G} H^g$. What can be said about the kernel and the complement of a
Frobenius group? In particular, which groups can be kernels?
complements? \hf {\sl V.\,P.\,Shunkov}\vs

Every group can be embedded into the kernel of a Frobenius group,
and every right-orderable group can be a complement in a Frobenius
group (V.\,V.\,Bludov, {\it Siberian Math.~J.}, {\bf 38}, no.\,6
(1997), 1054--1056) \emp

\bmp \textbf{6.54.} Are there infinite finitely
generated Frobenius groups? \hf{\sl V.\,P.\,Shunkov}\vs

Yes, such groups exist (A.\,I.\,Sozutov, {\it Siberian Math.~J.},
{\bf 35}, no.\,4 (1994), 795--801).
 \emp

\bmp \textbf{6.57.}
 A group $G$ is said to be ({\it conjugacy, $p\hs$-con\-ju\-gacy\/})
{\it biprimitively finite\/} if, for any finite subgroup $H$,
any two elements of prime order (any two conjugate elements of
prime order, of prime order $p$) in $N_G(H)/H$ generate a finite
subgroup. Do the elements of finite order in a (conjugacy)
biprimitively finite group $G$ form a subgroup (the periodic
part of $G$)?\hf {\sl V.\,P.\,Shunkov}\vs

Not always (A.\,A.\,Cherep, {\it Algebra and Logic}, {\bf 29}
(1987), 311--313).

\emp \bmp \textbf{6.58.}
 Are

 \mb{a)} the Al\"eshin $p\hs$-groups and

 \mb{b)} the 2-generator Golod
$p\hs$-groups

conjugacy biprimitively finite groups?\hfill {\sl
V.\,P.\,Shunkov}\vs

a) Not always (A.\,V.\,Rozhkov, {\it Math.
USSR--Sb.}, {\bf 57} (1987), 437--448).

b) Not always (A.\,V.\,Timofeenko,
{\it Algebra and Logic}, {\bf 24} (1985), 129--139).

\emp

\bmp \textbf{6.63.} An infinite group $G$ is called a
{\it monster of the first kind\/} if it has elements of order $> 2$
and for any such an element $a$ and for any proper subgroup $H$ of
$G$, there is an element $g$ in $G\setminus H$, such that $\left<
a,\, a^g\right> =G$. Classify the monsters of the first kind all of
whose proper subgroups are finite. \hf {\sl V.\,P.\,Shunkov}\vs

The centre of such a group coincides with the set of elements of
order $\leq 2$ (V.\,P.\,Shunkov, {\it Algebra and Logic}, {\bf 7}
(1968), no.\,1 (1970), 66--69). An infinite group all of whose
proper subgroups are finite is a monster of the first kind if its
centre coincides with the set of elements of order~$\leq 2$
(A.\,I.\,Sozutov, {\it Algebra and Logic}, {\bf 36}, no.\,5
(1997), 336--348). There are continuously many such groups
(A.\,Yu.\,Olshanskii, {\it Geometry of defining relations in
groups}, Kluwer, Dordrecht, 1991). \emp

\bmp \textbf{6.64.}
 A group $G$ is called a {\it monster of the second kind\/} if
it has elements of order $> 2$ and if for any such element $a$
and any proper subgroup $H$ of $G$ there exists an infinite subset
${\goth M}_{a,H}$ consisting of conjugates of $a$ by elements of
$G\setminus H$ such that $ \left< a,c\right> =G$ for all $ c\in
{\goth M}_{a,H}$. Do mixed monsters (that is, with elements of
both finite and infinite orders) of the
second kind exist? Do there exist torsion-free monsters of
the second kind?\hf {\sl V.\,P.\,Shunkov}\vs

Yes, such monsters exist, in both cases (A.\,Yu.\,Olshanskii, {\it The
geometry of defining relations in groups}, Kluwer, Dordrecht,
1991). \emp

\bmp \textbf{7.1.}
The free periodic groups $B(m,p)$ of prime exponent $p > 665$ are
known to possess many properties similar to those of absolutely
free groups (see S.\,I.\,Adian, {\it The Burnside Problem and
Identities in Groups}, Springer, Berlin, 1979). Is it true that
all normal subgroups of $B(m,p)$ are not free periodic groups?
\hfill {\sl S.\,I.\,Adian}\vs

Yes, it is true for all sufficiently large $p$
(A.\,Yu.\,Olshanskii, in: {\it Groups, rings, Lie and Hopf
algebras, Int. Workshop, Canada, 2001}, Dordrecht, Kluwer, 2003,
 179--187); this is also proved for all primes $p\geq 1003$ (V.\,S.\,Atabekyan, {\it Fund. Prikl. Mat.},
 {\bf 15}, no.\,1 (2009), 3--21 (Russian)). \emp
\markboth{\protect\vphantom{(y}{Archive of solved
problems (7th ed., 1980)}}{\protect\vphantom{(y}{Archive of solved
problems (7th ed., 1980)}}

\bmp \textbf{7.2.}
 Prove that the free periodic groups $B(m, n)$ of odd exponent
$n\geq 665$ with $ m \geq 2$ generators are non-amenable and
that random walks on these groups do not have the recurrence
property.\hf {\sl S.\,I.\,Adian}\vs

Both assertions are proved (S.\,I.\,Adian, {\it
Math. USSR--Izv.}, {\bf 21} (1982), 425--434).

\emp

\bmp \textbf{7.4.}
 Is a finitely generated group with quadratic growth almost
abelian?

\hf {\sl V.\,V.\,Belyaev}\vs

Yes, it is (M.\,Gromov, {\it Publ. Math. IHES}, {\bf 53} (1981), 53--73; \
J.\,A.\,Wolf, {\it Diff. Geometry}, {\bf 2} (1968), 421--446).
\emp

\bmp \textbf{7.6.} Describe the infinite simple locally finite groups
with a Chernikov Sylow 2-sub\-group. In particular, are such groups
the Chevalley groups over locally finite fields of odd
characteristic? \hf {\sl V.\,V.\,Belyaev, N.\,F.\,Sesekin}

\vs

No, as
every countably infinite locally finite $p$-group can be embedded as
a maximal $p$-subgroup of the simple Hall's universal group $U$, the unique countable
existentially closed group in the class of all locally finite groups (K.\,Hickin,
{\it Proc. London Math. Soc. (3)}, \textbf{ 52}, no.\,1 (1986), 53--72).
But if all Sylow 2-sub\-groups are Chernikov, then this is true mod CFSG
(O.\,H.\,Kegel, {\it Math. Z.}, \textbf{95} (1967), 169--195; \
V.\,V.\,Belyaev, in: {\it Investigations
in Group Theory,} Sverdlovsk, UNC AN SSSR, 1984, 39--50 (Russian);
\ A.\,V.\,Borovik, {\it Siberian Math.~J.}, \textbf{24}, no.\,6
(1983), 843--851; \ B.\,Hartley, G.\,Shute, {\it Quart. J.~Math.
Oxford (2)}, \textbf{35} (1984), 49--71; \ S.\,Thomas, {\it Arch.
Math.}, \textbf{41} (1983), 103--116).
When all Sylow 2-subgroups are Chernikov, the same result, without using CFSG,
 follows from (M.\,J.\,Larsen, R.\,Pink, {\it J.~Amer. Math. Soc.}, {\bf 24} (2011), 1105--1158 (also known as a preprint of 1998)).\emp

\bmp \textbf{7.7.}
 (Well-known problem). Is the group $G=\left< a,b\mid a^9=1,\;\;
ab= b^2a^2\right>$ finite? This group contains $F(2, 9)$, the only
Fibonacci group for which it is not yet known whether it is finite
or infinite. \hf {\sl R.\,G.\,Burns}\vs

No, it is infinite, since $F(2,9)$ is infinite (M.\,F.\,Newman,
{\it Arch. Math.}, {\bf 54}, no.\,3 (1990), 209--212).

\emp \bmp \textbf{7.8.}
 Suppose that $H$ is a normal subgroup of a group $G$, where $H$
and $G$ are subdirect products of the same $n$ groups
$G_1,\ldots ,G_n$. Does the nilpotency class of $G/H$ increase
with $n$?
\hf {\sl Yu.\,M.\,Gorchakov}\vs

Yes, it does (E.\,I.\,Khukhro, {\it Sibirsk. Mat. Zh.}, {\bf 23},
no.\,6 (1982), 178--180 (Russian)).

\emp \bmp \textbf{7.12.}
 Find all groups with a Hall $2'$-sub\-group.\hf {\sl
R.\,L.\,Griess}\vs

These are found mod CFSG (Z.\,Arad, M.\,B.\,Ward, {\it J.
Algebra}, {\bf 77} (1982), 234--246).

\emp \bmp \textbf{7.16.}
 If $H$ is a proper subgroup of the finite group $G$, is there
always an element of prime-power order not conjugate to an
element of $H$?\hf {\sl R.\,L.\,Griess}\vs

Yes, always (mod CFSG) (B.\,Fein, W.\,M.\,Kantor,
 M.\,Schacher, {\it J. Reine Angew. Math.}, {\bf 328} (1981), 39--57).
\emp

 \bmp \textbf{7.17.} Is
the number of maximal subgroups of the finite group $G$ at most
$|G| - 1$?

\makebox[15pt][r]{}{\it Editors' remarks:\/} This is proved for
$G$ soluble (G.\,E.\,Wall, {\it J.\,\,Austral. Math. Soc.}, {\bf
2} (1961--62), 35--59) and for symmetric groups ${\Bbb S} _n$ for
sufficiently large $n$ (M.\,Liebeck, A.\,Shalev, {\it
J.\,\,Combin. Theory, Ser.\,A}, \textbf{75} (1996), 341--352);
it is also proved (M.\,W.\,Liebeck, L.\,Pyber, A.\,Shalev, {\it
J.~Algebra}, \textbf{317} (2007), 184--197) that any finite group $G$
has at most $2C|G|^{3/2}$ maximal subgroups, where $C$ is an
absolute constant. \hfill {\sl R.\,Griess}
\vs

No, not always (R.\,Guralnick, F.\,L\"ubeck, L.\,Scott, T.\,Sprowl; see (C.\,P.\,Bendel, B.\,D.\,Boe, C.\,M.\,Drupieski, D.\,K.\,Nakano, B.\,J.\,Parshall, C.\,Pillen, C.\,B.\,Wright, in: {\it
Developments and retrospectives in Lie theory. Algebraic methods. Retrospective selected papers based on the presentations at the seminar ``Lie groups, Lie algebras and their representations'', 1991--2014}, Springer, Cham, 2014, 51--69). Infinitely many counterexamples were found in (F.\,L\"ubeck, \emph{Trans. Amer. Math. Soc.}, \textbf{373}, no.\,4 (2020), 2331--2347).

\emp

\bmp \textbf{7.22.} Suppose that a finite group $G$ is realized as the automorphism group of some torsion-free abelian group. Is it true that for every infinite cardinal $\frak m$ there exist $2^{\frak m}$ non-isomorphic torsion-free abelian groups of cardinality $\frak m$ whose automorphism groups are isomorphic to $G$?
\hf {\sl
S.\,F.\,Kozhukhov}\vs

Yes, this is true in the Zermelo--Frenkel system with axioms of choice and `weak diamond' (M.\,Dugas, R.\,G\"obel, \emph{Proc. London Math. Soc. (3)}, {\bf 45}, no.\,2 (1982), 319--336), or if $\frak m$ is smaller than the first measurable cardinal (V.\,A.\,Nikiforov, \emph{Mat. Zametki}, {\bf 39}, no.\,5 (1986), 641--646 (Russian)).
 \emp

\bmp \textbf{7.24.}
 We say that a group is {\it sparse\/} if the variety generated
by it has at most countably many subvarieties. Does there exist
a finitely generated
sparse group that has undecidable word problem?\hf {\sl
A.\,V.\,Kuznetsov}\vs

Yes, such groups do exist; for example, a free group of ${\frak
N}_3{\frak A}$. By (A.\,N.\,Krasil'nikov, {\it Math. USSR--Izv.},
{\bf 37} (1991), 539--553) every subvariety of ${\frak N}_3{\frak
A}$ has a finite basis for its laws; hence there are only
countably many of them. It has undecidable word problem by
(O.\,G.\,Kharlampovich, {\it Sov. Math.}, {\bf 32}, no.\,11
(1988), 136--140).

 \emp

\bmp {\bf 7.30.}
Which finite simple groups can be generated by three involutions,
two of which commute? \hfill {\sl V.\,D.\,Mazurov}\vs

The answer is known mod CFSG. For the alternating groups and
groups of Lie type see (Ya.\,N.\,Nuzhin, {\it Algebra and Logic},
{\bf 36}, no.\,4 (1997), 245--256). For sporadic groups
B.\,L.\,Abasheev, A.\,V.\,Ershov, N.\,S.\,Nevmerzhitskaya,
S.\,Norton, Ya.\,N.\,Nuzhin, A.\,V.\,Timofeenko have shown that
the groups $M_{11}$, $M_{22}$, $M_{23}$, and $McL$ cannot be
generated as required, while the others can; see more details in
(V.\,D.\,Mazurov, {\it Siberian Math.~J.}, {\bf 44}, no.\,1
(2003), 160--164).\emp

\bmp \textbf{7.36.} Is it true that every residually finite group in which every subgroup of finite index (including the group itself) is defined by a single defining relation is either free or isomorphic to the fundamental group of a compact surface?
 \hfill {\sl O.\,V.\,Mel'nikov}\vs

No; for example, let $H_n=\langle x,y\mid y^{-1}xy=x^n\rangle $, $n=2,3,\dots $; then every subgroup of finite index in $H_n$ is isomorphic to a group $H_m$ for some $m$ (V.\,A.\,Churkin, \emph{Abstracts of 8th All-USSR Symp. on Group Theory}, Kiev, 1982, 139--140 (Russian)).
\emp

 \bmp \textbf{7.37.} We
say that a profinite group is {\it strictly complete\/} if each of
its subgroups of finite index is open. It is known (B.\,Hartley,
{\it Math.~Z.}, {\bf 168}, no.\,1 (1979), 71--76) that finitely
generated profinite groups having a finite series with
pronilpotent factors are strictly complete. Is a profinite group
strictly complete if it is \ \

\mb{a)} finitely generated? \

\mb{b)} finitely generated and prosoluble?
\hfill {\sl O.\,V.\,Mel'nikov}\vs

a) Yes, it is (N.\,Nikolov,
D.\,Segal, {\it Ann. Math. (2)}, \textbf{165}, no.\,1 (2007),
171--273).

b) Yes, it is (D.\,Segal, {\it Proc. London Math. Soc. (3)}, {\bf 81}
(2000), 29--54).
\emp

\bmp
 \textbf{7.40.}
Describe the (lattice of) subgroups of a given classical matrix group
over a ring which contain the subgroup consisting of all matrices in
that group with coefficients in some subring (see
Ju.\,I.\,Merzljakov, {\it J.~Soviet Math.}, \textbf{1} (1973),
571--593).

\hfill {\sl Yu.\,I.\,Merzlyakov}
\vs

If the root system has single edges ($A_l$, $D_l$, $E_l$) and the larger ring is not quasi-algebraic over the smaller one, then it is shown that the lattice in question contains the lattice of subgroups of a free product containing one of the factors (A.\,V.\,Stepanov, \textit{J.~Algebra}, \textbf{324}, no.\,7 (2010), 1549--1557), which means that a reasonable description is unfeasible. In almsot all other cases a good description was obtained in (R.\,A.\,Schmidt, {\it Zap. Nauchn. Semin. Leningr. Otd. Mat. Inst. Steklova}, {\bf 94} (1979), 119--130 (Russian); \ Ya.\,N.\,Nuzhin, {\it Algebra Logic}, {\bf 22} (1983), 378--389; \ Ya.\,N.\,Nuzhin, A.\,V.\,Yakushevich, {\it Algebra Logic}, {\bf 39} (2000), 199--206; \
\ A.\,Stepanov, \textit{J.~Algebra}, \textbf{362} (2012), 12--29; \ Ya.\,N.\,Nuzhin, \textit{Siberian Math. J.}, {\bf 54}, no.\,1 (2013), 119--123; \ A.\,Stepanov, A.\,Bak, \textit{St. Petersburg Math. J.}, \textbf{28}, no.\,4 (2017), 47--61).

\emp

\bmp {\bf 7.42.}
A group $U$ is called an {\it $F_q $-group\/}
(where $q \in \pi (U)$) if, for each finite subgroup $K$ of $U$
and for any two elements $a,b$ of order $q$ in $T=N_U(K)/K$, there
exists $c \in T$ such that the group $\left< a,\, b^c\right>$ is
finite. A group $U$ is called an {\it $F ^{\displaystyle{*}}
$-group\/} if each subgroup $H$ of $U$ is an $F_q $-group for
every $q \in \pi (H)$ (V.\,P.\,Shunkov, 1977).

 \makebox[25pt][r]{a)} Is every primary $F
^{\displaystyle{*}} $-group satisfying the minimum condition for
subgroups almost abelian?

 \makebox[25pt][r]{b)} Does every $F
^{\displaystyle{*}} $-group satisfying the minimum condition for
(abelian) subgroups possess the radicable part?\hfill {\sl
A.\,N.\,Ostylovski\u{\i}}\vs

No. A counterexample to both questions is given by an infinite
group all of whose subgroups are conjugate and have prime order
(A.\,Yu.\,Olshanskii, {\it Math. USSR--Izv.}, {\bf 16}
(1981), 279--289). \emp

\bmp \textbf{7.44.} Does a normal subgroup $H$ in a finite group $G$
possess a complement in $G$ if each Sylow subgroups of $H$ is
a direct factor in some Sylow subgroup of $G$?
\hf {\sl V.\,I.\,Sergiyenko}\vs

Yes, it does, mod CFSG (W.\,Xiao, {\it J.~Pure Appl. Algebra},
{\bf 87}, no.\,1 (1993), 97--98). \emp

\bmp \textbf{7.48.} (Well-known problem). Suppose that, in a finite
group $G$, each two elements of the same order are conjugate. Is
then $|G| \leq 6$? \hf {\sl S.\,A.\,Syskin}\vs

Yes, it is, mod CFSG (P.\,Fitzpatrick, {\it Proc. Roy. Irish Acad. Sect. A}, {\bf 85}, no.\,1 (1985),
53--58); see also (W.\,Feit, G.\,Seitz, {\it Illinois J.
Math.}, {\bf 33}, no.\,1 (1988), 103--131) and (R.\,W.\,van~der~Waall, A.\,Bensa\"{\i}d, {\it Simon Stevin}, {\bf
65} (1991), 361--374).

 \emp

\bmp \textbf{7.53.} Let $p$ be a prime.
 The law $ x\cdot x^{\f}\cdots x^{\f^{p-1}}=1$ from the
definition of a splitting automorphism (see Archive, 1.10) gives
rise to a variety of groups with operators $\left< \f\right>$
consisting of all groups that admit a splitting automorphism of
order $p$. Does the analogue of Kostrikin's theorem hold for this
variety, that is, do the locally nilpotent groups in this variety
form a subvariety?\hf {\sl E.\,I.\,Khukhro}\vs

Yes, it does (E.\,I.\,Khukhro, {\it Math. USSR--Sb.}, {\bf 58}
(1987), 119--126).

\emp \bmp \textbf{7.57.}
 A set of generators of a finitely presented group $G$ that
consists of the least possible number $d(G)$ of generators is
called a {\it basis\/} for $G$. Let $r_M(G)$ be the least number
of relations necessary to define $G$ in the basis $M$, and $r(G)$
the minimum of $r_M(G)$ over all bases $M$ for $G$. Let $G_1$,
$G_2$ be any non-trivial groups.

\mb{b)} Is it true that $r_{M_1\cup M_2}(G_1*G_2)=
r_{M_1}G_1+r_{M_2}(G_2)$ for any bases $M_1$, $M_2$ of $G_1$,
$G_2$, respectively?

\mb{c)} Is it true that $r(G_1*G_2)=r(G_1)+r(G_2)$?\hf {\sl
V.\,A.\,Churkin}\vs

b) Not always; c) not always (C.\,Hog, M.\,Lusztig, W.\,Metzler, in:
{\it Presentation classes, 3-ma\-ni\-folds and free products
$($Lecture Notes in Math.}, {\bf 1167}), Springer, Berlin, 1985,
154--167).

\emp

\bmp \textbf{8.5.}
Prove that if $X$ is a finite group, $F$ is any
field, and $M$ is a non-trivial irre\-du\-cible $FX $-module
then $\frac{1}{|X|}\sum_{x\in X} \dim \, {\rm fix}\,(x) \leq
\frac{1}{2}\dim \, M$. \hfill {\sl M.\,R.\,Vaughan-Lee}
\vs

Proved, even with ${...}\leq
\frac{1}{p} \dim \, M$, where $p$ is the least prime divisor of $|G|$ \ (R.\,M.\,Guralnick,
A.\,Mar\'oti, {\it Adv. Math.}, {\bf 226}, no.\,1 (2011), 298--308).
\emp

\bmp \textbf{8.6.}
 (M.\,M.\,Day). Do the classes of amenable groups and elementary
groups coincide? The latter consists of groups that can be
obtained from commutative and finite groups by forming
 subgroups, factor-groups, extensions, and direct limits.
 \hf
{\sl R.\,I.\,Grigorchuk}\vs

No (R.\,I.\,Grigorchuk, {\it Soviet Math. Dokl.}, {\bf 28}
(1983), 23--26).

\emp
\markboth{\protect\vphantom{(y}{Archive of solved
problems (8th ed., 1982)}}{\protect\vphantom{(y}{Archive of solved
problems (8th ed., 1982)}}
\bmp \textbf{8.7.}
Does there exist a non-amenable finitely presented group which has
no free subgroups of rank 2? \hfill {\sl R.\,I.\,Grigorchuk}\vs

Yes, there does (A.\,Yu.\,Olshanskii, M.\,V.\,Sapir, {\it Publ.
Math. Inst. Hautes \'Etud. Sci.}, {\bf 96} (2002), 43--169).\emp

\bmp \textbf{8.8.}
 (D.\,V.\,Anosov). a) Does there exist a non-cyclic
finitely-generated group $G$ containing an element $a$ such that
every element of $G$ is conjugate to a power of $a$?

\hf {\sl
R.\,I.\,Grigorchuk}\vsh

Yes, there does (V.\,S.\,Guba, {\it Math. USSR--Izv.}, {\bf 29}
(1986), 233--277).

\emp

\bmp \textbf{8.10.} Is the group $G=\left< a,b\mid a^n=1,\;
ab=b^3a^3\right>$ finite or infinite for $n=7$, $n=9$, and $n=15$? All
other cases known, see Archive, 7.7.\hf {\sl D.\,L.\,Johnson}\vs

Infinite in all three cases. For $n=15$ by (D.\,J.\,Seal, {\it Proc. Roy. Soc. Edinburgh
(A)}, {\bf 92} (1982), 181--192). For $n=9$ (and $15$) by (M.\,I.\,Prishchepov, {\it Commun. Algebra}, {\bf 23}
(1995), 5095--5117). For $n=7$, because it contains the Fibonacci group $F(3,7)$ as an index~$7$ subgroup, as follows from Theorem 3.0 of (C.\,P.\,Chalk, {\it Commun. Algebra\/} {\bf 26}, no.\,5 (1998), 1511--1546) by standard technique for working with Fibonacci groups (G.\,Williams, {\it Letter of 6 October 2015}).
\emp

\bmp \textbf{8.12.}
b) Let $D0$ denote the class of finite groups of deficiency zero,
i.~e. having a presentation $\left< X\mid R\right>$ with
$|X|=|R|$. Are the central factors of nilpotent $D0\hskip0.1ex
$-groups 3-generated? \hfill {\sl D.\,L.\,Johnson,
E.\,F.\,Robertson}\vs

No, not always (G.\,Havas, E.\,F.\,Robertson, {\it Commun.
Algebra}, {\bf 24} (1996), 3483--3487). \emp

\bmp \textbf{8.13.} Let $G$ be a simple algebraic group over an
algebraically closed field of
characteristic $p$ and ${\goth G}$ be the Lie algebra of $G$.
Is the number of orbits of nilpotent elements of ${\goth G}$ under the
adjoint action of $G$ finite? This is known to be true if $p$
is not too small (i.~e. if $p$ is not a ``bad prime" for $G$).
\hf {\sl R.\,W.\,Carter}\vs

Yes, it is (D.\,F.\,Holt, N.\,Spaltenstein,
{\it J. Austral. Math. Soc. (A)}, {\bf 38} (1985), 330--350).
\emp

\bmp \textbf{8.14.} a) Assume a group $G$ is existentially closed in the
class~$L{\goth N}_p$ of all locally finite $p\hskip0.1ex
$-groups. Is it true that $G$ is characteristically simple? This is true for
 $G$ countable in $L{\goth N}_p$ (Berthold Maier, Freiburg); in
fact, up to isomorphism, there is
only one such countable locally finite $ p\hskip0.1ex $-group.
 \hf {\sl O.\,H.\,Kegel}\vs

Not always (S.\,R.\,Thomas, {\it
Arch. Math.}, {\bf 44} (1985), 98--109).
\emp

\bmp \textbf{8.17.}
 An {\it $\widetilde{RN}$-group\/} is one
whose every homomorphic image is an $RN$-group. Is the class of
$\widetilde{RN}$-groups closed under taking normal subgroups?\hf
{\sl Sh.\,S.\,Kemkhadze}\vs

No (J.\,S.\,Wilson, {\it Arch. Math.}, {\bf 25} (1974), 574--577).

\emp

\bmp \textbf{8.18.} Is every countably infinite
abelian group a verbal subgroup of some finitely generated (soluble)
relatively free group? \hf {\sl Yu.\,G.\,Kleiman}\vs

Yes, it is (A.\,Storozhev, {\it Commun. Algebra}, {\bf 22}, no.\,7
(1994), 2677--2701). \emp

\bmp \textbf{8.20.}
What is the cardinality of the set of all varieties covering an
abelian (nilpotent? Cross? hereditarily finitely based?) variety
of groups? The question is related to 4.46 and 4.73. \hfill {\sl
Yu.\,G.\,Kleiman}\vs

There are continually many varieties covering the variety $A$ of
abelian groups, as well as the variety $A_n$ of abelian groups of
sufficiently large odd exponent $n$ (P.\,A.\,Kozhevnikov, {\em On
varieties of groups of large odd exponent}, Dep.\,1612-V00,
VINITI, Moscow, 2000 (Russian); \ S.\,V.\,Ivanov,
A.\,M.\,Storozhev, {\it Contemp. Math.}, {\bf 360} (2004),
55--62). \emp

\bmp \textbf{8.22.} If $G$ is a (non-abelian) finite
group contained in a join ${\goth A} \vee {\goth B}$ of two varieties
${\goth A}, {\goth B}$ of groups, must there exist finite groups $A
\in {\goth A}$, $B \in {\goth B}$ such that $G$ is a section of the
direct product $A \times B$? \hf {\sl L.\,G.\,Kov\'acs}\vs

Not always (A.\,Storozhev, {\it Bull. Austral. Math. Soc.}, {\bf
51}, no.\,2 (1995), 287--290).

\emp

\bmp \textbf{8.26.} We call a variety {\it
passable\/} if there exists an unrefinable chain of its subvarieties
which is well-ordered by inclusion. For example, every variety
generated by its finite groups is passable --- this is an easy
consequence of (H.\,Neumann, {\it Varieties of groups}, Berlin et
al., Springer, 1967, Chapter~5). Do there exist non-passable
varieties of groups? \hf {\sl A.\,V.\,Kuznetsov}\vs

Yes, there do (M.\,I.\,Anokhin, {\it Moscow Univ. Math. Bull.},
{\bf 51}, no.\,1 (1996), 48--49).

\emp

\bmp \textbf{8.28.} Is the variety of groups finitely
based if it is generated by a finitely based quasivariety of groups?
\hf {\sl A.\,V.\,Kuznetsov}\vs

No, not always (M.\,I.\,Anokhin, {\it Sb. Math.}, {\bf 189},
no.\,7--8 (1998), 1115--1124). \emp

\bmp \textbf{8.31.}
Describe the finite groups in which every proper
subgroup has a complement in some larger subgroup. Among these
groups are, for example, $PSL_2(7)$ and all Sylow subgroups of
symmetric groups. \hfill {\sl V.\,M.\,Levchuk}

\vs
These groups are described independently in (V.\,M.\,Levchuk, A.\,G.\,Likharev, {\it Siberian Math.~J.}, {\bf 47}, no.\,4 (2006), 659--668) and in (V.\,N.\,Tyutyanov, {\it Proc. Gomel' State Univ.}, {\bf 3} (2006), 178--183 (in Russian)).
\emp

\bmp \textbf{8.32.} Suppose $G$ is a finitely
generated group such that, for any set $\pi$ of primes and any
subgroup $H$ of $G$, if $G/\left< H^G\right>$ is a finite $\pi
$-group then $|G:H|$ is a finite $\pi $-number. Is $G$ nilpotent?
This is true for finitely generated soluble groups. \hf {\sl
J.\,C.\,Lennox}\vs

Not always (V.\,N.\,Obraztsov, {\it J.~Austral. Math. Soc. (A)},
{\bf 61}, no.\,2 (1996), 267--288). \emp

\bmp \textbf{8.34.}
Let $G$ be a finite group. Is it true that
indecomposable projective ${\Bbb Z}G $-modu\-les are finitely
generated (and hence locally free)? \hfill {\sl
P.\,A.\,Linnell}
\vs

 Yes, it is; see Remark~2.13 and Corollary~4.2 in
(P.\,P{\v{r}}{\'{\i}}hoda, {\it Rend. Semin. Mat. Univ. Padova}, {\bf 123} (2010), 141--167).
\emp

\bmp \textbf{8.35.} Determine the conjugacy classes
of maximal subgroups in the sporadic simple groups \

\mb{a)} $F'_{24} $; \
\

\mb{b)} $F_2$.
 \hf{\sl V.\,D.\,Mazurov}\vs

a) They are determined mod CFSG (S.\,A.\,Linton, R.\,A.\,Wilson,
{\it Proc. London Math. Soc.}, {\bf 63}, no.\,1 (1991), 113--164).

b) They are determined mod CFSG (R.\,A.\,Wilson, {\it J.~Algebra},
{\bf 211} (1999), 1--14).

\emp

\bmp \textbf{8.37.}
 (R.\,Griess). \ a) Is $M_{11}$ a section in $O'N$?

 \mb{b)} The same question for 
 $M_{24}$ in $F_2$; \ $J_1$ in $F_1$
and in $F_2$; \ $J_2$ in $F_{23}$, $F_{24}'$, and $F_2$.

\hf {\sl
V.\,D.\,Mazurov}\vs

a) Yes (A.\,A.\,Ivanov, S.\,V.\,Shpektorov, {\it Abs. 18 All-Union
Algebra Conf.}, Part~1, Ki\-shi\-nev, 1985, 209 (Russian); \
S.\,Yoshiara, {\it J. Fac. Sci. Univ. Tokyo (IA)}, {\bf 32} (1985),
105--141). \ \ \

b) No (R.\,A.\,Wilson, {\it Bull. London Math.
Soc.}, {\bf 18} (1986), 349--350).

\emp \bmp \textbf{8.39.}
 b) Describe the irreducible subgroups of $SL_6(q)$.
\hf {\sl V.\,D.\,Mazurov}\vs

They are described mod CFSG (A.\,S.\,Kondratiev, {\it Algebra and
Logic}, {\bf 28} (1989), 122--138; \ P.\,Kleidman, {\it
Low-dimensional finite classical groups and their subgroups}, Harlow,
Essex, 1989).

\emp \bmp \textbf{8.46.}
 Describe the automorphisms of the symplectic group $Sp_{2n}$
over an arbitrary commutative ring. {\it Conjecture:} they are all
standard.\hf {\sl Yu.\,I.\,Merzlyakov}\vs

The conjecture was proved (V.\,M.\,Petechuk, {\it Algebra and
Logic}, {\bf 22} (1983), 397--405).

\emp \bmp \textbf{8.47.}
 Do there exist finitely presented soluble groups in which the
maximum condition for normal subgroups fails but all central sections
are finitely generated?

\hf {\sl Yu.\,I.\,Merzlyakov}\vsh

Yes (Yu.\,V.\,Sosnovski\u{\i}, {\it Math. Notes}, {\bf 36} (1984),
577--580).

\emp \bmp \textbf{8.48.}
 If a finite group $G$ can be written as the product of two
soluble subgroups of odd index, then is $G$ soluble?
\hf {\sl V.\,S.\,Monakhov}\vs

Yes, it is (mod CFSG) (L.\,S.\,Kazarin, {\it Ukrain. Math.~J.}, \textbf{43} (1991), 883--886).

\emp \bmp \textbf{8.49.}
 Let $G$ be a $p\hs$-group acting transitively as a permutation
group on a set $\W$, let $F$ be a field of characteristic $p$, and
regard $F\W$ as an $FG$-module. Then do the descending and
ascending Loewy series of $F\W$ coincide?\hf {\sl
P.\,M.\,Neumann}\vs

Yes, they do (J.\,L.\,Alperin, {\it Quart. J. Math. Oxford (2)},
{\bf 39}, no.\,154 (1988), 129--133).

\emp \bmp \textbf{8.53.}
 b) Let $n$ be a sufficiently large odd number. Is it true
that every non-cyclic subgroup of $B(m, n)$ has a subgroup
isomorphic to $B(2, n)$?\hfill {\sl
A.\,Yu.\,Olshanskii}\vs

Yes, it is true for odd $n\geq 1003$
(V.\,S.\,Atabekyan, {\it Izv. Math.}, {\bf 73}, no.\,5 (2009), 861--892).

\emp \bmp \textbf{8.56.}
 Let $X$ be a finite set and $f$ a mapping from the set of
subsets of $X$ to
the positive integers. Under the requirement that in a group
generated by $X$, every subgroup $\left< Y\right>$, \
$Y\subseteq X$, be nilpotent of class $\leq f(Y)$, is it true
that the free group $G_f$ relative to this condition is
torsion-free?\hf {\sl A.\,Yu.\,Olshanskii}\vs

Not always (V.\,V.\,Bludov, V.\,F.\,Kleimenov, E.\,V.\,Khlamov, {\it
Algebra and Logic}, {\bf 29} (1990), 95--96).

\emp

 \bmp \textbf{8.61.}
 Suppose that a locally compact group $G$ contains a subgroup that is topologically isomorphic to the additve group of the field of real numbers with natural topology. Is the space of all closed subgroups of $G$ connected in the Chabauty topology? \hf {\sl I.\,V.\,Protasov}\vs

No, not always: let $H$ be the group of matrices of the form
$\begin{pmatrix}
1&z_1&r\\
0&1&z_2\\
0&0&1\end{pmatrix}$, where $z_1,z_2\in \Z$ and $r\in {\Bbb R}$. Then $H$ is locally compact in the natural topology and contains a central subgroup topologically isomorphic to ${\Bbb R}$, but $L(H)$ is not connected (Yu.\,V.\,Tsybenko, \emph{Abstracts of 17th All-USSR Algebraic Conf., Part 1}, Minsk, 1983, 213 (Russian)).
\emp

\bmp \textbf{8.63.}
 Suppose that the space of all closed subgroups of a locally
compact group $G$ is $\s$-compact in the $E$-topo\-logy. Is
it true that the set of closed non-compact subgroups of $G$ is
at most countable?\hf {\sl I.\,V.\,Protasov}\vs

Yes, it is (A.\,G.\,Piskunov, {\it Ukrain. J. Math.}, {\bf 40}
(1988), 679--683).

\emp

\bmp \textbf{8.66.} Construct examples of residually finite groups which
would separate Shun\-kov's classes of groups with $(a,b)
$-finite\-ness condition, (weakly) conjugacy biprimitively finite
groups and (weakly) biprimitively finite groups (see, in particular,
6.57). Can one derive such examples from Golod's construction? \hf
{\sl A.\,I.\,Sozutov}\vs

Such examples are constructed (L.\,Hammoudi, {\it Nil-alg\`ebres
non-nilpotentes et groupes p\'eriodiques infinis $($Doctor
Thesis}), Strasbourg, 1996; \ A.\,V.\,Rozhkov, {\it Finiteness
conditions in groups of automorphisms of trees $($Doctor of Sci.
Thesis}), Krasnoyarsk, 1997; {\it Algebra and Logic}, {\bf 37},
no.\,3 (1998), 192--203). \emp

\bmp \textbf{8.68.}
Let $G=\left< a,b\mid r=1\right> $ where $r$ is
a cyclically reduced word that is not a proper power of any word in
$a,b$. If $G$ is residually finite, is $G_t=\left<a,b\mid
r^t=1\right> $, $t > 1$, residually finite?

\makebox[15pt][r]{}{\it Editors' comment:\/} An analogous fact was
proved for relative one-relator groups (S.\,J.\,Pride, {\it Proc.
Amer. Math. Soc.}, \textbf{136}, no.\,2 (2008), 377--386). \hfill
{\sl C.\,Y.\,Tang}
\vs

Yes, it is (D.\,T.\,Wise, {\it The structure of groups with a quasiconvex hierarchy},
Annals of Mathematics Studies, {\bf 209}, Princeton University Press, Princeton, NJ, 2021).
\emp

\bmp \textbf{8.69.}
Is every 1-relator group with non-trivial torsion
conjugacy separable?

\hfill {\sl C.\,Y.\,Tang}
\vsh

Yes, it is (A.\,Minasyan, P.\,Zalesskii, \textit{J.~Algebra}, {\bf 382} (2013), 39--45).

\emp

\bmp \textbf{8.70.} Let $A, B$ be polycyclic-by-finite groups. Let
$G=A*_H B$ where $H$ is cyclic. Is $G$ conjugacy separable? \hf {\sl
C.\,Y.\,Tang}\vs

Yes, it is (L.\,Ribes, D.\,Segal, P.\,A.\,Zalesskii, {\it
J.~London Math. Soc. (2)}, {\bf 57}, no.\,3 (1998), 609--628).
\emp

\bmp \textbf{8.71.}
 Is every countable conjugacy-separable group embeddable in a
2-generator conjugacy-separable group?\hf {\sl
C.\,Y.\,Tang}\vs

Yes, it is (V.\,A.\,Roman'kov, {\it Embedding theorems for residually
finite groups}, Preprint 84-515, Comp. Centre, Novosibirsk, 1984
(Russian)).

\emp

\bmp \textbf{8.73.}
 We say that a finite group $G$ {\it separates cyclic
subgroups\/} if, for any cyclic subgroups $A$ and $B$ of $G$, there
is a $g\in G$ such that $ A\cap B^g=1$. Is it true that $G$ separates
cyclic subgroups if and only if $G$ has no non-trivial cyclic
normal subgroups?

\hf {\sl J.\,G.\,Thompson}\vs

Not always. For $p_1=2$, $p_2=5$, $p_3=11$, $p_4=17$ let
$R_i$ be an elementary abelian group of order $p_i^2$ and $\f
_i$ a regular automorphism of order 3 of $R_i$, \ $i=1,\,2,\,3,\,4$.
In the direct product $R_1\langle \f _1\rangle \times R_2\langle
\f _2\rangle \times R_3\langle \f _3\rangle \times R_4\langle \f
_4\rangle$ let $G$ be the subgroup generated by all the $R_i$
and the elements $\f _2\f _3\f _4$ and $\f _1\f _3\f _4^{-1}$.
Then $G$ has no non-trivial cyclic
normal subgroups and every (cyclic) subgroup of order $2\cdot 5\cdot
11\cdot 17$ intersects any of its conjugates non-trivially.
(N.\,D.\,Podufalov, {\it Abstracts of
 the 9th All-Union Group Theory Symp.}, Moscow, 1984, 113--114
(Russian).) \vskip2ex
\emp

 \bmp \textbf{8.76.}
Give a realistic upper bound for the torsion-free rank of a
finitely generated nilpotent group in terms of the ranks of its
abelian subgroups. More precisely, for each integer $n$ let $f(n)$
be the largest integer $h$ such that there is a finitely generated
nilpotent group of torsion-free rank $h$ with the property that
all abelian subgroups have torsion-free rank at most $n$. It is
easy to see that $f(n)$ is bounded above by $n(n + 1)/2$. Describe
the behavior of $f(n)$ for large $n$. Is $f(n)$ bounded below by a
quadratic in $n$? \hfill {\sl John S.\,Wilson}\vs

Yes, it is (M.\,V.\,Milenteva, {\it J.~Group Theory}, {\bf 7},
no.\,3 (2004), 403--408). \emp

\bmp \textbf{8.80.} Let $G$ be a locally finite group containing a
maximal subgroup which is Cher\-ni\-kov. Is $G$ almost soluble?
\hfill{\sl B.\,Hartley}\vs

Yes, it is (B.\,Hartley, {\it Algebra and Logic}, {\bf 37}, no.\,1
(1998), 101--106).

\emp

\bmp \textbf{8.81.}
 Let $G$ be a finite $p\hs$-group admitting an automorphism $\alpha
$ of prime order $q$ with $|C_G(a)|\leq n$.

\mb{a)} If $ p = q$, then does $ G$ have a nilpotent subgroup of
class at most 2 and index bounded by a function of $n$?

\mb{b)} If $p\ne q$, then does $G$ have a nilpotent subgroup of
class bounded by a function of $q$ and index bounded by a
function of $n$ (and, possibly, $q$)? \hf {\sl B.\,Hartley
}\vs

a) Not necessarily (E.\,I.\,Khukhro, {\it Math. Notes}, {\bf 38}
(1985), 867--870).

b) Yes, it does (E.\,I.\,Khukhro, {\it Math. USSR--Sb.}, {\bf 71},
no.\,1 (1992), 51--63). \emp

\bmp \textbf{8.84.} We say that an automorphism $\varphi$ of the
group $G$ is a {\it pseudo-identity\/} if, for all $x \in G$,
there exists a finitely generated subgroup $K_x$ of $G$ such that
$x \in K_x$ and $\varphi |_{K_{\scriptstyle{x}}}$ is an
automorphism of $K_x$. Let $G$ be generated by subgroups $H, K$
and let $G$ be locally nilpotent. Let $\varphi \! : G \rightarrow
G$ be an endomorphism such that $\varphi |_H$ is pseudo-identity
of $H$ and $\varphi |_K$ is an automorphism of $K$. Does it follow
that $\varphi$ is an automorphism of $G$? It is known that
$\varphi$ is a pseudo-identity of $G$ if, additionally, $\varphi
|_K$ is pseudo-identity of $K$; it is also known that $\varphi$ is
an automorphism if, additionally, $K$ is normal in $G$. \hf {\sl
P.\,Hilton}\vs

No, not necessarily (A.\,V.\,Yagzhev, {\it Math. Notes}, {\bf 56},
no.\,5 (1994), 1205--1207). \emp

\bmp \textbf{8.87.} Find all hereditary local formations ${\goth F}$ of
finite groups satisfying the following condition: every finite
minimal non-${\goth F} $-group is biprimary. A finite group is said
to be {\it biprimary\/} if its order is divisible by precisely two
distinct primes. \hfill {\sl L.\,A.\,Shemetkov}\vs

They are found (V.\,N.\,Semenchuk, {\it Problems of Algebra: Proc.
of the Gomel' State Univ.}, no.\,1\,(15) (1999), 92--102
(Russian)). \emp

\bmp \textbf{9.2.} Is $G=\left<a,b\mid a^l=b^m=(ab)^n=1\right> $
conjugacy separable?\hf {\sl R.\,B.\,J.\,T.\,Allenby}\vs

Yes, it is (B.\,Fine, G.\,Rosenberger, {\it Contemp. Math.},
{\bf 109} (1990), 11--18).
\emp
\markboth{\protect\vphantom{(y}{Archive of solved
problems (9th ed., 1984)}}{\protect\vphantom{(y}{Archive of solved
problems (9th ed., 1984)}}

\bmp \textbf{9.3.} Suppose that a countable locally finite group $G$
contains no proper subgroups isomorphic to $G$ itself and suppose
that all Sylow subgroups of $G$ are finite. Does $G$ possess a
non-trivial finite normal subgroup? \hfill {\sl
V.\,V.\,Belyaev}\vs

 Yes, it does
(S.\,D.\,Bell, {\it Locally finite groups with \v{C}ernikov Sylow
subgroups}, ({\it Ph.\,D. Thesis\/}), University of Manchester,
1994). \emp

\bmp \textbf{9.8.}
Does there exist a finitely generated simple group
of intermediate growth?

\hfill {\sl R.\,I.\,Grigorchuk}
\vs

Yes, it does (V.\,Nekrashevych, {\it Ann. Math.}, {\bf 87}, no.\,3 (2018), 667--719).
\emp

\bmp \textbf{9.10.}
Do there exist
finitely generated groups different from ${\Bbb Z}/2{\Bbb Z}$ which
have precisely two conjugacy classes? \hfill {\sl
V.\,S.\,Guba}
\vs

Yes, there do
(D.\,Osin, {\it Ann. Math.}, {\bf 172}, no.\,1 (2010), 1--39).

\emp

\bmp \textbf{9.12.}
 Is there a soluble group with the following properties:
torsion-free of finite rank, not finitely generated, and having
a faithful irreducible representation over a finite field?\hf
{\sl D.\,I.\,Zaitsev}\vs

Yes, there is (A.\,V.\,Tushev, {\it Ukrain. Math.~J.}, {\bf 42}
(1990), 1233--1238).

\emp \bmp \textbf{9.16.}
 (Well-known problem). The {\it prime graph\/} of a finite group
$G$ is the graph with vertex set $\pi (G)$ and an edge joining
$p$ and $q$ if and only if $G$ has an element of order~$pq$.
Describe all finite Chevalley groups over a field of
characteristic $2$ whose prime graph is not connected and
describe the connected components.\hf {\sl
A.\,S.\,Kondratiev}\vs

This was done in (A.\,S.\,Kondratiev, {\it Math. USSR--Sb.}, {\bf 67}
(1990), 235--247).

\emp \bmp \textbf{9.17.}
 a) Let $G$ be a locally normal residually finite group. Can $G$
be embedded in a direct product of finite groups when
the factor group $G/[G, G]$ is a direct product of cyclic
groups?\hf {\sl L.\,A.\,Kurdachenko}\vs

Not always (L.\,A.\,Kurdachenko, {\it Math. Notes}, {\bf 39}
(1986), 273--279).

\emp

\bmp
 \textbf{9.18.} Let ${\goth
S}_{\displaystyle{*}}$ be the smallest
normal
Fitting
class. Are there
Fitting classes which are maximal in ${\goth S}_{\displaystyle{*}}$ (with
respect to
inclusion)? \hf {\sl H.\,Lausch}\vs

No, there are no such classes (N.\,T.\,Vorob'y\"ev, {\it Dokl.
Akad. Nauk Belorus. SSR}, {\bf 35}, no.\,6 (1991), 485--487
(Russian)). \emp

\bmp \textbf{9.19.} b) Let $n(X)$ denote the minimum of the indices of proper
subgroups of a group~$X$. A subgroup $A$ of a finite group $G$
is called {\it wide\/} if $A$ is a maximal element by inclusion of
the set $\{ X\mid X \mbox{ is a proper subgroup of}G\mbox{ and
}n(X)=n(G)\} $. Prove, without using CFSG, that
$n(F_1)=|F_1:2F_2|$, where $F_1$ and $F_2$ are the Fischer simple
groups and $2F_2$ is an extension of a group of order 2 by $F_2$.
\hf {\sl V.\,D.\,Mazurov}\vs

This was proved (S.\,V.\,Zharov, V.\,D.\,Mazurov, in: {\it
Matematicheskoye programmi\-ro\-va\-niye i prilozheniya},
Ekaterinburg, 1995, 96--97 (Russian)). \emp

\bmp \textbf{9.21.} Let $P$ be a maximal parabolic subgroup of the
smallest
index in a finite group $G$ of Lie type $ E_6$, $ E_7$, $ E_8$,
or $^2E_6$ and let $X$ be a subgroup such that $PX=G$. Is it
true that $X=G$? \hf {\sl V.\,D.\,Mazurov}\vs

Yes, it is true (mod CFSG) (C.\,Hering, M.\,W.\,Liebeck, J.\,Saxl,
{\it J.~Algebra}, {\bf 106}, no.\,2 (1987), 517--527). \emp

\bmp \textbf{9.25.}
Find an algorithm which recognizes, by an
equation $w(x_1,\ldots ,x_n)=1$ in a free group $F$ and by a list of
finitely generated subgroups $H_1,\ldots ,H_n$ of $ F$, whether there is a
solution of this equation satisfying the condition $x_1 \in
H_1$, $\ldots$ , $x_n \in H_n$.

\hfill
{\sl G.\,S.\,Makanin}

\vs
Such an algorithm is found in (V.\,Diekert, C.\,Guti{\'e}rrez, C.\,Hagenah,
{\it Inf. Comput.}, {\bf 202}, no.\,{2} (2005), 105--140).
\emp

\bmp \textbf{9.26.} a) Describe the finite groups of 2-local 3-rank 1
which have 3-rank at least 3.

 \hf {\sl A.\,A.\,Makhn\"ev}\vsh

Described in (A.\,A.\,Makhn\"ev, {\it Siber. Math.~J.}, {\bf 29}
(1988), 951--959). \emp

\bmp \textbf{9.27.} Let $M$ be a subgroup of a finite group $G$, \ $A$ an
abelian 2-subgroup of $M$, and suppose that $A^g$ is not
contained in $M$ for some $g$ from $G$. Determine the structure
of $G$ under the hypothesis that $\left< A,A^x\right> =G$
whenever the subgroup $A^x$, $x \in G$, is not contained in $M$.
\hf {\sl A.\,A.\,Makhn\"ev}\vs

It is described mod CFSG (V.\,I.\,Zenkov, {\it Algebra and Logic},
{\bf 35}, no.\,3 (1996), 160--163). \emp

\bmp \textbf{9.30.}
(Well-known problem). A finite set of reductions
$u_i \rightarrow v_i$ of words on a finite alphabet $\Sigma
=\Sigma ^{-1}$ is called a {\it group set of reductions\/} if
$\mbox{{\it length}}\, (u_i) > \mbox{{\it length}}\, (v_i)$ or
$\mbox{{\it length}}\, (u_i)=\mbox{{\it length}}\, (v_i)$ and $u_i
> v_i$ in the lexicographical ordering, and every word in $\Sigma$
can be reduced to the unique reduced form which does not depend on
the sequence of reductions. Do there exist group sets of
reductions satisfying the condition $\mbox{{\it length}}\, (v_i)
\leq 1$ for all $i$, which are different from 1) sets of trivial
reductions $x^{-\varepsilon }x^{\varepsilon} \rightarrow 1$,
$\varepsilon =\pm 1$,\,\, 2)~multiplication tables $xy \rightarrow
z$ of finite groups, and 3)~their finite unions? \hfill {\sl
Yu.\,I.\,Merzlyakov}\vs

Yes, there do (J.\,Avenhaus, K.\,Madlener, F.\,Otto, {\it Trans.
Amer. Math. Soc.}, {\bf 297} (1986), 427--443). \emp

\bmp \textbf{9.32.}
 What locally compact groups satisfy the following
condition: the product of any two closed subgroups is also a closed
subgroup? Abelian groups with this property were described in
(Yu.\,N.\,Mukhin, {\it Math. Notes}, \textbf{8} (1970), 755--760).
\hfill {\sl Yu.\,N.\,Mukhin}

\vs
 Such groups are described (W.\,Herfort, K.\,H.\,Hofmann, F.\,G.\,Russo, \emph{Adv. Math.}, \textbf{390} (2021), 107894).
\emp

\bmp \textbf{9.33.} (F.\,K\"ummich, H.\,Scheerer). If $H$ is a
closed
subgroup of a
connected locally-compact group $G$ such that $ \,\overline{\!
HX}=\,\overline{\! H}\,\overline{\! X}$ for every closed
subgroup $X$ of $G$, then is $H$ normal?\hf {\sl
Yu.\,N.\,Mukhin}\vs

Yes, it is (C.\,Scheiderer, {\it Monatsh. Math.}, {\bf 98} (1984), 75--81).
\emp

\bmp \textbf{9.34.} (S.\,K.\,Grosser, W.\,N.\,Herfort). Does there exist
an infinite compact $p\hskip0.1ex$-group in which the
centralizers of all elements are finite? \hf
{\sl Yu.\,N.\,Mukhin}\vs

No, it does not exist, in view of the connection of this problem
with the Restricted Burnside Problem (S.\,K.\,Grosser,
W.\,N.\,Herfort, {\it Trans. Amer. Math. Soc.}, {\bf 283}, no.\,1
(1984), 211--224) and because of the positive solution to the
latter (E.\,I.\,Zel'manov, {\it Math. USSR--Izv.}, {\bf 36}
(1991), 41--60; {\it Math. USSR--Sb.}, {\bf 72} (1992), 543--565).
\emp

\bmp \textbf{9.41.}
 a) Let $\W$ be a countably infinite set. For $k\geq 2$, we
define a {\it $k$-sec\-tion\/} of $\W$ to be a partition of $\W$
into a union of $k$ infinite subsets. Then does there exist a
transitive permutation group on $\W$ that is transitive on
$k$-sec\-tions but intransitive on ordered $k$-sec\-tions?\hf {\sl
P.\,M.\,Neumann}\vs

Yes, there does. Let $U$ be a non-principal ultrafilter in $
{\goth P}(\W )$ and let $G=\{ g\in {\rm Sym}(\W )\mid {\rm
Fix}(g)\in U\}$. It is not hard to prove that $G$ is transitive on
$k$-sec\-tions but not on ordered $k$-sec\-tions for any $k$ in
the range $ 2\leq k\leq \aleph _0$. (P.\,M.\,Neumann, {\it Letter
of October, 5, 1989\/}.)
\emp

\bmp \textbf{9.43.}
 a) The group $G$ described in the solution of 8.73 (Archive)
enables us to construct a~projective plane of order 3 in which
the lines are the elements of any con\-ju\-gacy class of subgroups of
order $2\cdot 5\cdot 7\cdot 11$ together with four
lines added in a natural way. In a~similar way, we can construct
a projective plane of
order $p^n$ for any prime $p$ and any positive integer $n$. Does
the resulting projective plane have the Galois property?\hf {\sl
N.\,D.\,Podufalov}\vs

Not always. The answer is affirmative for $n=1$. But for $n>1$ the
 planes in question may not have the Galois property. For example, if
there exists a near-field of $p^n$ elements, then among the
planes of order $p^n$
indicated there will definitely be some non-Desarguesian planes.
(N.\,D.\,Podufalov, {\it Letter of November, 24, 1986\/}.)

\emp \bmp \textbf{9.46.}
 Let $G$ be a locally compact group of countable weight and
$L(G)$ the space of all its closed subgroups equipped with the
$E$-topo\-logy. Then is $L(G)$ a $k$-space?

\hf {\sl
I.\,V.\,Protasov}\vs

Yes (I.\,V.\,Protasov, {\it Dokl. Akad. Nauk Ukr.\,SSR Ser.\,A}, {\bf
10} (1986), 64--66 (Russian)).

\emp

\bmp \textbf{9.48.} In a null-dimensional locally compact group, is
the set of all compact elements closed? \hf {\sl
I.\,V.\,Protasov}\vs

Yes, it is (G.\,A.\,Willis, {\it Math. Ann.}, {\bf 300} (1994),
341--363). \emp

\bmp \textbf{9.49.}
Let $G$ be a compact group of weight $> \omega
_2$. Is it true that the space of all closed subgroups of $G$ with
respect to $E$-topo\-logy is non-dyadic?

\hfill {\sl I.\,V.\,Protasov,
Yu.\,V.\,Tsybenko}

\vs  Yes, it is true (Yu.\,V.\,Tsybenko, \emph{Ukrain. Math.~J.}, \textbf{38} (1986), 542--545).
\emp

\bmp \textbf{9.50.} Is every 4-Engel group

\makebox[15pt][r]{}a) without elements of order 2 and 5
necessarily soluble?

\makebox[15pt][r]{}b) satisfying the identity $x^5=1$ necessarily
locally finite?

\makebox[15pt][r]{}c) (R.\,I.\,Grigorchuk) satisfying the identity
$x^8=1$ necessarily locally finite?

 \hf {\sl Yu.\,P.\,Razmyslov}\vs

a) Yes; this follows from the local nilpotency of 4-Engel groups
(G.\,Traustason, {\it Int. J. Algebra Comput.}, {\bf 15} (2005),
309--316; \ G.\,Havas, M.\,R.\,Vaughan-Lee, {\it Int. J. Algebra
Comput.} {\bf 15} (2005), 649--682) and the affirmative answer
for locally nilpotent groups in (A.\,Abdollahi, G.\,Traustason,
{\it Proc. Amer. Math. Soc.}, {\bf 130} (2002), 2827--2836).
\vs
b) Yes, it is (M.\,R.\,Vaughan-Lee, {\it Proc. London Math. Soc.
(3)}, {\bf 74} (1997), 306--334).
\vs
c) Yes; moreover, every 4-Engel 2-group is locally finite
(G.\,Traustason, {\it J.~Algebra}, {\bf 178}, no.\,2 (1995),
414--429).

\emp

\bmp \textbf{9.54.} If $G=HK$ is a soluble group and $ H$, $ K$ are
minimax groups, is it true that $G$ is a minimax group? \hf
{\sl D.\,J.\,S.\,Robinson}\vs

Yes, it is true (J.\,S.\,Wilson, {\it J.~Pure Appl. Algebra}, {\bf
53}, no.\,3 (1988), 297--318; Ya.\,P.\,Sysak, {\it Radical modules
over groups of finite rank}, Inst. of Math. Acad. Sci.
 Ukrain. SSR, Kiev, 1989 (Russian)).

\emp

\bmp \textbf{9.56.}
Find all finite groups with the property that the
tensor square of any ordinary irreducible character is
multiplicity free.
\hfill {\sl J.\,Saxl}
\vs

Such groups are shown to be soluble (L.\,S.\,Kazarin, E.\,I.\,Chankov,
{\it Sb. Math.}, {\bf 201}, no.\,5 (2010), 655--668), and abundance of examples shows that a classification of such soluble groups is not feasible. 
\emp

\bmp \textbf{9.58.}
 Can a product of non-local formations of finite groups be
local?

\hf {\sl A.\,N.\,Skiba,~L.\,A.\,Shemetkov}\vs

Yes, it can (V.\,A.\,Vedernikov, {\it Math. Notes}, {\bf 46} (1989),
910--913; \ N.\,T.\,Vorob'ev, in: {\it Proc. 7th Regional Sci.
Session Math.}, Zielona Goru, 1989, 79--86).

\emp

\bmp \textbf{9.59.}
(W.\,Gasch\"utz). Prove that the formation
generated by a finite group has finite lattice of subformations.
\hfill {\sl A.\,N.\,Skiba, L.\,A.\,Shemetkov}

\vs

Counterexamples are constructed (V.\,P.\,Burichenko, \textit{J.~Algebra}, {\bf 372}
(2012), 428--458).
\emp

\bmp \textbf{9.60.}
Let ${\goth F}$ and ${\goth H}$ be local formations
of finite groups and suppose that ${\goth F}$ is not contained in
${\goth H}$. Does ${\goth F}$ necessarily have at least one minimal
local non-${\goth H}$-sub\-for\-ma\-tion? \hfill {\sl \mbox{A.\,N.\,Skiba, L.\,A.\,Shemetkov}}

\vs

No, not necessarily (V.\,P.\,Burichenko, {\it Trudy Inst. Math. National Akad. Sci. Belarus'}, {\bf 21}, no.\,1 (2013),
15--24 (Russian)).
\emp

\bmp \textbf{9.62.} In any group $G$, the cosets of all of its normal
subgroups together with the empty set form the {\it block
lattice\/} $C(G)$ with respect to inclusion, which, for infinite
$|G|$, is subdirectly irreducible and, for finite $|G| \geq 3$, is
even simple (D.\,M.\,Smirnov, A.\,V.\,Reibol'd, {\it Algebra and
Logic}, {\bf 23}, no.\,6 (1984), 459--470). How large is the class
of such lattices? Is every finite lattice embeddable in the
lattice $C(G)$ for some finite group $G$? \hf {\sl
D.\,M.\,Smirnov}\vs

For each $g \in G$ the filter $\{ x \in C(G)\mid g \in x\} $ is
modular (and even arguesian); the variety generated by the block
 lattices of groups does not contain all finite
lattices (D.\,M.\,Smirnov, {\it Siberian Math.~J.}, {\bf 33},
no.\,4 (1992), 663--668). \emp

\bmp
 \textbf{9.63.} Is a finite group of
the form $G=ABA$
soluble
if $A$ is an
abelian subgroup and $B$ is a cyclic subgroup? \hf {\sl Ya.\,P.\,Sysak}\vs

Yes, it is soluble (mod CFSG) (D.\,L.\,Zagorin, L.\,S.\,Kazarin,
{\it Dokl. Math.}, {\bf 53}, no.\,2 (1996), 237--239).

\emp

\bmp \textbf{9.64.}
Is it true that, in a group of the form $ G=AB$,
every subgroup $N$ of $A \cap B$ which is subnormal both in $A$
and in $B$ is subnormal in $G$? The answer is affirmative in the
case of finite groups (H.\,Wielandt). \hfill {\sl
Ya.\,P.\,Sysak}

\vs

No, not always (C.\,Casolo, U.\,Dardano, {\it J. Group Theory}, {\bf 7} (2004), 507--520).
\emp

\bmp \textbf{9.67.} (A.\,Tarski). Let
$F_n$ be a free group of rank $ n$; is it true that $
Th(F_2)=Th(F_3)$?

\hfill {\sl A.\,D.\,Taimanov}

\vs

Yes, it is
(O.\,Kharlampovich, A.\,Myasnikov, {\it J.\,Algebra}, \textbf{302},
(2006), 451--552; \ Z.\,Sela, {\it Geom. Funct. Anal.}, \textbf{16} (2006), 
707--730).

\emp

\bmp \textbf{9.73.}
 Let ${\goth F}$ be any local formation of finite soluble groups
containing all finite nilpotent groups. Prove that $H\sp{{\frak
F}}K=KH\sp{{\frak F}}$ for any two
subnormal subgroups $H$ and $K$ of an arbitrary finite group
$G$.\hf {\sl L.\,A.\,Shemetkov}\vs

This has been proved, even without the hypothesis of solubility of
${\frak F}$ \ (S.\,F.\,Ka\-mor\-ni\-kov, {\it Dokl. Akad. Nauk
BSSR}, {\bf 33}, no.\,5 (1989), 396--399 (Russian)).

\emp

\bmp \textbf{9.74.} Find all local formations ${\goth F}$ of finite
groups such that every finite minimal non-${\goth F}$-group is
either a Shmidt group (that is, a non-nilpotent finite group all
of whose proper subgroups are nilpotent) or a group of prime
order. \hf {\sl L.\,A.\,Shemetkov}\vs

They are found (S.\,F.\,Kamornikov, {\it Siberian Math.~J.}, {\bf
35}, no.\,4 (1994), 713--721). \emp

\bmp \textbf{9.79.}
 (A.\,G.\,Kurosh). Is every group with the minimum
condition countable?

\hf {\sl V.\,P.\,Shunkov}\vsh

No (V.\,N.\,Obraztsov, {\it Math. USSR--Sb.}, {\bf
66} (1989), 541--553).

\emp \bmp \textbf{9.80.}
 Are the 2-elements of a group with the minimum condition
contained in its locally finite radical?\hf {\sl
V.\,P.\,Shunkov}\vs

Not necessarily (A.\,Yu.\,Olshanskii, {\it The geometry of
defining relations in groups}, Kluwer, Dordrecht, 1991).

\emp \bmp \textbf{9.81.}
 Does there exist a simple group with the minimum condition
possessing a non-trivial quasi-cyclic subgroup? \hfill {\sl
V.\,P.\,Shunkov}\vs

Yes, there does (V.\,N.\,Obraztsov, {\it Math. USSR--Sb.}, {\bf
66} (1989), 541--553).

\emp

\bmp \textbf{9.82.}
 An infinite group, all of whose proper subgroups are finite, is
called {\it quasi-finite}. Is it true that an element of a
quasi-finite group $G$ is central if and only if it is contained
in infinitely many subgroups of $G$?\hf {\sl
V.\,P.\,Shunkov}\vs

No (K.\,I.\,Lossov, Dep. no.\,5529-V89, VINITI, Moscow, 1988
(Russian)).

\emp

\bmp
 {\bf 10.1.} Let $p$ be a
prime number. Describe
the
groups of order
$
p^9$ of nilpotency class~2 which contain subgroups $X$ and $Y$
such
that $|X|=|Y|=p^3$ and any non-identity elements $x \in X$, $y
\in
Y$ do not commute. An answer to this question would yield a
description of the semifields of order $p^3$.
\hf {\sl S.\,N.\,Adamov,~A.\,N.\,Fomin}\vs

They are described (V.\,A.\,Antonov,
{\it Preprint}, Chelyabinsk, 1999 (Russian)).
\emp
\markboth{\protect\vphantom{(y}{Archive of solved
problems (10th ed., 1986)}}{\protect\vphantom{(y}{Archive of solved
problems (10th ed., 1986)}}

\bmp \textbf{10.6.} Is it true that in
an abelian group every non-discrete group topology can be
strengthened up to a non-discrete group topology such that the
group becomes a complete topological group? \hf
{\sl V.\,I.\,Arnautov}\vs

This is true for group topologies satisfying the first axiom of
countability (E.\,I.\,Marin, {\it in: Modules, Algebras, Topologies
$($Mat. Is\-sle\-do\-va\-niya}, {\bf 105}), Ki\-shi\-n\"ev, 1988,
105--119 (Russian)), but this may not be true in general, see
Archive,~12.2. \emp

\bmp \textbf{10.7.}
 Is it true that in a countable group $G$, any non-discrete
group topology satisfying the first axiom of countability can be
strengthened up to a non-discrete group topology such that $G$
becomes a complete topological group?\hf {\sl
V.\,I.\,Arnautov}\vs

Yes, it is (V.\,I.\,Arnautov, E.\,I.\,Kabanova, {\it Siberian
Math.~J.}, {\bf 31} (1990), 1--10).

\emp

\bmp \textbf{10.9.} Let $p$ be a prime
number and let $L_p$ denote the set of all quasivarieties each of
which is generated by a finite
$
p\hskip0.1ex $-group. Is $ L_p$ a sublattice of the lattice of all
quasivarieties of groups? \hf {\sl A.\,I.\,Budkin}\vs

No, not always (S.\,A.\,Shakhova, {\it Math. Notes}, {\bf 53},
no.\,3 (1993), 345--347).

\emp

\bmp \textbf{10.14.} a) Does every group satisfying the
minimum condition on subgroups satisfy the weak maximum condition
on subgroups?

 \mb{b)}
Does every group satisfying the maximum condition on subgroups
satisfy the weak minimum condition on subgroups? \hf {\sl
D.\,I.\,Zaitsev}\vs

a) No, not always (V.\,N.\,Obraztsov, {\it Math. USSR--Sb.}, {\bf
66}, no.\,2 (1990), 541--553).

b) No, not always (V.\,N.\,Obraztsov, {\it Siberian Math.~J.},
{\bf 32}, no.\,1 (1991), 79--84).

\emp

\bmp \textbf{10.23.} Is it true that extraction of roots in braid
groups is unique up to conjugation?

\hfill {\sl G.\,S.\,Makanin}
\vs

Yes, it is true (J.\,Gonz\'alez-Meneses, {\it Algebr. Geom. Topology}, {\bf 3} (2003), 1103--1118).
 \emp

\bmp \textbf{10.24.} A braid is said
to be {\it coloured\/} if its strings represent the identity
permutation. Is it true that links obtained by closing coloured
braids are equivalent if and only if the original braids are
conjugate in the braid group? \hf
{\sl G.\,S.\,Makanin}\vs

No. For non-oriented links: the braid words $\sigma _1^2$ and
$\sigma _1^{-2}$ both represent the Hopf link but are not
conjugate in the braid group (easy to see using the Burau
representation). For oriented links see Example~4 in (J.\,Birman,
{\it Braids, links, and mapping class groups $($Ann. Math. Stud.
Princeton}, {\bf 82}$)$, Princeton, NJ, 1975, p.\,100).
(V.\,Shpil\-rain, {\it Letters of May, 28, 1998 and January, 20,
2002\/}). \emp

\bmp \textbf{10.25.}
(Well-known problem). Does there exist an
algorithm which decides for a given automorphism of a free group
 whether this automorphism has a non-trivial fixed point? \hfill
{\sl G.\,S.\,Makanin}
\vs

Yes, it does (O.\,Bogopolski, O.\,Maslakova,
{\it Int. J. Algebra Comput.}, {\bf 26}, no.\,1 (2016), 29--67).
\emp

\bmp \textbf{10.26.} a)
Does there exist an algorithm which decides,
for given elements $a, b$ and an automorphism $ \varphi$ of a free
group, whether the equation $ax^{\varphi}=xb$ is soluble in this
group? This question seems to be useful for solving the problem of
equivalence of two knots.
 \hfill {\sl G.\,S.\,Makanin}

\vs

a) Yes, it does
(O.\,Bogopolski, A.\,Martino, O.\,Maslakova, E.\,Ventura, {\it Bull. London Math. Soc.}, {\bf 38}, no.\,5 (2006), 787--794); further generalizations and applications are in (O.\,Bogopolski, A.\,Martino, E.\,Ventura, {\it Trans. Amer. Math. Soc.}, {\bf 362} (2010), 2003--2036).
\emp

\bmp \textbf{10.28.}
Is it true that finite strongly regular graphs
with $\lambda =1$ have rank~3?

\hfill {\sl A.\,A.\,Makhn\"ev}

\vs

No, not always
(A.\,Cossidente, T.\,Penttila, {\it J.~London Math. Soc.}, {\bf 72} (2005), 731--741).
\emp

\bmp \textbf{10.30.} Does there exist a non-right-orderable group
which is residually finite $p\hskip0.1ex $-group for a finite set
of prime numbers $p$ containing at least two different primes? If
a group is residually finite $p\hskip0.1ex $-group for an infinite
set of primes $p$, then it admits a linear order
(A.\,H.\,Rhemtulla, {\it Proc. Amer. Math. Soc.}, {\bf 41}, no.\,1
(1973), 31--33).

\hfill {\sl N.\,Ya.\,Medvedev}\vsh

Yes, there does (P.\,A.\,Linnell, {\it J.~Algebra}, {\bf 248} (2002)
605--607). \emp

\bmp \textbf{10.33.} For $HNN
$-exten\-sions of the form $ G=\left< t, A\mid t^{-1}Bt=C, \; \,
\varphi \right> $, where $A$ is a finitely generated abelian group
and $\varphi \! : B \rightarrow C$ is an isomorphism of two of its
subgroups, find

\mb{a)} a criterion to be residually finite;

\mb{b)} a criterion to be Hopfian.\hf {\sl
Yu.\,I.\,Merzlyakov}\vs

a) Such a criterion is found (S.\,Andreadakis, E.\,Raptis,
D.\,Varsos, {\it Arch. Math.}, {\bf 50} (1988), 495--501.

b) Such a criterion is found (S.\,Andreadakis, E.\,Raptis,
D.\,Varsos, {\it Commun. Algebra}, {\bf 20} (1992), 1511--1533). \emp

\bmp \textbf{10.37.}
Suppose that $G$ is a finitely generated metabelian group all of
whose integral homology groups are finitely generated. Is it true
that $G$ is a group of finite rank? The answer is affirmative if
$G$ splits over the derived subgroup (J.\,R.\,J.\,Groves, {\it
Quart.~J.~Math.}, {\bf 33}, no.\,132 (1982), 405--420). \hfill
{\sl G.\,A.\,Noskov} \vs

Yes, it is (D.\,H.\,Kochloukova, {\it Groups St. Andrews 2001 in
Oxford}, Vol. II, Cambridge Univ. Press, 2003, 332--343). \emp

\bmp \textbf{10.41.}
 (Well-known problem). Let $\G$ be an almost polycyclic group
with no non-trivial finite normal subgroups, and let $k$ be a
field. The complete ring of quotients $Q(k\G )$ is a matrix ring
$M_n(D)$ over a skew field. Conjecture: $n$ is the least common
multiple of the orders of the finite subgroups of $\G$. An
equivalent formulation (M.\,Lorenz) is as follows: $\rho (G_0(k\G
))=\rho (G_0(k\G )_{{\cal F}})$, where $G_0(k\G )$ is the
Grothendieck group of the category of finitely-generated
$k\G$-modules, $G_0(k\G )_{{\cal F}}$ is the subgroup generated by
classes of modules induced from finite subgroups of $\G$, and
$\rho$ is the Goldie rank. There is a stronger conjecture:
$G_0(k\G )=G_0(k\G )_{{\cal F}}$.\hf {\sl G.\,A.\,Noskov}\vs

The strong conjecture $G_0(k\G )=G_0(k\G )_{{\cal F}}$ has been
proved (J.\,A.\,Moody, {\it Bull. Amer. Math. Soc.}, {\bf 17}
(1987), 113--116).

\emp \bmp \textbf{10.48.}
Let $V$ be a
vector space of finite dimension over a field of prime order.
A~sub\-set $R$ of $GL(V) \cup \{ 0\}$ is called {\it regular\/} if
$|R|=|V|$, $0,1 \in R$ and $vx \ne vy$ for any non-trivial vector
$v \in V$ and any distinct elements $x,y
\in
R$. It is obvious that $\tau ,\varepsilon ,\mu _g$ transform a
regular set into a regular one, where $x^{\tau}=x^{-1}$ for $x \ne 0$
and $0^{\tau}=0$, $x^{\varepsilon}=1 - x$, \ $x^{\mu _g}=xg^{-1}$ and
$ g$ is a non-zero element of the set being transformed. We say that
two regular subsets are {\it equivalent\/} if one can be obtained
from the other by a sequence of such transformations. \ \ \

\mb{a)} Study the equivalence classes of regular subsets. \ \ \

\mb{b)} Is every regular subset equivalent to a subgroup of
$GL(V)$ together with $0$?

\hf {\sl N.\,D.\,Podufalov}\vs

a) They were studied (N.\,D.\,Podufalov, {\it Algebra and Logic},
{\bf 30}, no.\,1 (1991), 62--69).

b) No. A regular set is closed with respect to multiplication if
and only if the corresponding $ (\g ,\g )$-tran\-si\-tive plane is
defined over a near-field. (N.\,D.\,Podufalov, {\it Letter of
February, 13, 1989.})

\emp \bmp \textbf{10.56.}
 Is the lattice of formations of finite nilpotent groups of
class $\leq 4$ distributive?

\hf {\sl A.\,N.\,Skiba}\vs

No. The lattice $L$ of all varieties of nilpotent 2-groups of
class $\leq 4$ is non-distributive (Yu.\,A.\,Belov, {\it Algebra
and Logic\/}, {\bf 9}, no.\,6 (1970), 371-374; \ R.\,A.\,Bryce,
{\it Philos. Trans. Roy. Soc. London (A)\/}, {\bf 266} (1970),
281--355, footnote on p.\,335). The mapping $ V\rightarrow V\cap
F$, where $F$ is the class of all finite groups, is an embedding
of $L$ into the lattice of formations. (L.\,Kov\'acs, {\it Letter
of May, 12, 1988.}) See also (L.\,A.\,Shemetkov, A.\,N.\,Skiba,
{\it Formations of algebraic systems}, Nauka, Moscow, 1989
(Russian)).

\emp \bmp \textbf{10.63.}
 Is there a doubly transitive permutation group in which the
stabilizer of a~point is infinite cyclic?\hf {\sl
Ya.\,P.\,Sysak}\vs

No (V.\,D.\,Mazurov, {\it Siberian Math.~J.}, {\bf 31} (1990),
615--617; \ Gy.\,K\'arolyi, S.\,J.\,Kov\'acs, P.\,P.\,P\'alfy,
{\it Aequationes Math.}, {\bf 39}, no.\,1 (1990), 161--166).

\emp

\bmp \textbf{10.66.} Is a group $G$ non-simple if it contains two
non-trivial subgroups $A$ and $B$ such that $AB \ne G$ and
$AB^g=B^gA$ for any $g \in G$? This is true if $G$ is finite
(O.\,H.\,Kegel, {\it Arch. Math.}, {\bf 12}, no.\,2 (1961),
90--93). \hf {\sl A.\,N.\,Fomin,~V.\,P.\,Shunkov}\vs

No; for example, for $G$ being any simple
linearly ordered group and $A$ and $B$ arbitrary proper
convex subgroups of $G$ \ (V.\,V.\,Bludov, {\it Letter of February,
12, 1997\/}). \emp

\bmp \textbf{10.68.} Suppose that a finite $ p\hskip0.1ex $-group $G$
admits an automorphism of order $p$ with exactly $p^m$ fixed points.
Is it true that $G$ has a subgroup whose index is bounded by a
function of $p$ and $m$ and whose nilpotency class is bounded by a
function of $m$ only? The results on $p\hskip0.1ex $-groups of
maximal class give an affirmative answer in the case $m=1$. \hfill
{\sl E.\,I.\,Khukhro}\vs

Yes, it is (Yu.\,Medvedev, {\it J. London Math. Soc. (2)}, {\bf
59}, no.\,3 (1999), 787--798). \emp

\bmp \textbf{10.69.}
Suppose that $S$ is a closed oriented surface of genus $g
> 1$ and $G=\pi _1S=G_1*_A G_2$, where $G_1 \ne A \ne G_2$ and
the subgroup $A$ is finitely generated (and hence free). Is this
decomposition geometrical, that is, do there exist connected
surfaces $S_1, S_2, T$ such that $S=S_1 \cup S_2$, $S_1 \cap
S_2=T$, $G_i=\pi _1S_i$, $A=\pi _1T$ and embeddings $S_i \subset
S$, $T \subset S$ induce, in a natural way, embeddings $G_i
\subset G$, $A \subset G$? This is true for $A={\Bbb Z}$. \hfill
{\sl H.\,Zieschang}\vs

No, not always (O.\,V.\,Bogopolski, {\it Geom. Dedicata},
{\bf 94} (2002), 63--89). \emp

 \bmp \textbf{10.72.}
 Prove that the formation of all finite $p\hs$-groups does not
decompose into a product of two non-trivial subformations.\hf
{\sl L.\,A.\,Shemetkov}\vs

This has been proved (A.\,N.\,Skiba, L.\,A.\,Shemetkov, {\it Dokl. Akad.
Nauk BSSR}, {\bf 33} (1989), 581--582 (Russian)).

\emp

\bmp \textbf{10.76.} Suppose that $G$ is a periodic group having an
infinite Sylow 2-sub\-group $S$ which is either elementary abelian or
 a Suzuki 2-group, and suppose that the normalizer $
N_G(S)$ is strongly embedded in $G$ and is a Frobenius group (see
6.55) with locally cyclic complement. Must $G$ be locally finite?
\hfill {\sl V.\,P.\,Shunkov}\vs

Yes, it must (A.\,I.\,Sozutov, {\it Algebra and Logic}, {\bf 39},
no.\,5 (2000),
 345--353; \
 A.\,I.\,So\-zu\-tov, N.\,M.\,Suchkov, {\it Math. Notes}, {\bf
68}, no.\,2 (2000), 237--247). \emp

 \bmp \textbf{11.2.}
 Classify the simple groups
that are isomorphic to the multiplicative groups of finite rings,
in particular, of the group rings of finite groups over finite
fields and over ${\Bbb Z}/n{\Bbb Z}$, $n \in {\Bbb Z}$. \hfill
{\sl R.\,Zh.\,Aleev}
\vs

Such a classification is obtained in (C.\,Davis, T.\,Occhipinti, \emph{J.~Pure Appl. Algebra}, \textbf{218}, no.\,4 (2014), 743--744).
\emp

\bmp \textbf{11.4.} Is it true that the lattice of
centralizers in a group is modular if it is a sub\-lattice of the
lattice of all subgroups? This is true for finite groups. \hf {\sl
V.\,A.\,Antonov}\vs

 No, not always
(V.\,N.\,Obraztsov, {\it J.~Algebra}, {\bf 199}, no.\,1 (1998),
337--343).

\emp
\markboth{\protect\vphantom{(y}{Archive of solved
problems (11th ed., 1990)}}{\protect\vphantom{(y}{Archive of solved
problems (11th ed., 1990)}}

\bmp \textbf{11.6.} Let $p$ be an odd prime.
Is it true that every finite $p\hskip0.2ex$-group possesses a set
of generators of equal orders? \hfill \raisebox{-0.6ex}{\sl
C.\,Bagi\'nski}

\vs

No, it is not true
(E.\,A.\,O'Brien, C.\,M.\,Scoppola, M.\,R.\,Vaughan-Lee, {\it
Proc. Amer. Math. Soc.}, \textbf{134}, no.\,12 (2006), 3457--3464).
\emp

\bmp \textbf{11.8.} a)
For a finite group $X$,
let $\chi _1(X)$ denote the totality of the degrees of all
irre\-du\-cible complex characters of $X$ with allowance for their
multiplicities. Suppose that $\chi _1(G)=\chi _1(H)$ for groups
$G$ and $H$. Clearly, then $|G|=|H|$. Is it true that
$H$ is simple if $G$ is simple?
\hfill {\sl Ya.\,G.\,Berkovich}

\vs
a) Yes, it is (H.\,P.\,Tong-Viet, {\it J.~Algebra}, {\bf 357} (2012), 61--68;
 {\it Monatsh. Math.}, {\bf 166}, no.\,3--4 (2012), 559--577;
 {\it Algebr. Represent. Theory}, {\bf 15}, no.\,2 (2012), 379--389).

\emp

\bmp \textbf{11.10.} (R.\,C.\,Lyndon). b) Is it true
that $a \ne 1$ in $G=\langle a,x\mid a^5=1,\; \,
a^{x^2}=[a,a^x]\rangle $?

\hf {\sl V.\,V.\,Bludov}\vs

Yes, it is true, since the mapping $a\rightarrow (1\, 3\, 5\, 2\,
4)$, $x\rightarrow (1\, 2\, 4\, 3\, 5)$ can be extended to a
homomorphism of the group $G$ onto the alternating group ${\Bbb
A}_5$ (D.\,N.\,Azarov, in: {\it
 Algebraicheskije Systemy}, Ivanovo, 1991, 4--5 (Russian)).

\emp

\bmp \textbf{11.11.} The well-known
Baer--Suzuki theorem
 states that if every two conjugates of an element $a$ of
a finite group $G$ generate a finite $p\hskip0.2ex$-sub\-group,
then $a$ is contained in a normal $p\hskip0.2ex$-sub\-group.

\makebox[25pt][r]{a)} Does
such a theorem hold in the class of periodic groups?
The case $p=2$ is of particular interest.

\makebox[25pt][r]{b)} Does such a theorem hold in the class of binary finite groups?

 \hfill {\sl A.\,V.\,Borovik}\vs

a) No, it does not hold for $p=2$ (V.\,D.\,Mazurov, A.\,Yu.\,Olshanskii, A.\,I.\,Sozutov, \textit{Algebra and Logic},
{\bf 54}, no.\,2 (2015), 161--166). 
\vs

b) Yes, such a theorem does hold for binary finite groups. By the
Baer--Suzuki theorem $\left< a_1,\,b\right>$ is a finite
$p\hskip0.1ex $-group for any element $a_1\in a^{G}=\{ a^g\mid
g\in G\}$ and for any $p\hskip0.1ex $-element $b\in G$. Now
induction on $n$ yields that the product $a_1\cdots a_n$ is a
$p\hskip0.1ex $-element for any $a_1,\ldots ,a_n\in a^G$.
(A.\,I.\,Sozutov, {\it Siberian Math. J.}, {\bf 41}, no.\,3
(2000), 561--562.)

\emp

 \bmp \textbf{11.13.}
 Suppose that $G$ is a
periodic group with an involution $ i$ such that $i^g \cdot i$
has odd order for any $g \in G$. Is it true that the image of
$i$ in $G/O(G)$ belongs to the centre of $G/O(G)$? \hfill {\sl A.\,V.\,Borovik}
\vs

No, not always (E.\,B.\,Durakov, A.\,I.\,Sozutov, {\it Algebra Logic}, {\bf 52}, no.\,5 (2013), 422--425).
\emp

 \bmp \textbf{11.20.}
Suppose we have $[a,b]=[c,d]$ in an absolutely free group, where
$a, b, [a,b]$ are basic commutators (in some fixed free
generators). If $c$ and $d$ are arbitrary (proper) commutators,
does it follow that $a=c$ and $b=d$? \hfill
{\sl A.\,Gaglione, D.\,Spellman}\vs

No, not always.
For example, if $x_1,\, x_2,\, x_3$ are free generators, \ $
x_1<x_2< x_3$, \linebreak $a=[[x_2,x_1],x_3]$, \ $b=[x_2,x_1]$,
\ $c=b^{-1}$, \ $d=b^{x_3b}$ (V.\,G.\,Bardakov, {\it Abstracts of the IIIrd
Intern. Conf. on Algebra}, Krasnoyarsk, 1993, p.~33 (Russian)).

\emp

\bmp \textbf{11.21.} Let ${\goth N}_p$ denote the
formation of all finite $p\hskip0.1ex$-groups, for a given prime
number~$p$. Is it true that, for every subformation ${\goth F}$ of
${\goth N}_p$, there exists a variety ${\goth M}$ such that ${\goth
F}={\goth N}_p \cap {\goth M}$? \hf {\sl A.\,F.\,Vasil'yev}\vs

 No, it is not true.
There is a natural one-to-one correspondence between formations of
finite $p\hs$-groups and varieties of pro-$p\hs$-groups: for every
variety ${\goth V}$ of pro-$p\hs$-groups the class of all finite
groups in ${\goth V}$ is a formation of $p\hs$-groups and every
formation of $p\hs$-groups arises in this way. There are
continuously many varieties of nilpotent pro-$p\hs$-groups of
class at most 6 (A.\,N.\,Zubkov, {\it Siberian Math.~J.}, {\bf
29}, no.\,3 (1988), 491--494) and only countably many varieties of
nilpotent groups of class at most 6. (A.\,N.\,Krasil'nikov, {\it
Letter of July, 17th, 1998}.)

\emp

\bmp \textbf{11.24.} A Fitting class ${\goth F}$ is
said to be {\it local\/} if there exists a group function $f$ (for
definition see (L.\,A.\,Shemetkov, {\it Formations of Finite Groups},
Moscow, Nauka, 1978 (Russian)) such that $f(p)$ is a Fitting class
for every prime number $p$ and ${\goth F}={\goth G}_{\pi ( {\goth
F})} \bigcap \left(\bigcap_{p\in \pi ({\goth F})} \,f(p){\goth
N}_p{\goth G}_{p'}\right) $. Is every hereditary Fitting class of
finite groups local? \hf {\sl N.\,T.\,Vorob'\"ev}\vs

 No, not every
(L.\,A.\,Shemetkov, A.\,F.\,Vasil'yev, {\it Abstracts of the Conf.
of Mathematicians of Belarus', Part 1}, Grodno, 1992, p.~56
(Russian); \ S.\,F.\,Kamornikov, {\it Math. Notes}, {\bf 55},
no.\,6 (1994),
 586--588).

\emp

\bmp \textbf{11.25.} a) Does there exist a local
product (different from the class of all finite groups and from the
class of all finite soluble groups) of Fitting classes each of which
is not local and is not a formation? See the definition of the {\it
product\/} of Fitting classes in (N.\,T.\,Vorob'ev, {\it Math.
Notes}, {\bf 43}, No\,\,1--2 (1988), 91--94). \hf {\sl
N.\,T.\,Vorob'\"ev}\vs

Yes, there does (N.\,T.\,Vorob'ev, A.\,N.\,Skiba, {\it
Problems in Algebra}, {\bf 8}, Gomel', 1995, 55--58 (Russian)).

\emp

\bmp \textbf{11.26.} Does there exist a group which is not isomorphic to
outer auto\-mor\-phism group of a metabelian group with trivial
center? \hfill {\sl R.\,G\"obel}\vs

No, given any group $G$ there is a metabelian group $M$ with
trivial center such that ${\rm Out}\, M \cong G$ (R.\,G\"obel,
A.\,Paras, {\it J.~Pure Appl. Algebra}, {\bf 149}, no.\,3 (2000)
251--266), and if $G$ is finite or countable then $M$ above can be
chosen countable (R.\,G\"obel, A.\,Paras, in: {\it Abelian Groups
and Modules, Proc. Int. Conf. Dublin, 1998}, Birkh\"auser, Basel,
1999, 309--317).

\emp

\bmp \textbf{11.27.} What are the minimum numbers of generators for
groups $G$ satisfying $S \leq G \leq {\rm Aut}\, S$ where $S$ is a
finite simple non-abelian group? \hf {\sl K.\,Gruenberg}\vs

 The number $d(G)$
is found (mod CFSG) for every such a group
 $G$. In particular, $d(G)=\max
\{ 2,\, d(G/S)\}$ and $d(G)\leq 3$ (F.\,\,Dalla\,\,Volta,
A.\,Lucchini, {\it J.~Algebra}, {\bf 178}, no.\,1 (1995),
194--223).

\emp

\bmp \textbf{11.29.} f) Let $F$ be a free group and ${\goth f}={\Bbb
Z}F(F - 1)$ the augmentation ideal of the integral group ring
${\Bbb Z}F$. For any normal subgroup $R$ of $F$ define the
corresponding ideal ${\goth r}={\Bbb Z}F(R - 1)=\, _{{\rm id}} (r
- 1\mid r \in R)$. One may identify, for instance, $F \cap (1 +
{\goth r} {\goth f})=R'$, where $F$ is naturally imbedded into
${\Bbb Z}F$ and $1 + {\goth r} {\goth f}=\{ 1 + a\mid a \in {\goth
r}{\goth f}\} $. Is the quotient group $(F \cap (1 + {\goth r} +
{\goth f}^n))/R\cdot \gamma _n(F)$ always abelian?
 \hf {\sl N.\,D.\,Gupta}\vs

 Yes, it is (N.\,D.\,Gupta, Yu.\,V.\,Kuz'min, {\it J.~Pure Appl. Algebra},
{\bf 78}, no.\,1 (1992), 165--172). \emp

\bmp \textbf{11.33.} a) Let $G(q)$ be a simple Chevalley group over a
field of order $q$. Prove that there exists $m$ such that the
restriction of every non-one-dimensional complex representation of
$G(q^m)$ to $G(q)$ contains all irreducible representations of
$G(q)$ as composition factors.\hf {\sl A.\,E.\,Zalesskii}\vs

This is proved in (D.\,Gluck, {\it J.~Algebra}, {\bf 155}, no.\,2
(1993), 221--237).

\emp

\bmp \textbf{11.35.} Suppose that $H$ is a finite linear group over
${\Bbb C}$ and $h$ is an element of $H$ of prime order $p$ which
is not contained in any abelian normal subgroup. Is it true that
$h$ has at least $(p - 1)/2$ different eigenvalues? \hf {\sl
A.\,E.\,Zalesskii}\vs

 Yes, it is
(G.\,R.\,Robinson, {\it J.~Algebra}, {\bf 178}, no.\,2 (1995),
635--642).

\emp

 \bmp \textbf{11.36.} e) Let $G=B(m,n)$ be the
free Burnside group of rank $m$ and of odd exponent $n \gg 1$. Is it true that
all zero divisors in the group ring ${\Bbb
Z}G$ are trivial? which means that if $ab=0$ then $a=a_1c$,
$b=db_1$ where $a_1, c, b_1, d \in {\Bbb Z}G$, \ $cd=0$, and the
set \ ${\rm
supp}\, c \cup {\rm supp}\, d$ \ is contained in a cyclic subgroup
of $G$. \hfill {\sl S.\,V.\,Ivanov}
\vs

e) No, there are nontrivial zero divisors (S.\,V.\,Ivanov, R.\,Mikhailov, {\it Canad. Math. Bull.}, {\bf 57} (2014), 326--334).
\emp

\bmp \textbf{11.37.} b) Can the free Burnside group $ B(m,n)$, for
any $m$ and $ n=2^l \gg 1$, be given by defining relations of the
form $v^n=1$ such that for any natural divisor $d$ of $n$ distinct
from $n$ the element $v^d$ is not trivial in $B(m,n)$? \hf {\sl
S.\,V.\,Ivanov}\vs

Yes, it can (S.\,V.\,Ivanov, {\it Int. J.~Algebra Comput.}, {\bf
4}, no.\,1--2 (1994), 1--308).

\emp

\bmp \textbf{11.42.} Does there exist
a torsion-free group having exactly three conjugacy classes and
containing a subgroup of index 2? \hfill {\sl
A.\,V.\,Izosov}

\vs

Yes, there does
(A.\,Minasyan, {\it Comment. Math. Helv.}, \textbf{84} (2009),
259--296).

\emp

\bmp \textbf{11.43.} For a finite group $X$, we denote by $k(X)$ the
number of its conjugacy classes. Is it true that $k(AB) \leq
k(A)k(B)$? \hf {\sl L.\,S.\,Kazarin}\vs

 No, it is not true in general:
let $G=\langle a,b\mid a^{30}=b^2=1,\, a^b=a^{-1}\rangle\cong
D_{60}$ be the dihedral group of order $60$, and let $A=\langle
a^{10}, ba\rangle\cong D_6$ and $B=\langle a^6, b\rangle\cong
D_{10}$. Then $G=AB$, but $G$ has $18$ conjugacy classes and $A$
and $B$ only $3$ and $4$, respectively. From (P.\,Gallagher, {\it
Math. Z.}, {\bf 118} (1970), 175--179) a positive answer follows
if $A$ or $B$ is normal. The problem remains open in the case
where $A$ and $B$ have coprime orders, see new problem 14.44.
(J.\,Sangroniz, {\it Letter of December, 17,
 1998}.)

\emp

\bmp \textbf{11.46.} b) Does there exist a finite
$p\hskip0.1ex$-group of nilpotency class greater than 2, with
${\rm Aut}\, G={\rm Aut}_c G\cdot {\rm Inn}\, G$, where ${\rm
Aut}_c G$ is the group of central automorphisms of $G$?

\makebox[25pt][r]{c)}
Does there exist a 2-Engel finite
$p\hskip0.2ex$-group $G$ of nilpotency class greater than 2 such
that ${\rm Aut}\, G={\rm Aut}_c G\cdot {\rm Inn}\, G$, where ${\rm
Aut}_c G$ is the group of central automorphisms of~$G$?
 \hf {\sl A.\,Caranti}\vs

b) Yes, there does (I.\,Malinowska, {\it Rend. Sem. Mat.
Padova}, {\bf 91} (1994), 265--271).

c) Yes, it does (A.\,Abdollahi, A.\,Faghihi, S.\,A.\,Linton,
E.\,A.\,O'Brien,
{\it Arch. Math. (Basel)}, {\bf 95}, no.\,1 (2010), 1--7).
\emp

\bmp \textbf{11.47.} Let ${\cal L}_d$ be the
homogeneous component of degree $d$ in a free Lie algebra ${\cal L}$
of rank 2 over the field of order 2. What is the dimension of the
fixed point space in ${\cal L}_d$ for the automorphism of ${\cal L}$
which interchanges two elements of a free generating set of ${\cal
L}$? \hf {\sl L.\,G.\,Kov\'acs}\vs

 It is found; R.\,M.\,Bryant and R.\,St\"ohr ({\it Arch. Math.}, {\bf
67}, no.\,4 (1996), 281--289) confirmed the conjecture from
(M.\,W.\,Short, {\it Commun. Algebra}, {\bf 23}, no.\,8 (1995),
3051--3057).

\emp

 \bmp \textbf{11.52.} (Well-known problem). A
permutation group on a set $\Omega $ is called {\it sharply doubly
transitive\/} if for any two pairs $(\alpha , \beta )$ and
$(\gamma , \delta )$ of elements of $\Omega$ such that $\alpha \ne
\beta$ and $\gamma \ne \delta $, there is exactly one element of
the group taking $\alpha$ to $\gamma$ and $\beta$ to~$\delta$.
Does every sharply doubly transitive group possess a non-trivial
abelian normal subgroup? A positive answer is well known for
finite groups. \hfill {\sl V.\,D.\,Mazurov}
\vs

No, not every (E.\,Rips, Y.\,Segev, K.\,Tent, {\it J.~Europ. Math. Soc.}, {\bf 19}, no.\,10 (2017), 2895--2910).
\emp

\bmp \textbf{11.53.} (P.\,Kleidman). Do the sporadic
simple groups of Rudvalis $ Ru$, Mathieu $M_{22}$, and Higman--Sims
$HS$ embed into the simple group $E_7(5)$? \hf {\sl
V.\,D.\,Mazurov}\vs

 Yes, they do
(P.\,B.\,Kleidman, R.\,A.\,Wilson, {\it J.~Algebra}, {\bf 157},
no.\,2 (1993), 316--330).

\emp

\bmp \textbf{11.54.} Is it true that in the group
of coloured braids only the identity braid is a~conjugate to its
inverse? (For definition see Archive, 10.24.) \hfill {\sl
G.\,S.\,Makanin}\vs

 {\bf 11.55.} Is it true that extraction of roots in the
group of coloured braids is uniquely determined? \hf {\sl
G.\,S.\,Makanin}\vs

 The answers to both 11.54 and 11.55 are affirmative,
 since the group $K_n$ of coloured braids is
embeddable into the group of those automorphisms of the free group
$F_n$ that act trivially modulo the derived subgroup of $F_n$,
and this group is residually in the class of torsion-free
nilpotent groups, for which the corresponding assertions are true
(V.\,A.\,Roman'kov, {\it Letter of October, 3, 1990\/}).
See also (V.\,G.\,Bardakov,
 {\it Russ. Acad. Sci. Sbornik Math.}, {\bf 76}, no.\,1
 (1993), 123--153).
 \emp

\bmp \textbf{11.57.} An {\it upper composition
factor\/} of a group $G$ is a composition factor of some finite
quotient of $ G$. Is there any restriction on the set of non-abelian
upper composition factors of a finitely generated group? \hfill {\sl
\,A.\,Mann, \,D.\,Segal}\vs

There are no restrictions: any set of non-abelian finite simple
groups can be the set of upper composition factors of a 63-generator
group; if we allow the group to have also abelian upper composition
factors, the number of generators can be reduced to~3 \ (D.\,Segal,
 {\it Proc. London Math. Soc. (3)}, {\bf 82} (2001), 597--613).
\emp

\bmp \textbf{11.64.} Let $\pi (G)$ denote the set of
prime divisors of the order of a finite group~$G$. Are there only
finitely many finite simple groups $G$, different from alternating
groups, which have a proper subgroup $H$ such that $\pi (H)=\pi (G)$?
\hf {\sl V.\,S.\,Monakhov}\vs

 No, there are
infinitely many such groups. If
 $G=S_{4k}(2^s)$ and $H\cong \Omega ^{-}_{4k}(2^s)$, then $\pi
 (G)=\pi (H)$ (V.\,I.\,Zenkov, {\it
Letter of March, 10, 1994}.)

\emp

\bmp \textbf{11.68.}
 Can every fully ordered
group be embedded in a fully ordered group (continuing the given
order) with only 3 classes of conjugate elements? \hfill
 {\sl B.\,Neumann}
\vs

No, not every
(V.\,V.\,Bludov, {\it Algebra and Logic}, \textbf{44}, no.\,6
(2005), 370--380). \emp 

\bmp \textbf{11.70.} b) Let $F$ be an infinite
field or a skew-field.
What conditions on the subring $R$ of $F$
will ensure that $PGL(d + 1,\, R)$ is flag-transitive on the
projective $d$-space $PG(d,F)$?
\hfill {\sl P.\,M.\,Neumann, C.\,E.\,Praeger}
\vs

b) Such necessary and sufficient conditions are given in (S.\,A.\,Zyubin, {\it Siberian Electron. Math. Rep.}, {\bf 11} (2014), 64--69 (Russian)).
\emp

\bmp \textbf{11.74.} Let $G$ be a non-elementary
hyperbolic group and let $G^n$ be the subgroup generated by the $n$th
powers of the elements of $G$.

\mb{a)} (M.\,Gromov). Is it true that $G/G^n$ is
infinite for some $n = n(G)$?

\mb{b)} Is it true that $\bigcap\limits_{n=1}^{\infty} G^n=\{ 1\}
$? \hf {\sl A.\,Yu.\,Olshanskii}\vs

a) Yes, it is; b) Yes, it is (S.\,V.\,Ivanov,
A.\,Yu.\,Olshanskii, {\it Trans. Amer. Math. Soc.}, {\bf 348}, N
6 (1996), 2091--2138).

\emp

\bmp \textbf{11.75.} Let us consider the class of groups with $n$
generators and $ m$ relators. A~subclass of this class is called
{\it dense\/} if the ratio of the number of presentations of the
form $\left< a_1,\ldots ,a_n\mid R_1,\ldots ,R_m\right>$ (where
$|R_i|=d_i$) for groups from this subclass to the number of all
such presentations converges to 1 when $d_1 + \cdots + d_m$ tends
to infinity. Prove that for every $k < m$ and for any $n$ the
subclass of groups all of whose $k$-gene\-ra\-tor subgroups are
free is dense. \hfill {\sl A.\,Yu.\,Olshanskii}\vs

This is proved (G.\,N.\,Arzhantseva, A.\,Yu.\,Olshanskii,
{\it Math. Notes}, {\bf 59}, no.\,4 (1996), 350--355). \emp

\bmp \textbf{11.79.} Let $G$ be a finite group of automorphisms of an
infinite field $F$ of characteristic~$p$. Taking integral powers of
the elements of $F$ and the action of~$G$ define the action of the
group ring ${\Bbb Z}G$ of $G$ on the multiplicative group of~$F$. Is
it true that any subfield of $F$ that contains the images of all
elements of $F$ under the action of some fixed element of ${\Bbb
Z}G\setminus p{\Bbb Z}G$ contains infinitely many $ G$-inva\-ri\-ant
elements of~$F$?

\hf {\sl K.\,N.\,Ponomar\"ev}\vsh

 Yes
(K.\,N.\,Ponomar\"ev, {\it Siber. Math.~J.}, {\bf 33} (1992),
1094--1099). \emp

 \bmp \textbf{11.82.} Let $R$ be the normal
closure of an element $r$ in a free group $F$ with the natural length
function and suppose that $s$ is an element of minimal length in $R$.
Is it true that $s$ is conjugate to one of the following elements:
$r$, $r^{-1}$, $[r,f]$, $[r^{-1},f]$ for some $f \in F$? \hfill {\sl
V.\,N.\,Remeslennikov}\vs

No, not always (J.\,McCool, {\it Glasgow Math.~J.}, {\bf 43},
no.\,1 (2001), 123--124). \emp

 \bmp \textbf{11.88.}
We define the {\it length\/} $l(g)$ of an Engel element $g$ of a
group $G$ to be the smallest number $l$ such that $[h,g;\, l]=1$
for all $h \in G$. Here $[h,g;\, 1]=[h,g]$ and $[h,g;\, i +
1]=[[h,g;\, i],g]$. Does there exist a polynomial
function $\varphi (x,y)$ such that $l(uv) \leq \varphi
(l(u),l(v))$? Up to now, it is unknown whether a product of Engel
elements is again an Engel element.
\hfill
{\sl V.\,A.\,Roman'kov}

\vs

No, it does not (L.\,V.\,Dolbak, {\it
Siberian Math.~J.}, \textbf{47}, no.\,1 (2006), 55--57). 
\emp

 \bmp \textbf{11.91.}
Prove that a hereditary formation ${\goth F}$ of finite soluble
groups is local if every finite soluble non-simple minimal
non-${\goth F}$-group is a Shmidt group (that is, a non-nil\-po\-tent
finite group all of whose proper subgroups are nilpotent). \hfill{\sl
V.\,N.\,Semenchuk}\vs

This is proved (A.\,N.\,Skiba, {\it Dokl. Akad. Nauk Belorus.
SSR}, {\bf 34}, no.\,11 (1990), 982--985 (Russian)). \emp

\bmp \textbf{11.93.} Is the variety of all lattices generated by the
block lattices (see Archive, 9.62) of finite groups? \hf {\sl
D.\,M.\,Smirnov}\vs

 No, it is not
(D.\,M.\,Smirnov, {\it Siberian Math.~J.}, {\bf 33}, no.\,4
(1992), 663--668).

\emp

\bmp \textbf{11.94.}
 Describe all {\it simply
reducible\/} groups, that is, groups such that all their
characters are real and the tensor product of any two irreducible
representations contains no multiple components. This question is
interesting for physicists. Is every finite simply reducible group
soluble? \hfill {\sl S.\,P.\,Strunkov}

\vs

Yes, it is
(L.\,S.\,Kazarin, E.\,I.\,Chankov,
{\it Sb. Math.}, {\bf 201}, no.\,5 (2010), 655--668).
\emp

\bmp \textbf{11.97.} Are there only finitely many
finite simple groups with a given set of all different values of
irreducible characters on a single element? \hfill {\sl
S.\,P.\,Strunkov}\vs

 No; all
complex irreducible characters of the groups $L_2(2^m)$, $m\geq
2$, take the values $0,\, \pm1$ on involutions (V.\,D.\,Mazurov).

\emp

\bmp \textbf{11.98.} b) Let $r$ be the number of conjugacy classes of
elements in a finite (simple) group~$G$. Is it true that $|G|\leq
{\rm exp}(r)$? \hfill {\sl S.\,P.\,Strunkov}\vs

No, this is not true: the group $M_{22}$ has order 443520 and
contains 12 conjugacy classes (T.\,Plunkett, {\it Letter of May, 9,
2000\/}).

\emp

\bmp \textbf{11.101.}
Does there
exist a Golod group (see 9.76) with finite centre?

\hfill {\sl A.\,V.\,Timofeenko}

\vs

Yes, there is a Golod group with trivial centre
(V.\,A.\,Sereda, A.\,I.\,Sozutov, {\it Algebra and Logic}, {\bf
45}, no.\,2 (2006), 134--138).
\emp

 \bmp \textbf{11.103.} Is a
2-group satisfying the minimum condition for centralizers necessarily
locally finite? \hfill {\sl John S.\,Wilson}\vs

Yes, it is (F.\,O.\,Wagner, {\it J.~Algebra}, {\bf 217}, no.\,2
(1999),
 448--460). \emp

\bmp
\textbf{11.104.} Let $G$ be a finite
group of order $p^a\cdot q^b\cdots\,$, where $p, q, \ldots $ are
distinct primes. Introduce distinct variables $x_p, x_q, \ldots $
corresponding to $p, q, \ldots $~. Define functions $f, \varphi$
from the lattice of subgroups of $G$ to the polynomial ring ${\Bbb
Z}[x_p,x_q,\ldots ]$ as follows:

\vskip0.5ex
\makebox[15pt][r]{}(1) if $H$ has order $p^{\alpha}\cdot q^{\beta
}\cdots\,$, then $f(H)=x_p^{\alpha}\cdot x_q^{\beta}\cdots\,$;

\vskip0.5ex
\makebox[15pt][r]{}(2) for all $H \leq G$, we have
$\sum\limits_{K\leq H}\varphi (K)=f(H)$.

\vskip0.5ex Then $f(G), \varphi (G)$ may be called the {\it
order\/} and {\it Eulerian polynomials\/} of~$G$. Substituting
$p^m, q^m, \ldots $ for $x_p, x_q, \ldots $ in these polynomials
we get the $m$th power of the order of~$G$ and the number of
ordered $m$-tuples of elements that generate $G$ respectively.

\makebox[15pt][r]{}It is known that if $G$ is
$p\hskip0.2ex$-solvable, then $\varphi (G)$ is a product of a
polynomial in~$x_p$ and a polynomial in the remaining variables.
Consequently, if $G$ is solvable, $\varphi (G)$ is the product of
a polynomial in~$x_p$ by a polynomial in~$x_q$ by~$\ldots $~. Are
the converses of these statements true \ a) for solvable groups? \
b) for $p\hskip0.2ex$-solvable groups? \hfill {\sl G.\,E.\,Wall}

\vs

Yes, the converses are
true:
 a) for~solvable groups (E.\,Detomi, A.\,Lucchini, {\it J.~London
Math. Soc. (2)}, \textbf{70} (2004), 165--181); \ b)~for
$p\hskip0.2ex$-solvable groups (E.\,Damian, A.\,Lucchini, {\it
Commun. Algebra}, \textbf{35} (2007), 3451--3472).

 \emp

\bmp \textbf{11.105.} a) Let ${\goth V}$ be a variety of groups. Its
relatively free group of given rank has a presentation $F/N$, where
$F$ is absolutely free of the same rank and $N$ fully invariant in
$F$. The associated Lie ring ${\cal L}(F/N)$ has a presentation ${
L}/J$, where ${ L}$ is the free Lie ring of the same rank and $J$ an
ideal of ${ L}$. Is $J$ always fully invariant in ${ L}$?

 \hfill {\sl G.\,E.\,Wall}\vs

No, not always; the ideal $J$ is not fully invariant for
$F/(F^2)^4$, that is, for ${\goth V}={\goth B}_4{\goth B}_2$
(D.\,Groves, {\it J.~Algebra}, {\bf 211}, no.\,1 (1999), 15--25).
\emp

\bmp \textbf{11.106.} Can every periodic group be embedded in a simple
periodic group?

 \hf {\sl R.\,Phillips}\vsh

Yes (A.\,Yu.\,Olshanskii, {\it Ukrain. Math.~J.}, {\bf 44} (1992),
761--763). \emp

\bmp \textbf{11.108.} Is every locally finite simple group absolutely
simple? A group $G$ is said to be {\it absolutely simple} if the only
composition series of $G$ is $\{ 1, \, G\} $. For equivalent
formulations see (R.\,E.\,Phillips, {\it Rend. Sem. Mat. Univ.
Padova}, {\bf 79} (1988), 213--220). \hfill {\sl R.\,Phillips}\vs

No, not every (U.\,Meierfrankenfeld, in: {\it Proc. Int. Conf. Finite
and Locally Finite Groups, Istanbul, 1994}, Kl\"uwer, 1995,
189--212).

\emp

\bmp \textbf{11.110.} Is it possible to embed $SL(2,{\Bbb Q})$ in the
multiplicative group of some division ring? \hf {\sl
B.\,Hartley}\vs

No, it is not (W.\,Dicks, \,B.\,Hartley, {\it Commun. Algebra},
{\bf 19} (1991), 1919--1943). \emp

\bmp \textbf{11.112.}
b) Let $L=L(K(p))$ be the associated Lie ring of a free countably
generated group $K(p)$ of the Kostrikin variety of locally finite
groups of a given prime exponent $p$. Is it true that all
identities of $L$ follow from multilinear identities of $L$?
\hfill {\sl E.\,I.\,Khukhro}\vs

No, there are two relations of weight 29 which hold in $L(K(7))$
(in two generators, of multiweights $(14,15)$ and $(15,14)$) that
are not consequences of multilinear relations (E.\,O'Brien,
M.\,R.\,Vaughan-Lee, {\it Int. J. Algebra Comput.}, {\bf 12}
(2002), 575--592; \ M.\,F.\,Newman, M.\,R.\,Vaughan-Lee, {\it
Electron. Res. Announc. Amer. Math. Soc.}, {\bf 4}, no.\,1 (1998),
1--3). \emp

\bmp \textbf{11.126.}
Do there exist a constant $h$ and a function $f$ with the
following property: if a finite soluble group $G$ admits an
automorphism $\varphi$ of order 4 such that $|C_G(\varphi )| \leq
m$, then $G$ has a normal series
 $1\leq M\leq N\leq G$ such that the index
$|G:N|$ does not exceed $f(m)$, the group $N/M$ is nilpotent of
class $\leq 2$, and the group $M$ is nilpotent of class $\leq h$?
\hfill {\sl P.\,V.\,Shumyatski\u{\i}}

\vs

Yes, there do (N.\,Yu.\,Makarenko, E.\,I.\,Khukhro, {\it Algebra and Logic},
\textbf{45} (2006), 326--343).

\emp

\bmp \textbf{11.128.} A group $G$ is said to be a {\it $K$-group\/} if
for every subgroup $A \leq G$ there exists a subgroup $B \leq G$
such that $A \cap B=1$ and $\left< A,B\right> =G$. Is it true that
normal sub\-groups of $K$-groups are also $K$-groups? \hf {\sl
M.\,Emaldi}\vs

 No, it is not
(V.\,N.\,Obraztsov, {\it J.~Austral. Math. Soc. (Ser.~A)}, {\bf
61}, no.\,2 (1996), 267--288).

\emp

\bmp \textbf{12.1.} b) (A.\,A.\,Bovdi). H.\,Bass ({\it Topology},
{\bf 4}, no.\,4 (1966), 391--400) has constructed explicitly a
proper subgroup of finite index in the group of units of the
integer group ring of a finite cyclic group. Construct analogous
subgroups for finite abelian non-cyclic groups. \hfill {\sl
R.\,Zh.\,Aleev}\vs

They are constructed: for $p\hs$-groups,
 $p\ne 2$, in
 (K.\,Hoechsmann, J.\,Ritter, {\it J.~Pure Appl. Algebra}, {\bf
68}, no.\,3 (1990), 325--333); for the general case of groups of
central units of integral group rings of arbitrary (not
necessarily abelian) finite groups in (R.\,Zh.\,Aleev, {\it Mat.
Tr. Novosibirsk Inst. Math.}, {\bf 3}, no.\,1 (2000), 3--37
(Russian)). \emp
\markboth{\protect\vphantom{(y}{Archive of solved
problems (12th ed., 1992)}}{\protect\vphantom{(y}{Archive of solved
problems (12th ed., 1992)}}

\bmp \textbf{12.2.} Is it true that every non-discrete group topology
(in an abelian group) can be strengthened up to a maximal complete
group topology? \hf {\sl V.\,I.\,Arnautov}\vs

 No, not always, since
under the assumption of the Continuum Hypothesis there exist
non-complete maximal topologies (V.\,I.\,Arnautov,
E.\,G.\,Zelenyuk, {\it Ukrain. Math.~J.}, {\bf 43}, no.\,1 (1991),
15--20). \emp

 \bmp \textbf{12.5.} Does there exist a countable
non-trivial filter in the lattice of quasivarieties of metabelian
torsion-free groups? \hf {\sl A.\,I.\,Budkin}\vs

 No, there is no such filter
(S.\,V.\,Lenyuk, {\it Siberian Math.~J.}, {\bf 39}, no.\,1 (1998),
57--62).

\emp

\bmp \textbf{12.7.} Is it true that every radical hereditary
formation of finite groups is a com\-position one? \hfill {\sl
A.\,F.\,Vasil'ev}\vs

 No, it is not
(S.\,F.\,Kamornikov, {\it Math.
 Notes}, {\bf 55}, no.\,6 (1994), 586--588).

\emp

\bmp \textbf{12.10.} (P.\,Neumann). Can the free group on two
generators be embedded in ${\rm Sym}\,({\Bbb N})$ so that the
image of every non-identity element has only a finite number of
orbits?

\hf {\sl A.\,M.\,W.\,Glass}\vs

 Yes, it can
(H.\,D.\,Macpherson, in: {\it Ordered groups and infinite permutation
groups, Partially based on Conf. Luminy, France, 1993\/} ({\it
Mathematics and its Applications}, {\bf 354}), Kluwer, Dordrecht,
1996, 221--230).

\emp

\bmp \textbf{12.14.}
If $T$ is a
countable theory, does there exist a model ${\frak A}$ of $T$ such
that the theory of ${\rm Aut}\, ({\frak A})$ is undecidable?
\hfill {\sl M.\,Giraudet, A.\,M.\,W.\,Glass}
\vs

Yes, moreover, every first
order theory having infinite models has a model whose automorphism
group has undecidable existential theory (V.\,V.\,Bludov,
M.\,Giraudet, A.\,M.\,W.\,Glass, G.\,Sabbagh, in {\it Models,
Modules and Abelian Groups}, de Gruyter, Berlin, 2008, 325--328).
\emp

\bmp \textbf{12.22.}  Let $\Delta (G)$ be the augmentation ideal of
the integer group ring of an arbitrary group $G$. Then $D_n(G)=G \cap
(1 + \Delta ^n(G))$ contains the $n $th lower central subgroup
$\gamma _n(G)$ of $G$.

\makebox[25pt][r]{a)} Is it true that $D_n(G)/\gamma_n(G)$ is
central in $G/\gamma _n(G)$?

\makebox[25pt][r]{b)} Is it true that $D_n(G)/\gamma _n(G)$ has exponent dividing 2?

 \hf {\sl N.\,D.\,Gupta,~Yu.\,V.\,Kuz'min}\vs

 a) No, not always; moreover,
$D_n(G)/\gamma_n(G)$ need not be contained in any term of the
upper central series of $G/\gamma _n(G)$ with fixed number
(N.\,D.\,Gupta, Yu.\,V.\,Kuz'min, {\it J.~Pure Appl. Algebra},
{\bf 104}, no.\,1 (1995), 191--197).
\vs

b) No, not always (L.\,Bartholdi, R.\,Mikhailov, {\it J.~Topology}, {\bf 16}, no.\,2 (2023), 822--853).
\emp

\bmp \textbf{12.24.} Given a ring $R$ with identity, the automorphisms
of $R[[x]]$ sending $x$ to $x\left( 1 + \sum_{i=1}^{\infty}
a_ix^i\right) $, $a_i\in R$, form a group $N(R)$. We know that
$N({\Bbb Z})$ contains a copy of the free group $F_2$ of rank 2 and,
from work of A.\,Weiss, that $N({\Bbb Z}/p{\Bbb Z})$ contains a copy
of every finite $p\hskip0.1ex $-group (but not of ${\Bbb
Z}_{p^{\infty}})$, $p$ a prime. Does $N({\Bbb Z}/p{\Bbb Z})$ contain
a copy of~$F_2$? \hf {\sl D.\,L.\,Johnson}\vs

 Yes, it does (R.\,Camina,
{\it J.~Algebra}, {\bf 196}, no.\,1 (1997), 101--113).
I.\,B.\,Fesenko noted that this fact might have been known to
specialists in the theory of fields of norms back in 1985
(J.-P.\,Wintenberger, J.-M.\,Fontaine, F.\,Laubie); however a
proof based on this theory first appeared only in
(I.\,B.\,Fesenko, {\it Preprint}, Nottingham Univ., 1998).

\emp

\bmp \textbf{12.25.} Let $G$ be a finite group acting irreducibly on a
vector space $V$. An orbit $\alpha ^G$ for $\alpha \in V$ is said to
be {\it $p$-regu\-lar\/} if the stabilizer of $\alpha$ in $G$ is a
$p' $-sub\-group. Does $G$ have a regular orbit on $V$ if it has a
$p\hskip0.1ex $-regu\-lar orbit for every prime $p$? \hfill {\sl
Jiping~Zhang}\vs

 No, not always.
For example, let $H={\Bbb A}_4\wr {\Bbb Z}_5$, \ $V=O_2(H)$, and $G$
be a complement to $V$ in $H$, $G$ acting on
 $V$ by conjugation.
(V.\,I.\,Zenkov, {\it Letter of February, 11, 1994}.)

\emp

\bmp \textbf{12.26.} (Shi Shengming). Is it true that a finite
$p\hskip0.1ex $-soluble group $G$ has a $p\hskip0.1ex $-block of
defect zero if and only if there exists an element $x \in O_{p'}(G)$
such that $C_G(x)$ is a $p' $-sub\-group? \hf {\sl Jiping Zhang}\vs

 No, it is not. The
group $3^2 . GL_2(3)$ has a 2-block of defect $0$, but the
centralizer of every element from $O(G)$ has even order
(V.\,I.\,Zenkov, in: {\it Trudy Inst. Matem. i Mekh. UrO RAN}, {\bf
3} (1995), 36--40 (Russian)).

\emp

\bmp \textbf{12.31.} For relatively free groups $G$, prove or disprove
the following conjecture of P.\,Hall: if a word $v$ takes only
finitely many values on $G$ then the verbal subgroup $vG$ is finite.
\hfill {\sl S.\,V.\,Ivanov}

The conjecture is
 disproved (A.\,Storozhev, {\it Proc. Amer. Math. Soc.}, {\bf
 124}, no.\,10 (1996), 2953--2954).
\emp

\bmp \textbf{12.32.}
Prove an analogue of Higman's theorem for the
Burnside variety ${\frak B}_n$ of groups of odd exponent $n \gg
1$, that is, prove that every recursively presented group of
exponent $n$ can be embedded in a finitely presented (in~${\frak
B}_n$) group of exponent $n$.

\hfill
{\sl S.\,V.\,Ivanov}

 It is proved (A.\,Yu.\,Olshanskii, {\it J.~Algebra}, \textbf{560} (2020), 960--1052).
\emp

\bmp \textbf{12.36.} Let $p$ be a prime, $V$ an $n $-dimen\-sional
vector space over the field of $p$~elements, and let $G$ be a
subgroup of $GL(V)$. Let $S=S[V ^{\displaystyle{*}} ]$ be the
symmetric algebra on $V ^{\displaystyle{*}} $, the dual of $V$. Let
$T=S^G$ be the ring of invariants and let $b_m$ be the dimension of
the homogeneous component of degree $m$. Then the Poincar\'e series
$\sum\limits_{m\geq 0}b_mt^m$ is a rational function with a Laurent
power series expansion $\sum\limits_{i\geq - n}a_i(1 - t)^i$ about
$t=1$ where $a_{- n}=\frac{\displaystyle{1}}{\displaystyle{|G|}}$.

\makebox[15pt][r]{}{\it Conjecture:\/} $a_{- n+1}=
\frac{\displaystyle{r}}{\displaystyle{2|G|}}$ where
$r=\sum\limits_{W} ((p - 1)\alpha _W + s_W - 1)$, the sum is taken
over all maximal subspaces $W$ of $V$, and $\alpha _W$, $s_W$ are
defined by $|G_W|=p^{\alpha _W}\cdot s_W$ where $p\nmid s_W$ and
$G_W$ denotes the pointwise stabilizer of $W$.

\hfill {\sl D.\,Carlisle, P.\,H.\,Kropholler}\vs

The conjecture is proved (D.\,J.\,Benson,
 W.\,W.\,Crawley-Boevey, {\it Bull. London Math. Soc.}, {\bf
27}, no.\,5 (1995),
 435--440); another proof based on the Grothendieck--Riemann--Roch
Theorem was later obtained in (A.\,Neeman, {\it Comment. Math.
Helv.}, {\bf 70}, no.\,3 (1995), 339--349).
\emp

\bmp \textbf{12.38.}
(J.\,G.\,Thompson). For a finite group $G$, we
denote by $N(G)$ the set of all orders of the conjugacy classes of
$G$. Is it true that if $G$ is a finite non-abelian simple group,
$H$ a finite group with trivial centre and $N(G)=N(H)$, then $G$
and $H$ are isomorphic? \hfill {\sl
A.\,S.\,Kondratiev, W.\,J.\,Shi}

\vs Yes, it is true. The final step of the proof is in the paper (I.\,B.\,Gorshkov, {\it Commun. Algebra}, {\bf 47}, no.\,12 (2019), 5192--5206),
which contains references to the previous steps by M.\,Ahanjideh, N.\,Ahanjideh, S.\,H.\,Alavi, G.\,Y.\,Chen, A.\,Daneshkhah, M.\,R.\,Darafsheh, I.\,B.\,Gorshkov, A.\,Iranmanesh, I.\,Kaygorodov, Behn.\,Khosravi, Behr.\,Khosravi, A.\,Kukharev, W.\,Shi, A.\,Shlepkin, A.\,V.\,Vasil'ev, L.\,Wang, M.\,Xu.
\emp

\bmp \textbf{12.39.} (W.\,J.\,Shi). Must a finite group
and a finite simple group be isomorphic if they have equal orders
and the same set of orders of elements?
\hf {\sl
A.\,S.\,Kondratiev}

\vs

Yes, they must
(M.\,C.\,Xu, W.\,J.\,Shi, {\it Algebra Colloq.}, {\bf
10} (2003), 427--443; A.V.\,Vasil'ev, M.A.\,Grechkoseeva, V.D.\,Mazurov, {\it
Algebra Logic}, \textbf{48} (2009), 385--409). 
\emp

\bmp \textbf{12.42.}
Describe the automorphisms of the Sylow $p\hskip0.2ex $-sub\-group
of a Chevalley group of normal type over ${\Bbb Z}/p^m{\Bbb Z} $,
$m \geq 2$, where $p$ is a prime. \hfill {\sl
V.\,M.\,Levchuk}

\vs

Described in
(S.\,G.\,Kolesnikov, {\it Algebra and Logic}, \textbf{43}, no.\,1
(2004), 17--33); {\it Izv. Gomel' Univ.}, \textbf{36}, no.\,3 (2006),
137--146 (Russian); {\it J.~Math. Sci., New York}, \textbf{152}
(2008), 220--246).

\emp
 \bmp \textbf{12.44.}
(P.\,Hall). Is there a non-trivial group which is isomorphic with
every proper extension of itself by itself? \hfill {\sl
J.\,C.\,Lennox}\vs

Yes, such groups exist of cardinality any regular cardinal
(R.\,G\"obel, S.\,Shelah, {\it Math. Proc. Cambridge Philos.
Society (2)}, {\bf 134}, no.\,1 (2003), 23--31). \emp

\bmp \textbf{12.45.} (P.\,Hall). Must a non-trivial group, which is
isomorphic to each of its non-trivial normal subgroups, be either
free of infinite rank, simple, or infinite cyclic? (Lennox, Smith
and Wiegold, 1992, have shown that a finitely generated group of
this kind which has a proper normal subgroup of finite index is
infinite cyclic.) \hf {\sl J.\,C.\,Lennox}\vs

 No (V.\,N.\,Obraztsov, {\it Proc. London Math. Soc.}, {\bf 75}
 (1997), 79--98). \emp

\bmp \textbf{12.46.} Let $F$ be the nonabelian free group on two
generators $x,y$. For $a,b \in {\Bbb C}$, $|a|=|b|=1$, let $\vartheta
_{a,b}$ be the automorphism of ${\Bbb C}F$ defined by $\vartheta
_{a,b}(x)=ax$, $\vartheta _{a,b}(y)=by$. Given $0 \ne \alpha \in
{\Bbb C}F$, can we always find $a,b \in {\Bbb C}\setminus \{ 1\} $
with $|a|=|b|=1$ such that $\alpha {\Bbb C}F \cap \vartheta
_{a,b}(\alpha ){\Bbb C}F \ne 0$? \hf {\sl P.\,A.\,Linnell}\vs

 No, not always (A.\,V.\,Tushev, {\it Ukrain. Mat.~J.}, {\bf
47}, no.\,4 (1995), 571--572). \emp

\bmp \textbf{12.47.}
Let $k$ be a field, let $p$ be a prime, and let
$G$ be the Wreath product ${\Bbb Z}_p\wr {\Bbb Z}$ (so the base
group has exponent $p$). Does $kG$ have a classical quotient ring?
(i.~e. do the non-zero-divisors of $kG$ form an Ore set?)
 \hfill {\sl P.\,A.\,Linnell}\vs

No, it does not (P.\,A.\,Linnell, W.\,L\"uck, T.\,Schick, in: {\it
High-dimensional manifold topology. Proc. of the school, ICTP
(Trieste, 2001)}, World Scientific, River Edge, NJ, 2003,
315--321).

\emp

\bmp \textbf{12.49.} Construct all non-split extensions of elementary
abelian 2-groups $V$ by $H=PSL_2(q)$ for which $H$ acts irreducibly
on $V$. \hfill {\sl V.\,D.\,Mazurov}\vs

They are constructed (V.\,P.\,Burichenko, {\it Algebra and Logic},
{\bf 39} (2000),
 160--183).
\emp

\bmp \textbf{12.65.} Let ${\cal P}=(P_0,P_1,P_2)$ be a parabolic system
in a finite group $G$, belonging to the $C_3$ Coxeter diagram \
$\parbox[t]{2ex}{ $ \circ $ \\ \mathstrut $\hskip0.2ex{}^1$}
\rule[0.5ex]{8ex}{0.1ex}\, \;
\parbox[t]{2ex}{ $ \circ $ \\ \mathstrut $\hskip0.2ex{}^2$}
\raisebox{0.4ex}{\underline{\rule[0.5ex]{8ex}{0.1ex}}}\, \;
\parbox[t]{2ex}{ $ \circ $ \\ \mathstrut $\hskip0.2ex{}^3$}
$ \ and let the Borel subgroup have index at least 3 in $P_0$ and
$P_1$. It is known that, if we furthermore assume that the chamber
system of ${\cal P}$ is geometric and that the projective planes
arising as $\{ 0,1\} $-re\-si\-dues from ${\cal P}$ are
desarguesian, then either $G^{(\infty )}$ is a Chevalley group of
type $C_3$ or $B_3$ or $G=A_7$. Can we obtain the same conclusion
in general, without assuming the previous two hypothesis? \hf {\sl
A.\,Pasini}\vs

 Yes, we can
(S.\,Yoshiara, {\it J.~Algebraic Combin.}, {\bf 5}, no.\,3 (1996),
251--284). \emp

\bmp \textbf{12.70.} Let $p$ be a
prime number, $F$ a free pro-$p\hskip0.1ex $-group of finite rank,
and $\Delta \ne 1$ an automorphism of $F$ whose order is a power of
$p$. Is the rank of ${\rm Fix}_F(\Delta )=\{ x \in F\mid \Delta
(x)=x\}$ finite?
 If the order of $\Delta$ is prime to $p$, then
${\rm Fix}_F(\Delta )$ has infinite rank (W.\,Herfort, L.\,Ribes,
{\it Proc. Amer. Math. Soc.}, {\bf 108} (1990), 287--295). \hfill
{\sl L.\,Ribes}\vs

Yes, it is finite (C.\,Scheiderer, {\it Proc. Amer. Math. Soc.},
{\bf 127}, no.\,3 (1999),
 695--700).
Moreover, in (W.\,N.\,Herfort, L.\,Ribes, P.\,A.\,Zalesskii, {\it
Forum Math.}, {\bf 11}, no.\,1 (1999), 49--61) it is proved that
if $P$ is a finite $p\hs$-group of automorphisms of a free pro-$p$
group $F$ of arbitrary rank, then the group of fixed points of $P$
is a free factor of $F$ and the latter result has been extended in
(P.\,A.\,Zalesskii, {\it J.~Reine Angew. Math.}, {\bf 572} (2004),
97--110) to the situation of a free pro-$p$ group of arbitrary
rank. \emp

\bmp \textbf{12.71.} Let $d(G)$ denote the smallest cardinality of a
generating set of the group~$G$. Let $A$ and $B$ be finite groups.
Is there a finite group $G$ such that $A,B \leq G$, $G=\left<
A,B\right>$, and $d(G)=d(A) + d(B)$? The corresponding question
has negative answer in the class of solvable groups
(L.\,G.\,Kov\'acs, H.-S.\,Sim, {\it Indag. Math.}, {\bf 2}, no.\,2
(1991), 229--232).

\hfill {\sl L.\,Ribes}\vsh

No, not always (A.\,Lucchini, {\it J. Group Theory}, {\bf 4},
no.\,1 (2001), 53--58). \emp

\bmp \textbf{12.74.} Let ${\frak F}$ be a non-primary one-generator
composition formation of finite groups. Is it true that if ${\frak
F}={\frak M}{\frak H}$ and the formations ${\frak M}$ and ${\frak H}$
are non-trivial, then ${\frak M}$ is a composition formation? \hfill
{\sl A.\,N.\,Skiba}\vs

This is true if ${\frak F}\ne {\frak H}$ (A.\,N.\,Skiba, {\it
Algebra of formations}, Belaruskaja Navuka, Minsk, 1997, p.\,144
(Russian)). In general the answer is negative (W.\,Guo, {\it
Commun. Algebra}, {\bf 28}, no.\,10 (2000), 4767--4782). \emp

\bmp \textbf{12.76.} Is every group generated by a set of
3-transpositions locally finite? A {\it set of $3$-transpositions\/}
is, by definition, a normal set of involutions such that the orders
of their pairwise products are at most 3. \hf {\sl
A.\,I.\,Sozutov}\vs

 Yes, it is
(H.\,Cuypers, J.\,I.\,Hall, in: {\it Groups, combinatorics and
geometry. Proc. L.\,M.\,S. Durham symp., July 5--15, 1990 $($London
Math. Soc. Lecture Note Ser.}, {\bf 165}), Cambridge Univ. Press,
1992, 121--138).

\emp

\bmp \textbf{12.77.} (Well-known problem). Does the order (if it is
greater than $p^2$) of a finite non-cyclic $p\hskip0.2ex $-group
divide the order of its automorphism group? \hfill
{\sl A.\,I.\,Starostin}
\vs

No, not always (J.\,Gonz\'alez-S\'anchez, A.\,Jaikin-Zapirain, {\it Forum Math. Sigma}, {\bf 3}, Article ID e7, 11 p., electronic only (2015)).
\emp

\bmp \textbf{12.78.} (M.\,J.\,Curran). a) Does there exist a group of
order $p^6$ (here $p$ is a prime number), whose automorphism group
has also order $p^6$?

 \mb{b)} Is it true that for $p \equiv 1\, ({\rm
mod}\, 3)$, the smallest order of a $p\hskip0.1ex $-group that is
the automorphism group of a $p\hskip0.1ex $-group is $p^7$?

 \mb{c)} The same question for $p=3$ with $3^9$
replacing $p^7$. \hf
{\sl A.\,I.\,Starostin}\vs

 a) No, there is no such a
group; b)~yes, it is true; c)~no, it is not true (S.\,Yu, G.\,Ban,
J.\,Zhang, {\it Algebra Colloq.}, {\bf 3}, no.\,2 (1996),
97--106).
\emp

\bmp \textbf{12.79.} Suppose that $a$ and $b$ are two elements of a
finite group $G$ such that the function 
$\varphi
(g)=1^G(g) - 1_{\left< a\right>}^G(g) - 1_{\left< b\right>}^G(g)
- 1_{\left< ab\right>}^G(g) + 2$ 
 is a character of $G$. Is it
true that $G=\left< a,b\right>$? The converse statement is true.
\hfill {\sl S.\,P.\,Strunkov}

\vs
 No, not always: a counterexample is given by $G=A_4$, $a=b=(123)$. (S.\,V.\,Skresanov, {\it Algebra Logic}, {\bf 58}, no.\,3 (2019),
249--253).
\emp

\bmp \textbf{12.80.} (K.\,W.\,Roggenkamp). a)
 Is it true that the number of
$p\hskip0.1ex $-blocks of defect~0 of a finite group $G$ is equal
to the number of the conjugacy classes of elements $g \in G$ such
that the number of solutions of the equation $[x,y]=g$ in~$G$ is
not divisible by~$p$?

 \mb{b)} The same question in the case of $G$ being a simple group.
 \hf {\sl S.\,P.\,Strunkov}\vs

a) No it is not true; for example, if $p=2$ and $G={\Bbb Z} _3\times
{\Bbb S} _3$ (L.\,Barker, {\it Letter of May, 27, 1996}).\vs

b) No, not always. For example in the alternating group $A_5$ for $p=2$, the number of 2-blocks of defect 0 is 1, but there are 3 conjugacy classes of elements $g$ such that the number of solutions $[x,y] = g$ is odd, namley, $g = (1,2,3)$, $(1,2,3,4,5)$, and $(1,3,5,2,4)$. (B.\,Sambale, \textit{Letter of 5 May 2021}.)
\emp

\bmp \textbf{12.82.} Find all pairs
$(n,r)$ such that the symmetric group ${\Bbb S}_n$ contains a maximal
subgroup isomorphic to ${\Bbb S}_r$. \hfill {\sl
V.\,I.\,Sushchanski\u{\i}}

\vs

They are found 
(B.\,Newton, B.\,Benesh, {\it J.~Algebra}, \textbf{304}, no.\,2
(2006), 1108--1113).

\emp

\bmp \textbf{12.84.} (Well-known problem). Is it true that if there
exist two non-iso\-mor\-phic groups with the given set of orders
of the elements, then there are infinitely many groups with this
set of orders of the elements? \hf {\sl S.\,A.\,Syskin}\vs

 No, it is not
(V.\,D.\,Mazurov, {\it Algebra and Logic}, {\bf 33}, no.\,1
(1994), 49--55).

\emp

\bmp \textbf{12.90.} Let $G$ be a finitely generated soluble minimax
group and let $H$ be a finitely generated residually finite group
which has precisely the same finite images as $G$. Must $H$ be a
minimax group? \hf {\sl John S.\,Wilson}\vs

 Yes, it must
 (A.\,Mann, D.\,Segal, {\it Proc.
London Math. Soc. (3)}, {\bf 61} (1990), 529--545; see also
A.\,V.\,Tushev, {\it Math. Notes}, {\bf 56}, no.\,5 (1994),
1190--1192).

\emp

\bmp \textbf{12.91.} Every metabelian group
belonging to a Fitting class of (finite) supersoluble groups is
nilpotent. Does the following generalization also hold: Every
group belonging to a Fitting class of supersoluble groups has only
central minimal normal subgroups? \hf
{\sl H.\,Heineken}\vs

 No, it does not
(H.\,Heineken, {\it Rend. Sem. Mat. Univ. Padova}, {\bf 98}
(1997), 241--251).

\emp

\bmp \textbf{12.93.} Let $N \rightarrowtail G \twoheadrightarrow Q$
be an extension of nilpotent groups, with $Q$ finitely generated,
which splits at every prime. Does the extension split? This is
known to be true if $N$ is finite or commutative. \hf {\sl
P.\,Hilton}\vs

 No, not always (K.\,Lorensen,
{\it Math. Proc. Cambridge Philos. Soc.}, {\bf 123}, no.\,2
(1998), 213--215). \emp

\bmp \textbf{12.94.}
Let $G$ be a finitely generated
pro-$p\hskip0.2ex $-group not involving the wreath product $C_p\wr
{\Bbb Z}_p$ as a closed section (where $C_p$ is a cyclic group of
order $p$ and ${\Bbb Z}_p$ is the group of $p\hskip0.2ex $-adic
integers). Does it follow that $G$ is $p\hskip0.2ex $-adic
analytic? \hfill {\sl A.\,Shalev}

\vs  No, it does not (A.\,Jaikin-Zapirain, B.\,Klopsch, \emph{J.~London Math. Soc. (2)}, \textbf{76}, no.\,2 (2007), 365--383).
\emp

\bmp \textbf{12.96.}
Find a
non-empty Fitting class ${\frak F}$ and a non-soluble finite group
$G$ such that $G$ has no ${\frak F} $-injectors. \hfill
 {\sl L.\,A.\,Shemetkov}

\vs

Found by E.\,Salomon
(Mainz Univ., unpublished); his example is presented in \S\,7.1 of
(A.\,Ballester-Bolinches, L.\,M.\,Ezquerro, {\it Classes of finite
groups}, Springer, 2006).

\emp

\bmp \textbf{12.97.} Let ${\frak F}$ be the formation of all finite
groups all of whose composition factors are isomorphic to some
fixed simple non-abelian group $T$. Prove that ${\frak F}$ is
indecomposable into a product of two non-trivial subformations.
\hf {\sl L.\,A.\,Shemetkov}\vs

 It is proved (O.\,V.\,Mel'nikov,
 {\it Problems in Algebra}, {\bf 9}, Gomel', 1996, 42--47 (Russian)).

\emp

\bmp \textbf{12.98.} Let $F$ be a free group of finite rank, $R$ its
recursively defined normal sub\-group. Is it true that

 \mb{a)} the word problem for
$F/R$ is soluble if and only if it is soluble for $F/[R,R]$?

 \mb{b)} the conjugacy problem for
$F/R$ is soluble if and only if it is soluble for $F/[R,R]$?

 \mb{c)} the conjugacy problem
for $F/[R,R]$ is soluble if the word problem is soluble for
$F/[R,R]$? \hf {\sl V.\,E.\,Shpil'rain}\vs

 a) Yes, it is; b)
No, it is not; c) No, it is not (M.\,I.\,Anokhin, {\it Math.
Notes}, {\bf 61}, no.\,1--2 (1997), 3--8).

\emp

\bmp

 {\bf 12.102.} Is
every proper factor-group of a group of Golod (see 9.76) residually
finite?

\hfill {\sl V.\,P.\,Shunkov}\vsh

No, not every (L.\,Hammoudi, {\it Algebra Colloq.}, {\bf 5} (1998),
371--376).

\emp

\markboth{\protect\vphantom{(y}{Archive of solved
problems (13th ed., 1995)}}{\protect\vphantom{(y}{Archive of solved
problems (13th ed., 1995)}}

\bmp \textbf{13.10.} Is there a function
$f\! : {\Bbb N} \rightarrow {\Bbb N}$ such that, for every soluble
group $G$ of derived length $k$ generated by a set $A$, the
validity of the identity
 $x^4=1$ on each subgroup generated by at most $f(k)$
elements of $A$ implies that $G$ is a group of exponent 4?

\hfill {\sl V.\,V.\,Bludov}

No, there is not (G.\,S.\,Deryabina, A.\,N.\,Krasil'nikov, {\it
Siberian Math. J.}, {\bf 44}, no.\,1 (2003), 58--60).

 \emp

\bmp \textbf{ 13.11.}
Is a torsion-free
group almost polycyclic if it has a finite set of generators
$a_1,\ldots , a_n$ such that every element of the group has a
unique presentation in the form $a_1^{k_1}\cdots\, a_n^{k_n}$,
where $k_1,\ldots , k_n\in {\Bbb Z}$?\hfill {\sl
V.\,V.\,Bludov}

\vs

No, not always
(A.\,Muranov, {\it Trans. Amer. Math. Soc.}, \textbf{359} (2007),
3609--3645).

 \emp

\bmp \textbf{13.16.} Is every locally nilpotent group with minimum
condition on centralizers hyper\-central?\hfill {\sl
F.\,O\,Wagner}\vs

 Yes, it is (V.\,V.\,Bludov,
{\it Algebra and Logic}, {\bf 37}, no.\,3 (1998), 151--156).
 \emp

\bmp \textbf{13.18.}
Let $F$ be a finitely generated non-Abelian free group and let $G$
be the Cartesian (unrestricted) product of countable infinity of
copies of $F$. Must the Abelianization $G/G'$ of $G$ be
torsion-free?\hfill {\sl A.\,M.\,Gaglione,~D.\,Spellman}\vs

No, it may contain elements of order~$2$ \ (O.\,Kharlampovich,
A.\,Myasnikov, in: {\it Knots, braids, and mapping class
groups\,---\,papers dedicated to Joan S.~Birman (New York, 1998)},
Amer. Math. Soc., Providence, RI, 2001, 77--83).
 \emp

\bmp \textbf{13.21.} a) Is there an infinite finitely generated
residually finite
 $p\hskip0.1ex $-group, in which the order $|g|$ of an arbitrary
element $g$ does not exceed $f(\delta (g))$, where $\delta (g)$ is
the length of $g$ with respect to a fixed set of generators and
 $f(n)$ is a function growing at $n\rightarrow \infty$ slower than any
power function $n^{\lambda}$, $\lambda >0$? \hfill {\sl
R.\,I.\,Grigorchuk}\vs

Yes, there is (A.\,V.\,Rozhkov, {\it Dokl. Math.}, {\bf 58},
no.\,2 (1998), 234--237).
 \emp

\bmp \textbf{13.26.}
Is it true that a countable topological group of
exponent 2 with unique free ultrafilter converging to the identity
has a basis of neighborhoods of the identity consisting of
subgroups?\hfill \raisebox{0ex}{\sl E.\,G.\,Zelenyuk,
I.\,V.\,Protasov}

\vs

Assuming Martin's Axiom, the answer is ``No'' (Ye.\,Zelenyuk, {\it
Adv. Math.}, {\bf 229}, no.\,4 (2012), 2415--2426).

 \emp

\bmp \textbf{13.28.}
(D.\,M.\,Evans). A permutation group on an
infinite set is {\it cofinitary\/} if its non-iden\-tity elements
fix only finitely many points. Is it true that a closed cofinitary
per\-mu\-tation group is locally compact (in the topology of
pointwise convergence)?

\hfill {\sl P.\,J.\,Cameron}\vsh

No, not always (G.\,Hjorth, {\it J.~Algebra}, {\bf 200}, no.\,2
(1998), 439--448).
 \emp

\bmp \textbf{13.29.}
Given an infinite set $\Omega$, define an
algebra $A$ (the {\it reduced incidence algebra of finite
subsets\/}) as follows. Let $V_n$ be the set of functions from the
set of $n $-element subsets of $\Omega$ to the rationals $\Bbb Q$.
Now let $A=\bigoplus V_n$, with multiplication as follows: for
$f\in V_n$, $g\in V_m$, and $|X|=m+n$, let $(fg)(X)=\sum
f(Y)g(X\setminus Y),$ where the sum is over the $n $-element
subsets $Y$ of $X$. If $G$ is a permutation group on $\Omega$, let
$A^G$ be the algebra of $G $-inva\-ri\-ants in $A$.

\makebox[15pt][r]{}
{\it Conjecture:\/} If $G$ has no finite orbits
on $\Omega$, then $A^G$ is an integral domain.

\hfill {\sl P.\,J.\,Cameron}
\vs

Conjecture is proved (M.\,Pouzet, \textit{Theor. Inf. App.}, \textbf{42}, no.\,1 (2008), 83--103).
 \emp

\bmp \textbf{13.33.} (F.\,Gross).
Is a normal subgroup of a finite $D_{\pi} $-group (see Archive, 3.62)
always a $D_{\pi} $-group? \hfill {\sl
V.\,D.\,Mazurov}
\vs

Yes, it is mod CFSG
 (E.\,P.\,Vdovin, D.\,O.\,Revin, {\it in: Ischia Group Theory 2004 $($Contemp.
Math.}, \textbf{402}), Amer. Math. Soc., 2006, 229--263).
 \emp

\bmp \textbf{13.34.} (I.\,D.\,Macdonald). If the identity $[x,\,
y]^n=1$ holds on a group, is the derived subgroup of the group
periodic? \hfill {\sl V.\,D.\,Mazurov}\vs

No, not always (G.\,S.\,Deryabina, P.\,A.\,Kozhevnikov, {\it
Commun. Algebra}, {\bf 27}, no.\,9 (1999), 4525--4530; \
S.\,I.\,Adyan, {\it Dokl.
 Math.}, {\bf 62}, no.\,2 (2000), 174--176).

 \emp

\bmp \textbf{13.36.} b) (A.\,Lubotzky,
 A.\,Shalev). For a finitely generated pro-$p\hskip0.2ex
$-group $G$ set $a_n(G)={\rm dim}_{\,{\Bbb F}_p}I^n/I^{n+1}$,
where $I$ is the augmentation ideal of the group ring
 ${\Bbb F}_p[[G]]$.
We define the {\it growth\/} of $G$ to be the growth of the
sequence
 $\{ a_n(G)\} _{n\in {\Bbb N}}$.
 Is the growth of $G$
exponential if $G$ contains a finitely generated closed
subgroup of exponential growth?
\hfill {\sl O.\,V.\,Mel'nikov}

\vs

b) Not necessarily. Let $G$ be the Nottingham
group over ${\Bbb F}_p$. The growth of $c_n(G):=\log_p |G:\omega_n(G)|$,
where $\{\omega_n(G)\}$ is the Zassenhaus filtration, is linear,
so by Quillen's theorem the growth of $G$ is
subexponential. But $G$ has a 2-generator free pro-$p$ subgroup
(R.\,Camina, {\it J.~Algebra}, \textbf{ 196} (1997), 101--113) with
exponential growth. (M.\,Ershov, {\it Letter of 19.10.2009}.)

 \emp

\bmp \textbf{13.38.} Let $G=F/R$ be a pro-$p\hskip0.1ex $-group with one
defining relation, where
 $R$ is the normal subgroup of a free pro-$p\hskip0.1ex $-group
$F$ generated by a single element $r\in F^p[F,\, F]$.

\makebox[25pt][r]{a)} Suppose that $r=t^p$ for some $t\in F$; can $G$
contain a Demushkin group as a subgroup?

\makebox[25pt][r]{b)} Do there exist two pro-$p\hskip0.1ex $-groups
$G_1\supset G_2$ with one defining relation, where $G_1$ has elements
of finite order, while the subgroup $G_2$ is torsion-free?\hfill {\sl
O.\,V.\,Mel'nikov}\vs

a) Yes, it can. b) Yes,
 they exist. The pro-$p\hs$-group with presentation
$\left< a,b\mid [a,b]^p=1\right>$ contains the Demushkin group
$\left< x_1,\ldots
 ,x_{2g}\mid \prod_{i=1}^{g}[x_{2i-1},x_{2i}]\,=1\right>$ for some
 $g>1$. (O.\,V.\,Mel'nikov, {\it Letter of February, 28, 1999\/}).
 \emp

\bmp
\textbf{13.39.} Let $A$ be an associative ring with unity and with
torsion-free additive group, and let $F^A$ be the tensor product
of a free group $F$ by $A$ (A.\,G.\,Myasnikov,
V.\,N.\,Re\-mes\-len\-ni\-kov, {\it Siberian Math.~J.}, \textbf{35},
no.\,5 (1994), 986--996);
 then $F^A$ is a {\it free exponential group over $A$}; in
(A.\,G.\,Myasnikov, V.\,N.\,Remeslennikov, {\it
Int. J.~Algebra Comput.}, \textbf{6} (1996), 687--711), it is shown
how to construct
 $F^A$ in terms of free products with amalgamation.

\makebox[25pt][r]{a)}
(G.\,Baumslag). Is $F^{A}$ residually
nilpotent torsion-free?

\makebox[25pt][r]{c)}
(G.\,Baumslag). Is the Magnus homomorphism
of $F^{\Bbb Q}$ into the group of power series over the rational
number field $\Bbb Q$ faithful or not?

\hfill {\sl
A.\,G.\,Myasnikov,
V.\,N.\,Remeslennikov}

\vs

a) Yes, it is (A.\,Jaikin-Zapirain, \emph{Ann. Sci. \'Ec. Norm. Sup\'er. (4)}, {\bf 57}, no.\,4 (2024), 1101--1133).
\vs
 c) Yes, it is faithful (A.\,Jaikin-Zapirain, \emph{Ann. Sci. \'Ec. Norm. Sup\'er. (4)}, {\bf 57}, no.\,4 (2024), 1101--1133).
\emp

\bmp \textbf{13.40.} A group $G$ is said to be $\omega$-{\it
resi\-du\-al\-ly free\/} if, for every finite set of non-trivial
elements of $G$, there is a homomorphism of~$G$ into a free group
such that the images of all these elements are non-trivial. Is
every finitely-generated $\omega$-resi\-du\-al\-ly free group
embeddable in a free ${\Bbb Z} [x]$-group? \hfill {\sl
A.\,G.\,Myasnikov,~V.\,N.\,Remeslennikov}\vs

 Yes, it is (O.\,Kharlampovich,
A.\,Myasnikov, {\it J.\,Algebra}, {\bf 200} (1998), 472--570).

 \emp

\bmp \textbf{13.45.}
Every infinite group $G$ of regular cardinality
 ${\frak m}$ can be partitioned into two subsets
$G=A_1\cup A_2$ so that
 $A_1F\ne G$ and $A_2F\ne G$ for every subset
$F\subset
G$ of cardinality less than
 ${\frak m}$. Is this statement true for groups of singular
cardinality?

\hfill {\sl I.\,V.\,Protasov}
 \vs

 No, not always. The answer depends on the algebraic structure of $G$. In particular, this is true for a free group, but
the statement does not hold for every Abelian group $G$ of singular cardinality (I.\,Protasov, S.\,Slobodianiuk, {\it Quest. Answers Gen. Topology}, {\bf 33}, no.\,2 (2015), 61--70).
\emp

\bmp \textbf{13.46.} Can every uncountable abelian group of finite odd
exponent be partitioned into two subsets so that neither of them
contains cosets of infinite subgroups? Among countable abelian
groups, such partitions exist for groups with finitely many
involutions. \hfill {\sl I.\,V.\,Protasov}\vs

Yes, it can (E.\,G.\,Zelenyuk, {\it Math. Notes}, {\bf 67} (2000),
599--602).
 \emp

\bmp \textbf{13.47.} Can every countable abelian group with finitely
many involutions be partitioned into two subsets that are dense in
every group topology? \hfill {\sl I.\,V.\,Protasov}\vs

 Yes, it can
(E.\,G.\,Zelenyuk, {\it Ukrain. Math.~J.}, {\bf 51}, no.\,1
(1999), 44--50).
 \emp

\bmp
\textbf{13.50.}
Let ${\frak F}$ be a local Fitting class. Is it true
that there are no maximal elements in the partially ordered by
inclusion set of the Fitting classes contained in
 ${\frak F}$ and distinct from
 ${\frak F}$? \hfill {\sl A.\,N.\,Skiba}

\vs

No, there may exist
maximal elements (N.\,V.\,Savel'eva, N.\,T.\,Vorob'ev, \textit{Siberian Math. J.}, {\bf 49}, no.\,6
(2008), 1124--1130). 
 \emp

\bmp \textbf{13.54.} b) Is it true that every group is embeddable in
the kernel of some Frobenius group (see Archive, 6.53)?
 \hf {\sl A.\,I.\,Sozutov}\vs

 Yes, it is true (V.\,V.\,Bludov, {\it Siberian Math.~J.}, {\bf 38},
no.\,6 (1997), 1054--1056). \emp

\bmp \textbf{13.56.}
(A.\,Shalev). Let $G$ be a finite $p\hskip0.1ex $-group of
sectional rank $r$, and $\f$ an auto\-mor\-phism of~$G$ having
exactly $m$ fixed points. Is the derived length of $G$ bounded by
a function depending on $r$ and $m$ only?\hfill {\sl
E.\,I.\,Khukhro}\vs

Yes, it is (A.\,Jaikin-Zapirain, {\it Israel J.~Math.}, {\bf 129}
(2002), 209--220).

 \emp

\bmp \textbf{13.58.}
Let $\f$ be
an automorphism of prime order
 $p$ of a nilpotent (periodic) group~$G$ such that
$C_G(\f )$ is a group of finite sectional rank
 $r$.
Does $G$ possess a normal subgroup $N$ which is nilpotent of class
bounded by a function of $p$ only and is such that
 $G/N$ is a group of finite sectional rank bounded in terms
of
 $r$ and $p$?

\makebox[15pt][r]{}This was proved for $p=2$ in (P.\,Shumyatsky, {\it
Arch. Math.}, \textbf{71} (1998), 425--432).

\hfill {\sl E.\,I.\,Khukhro}

\vs

Yes, it does
(E.\,I.\,Khukhro, \textit{J.~London Math. Soc.}, \textbf{77} (2008),
130--148).

 \emp

\bmp \textbf{13.61.} We call a metric space {\it narrow\/} if it is
quasiisometric to a subset of the real line, and {\it wide\/}
otherwise. Let $G$ be a group with the finite set of generators
$A$, and let $\Gamma =\Gamma (G,A)$ be the Cayley graph of $G$
with the natural metric. Suppose that, after deleting any narrow
subset $L$ from $\Gamma$, at most two connected components of the
graph $\Gamma \setminus L$ can be wide, and there exists at least
one such a subset $L$ yielding exactly two wide components in
$\Gamma \setminus L$. Is it true that $\Gamma$ is quasiisometric
to an Euclidean or a hyperbolic plane? \hfill {\sl
V.\,A.\,Churkin}\vs

No, it is not true in general (O.\,V.\,Bogopol'skii, {\it
Preprint}, Novosibirsk, 1998 (Russian)). See also new problem 14.98.
 \emp

\bmp \textbf{13.62.} Let $U$ and $V$ be non-cyclic subgroups of a
free group. Does the inclusion
 $[U,\,
U]\leq [V,\, V]$ imply that $U\leq V$? \hf {\sl
V.\,P.\,Shaptala}\vs

 No, it does not
(M.\,J.\,Dunwoody, {\it Arch. Math.}, {\bf 16} (1965), 153--157).
 \emp

\bmp \textbf{13.63.}
Let $\pi _e(G)$ denote the set of orders of elements of a group
$G$. For $\Gamma \subseteq {\Bbb N}$ let $h(\Gamma )$ denote the
number of non-isomorphic finite groups $G$ with $\pi
_e(G)=\Gamma$. Is there a number $k$ such that, for every
$\Gamma$, either $h(\Gamma )\leq k$, or $h(\Gamma )=\infty$?
\hfill {\sl W.\,J.\,Shi}\vs

No, there is no such number: $h(\pi _e(L_3 (7^{3^r})))=r+1$ for any $r\geq 0$
(A.\,V.\,Zavarnitsine, {\it J.~Group Theory}, {\bf 7}, no.\,1
(2004), 81--97).\emp

\bmp \textbf{13.66.} Let $F$ be a \

\mb{a)} free; \

\mb{b)} free metabelian \

group of finite rank. Let $M$ denote the set of all
endo\-mor\-phisms of
 $F$ with non-cyclic images.
Can one choose two elements $g,\, h\in F$ such that, for every $\f
,\, \psi \in M$, equalities $\f (g)=\psi (g)$ and $\f (h)=\psi (h)$
imply that $\f =\psi$, that is, the endomorphisms in
 $M$ are uniquely determined by their values at
 $g$ and $h$?
\hf {\sl V.\,E.\,Shpil'rain}\vs

a) Yes, one can (D.\,Lee, {\it J.~Algebra}, {\bf 247} (2002),
no.\,2, 509--540).

 b) No, not always (E.\,I.\,Timoshenko, {\it Math. Notes},
{\bf 62}, no.\,5--6 (1998), 767--770).
 \emp

\bmp \textbf{14.1.}
Suppose that $G$ is a finite group with no
non-trivial normal subgroups of odd order, and $\varphi$ is its
$2$-automorphism centralizing a Sylow $2$-sub\-group of $G$. Is it
true that $\varphi ^2$ is an inner automorphism of $G$?
 \hfill {\sl R.\,Zh.\,Aleev}

\vs Yes, it is true (G.\,Glauberman, 
{\it Math.~Z.}, {\bf 107} (1968), 1--20).
\emp

\bmp
\textbf{14.10.}
a)
(Well-known problem). It is known that any
recursively presented group embeds in a finitely presented group
(G.\,Higman, {\it Proc. Royal Soc. London Ser. A}, \textbf{262} (1961),
455--475).
 Find an explicit and {\it \lq\lq
natural\rq\rq} \ finitely presented group $\Gamma$ and an
embedding of the additive group of the rationals ${\Bbb Q}$ in
$\Gamma$.

\makebox[15pt][r]{}b)
Find an explicit embedding of $\Bbb Q$ in
a finitely generated group; such a group exists by Theorem~IV in
(G.\,Higman, B.\,H.\,Neumann, H.\,Neumann, {\it J.\,\,London Math.
Soc.}, \textbf{24} (1949), 247--254). \hfill {\sl
P.\,\,de la Harpe}

\vs

a) Such an embedding is found  (J.\,Belk, J.\,Hyde, F.\,Matucci, {\it Bull. Amer. Math. Soc.}, {\bf 59}, no.\,4 (2022), 561--567; \ V.\,H.\,Mikaelian,  {\it  J.~Group Theory}, {\bf 29}, no.\,2 (2026),  371--417).
\vs

b) Such an embedding is found
(V.\,H.\,Mikaelian, {\it Int. J. Math. Math. Sci.}, \textbf{
2005}, no.\,13 (2005), 2119--2123); another embedding can be derived from the proof in (P.\,Hall, {\it  Proc. London Math. Soc. (3)}, {\bf 9}  (1959), 595--622).
\emp

\markboth{\protect\vphantom{(y}{Archive of solved
problems (14th ed., 1999)}}{\protect\vphantom{(y}{Archive of solved
problems (14th ed., 1999)}}

\bmp \textbf{14.13.} a)
By definition the {\it commutator length\/} of
an element $z$ of the derived subgroup of a group $G$ is the least
possible number of commutators from $G$ whose product is equal to
$z$. Does
there exist a simple group on which the commutator length is not
bounded?

\mb{b)}
Does there exist a finitely presented simple
group on which the commutator length is not bounded?
\hfill
 {\sl V.\,G.\,Bardakov}

\vs

a) Simple groups with this
property were constructed in (J.\,Barge, \'{E}.\,Ghys, {\it Math.
Ann.}, \textbf{294} (1992), 235--265), and finitely generated simple ones
in (A.\,Muranov, {\it Int. J. Algebra Comput.}, \textbf{17} (2007),
607--659).\vs

b) Yes, there does (P.-E.\,Caprace, K.\,Fujiwara, {\it Geom. Funct. Anal.}, {\bf 19} (2010), 1296--1319).
 \emp

\bmp \textbf{14.25.} Let $qG$ denote the quasivariety generated by a
group~$G$. Is it true that there exists a finitely generated group
$G$ such that the set of proper maximal subquasivarieties in $qG$
is infinite? \hfill {\sl A.\,I.\,Budkin}\vs

Yes, it is (V.\,V.\,Bludov, {\it Algebra and Logic}, {\bf 41},
no.\,1 (2002), 1--8).

 \emp

 \bmp \textbf{14.27.}
Let $\Gamma$ be a group generated by a finite set $S$. Assume that
there exists a nested sequence $F_1\subset F_2\subset\cdots$ of
finite subsets of $\Gamma$ such that \ (i)\,\, $F_k\ne F_{k+1}$
for all $k\geq 1$, \ \ (ii)\,\, $\displaystyle{\Gamma=\bigcup_{k
\geq 1}F_k}$, \ \ (iii)\,\,
$\displaystyle{\lim_{k\to\infty}\frac{\vert\partial F_k\vert}
{\vert F_k\vert}=0}$, where, by definition, \linebreak $\partial
F_k = \{\,\gamma\in\Gamma \setminus F_k \mid \text{there exists}\;
s\in S \text{ such that}\;\gamma s\in F_k\,\}$, \ and \ (iv)\,\,
there exist constants $c\geq 0,\,\,d\geq 1$ such that $\vert
F_k\vert\le c k^d$ for all $k\geq 1$. Does it follow that $\Gamma$
has polynomial growth? \hfill {\sl A.\,G.\,Vaillant}\vs

No, not always. These properties are enjoyed by every group of
intermediate growth (V.\,G.\,Bardakov, {\it Algebra and Logic},
{\bf 40}, no.\,1 (2001), 12--16). 
 \emp

\bmp \textbf{14.32.}
Extending the classical
definition of formations, let us define a {\it formation\/} of
(not necessarily finite) groups as a nonempty class of groups
closed under taking homomorphic images and subdirect products with
finitely many factors. Must every first order axiomatizable
formation of groups be a variety?
\hf {\sl \mbox{A.\,M.\,Gaglione, D.\,Spellman}}

\vs
No, every variety of
groups that contains a finite nonsolvable member contains an
axiomatic subformation that is not a variety (K.\,A.\,Kearnes,
{\it J.~Group Theory}, {\bf 13}, No. 2 (2010), 233--241).
 \emp

\bmp \textbf{14.33.}
Does there exist a finitely presented
pro-$p\hs$-group ($p$ being a prime) which contains an isomorphic
copy of every countably based pro-$p\hs$-group?

\hfill {\sl R.\,I.\,Grigorchuk}

\vs

Yes, there does:
the Nottingham group over ${\Bbb F} _p$ for $p>2$, which
was shown to be finitely presented in
(M.\,V.\,Ershov, {\it J. London Math. Soc.}, \textbf{71} (2005),
362--378) (M.\,Ershov, {\it Letter
of 19.10.2009}).

 \emp

\bmp \textbf{14.34.}
By definition, a locally compact group has the
{\it Kazhdan $T$-property\/} if the trivial representation is an
isolated point in the natural topological space of unitary
representations of the group. Does there exist a profinite group
with two dense discrete subgroups one of which is amenable, and
the other has the Kazhdan $T$-property?

\makebox[15pt][r]{}See (A.\,Lubotzky, {\it
Discrete groups, expanding graphs and invariant measures\/} ({\it
Progress in Mathematics, Boston, Mass.}, \textbf{125}),
Birkh\"{a}user, Basel, 1994) for further motivation.
\hfill {\sl
R.\,I.\,Grigorchuk,~A.\,Lubotzky}

\vs

Yes, there does
(M.\,Kassabov, {\it Invent. Math.}, \textbf{170} (2007),
297--326; \ M.\,Ershov, A.\,Jaikin-Zapirain, {\it
Invent. Math.}, \textbf{179}, 303--347).

 \emp

\bmp \textbf{14.48.} If an equation over a free group
 $F$ has no solution in
$F$, is there a finite quotient of $F$ in which the equation has no
solution? \hfill {\sl L.\,Comerford}\vs

No, not always (T.\,Coulbois, A.\,Khelif, {\it Proc. Amer. Math.
Soc.}, {\bf 127}, no.\,4 (1999), 963--965).
 \emp

\bmp \textbf{14.49.}
 Is $SL_3({\Bbb Z} )$ a factor group of the modular group $
PSL_2({\Bbb Z} )$? Since the latter is isomorphic to the free product
of two cyclic groups of orders 2 and 3,
 the
question asks if $SL_3({\Bbb Z} )$ can be generated by two
elements of orders $2$ and $3$. \hfill {\sl M.\,Conder}\vs

No, it is not (M.\,C.\,Tamburini, P.\,Zucca, {\it Atti Accad. Naz.
Lincei Cl. Sci. Fis. Mat. Natur. Rend. Lincei (9) Mat. Appl.},
{\bf 11}, no.\,1 (2000), 5--7; \ Ya.\,N.\,Nuzhin, {\it Math.
Notes}, {\bf
70}, no.\,1--2 (2001), 71--78). M.\,Conder has shown, however,
that $SL_3({\Bbb Z})$ has a subgroup of index 57 that is a factor
group of the modular group $PSL_2({\Bbb Z})$. Also it has been
shown in (M.\,C.\,Tamburini, J.\,S.\,Wilson, N.\,Gavioli, {\it
J.~Algebra}, {\bf 168} (1994), 353--370) that $SL_d({\Bbb Z})$ is
a factor group of the modular group $PSL_2({\Bbb Z})$ for all $d
\geq 28$.
 \emp

\bmp \textbf{14.50.} (Z.\,I.\,Borevich). A subgroup $A$ of a group $G$
is said to be {\it paranormal\/} (respectively, {\it polynormal\/})
if $A^x\leq \left< A^u \mid u\in \left< A, A^x\right> \right> $
(respectively, $A^x\leq \left< A^u \mid u\in A^{\left< x
\right>}\right> $)
 for any
 $x\in G$.
Is every polynormal subgroup of a finite group paranormal?

\hfill {\sl A.\,S.\,Kondratiev}\vsh

Not always (V.\,I.\,Mysovskikh, {\it Dokl. Math.}, {\bf 60} (1999),
71--72).
 \emp

\bmp \textbf{14.52.} It is known that if a finitely generated group is
residually torsion-free nilpotent, then the group is residually
finite $p\hs$-group, for every prime $p$. Is the converse true?

 \hfill {\sl Yu.\,V.\,Kuz'min}\vs

Not always (B.\,Hartley, in: {\it Symposia Mathematica, Bologna,
Vol.~XVII, Convegno sui Gruppi Infiniti, INDAM, Roma, 1973}, Academic
Press, London, 1976, 225--234).
 \emp

\bmp
\textbf{14.55.} a)
 Prove that the Nottingham group $J=N({\Bbb Z}/p{\Bbb Z})$ (as
defined in Archive, 12.24) is finitely presented
 for $p>2$. \hfill {\sl
C.\,R.\,Leedham-Green}

\vs

 Proved for $p>2$
(M.\,V.\,Ershov, {\it J. London Math. Soc.}, \textbf{71} (2005),
362--378).

 \emp

\bmp \textbf{14.58.} b) Suppose that $A$ is a periodic group of regular
automorphisms of an abelian group.
 Is $A$ finite if $A$ is generated by elements of order $3$? \hfill
{\sl V.\,D.\,Mazurov}\vs

Yes, it is (A.\,Kh.\,Zhurtov, {\it Siberian Math. J.}, {\bf 41},
no.\,2 (2000), 268--275).
 \emp

\bmp \textbf{14.60.}
Suppose that $H$ is a non-trivial normal
subgroup of a finite group~$G$ such that the factor-group
 $G/H$ is isomorphic to one of the simple groups $L_n(q)$, $n\geq 3$.
Is it true that $G$ has an element whose order is distinct from
the order of any element in $G/H$? \hfill {\sl
V.\,D.\,Mazurov}
\vs

Yes, it is true: for $n\ne 4$ proved in (A.\,V.\,Zavarnitsine, {\it Siberian Math.~J.}, \textbf{49} (2008),
246--256), and for $n=4$ in (M.\,A.\,Grechkoseeva, S.\,V.\,Skresanov, {\it Siberian Electron. Math. Rep.}, {\bf 17} (2020), 585--589). {\it Editors' comment}: previous claim that it is not true for $n=4$ was erroneous. %
 \emp

\bmp \textbf{14.62.}
Suppose that $H$ is a non-soluble normal subgroup of a finite
group~$G$. Does there always exist a maximal soluble subgroup $S$
of $H$ such that $G=H\cdot N_G(S)$?

 \hfill {\sl V.\,S.\,Monakhov}\vs

Yes, it always exists mod CFSG (V.\,I.\,Zenkov, V.\,S.\,Monakhov,
D.\,O.\,Revin, {\it Algebra and Logic}, {\bf 43}, no.\,2 (2004),
102--108). \emp

\bmp \textbf{14.63.}
What are the
composition factors of non-soluble finite groups all of whose
normalizers of Sylow subgroups are $2$-nilpotent, in particular,
supersoluble?

\hfill {\sl V.\,S.\,Monakhov}

 \vs

Described
(L.\,S.\,Kazarin, A.\,A.\,Volochkov, {\it Mathematics in
Yaroslavl' Univ.: Coll. of Surveys to 30-th Anniv. of Math. Fac.},
Yaroslavl', 2006, 243--255 (Russian)).

 \emp

\bmp \textbf{14.66.} (Well-known problem). Let $G$ be a finite soluble
group, $\pi(G)$ the set of primes dividing the order of $G$, and
$\nu(G)$ the maximum number of primes dividing the order of some
element. Does there exist a linear bound for $|\pi(G)|$ in terms
of~$\nu(G)$?

\hfill {\sl A.\,Moret\'o}\vsh

Yes, it does (T.\,M.\,Keller, {\it J.~Algebra}, {\bf 178}, no.\,2
(1995),
 643--652).
 \emp

\bmp \textbf{14.71.} Consider a free group $F$ of finite rank and an
arbitrary group $G$.
 Define the {\it $G$-closure\/} ${\rm cl}_G(T)$ of any subset
$T\subseteq F$ as the intersection of the kernels of all those
homomorphisms $\mu:F \rightarrow G$ of $F$ into $G$ that vanish on
$T$: \ $ {\rm cl}_G(T) = \bigcap \,\{ {\rm Ker}\,\mu \,\mid \,\mu :F
\rightarrow G;\;\, T\subseteq {\rm Ker}\,\mu \}$. \
 Groups $G$ and $H$ are
called {\it geometrically equivalent\/} if for every free group $F$
and every subset $T\subset F$ the $G$- and $H$-closures of $T$
coincide: ${\rm cl}_G(T)={\rm cl}_H(T)$. It is easy to see that if
$G$ and $H$ are geometrically equivalent then they have the same
quasiidentities. Is it true that if two groups have the same
quasiidentities then they are geometrically equivalent? This is true
for nilpotent groups.\hfill {\sl B.\,I.\,Plotkin}\vs

No, not always (V.\,V.\,Bludov, {\it Abstracts of the 7th Int. Conf.
Groups and
Group Rings}, Supras\'l, Poland, 1999, p.\,6; \ R.\,G\"obel,
S.\,Shelah, {\it Proc. Amer. Math. Soc.}, {\bf 130} (2002),
673--674 (electronic); \ A.\,G.\,Myasnikov, V.\,N.\,Remeslennikov,
{\it J.~Algebra}, {\bf 234}, no.\,1 (2000), 225--276). The latter
paper contains also necessary and sufficient conditions for
geometrical equivalence.\emp

\bmp \textbf{14.77.} Let $p$ be a prime number and $X$ a finite set of
powers of $p$ containing~$1$. Is it true that $X$ is the set of all
lengths of the conjugacy classes of some finite $p\hs$-group?

\hfill {\sl J.\,Sangroniz}\vs

Yes, it is (J.\,Cossey, T.\,Hawkes, {\it Proc. Amer. Math. Soc.},
{\bf 128}, no.\,1 (2000), 49--51).

 \emp

\bmp \textbf{14.80.}
Is the lattice of all
totally local formations of finite groups modular? The definition
of a {\it totally local\/} formation see in (L.\,A.\,Shemetkov,
A.\,N.\,Skiba, {\it Formatsii algebraicheskikh system}, Mos\-cow,
Nauka, 1989 (Russian)).

\hfill {\sl A.\,N.\,Skiba, L.\,A.\,Shemetkov}

\vs

Yes, even distributive
(V.\,G.\,Safonov, {\it Commun. Algebra}, \textbf{35} (2007),
3495--3502).

 \emp

\bmp
\textbf{14.82.}
(Well-known problem).
Describe the finite simple groups in which
every element is a product of two involutions.
 \hfill {\sl A.\,I.\,Sozutov}

 \vs
Described in (E.\,P.\,Vdovin, A.\,A.\,Gal't, {\it Siberian Math. J.}, {\bf 51} (2010), 610--615).
 \emp

\bmp \textbf{14.86.} Does there exist an infinite locally nilpotent
$p\hs$-group that is equal to its commutator subgroup and in which
every proper subgroup is nilpotent? \hfill {\sl J.\,Wiegold}\vs

No, it does not exist (A.\,O.\,Asar, {\it J.~Lond. Math. Soc.
(2)}, {\bf 61}, no.\,2 (2000), 412--422).

 \emp

\bmp \textbf{14.88.}
We say that an element $u$ of a group $G$ is a {\it test
element\/} if for any endomorphism $\f$ of $G$ the equality $\f
(u)=u$ implies that $\f$ is an automorphism of~$G$. Does a free
soluble group of rank $2$ and derived length $d>2$ have any test
elements?

\hfill
 {\sl B.\,Fine, V.\,Shpilrain}

\vs

Yes, it does
(E.\,I.\,Timoshenko, {\it Algebra and Logic}, \textbf{45}, no.\,4
(2006), 254--260.)

 \emp

\bmp \textbf{14.92.}
(I.\,D.\,Macdonald). Every finite $p\hs$-group
has at least $p-1$ conjugacy classes of maximum size.
I.\,D.\,Macdonald ({\it Proc. Edinburgh Math. Soc.}, \textbf{26}
(1983), 233--239) constructed groups of order $2^n$ for any $n\geq
7$ with just one conjugacy class of maximum size. Are there any
examples with exactly $p-1$ conjugacy classes of maximum size for
odd $p$?
\hfill {\sl G.\,Fern\'andez--Alcober}

\vs

Yes, there are, for $p=3$
(A.\,Jaikin-Zapirain, M.\,F.\,Newman, E.\,A.\,O'Brien, {\it
Israel J. Math.}, \textbf{172}, no.\,1 (2009), 119--123).

 \emp

\bmp \textbf{14.96.} Suppose that a finite $p\hs$-group $P$ admits an
automorphism of order $p^n$ having exactly $p^m$ fixed points. By
(E.\,I.\,Khukhro, {\it Russ. Acad. Sci. Sbornik Math.}, {\bf 80}
(1995), 435--444)
 then
$P$ has a subgroup of index bounded in terms of
 $p$, $n$ and $m$ which is soluble of derived length bounded in terms
of $p^n$. Is it true that
 $P$ has also a subgroup of index bounded in terms of
 $p$, $n$ and $m$ which is soluble of derived length bounded in
terms of $m$? There are positive answers in the cases of $m=1$
(S.\,McKay, {\it Quart. J.~Math. Oxford, Ser. (2)}, {\bf 38}
(1987), 489--502; \ I.\,Kiming, {\it Math. Scand.}, {\bf 62}
(1988), 153--172) and $n=1$ (Yu.\,A.\,Medvedev, see
Archive,~10.68). \hfill {\sl E.\,I.\,Khukhro}\vs

Yes, it is true (A.\,Jaikin-Zapirain, {\it Adv. Math.}, {\bf
153}, no.\,2 (2000),
 391--402).
 \emp

\bmp \textbf{14.99.} b)
A formation ${\frak F}$ of finite groups is
called {\it superradical\/} if it is $S_{n}$-closed and contains
every finite group of the form $G=AB$ where $A$ and $B$ are $
{\frak F}$-subnormal ${\frak \ F}$-sub\-groups.
Prove that every $S$-closed superradical
formation is a solubly saturated formation.
 \hfill
{\sl L.\,A.\,Shemetkov}
\vs

b) A counterexample is constructed (S.\,Yi, S.\,F.\,Kamornikov, 
{\it Siberian Math. J.}, {\bf 57}, no.\,2 (2016), 260--264).
 \emp

\bmp \textbf{14.102.} c)
(V.\,Lin). Let $B_n$ be the braid group
on
 $n$
strings, and let $n>4$.
 Does
$B_n$ have proper non-abelian torsion-free factor-groups?

\makebox[15pt][r]{}{\it Comment of 2001:\/} S.\,P.\,Humphries, ({\it Int. J. Algebra
Comput.}, \textbf{11}, no.\,3 (2001), 363--373) has constructed a
representation of $B_n$ which is shown to provide torsion-free
non-abelian factor groups of $B_n$ as well as of $[B_n, B_n]$ for
$n < 7$. It is likely that the same representation should work for
other values of $n$ as well.
\hfill {\sl V.\,Shpilrain}

\vs

 Yes, it does; moreover,
every braid group $B_n$ is residually torsion-free
nilpotent-by-finite (P.\,Linnell, T.\,Schick, {\it J. Amer. Math.
Soc.}, \textbf{20}, no.\,4 (2007), 1003--1051).
\emp

\bmp \textbf{14.103.} Let $H$ be a proper subgroup of a group $G$ and
let elements $a,b\in H$ have distinct prime orders~$p, q$. Suppose
that, for every $g\in G\setminus H$, the subgroup $\langle
a,b^g\rangle$ is a finite Frobenius group with complement of
order~$pq$. Does the subgroup generated by the union of the kernels
of all Frobenius subgroups of $G$ with complement $\langle a\rangle$
intersect $\langle a\rangle$ trivially? The case where all groups
$\langle a,b^g\rangle$, $g\in G$, are finite is of special interest.
 \hfill {\sl V.\,P.\,Shunkov}\vs

No, in general case not necessarily. As a counter-example one can
take $G=\langle x,y,z,a,b\mid
x^7=y^7=[x,y,y]=[x,y,x]=a^3=b^2=[a,b]=1,\;\; x^a=x^2,\;\;
y^a=y^2,\;\; x^b=x^{-1},\;\; y^b=y^{-1} \rangle$ with $H=\langle
a,b\rangle$. Analogous examples exist also for $p=2$ and every odd
prime $q$. \ (A.\,I.\,Sozutov, {\it Letter of 2002.\/})

\emp

\bmp \textbf{14.104.} An infinite group $G$ is called a {\it monster
of the first kind\/} if it has elements of order $> 2$ and for any
such an element $a$ and for any proper subgroup $H$ of $G$, there
is an element $g$ in $G\setminus H$, such that $\left< a,\,
a^g\right> =G$. A.\,Yu.\,Olshanskii showed that there are
continuously many monsters of the first kind (see~Archive,~6.63).
Does there exist, for any such a monster, a torsion-free group
which is a central extension of a cyclic group by the given
monster? \hfill {\sl V.\,P.\,Shunkov}\vs

No, not always, since for such an extension to exist it is necessary
that every finite subgroup of the monster is cyclic, but this is not
always true (A.\,I.\,Sozutov, {\it Letter of November, 20, 2001\/}).
\emp

\bmp \textbf{15.4.}
Is it true that large growth implies non-amenability? More
precisely, consider a number $\epsilon > 0$, an integer $k \geq
2$, a group $\Gamma$ generated by a set $S$ of $k$ elements, and
the corresponding exponential growth rate $\omega(\Gamma,S)$
defined as in~14.7. For $\epsilon$ small enough, does the
inequality \ $\omega(\Gamma,S) \geq 2k-1-\epsilon$ \ imply that
$\Gamma$ is non-amenable?

\makebox[15pt][r]{}If $\omega(\Gamma,S) = 2k-1$, it is easy to
show that $\Gamma$ is free on $S$, and in particular non-amenable;
see Section 2 in (R.\,I.\,Grigorchuk,
 P.~de~la~Harpe, {\it J.~Dynam. Control Syst.}, {\bf 3}, no.\,1
(1997), 51--89). \hfill {\sl P.~de~la~Harpe}\vs

No, it does not. Counterexamples can be found even in the classes
of abelian-by-nilpotent and metabelian-by-finite groups.
(G.\,N.\,Arzhantseva, V.\,S.\,Guba, L.\,Guyot, {\it J.~Group
Theory}, {\bf 8} (2005), 389--394). \emp

\markboth{\protect\vphantom{(y}{Archive of solved
problems (15th ed., 2002)}}{\protect\vphantom{(y}{Archive of solved
problems (15th ed., 2002)}}

 \bmp \textbf{15.5.}
 (Well-known problem). Does there
exist an infinite finitely generated group which is simple and
amenable? \hfill {\sl P.~de~la~Harpe}
\vs

Yes, it does (K.\,Juschenko, N.\,Monod, {\it
Ann. Math.}, \textbf{178}, no.\,2 (2013), 775--787).
Another example
was later constructed in (V.\,Nekrashevych, {\it Ann. Math.}, {\bf 187}, no.\,3 (2018), 667--719).
\emp

 \bmp \textbf{15.6.}
(Well-known problem). Is it true that Golod $p\hs$-groups are
non-amenable?

\makebox[15pt][r]{}These infinite finitely generated torsion
groups are defined in (E.\,S.\,Golod, {\it Amer. Math. Soc.
Transl. (2)}, \textbf{48} (1965), 103--106). \hfill
{\sl P.~de~la~Harpe}

\vs
Yes, moreover, such groups
have infinite quotients with Kazhdan's property~(T) and are uniformly non-amenable
(M.\,Ershov, A.\,Jaikin-Zapirain, \emph{Proc. London Math. Soc. (3)}, {\bf 102}, no.\,4 (2011), 599--636).
 \emp

\bmp \textbf{15.7.}
(Well-known problem). Is it true that the reduced
$C^*$-alge\-bra of a countable group without amenable normal
subgroups distinct from $\{ 1\}$ is always a simple
$C^*$-alge\-bra with unique trace? See (M.\,B.\,Bekka,
P.~de~la~Harpe, {\it Expos. Math.}, \textbf{18}, no.\,3 (2000),
215--230).

\makebox[15pt][r]{}{\it Comment of 2009\/}: This was proved for
linear groups (T.\,Poznansky, \url{https://arxiv.org/pdf/0812.2486.pdf}).
\hfill {\sl P.~de~la~Harpe}
\vs

Yes, it is true (E.\,Breuillard, M.\,Kalantar, M.\,Kennedy, N.\,Ozawa, {\it Publ. Math. IH\'ES}, {\bf 126}, no.\,1 (2017), 35--71). If $\Gamma$ is assumed to be linear, then $C^*_{\text{red}}(\Gamma)$ is simple if and only if it has a unique trace (Theorem 1.6 ibid., and also Poznansky's 2009 preprint). However, there are infinite countable groups $\Gamma$ having the following two properties: the only amenable normal subgroup is $\{1\}$, and $C^*_{\text{red}}(\Gamma)$ is not simple, as shown in (A.\,Le~Boudec, {\it Invent. Math.}, \textbf{209}, no.\,1 (2017), 159--174).
\emp

\bmp \textbf{15.8.} a)
(S.\,M.\,Ulam). Consider the usual compact group $SO(3)$ of all
rotations of a
 3-dimensional Euclidean space, and let $G$ denote this group viewed
 as a discrete group. Can $G$ act non-trivially on a countable
set?
\hfill {\sl
P.~de~la~Harpe}\vs

Yes, it can (S.\,Thomas, {\it J.~Group Theory}, {\bf 2} (1999),
401--434); another proof is in (Yu.\,L.\,Ershov, V.\,A.\,Churkin,
{\it Dokl. Math.}, {\bf 70}, no.\,3 (2004), 896--898). \emp

\bmp \textbf{15.10.}
 (Yu.\,I.\,Merzlyakov). Is the group of all automorphisms of the free
group $F_n$ that act trivially on $F_n/[F_n,F_n]$ linear for
$n\geq 3$? \hfill {\sl V.\,G.\,Bardakov}

\vs

No, it is not: for $n \geq
5$ (A.\,Pettet, {\it Cohomology of some subgroups of the
automorphism group of a free group}, Ph.D. Thesis, 2006), and for
$n \geq 3$ (V.\,G.\,Bardakov, R.\,Mikhailov, {\it Commun.
Algebra}, \textbf{36}, no.\,4 (2008), 1489--1499).

\emp

\bmp \textbf{15.14.}
Do there exist finitely generated branch groups
(see 15.12)

\makebox[15pt][r]{}a) that are non-amenable?

\makebox[15pt][r]{}c) that contain $F_2$?

\makebox[15pt][r]{}d) that have exponential growth? \hfill {\sl
L.\,Bartholdi, R.\,I.\,Grigorchuk, Z.\,{\v S}uni\'k}\vs

a), c), d) Such groups do exist (S.\,Sidki, J.\,S.\,Wilson, {\it
Arch. Math.}, {\bf 80}, no.\,5 (2003), 458--463).
\emp

\bmp
\textbf{15.15.}
Is every maximal subgroup of a finitely generated
branch group necessarily of finite index?
\hfill
{\sl L.\,Bartholdi, R.\,I.\,Grigorchuk, Z.\,{\v
S}uni\'k}

\vs
No, not necessarily (I.\,V.\,Bondarenko, \emph{Arch. Math. (Basel)}, \textbf{95}, no.\,4 (2010), 301--308).
\emp

 \bmp
\textbf{15.17.}
An infinite group is {\it just infinite\/} if all of its
proper quotients are finite. Is every finitely generated just
infinite group of intermediate growth necessarily a branch group?
\hfill {\sl L.\,Bartholdi, R.\,I.\,Grigorchuk, Z.\,{\v
S}uni\'k}
\vs

No, not every; there exist simple finitely generated groups of intermediate growth (V.\,Nekrashevych, {\it Ann. Math.}, {\bf 87}, no.\,3 (2018), 667--719).
\emp

\bmp \textbf{15.18.}
A group is {\it hereditarily just infinite\/} if it is
residually finite and all of its non-trivial normal subgroups are
just infinite.

\makebox[15pt][r]{}a) Do there exist finitely generated
hereditarily just infinite torsion groups? (It is believed there
are none.)

\makebox[15pt][r]{}b) Is every finitely generated hereditarily
just infinite group necessarily linear?

\makebox[15pt][r]{}A positive answer to the question b) would
imply a negative answer to a).

\hfill {\sl L.\,Bartholdi, R.\,I.\,Grigorchuk,
Z.\,{\v S}uni\'k}

\vs

 a) Yes, such groups exist (M.\,Ershov, A.\,Jaikin-Zapirain, \emph{J.~Reine Angew. Math.}, \textbf{677} (2013), 71--134).
\vs

b) No, not every (M.\,Ershov, A.\,Jaikin-Zapirain, \emph{J.~Reine Angew. Math.}, \textbf{677} (2013), 71--134).
\emp

\bmp \textbf{15.24.}
Suppose that a finite $p\hs$-group $G$ has a subgroup of exponent
$p$ and order $p^n$. Is it true that if $p$ is sufficiently
large relative to $n$,
then $G$
 contains a normal subgroup of exponent $p$ and order $p^n$?

\makebox[15pt][r]{}J.\,L.\,Alperin and G.\,Glauberman ({\it
J.~Algebra}, \textbf{203}, no.\,2 (1998), 533--566) proved that if a
finite $p\hs$-group
 contains an elementary abelian subgroup of order $p^n$, then it
 contains a normal elementary abelian subgroup of the same order
 provided $p>4n-7$, and the analogue for arbitrary
abelian subgroups is proved in (G.\,Glauberman, {\it J.~Algebra},
\textbf{272} (2004), 128--153).
 \hfill {\sl Ya.\,G.\,Berkovich}
\vs

Yes, it is true if $p>n$.
By induction, $G$ contains a subgroup $H$ of exponent $p$ and
order $p^n$ which is normal in a maximal subgroup $M$ of $G$. Then
$H\leq \zeta _{p-1}(M)$. The elements of order $p$ of $\zeta
_{p-1}(M)$ constitute a normal subgroup of $G$, which
contains~$H$. (A.\,Mann, {\it Letter of 1 October 2002}.)

\emp

\bmp \textbf{15.25.} A finite group $G$
is said to be {\it rational\/} if every irreducible character of $G$
takes only rational values. Are the Sylow $2$-sub\-groups of the
symmetric groups ${\Bbb S} _{2^n}$ rational? \hfill
{\sl Ya.\,G.\,Berkovich}

\vs

Yes, they are. A~Sylow
2-subgroup $T_n$ of ${\Bbb S} _{2^n}$ is the wreath product of the
one for ${\Bbb S} _{2^{n-1}}$ with the group $C$ of order 2. If
$T_{n}$ is rational, then the wreath product of $T_{n}$ with $C$
is also rational by Corollary~70 in (D.\,Kletzing {\it Structure
and representations of $Q$-groups} (Lecture Notes in Math., {\bf
1084}), Springer, Berlin, 1984). The result follows by induction.
(A.\,Mann, {\it Letter of 1 October 2002.\/})\emp

\bmp \textbf{15.27.}
Is it possible
that $\text{Aut}\,G\cong\text{Aut}\,H$ for a finite $p\hs$-group
$G$ of order $>2$ and a proper subgroup $H<G$? \hfill
{\sl Ya.\,G.\,Berkovich}

\vs

Yes, it is possible
(T.\,Li, {\it Arch. Math.}, \textbf{92} (2009), 287--290) for
$|G|=2|H|=32$.

\emp

\bmp
 \textbf{15.33.}
 Suppose that all $2$-gene\-ra\-tor subgroups of a finite $2$-group
$G$ are metacyclic. Is the derived length of $G$ bounded? This is
true for finite $p\hs$-groups if $p\ne 2$, see 1.1.8 in
(M.\,Suzuki, {\it Structure of a group and the structure of its
lattice of subgroups}, Springer, Berlin, 1956). \hfill
{\sl Ya.\,G.\,Berkovich}

\vs
Yes, it is at most 2 (E.\,Crestani, F.\,Menegazzo, \emph{J.~Group Theory}, {\bf 15}, no.\,3 (2012), 359--383).
\emp

\bmp \textbf{15.34.}
Is any free product of linearly ordered
 groups with an amalgamated subgroup right-orderable?
\hfill {\sl V.\,V.\,Bludov}

\vs

Yes, it is
(V.\,V.\,Bludov, A.\,M.\,W.\,Glass, {\it Math. Proc. Cambridge
Philos. Soc.}, \textbf{146} (2009), 591--601).
\emp

\bmp \textbf{15.35.} Let $F$ be the free group of
finite rank $ r$ with basis $\{ x_1,\ldots, x_r\}$. Is it true that
there exists a number $C=C(r)$ such that any reduced word
 of length $n>1$ in the $x_i$ lies outside some subgroup of $F$
 of index at most $C
\log n$? \hfill {\sl O.\,V.\,Bogopolski}
\vs

No, it is not true (K.\,Bou-Rabee, D.\,B.\,McReynolds, {\it Bull. London Math. Soc.}, {\bf 43}, no.\,6 (2011), 1059--1068).
\emp

\bmp
 \textbf{15.38.}
 Does there exist a non-local hereditary composition
formation
 $\frak {F} $ of finite groups such that the set of all
 $\frak {F} $-sub\-nor\-mal subgroups is a sublattice of the
subgroup lattice in any finite group?\hfill \raisebox{0ex}{\sl
A.\,F.\,Vasil'ev, S.\,F.\,Kamornikov}

\vs
No, it does not
(S.\,F.\,Kamornikov, {\it Dokl. Math.}, {\bf 81}, no.\,1 (2010), 97--100).
\emp

\bmp \textbf{15.39.}
Axiomatizing the basic properties of subnormal subgroups, we say that
a
 functor $\tau$ associating with every finite group
 $G$ some non-empty set $\tau (G) $ of its subgroups
is an {\it $ETP$-func\-tor\/} if

\makebox[15pt][r]{}1) $ \tau (A) ^ {\varphi} \subseteq \tau (B ) $
and $\tau (B) ^ {\varphi^ {-1}} \subseteq\tau (A) $ for any
epimorphism
 $\varphi : A\longmapsto B$, as well as $\{ H\cap R\mid
R\in\tau
 (G) \} \subseteq\tau (H) $
 for any subgroup $H\leq G$;

\makebox[15pt][r]{}2) $\tau (H) \subseteq\tau (G) $ for any subgroup
$H\in\tau
(G) $;

\makebox[15pt][r]{}3) $\tau (G) $ is a sublattice of the lattice of all
subgroups of
$G$.

\makebox[15pt][r]{}Let $\tau$ be an $ETP$-func\-tor. Does there exist a
hereditary
formation $\frak {F} $ such that $\tau (G) $ coincides with the set
of all
 $\frak {F} $-sub\-nor\-mal subgroups in any finite group $G$?
 This is true for finite soluble groups (A.\,F.\,Vasil'ev,
S.\,F.\,Kamornikov, {\it Siberian Math.~J.}, \textbf{42}, no.\,1
(2001), 25--33). \hfill {\sl A.\,F.\,Vasil'ev,
S.\,F.\,Kamornikov}

\vs

Not always (S.\,F.\,Kamornikov, {\it Math. Notes}, {\bf 89}, no.\,3--4 (2011), 340--348).

\emp

\bmp \textbf{15.43.}
Let $G$ be a
finite group of order $n$.

\makebox[15pt][r]{}a) Is it true that $|{\rm Aut}\,G|\ge\varphi(n)$
where $\varphi$ is Euler's
function?

 \makebox[15pt][r]{}b) Is it true that $G$ is cyclic if $|{\rm
Aut}\,G|=\varphi(n)$? \hfill {\sl M.\,Deaconescu}

 \vs

Both questions have
negative answers; moreover, $|{\rm Aut}\, G|/\varphi(|G|)$ can be
made arbitrarily small (J.\,N.\,Bray, R.\,A.\,Wilson, {\it Bull.
London Math. Soc.}, \textbf{37}, no.\,3 (2005), 381--385). \emp

\bmp \textbf{15.44.}
b) Is the ring of invariants $K[M(n)^m]^{GL(n)}$ Cohen--Macaulay
in all characteristics? (Here $M(n)^m$ is the direct sum of $m$
copies of the space of $n\times n$ matrices.)

\hfill {\sl A.\,N.\,Zubkov}\vsh

Yes, it is (M.\,Hashimoto, {\it Math. Z.}, {\bf 236} (2001),
605--623). \emp

\bmp
 \textbf{15.49.}
 A group
 $G$ is a {\it unique product group}\/ if, for any finite
nonempty subsets $X,Y$ of $G$, there is an element of $G$ which can
be written in exactly one way in the form $xy$ with $x \in X$ and $y
\in Y$. Does there exist a unique product group which is not
left-orderable? \hfill {\sl P.\,Linnell}

\vs Yes, there does (N.\,Dunfield, {\it Appendix~B\,} in S.\,Kionke, J.\,Raimbault,
{\it Doc. Math.}, {\bf 21} (2016), 873--915).
\emp

\bmp \textbf{15.60.}
Is it true that any finitely generated $p\sp\prime$-iso\-la\-ted
subgroup of a free group is separable in the class of finite
$p\hs$-groups? It is easy to see that this is true for cyclic
subgroups. \hfill {\sl D.\,I.\,Moldavanski\u{\i}}\vs

No, it is not (V.\,G.\,Bardakov, {\it Siberian Math. J.}, {\bf
45}, no.\,3 (2004), 416--419).\emp

 \bmp \textbf{15.62.}
Given an ordinary irreducible character $\chi$ of a finite group
 $G$ write $p^{e_p(\chi)}$ to denote the $p\hs$-part of $\chi(1)$ and
put $e_p(G)=\max\{e_p(\chi)\mid\chi\in{\rm Irr}(G)\}$. Suppose
that $P$
 is a Sylow $p\hs$-sub\-group of a group $G$. Is it true that $e_p(P)$
is bounded above by a function of $e_p(G)$?

\makebox[15pt][r]{}{\it Comment of 2005:} An affirmative answer in
the case of solvable groups has been given in (A.\,Moret\'o,
T.\,R.\,Wolf, {\it Adv. Math.}, \textbf{184} (2004), 18--36).\hfill
{\sl A.\,Moret\'o}
\vs

Yes, it is true (Yong Yang, Guohua Qian, {\it Adv. Math.}, {\bf 328} (2018), 356--366).
\emp

 \bmp \textbf{15.63.}
b) Let $F_n$ be the free group of finite rank $n$ on the free
generators $x_1,\ldots ,x_n$. An element $u \in F_n$ is called
{\it positive\/} if $u$ belongs to the semigroup generated by the
$x_i$. An element $u \in F_n$ is called {\it potentially
positive\/} if $\alpha (u)$ is positive for some automorphism
$\alpha$ of~$F_n$. Finally, $u \in F_n$ is called {\it stably
potentially positive\/} if it is potentially positive as an
element of $F_m$ for some $m \geq n$. Are there stably potentially
positive elements that are not potentially positive? \hfill {\sl
A.\,G.\,Myasnikov,~V.\,E.\,Shpilrain}\vs

No, there are none (A.\,Clark, R.\,Goldstein, {\it Commun.
Algebra}, {\bf 33}, no.\,11 (2005), 4097--4104). \emp

 \bmp \textbf{15.75.}
a) Does there exist a sequence of identities in two variables
$u_1=1$, $u_2=1,\ldots $ with the following properties: \ 1)~each
of these identities
 implies the next one, and \ 2)~an arbitrary finite group is
soluble if and only if it satisfies one of the identities $u_n =
1$? \hfill {\sl B.\,I.\,Plotkin}\vs

Yes, such sequences exist (T.\,Bandman, F.\,Grunewald,
G.-M.\,Greuel, B.\,Ku\-nyav\-skii, G.\,Pfister, E.\,Plotkin, {\it
Compositio Math.}, {\bf 142} (2006), 734--764; \
J.\,N.\,Bray, J.\,S.\,Wil\-son, R.\,A.\,Wilson, {\it Bull. Lond.
Math. Soc.}, {\bf 37} (2005) 179--186; \ E.\,Ribnere, {\it Monatsh. Math.}, {\bf 157} (2009), 387--401).

\emp

\bmp \textbf{15.76.} a)
If $\Theta$ is a variety of groups, then let $\Theta ^0$ denote the
category of all free groups of finite rank in $\Theta$. It is proved
(G.\,Mashevitzky, B.\,Plotkin, E.\,Plotkin {\it J.~Algebra}, {\bf
282} (2004), 490--512) that if $\Theta$ is the variety of all groups,
then every automorphism of the category $\Theta ^0$ is an inner one.
The same is true if $\Theta$ is the variety of all abelian groups. Is
this true for the variety of
nilpotent groups of class~2?

\makebox[15pt][r]{}An automorphism $\varphi$ of a category is called
{\it inner\/} if it is isomorphic to the identity automorphism. Let
$s:1\to\varphi$ be a function defining this isomorphism. Then for
every object $A$ we have an isomorphism $s_A: A\to \varphi(A)$ and
for any morphism of objects $\mu: A\to B$ we have
$\varphi(\mu)=s_B\mu s_A^{-1}$.
 \hfill {\sl B.\,I.\,Plotkin}\vs

Yes, it is, even for
the variety of nilpotent groups of any class $n$ (A.\,Tsurkov,
{\it Int. J. Algebra Comput.}, \textbf{17} (2007), 1273--1281).
\emp

\bmp \textbf{15.79.}
 Does there exist a
Hausdorff group topology on ${\Bbb Z}$ such that the sequence
\linebreak $ \{ 2^n+3^n \}$ converges to zero? \hfill
{\sl I.\,V.\,Protasov}

\vs  Yes, there does,
as follows from Theorem~3 in (I.\,Z.\,Ruzsa, \emph{Proc. Conf. Number theory (Budapest, 1987). Vol.\,{\rm I}: Elementary and analytic}, \emph{Colloq. Math. Soc. J\'anos Bolyai}, \textbf{51}, North-Holland, Amsterdam, 1990, 473--504). Another solution is in (S.\,V.\,Skresanov,
\textit{Siberian Math.~J.}, {\bf 61}, no.\,3
 (2020), 542--544).
\emp

\bmp \textbf{15.81.}
Let $G$ be a finite non-supersoluble group. Is it true that $G$
has a non-cyclic Sylow subgroup $P$ such that some maximal
subgroup of $P$ has no proper complement in $G$? \hfill {\sl
A.\,N.\,Skiba}\vs

Yes, it is (Wang, Yanming, Wei, Huaquan, {\it Sci. China Ser. A},
{\bf 47} (2004), no.\,1, 96--103). \emp

 \bmp \textbf{15.86.}
A group $G$ is called {\it discriminating\/} if for any finite set
of nontrivial elements of the direct square $G\times G$ there is a
homomorphism $G\times G\rightarrow G$ which does not annihilate
any of them (G.\,Baumslag, A.\,G.\,Myasnikov,
V.\,N.\,Remeslennikov). A group $G$ is called {\it squarelike\/}
if $G$ is universally equivalent (in the sense of first order
logic) to a discriminating group (B.\,Fine, A.\,M.\,Gaglione,
A.\,G.\,Myasnikov, D.\,Spellman). Must every squarelike group be
elementarily equivalent to a discriminating group?

\hf {\sl D.\,Spellman}
\vs

Yes, it must
(O.\,Belegradek, {\it J. Group Theory}, \textbf{7}, no.\,4 (2004),
521--532; \ B.\,Fine, A.\,M.\,Gaglione, D.\,Spellman, {\it Archiv
Math.}, \textbf{83}, no.\,2 (2004), 106--112). \emp

\bmp \textbf{16.1.} a) Let $G$ be a finite non-abelian group, and
$Z(G)$ its centre. One can associate a graph $\Gamma_G$ with $G$
as follows: take $G\backslash Z(G)$ as vertices of $\Gamma_G$ and
join two vertices $x$ and $y$ if $xy\not=yx$. Let $H$ be a finite
non-abelian group such that $\Gamma_G \cong \Gamma_H$.
If $H$ is simple, is it true that $G\cong
H$? \ {\it Comment of 2009:\/}
This is true for groups with disconnected prime graphs (L.\,Wang,
W.\,Shi, {\it Commun. Algebra}, \textbf{36} (2008), 523--528).

 \hf
{\sl A.\,Abdollahi, S.\,Akbari, H.\,R.\,Maimani}
\vs

a) Yes, it is mod CFSG
(Ch.\,Khan', G.\,Ch\`en', S.\,Go, \emph{
 Siberian Math. J.}, {\bf 49}, no.\,6 (2008), 1138--1146, for sporadic simple groups; \ A.\,Abdollahi, H.\,Shahverdi, {\it J. Algebra}, {\bf 357} (2012), 203--207, for alternating groups; \ R.\,M.\,Solomon, A.\,J.\,Woldar, \emph{J.~Group Theory}, {\bf 16}, no.\,6 (2013), 793--824, for simple groups of Lie type).
\emp

 \bmp \textbf{16.2.}
 A group $G$ is {\it subgroup-separable\/} if for
any subgroup $H\leq G$ and element $x\in G\setminus H$ there is a
homomorphism to a finite group $f: G\rightarrow F$ such that
$f(x)\not\in f(H)$. Is it true that a finitely generated solvable
group is locally subgroup-separable if and only if it does not
contain a solvable Baumslag--Solitar group? Solvable
Baumslag--Solitar groups are $BS(1,n)= \left< a, b \mid
bab^{-1}=a^n\right>$ for $n>1$.

\makebox[15pt][r]{}Background: It is known that a finitely
generated solvable group is subgroup-separable if and only if it
is polycyclic (R.\,C.\,Alperin, in: {\it Groups--Korea'98
(Pusan)}, de~Gruyter, Berlin, 2000, 1--5).
 \hf {\sl R.\,C.\,Alperin}

\vs
No, it is not (J.\,O.\,Button, {\it Ricerche di Matematica}, {\bf 61}, no.\,1 (2012), 139--145).
\emp

\markboth{\protect\vphantom{(y}{Archive of solved
problems (16th ed., 2006)}}{\protect\vphantom{(y}{Archive of solved
problems (16th ed., 2006)}}

\bmp \textbf{16.8.}
The {\it width\/} $w(G')$ of the derived
subgroup $G'$ of a finite non-abelian group $G$ is the smallest
positive integer $m$ such that every element of $G'$ is a product
of $\leq m$ commutators. Is it true that the maximum value of the
ratio $w(G')/|G|$ is $1/6$ (attained at the symmetric group ${\Bbb
S}_3$)? \hf {\sl V.\,G.\,Bardakov}

\vs

Yes, it is (T.\,Bonner,
{\it J.\,Algebra}, \textbf{320} (2008), 3165--3171). \emp

\bmp \textbf{16.12.}
Given a finite $p\hs$-group $G$, we define a
{\it $\Phi$-extension\/} of $G$ as any finite $p\hs$-group $H$
containing a normal subgroup $N$ of order $p$ such that $H/N\cong
G$ and $N\leq\Phi(H)$. Is it true that for every finite
$p\hs$-group $G$ there exists an infinite sequence
$G=G_1,\,G_2,\ldots $ such that $G_{i+1}$ is a $\Phi$-extension of
$G_i$ for all $i=1,\,2,\ldots$?\hf {\sl
Ya.\,G.\,Berkovich}

\vs

Yes, it is
(S.\,F.\,Kamornikov, {\it Izv. Gomel' Univ.}, \textbf{5} (2008),
200--201 (Russian)). \emp

\bmp \textbf{16.13.}
Does there exist a finite $p$-group $G$
all of whose maximal subgroups $H$ are special, that is, satisfy
$Z(H)=[H,H]=\Phi(H)$?
\hf {\sl Ya.\,G.\,Berkovich}
\vs

Yes, for all primes there are groups of arbitrarily large size with this property (J.\,Cossey, \textit{Bull. Austral. Math. Soc.}, {\bf 89} (2014), 415--419); an example for $p=2$ given by a Sylow $2$-subgroup of $L_3(4)$ was independently presented by V.\,I.\,Zenkov at Mal'cev Meeting--2014, 10--14 November, 2014, Novosibirsk.
\emp

\bmp \textbf{16.15.} An element $g$ of a group $G$ is an {\it Engel
element\/} if for every $h\in G$ there exists $k$ such that
$[h,g,\ldots ,g]=1$, where $g$ occurs $k$ times; if there is such
$k$ independent of $h$, then $g$ is said to be {\it boundedly
Engel}.

\makebox[25pt][r]{a)}
(B.\,I.\,Plotkin). Does the set of boundedly
Engel elements of a group form a subgroup?
 \hf {\sl V.\,V.\,Bludov}

\vs

a) No, not always (A.\,I.\,Sozutov,
{\it Siberian Math.~J.}, {\bf 60}, no.\,6 (2019), 1099--1100).
\emp

\bmp \textbf{16.17.}
Is it true that a non-abelian simple group cannot contain Engel
elements other than the identity element?\hf {\sl
V.\,V.\,Bludov}

\vs

No, it can: since an
involution is an Engel element in any 2-group, counterexamples are
provided by non-abelian simple 2-groups; see Archive, 4.74a).
(V.\,V.\,Bludov, {\it Letter of 30 March 2006}.)
\emp

\bmp \textbf{16.24.}
The {\it spectrum\/} of a finite group is the
set of orders of its elements. Does there exist a finite group
$G$ whose spectrum coincides with the spectrum of a finite simple
exceptional group $L$ of Lie type, but $G$ is not isomorphic to
$L$? \hf {\sl A.\,V.\,Vasil'ev}
\vs

Yes, for example, for $L={}^3D_4(2)$ (V.\,D.\,Mazurov, 
{\it Algebra Logic}, {\bf 52}, no.\,5 (2013), 400--403). Further examples are given in
M.\,A.\,Grechkoseeva, M.\,A.\,Zvezdina, {\it J. Algebra Appl.}, {\bf 15}, no.\,9 (2016), Article ID 1650168, 13 pp.).
\emp

\bmp \textbf{16.25.}
Do there exist three pairwise non-isomorphic
finite non-abelian simple groups with the same spectrum? \hf
{\sl A.\,V.\,Vasil'ev}

\vs
No (A.\,A.\,Buturlakin, {\it Sibirsk. Electron. Math. Reports}, {\bf 7} (2010), 111--114).
\emp

\bmp
 \textbf{16.27.}
 Suppose that a finite group $G$ has the
same spectrum as an alternating group. Is it true that $G$ has at
most one non-abelian composition factor?

\makebox[15pt][r]{}For any finite simple group other than
alternating
 the answer to the corresponding question is
 affirmative (mod CFSG).
 \hf \raisebox{0ex}{\sl A.\,V.\,Vasil'ev and
V.\,D.\,Mazurov}

\vs
Yes, it is (I.\,B.\,Gorshkov, {\it Algebra and Logic}, {\bf 52}, no.\,1 (2013), 41--45). 
\emp

\bmp \textbf{16.31.}
Suppose that a group $G$ has a composition
series and let ${\frak F}(G)$ be the formation generated by $G$.
Is the set of all subformations of ${\frak F}(G)$ finite?
\hf {\sl V.\,A.\,Vedernikov}

\vs
No, not always, even for $G$ finite (V.\,P.\,Burichenko, \textit{J.~Algebra}, {\bf 372}
(2012), 428--458).
\emp

 \bmp \textbf{16.35.}
Is every finitely presented soluble group
 nilpotent-by-nilpotent-by-finite?

 \hf{\sl
J.\,R.\,J.\,Groves}

\vs

No: for a prime $p$ the
group $G=\langle a,b,c,d\mid b^a =b^p,\;\, c^a=c^p,\;\, d^a=d,\;\,
c^b=c,\break d^c=dd^b,\;\,[d,d^b]=1,\;\, d^p=1\rangle $ is soluble
of derived length~3, but is not metanilpotent-by-finite. (V.\,V.\,Bludov, {\it Letter of 27 April 2007}.)
\emp

\bmp \textbf{16.36.}
 (Well-known problem).
 We call a finite group {\it rational\/} if all of its
ordinary characters are rational-valued. Is every Sylow 2-subgroup
of a rational group also a rational group? \hf {\sl M.\,R.\,Darafsheh}

\vs
No, not always
(I.\,M.\,Isaacs, G.\,Navarro, \textit{Math. Z.}, {\bf 272} (2012), 937--945).
\emp

 \bmp \textbf{16.37.}
 Let $G$ be a solvable rational finite group with an
extra-special Sylow 2-sub\-group. Is it true that either $G$ is a
2-nilpotent group, or there is a normal subgroup $E$ of $G$ such
that $G/E$ is an extension of a normal 3-subgroup by an elementary
abelian 2-group?\hf {\sl M.\,R.\,Darafsheh}
\vs

Yes, it is (M.\,R.\,Darafsheh, H.\,Sharifi, {\it Extracta Math.},
{\bf 22}, no.\,1 (2007), 83--91).
\emp

\bmp \textbf{16.42.}
Is a topological Abelian group $(G,\tau)$
compact if every group topology $\tau '\subseteq\tau$ on $G$ is
complete? (The answer is yes if every continuous homomorphic image
of $(G,\tau)$ is complete.) \hf {\sl
E.\,G.\,Zelenyuk}

\vs  Yes, it is (T.\,Banakh, \textit{Topology Appl.}, \textbf{271} (2020), Article ID 106983, 17 p.).
\emp

\bmp \textbf{16.43.}
Is there a
partition of the group $\bigoplus_{\omega_1}({\mathbb Z}/3{\mathbb
Z})$ into three subsets whose complements do not contain cosets
modulo infinite subgroups? (There is a partition into two such
subsets.)
 \hf {\sl E.\,G.\,Zelenyuk}
\vs

Yes, moreover, every
infinite abelian group with finitely many involutions can be
partitioned into infinitely many subsets such that every coset
modulo an infinite subgroup meets each subset of the partition
(Y.\,Zelenyuk, {\it J.~Combin. Theory (A)}, \textbf{115} (2008),
331--339). \emp

\bmp \textbf{16.49.}
Is it true that a free product of groups
without generalized torsion is a group without generalized
torsion?\hf {\sl V.\,M.\,Kopytov,
N.\,Ya.\,Medvedev}

\vs  Yes, it is true; moreover, the generalized torsion in a free product of torsion-free groups is conjugate to a generalized torsion of one of its factor groups (T.\,Ito, K.\,Motegi, M.\,Teragaito,
\textit{Proc. Amer. Math. Soc.}, {\bf 147}, no.\,11 (2019), 4999--5008).
\emp

\bmp \textbf{16.50.}
Do there exist simple finitely generated
right-orderable groups?

\makebox[15pt][r]{}There exist finitely generated right-orderable
groups coinciding with the derived subgroup (G.\,M.\,Bergman).\hf
{\sl V.\,M.\,Kopytov, N.\,Ya.\,Medvedev}

\vs  Yes, such groups do exist (J.\,Hyde, Y.\,Lodha, {\it Invent. Math.}, {\bf 218} (2019), 83--112).
\emp

\bmp \textbf{16.52.}
Is every finitely presented elementary amenable group
solvable-by-finite?

\hf {\sl P.\,Linnell, T.\,Schick}

\vs
No, not every. I.\,Belegradek and Y.\,Cornulier have pointed out that the
groups in (C.\,H.\,Houghton, {\it Arch. Math. (Basel)}, {\bf 31}, no.\,3 (1978/79), 254--258) are finitely presented as shown in (K.\,S.\,Brown, {\it J.~Pure Appl. Algebra}, {\bf 44}, no.\,1--3 (1987), 45--75) and elementary amenable, but not virtually solvable; \url{http://mathoverflow.net/questions/107996}.
\emp

\bmp \textbf{16.54.} We say that a group $G$ acts {\it freely\/}
on a group $V$ if $vg\ne v$ for any nontrivial elements $g\in
G,\; v\in V$. Is it true that a group $G$ that can act
 freely on a non-trivial abelian group is embeddable in the
multiplicative group of some skew-field?

\hf {\sl V.\,D.\,Mazurov}
\vs

No, it is not. For example, the group $2.A_5.2$ with a quaternion Sylow 2-subgroup can act freely on an elementary abelian group of order $7^4$ but is not embeddable in the multiplicative group of any skew-field by Theorem~7 in (S.\,A.\,Amitsur, {\it Trans. Amer. Math. Soc.}, {\bf 80}, no.\,2 (1955), 361--386). (D.\,Nedrenko, {\it Letter of 20 January 2014}.)
\emp

\bmp \textbf{16.55.}
(Well-known
problem). Let $V$ be a faithful absolutely irreducible module for
a finite group $G$. Is it true that $\dim H^1(G,V)\leq 2$? \hf
 {\sl V.\,D.\,Mazurov}

\vs

No, it is not
(L.\,L.\,Scott, {\it J.~Algebra}, \textbf{260} (2003), 416--425; see
also J.\,N.\,Bray, R.\,A.\,Wilson, {\it J.~Group Theory}, \textbf{11}
(2008), 669--673).
 \emp

\bmp \textbf{16.57.}
Suppose that $\omega(G)=\omega(L_2(7))=\{1,2,3,4,7\}$. Is $G\cong
L_2(7)$? This is true for finite $G$. \hf {\sl
V.\,D.\,Mazurov}

\vs

Yes, it is
(D.\,V.\,Lytkina, A.\,A.\,Kuznetsov, {\it Siberian Electr. Math.
Rep.}, \textbf{4} (2007), 136--140; \url{http://semr.math.nsc.ru}).

 \emp

\bmp \textbf{16.58.} 
 Is $SU_2({\Bbb C})$ the only group
that has just one irreducible complex representation of dimension
$n$ for each $n = 1,\,2,\ldots $?

\makebox[15pt][r]{}(If $R[n]$ is the $n$-dimensional irreducible
complex representation of $SU_2({\Bbb C})$, then $R[2]$ is the natural
two-dimensional representation, and $R[2]\otimes R[n] = R[n-1] +
R[n+1]$ for $n > 1$.)
 \hf {\sl J.\,McKay}
 \vs

No, it is not. It is known that ${\mathbb C}$ has infinitely many (discrete) automorphisms. For $\varphi\in{\rm Aut}({\mathbb C})$ and a matrix $x$, let $x^\varphi$ denote the matrix obtained from $x$ by applying $\varphi$ to each element. Then $T_\varphi:x\mapsto x^\varphi$ is an irreducible
representation of $SU_2({\mathbb C})$. It is easy to show that among these representations there are infinitely many pairwise non-equivalent ones.
Therefore the group $SU_2({\mathbb C})$ itself does not satisfy the condition of the problem. (This observation belongs to von~Neumann.)
One can show that any group satisfying the condition of the problem is isomorphic to a subgroup of $SL_2({\mathbb Q})$ satisfying the condition of Problem~15.57. On the other hand, E.\,Cartan proved that a unitary representation of a simple compact Lie group is always continuous. Hence, if we restrict ourselves to the \textit{unitary} representations, then $SU_2({\mathbb C})$ does satisfy the condition of the problem. (V.\,P.\,Burichenko, \emph{Letter of 16 July 2013}.)
\emp

\bmp \textbf{16.59.}
Given a finite group $K$, does there exist a
finite group $G$ such that $K\cong {\rm Out\,}G={\rm Aut\,} G/{\rm
Inn\,}G$? (It is known that an infinite group $G$ exists with this
property.)\hf {\sl D.\,MacHale}

\vs
 Yes, it does (Y.\,Cornulier, \url{https://mathoverflow.net/questions/372480/is-every-finite-group-the-outer-automorphism-group-of-a-finite-group/372563}). \emp

 \bmp \textbf{16.61.}
A subgroup $H$ of a group $G$ is {\it fully invariant\/} if
$\vartheta (H)\leq H$ for every endomorphism $\vartheta$ of~$G$.
Let $G$ be a finite group such that $G$ has a fully invariant
subgroup of order $d$ for every $d$ dividing $|G|$. Must $G$ be
cyclic? \hf {\sl D.\,MacHale}

\vs

No: take $G=C_p\times
C_{p^2}$, where $C_p,C_{p^2}$ are cyclic $p$-groups of orders
$p,\;p^2$ (A.\,Abdol\-lahi, {\it Letter of 4 March 2009}). \emp

\bmp \textbf{16.62.}
Let $G$ be a group such that every $\alpha
\in {\rm Aut\,} G$ fixes every conjugacy class of $G$ (setwise).
Must ${\rm Aut\,} G ={\rm Inn\,}G$? \hf {\sl
D.\,MacHale}

\vs

Not necessarily: the paper
(H.\,Heineken, {\it Arch. Math. $($Basel\/$)$}, \textbf{33}
(1979/80), no.\,6, 497--503), in particular, produces finite
$p$-groups all of whose automorphisms preserve all the conjugacy
classes, while by Gasch\"utz' theorem every finite $p$-group has
outer automorphisms (A.\,Mann, {\it Letter of 20 April 2006\/}).
\emp

 \bmp \textbf{16.65.}
Does there exist a finitely presented residually torsion-free
nilpotent group with a free presentation $G=F/R$ such that the
group $F/[F,R]$ is not residually nilpotent?
\hf {\sl R.\,Mikhailov,
I.\,B.\,S.\,Passi}

\vs

Yes, there does
(V.\,V.\,Bludov, {\it J. Group Theory}, \textbf{12}, no.\,4 (2009),
579--590). \emp

 \bmp \textbf{16.67.}
{\it Conjecture:\/} Given any integer $k$, there exists an integer
$n_0=n_0(k)$ such that if $n\geq n_0$ then the symmetric group of
degree $n$ has at least $k$ different ordinary irreducible
characters of equal degrees. \hf {\sl
A.\,Moret\'o}

\vs

This is proved
(D.\,Craven, {\it Proc. London Math. Soc.}, \textbf{96} (2008),
26--50).
\emp

\bmp \textbf{16.71.}
Is the elementary theory of a torsion-free
hyperbolic group decidable?

\hf {\sl
A.\,G.\,Myasnikov, O.\,G.\,Kharlampovich}
\vs

Yes, it is (O.\,Kharlampovich, A.\,Myasnikov, \url{https://arxiv.org/pdf/1303.0760.pdf}).
\emp

\bmp \textbf{16.73.} a) Let $G$ be a group generated by a finite set
$S$, and let $l(g)$ denote the word length function of $g\in G$
with respect to $S$. The group $G$ is said to be {\it
contracting\/} if there exist a faithful action of $G$ on the set
$X^*$ of finite words over a finite alphabet $X$ and constants
$0<\lambda <1$ and $C>0$ such that for every $g\in G$ and $x\in X$
there exist $h\in G$ and $y\in X$ such that $l(h)<\lambda l(g)+C$
and $g(xw)=yh(w)$ for all $w\in X^*$.

\makebox[15pt][r]{}Can a contracting group
have a non-abelian free subgroup?
 \hf {\sl V.\,V.\,Nekrashevich}
\vs

a) No, it cannot
(V.\,V.\,Nekrashevich, {\it Groups Geom. Dynam.}, {\bf 4}, no.\,4 (2010), 847--862).
\emp
\bmp \textbf{16.74.} a)
Let $G=\left<\alpha,\beta\right>$ be the group generated by the
following two permutations of ${\Bbb Z}$: \ $\alpha (n)=n+1$; \
$\beta (0)=0$, \ $\beta (2^km)=2^k(m+2)$, where $m$ is odd and $n$
is a positive integer. Is $G$ amenable?
\hf {\sl V.\,V.\,Nekrashevich}

\vs
a) Yes, it is
(G.\,Amir, O.\,Angel, B.\,Vir\'ag, {\it J. Eur. Math. Soc.}, {\bf 15}, no.\,3 (2013), 705--730).
\emp

\bmp \textbf{16.75.}
Can a non-abelian one-relator group be the
group of all automorphisms of some group?
 \hf {\sl M.\,V.\,Neshchadim}

\vs

Yes, it can: in
(D.\,J.\,Collins, {\it Proc. London Math. Soc. (3)}, \textbf{36}
(1978), no.\,3, 480--493) it was proved that for integers $r,s$
such that $(r,s)=\nobreak 1$, $|r|\ne |s|$, and $r-s$ is even, the
group $G=\langle a,b\mid a^{-1}b^ra=b^s\rangle$ is
 isomorphic to ${\rm Aut}({\rm Aut} (G))$ (V.\,A.\,Churkin, {\it
Letter of 5 April 2006\/}).
\emp

\bmp \textbf{16.79.} Is it true that in any finitely generated
$AT$-group over a sequence of cyclic groups of uniformly bounded
orders all Sylow subgroups are locally finite? For the definition
of an {\it $AT $-group\/} see (A.\,V.\,Rozhkov, {\it Math. Notes},
\textbf{40} (1986), 827--836)

\hf {\sl A.\,V.\,Rozhkov}

\vs
No, it is not true (A.\,V.\,Rozhkov, in {\it Group theory and its applications}, Proc. XXII School-Conf. on Group Theory, Kuban' Univ., Krasnodar, 2018, 126--131 (Russian)).
\emp

\bmp \textbf{16.82.}
Let $\cal{X}$ be a non-empty class of finite groups closed under
taking homomorphic images, subgroups, and direct products. With
every group $G\in\cal{X}$ we associate some set $\tau(G)$ of
subgroups of $G$. We say that $\tau$ is a {\it subgroup functor
on\/} $\cal{X}$ if:

\makebox[25pt][r]{1)} $G\in\tau(G)$ for all $G\in\cal{X}$, and

\makebox[25pt][r]{2)} for each epimorphism $\varphi :A\mapsto B$,
 where $A,B\in\cal{X}$, and for any $H\in\tau(A)$ and

 \makebox[25pt][r]{} $T\in\tau(B)$ we have $H^{\varphi}\in\tau(B)$ and
$T^{\varphi^{-1}}\in\tau(A)$.

A subgroup functor $\tau$ is {\it closed\/} if for each group
$G\in{\cal X}$ and for every subgroup $H\in {\cal X} \cap \tau
(G)$ we have $\tau (H)\subseteq \tau (G)$. The set $F(\cal {X})$
consisting of all closed subgroup functors on $\cal X$ is a
lattice (in which $\tau _1 \leq \tau _2$ if and only if $\tau_1
(G) \subseteq \tau_2 (G)$ for every group $G\in\cal X$). It is
known that $F(\cal {X})$ is a chain if and only if $\cal X$ a
class of $p$-groups for some prime $p$ (Theorem~1.5.17 in
S.\,F.\,Kamornikov and M.\,V.\,Sel'kin, {\it Subgroups functors
and classes of finite groups}, Belaruskaya Navuka, Minsk, 2001
(Russian)).

\makebox[15pt][r]{}Is there a non-nilpotent class $\cal X$ such
that the width of the lattice $F(\cal {X})$ is at most $|\pi(\cal
X)|$ where $\pi (\cal X)$ is the set of all prime divisors of the
orders of the groups in~$\cal X$?
\hf {\sl A.\,N.\,Skiba}

\vs
No such classes exist
(S.\,F.\,Kamornikov, {\it Siberian Math. J.}, {\bf 51}, no.\,5 
(2010), 824--829).
\emp

\bmp \textbf{16.83.} a) Let $E_n$ be a free locally nilpotent $n$-Engel
group on countably many generators, and let $\pi(E_n)$ be the set
of prime divisors of the orders of elements of the periodic part
of $E_n$. It is known that $2,\,3,\,5\in \pi(E_4)$.

\makebox[15pt][r]{}Does there exist $n$ for which $ 7\in
\pi(E_n)$?
\hf {\sl
Yu.\,V.\,Sosnovski\u{\i}}
\vs

a) Yes, $n=6$, since the 2-generator free nilpotent
 6-Engel group has an element of order 7
(W.\,Nickel, \emph{J.~Austral. Math. Soc., Ser. A}, {\bf 67}, no.\,2 (1999), 214--222;
\url{http://www.mathematik.tu-darmstadt.de/~nickel/Engel/Engel.html})
 (A.\,Abdollahi, {\it Letter of 19 August 2011}).
\emp

\bmp \textbf{16.85.} Suppose that groups $G,\,H$ act faithfully on a
regular rooted tree by finite-state automorphisms. Can their free
product $G*H$ act faithfully on a regular rooted tree by finite
state automorphisms?\hf {\sl
V.\,I.\,Sushchanski\u{\i}}

\vs  Yes, it can (M.\,Fedorova, A.\,Oliynyk, {\it Algebra Discrete Math.}, {\bf 23}, no.\,2 (2017), 230--236).
\emp

 \bmp \textbf{16.86.} Does the group of all finite-state
automorphisms of a regular rooted tree possess an irreducible
system of generators?\hf {\sl
V.\,I.\,Sushchanski\u{\i}}

\vs
 Yes, it does (Ya.\,Lavrenyuk, {\it Geometriae Dedicata}, {\bf 183}, no.\,1 (2016), 59--67).
\emp

\bmp \textbf{16.101.}
Do there exist uncountably many infinite 2-groups that are
quotients of the group $\left<x, y \mid x^2 = y^4 = (xy)^8 =
1\right>$? There certainly exists one, namely the subgroup of
finite index in Grigorchuk's first group generated by $b$ and
$ad$; see (R.\,I.\,Grigorchuk, {\it Functional Anal. Appl.},
\textbf{14} (1980), 41--43).
 \hf {\sl J.\,Wiegold}

\vs

Yes, there do
(R.\,I.\,Grigorchuk, {\it Algebra Discrete Math.}, {\bf 2009}, no.\,4 (2009), 78--96 (subm. 25 November 2009); \
A.\,Minasyan, A.\,Yu.\,Olshanskii, D.\,Sonkin, {\it Groups Geom.
Dynam.}, \textbf{3}, no.\,3 (2009) 423--452 (subm. 21 April 2008)). \emp

\bmp \textbf{16.103.}
Is there a rank analogue of the
Leedham-Green--McKay--Shepherd
theorem on $p\hs$-groups of maximal class? More precisely, suppose
that $P$ is a 2-generator finite $p\hs$-group whose lower central
quotients $\gamma _i(P)/\gamma _{i+1}(P)$ are cyclic for all
$i\geq 2$. Is it true that $P$ contains a normal subgroup $N$ of
nilpotency class $\leq 2$ such that the rank of $P/N$ is bounded
in terms of $p$ only?\hf {\sl E.\,I.\,Khukhro}

\vs
No, moreover, there are no functions $d(p)$ and
$r(p)$ such that a group with these properties would necessarily
have a normal subgroup of derived length ${}\leq d(p)$ with quotient of rank
${}\leq r(p)$ (E.\,I.\,Khukhro, {\it Siber. Math. J.}, {\bf 54}, no.\,1 (2013),
174--184). 
\emp

 \bmp \textbf{16.104.}
If \( G \) is a finite group, then every element \( a \) of the
rational group algebra \( {\mathbb Q}[G] \) has a unique Jordan
decomposition \( a =a _{s}+a _{n} \), where \( a _{n}\in {\mathbb
Q}[G] \) is nilpotent, \( a _{s}\in {\mathbb Q}[G] \) is
semisimple over \( \mathbb Q \), and \( a _{s}a _{n}=a _{n}a _{s}
\). The integral group ring \( {\mathbb Z}[G] \) is said to have
the \emph{additive Jordan decomposition\/} property (AJD) if \( a
_{s},\, a _{n} \in {\mathbb Z}[G] \) for every
 \( a \in {\mathbb Z}[G] \). If \( a \in {\mathbb Q}[G] \) is invertible, then
 \( a _{s} \) is also invertible and so \( a =a _{s}a _{u} \)
with \( a _{u}=1+a _{s}^{-1}a _{n} \) unipotent and \( a _{s}a
_{u}=a _{u}a _{s} \). Such a decomposition is again unique. We say
that \( {\mathbb Z}[G] \) has \emph{multiplicative Jordan
decomposition\/} property (MJD) if \( a _{s},\, a _{u} \in
{\mathbb Z}[G] \) for every invertible
 \( a \in {\mathbb Z}[G] \). See the survey (A.\,W.\,Hales, I.\,B.\,S.\,Passi,
 in:
 {\it Algebra, Some Recent Advances}, Birkh\"{a}user, Basel, 1999, 75--87).

\makebox[15pt][r]{}Is it true that there are only finitely many
isomorphism classes of finite $2$-groups $G$ such that ${\mathbb
Z}[G]$ has MJD but not AJD?
 \hf {\sl A.\,W.\,Hales, I.\,B.\,S.\,Passi}

\vs
 Yes, it is true (A.\,W.\,Hales, I.\,B.\,S.\,Passi, L.\,E.\,Wilson, \emph{J.~Algebra}, \textbf{316}, no.\,1 (2007), 109--132; \textbf{371} (2012), 665--666).
\emp

\bmp \textbf{16.106.}
Let $\pi _e(G)$ denote the set of orders of
elements of a group $G$, and $h(\Gamma )$ the number of
non-isomorphic finite groups $G$ with $\pi _e(G)=\Gamma$. Do there
exist two finite groups $G_1,\,G_2$
such that $\pi _e(G_1)=\pi _e(G_2)$, \ $h(\pi _e(G_1)) < \infty$,
and neither of the two groups $G_1,\,G_2$ is isomorphic to a
subgroup or a quotient of a normal subgroup of the other? \hf
 {\sl W.\,J.\,Shi}

\vs

Yes, there do: for example, $G_1=L_{15}(2^{60}).3$ and
$G_2=L_{15}(2^{60}).5$ (M.\,A.\,Grech\-ko\-se\-eva, {\it Algebra and Logic}, \textbf{47}, no.\,4
(2008), 229--241). 
\emp

\bmp \textbf{16.107.}
Is
it true that almost every alternating group
${A}_n$ is uniquely determined in the class of finite groups by
its set of element orders, i.\,e. that $h(\pi _e({ A}_n))=1$ for
all large enough $n$? \hf {\sl W.\,J.\,Shi}

\vs

Yes, it is (I.\,B.\,Gorshkov, {\it Algebra and Logic}, {\bf 52}, no.\,1 (2013), 41--45).
\emp

\bmp \textbf{16.109.}
Is there a polynomial time algorithm for
solving the word problem in the group ${\rm Aut}\, F_n$ (with
respect to some particular finite presentation), where $F_n$ is
the free group of rank $n \geq 2$?
 \hf {\sl V.\,E.\,Shpilrain}

\vs

Yes, there is
(S.\,Schleimer, {\it Comment. Math. Helv.}, \textbf{83} (2008),
741--765).
\emp

\markboth{\protect\vphantom{(y}{Archive of solved
problems (17th ed., 2010)}}{\protect\vphantom{(y}{Archive of solved
problems (17th ed., 2010)}}

 \bmp \textbf{17.2.}
 (P.\,Schmid). Does there exist a finite non-abelian $p$-group $G$ such that $H^1({G}/{\Phi(G)},\,Z(\Phi(G) ) )=0$?
\hf {\sl A.\,Abdollahi}

\vs
Yes, it does (A.\,Abdollahi, {\it J. Algebra}, {\bf 342} (2011), 154--160).
 \emp

\bmp \textbf{17.3.}
Let $G$ be a group in which every $4$-element
subset contains two elements generating a nilpotent subgroup. Is
it true that every $2$-generated subgroup of $G$ is nilpotent? \hf
{\sl A.\,Abdollahi}
\vs

 No, not always (A.\,I.\,Sozutov,
{\it Siberian Math.~J.}, {\bf 60}, no.\,6 (2019), 1099--1100).
\emp

 \bmp \textbf{17.12.}
 Are there functions
$e,c:\mathbb{N}\rightarrow \mathbb{N}$ such that if in a nilpotent
group $G$ a normal subgroup $H$ consists of right $n$-Engel
elements of $G$, then $H^{e(n)}\leq \zeta_{c(n)}(G)$?

\hf {\sl A.\,Abdollahi}

\vs

Yes, there are
(P.\,G.\,Crosby, G.\,Traustason, {\it J. Algebra}, {\bf 324} (2010), 875--883).
 \emp

\bmp \textbf{17.14.}
Can the braid group $B_n$ for $n \geq 4$ be embedded into the
automorphism group $Aut (F_{n-2})$ of a free group $F_{n-2}$ of
rank~$n-2$?

 \makebox[15pt][r]{}Artin's theorem implies that $B_n$ can be
embedded into $Aut (F_n)$, and an embedding of $B_n$, $n \geq 3$,
into $Aut (F_{n-1})$ was constructed in (B.\,Perron,
J.\,P.\,Vannier, {\it Math. Ann.}, \textbf{306} (1996), 231--245).
 \hf {\sl V.\,G.\,Bardakov}
\vs

No, it cannot for $n=4$ (V.\,G.\,Bardakov, P.\,Bellingeri, in {\it Knot theory and its applications. ICTS program knot theory and its applications (KTH-2013), IISER Mohali, India, December 10--20, 2013}, Amer. Math. Soc., {\it Contemporary Mathematics}, {\bf 670} (2016), 285--298).
 \emp

\bmp \textbf{17.17.}
If a finitely generated group $G$ has
$n<\infty$ maximal subgroups, must $G$ be finite? In particular,
what if $n=3$? \hf {\sl G.\,M.\,Bergman}
\vs

No, it need not. An example of an infinite group with 3 maximal subgroups is given by a 2-generated 2-LERF group in \S\,7 of (M.\,Ershov, A.\,Jaikin-Zapirain, {\it J. Reine Angew. Math.} {\bf 677} (2013), 71--134) as
its maximal subgroups have index~2.
 \emp

\bmp \textbf{17.19.}
If $F$ is a free group of finite rank, $R$ a
retract of $F$, and $H$ a subgroup of $F$ of finite rank, must
$H\cap R$ be a retract of $H$\/?\hf {\sl
G.\,M.\,Bergman}

\vs
 No, it need not (I.\,Snopce, S.\,Tanushevski, P.\,Zalesskii, {\it Int. Math. Res. Notices}, {\bf 2022}, no.\,11 (2022), 8280--8294). \emp

\bmp \textbf{17.20.}
If $M$ is a real manifold with nonempty
boundary, and $G$ the group of self-homeo\-mor\-phisms of $M$
which fix the boundary pointwise, is $G$ right-orderable?

\hfill {\sl G.\,M.\,Bergman}

\vs  No, not always (J.\,Hyde,
{\it Ann. of Math. (2)}, {\bf 190}, no.\,2 (2019), 657--661).
 \emp

\bmp \textbf{17.21.}
 a) If $A, B, C$ are torsion-free abelian
groups with $A \cong A\oplus B\oplus C$, must $A\cong A\oplus B$?

\makebox[15pt][r]{}b) What if, furthermore, $B\cong C$? \hf
\raisebox{0ex}{\sl G.\,M.\,Bergman}

\vs
No, not necessarily, even in case b). See
(A.\,L.\,S.\,Corner, in {\it Proc. Colloq. Abelian
Groups (Tihany, 1963)}, Budapest, 1964, 43--48); I also have an example with
$B$ and $C$ of rank 1 as torsion-free abelian groups. (G.\,M.\,Bergman, {\it Letter of 20 February 2011}.)

 \emp

\bmp \textbf{17.24.}
(A.\,Blass, J.\,Irwin, G.\,Schlitt). Does
$\mathbb{Z}^\omega$ have a subgroup whose dual is free of
uncountable rank? \hf {\sl G.\,M.\,Bergman}

\vs
Yes, it does
(G.\,M.\,Bergman, {\it Portugaliae Math.}, {\bf 69} (2012) 69--84).
\emp

\bmp \textbf{17.28.}
 Is there a soluble right-orderable group with insoluble word problem?

 \hf {\sl V.\,V.\,Bludov, A.\,M.\,W.\,Glass}

 \vs Yes, there is (A.\,Darbinyan, \emph{J.~Symbolic Logic}, \textbf{85}, no.\,4 (2020), 1588--1598). \emp

\bmp \textbf{17.29.}
Construct a finitely presented orderable group
with insoluble word problem.

\hf \raisebox{0ex}{\sl
V.\,V.\,Bludov, A.\,M.\,W.\,Glass}

\vs
Such a group is constructed (V.\,V.\,Bludov, A.\,M.\,W.\,Glass, {\it
Bull. London Math. Soc.}, {\bf 44} (2012), 85--98).
 \emp

\bmp \textbf{17.35.}
Suppose we have a finite two-colourable
triangulation of the sphere, with triangles each coloured either
black or white so that no pair of triangles with the same colour
share an edge. On each vertex we write a generator, and we assume
the generators commute. We use these generators to generate an
abelian group $G_W$ with relations stating that the generators
around each white triangle add to zero. Doing the same thing with
the black triangles, we generate a group $G_B$.

\makebox[15pt][r]{}{\it Conjecture:}
$G_W$ is isomorphic to $G_B$. \hf {\sl
 I.\,M.\,Wanless, N.\,J.\,Cavenagh}
\vs

Conjecture is proved
(S.\,R.\,Blackburn, T.\,A.\,McCourt,
\textit{Combinatorica}, {\bf 34}, no.\,5 (2014), 527--546).
\emp

\bmp \textbf{17.36.}
Two groups are called {\it isospectral\/} if
they have the same set of element orders. Find all finite
non-abelian simple groups $G$ for which there is a finite group
$H$ isospectral to $G$ and containing a proper normal subgroup
isomorphic to $G$. For every simple group $G$ determine all groups
$H$ satisfying this condition.

\makebox[15pt][r]{}It is easy to show that a group $H$ must
satisfy the condition $G<H\leqslant\operatorname{Aut}G$.

\hf {\sl A.\,V.\,Vasil'ev}
\vs

All such groups are determined: for exceptional groups of Lie type in (M.\,A.\,Zvezdina, {\it Algebra Logic},
{\bf 55}, no.\,5 (2016), 354--366);
for the other simple groups in (M.\,A.\,Grechkoseeva,
{\it Siberian Math. J.}, {\bf 59}, no.\,4 (2018), 623--640).
\emp

\bmp \textbf{17.40.}
Let $N$ be a nilpotent subgroup of a finite group $G$. Do there
always exist elements $x,y\in G$ such that $N\cap N^x\cap N^y\leq
F(G)$?
\hf {\sl E.\,P.\,Vdovin}

\vs Yes, such elements always exist, mod CFSG (V.\,I.\,Zenkov, {\it Siberian Math.~J.}, {\bf 62}, no.\,4 (2021), 621--637).
 \emp

\bmp \textbf{17.44.}
Let $\pi$ be a set of primes. A finite group is called a {\it $C_\pi$-group\/} if it possesses exactly one class of conjugate Hall
$\pi$-subgroups. A finite group is called a {\it $D_{\pi} $-group\/} if any two of its maximal $\pi
 $-sub\-groups are conjugate.

\makebox[25pt][r]{a)} In a $C_\pi$-group, is an overgroup of a Hall $\pi$-subgroup always a $C_\pi$-group?

\makebox[15pt][r]{}An affirmative answer in the case $2\not\in\pi$ follows
mod CFSG from (F.\,Gross, {\it Bull. London Math. Soc.}, \textbf{19},
no.\,4 (1987), 311--319).

\makebox[25pt][r]{b)} In a $D_\pi$-group, is an overgroup of a Hall $\pi$-subgroup always a
 $D_\pi$-group?

 \hf\raisebox{0ex}{\sl E.\,P.\,Vdovin, D.\,O.\,Revin}

\vs

a) Yes, it is
(E.\,P.\,Vdovin, D.\,O.\,Revin, {\it Siberian Math. J.}, {\bf 54}, no.\,1 (2013), 22--28).\vs

b) Yes, it is (N.\,Ch.\,Man\-zaeva, \url{https://arxiv.org/pdf/1504.03137.pdf}) and (E.\,P.\,Vdovin, N.\,Ch.\,Manzaeva, D.\,O.\,Revin, {\it Sb. Math.}, {\bf 211}, no.\,3 (2020), 309--335).

 \emp

\bmp \textbf{17.45.}
A subgroup $H$ of a group $G$ is called {\em pronormal\/} if $H$
and $H^g$ are conjugate in $\langle H,H^g\rangle$ for every $g\in
G$. We say that $H$ is {\em strongly pronormal\/} if $L^g$ is
conjugate to a subgroup of $H$ in $\langle H,L^g\rangle$ for every
$L\leq H$ and $g\in G$.

\makebox[25pt][r]{a)} In a finite simple group, are Hall subgroups always pronormal?

\makebox[15pt][r]{}An affirmative answer is known modulo
CFSG for Hall subgroups of odd order (F.\,Gross, {\it Bull. London
Math. Soc.}, \textbf{19}, no.\,4 (1987), 311--319).

\makebox[25pt][r]{b)} In a finite simple group, are Hall subgroups always strongly pronormal?

\makebox[15pt][r]{}Notice that there
exist finite (non-simple) groups with a non-pronormal Hall
subgroup. Hall subgroups with a Sylow tower are known to be
pronormal.

 \hf \raisebox{0ex}{\sl E.\,P.\,Vdovin, D.\,O.\,Revin}

\vs

a) Yes, they are (E.\,P.\,Vdovin, D.\,O.\,Revin, {\it Siberian Math. J.}, 
{\bf 53}, no.\,3 (2012), 419--430).\vs

b) No, not always (M.\,N.\,Nesterov, {\it Siberian Math. J.}, {\bf 58},\ no.\,1 (2017), 
 128--133).
 \emp

 \bmp \textbf{17.46.}
 Let $G$ be a finite $p$-group in which every
2-generator subgroup has cyclic derived subgroup. Is the derived
length of $G$ bounded?

\makebox[16pt][r]{}If $p\ne 2$, then $G^{(2)}=1$ (J.\,Alperin),
but for $p=2$ there are examples with $G^{(2)}$ cyclic and
elementary abelian of arbitrary order.
 \hf {\sl B.\,M.\,Veretennikov}
 \vs

 Yes, it is at most 4 (B.\,Wilkens, {\it J. Group Theory}, {\bf 17}, no.\,1 (2014), 151--174).
 \emp
\bmp \textbf{17.50.}
Is it true that for every finite group $G$,
there is a finite group $F$ and a surjective homomorphism $f: F
\to G$ such that for each nontrivial subgroup $H$ of $F$, the
restriction $f|_H$ is not injective?

\makebox[15pt][r]{}It is known that for every finite group $G$
there is a finite group $F$ and a surjective homomorphism $f: F
\to G$ such that for each subgroup $H$ of $F$ the restriction
$f|_H$ is not bijective. \hf {\sl E.\,G.\,Zelenyuk}

\vs
Yes, it is
(V.\,P.\,Burichenko, {\it Math. Notes}, {\bf 92}, no.\,3 (2012), 327--332).
 \emp

\bmp \textbf{17.55.}
Does there exist an absolute constant $k$ such that for any
prefrattini subgroup $H$ in any finite soluble group $G$ there
exist $k$ conjugates of $H$ whose intersection is~$\Phi (G)$?
\hf {\sl S.\,F.\,Kamornikov}
\vs

Yes, such a constant exists and is equal to 3 (S.\,F.\,Kamornikov, {\it Int. J. Group Theory}, {\bf 6}, no.\,2 (2017), 1--5).
 \emp

 \bmp
\textbf{17.67.}
(H.\,Zassenhaus). {\it Conjecture\/}: Every
invertible element of finite order
of the integral group ring ${\mathbb Z}G$ of a
finite group $G$ is conjugate in the rational group ring ${\mathbb
Q}G$ to an element of $\pm G$. \hf {\sl
V.\,D.\,Mazurov}
\vs

A counterexample was constructed (F.\,Eisele, L.\,Margolis,
{\it Adv. Math.} {\bf 339} (2018), 599--641).
\emp

 \bmp \textbf{17.71.} a) Let $\alpha$ be a fixed-point-free automorphism
of prime order $p$ of a periodic group~$G$.
Is it true that $G$ does not contain a
non-trivial $p$-element?
\hfill {\sl V.\,D.\,Mazurov}
\vs

a) No, it may contain a non-trivial $p$-element (E.\,B.\,Durakov, A.\,I.\,Sozutov, {\it Algebra Logic}, {\bf 52}, no.\,5 (2013), 422--425).
\emp

\bmp \textbf{17.72.} a)
Let $AB$ be a Frobenius group with kernel $A$
and complement $B$. Suppose that $AB$ acts on a finite group $G$
so that $GA$ is also a Frobenius group with kernel $G$ and
complement $A$. Is the nilpotency class of $G$ bounded in
terms of $|B|$ and the 
class of $C_G(B)$?
\hf {\sl V.\,D.\,Mazurov}

\vs
a) Yes, it is (N.\,Yu.\,Makarenko, P.\,Shumyatsky,
\emph{Proc. Amer. Math. Soc.}, {\bf 138} (2010), 3425--3436).
\emp

\bmp \textbf{17.73.} Let $G$ be a finite simple group of Lie type
defined over a field of characteristic~$p$, and $V$ an absolutely
irreducible $G$-module over a field of the same characteristic.
Is it true that in the following cases the split extension of $V$ by $G$ must contain an element whose order is distinct from the order of any element
of~$G$?

\makebox[25pt][r]{a)} $G=U_4(q)$;

\makebox[25pt][r]{b)} $G=S_{2n}(q), n\geq 3$; \

\makebox[25pt][r]{c)} $G=O_{2n+1}(q), n\geq 3$; \

\makebox[25pt][r]{d)} $G=O_{2n}^+(q), n\geq 4$; \

\makebox[25pt][r]{e)} $G=O_{2n}^-(q), n\geq 4$; \

\makebox[25pt][r]{f)} $G=~^3D_{4}(q)$, $q\ne 2$; \

\makebox[25pt][r]{g)} $G=E_6(q)$; \

\makebox[25pt][r]{h)}~$G=~^2E_6(q)$; \

\makebox[25pt][r]{i)} $G=E_7(q)$; \

\makebox[25pt][r]{j)} $G=G_2(2^m)$. \
\hf {\sl V.\,D.\,Mazurov}
\vs

Yes, it is true: a) (M.\,A.\,Grechkoseeva, S.\,V.\,Skresanov, {\it Siberian Electron. Math. Rep.}, {\bf 17} (2020), 585--589); \  b)--i) (M.\,A.\,Grechkoseeva, {\it J.~Algebra Appl.}, {\bf 14}, no.\,4 (2015), Article ID 1550056); \  j) (A.\,V.\,Vasil'ev, A.\,M.\,Staroletov, {\it Algebra Logic}, {\bf 52}, no.\,1 (2013), 1--14).
 \emp

\bmp \textbf{17.74.}
Let $G$ be a finite simple group of Lie type
defined over a field of characteristic~$p$ whose Lie rank is at
least three, and $V$ an absolutely irreducible $G$-module over a
field of characteristic that does not divide $p$. It is true that
the split extension of $V$ by $G$ must contain an element whose
order is distinct from the order of any element of~$G$? The case
of $G=U_n(p^m)$ is of special interest.
\hf {\sl V.\,D.\,Mazurov}
\vs

Yes, it is true: for linear groups (A.\,V.\,Zavarnitsine, {\it Siberian Math.~J.}, {\bf 49} (2008), 246--256); for other classical groups (M.\,A.\,Grechkoseeva, {\it J.~Algebra}, {\bf 339}, no.\,1 (2011), 304--319); for exceptional gorups (M.\,A.\,Grechkoseeva, {\it J. Algebra Appl.}, {\bf 14}, no.\,4 (2015), Article ID 1550056).
 \emp

 \bmp \textbf{17.77.}
Let $k$ be a positive integer such that there is an insoluble
finite group with exactly $k$ conjugacy classes. Is it true that a
finite group of maximal order with exactly $k$ conjugacy classes
is insoluble?\hf \raisebox{0ex}{\sl R.\,Heffernan, D.\,MacHale}

\vs
No, it is not; see Table 1 in (A.\,Vera-Lopez, J.\,Sangroniz, \emph{Math. Nachr.}, {\bf 280}, no.\,5--6 (2007), 676--694).
 \emp

 \bmp \textbf{17.79.}
Does there exist a finitely generated torsion
group of unbounded exponent generating a proper variety? \hf
{\sl O.\,Macedo\'nska}

\vs
Yes, moreover, there is a continuum of such groups
(V.\,S.\,Atabekyan, {\it Infinite simple groups satisfying an identity}, Dep. VINITI no.\,5381-V86, Moscow, 1986 (Russian)).
Another example was suggested by D.\,Osin in a letter of 31 August 2013: a free group $G$ in the variety $\frak M$ defined by the law $x^ny=yx^n$ is a central extension of a free Burnside group of exponent $n$ such that the centre of $G$ is a free
abelian group of countable rank (I.\,S.\,Ashmanov, A.\,Yu.\,Olshanskii, {\it Izv. Vyssh. Uchebn. Zaved. Mat.}, {\bf 1985}, no.\,11 (1985), 48--60 (Russian)). Let $x_1, x_2, \dots $ be a basis in $Z(G)$. By adding to $G$ the relations $x_i^i=1$ for all $i$ we obtain a periodic group of unbounded exponent generating a proper variety (contained in $\frak M$).
\emp

\bmp \textbf{17.82.}
Is it true that in every finitely presented
group the intersection of the derived series has trivial
abelianization? \hf {\sl R.\,Mikhailov}

\vs

No, not always (V.\,A.\,Roman'kov, {\it Siberian Math. J.}, {\bf 52}, no.\,2 (2011), 348--351).
 \emp

\bmp \textbf{17.86.} (Simplest questions related to the
 Whitehead asphericity conjecture). Let $\langle x_1,x_2,x_3\mid
r_1,r_2,r_3\rangle$ be a presentation of the trivial group.

\makebox[25pt][r]{b)} Let $F=F(x_1,x_2,x_3)$ and $R_i=\langle
r_i\rangle^F$. Is it true that the group $F/[R_1,R_2]$ is
residually soluble?
 \hf
{\sl R.\,Mikhailov}
\vs

b) No, not always (V.\,A.\,Roman'kov, {\it Siberian Math. J.}, {\bf 52}, no.\,2 (2011), 348--351).
 \emp

\bmp \textbf{17.92.}
What are the non-abelian composition factors of
a finite non-soluble group all of whose maximal subgroups are Hall
subgroups?
\hf {\sl V.\,S.\,Monakhov}

\vs
Answer: $L_2(7)$, $L_2(11)$, $L_5(2)$
(N.\,V.\,Maslova, {\it Siberian Math. J.}, {\bf 53}, no.\,5 (2012), 853--861).
 \emp

\bmp \textbf{17.108.} Is the group $\langle a,b,t \mid a^t=ab,\;\,
b^t=ba\rangle$ linear?

\makebox[15pt][r]{}If not, this would be an easy example of a
non-linear hyperbolic group.
 \hf
{\sl M.\,Sapir}

\vs  Yes, this group is linear. Indeed, as noticed by M.\,Sapir, this group is
the mapping torus of an irreducible atoroidal self-monomorphism of a free group; thus it is virtually special, and hence $\Z$-linear by Theorem~B in (M.\,F.\,Hagen, D.\,T.\,Wise, {\it Duke Math.~J.}, {\bf 165}, no.\,9 (2016), 1753--1813). (M.\,F.\,Hagen, {\it Letter of 6 August 2018}.)
\emp

\bmp \textbf{17.111.}
Let $G$ be a finite group, and $p$ a prime
divisor of~$|G|$. Suppose that every maximal subgroup of a Sylow
$p$-subgroup of $G$ has a $p$-soluble supplement in $G$. Must $G$
be $p$-soluble? \hf {\sl A.\,N.\,Skiba}

\vs
Yes, it must (GuoHua Qian, {\it Science China Mathematics}, {\bf 56}, no.\,5 (2013), 1015--1018). \emp

\bmp \textbf{17.116.}
Let $n(G)$ be the maximum of positive integers
$n$ such that the $n$-th direct power of a finite simple group $G$
is 2-generated. Is it true that $n(G)\geq \sqrt{|G|}$?

\hf\raisebox{0ex}{\sl A.\,Erfanian, J.\,Wiegold}

\vs

Yes, it is, even with $n(G)>2\sqrt{|G|}$
(A.\,Mar\'oti, M.\,C.\,Tamburini, {\it Commun. Algebra}, {\bf 41}, no.\,1 (2013), 34--49).
\emp

\markboth{\protect\vphantom{(y}{Archive of solved
problems (18th ed., 2014)}}{\protect\vphantom{(y}{Archive of solved
problems (18th ed., 2014)}}

\bmp
 \textbf{18.9.} Does there exist a subgroup-closed saturated formation $\mathfrak F$ of finite groups
properly contained in $\mathfrak{E}_{\pi}$, where $\pi = char(\mathcal F)$, satisfying the following property: if $G \in \mathfrak F$, then there exists a prime $p$ (depending on the group $G$) such that the wreath product $C_p\wr G$ belongs to $\mathfrak F$, where $C_p$ is the cyclic group of order $p$?

\hfill {\sl A.\,Ballester-Bolinches}
\vs

Yes, there does
(S.\,F.\,Kamornikov, {\it Trudy. Inst. Mat. (Minsk)}, {\bf 24}, no.\,1 (2016), 30--33).
\emp

\bmp
 \textbf{18.15.} For an automorphism $\varphi \in Aut(G)$ of a group $G$, let $[e]_{\varphi} =\{g^{-1}g^{\varphi}\mid g\in G\}$.
Is it true that if a group $G$ has trivial center, then there is an inner automorphism $\varphi$ such that
$[e]_{\varphi}$ is not a subgroup?

\hfill {\sl V.\,G.\,Bardakov, M.\,V.\,Neshchadim, T.\,R.\,Nasybullov}
\vs

No, not always; see V.\,V.\,Bludov's Example~2 in (D.\,Gon\c{c}alves, T.\,Nasybullov, {\it Commun. Algebra}, {\bf 47}, no.\,3 (2019), 930--944).
 \emp

\bmp
 \textbf{18.17.}
 Is there a torsion-free group which is finitely presented in the quasi-variety of torsion-free groups but not finitely presentable in the variety of
all groups?

\hfill {\sl O.\,V.\,Belegradek}
\vs

Yes, there is (A.\,I.\,Budkin, {\it Siberian Electron. Math. Rep.}, {\bf 14} (2017), 937--945, \url{http://semr.math.nsc.ru}).
\emp

\bmp \textbf{18.23.}
The normal covering number of the symmetric group $S_n$ of degree $n$ is the minimum number $\gamma(S_n)$ of proper subgroups $H_1,\ldots,H_{\gamma(S_n)}$ of $S_n$
such that every element of $S_n$ is conjugate to an element of $H_i$, for some $i=1,\dots, \gamma(S_n)$.
Write $n=p_1^{\alpha_1}\cdots p_r^{\alpha_r}$ for primes $p_1<\cdots<p_r$ and positive integers $\alpha_1,\ldots,\alpha_r$.

 \makebox[15pt][r]{}{\it Conjecture\/}:
\[
\gamma(S_n)=\left\{
\begin{array}{lcl}
\frac{n}{2}\left(1-\frac{1}{p_1}\right)&&\textrm{if }r=1\,\,\textrm{and}\,\,\alpha_1=1\\
\frac{n}{2}\left(1-\frac{1}{p_1}\right)+1&&\textrm{if }r=1\,\,\textrm{and}\,\,\alpha_1\geq 2\\
\frac{n}{2}\left(1-\frac{1}{p_1}\right)\left(1-\frac{1}{p_2}\right)+1&&\textrm{if }r=2\,\,\textrm{and}\,\,\alpha_1+\alpha_2=2\\
\frac{n}{2}\left(1-\frac{1}{p_1}\right)\left(1-\frac{1}{p_2}\right)+2&&\textrm{if }r\geq 2\,\,\textrm{and}\,\,\alpha_1+\cdots+\alpha_r\geq 3\\
\end{array}
\right.
\]
This is the strongest form of the conjecture. We would be also interested in a proof that this holds for $n$ sufficiently large. The result for $r\leq 2$, which includes the first three cases above, is proved for $n$ odd (D.\,Bubboloni, C.\,E.\,Praeger, {\it J. Combin. Theory (A)}, {\bf 118} (2011), 2000--2024). When $r\geq3$ we know that $cn\leq \gamma(S_n)\leq \frac{2}{3}n$ for some positive constant $c$ (D.\,Bubboloni, C.\,E.\,Praeger, P.\,Spiga, {\it J. Algebra}, {\bf 390} (2013) 199--215). We showed that the conjectured value for $\gamma(S_n)$ is an upper bound, by constructing a normal covering for $S_n$ with this number of conjugacy classes of maximal subgroups, and gave further evidence for the truth of the conjecture in other cases (D.\,Bubboloni, C.\,E.\,Praeger, P.\,Spiga, {\it Int. J. Group Theory}, {\bf 3}, no.\,2 (2014), 57--75).
\hf {\sl D.\,Bubboloni, C.\,E.\,Praeger, P.\,Spiga}

\vs
 A.\,Mar\'oti showed that the conjecture is incorrect for odd~$n$; his counterexample is presented in \S\,2 of (D.\,Bubboloni, C.\,E.\,Praeger, P.\,Spiga, {\it Monatsh. Math.}, {\bf 191} (2020), 229--247).
 \emp

\bmp
 \textbf{18.30.}
 A subgroup $H$ of a group $G$ is called \emph{$\mathbb{P}$-subnormal} in $G$ if either $H=G$, or there is a chain of subgroups $H_0\subset H_1\subset \cdots\subset H_{n-1}\subset H_{n} =G$ such that $|H_i : H_{i-1}|$ is a prime for all $i=1,\ldots, n$. Must a finite group be soluble if every Shmidt subgroup of it is $\mathbb{P}$-subnormal? \hfill {\sl A.\,F.\,Vasil'ev, T.\,I.\,Vasil'eva, V.\,N.\,Tyutyanov}
\vs

Yes, it must (V.\,N.\,Tyutyanov, {\it Probl. Fiz. Mat. Tekhn.}, {\bf 2015}, no.\,1 (22), 88--91 (Russian)).
\emp

\bmp \textbf{18.31.}
Let $\pi$ be a set of primes. Is it true that in any
$D_\pi$-group $G$ (see Archive, 3.62) there are three Hall $\pi$-subgroups whose intersection coincides with $O_\pi(G)$?

\hf {\sl E.\,P.\,Vdovin, D.\,O.\,Revin}

\vs No, it is not true, for example, for $G$ being an extension of $L_2(27)$ by a field automorphism of order $3$, in which $H$ is an extension of a Borel subgroup of $L_2(27)$ by a field automorphism of order $3$ (V.\,I.\,Zenkov, {\it Siberian Math.~J.}, {\bf 63}, no.\,4 (2022), 720--722).
 \emp

\bmp \textbf{18.32.}
Is every Hall subgroup of a finite group pronormal in its normal closure?

\hf {\sl E.\,P.\,Vdovin, D.\,O.\,Revin}
 \vs

 No, not always (M.\,N.\,Nesterov, {\it Siberian Math. J.}, {\bf 58},\ no.\,1 (2017), 
 128--133).
 \emp

\bmp \textbf{18.33.}
A group in which the derived subgroup of every 2-generated subgroup is cyclic is called an \emph{Alperin group}. Is there a bound for the derived length of finite Alperin groups?

\makebox[15pt][r]{}G.\,Higman proved that finite Alperin groups are soluble, and finite Alperin $p$-groups have bounded derived length, see Archive 17.46.
 \hfill {\sl B.\,M.\,Veretennikov}
 \vs

 Yes, there is: it is at most 6. Indeed, in (P.\,Longobardi, M.\,Maj, H.\,Smith, {\it Rend. Semin. Mat. Univ. Padova}, {\bf 115} (2006), 29--40) it was proved that finite Alperin groups of odd order are metabelian, and any finite Alperin group is supersoluble, so the elements of odd order form a metabelian normal subgroup, while finite Alperin 2-groups have derived length at most 4 by (B.\,Wilkens, {\it J. Group Theory}, {\bf 17}, no.\,1 (2014), 151--174).
\emp

\bmp \textbf{18.49.}
Let \(n \in \mathbb{N}\). Is it true that for any \(a, b, c \in \mathbb{N}\) satisfying
\(1 < a, b, c \leq n-2\) the symmetric group \({\rm S}_n\) has elements of order \(a\)
and \(b\) whose product has order \(c\)?

\hf {\sl
S.\,Kohl}

\vs Yes, it is
(G.\,A.\,Miller, \emph{Amer. J. Math.}, \textbf{22}, no.\,2 (1900), 185--190).
Another solution is in (J.\,K\"onig, {\it Eur. J. Comb.}, {\bf 57} (2016), 50--56).
\emp

 \bmp \textbf{18.52.}
 Is every finite simple group generated by two elements of prime-power orders $m$, $n$? (Here numbers $m,n$ may depend on the group.) The work of many authors shows that it remains to verify this property for a small number of finite simple groups.

 \hfill {\sl J.\,Krempa}
\vs

Yes, it is (mod CFSG); moreover, every finite simple group is generated by an involution and
an element of prime order (C.\,S.\,H.\,King, {\it J. Algebra}, {\bf 478} (2017), 153--173).
 \emp

\bmp \textbf{18.57.} Let $G$ be a finite 2-group generated by involutions in which
$[x,u,u]=1$ for every $x\in G$ and every involution $u\in G$.
Is the derived length of $G$ bounded?

\hf{\sl D.\,V.\,Lytkina}

\vs

No, it is not (A.\,Abdollahi,
 {\it J.~Group Theory}, {\bf 18}, no.\,1 (2015), 111--114).
 \emp

\bmp \textbf{18.73.}
b) Does every finitely generated solvable group of derived length $l\ge 2$ embed
 into some $k$-generated $(l+1)$-solvable group, where $k=k(l)$?
\hf {\sl A.\,Yu.\,Olshanskii}

\vs b) It is proved that any finitely generated solvable group of derived length $l$ can be embedded in a $4$-generated solvable group of derived length $l + 1$ (V.\,A.\,Roman'kov, {\it Proc. Amer. Math. Soc.}, \textbf{149} (2021), 4133--4143).
 \emp

\bmp \textbf{18.75.}
 Does every finite solvable group $G$ have the following property: there is a number
$d=d(G)$ such that $G$ is a homomorphic image of every group with $d$ generators and one relation?

\makebox[15pt][r]{}This property holds for finite nilpotent groups and does not hold for every
non-solvable finite group; see (S.\,A.\,Za\u{\i}tsev,
{\it Moscow Univ. Math. Bull.}, {\bf 52}, no.\,4 (1997), 42--44). 
\hf {\sl A.\,Yu.\,Olshanskii}
\vs

Yes, it does (N.\,Nikolov, D.\,Segal, {\it Bull. London Math. Soc.}, {\bf 39}, no.\,2 (2007), 209--213).
 \emp

\bmp \textbf{18.82.} Is there a function $f:{\Bbb N} \to {\Bbb N}$ such that for any prime $p$, if $p^{f(n)}$ divides the order of a finite group $G$, then $p^n$ divides the order of ${\rm Aut}\,G$?
\hf {\sl R.\,M.\,Patne}
\vs

Yes, there is. This was established in (J.\,A.\,Green, {\it Proc. Roy. Soc. London Ser.~A.}, {\bf 237} (1956), 574--581), which improved a previous result with a function depending on $p$ in (W.\,Ledermann, B.\,H.\,Neumann, {\it Proc. Roy. Soc. London. Ser.~A}, {\bf 235} (1956), 235--246).
\emp

\bmp
 \textbf{18.86.}
 Is the group $G = \langle a, b \mid [[a,b],b]=1\rangle$, which is isomorphic to the group of all unitriangular automorphisms of the free group of rank 3, linear?\hfill {\sl V.\,A.\,Roman'kov}
 \vs

 Yes, it is linear, since it embeds in the holomorph ${\rm Hol}\,F_2$ of the free group $F_2$, which can be seen by adding a new generator $c=[a,b]$, so that $G = \langle a, b, c \mid a^b=ac,\; c^b=c\rangle$, while ${\rm Hol}\,F_2$ was shown to be linear in Corollary~3 in (V.\,G.\,Bardakov, O.\,V.\,Bryukhanov, {\it Vestnik Novosibirsk Univ. Ser. Mat. Mekh. Inf.}, {\bf 7}, no.\,3 (2007), 45--58 (Russian)) (O.\,V.\,Bryukhanov, {\it Letter of 27 January 2014}). Another proof can be found in (V.\,A.\,Roman'kov, {\it J.~Siberian Federal Univ. Math. Phys.}, {\bf 6}, no.\,4 (2013), 516--520).
\emp

\bmp
 \textbf{18.91.}
 A subgroup $H$ of a group $G$ is said to be \emph{propermutable} in $G$
if there is a subgroup $B\leq G$ such that $G=N_{G}(H)B$ and $H$
permutes with every subgroup of~$B$.

 \makebox[15pt][r]{}a) Is there a finite group $G$ with subgroups $A\leq B\leq G$ such that $A$ is propermutable in $G$ but $A$ is not propermutable in~$B$?

\makebox[15pt][r]{}b) Let $P$ be a non-abelian Sylow 2-subgroup of a finite group $G$ with $|P|=2^{n}$. Suppose that there is an integer $k$ such that $1 < k < n$ and every subgroup of $P$ of order $2^{k}$ is propermutable in $G$, and also, in the case of $k=1$, every cyclic subgroup of $P$ of order 4 is propermutable in $G$. Is it true that then $G$ is 2-nilpotent?

\hfill {\sl A.\,N.\,Skiba}
 \vs

 a) Yes, there is
(A.\,A.\,Pypka, D.\,Yu.\,Storozhenko,
{\it Dopov. Nac. Akad. Nauk Ukrain.}, {\bf 2017}, no.\,7, 18--20 (Ukrainian)).
\vs

b) Yes, it is (Kh.\,A.\,Al-Sharo, Finite groups with given systems of weakly $S$-pro\-per\-mu\-table subgroups,
{\it J.~Group Theory}, {\bf 19}, no.\,5 (2016), 871--887).
\emp

\bmp
 \textbf{18.94.}
 Let $G$ be a group without involutions, $a$ an element of it that is not a square of any element of $G$, and $n$ an odd positive integer.
Is it true that the quotient $G/\langle (a^{n})^{G}\rangle $ does not contain involutions?
\hfill {\sl A.\,I.\,Sozutov}
 \vs

 No, not always, as shown by an example of V.\,I.\,Trofimov; another example: the quotient of $G=\langle a,b\mid a^2=b^2\rangle$ by $\langle a^G\rangle$ has order 2.
\emp

\bmp
 \textbf{18.95.}
 Suppose that a group $G=AB$ is a product of an abelian subgroup $A$ and a locally quaternion group $B$ (that is, $B$ is a union of an increasing chain of finite generalized quaternion groups). Is $G$ soluble?
\hfill {\sl A.\,I.\,Sozutov}
 \vs

 Yes, it is (B.\,Amberg, Ya.\,Sysak, On products of groups with abelian subgroups of small index, {\it J.~Group Theory}, {\bf 20}, no.\,6 (2017), 1061--1072).
\emp

\bmp \textbf{18.112.}
Is it true that the orders of all elements of a finite group $G$ are powers of primes if, for every divisor $d$ ($d>1$) of $|G|$ and for every subgroup $H$ of $G$
of order coprime to $d$, the order $|H|$
divides the number of elements of $G$ of order $d$?

\makebox[15pt][r]{}The converse is true (W.\,J.\,Shi, {\it Math. Forum (Vladikavkaz)}, {\bf 6} (2012), 152--154).

\hf {\sl W.\,J.\,Shi}
\vs

Yes, it is true (A.\,A.\,Buturlakin, R.\,Shen, W.\,Shi, {\it Siberian Math. J.}, {\bf 58}, no.\,3 (2017), 405--407).
 \emp

\bmp
\textbf{18.115.}
 Let $G$ be a finite simple group, and let $X, Y$ be isomorphic simple maximal subgroups of $G$. Are $X$ and $Y$ conjugate in ${\rm Aut}\,G$?\hfill {\sl P.\,Schmid}

\vs  No, not always. For
example, there are two ${\rm Aut}\,G$-classes of $M_{12}$ in $E_6(5)$ (P.\,B.\,Kleidman, R.\,A.\,Wilson, \emph{J.~London Math.
Soc.}, \textbf{42} (1990), 555--561); other counterexamples in $E_6(q)$ for certain $q$
include two classes of $PSL_2(11)$, and two of $PSL_2(19)$ (D.\,A.\,Craven, {\it Invent. Math.}, {\bf 234}, no.\,2 (2023), 637--719).
 (D.\,A.\,Craven \textit{Letter of 8 September 2021}.)
\emp

\markboth{\protect\vphantom{(y}{Archive of solved
problems (19th ed., 2018)}}{\protect\vphantom{(y}{Archive of solved
problems (19th ed., 2018)}}

\bmp \textbf{19.23.}
For a group $G$, let $\operatorname{Tor}_{1}(G)$ be the normal closure of all torsion elements of $G$, and then by induction let $\operatorname{Tor}_{i+1}(G)$ be the inverse image of $\operatorname{Tor}_{1}(G/\operatorname{Tor}_{i}(G))$. The torsion length of $G$ is defined to be either the least positive integer $l$ such that $G/\operatorname{Tor}_{l}(G)$ is torsion-free, or $\omega$ if no such integer exists (since $G/\bigcup \operatorname{Tor}_{i}(G)$ is always torsion-free).

\makebox[15pt][r]{}Does there exists a finitely generated, or even finitely presented, soluble group with torsion length greater than $2$?
\hfill {\sl M.\,Chiodo, R.\,Vyas}

\vs Yes, such groups exist (I.\,J.\,Leary, A.\,Minasyan, {\it J.~Group Theory}, {\bf 26}, no.\,4 (2023), 741--750).
\emp

 \bmp \textbf{19.24.}
For a group $G$, let $\operatorname{Tor}(G)$ be the normal closure of all torsion elements of~$G$. Does there exist a finitely presented group $G$ such that $G/\operatorname{Tor}(G)$ is not finitely presented? Such a group must necessarily be non-hyperbolic.
\hfill {\sl M.\,Chiodo, R.\,Vyas}

\vs  Yes, such groups exist: one soluble example and another virtually torsion-free are constructed in (I.\,J.\,Leary, A.\,Minasyan, {\it J.~Group Theory}, {\bf 26}, no.\,4 (2023), 741--750).
\emp

\bmp \textbf{19.35.}
 Let $G$ be a finite group of order $n$. Is it true that for every factorization $n=a_1 \cdots a_k$ there exist subsets $A_1,\dots ,A_k$ such that
$|A_1|=a_1,\dots ,|A_k|=a_k$ and $G=A_1 \cdots A_k$?
\hf {\sl M.\,H.\,Hooshmand}

\vs No, it is not true. A counterexample with $k=3$ is given by the alternating group on 4 letters $G={A}_4$ and $(a_1, a_2, a_3) =
(2,3,2)$. (G.\,M.\,Bergman, {\it Letter of 19 December 2019}, \url{https://math.berkeley.edu/~gbergman/papers/gp_factzn.pdf}.) Counterexamples for all $k\geq 3$ are constructed in (R.\,McCulloch, {\it Preprint}, 2026, \url{https://arxiv.org/abs/2607.20569}).
\emp

\bmp \textbf{19.36.}
Let $G$ be a periodic group and let $\mathscr{I}=\{ x \in G \mid x^{2}=1 \not =x \}$ be the set of its involutions. Let $D$ be a non-empty set of odd integers greater than $1$; then $G$ is called a group with $D$-involutions if $G=\langle \mathscr{I} \rangle$ and for $x,y \in \mathscr{I}$ the order of $xy$ is in the set $\{1,2\} \cup D$ and all these values actually occur.
It is clear that if $G$ is a group with $D$-involutions, then $\mathscr{I}$ is a single conjugacy class.

 \makebox[15pt][r]{} {\it Conjecture:} If $G$ is a group with $\{3,5\}$-involutions, then $G \simeq A_{5}$ or $G \simeq PSU(3,4)$.
\hfill {\sl E.\,Jabara}

 \vs The conjecture is proved, even without using the hypothesis
that the group $G$ is periodic (E.\,Bettio, \emph{J.~Group Theory}, \textbf{24} (2021), 1055--1067).
\emp

\bmp \textbf{19.37.}
{a)} Does there exist an absolute constant $k$ such that for any nilpotent injector $H$ of an arbitrary finite group $G$ there are $k$ conjugates of $H$ the intersection of which is equal to the Fitting subgroup $F(G)$ of $G$?

 \makebox[25pt][r]{b)} Can one choose $k=3$ as such a constant? This is true for finite soluble groups (D.\,S.\,Passman,
 {\it Trans. Amer. Math. Soc.}, {\bf 123}, no.\,1 (1966), 99--111; \ A.\,Mann,
 {\it Proc. Amer. Math. Soc.}, {\bf 53}, no.\,1 (1975), 262--264).
\hfill {\sl S.\,F.\,Kamornikov}

\vs  Yes, one can choose $k=3$ as such a constant, mod CFSG (V.\,I.\,Zenkov, {\it Siberian Math.~J.}, {\bf 62}, no.\,4 (2021), 621--637).
\emp

\bmp \textbf{19.40.}
Does Thompson's group $F$ (see 12.20) have the Howson property, that is, is the intersection of any two finitely generated subgroups of $F$ finitely generated?

\hfill {\sl I.\,Kapovich}

\vs
 No, it does not. Otherwise any subgroup of $F$ would have this property. But the wreath product $\Z\wr \Z$, which does not have the Howson property
(A.\,S.\,Kirkinski\u{\i}, {\it Algebra Logic}, {\bf 20}, no.\,1 (1981), 24--36), can be embedded in $F$ (S.\,Cleary, {\it Pacific J.~Math.}, {\bf 228}, no.\,1 (2006), 53--61). (D.\,Robertson, {\it Letter of 24 April 2018}.)
 \emp

\bmp \textbf{19.49.}
A {\it skew brace} is a set $B$ equipped with two operations $+$ and $\cdot$ such that $(B,+)$ is an additively written (but not necessarily abelian) group, $(B,\cdot )$ is a multiplicatively written group, and $a\cdot (b+c)=ab-a+ac$ for any $a,b,c\in B$.

 \makebox[15pt][r]{}Let $A$ be a skew brace with left-orderable multiplicative group. Is the additive group of $A$ left-orderable?
\hfill {\sl V.\,Lebed, L.\,Vendramin}

\vs
No, not always (T.\,Nasybullov,
{\it J.~Algebra}, {\bf 540} (2019), 156--167).
\emp

\bmp \textbf{19.50.} A finite graph is said to be integral if all eigenvalues of its adjacency matrix are integers.

 \makebox[25pt][r]{a)} Let $G$ be a finite group generated by a normal subset $R$ consisting of involutions. Is it true that the Cayley graph $Cay(G,R)$ is integral?

 \makebox[25pt][r]{b)} Let $A_n$ be the alternating group of degree $n$, let $S=\{(123),\, (124),\dots,(12n)\}$ and $R=S\cup S^{-1}$. Is it true that the Cayley graph $Cay(A_n,R)$ is integral?

\hfill {\sl D.\,V.\,Lytkina}

\vs  a) Yes, it is true (D.\,O.\,Revin, {\it Letter of 21 April 2018\,}; see also the reference for part (b) below; \ A.\,Abdollahi, {\it Letter of 3 May 2018}\,). Both proofs suggested are based on character theory; here is the second one. It suffices to show that the eigenvalues of $Cay(G,R)$ are rational, since the eigenvalues of a simple graph are algebraic integers.
It is known that every eigenvalue of $Cay(G,R)$ has the form $\theta_\chi = \frac{1}{\chi(1)} \sum_{r\in R} \chi(r)$ for some complex irreducible character $\chi$ of $G$ (implicit on pages 175--177 in P.\,Diaconis, M.\,Shahshahani, {\it Z.~Wahrscheinlichkeitstheorie Verw. Gebiete}, {\bf 57} (1981), 159--179, see also Theorem~9 in M.\,R.\,Murty, {\it J. Ramanujan Math. Soc.}, {\bf 18}, no.\,1 (2003) 1--20). Since the value of any complex character on an involution is an integer, it follows that $\theta_\chi$ is rational.

\vs  b) Yes, it is true
(W.\,Guo, D.\,V.\,Lytkina, V.\,D.\,Mazurov, D.\,O.\,Revin,
{\it Algebra Logic}, {\bf 58}, no.\,4 (2019), 297--305).
\emp

\bmp \textbf{19.55.}
Suppose that in a finite group $G$ every maximal subgroup $M$ is supersoluble whenever $\pi (M)=\pi (G)$, where $\pi (G)$ is the set of all prime divisors of the order of $G$.

\makebox[25pt][r]{a)} What are the non-abelian composition factors of $G$?

\makebox[25pt][r]{b)} Determine the exact upper bounds for the nilpotency length, the rank, and the $p$-length of $G$ if $G$ is soluble.
\hfill {\sl V.\,S.\,Monakhov}

\vs

 a) Every nonabelian finite simple group can occur as a composition factor of $G$ (A.\,Moret\'o, {\it Monatsh. Math.}, \textbf{195}, no.\,3 (2021), 497--500).
\vs

b) There is not any bound for the nilpotency length or the rank, but the \text{$p$-length} is at most 1 for every prime $p$ (A.\,Moret\'o, {\it Monatsh. Math.}, \textbf{195}, no.\,3 (2021), 497--500).
\emp

\bmp \textbf{19.67.}
 Let $G\le {\rm Sym}(\Omega)$, where $\Omega$ is finite. The $2$-closure $G^{(2)}$ of the group $G$ is defined to be the largest subgroup of ${\rm Sym}(\Omega)$ containing~$G$ which has the same orbits as~$G$ in the induced action on~$\Omega\times\Omega$. Is it true that if $G$ is solvable, then every composition factor of $G^{(2)}$ is either a cyclic or an alternating group?

 \hfill {\sl I.\,Ponomarenko}

 \vs  No, it is not true (S.\,V.\,Skresanov,
 {\it Algebra Logic}, {\bf 58}, no.\,3 (2019), 249--253).
\emp

\bmp \textbf{19.75.}
Let $P$ be a finite 2-group of exponent $2^e$ such that the rank of every abelian subgroup is at most $r$. Is it true that $|P|\le 2^{r(e+1)}$?
This bound would be sharp (for a direct product of quaternion groups).
\hf {\sl B.\,Sambale}

\vs  No, it is not true, as follows from the examples constructed in (A.\,Yu.\,Olshanskii, {\it Math. Notes}, {\bf 23} (1978), 183--185) (A.\,Mann, {\it Letter of 23 April 2018}).
\emp

\bmp \textbf{19.80.}
For a periodic group $G$, let $\pi _e(G)$ denote the set of orders of elements of~$G$. A periodic group $G$ is said to be an {\it $OC_n $-group\/} if $\pi _e(G)=\{ 1,\, 2, \ldots ,\, n\}$.
Is it true that every $OC_7 $-group is isomorphic to the alternating group $A_7$?

 \makebox[15pt][r]{}For finite groups, the answer is affirmative.
\hfill {\sl W.\,J.\,Shi}

\vs  Yes, it is true (E.\,Jabara, A.\,S.\,Mamontov, \emph{Siberian Math.~J.}, \textbf{62}, no.\,1 (2021), 94--104).
\emp

\bmp \textbf{19.81.}
(Well-known problem). Is the conjugacy problem in the braid group $B_n$ in the class NP (that is, decidable in nondeterministic polynomial time with respect to the maximum of the lengths $|u|$, $|v|$, where $u, v$ are the input braid words)? A~stronger question: given two conjugate elements of $B_n$ represented by braid words of lengths $\leq m$, is there a conjugator whose length is bounded by a polynomial function of $m$?

\hfill {\sl V.\,Shpilrain}

\vs Yes, there is such a conjugator, even for any mapping class groups (J.\,Tao,
\emph{Geom. Funct. Anal.}, \textbf{23}, no.\,1 (2013), 415--466).
\emp

\bmp \textbf{19.84.}
Let $\Bbb{P}$ be the set of all primes, and let
 $\sigma =\{\sigma_{i} \mid i\in I\}$ be some partition of
$\Bbb{P}$ into disjoint subsets. A finite group $G$ is said to be \emph{$\sigma$-primary} if $G$ is a $\sigma_{i}$-group for some $i$; \emph{$\sigma$-nilpotent} if $G$ is a direct product of $\sigma$-primary groups;
 \emph{$\sigma$-soluble} if every chief factor of $G$ is $\sigma$-primary.
A subgroup $A$ of a finite group $G$ is said to be \emph{${\sigma}$-sub\-normal in $G$} if there is a chain $A=A_{0} \leq A_{1} \leq \cdots \leq
A_{n}=G$ such that for every $i$ either $A_{i-1} \trianglelefteq A_{i}$ or
$A_{i}/(A_{i-1})_{A_{i}}$ is ${\sigma}$-primary, where $(A_{i-1})_{A_{i}}$ is the largest normal subgroup of $A_{i}$ contained in $A_{i-1}$.

 \makebox[15pt][r]{}Suppose that a subgroup $A$ of a finite group $G$ is
${\sigma}$-subnormal in $\langle A, A^{x} \rangle$ for all $x\in G$. Is it true that then $A$ is ${\sigma}$-subnormal in $G$?

\makebox[15pt][r]{}An affirmative answer is known if $\sigma =\{\{2\}, \{3 \}, \ldots \}$ (Wielandt).
\hfill {\sl A.\,N.\,Skiba}

\vs
 No, not always. For example, in $G = S_5$ with partition $\sigma = \{2, 3\}\cup \{5\}$ the subgroup $A =\langle (12) \rangle$ is $\sigma$-subnormal in $\langle A, A^x\rangle$ for every $x\in G$, but it is not $\sigma$-subnormal in~$G$. (V.\,N.\,Tyutyanov, {\it Letter of 28 August 2019}.)
\emp

\bmp \textbf{19.85.}
Suppose that every Schmidt subgroup of a finite group $G$ is ${\sigma}$-subnormal in $ G$ (see Archive 19.84). Is it true that then there is a normal $\sigma$-nilpotent subgroup $N\leq G$ such that $G/N$ is cyclic?

\makebox[15pt][r]{}An affirmative answer is known if $\sigma =\{\{2\}, \{3 \}, \ldots \}$.
\hfill {\sl A.\,N.\,Skiba}

\vs  Yes, it is true (X.\,Yi, S.\,F.\,Kamornikov, \emph{J.~Algebra}, \textbf{560} (2020), 181--191).
\emp

\bmp \textbf{19.87.}
Suppose that for every Sylow subgroup $P$ of a finite group $G$ and every maximal subgroup $V$ of $P$ there is a
${\sigma}$-soluble subgroup $T$ such that $VT=G$. Is it true that then $G$ is ${\sigma}$-soluble?
\hfill {\sl A.\,N.\,Skiba}

\vs
Yes, it is true (A.-M.\,Liu, W.\,Guo, I.\,N.\,Safonova, A.\,N.\,Skiba, {\it J.~Algebra}, {\bf 585} (2021), 280--293); another solution using CFSG is in (S.\,F.\,Kamornikov, V.\,N.\,Tyutyanov, \textit{Trudy Inst. Mat. Mekh. UrO RAN}, \textbf{27}, no.\,1 (2021), 98--102 (Russian); \ X.\,Yi, S.\,F.\,Kamornikov, V.\,N.\,Tyutyanov, \textit{Probl. Physics, Math. Techn.}, \textbf{46}, no.\,1 (2021), 50--53).
\emp

\bmp \textbf{19.88.}
Suppose that for every Sylow subgroup $P$ of a finite group $G$ and every maximal subgroup $V$ of $P$ there is a
${\sigma}$-nilpotent subgroup $T$ such that $VT=G$. Is it true that then $G$ is ${\sigma}$-nilpotent?
\hfill {\sl A.\,N.\,Skiba}

\vs  Yes, it is true (A.-M.\,Liu, W.\,Guo, I.\,N.\,Safonova, A.\,N.\,Skiba, {\it J.~Algebra}, {\bf 585} (2021), 280--293); another solution using CFSG is in (S.\,F.\,Kamornikov, V.\,N.\,Tyutyanov, \textit{Trudy Inst. Mat. Mekh. UrO RAN}, \textbf{27}, no.\,1 (2021), 98--102 (Russian); \ X.\,Yi, S.\,F.\,Kamornikov, V.\,N.\,Tyutyanov, \textit{Probl. Physics, Math. Techn.}, \textbf{46}, no.\,1 (2021), 50--53).
\emp

\bmp \textbf{19.90.}
A {\it skew brace} is a set $B$ equipped with two operations $+$ and $\cdot$ such that $(B,+)$ is an additively written (but not necessarily abelian) group, $(B,\cdot )$ is a multiplicatively written group, and $a\cdot (b+c)=ab-a+ac$ for any $a,b,c\in B$.

 \makebox[25pt][r]{a)} Is there a skew brace with soluble additive group but non-soluble multiplicative group?

 \makebox[25pt][r]{d)}
 Is there a finite skew brace with non-soluble additive group but nilpotent multiplicative group? \hfill \mbox{\sl A.\,Smoktunowicz, L.\,Vendramin}

\vs
a) Yes, there is (T.\,Nasybullov, {\it J.~Algebra}, {\bf 540} (2019), 156--167).
\vs

d) No, there is not (C.\,Tsang, Q.\,Chao, \textit{ Int. J. Algebra Comput.}, \textbf{30}, no.\,2 (2020), 253--265).
\emp

\bmp \textbf{19.98.}
A connected graph $\Sigma$ is a {\it symmetrical extension} of a graph
$\Gamma$ by a graph $\Delta$ if there exist a vertex-transitive group $G$ of automorphisms of $\Sigma$ and an imprimitivity system $\sigma$ of $G$ on the set of vertices of $\Sigma$ such that the quotient graph $\Sigma/\sigma$ is
isomorphic to $\Gamma$ and blocks of $\sigma$ generate in $\Sigma$
subgraphs isomorphic to $\Delta$.

 \makebox[15pt][r]{}{(a)} Let $\Gamma$ be a locally finite Cayley graph of a finitely presented group, and $\Delta$ a~finite
graph. Are there only finitely many symmetrical extensions of $\Gamma$
by $\Delta$?

 \makebox[15pt][r]{}{(b)} Let $\Gamma$ be a locally finite graph which has the property of $k$-contractibility for some positive integer $k$ (see the definition
in (V.\,I.\,Trofimov, {\it Proc. Steklov Inst. Math.},
\textbf{279}, suppl. 1 (2012), 107--112); note that any $\Gamma$ from
(a) is such a graph) and let $\Delta$ be a finite graph. Are there only
finitely many symmetrical extensions of $\Gamma$ by~$\Delta$?
\hfill {\sl V.\,I.\,Trofimov}

\vs  (a) No, not always, as follows from the construction in the proof of Theorem~H in (M.\,de\,la\,Salle, R.\,Tessera, \emph{J.~Topology}, \textbf{12} (2019), 705--743).

\vs (b) No, not always as follows from the construction in the proof of Theorem~H in (M.\,de\,la\,Salle, R.\,Tessera, \emph{J.~Topology}, \textbf{12} (2019), 705--743).
\emp

\bmp \textbf{19.101.}
The maximum length of a chain of nested centralizers of a group is called its $c$-dimension. Let $G$ be a locally finite group of finite $c$-dimension $k$, and let $S$ be the preimage in $G$ of the socle of $G/R$, where $R$ is the locally solvable radical of~$G$. Is it true that the factor group $G/S$ contains an abelian subgroup of index bounded by a function of~$k$?
\hfill {\sl A.\,V.\,Vasil'ev}

\vs  Yes, it is true (A.\,A.\,Buturlakin,
{\it J. Algebra Appl.}, {\bf 18}, no.\,12 (2019), 1950223, 12~pp).
\emp

\bmp \textbf{19.109.}
A subgroup $H$ of a finite group $G$ is called pronormal if for any $g\in G$ the subgroups $H$ and $H^g$ are conjugate by an element of $\langle H,H^g\rangle$. A maximal subgroup of a maximal subgroup is called second maximal.
Is it true that in a non-abelian finite simple group $G$ all maximal subgroups are Hall subgroups if and only if every second maximal subgroup of $G$ is pronormal in $G$?
\hfill {\sl V.\,I.\,Zenkov}

\vs  No, not always, for example, in $SL_2(2^{11})$ every second maximal subgroup is pronormal (V.\,N.\,Tyutyanov, {\it Letter of 23 November 2018}).
\emp

\newpage

\columnsep=0.7cm
\pagestyle{myheadings}
\markboth{\protect\vphantom{(y} Index of
Names (\textbf{A:} Archive)}{\protect\vphantom{(y} Index of
Names (\textbf{A:} Archive)}
\thispagestyle{headings}

\twocolumn[~
\vskip1ex \centerline {\Large \textbf{Index of names (\textbf{A:} Archive)}}
\phantomsection\label{index}
\vskip4ex

]

\textbf{A}basheev B.\,L.~ \textbf{A:}~{\it 7.30}

 Abdollahi A.~ {\it 14.53},\, 15.1,\, 16.1,\, {\it 16.1},\lb17.1,\, 17.4--11,\, {\it 20.1},\, 21.1,\, \textbf{A:}~{\it 9.50},\lb{\it 11.46},\, 16.1,\, {\it 16.61},\, {\it 16.83},\, 17.2--3,\lb{\it 17.2},\, 17.12,\, {\it 18.57},\, {\it 19.50}

 Adamov S.\,N.~ \textbf{A:}~10.1

Adeleke S.\,A.~ 14.8

 Adyan S.\,I.~ 4.2,\, 4.5,\, 7.3, \textbf{A:}~{\it
1.24},\, {\it 1.63},\lb{\it 1.82,\, 2.13,\,
2.15,\, 2.39,\, 3.9},\, 4.1, \lb4.3--5,\, 7.1--2,\,
{\it 7.2},\, {\it 13.34}

 Agalakov S.\,A.~ \textbf{A:}~{\it 2.19}

Agrachev A.~ 21.123

Ahanjideh M.~ \textbf{A:}~{\it 12.38}

Ahanjideh N.~ \textbf{A:}~{\it 12.38}

 Akbari S.~ 16.1,\, \textbf{A:}~16.1

Akhlaghi Z.~ {\it 20.78},\, {\it 20.79}

Akhmedov N.\,S.~ 4.42

Alavi S.\,H.~ \textbf{A:}~{\it 12.38}

 Aleev R.\,Zh.~ 11.1,\, 11.3,\, 12.1,\, 13.1,\lb14.2--3,\, \textbf{A:}~11.2,\, {\it 12.1},\, 14.1

Allenby R.\,B.\,J.\,T.~ 9.1,\, \textbf{A:}~9.2

 Alperin J.~ 4.24,\, 16.33,\,
\textbf{A:}~{\it 3.64}, {\it 4.76},\lb{\it 8.49},\, 15.24,\, 17.46

 Alperin R.~ 13.39,\, \textbf{A:}~16.2

Al-Sharo Kh.\,A.~ \textbf{A:}~{\it 18.91}

Amberg B. 7.56,\, 13.27,\,  18.2--5,\, \textbf{A:}~{\it 18.95}

Amelio M.~ {\it 10.64}

 Amir G.~ \textbf{A:}~{\it 16.74}

 Amiri M.~ {\it 18.1},\, 21.2

Amitsur S.\,A.~ \textbf{A:}~{\it 16.54}

Anagnostopoulou-Merkouri M.~ 21.3--4

Anashin V.\,S.~ 10.2--5

Andreadakis S.~ \textbf{A:}~{\it 10.33}

 Angel O.~ \textbf{A:}~{\it 16.74}

 Anokhin M.\,I.~ 13.2,\, \textbf{A:}~{\it 6.14},\, {\it
8.26},\, {\it 8.28}, \lb{\it 12.98}

 Anosov D.\,V.~ 8.8,\, \textbf{A:}~8.8

Antol\'{i}n Y.~ {\it 8.11}

Antonov V.\,A.~ 6.1--2,\, \textbf{A:}~{\it 10.1,\,} 11.4

Apps A.\,B.~ 21.39

 Arad Z.~ 16.3--4,\, {\it 16.3},\, \textbf{A:}~{\it 3.39,\, 7.12}

Arezoomand M.~ , {\it 20.1},\, 20.1--3

 Arnautov V.\,I.~ 10.8,\, 14.4--5,\, {\it 14.4},\, {\it 14.5},\lb\textbf{A:} 10.6,\,
{\it 10.6},\, 10.7,\, {\it 10.7},\,
12.2,\lb{\it 12.2}

Artamonov V.\,A.~ 4.6--7,\, 11.5,\, 13.3

Artemovich O.\,D.~ \textbf{A:}~{\it 4.21\/}

Artin E.~ \textbf{A:}~1.11,\, 17.14

 Asar A.\,O.~ 16.5--6,\, 17.13,\, 21.5--6,\lb\textbf{A:}~{\it 4.36,\, 14.86}

 Aronszajn N.~ 2.48,\, 18.71

 Arzhantseva G.\,N.~ {\it 14.6},\, 18.6,\, \textbf{A:}~{\it 15.4},\lb{\it 11.75}

 Aschbacher M.~ 11.3,\, 11.60,\, 19.91,\lb\textbf{A:}~{\it 2.37,\,
5.41}

Ashmanov I.\,S.~ \textbf{A:}~{\it 2.45},\, {\it 17.79}

Aslanyan H.\,T.~ {\it 17.70}

Atabekyan V.\,S.~ {\it 11.36},\, {\it 17.70},\, {\it 17.71},\, \lb 18.7,\, \textbf{A:}~{\it 1.83},\, {\it 7.1},\, {\it 8.53},\, {\it 17.79}

Auslander L.~ \textbf{A:}~{\it 1.34}

 Avenhaus J.~ \textbf{A:}~{\it 9.30}

Azarian M.\,K.~ 19.1--6

Azarov D.\,N.~ \textbf{A:}~{\it 11.10}

~

 \textbf{B}abai L.~ 17.41

B{\"a}chle A.~ 20.4--6

Bachmuth S.~ \textbf{A:}~{\it 5.49}

 Baer R.~ 4.17--19,\, 16.77,\, 19.13,\, 20.93,\lb\textbf{A:}~{\it
2.88},\, 4.16,\, 11.11,\, {\it 11.11}

Bagayoko V.~ {\it 11.67}

 Bagi\'nski C.~ \textbf{A:}~11.6

Bahri A.~ {\it 20.79}

Bak A.~ \textbf{A:}~{\it 7.40}

Bakhturin Yu.\,A.~ \textbf{A:}~{\it 2.41}

Ballester-Bolinches A.~ {\it 14.28},\, {\it 17.39},\lb{\it 17.56},\, 18.8,\, 18.10,\, 21.7,\, \textbf{A:}~12.96,\lb18.9

Ban G.~ \textbf{A:}~{\it 12.78}

Banach S.~ 10.12,\, 19.59

Banakh T.~ \textbf{A:}~{\it 16.42}

 Bandman T.~ \textbf{A:}~{\it 15.75}

 Baranov D.\,V.~ 18.46

Barbosa-Fernandes R.~ {\it 20.72}

 Bardakov V.\,G.\,\, 14.11--12,\, 14.14--17,\lb15.9,\, 15.11,\, 16.7,\, 16.9--10,\, 17.15--16,\lb{\it 17.14},\, 18.11--14,\, 19.7--9,\, 20.124--126,\lb21.8\, \textbf{A:}~{\it 11.20},\, {\it 11.54},\, {\it 11.55},\, 14.13,\lb{\it 14.27},\, 15.10,\, {\it 15.10},\, {\it 15.60},\, 16.8,\lb17.14,\, 18.15,\, {\it 18.86}

Barge J.~ \textbf{A:}~{\it 14.13}

Barker L.~ \textbf{A:}~{\it 12.80}

Barnea Y. 19.10,\, {\it 19.10},\, 20.102,\, 21.9--10

Bartholdi L.~ 15.12--14,\, 15.16,\, 15.19,\lb\textbf{A:}~{\it 12.22},\, 15.14--15,\, 15.17--18

 Bass H.~ 10.40,\, 12.1,\, 12.9,\, 13.39,\, 19.27,\lb\textbf{A:}~12.1

 Bauer A.~ 20.33

Baumann B.~ 13.4,\, \textbf{A:}~{\it 2.1,\, 3.27\/}

Baumeister B.~ {\it 19.11}

 Baumslag B.~ \textbf{A:}~{\it 4.39}

 Baumslag G.~ 4.9,\, 4.11,\, 5.16,
\,6.5,\, 8.2,\lb13.39,\, 14.18--22,\,
16.65,\, 17.90,\, 17.99,\lb21.47,\,  \textbf{A:}~2.72,\, {\it
3.15},\, 4.8,\, 4.10,\lb13.39,\, 15.86,\, 16.2

Bayes A.\,J.~ 20.36

Baykalov A.\,A.~ {\it 17.41}

 Bazhenova G.\,A.~ 18.85

Beek van M.~ {\it 11.34}

Bekka M.~ 14.6,\, \textbf{A:}~15.7

Belegradek O.~ 18.16,\, 18.18--19,\, \textbf{A:}~{\it 4.20},\lb{\it 15.86},\, 18.17

Belegradek I.~ \textbf{A:}~{\it 16.52}

 Beletski\u{i} P.\,M.~ \textbf{A:}~{\it 4.38}

Belk J.~ 20.7,\, \textbf{A:}~{\it 14.10}

Bell S.\,D.~ 21.66,\, \textbf{A:}~{\it 9.3}

 Bellingeri P.~ \textbf{A:}~{\it 17.14}

Belonogov V.\,A.~ 18.20,\, \textbf{A:}~2.1

Belousov I.\,N.~ {\it 20.123}

Belousov V.\,D.~ \textbf{A:}~2.2

 Belov A.\,Ya.~ \textbf{A:}~{\it 2.39--40}

Belov Yu.\,A.~ \textbf{A:}~{\it 3.53},\, {\it 10.56}

Beltr\'an A.~  {\it 19.57},\, {\it 19.58}

Belyaev V.\,V.~ 5.1,\, 7.5,\,
 11.109,\, 13.5--9,\lb15.20--23,\, \textbf{A:}~{\it 1.75},\,
{\it 4.35},\, 5.1--2,\, 7.4,\lb7.6,\,
{\it 7.6},\, 9.3

Bencs\'ath K.~ 11.7

 Benesh B.~ \textbf{A:}~{\it 12.82}

 Ben-Ezra D.\,E.-C.~ {\it 14.42}

 Bensa\"{\i}d A.~ \textbf{A:}~{\it 7.48}

 Benson D.\,J.~ \textbf{A:}~{\it 12.36}

Berger M.\,A.~ 21.93

Berger T.~ 5.30,\, 15.56

 Bergman G.\,M.~ 6.39,\, 16.88--89,\, 17.18,\lb17.22--23,\, 17.25,\, 18.21--22,\, 19.110,\lb20.8--15,\, 21.11--18,\, \textbf{A:}~{\it 2.29},\, {\it 2.71},\lb5.3,\, 16.50,\, 17.17,\, 17.19--21,\, {\it 17.21},\lb17.24,\, {\it 17.24},\, {\it 19.35}

 Berkovich Ya.\,G.~ 4.13,\, 11.8--9,\, 15.26,\lb15.28--32,\, 16.11,\, 16.14,\, 21.51,\, \textbf{A:}~2.3,\lb{\it
3.28},\, 4.12,\, 11.8,\, 15.24--25,\, 15.27,\lb15.33,\, 16.12--13

Berlatto A.~ 21.42

Berman S.\,D.~ 14.2,\, \textbf{A:}~3.2

 Berrick A.\,J.~ {\it 10.40}

Bertram E.~ 19.11

Bessenrodt C.~ {\it 21.59}

Bestvina M.~ {\it 14.15},\, \textbf{A:}~{\it 6.35}

Bettio E.~ \textbf{A:}~{\it 19.36}

Beyer\,de\,Ryke B.~ {\it 21.88}

Bezerra M.~ 19.12

Bican L.~ \textbf{A:}~{\it 3.30},\, {\it 3.31}

Bieri R.~ 6.3,\, 6.5,\, \textbf{A:}~6.4,\, 6.35

Bigelow S.~ 17.16

 Bildanov R.\,R.~ {\it 20.37}

Birman J.\,S.~ 14.12,\, \textbf{A:}~10.24

Blackburn S.\,R.~ \textbf{A:}~{\it 17.35}

Blass A.~ \textbf{A:}~17.24

Bloshchitsyn V.\,Ya.~ \textbf{A:}~{\it 5.34}

 Bludov V.\,V.~ 11.10,
\,13.57,\,
 16.15--16,\lb16.18--19,
 \,17.26--27,\, 17.30--31,\, 21.79,\lb\textbf{A:}~{\it
1.47},\, {\it 1.49},\, {\it 1.52},\, {\it 1.60},\, {\it 3.17},\, 5.4,\lb{\it 5.4},\, {\it 5.23},\, {\it
6.53},\, {\it 8.56}, {\it 10.66},\,
11.10,\lb{\it 11.68},\, {\it 12.14},\, 13.10--11,\, {\it 13.16},\, {\it
13.54},\lb{\it 14.25},\, {\it 14.71},\, 15.34,\, {\it 15.34},\, 16.15,\lb16.17,{\it 16.17},\, {\it 16.35,\, 16.65},\, 17.28--29,\lb{\it 17.29},\, {\it 18.15}

Bogachev N.\,V.~ {\it 17.94},\, {\it 18.87}

Bogan Yu.\,A.~ \textbf{A:}~1.1,\, 2.4

 Bogopolski O.\,V.~ 13.12,\, 14.23,\, 14.24,\lb17.32,\, \textbf{A:}~{\it 10.25},\, {\it 10.26},\, {\it 10.69},\, {\it 13.61},\lb15.35

 Bokut' L.\,A.~ 1.3,\, 1.5--6,\, 3.3,\, \textbf{A:}~1.2,
\lb1.4,\, {\it 1.4},\, 3.4

Bondarenko I.\,V.~ \textbf{A:}~{\it 15.15}

 Bonner T.~ \textbf{A:}~{\it 16.8}

Bonvallet K.~ 12.29

Boone W.~ 20.7

Borcherds R.\,E.~ 15.55

Borel A.~ 19.79,\, \textbf{A:}~2.51

 Borevich Z.\,I.~ \textbf{A:}~14.50

 Borovik A.\,V.~ 11.12,\,
11.109,\, 19.13,\lb20.16,\, 21.19--21,\, \textbf{A:}~{\it
1.75},\, {\it 7.6},\, 11.11,\lb11.13

 Borshchev A.\,V.~ \textbf{A:}~{\it 3.33}

 Bou-Rabee K.~ \textbf{A:}~{\it 15.35}

Bousfield A.\,K.~ 17.89

 Bovdi A.\,A.~ \textbf{A:}~12.1

Bowtell A.~ \textbf{A:}~{\it 1.4}

Brachter J.~  21.24

Bradford H.~ 21.22

Brady N.~ 21.126,\, \textbf{A:}~{\it 6.35}

Brameret M.-P.~ \textbf{A:}~{\it 2.47}

Brandl R.~ 11.14--18,\,
12.3--4,\, 13.13

Brauer R.~ 9.23,\, 11.98,\, 20.93,\, 21.121,\lb\textbf{A:}~{\it 3.64},\, 4.14

 Bray J.~ \textbf{A:}~{\it 15.43},\, {\it 15.75},\, {\it 16.55}

Brenner J.\,L.~ 10.32

Brescia M.~  {\it 6.30},\, {\it 11.9},\, {\it 13.23},\, {\it 19.48}

Breuer T.~ 21.25

 Breuillard E.~ 18.108,\, \textbf{A:}~{\it 15.7}

 Bridson M.~ {\it 12.41},\, 16.90,\, 21.125,\, 21.145

 Brin M.~ 21.144

 Brodski\u{\i} S.\,D.~ 11.62,\, \textbf{A:}~{\it 1.2},\, {\it
4.10}

Bromberg K.~ {\it 14.15}

Brown K.\,S.~ \textbf{A:}~{\it 6.4},\, {\it 16.52}

Bruhat F.~ \textbf{A:}~2.52,\, {\it 2.52}

Bruijn de N.\,G.~ 17.22

 Brumberg N.\,R.~ 3.5

 Brunner A.\,M.~ 15.19,\, 15.92

 Bryant R.~ 10.5,\, 16.92,\, \textbf{A:}~{\it 11.47}

Bryce R.\,A.~ \textbf{A:}~{\it 10.56}

 Bryukhanov O.\,V.~ 17.16,\, \textbf{A:}~{\it 18.86}

Bryukhanova E.\,G.~ {\it 7.31}

 Bubboloni D.~ 20.17,\, 21.23--24,\, \textbf{A:}~18.23

 Budkin A.\,I.~
10.10,\, 13.14--15,\, 14.26,\lb15.36,\, 16.20,\, 17.33--34,\, 19.14--15,\lb20.18, \textbf{A:}~{\it 2.40},\, {\it 3.15},\, {\it
3.52},\,{\it 5.51},\, 10.9,\lb12.5,\, 14.25,\, {\it 18.17}

Bugeaud Y.~ 13.65

 Burger M.~ 18.24,\, \textbf{A:}~{\it 4.45}

Burichenko V.\,P.~ 18.25--28,\, \textbf{A:}~{\it 9.59--60},\lb{\it 12.49},\, {\it 16.31},\, {\it 16.58},\, {\it 17.50}

Burness T.\,C.~ {\it 15.40},\, {\it 17.41},\, {\it 17.42},\lb{\it 20.121},\, 21.3--4,\, 21.25--29,\, 21.38

 Burns R.\,G.~ 2.82,\, 15.92,\, \textbf{A:}~7.7

Burnside W.~ 15.2,\, 20.52,\, 20.93

Busetto G.~ 12.6

 Button J.\,O.~ 19.16--17,\, 21.30,\, \textbf{A:}~{\it 16.2}

 Buttsworth R.\,N.~ 16.51,\, \textbf{A:}~{\it 2.25}

 Buturlakin~A.\,A.~ {\it 3.12},\, {\it 15.53},\, {\it 18.38},\lb\textbf{A:}~{\it 16.25},\, {\it 18.112},\, {\it 19.101}

Byott N.~ 21.31

~

 \textbf{C}abanes M.~ {\it 8.51}

Cai T.\,W.~ {\it 16.48}

Calegari D.~ 17.107,\, 18.24,\, 18.39--42,\lb21.127

Cameron P.\,J.~ 9.69--70,\, 11.45,\, 13.30--31,\linebreak
\makebox[3ex][r]{}16.45,\, 19.18--20,\,  21.24,\, 21.32,\, 21.95,\lb\textbf{A:}~13.28--29

Camina A.~ 14.46,\, {\it 14.46}

Camina R.~
\textbf{A:}~{\it 12.24},\, {\it 13.36}

Cannon J.~ 19.42

Cannonito F.\,B.~ 5.15--16,\, 7.19,\, 7.21

Cant A.~ {\it 14.64}

Cantat S.~ 19.78

Cao Zhenfu~ 13.65

Capdeboscq I.~ {\it 21.142}

 Caprace P.-E.~ 18.70,\, 19.21--22,\, 19.70,\lb19.73,\, 20.19,\, \textbf{A:}~{\it 14.13}

Caranti A.~ 11.46,\, {\it 15.32}, 20.108,\lb \textbf{A:}~11.46

 Carlisle D.~ \textbf{A:}~12.36

 Carmichael R.\,D.~ 18.90

Carette M.~ {\it 13.12}

 Cartan E.~ \textbf{A:}~{\it 16.58}

Carter R.\,W.~ \textbf{A:}~8.13

Casolo C.~ {\it 14.65},\, \textbf{A:}~{\it 9.64}

Castel F.~ {\it 14.102}

Castellano I.~ {\it 19.22},\, 21.33,\, 21.37

Cat G.~ {\it 15.26}

Cavenagh N.\,J.~ \textbf{A:}~17.35

 Chalk C.\,P.~ 20.20,\, \textbf{A:}~{\it 8.10}

 Chang  A.~ {\it 21.137}

Chankov E.\,I.~ \textbf{A:}~{\it 9.56},\, {\it 11.94}

Chao Q.~ \textbf{A:}~{\it 19.90}

Charin V.\,S.~ 8.60,\, 8.86,\, \textbf{A:}~3.58

Chase S.~ \textbf{A:}~2.4

 Chatterji I.~ {\it 10.40},\, 18.111

 Chechin S.\,A.~ \textbf{A:}~6.49,\, {\it 6.49}

Chen J.~ {\it 4.13},\, {\it 21.51}

Chen  L.~ {\it 5.52},\, {\it 14.102}

Chen G.\,Y.~ \textbf{A:}~{\it 12.38},\, {\it 16.1}

Chen W.\,Y.~ 21.82

Cherep A.\,A.~ \textbf{A:}~{\it 6.57}

Cherepanov E.\,A.~ {\it 11.36}

Chermak A.\,L.~ 9.72

 Chernikov N.\,S.~ 6.48,\, 7.54--56,\, 11.114,\linebreak
\makebox[3ex][r]{}13.60,\, 14.97,\, 15.97,\, 16.105,\, 17.120,\lb18.3,\,
\textbf{A:}~{\it 2.69},\, {\it 5.19},\, {\it
5.62}

 Chernikov S.\,N.~ 2.80,\, 15.97,\, \textbf{A:}~5.62

Chevalley C.~ 2.56,\, 3.42

Chiodo M.~ \textbf{A:}~19.23--24

Chou C.~ 8.9

 Chunikhin S.\,A.~ \textbf{A:}~5.63

 Churkin V.\,A.~ 7.57,\, 10.71,\, 14.98,\, 18.36,\,
 \lb\textbf{A:}~{\it 2.31},\, {\it 3.11},\, {\it 7.36},\,
7.57,\, 13.61,\lb{\it 15.8},\, {\it 16.75}

Ciobanu L.~ 20.41

Clare F.~ \textbf{A:}~{\it 1.68}

 Clark A.~ \textbf{A:}~{\it 15.63}

 Clay A.~ {\it 16.48},\, {\it 16.51},\, 21.34

Cleary S.~ \textbf{A:}~{\it 19.40}

 Cohen A.\,M.~ 17.16

Cohen D.\,E.~ \textbf{A:}~{\it 1.85},\, 5.28

Cohn P.\,M.~ \textbf{A:}~{\it 2.38}

Collatz L.~ 18.47

Collins D.\,J.~ 10.70,\, \textbf{A:}~5.21--22,\, {\it 16.75}

 Comerford L.~ \textbf{A:}~14.48

Comfort W.\,W.~ 13.48

 Conder M.~ 14.73,\, 16.46,\, 20.21,\, \textbf{A:}~14.49,\lb{\it
14.49}

 Conrad J.~ 17.62

 Conrad P.~ 16.47

Contreras Rojas Y.~ 21.35

Conversano A.~ {\it 15.8}

Cook G.\,C.~ 21.36--37

Cooper C.\,D.\,H.~ 6.47

 Corner A.\,L.\,S.~ \textbf{A:}~{\it 17.21}

 Cornulier Y.~ 20.42,\, \textbf{A:}~{\it 16.52},\, {\it 16.59}

 Corson S.\,M.~ {\it 9.41},\, {\it 9.42},\, {\it 15.70},\, {\it 15.80},\lb{\it 18.36},\, {\it 19.74},\, {\it 20.10},\, {\it 20.13},\, {\it 20.40},\lb {\it 21.12},\, {\it 21.15},\,  21.38,\, \textbf{A:}~{\it 3.22}

Cosme-Ll\'opez E.~ {\it 4.13}

 Cossey J.~ {\it 17.39},\,
\textbf{A:}~{\it 5.17},\, {\it 14.77},\, {\it 16.13}

 Cossidente A.~ \textbf{A:}~{\it 10.28}

Coulbois T.~ \textbf{A:}~{\it 14.48}

Coulon R.~ {\it 6.45},\, {\it 17.62}

Cowling M.~ 14.6

 Coxeter H.\,S.\,M.~ 18.90

 Craven D.~ \textbf{A:}~{\it 16.67},\, {\it 18.115}

 Crawley-Boevey W.\,W.~ \textbf{A:}~{\it 12.36}

 Crestani E.~ \textbf{A:}~{\it 15.33}

 Crosby P.\,G.~ \textbf{A:}~{\it 17.12}

 Curran M.\,J.~ \textbf{A:}~12.78

Curtin B.~ 19.25,

 Curtis C.~ 16.39

 Cutolo G.~ 18.54--56,\, 21.98

Cuypers H.~ \textbf{A:}~{\it 12.76}

Czerniakiewicz A.\,J.~ {\it 3.3}

~

 \textbf{D}ade E.\,C.~ 18.102

Dahmani F.~ {\it 19.41}

 Dalla Volta F.~ \textbf{A:}~{\it 11.27}

Damian E.~ \textbf{A:}~{\it 11.104}

Daneshkhah A.~ \textbf{A:}~{\it 12.38}

Dantas A.\,C.~ {\it 15.19},\, 21.39--42

 Darafsheh M.\,R.~ \textbf{A:}~{\it 12.38},\, 16.36--37,\linebreak
\makebox[3ex][r]{}{\it 16.37}

Darbinyan A.~ \textbf{A:}~{\it 5.21},\, {\it 17.28}

Dardano U.~ \textbf{A:}~{\it 9.64}

Dark R.\,S.~ 17.47,\, \textbf{A:}~{\it 2.21}

Davis C.~ \textbf{A:}~{\it 11.2}

Day M.\,M.~ \textbf{A:}~8.6

Deaconescu M.~ 21.43,\, \textbf{A:}~15.43

De~Bruijn N.\,G.~ 17.22

DeCaro J.~ {\it 21.87}

De Falco M.~ {\it 16.38}

 De Giovanni F.~ 13.22--23,\, 14.36,\, 16.38

De Groot~ \textbf{A:}~2.36

De la Harpe P.~ 14.6--10,\, 15.8,\, \textbf{A:}~{\it 5.2},\lb14.10,\, 15.4--8

De la Salle M.~ {\it 12.20},\, {\it 15.42},\, \textbf{A:}~{\it 19.98}

Dekimpe K.~ 19.28

De Melo E.~ 21.39--40

Dennis R.\,K.~ 14.15

 Deryabina G.\,S.~ \textbf{A:}~{\it 1.81},\, {\it 13.10},\,
{\it 13.34}

 Detomi E.~ {\it 16.53},\, 20.49,\, 21.122,\lb\textbf{A:}~{\it 11.104}

 Diaconis P.~ \textbf{A:}~{\it 19.50}

Dicks W.~ 19.26,\, \textbf{A:}~{\it 5.28},\, {\it 11.110}

Diekert V.~ {\it 16.9},\,  \textbf{A:}~{\it 9.25}

Dietrich H.~ {\it 14.64}

 Digne F.~ 17.16

Di Martino L.~ {\it 11.32}

Di Matteo M.~ {\it 20.109}

Di~Siena B.\,G.~ {\it 15.50},\, {\it 20.56},\, {\it 20.57},\lb {\it 20.73},\, {\it 21.96}

Dixmier J.~ 20.62

 Dixon J.\,D.~ 6.10--11,\,
18.36

 Dixon M.~ {\it 6.9},\, 18.37

Doerk K.~ 14.30

Doktorov I.\,P.~ \textbf{A:}~4.27

 Dolbak L.\,V.~ \textbf{A:}~{\it 11.88}

Dolfi S.~ {\it 14.65}

Dolorfino M.~ {\it 20.79}

Dona D.~ {\it 20.23}

 Donkin S.~ 15.44,\, 16.39

Donoven C.~ 21.38

Doorn van W.~ {\it 3.46},\, {\it 18.50},\, {\it 19.25},\lb{\it 20.125},\, {\it 21.8},\, {\it 21.24},\, {\it 21.147},\, {\it 21.150}

Dovhyi, S.~ {\it 12.51}

Droste M.~ 10.32

Drutu C.~ 19.16

Dudkin F.\,A.~ 19.27

 Dugas M.~ 10.53,\, \textbf{A:}~{\it 7.22}

 Dunfield N.\,M.~ 17.107,\, \textbf{A:}~{\it 15.49}

Dunwoody M.\,J.~ \textbf{A:}~{\it 5.28},\, {\it 13.62}

Durakov E.\,B.~ \textbf{A:}~{\it 11.13},\, {\it 17.71}

Dwyer W.\,G.~ 19.92

Dyer M.\,N.~ \textbf{A:}~{\it 4.54}

~

\textbf{E}berhard S.~ {\it 20.17},\, 21.44

Edmunds C.\,C.~ 14.14

Eggert N.~ 17.52

Egorychev G.\,P.~ \textbf{A:}~{\it 2.18}

Eick B.~ {\it 14.64}

\`Eidinov M.\,I.~ \textbf{A:}~2.88

Eisele F.~ \textbf{A:}~{\it 17.67}

Emaldi M.~ \textbf{A:}~11.128

 Epstein D.\,B.\,A.~ 18.13

Epstein I.~ 20.63

 Erfanian A.~ \textbf{A:}~17.116

Ershov A.\,V.~ \textbf{A:}~{\it 7.30}

 Ershov Yu.\,L.~ 1.20,\, 11.66,\, {\bf
 A:}~1.18--19, \lb3.13,\, {\it 15.8}

Ershov M.~
{\it 14.93},\, {\it 18.89},\,
18.107--108,\lb{\it 19.10},\, \textbf{A:}~{\it 13.36},\, {\it 14.33--34},\, {\it 14.55},\lb{\it 15.6},\, {\it 15.18},\, {\it 17.17}

Esteban-Romero R.~ {\it 17.56}

Evans D.\,M.~ \textbf{A:}~13.28

Evans R.\,J.~ 10.32

Ezquerro L.\,M.~ {\it 14.28},\, \textbf{A:}~12.96

~

 \textbf{F}aghihi A.~ \textbf{A:}~{\it 11.46}

Falco, De, M.~ {\it 16.38}

 Fang X.\,G.~ {\it 11.80}

Farbman S.\,P.~ 17.25

J.\,Fari\~{n}a-Asategui J.~ {\it 19.22}

Faulkner J.~ 19.92

F\"edorov A.\,N.~ \textbf{A:}~{\it 5.50}

F\"edorova M.~ \textbf{A:}~{\it 16.85}

Feighn M.~ \textbf{A:}~{\it 4.8}

Fein B.~ \textbf{A:}~{\it 7.16}

Feit W.~ 4.65,\, 6.10,\, 18.97,\, 21.121,\lb\textbf{A:}~{\it 7.48}

 Felsch V.~ 18.90

Felsch W.~ 14.95

Fel'shtyn A.\,L.~ 19.28--29

Felzenbaum A.~ 21.93

Feng Z.~ {\it 13.43}

Fern\'andez-Alcober~G.\, 14.89--91,\, 21.104,\lb\textbf{A:}\,14.92

Fern\'{a}ndez-Ayala \'{O}.~ {\it 17.113}

Ferrara M.~ {\it 20.109}

Ferudun A.~ {\it 13.19}

Fesenko I.\,B.~ 14.93--94,\, \textbf{A:}~{\it 12.24}

Fine B.~ 20.22,\, 20.41,\, 20.85,\, \textbf{A:}~{\it 9.2},\lb14.88,\, 15.86,\, {\it 15.86}

Fine Y.~ \textbf{A:} {\it 3.51}

Fitzpatrick P.~ \textbf{A:}~{\it 2.43},\, {\it 7.48}

 Flavell P.~ 18.106

Flores R.\,J.~ 19.91

 Fomin A.\,N.~ 7.50--52,\, 9.69--71,\, 10.65, \linebreak
\makebox[3ex][r]{}14.67,\, \textbf{A:}~10.1,\, 10.66

Fon-Der-Flaass D.\,G.~ 8.77

Fong P.~ \textbf{A:}~{\it 4.76}

 Fontaine J.-M.~ \textbf{A:}~12.24

Formanek E.~ 12.92,\, \textbf{A:}~{\it 2.16}

Foote R.\,M.~ 19.91

Foroudi Ghasemabadi M.~ 19.30

 F\"orster P.~ 18.8

Fournier-Facio F.~ {\it 6.45},\, {\it 15.45},\, {\it 17.62},\lb{\it 21.14},\, 21.22,\, 21.45--48

Fraenkel A.~ 21.93

Frid A.\,E.~ {\it 16.9}

Fridman \`E.\,I.~ \textbf{A:}~{\it 3.21}

Frigerio R.~ {\it 18.111}

Frobenius G.~ 12.92,\, 20.86

Fuchs L.~ 1.35,\, 2.25--26

 Fujiwara K.~ {\it 14.15},\, 19.31,\, \textbf{A:}~{\it 14.13}

Fulman J.~ {\it 7.27},\,  15.65

Fumagalli F.~ {\it 18.110},\, 21.23--24

Fusari M.~ {\it 11.71},\, 21.28

~

 \textbf{G}aglione A.~ 11.19,\, 12.8--9,\, 13.39,\,
18.35,\lb19.32,\, \textbf{A:}~11.20,\, 13.18,\, 13.39,\, 14.32,\lb15.86,\, {\it 15.86}

Gallagher P.~ \textbf{A:}~{\it 11.43}

 Gal't A.\,A.~ {\it 17.37},\, {\it 17.112},\, \textbf{A:}~{\it 14.82}

Gamkrelidze R.~ 21.123

Garciarena M.\,E.~ {\it 19.22},\, {\it 20.101}

Garonzi M.~ 20.23

Garside F.\,A.~ \textbf{A:}~{\it 1.11}

Gasch\"utz W.~ 6.1,\, 9.59,\, \textbf{A:}~{\it 16.62}

 Gavioli N.~ 18.34,\, \textbf{A:}~{\it 14.49}

 Genevois A.~  {\it 14.15},\, 21.49

Geoghegan R.~ \textbf{A:}~{\it 6.4}

 Gerstenhaber M.~ {\it 17.94},\, 18.87

 Ghate S.\,H.~ 17.47

 Gheri P.~ 20.86

Ghys \'{E}.~ \textbf{A:}~{\it 14.13}

Giannelli E.~ 19.33,\, {\it 20.78}

 Gill N.~ {\it 14.46}

 Giovanni de F.~ 13.22--23,\, 14.36,\, 16.38

Giraudet M.~ \textbf{A:}~12.14,\ {\it 12.14}

Giudici M.~ 21.29

Glasner Y.~ 20.24

 Glass A.\,M.\,W.~ 12.11--13,\,
14.8,\lb17.26--27,\, 21.79,\, \textbf{A:}~{\it 5.23},\, 12.10,\lb12.14,\, {\it 12.14},\, {\it 15.34},\, 17.28--29,\, {\it 17.29}

 Glauberman G.~ 16.33,\, 20.93,\, 21.50--51,\lb\textbf{A:}~{\it
1.80},\, 4.21--23,\, 6.7,\, {\it 14.1},\, 15.24

Gleason~ 6.11

 Gluck D.~ 9.23,\, 20.25,\, \textbf{A:}~{\it 11.33}

Glushkov V.\,M.~ 6.11

 Go S.~ \textbf{A:}~{\it 16.1}

 G\"obel R.~ 10.54,\, \textbf{A:}~{\it 7.22},\, 11.26,\,
{\it 11.26},\lb{\it 12.44}, \,{\it 14.71}

Gogniat L.~ {\it 11.71}

Goldschmidt D.\,M.~ 4.24,\, 19.91,\, \textbf{A:}~4.24

 Goldstein R.~ {\it 15.63},\, \textbf{A:}~{\it 15.63}

Golikova E.\,A.~ 12.15

Golod E.\,S.~ 9.76,\, 21.132,\, \textbf{A:}~{\it 2.64},\, 8.66,\lb15.6

 Golovin O.\,N.~ 2.9--10,\, 12.30,\, \textbf{A:}~2.7,\, 2.8

Gon\c{c}alves D.~ 19.28,\, 19.34,\, \textbf{A:}~{\it 18.15}

Goncharov S.\,S.~ 12.16--17

Gong L.~ {\it 4.13},\, {\it 20.126},\,  {\it 21.51}

Gong  T.~ {\it 20.2},\, {\it 21.142}

 Gonz\'alez-Meneses J.~ \textbf{A:}~{\it 10.23}

 Gonz\'alez-S\'anchez J.~ 20.26,\, \textbf{A:}~{\it 12.77}

 Goodman A.\,J.~ 17.41

Gorbunov V.\,A.~ 9.5--6,\, \textbf{A:}~{\it 3.15}

 Gorchakov Yu.\,M.\,\, 3.12,\, 6.9,\, 10.16,\lb13.19,\, \textbf{A:}~1.9,\, {\it 1.9},\, 1.10,\,
{\it 2.18},\, {\it 2.43, \lb2.44},\,
3.9--11,\, 6.8,\, 7.8

Gorenstein D.~ \textbf{A:}~{\it 3.6},\, {\it 3.64}

 Gorshkov I.\,B.~ {\it 12.38},\, 20.27--30,\, {\it 20.58},\lb21.52--53,\, \textbf{A:}~{\it 15.53},\, {\it 16.27},\, {\it 16.107}

 Goryachenko V.\,A.~ {\it 20.37}

Goryushkin A.\,P.~ \textbf{A:}~{\it 3.11},\, {\it 4.59}

Grazian V.~ {\it 16.1},\, 21.35

 Grechkoseeva M.\,A.~ {\it 15.53},\, 20.31,\lb\textbf{A:}~{\it 12.39},\, {\it 14.60},\, {\it 16.24},\, {\it 16.106},\lb{\it 17.36},\, {\it 17.73},\, {\it 17.74}

 Green B.~ 18.108,\, 20.32

 Green J.\,A.~ \textbf{A:}~{\it 18.82}

 Greendlinger M.\,D.~ 1.12,\, \textbf{A:}~1.11,\lb1.13--14,\, 1.17

 Greuel G.-M.~ \textbf{A:}~{\it 15.75}

Griess R.~ 7.15 ~ \textbf{A:}~7.12,\, 7.16,\, 7.17,\lb8.37

 Grigorchuk R.\,I.~ 8.8--9,\,
9.7,\, 9.9,\lb10.11--13,\, 12.19--20,\, 13.20--21,\,
14.33,\lb15.12--14,\, 15.16,\, 15.19,\, 15.78,\, 16.34,\lb16.84,\, \textbf{A:}\,{\it 4.5},\, {\it 5.2},\, 8.6--8,\, {\it
8.6},\, 9.8,\lb9.50,\, 13.21,\, 14.34,\, 15.4,\, 15.14--15,\lb15.17--18,\, 16.101, {\it 16.101}

Griffith P.\,A.~ \textbf{A:}~{\it 2.4}

Grishkov A.\,N.~ 21.54

Grochow J.~ 20.33

Gromov M.~ 14.7,\, 18.12,\, 21.86,\, \textbf{A:}~{\it 5.68},\lb{\it 7.4},\, 11.74

Groot de~ \textbf{A:}~2.36

Gross F.~ 5.30,\, 7.31,\,
\textbf{A:}~3.26,\, 13.33,\lb17.44--45

Grosser S.\,K.~ \textbf{A:}~9.34,\, {\it 9.34}

Grossman J.\,W.~ \textbf{A:}~6.36

 Groves D.~ {\it 14.19},\, {\it 14.20},\, \textbf{A:}~{\it 11.105}

 Groves J.\,R.\,J.~ 17.30,\, \textbf{A:}~10.37,\, 16.35

 Gruenberg K.\,W.~ 11.28,\, 18.13,\, 20.89,\lb\textbf{A:}~5.8--11,\,
11.27

 Grunewald F.\,J.~ \textbf{A:}~{\it 4.53},\, {\it 15.75}

Grytczuk A.~ {\it 8.30}

 Guba V.\,S.~ 15.42,\, \textbf{A:}~{\it
1.1},\, {\it 6.44},\, {\it 8.8},\lb9.10,\, {\it 15.4}

Gudivok P.\,M.~ 12.21

G\"{u}lo\u{g}lu \.{I}~ 21.55--56

 Gunhouse S.\,V.~ {\it 12.11}

 Guo W.~ 20.34,\, 21.57,\, \textbf{A:}~{\it 12.74},\, {\it 19.50},\lb{\it 19.87},\, {\it 19.88}

 Gupta N.\,D.~ 11.17,\, 11.29,\, 12.4,\, 14.35,\lb16.66,\, 16.96,\, \textbf{A:}~11.29,\, {\it
11.29},\, 12.22,\lb{\it 12.22}

 Guralnick R.\,M.~ {\it 7.27},\, 17.66,\, 18.65,\lb 21.25,\, 21.38,\, 21.110,\, \textbf{A:}~{\it 6.6},\, {\it 7.17},\lb {\it 8.5}

Gurchenkov S.\,A.~ \textbf{A:}~{\it 1.61}

Guti{\'e}rrez C.~ \textbf{A:}~{\it 9.25}

Gutsan M.~ 12.23

 Guyot L.~ \textbf{A:}~{\it 15.4}

Gvaramiya A.\,A.~ 9.4,\, \textbf{A:}~2.2,\, {\it 2.2}

Gvozdev R.\,I.~ {\it 14.69}

~

\textbf{H}adjievangelou A.~ {\it 20.106}

Hagen M.\,F.~ 20.35,\, {\it 17.108}

Hagenah C.~ \textbf{A:}~{\it 9.25}

Haglund F.~ 20.60,\, 21.124

Halasi Z.~ 19.12

 Hales A.\,W.~ \textbf{A:}~16.104, {\it 16.104}

Hall J.\,I.~ \textbf{A:}~{\it 12.76}

 Hall M.~ 16.7

 Hall P.~ 2.45,\, 4.6,\, {\it 7.6},\, 10.16,\, {\it 14.10},\lb 14.89,\, 17.101,\,  {\it 20.120},\, {\bf
A:}~2.45,\lb 4.59--60,\, 12.31,\, 12.44--45

Hamidoune Y.\,O.~ 19.26

 Hammoudi L.~ 9.76,\, \textbf{A:}~{\it 8.66},\, {\it 12.102}

Handel M.~ \textbf{A:}~{\it 4.8}

 Harada K.~ 18.64

Harpe\,de\,la\,P.~
14.6--10,\, 15.8,\, \textbf{A:}~{\it 5.2},\lb14.10,\, 15.4--8

Harper S.~ {\it 11.71},\, 21.38

 Hartley B.~ 5.59,\, 8.78--79,\, 11.22,\, 11.49,\linebreak
\makebox[3ex][r]{}11.109,\, 11.111,\, 11.113,\,
{\it 12.27},\,
12.29,\lb13.6--9,\, 15.20--21,\, 17.81,\, 17.119,
\lb20.16,\, 21.65--67,\, \textbf{A:}~{\it
1.75},\, {\it 3.14},\, {\it 4.76},\lb{\it 5.18},\, 5.60--61,\, {\it 5.60},\, {\it 6.20},\, {\it
7.6},\,
7.37,\lb8.80--81,\, {\it
8.80},\, 11.110,\, {\it 11.110},\, {\it
14.52}

 Hashimoto M.~ {\it 15.44},\, \textbf{A:}~{\it 15.44}

Hassani A.~ 15.95

Hatem O.~ {\it 17.76}

Hauck P.~ 14.30

 Havas G.~ {\it 11.80},\, 20.36,\, \textbf{A:}~{\it 8.12},\, {\it 9.50}

 Hawkes T.~ \textbf{A:}~{\it 14.77}

 Heffernan R.~ \textbf{A:}~17.77

Heineken H.~ 11.17,\,
\textbf{A:}~{\it 2.63},\, 12.91,\lb{\it 12.91},\, {\it 16.62}

Helling H.~ 8.82--83

 Herfort W.\,N.~ {\it 20.116},\, \textbf{A:}~{\it 3.39},\, {\it 9.32},\lb9.34,\, {\it
9.34},\, 12.70,\, {\it 12.70}

Hering C.~ \textbf{A:}~{\it 9.21}

Herrman C.~ 12.75

Hertzig D.~ \textbf{A:}~1.72

Hertweck M.~ {\it 20.5}

Herzog M.~ 21.93

Hewitt P.\,R.~ 13.59

Hickin K.\,K.~ 8.79,\, {\it 12.11},\, \textbf{A:}~{\it 7.6}

 Higman G.~ 6.21,\, 17.101,\, 20.7,\, 20.33,\lb\textbf{A:}~{\it 1.14},\, {\it 4.76},\, 5.22,\, 12.32,\, 14.10,\lb18.33

 Hilbert D.~ 6.11,\, 14.22

Hilton P.~ \textbf{A:}~8.84,\, 12.93

Hirzebruch F.~ 15.56

 Hjorth G.~ \textbf{A:}~{\it 13.28}

 Hoechsmann K.~ \textbf{A:}~{\it 12.1}

Hofmann K.\,H.~ \textbf{A:}~{\it 9.32}

Hog C.~ \textbf{A:}~{\it 7.57}

H\"ohnke H.-J.~ 12.92

Holland W.\,C.~ \textbf{A:}~{\it 3.19}

Holt D.~ 7.52,\, \textbf{A:}~{\it 5.61}

Hooshmand M.\,H.~ 20.37,\, {\it 20.37},\, 21.58,\lb{\it 21.58},\, \textbf{A:}~19.35

Horbez C.~ {\it 14.15}

 Houghton C.\,H.~ \textbf{A:}~{\it 16.52}

Howie J.~ 15.94,\, \textbf{A:}~{\it 5.53}

Huang H.\,Y.~ {\it 15.40},\, {\it 20.121}

Huber M.~ {\it 11.45}

Hughes S.~ {\it 19.41}

Hull M.~ {\it 14.20}

 Humphreys J.\,E.~ 16.39

 Humphries S.\,P.~ \textbf{A:}~{\it 14.102}

Huppert B.~ 15.71,\, 18.99

Hyde J.~ \textbf{A:}~{\it 14.10},\, {\it 16.50},\, {\it 17.20}

~

 \textbf{I}gnatov Yu.\,A.~ \textbf{A:}~{\it 4.41}

 Ihara I.~ 5.33

Ilenko K.\,A.~ {\it 17.103}

Ingrosso E.~ {\it 6.30},\, {\it 11.9},\, {\it 13.23},\, {\it 20.33}

Ionin V~ {\it 14.85},\, {\it 16.11},\, {\it 17.32},\, {\it 17.47},\lb{\it 19.94}

Iranmanesh A.~ 19.30,\, 20.38,\, {\it 20.38},\lb21.59,\, \textbf{A:}~{\it 12.38}

Iranmanesh M.\,A.~  {\it 20.1},\, 20.1--3

Irwin J.~ \textbf{A:}~17.24

Isaacs I.\,M.~ 15.2--3,\, 17.1,\, 17.119,\, 18.1,\lb{\it 20.78},\, 21.65--66,\, \textbf{A:}~{\it 2.29},\, {\it 16.36}

Iskra A.\,L.~ {\it 18.48},\, 21.8,\, {\it 21.75}

Ito N.~ \textbf{A:}~5.17

Ito T.~ \textbf{A:}~{\it 16.49}

Ivanov A.\,A.~ \textbf{A:}~{\it 8.37}

Ivanov S.\,G.~ \textbf{A:}~2.13

 Ivanov S.\,V.~
 11.36,\, 11.37--40,\, 12.30,\lb\textbf{A:}~{\it 1.13,\, 1.81,\, 2.8},\, {\it
 2.45},\, {\it
 2.72},\lb{\it 4.1},\, {\it 4.45},\, {\it 4.46},\, {\it 4.73},\, {\it 4.74},\, {\it 6.16},\lb{\it 8.20},\, {\it 11.36},\, 11.37,\, {\it 11.37},\, {\it 11.74},\lb12.31,\, 12.32

 Iwahori N.~ 2.52

Izosov A.\,V.~ \textbf{A:}~11.42

~

 \textbf{J}abara E.~ 17.128,\, {\it 18.102},\,
18.120,\, 19.28,\lb20.106,\, \textbf{A:}~19.36,\, {\it 19.80}

 Jackson D.~ 11.19

 Jaikin-Zapirain A.~ 14.89,\, 14.95,\, 15.93,\lb16.102,\,
 17.109,\, 17.110,\, 18.107--109,\lb{\it 19.102},\, {\it 19.103},\, {\it 19.104},\, {\it 19.105},\lb20.26,\, 20.39,\, 20.67,\, \textbf{A:}~{\it 12.77},\, {\it 12.94},\lb{\it 13.39},\, {\it 13.56},\, {\it 14.34},\, {\it 14.92},\, {\it 14.96},\lb{\it 15.6},\, {\it 15.18},\, {\it 17.17}

 Jaligot E.~ {\it 10.49},\, 17.48

Janko Z.~ \textbf{A:}~4.77

 Jantzen J.\,C.~ 16.39

Januszkiewicz T.~ {\it 15.45}

Jayadevan Achyuth~ {\it 2.81},\, {\it 3.12},\, {\it 4.55},\lb {\it 5.38},\, {\it 6.9},\, {\it 7.15},\, {\it 7.38},\, {\it 8.16},\, {\it 9.4},\, {\it 9.41},\lb {\it 9.47},\, {\it 10.4},\, {\it 11.5},\, {\it 13.51},\, {\it 14.22},\, {\it 14.67},\lb {\it 15.32},\, {\it 15.77},\, {\it 15.89},\,  {\it 16.28},\, {\it 17.34},\, {\it 18.11},\lb {\it 18.76},\, {\it 21.106},\, {\it 21.121}

Ji H.~  {\it 18.16}

Jockusch C.\,G., jr.~ \textbf{A:}~{\it 4.20}

Johnson D.\,L.~ 8.11--12,\, \textbf{A:}~{\it 6.23},\,
8.10, \lb8.12,\, 12.24

Johnson K.\,W.~ 12.92

Johnston D.~ {\it 4.55},\, 21.60

Jonah D.~ 21.51

Jones G.\,A.~ \textbf{A:}~{\it 1.21},\, {\it 3.10}

Jones N.\,A.~  {\it 18.44}

 Jonsson B.~ 9.66,\, 12.75

Judin E.~ {\it 3.46},\, {\it 18.50},\, {\it 19.25},\, {\it 20.125},\lb{\it 21.8},\, {\it 21.24},\, {\it 21.147},\, {\it 21.150}

 Jung R.~ 15.56

 Juschenko K.~ \textbf{A:}~{\it 15.5}

~

\textbf{K}abanov V.\,V.~ 5.14,\, \textbf{A:}~5.12--13,\, {\it 6.13}

Kabanova E.\,I.~ \textbf{A:}~{\it 10.7}

Kabenyuk M.\,I.~ {\it 21.18},\, {\it 21.58},\, \textbf{A:}~{\it 3.58}

Ka\"{i}chouh A.~ {\it 18.36}

Kaja E.~   21.24

Kautsky J.~ 20.36

Kaygorodov I.~ \textbf{A:}~{\it 12.38},\,

Kalantar M.~ \textbf{A:}~{\it 15.7}

Kalithasan S.~ {\it 21.32}

 Kallen van der W.~ 10.38

Kaljulaid U.\,\`E.~ 14.45

 Kamornikov S.\,F.~ {\it 9.75},\,
 12.35,\, 14.47,\lb{\it 17.39},\, 17.54,\, {\it 17.56},\, 19.38--39,\, {\it 19.38},\lb\textbf{A:}~{\it
9.73},\, {\it 9.74},\, {\it 11.24},\, {\it 12.7},\, {\it 14.99},\lb15.38--39,\, {\it 15.38},\, {\it 15.39},\, {\it 16.12},\, 16.82,\lb{\it 16.82},\, 17.55,\, {\it 17.55},\, {\it 18.9},\, 19.37,\lb{\it 19.85},\, {\it 19.87},\, {\it 19.88}

Kantor W.\,M.~ \textbf{A:}~{\it 7.16}

Kaplansky I.~ 12.29,\, 18.80

 Kapovich I.~ 16.110,\, 18.43,\, 19.41--42,\lb21.61--62,\, \textbf{A:}~19.40

Kar A.~ {\it 14.6}

Karasev G.\,A.~ \textbf{A:}~2.15

 Kargapolov M.\,I.\,\, 1.27--28,\, 1.31,\, 1.33, \linebreak
\makebox[3ex][r]{}1.35,\, 3.16,\, 4.30--31,\, 4.33--34,\, 4.42,
\lb8.15,\, \textbf{A:}~1.21,\, 1.23--26,\,
1.29--30,\, 1.32, \lb1.34--35,\, {\it 1.48},\,
{\it 1.56},\, 2.16--20,\, {\it 2.16}, \linebreak
\makebox[3ex][r]{}{\it 2.22}, {\it 2.46},\, 3.15,\, {\it 3.18},\,
4.32

 K\'arolyi Gy.~ 12.62,\, \textbf{A:}~{\it 10.63}

 Kassabov M.~ 18.108,\, \textbf{A:}~{\it 14.34}

Kaya B.~ 20.40,\, {\it 20.40}

 Kazarin L.\,S.~ 11.44,\, 12.33--34,\, 13.27,\lb14.43--44,\, 17.52--53,\, {\it 19.100},\, \textbf{A:}~{\it 4.27},\lb{\it
5.63},\, {\it 5.64},\, {\it 6.37},\, {\it
8.48},\, {\it 9.56},\, {\it 9.63},\lb11.43,\, {\it 11.94},\, {\it 14.63}

 Kazhdan D.~ \textbf{A:}~{\it
2.50}

Kearnes K.~ \textbf{A:}~{\it 14.32}

Kechris A.\,S.~  {\it 21.15}

 Kegel O.\,H.~ 8.14,\, 14.28,
 14.43,
\lb\textbf{A:}~4.35--37,\, {\it 4.37},\,
5.17--19,\, {\it 5.19},\lb6.37,\, {\it 7.6},\, 8.14,\, 10.66

 Keller T.\,M.~ {\it 18.44},\, 18.44--45,\, \textbf{A:}~{\it 14.66}

 Kemer A.\,R.~ \textbf{A:}~{\it 2.39}

 Kemkhadze Sh.\,S.~ 1.40,\, 2.22,\, 8.15--16, \linebreak
\makebox[3ex][r]{}\textbf{A:}~1.36--39,\, 2.21--23,\, 8.17

Kennedy M.~ \textbf{A:}~{\it 15.7}

Kerov S.\,V.~ 12.28

 Kervaire M.~ 18.87

Kessar R.~ {\it 13.43}

 Khan' Ch.~ \textbf{A:}~{\it 16.1}

Kharchenko V.\,K.~ 13.57

 Kharlampovich O.\,G.~ {\it 2.81},\,
14.22,\lb{\it 14.23},\, 16.70,\, 16.72,\, {\it 16.70},\, {\it 17.32},\lb  17.117,\, \textbf{A:}~{\it 1.19},\, {\it
1.25},\, {\it 1.29},\, {\it 1.30},\lb {\it 4.3},\, {\it 7.24},\, {\it 9.67},\, {\it 13.18},\, {\it
13.40},\, 16.71,\lb {\it 16.71}

Khayaty M.~ 15.95

 Khelif A.~ \textbf{A:}~{\it 14.48}

Khlamov E.\,V.~ \textbf{A:}~{\it 8.56}

Khosravi Behn.~ \textbf{A:}~{\it 12.38}

Khosravi Behr.~ {\it 20.79},\, \textbf{A:}~{\it 12.38}

Khukhro E.\,I.~ 6.47,\, 8.85,\, 10.67,\lb11.112--113,\, {\it 13.8},\, 13.57,\,
15.95--96,\lb17.118--119,\, 18.66--67,\, 18.110,\, {\it 19.43}, \lb19.43--44,\, 21.63--67,\, \textbf{A:}~{\it 1.10},\, {\it
3.26},\lb{\it 4.76},\, {\it 5.56},\, {\it
6.19},\, {\it 7.8},\, 7.53,\, {\it 7.53},\lb{\it 8.81},\, 10.68,\, 11.112,\, {\it 11.126},\, 13.56,\lb13.58,\, {\it 13.58},\, 14.96,\, 16.103,\, {\it 16.103}

Khusainov A.\,A.~ \textbf{A:}~{\it 2.35}

Kida M.~ 21.68

Kielak D.~ {\it 15.42},\, {\it 15.92},\, 20.41,\, 21.69--72

Kim S.-H.~ {\it 12.20},\, {\it 15.42}

Kim Y.~ {\it 2.82}

 Kiming I.~ \textbf{A:}~14.96

Kimmerle W.~ 20.6,\, {\it 20.6}

 King C.\,S.\,H.~ \textbf{A:}~{\it 18.52}

King J.~ 21.78

Kinyon M.~
{\it 4.75},\, {\it 9.4},\, {\it 10.62},\, {\it 12.40},\lb{\it 11.116},\, {\it 13.42},\, {\it 14.72},\, {\it 15.92},\, {\it 16.20},\lb  {\it 16.38},\,  {\it 16.87},\, {\it 17.33},\, {\it  18.18},\, {\it 18.92},\lb {\it 19.61},\, {\it 19.62},\,  {\it 19.93},\, {\it 20.1},\, {\it 20.90},\lb {\it 20.108},\, {\it 20.122},\, {\it 21.60}

Kionke S.~ {\it 14.53},\, {\it 15.14},\, 20.42--43,\lb\textbf{A:}~{\it 15.49}

Kiralis G.~  {\it 18.89}

Kirkinski\u{\i}
A.\,S.~ \textbf{A:}~{\it 19.40}

Kizmaz M.\,Y.~ {\it 4.13}

 Kleidman P.\,B.~ \textbf{A:}~{\it 5.41},\, {\it 8.39},\, 11.53,\lb{\it 11.53},\, {\it 18.115}

 Kleiman Yu.\,G.~ {\it
2.45},\, 8.19--20,\, \textbf{A:}~{\it 2.40},\, \lb{\it 3.4},\, {\it 4.47},\,
{\it 4.64},\, {\it 5.43},\, {\it 6.17},\, 8.18

Kleim\"enov V.\,F.~ \textbf{A:}~{\it 1.53},\, {\it 8.56}

Klein A.~ \textbf{A:}~{\it 1.4}

 Kletzing D.~ \textbf{A:}~{\it 15.25}

Klopsch B.~ 21.122,\, \textbf{A:}~{\it 12.94}

 Klyachko Ant.\,A.~ {\it 17.94},\, 18.46

 Kneser M.~ \textbf{A:}~3.42

 Knyagina V.\,N.~ 17.56

Koch H.~ 5.26--27

 Koch-Hyde L.~  {\it 10.13},\, {\it 20.88}

 Kochloukova D.\,H.~ \textbf{A:}~{\it 10.37}

 Kohl S.~ 17.57--61,\, 18.47--48,\,  {\it 18.50},\lb18.50--51,\, 19.45--46,\, 20.44--45,\lb21.73--75,\, {\it 21.75},\, \textbf{A:}~18.49

Ko\u{\i}baev V.\,A.~ 19.47--48,\, {\it 19.48},\lb 21.76--77

 Kokorin A.\,I.~ 1.46,\, 1.51,\, 1.54--55,
\lb2.24--26,\, 2.28,\, 3.20,\,
\textbf{A:}~1.41--45, \lb{\it
1.41,\, 1.42,\, 1.43},\, 1.47--50,\, {\it 1.48},\lb1.52--53,\, 1.56,\, {\it 1.56},\, 1.60--61,\,
2.25,\lb3.13,\, 3.17--19,\, 3.21,\,
5.20

Kolchin E.~ \textbf{A:}~2.62

 Kolesnikov S.\,G.~ {\it 8.3},\, \textbf{A:}~{\it 12.42}

Kolmakov Yu.\,A.~ \textbf{A:}~{\it 5.56}

 Kondratiev A.\,S.~ 9.15,\, 10.15,\, 12.37,\lb12.40,\, \textbf{A:}~{\it 4.24},\, 6.13,\, {\it 6.13},\, {\it
8.39},\lb9.16,\, {\it 9.16},\, 12.38--39,\, 14.50

 K\"onig J.~ \textbf{A:}~{\it 18.49}

Konovalov A.~ {\it 20.6}

Kontorovich P.\,G.~ \textbf{A:}~1.63--64,\, 2.29--30

 Kontsevich M.~ 17.84

Konvisser M.~ 21.51

 Koolen J.\,H.~ 20.54

 Kopytov V.\,M.~
5.25,\, 8.24,\, 13.32, \linebreak
\makebox[3ex][r]{}16.47--48,\, 16.51,\, 21.147--149,\lb\textbf{A:}~{\it 1.35},\, {\it
1.41,\, 1.42,\, 1.43},\, {\it 1.48},\, {\it 1.56},\lb2.31,\, 5.23--24,\, 16.49--50

Kordek, K.~ {\it 14.102}

 Kostrikin A.\,I.~ 11.48,\, 11.112,\, 18.22,\lb\textbf{A:}~{\it
5.56},\, 7.53

 Kov\'acs L.\,G.~ 8.21,\, 8.23,\, 21.87,\, \textbf{A:}~{\it 2.43},\lb8.22,\, {\it 10.56},\, 11.47,\, 12.71

 Kov\'acs S.\,J.~ 12.62,\, \textbf{A:}~{\it 10.63}

 Kozhevnikov P.\,A.~ \textbf{A:}~{\it 4.46},\, {\it 4.73},\, {\it 6.15},
\lb{\it 6.16},\, {\it 8.20,\,
13.34}

 Kozhukhov S.\,F.~ \textbf{A:}~7.22

Kozlov G.\,T.~ \textbf{A:}~{\it 3.13},\, {\it 3.21}

 Krammer D.~ 17.16

 Krasil'nikov A.\,N.~ 14.51,\, \textbf{A:}~{\it 2.43},\,
{\it 7.24}, \lb{\it 11.21},\, {\it
13.10}

Kravchenko A.\,A.~ \textbf{A:}~{\it 3.31}

Kravchuk M.\,I.~ \textbf{A:}~2.3

 Krempa J.~ \textbf{A:}~18.52

Kriebel A.~ {\it 16.45}

 Kropholler P.\,H.~ 12.41,\, {\it 12.41},\, 15.45,\lb{\it 15.45},\, {\it 20.118},\, 21.36,\, 21.78,\, \textbf{A:}~{\it 5.8},\lb12.36

Krsti\v{c} S.~ {\it 18.89}

 Kruse R.\,L.~ \textbf{A:}~2.41

Kudlinska M.~ {\it 19.41}

Kukharev A.~ \textbf{A:}~{\it 12.38},\,

Kukin G.\,P.~ \textbf{A:}~{\it 5.46}

 Kulikov L.\,Ya.~ 1.65,\, 1.67,\, \textbf{A:}~1.66

 K\"ulshammer B.~ \textbf{A:}~{\it 4.14}

K\"ummich F.~ \textbf{A:}~9.33

Kun G.~ {\it 21.86}

 Kunyavski\u{\i} B.~ \textbf{A:}~{\it 15.75}

Kurdachenko L.\,A.~ 8.29,\, 9.17,\lb10.16--18,\, 11.49,\, \textbf{A:}~{\it
6.12},\,
9.17,\, {\it 9.17}

 Kurmazov R.\,K.~ {\it 15.40}

Kurosh A.\,G.~ 7.25,\, \textbf{A:}~9.79

 Kurzweil H.~ 5.30,\, 16.98,\, 21.66

 Kuz'min Yu.\,V.~ 16.66,\, 17.81,\, \textbf{A:}~{\it 11.29},\lb12.22,\, {\it 12.22},\, 14.52

 Kuz'minov V.\,I.~ \textbf{A:}~2.33,\, 2.35--36,\, {\it
2.36}, \lb3.22

 Kuznetsov A.\,A.~ \textbf{A:}~{\it 16.57}

 Kuznetsov A.\,V.~ 7.23,\, 7.25--26,\,
8.25, \lb8.27,\, \textbf{A:}~6.14--17,\,
7.24,\,
8.26,\, 8.28

Kuznetsov E.\,A.~ \textbf{A:}~{\it 4.70}

 Kuzucuo\=glu M.~ 15.20,\, 18.53,\, 20.40,\lb{\it 20.40},\, \textbf{A:}~{\it 5.1},\, {\it
5.18}

~

 \textbf{L}abute J.\,P.~ 5.27

Lackenby M.~ {\it 21.10}

Ladisch F.~ 18.1

Lady E.\,L.~ 11.50--51

 Lam C.\,H.~ {\it 15.55}

 Larsen M.\,J.~ {\it 12.95},\, {\it 20.74},\, 21.27,\lb21.121,\, \textbf{A:}~{\it 1.75,\, 7.6}

Laubie F.~ \textbf{A:}~12.24

 Laudenbach F.~ 18.87

Lauri J.~ 20.71

Lausch H.~ 8.30,\, \textbf{A:}~9.18

Lavrenyuk Ya.~ \textbf{A:}~{\it 16.86}

Leary I.\,J~ {\it 15.45},\, \textbf{A:}~{\it 19.23},\, {\it 19.24}

Lebed V.~ \textbf{A:}~19.49

Le~Boudec A.~ {\it 19.10},\, 20.46--47,\, \textbf{A:}~{\it 15.7}

Ledermann W.~ \textbf{A:}~{\it 18.82}

Lee D.~ {\it 13.66},\, {\it 16.110}

Lee M.~  {\it 17.75}

 Leedham-Green C.\,R.~ 14.55--57,\, 14.95, \linebreak
\makebox[3ex][r]{}15.48,\, 17.63--65,\, \textbf{A:}~16.103

Leemans D.~ {\it 11.80}

Lehner F.~ {\it 19.89}

Leinen F.~ {\it 12.27},\, {\it 12.28},\, {\it 18.110}

Le Ma\^{i}tre F.~ {\it 18.36}

 Lennox J.\,C.~ 11.31,\, 12.43,\, \textbf{A:}~8.32,
\lb12.44--45

Lenyuk S.\,V.~ \textbf{A:}~{\it 12.5}

Levai L.~ 14.53

 Levchuk V.\,M.~ 7.27--28,\, {\it 8.31},\, 10.19--20,\lb15.46,\, 18.68,\, \textbf{A:}~6.18,\,
6.19,\, {\it 6.34},\lb8.31,\, 12.42

 Levich E.\,M.~ {\it 3.47},\, \textbf{A:}~{\it 1.69}

 Levin F.~ 16.96

Lewis M.~ 20.48

Li B.~ {\it 21.51}

Li C.\,H.~ {\it 7.31},\, {\it 13.43},\, 14.73,\, 15.47

Li J.~ {\it 10.34},\,   {\it 20.87}

Li T.~ \textbf{A:}~{\it 15.27}

Li Y.~ {\it 21.134}

Liebeck M.\,W.~ 14.54,\, 18.117,\, \textbf{A:}~{\it
5.41},\lb{\it 7.17},\, {\it 9.21}

 Likharev A.\,G.~ 18.68,\, \textbf{A:}~{\it 8.31}

Liman F.\,N.~ 8.33

 Limanski\u{\i} V.\,V.~ \textbf{A:}~{\it 2.10}

Lin C.~ {\it 14.81}

Lin V.~ 14.102,\, \textbf{A:}~14.102

 Linnell P.\,A.~ 11.28,\,
{\it 16.51},\,
18.13,\lb\textbf{A:}~{\it 2.33},\, {\it 4.54},\, {\it 5.9,\, 5.10},\, 8.34,\, {\it 10.30},\lb12.46--47,\, {\it 12.47},\, {\it 14.102},\, 15.49,\, 16.52

Linton S.\,A.~ \textbf{A:}~{\it 8.35},\, {\it 11.46}

Lisi F.~ 21.26,\, 21.114

 Lisitsa A.\,P.~ 18.89

Liu A.-M. \textbf{A:}~{\it 19.87},\, {\it 19.88}

 Liu W.~ {\it 11.45},\, {\it 14.46}

Liu Y.~ {\it 13.43},\, {\it 20.79}

Lodha Y.~ {\it 5.52},\, {\it 12.20},\, {\it 15.42},\, 21.79--80,\lb \textbf{A:}~{\it 16.50}

Longobardi P.~ 15.1,\, {\it 20.40},\, \textbf{A:}~{\it 18.33}

Lopatin A.\,A.~ {\it 14.37}

Lorensen K.~ \textbf{A:}~{\it 12.93}

Lorenz M.~ \textbf{A:}~10.41

Lossov K.\,I.~ \textbf{A:}~{\it 2.75},\, {\it 4.45},\, {\it 9.82}

 L\"ubeck F.~ \textbf{A:}~{\it 7.17}

Lubotzky A.~ 18.108,\, 21.81--86,\, \textbf{A:}~{\it 4.67},\lb13.36,\, 14.34

 Lucchini A.~ 16.53,\, {\it 16.53},\, {\it 17.125},\lb20.49,\, 21.87,\, \textbf{A:}~{\it
11.27},\, {\it 11.104},\, {\it 12.71}

L\"uck W.~ \textbf{A:}~{\it 12.47}

Lusztig G.~ \textbf{A:}~{\it 2.51},\, {\it 7.57}

L'vov I.\,V.~ \textbf{A:}~{\it 2.41},\, {\it 5.56}

Lyndon R.~ 11.10,\, \textbf{A:}~5.29,\, 11.10

Lys\"enok I.\,G.~ \textbf{A:}~{\it 4.1},\, {\it 4.74}

Lytkina D.\,V.~ 18.58--61,\, {\it 19.53},\, {\it 20.31},\lb 20.50,\, {\it 20.96},\, \textbf{A:}~{\it 16.57},\, 18.57,\, 19.50,\lb {\it 19.50}

~

\textbf{M}acbeath A.\,M.~ 4.42,\, \textbf{A:}~4.39

 Macdonald I.\,D.~ \textbf{A:}~13.34, 14.92

 Macedo\'nska O.~ {\it 2.82},\, 16.64,\, 17.78,\lb17.80,\, 18.69,\, \textbf{A:}~17.79

 MacHale D. 16.60,\, 16.63,\, 17.76, \lb20.51--52,\, 21.88--89,\, \textbf{A:}~16.59,\lb16.61--62,\, 17.77

Macintyre A.~ 5.21

MacManus J.~ {\it 14.98}

Macpherson H.\,D.~ 12.63,\,  {\it  18.18},\, \textbf{A:}~{\it 12.10}

Madanha S.\,Y.~ {\it 20.78}

Mader A.~ 10.55,\, 11.50--51

 Madlener K.~ \textbf{A:}~{\it 9.30}

 Magaard K.~ {\it 9.23}

Magidin A.~ 19.110

Magnus W.~ 1.12,\, 9.29,\, 15.102

Maier B.~ \textbf{A:}~{\it 8.14}

 Maimani H.\,R.~ 16.1,\, \textbf{A:}~16.1

Maj M.~ 15.1,\, {\it 20.40},\, \textbf{A:}~{\it 18.33}

 Makanin G.\,S.~ 6.24,\, 10.26,\, 12.50,\lb\textbf{A:}~{\it 1.11},\, {\it
1.18},\, 6.22--23,\, 6.25,\, 9.25,\lb10.23,\, 10.24,\, 10.25,\, 10.26,\,
11.54--55

Makanina T.\,A.~ 6.24

Makar-Limanov L.\,G.~ {\it 3.3}

Makarenko N.\,Yu.~ 18.66--67,\, \textbf{A:}~{\it 11.126},\lb{\it 17.72}

 Makhn\"ev A.\,A.~ 5.14,\, 6.26,\, 6.28,\, 7.33, \linebreak
\makebox[3ex][r]{}7.34,\, 8.40--43,\, 9.26,\, 9.28,\, 10.27,\lb10.29,\,
11.58--60,\, 12.57--60,\, 14.61,\lb20.53--54,\, 21.90,\, \textbf{A:}~{\it 5.13},\, 9.26--27,\lb{\it
9.26},\, 10.28

 Mal'cev A.\,I.~ 1.6,\, 1.33,\, 1.35,\, 2.40,\, 2.42,\lb{\it 2.82},\, 3.47,\, 18.18,\, \textbf{A:}~1.4,\, 1.19,\,
1.35,\lb2.8,\, 2.38--41,\, 2.43--44,\,
{\it 2.64},\, 3.19

Malekan M.\,S.~ {\it 14.53},\, 21.1

Malinowska I.~ \textbf{A:}~{\it 11.46}

 Malle G.~ {\it 9.23},\, {\it 13.43},\, 14.69,\, 18.65,\lb 21.91--92

Malykhin V.\,I.~ 13.49

Mamontov A.\,S.~ \textbf{A:}~{\it 19.80}

Mangahas J.~ {\it 19.31}

 Mann A.~ 11.56,\, 12.55--56,\, 15.29,\, 15.31, \linebreak
\makebox[3ex][r]{}15.95,\, 19.51,\, 20.42,\, \textbf{A:}~{\it
4.67},\, 11.57,\lb{\it 12.90},\, {\it 15.24},\, {\it 15.25},\, {\it 16.62},\, 19.37,\lb{\it 19.75}

Manna P.~  21.24

 Manzaeva N.\,Ch.~ {\it 17.43},\, \textbf{A:}~{\it 17.44}

Marchionna B.~ {\it 19.22}

Margalit D.~ {\it 14.102}

Margolis L.~ 20.5--6,\, 21.93,\, \textbf{A:}~{\it 17.67}

Margulis G.\,A.~ \textbf{A:}~{\it 2.49,\, 2.50},\, {\it 2.55}

Marin E.\,I.~ \textbf{A:}~{\it 10.6}

 Mar\'oti A.~ {\it 19.11},\, 19.12,\, 20.23,\, \textbf{A:}~{\it 8.5},\lb{\it 17.116},\, {\it 18.23}

Martin L.~ {\it 20.79}

 Martino A.~ \textbf{A:}~{\it 10.26}

Martino Di L.~ {\it 11.32}

Mart\'{i} Salvador M.~ {\it 20.52}

 Mashevitzky G.~ 15.76,\, \textbf{A:}~15.76

 Maslakova O.~ \textbf{A:}~{\it 10.25},\, {\it 10.26}

 Maslova N.\,V.~ {\it 17.103},\, 18.68,\, 19.52,\lb21.94--95,\, \textbf{A:}~{\it 17.92}

Matsumoto M.~ \textbf{A:}~2.52

Mattarei S.~ 20.55

 Matte Bon N.~ 20.47

Matteo Di M.~ {\it 20.109}

Matucci F.~ \textbf{A:}~{\it 14.10}

Matveev G.\,V.~ \textbf{A:}~{\it 2.16}

Mazurov V.\,D.~ 5.30,\, 6.21,
 7.31,\,
 {\it 7.30},\lb8.35,\, 9.19,\, 9.23--24,\, 12.48,\, 12.62,\lb{\it 13.57},\,
14.58--59,\, 15.50--54,\, 15.98,\lb16.56,\, 17.66,\, 17.68--72,\, 17.75,\lb17.119,\, 18.62--65,\, 19.53,\, {\it 20.31},\lb {\it 19.53},\,\, 20.56--58,\, {\it 20.96},\, 21.96,\lb \textbf{A:}~{\it
2.3},\, {\it 2.30},\, 2.37,\, {\it 3.6},\, 3.26--28,\lb {\it 3.61},\, 4.38,\, {\it 4.77},\, 5.32,\, 6.20,\, 7.30,\lb 8.35,\, 8.37,\, 8.39,\, 9.19,\, {\it
9.19},\, 9.21,\lb {\it 10.63},\, {\it 11.11},\, 11.52--53,\, {\it 11.97},\, {\it 12.39},\lb 12.49,\, {\it 12.84},\, 13.33--34,\, 14.58,\, 14.60,\lb {\it 16.24},\, 16.27,\, 16.54--55,\, 16.57,\, 17.67,\lb 17.71--74,\, {\it 19.50}

McCleary S.~ 12.13,\, 12.51--54

 McCool J.~ {\it 18.89},\,
 \textbf{A:}~{\it 11.82}

McCourt T.\,A.~ \textbf{A:}~{\it 17.35}

McCulloch R.~ {\it 21.43},\, {\it 21.97},\, \textbf{A:}~{\it 19.35}

McCullough D.~ {\it 10.70}

 McDermott J.\,P.~ 17.61

 McKay J.~ 8.51,\, 15.55--57,\, \textbf{A:}~16.58

 McKay S.~ \textbf{A:}~14.96,\, 16.103

 McReynolds D.\,B.~ \textbf{A:}~{\it 15.35}

 Medvedev N.\,Ya.~ 16.47--48,\, 16.51,\lb21.147--149,\,\textbf{A:}~{\it 1.60, \, 5.20},\, {\it 5.24},\lb10.30,\, 16.49--50

 Medvedev Yu.\,A.~ {\it 2.82},\, \textbf{A:}~{\it 10.68},\,
14.96

Megibben C.~ 10.51

Megrelishvili M.~  21.103

Mehatari R.~  21.24

Meierfrankenfeld U.~ \textbf{A:}~{\it 11.108}

Meixner T.~ \textbf{A:}~{\it 4.76}

Mekler A.\,H.~ \textbf{A:}~{\it 3.19},\, {\it 5.40}

Mellon C.~ {\it 20.17}

Mel'nikov O.\,V.\,\, 6.29--32,\,
 7.35,\, 7.38,
\lb11.61,\, 13.35--37,\, 15.58--59,\, {\bf
A:}~{\it 3.40}, \lb6.31,\, 7.36--37,\, {\it
12.97},\, 13.36,\, 13.38,\lb{\it 13.38}

Melo de E.~ 21.39--40

Menal P.~ 11.22

Mendelssohn N.\,S.~ 12.4

Mendonc L.~ {\it 20.90}

 Menegazzo F.~ \textbf{A:}~{\it 15.33}

Meng H.~ {\it 17.39}

Mennicke J.~ 5.33,\, 7.39,\, 8.44

 Merzlyakov Yu.\,I.~ 2.45,\, 4.40,\, 4.42,\, 5.35,\lb7.41,\, 8.3,\, 8.45,\, 9.29--31, 10.31--32,\lb14.11--12,\, 14.16,\,
15.83--84, \, \textbf{A:}~{\it 1.38},\lb1.68--70,\, {\it 1.73},\, {\it 2.13},\, {\it 2.22}, 2.45--46,\lb{\it 2.46},\, {\it 3.15},\, 3.29,\,
4.41,\, {\it 4.60}, \, 5.34,\lb6.34--36,\, 7.40,\, 8.46--47,\, 10.33,\, 15.10

Metzler W.~ \textbf{A:}~{\it 7.57}

Mignotte M.~ 13.65

 Mikaelian V.\,H.~ \textbf{A:}~{\it 14.10}

 Mikhailov R.~ 16.66,\, 17.81,\, 17.83--90,\lb\textbf{A:}~{\it 11.36},\, {\it 12.22},\, {\it 15.10},\, 16.65,\, 17.82,\lb17.86

 Mikhailova K.\,A.~ 16.7

Mikhal\"ev A.\,V.~ 10.8

 Milenteva M.\,V.~ \textbf{A:}~{\it 8.76}

Mill van J.~ 13.48

 Miller A.~ {\it 10.70}

Miller C.\,F., III~ 5.16,\, 21.47

Miller G.\,A.~ \textbf{A:}~{\it 18.49}

Milnor J.~ \textbf{A:}~4.5,\, 5.68

Minasyan A.~ {\it 14.6},\, 18.70, 20.59--61,\lb\textbf{A:}~{\it 8.69},\, {\it 11.42},\, {\it 16.101},\, {\it 19.23},\, {\it 19.24}

Mines R.~ 10.52

 Minkowski H.~ 4.42

Mishina A.\,P.~ \textbf{A:}~2.47,\, 3.30--31

 Mislin G.~ {\it 10.40},\, 18.111

Miyamoto M.~ {\it 15.2}

Mochizuki H.\,Y.~ \textbf{A:}~{\it 2.62},\, {\it 5.49}

Moede T.~ {\it 14.64}

 Moldavanski\u{\i} D.\,I.~ 3.34,\, 11.62--63,\lb\textbf{A:}~3.33,\, {\it 3.33},\, 15.60

Moldenhauer A.~ 20.85

Molinier R.~ 21.97

 Monakhov V.\,S.~ 10.34,\,
15.61,\, 17.56, \lb17.91,\, 19.54,\, 19.56--58,\,
\textbf{A:}~6.37,\, 8.48,\lb11.64, \, 14.62--63,\, {\it 14.62},\, 17.92,\, 19.55

Monetta C.~ {\it 16.1},\, 21.35

 Monod N.~ {\it 15.8},\, 19.73,\, 20.24,\, 20.62--64,\lb\textbf{A:}~{\it 15.5}

Montgomery D.~ 6.11

 Monti V.~ 18.34

Monticone P.~ {\it 3.46},\, {\it 18.50},\, {\it 19.25},\lb{\it 20.125},\, {\it 21.5},\, {\it 21.8},\, {\it 21.24},\, {\it 21.147},\lb{\it 21.150}

Moody J.\,A.~ \textbf{A:}~{\it 10.41}

 Moret\'o A.~ 14.64--65,\, 20.65,\, \textbf{A:}~14.66,\lb15.62,\, {\it 15.62},\, 16.67,\, {\it 19.55}

Morigi M.~ 15.11,\, {\it 17.125},\, 21.98

Morley L.~ 14.8

Morrison D.~ {\it 3.46},\, {\it 18.50},\, {\it 19.25},\lb{\it 20.125},\, {\it 21.8},\, {\it 21.24},\, {\it 21.147},\, {\it 21.150}

Moskalenko A.\,I.~ \textbf{A:}~{\it 3.35}

Motegi K.~ \textbf{A:}~{\it 16.49}

Mozes S.~ \textbf{A:}~{\it 4.45}

Mukhamet'yanov I.\,T.~ 14.67

 Mukhin Yu.\,N.\,\, 3.38,\, 5.36,\, 9.35,
\lb9.36--38,\, \textbf{A:}~3.35,\, 3.37,\,
9.32--34

 Muktibodh A.\,S.~ 17.47

Muliarchyk K.~ {\it 8.74},\, {\it 9.4},\, {\it 10.4},\, {\it 10.46},\lb{\it 11.23},\, {\it 12.15},\, {\it 12.51},\, {\it 14.67},\, {\it 16.53},\lb{\it 18.59},\, {\it 18.60},\, {\it 19.20},\, {\it 20.30},\, {\it 20.75},\lb{\it 20.76},\, {\it 21.99},\, {\it 21.137},\, {\it 21.148}

 M\"uller P.~ {\it 8.75},\, 21.99

Muranov A.~ \textbf{A:}~{\it 13.11},\, {\it 14.13}

Murashka V.\,I.~ 20.66

 Murski\u{\i} V.\,L.~ \textbf{A:}~{\it 2.40}

 Murty M.\,R.~ \textbf{A:}~{\it 19.50}

 Muzychuk M.~ {\it 16.3}

 Myasnikov A.\,G.~ 11.65--66,\,
13.39,\lb13.41--42,\,
 14.18--22,\, 15.63,\, 16.70, \lb16.72,\, {\it 16.70},\, 18.72,\, 19.64,\, \textbf{A:}~{\it 1.19}, \lb{\it
 1.29,\, 9.67,\, 13.18},\,
 13.39--40,\, {\it 13.40}, \lb{\it 14.71},\, 15.63,\,
 15.86,\, 16.71,\, {\it 16.71}

 Mycielski J.~ 16.68--69,\, 17.93--94,\, 18.71,\lb19.59

 Mysovskikh V.\,I.~ \textbf{A:}~{\it 14.50}

~

 \textbf{N}agrebetski\u{\i} V.\,T.~ 4.44,\, \textbf{A:}~2.29,\,
3.39

Narzullaev Kh.\,N.~ 4.42

 Nasybullov T.\,R.~ 18.14,\, 19.34,\, \textbf{A:}~18.15,\lb{\it 18.15},\, {\it 19.49},\, {\it 19.90}

Natterer S.~ {\it 16.9}

Navarro G.~ {\it 9.23},\, 15.2,\, 20.67,\, 21.92,\lb  21.100,\, \textbf{A:}~{\it 16.36}

Navas A.~ {\it 16.51},\, 21.79--80

Nedorezov T.\,L.~ {\it 13.57}

 Nedrenko D.~ \textbf{A:}~{\it 16.54}

 Neeman A.~ \textbf{A:}~{\it 12.36}

 Nekrashevich V.\,V.~ 15.19,\, 16.34, \lb16.73--74,\, 16.84,\, \textbf{A:}~{\it 9.8},\, {\it 15.5},\, {\it 15.17}, \lb16.73--74,\, {\it 16.73}

 Neroslavski\u{\i} O.\,M.~ \textbf{A:}~{\it 2.65}

 Neshchadim M.\,V.~ 14.68,\, 18.13--14,\lb20.126,\, \textbf{A:}~16.75,\, 18.15

 Nesterov M.\,N.~
\textbf{A:}~{\it 17.45},\, {\it 18.32}

Neukirch J.~ 15.58

 Neub\"user J.~ 14.95,\, 18.90

 Neumann B.~ 11.67,\, 17.101,\, 18.21,\lb19.65,\, 21.81,\, \textbf{A:}~1.14,\, 3.18,\, 11.68,\lb14.10,\, {\it 18.82}

Neumann H.~ 7.58,\, 17.101,\, 18.21--22,\lb19.65,\, 21.81,\,
\textbf{A:}~2.39,\, 8.26,\, 14.10

 Neumann von J.~ \textbf{A:}~{\it 16.58}

 Neumann P.\,M.~ 3.32,\,
5.38--39,\, 6.38,\lb8.50,\, 9.39--42,\, 11.69--71,\, 12.62--63,\lb14.46,\, 14.95,\, 15.64--66,\, 17.95--98,\lb19.60,\, 19.65,\, 20.30,\,
{\it 20.101},\,
 \textbf{A:}~4.39,\lb4.45, \, 5.40,\, 6.38,\, 8.49,\,
9.41,\, {\it 9.41},\lb11.70, \, 12.10

Nevmerzhitskaya N.\,S.~ \textbf{A:}~{\it 7.30}

 Newman M.\,F.~ 14.64,\, 15.41,\, 20.36, \lb\textbf{A:}~{\it 7.7},\, {\it 11.112},\, {\it 14.92}

 Newton B.~ \textbf{A:}~{\it 12.82}

Nguyen H.\,N.~  {\it 21.59}

Nickel W.~ \textbf{A:}~{\it 16.83}

Nicotera C.~ {\it 18.14},\, 21.98

Nies A.~ 21.120

Nikiforov V.\,A.~ \textbf{A:}~{\it 7.22}

Nikolov N.~
18.108,\, \textbf{A:}~{\it 7.37},\, {\it 18.75}

Nilradical AI~ {\it 16.68},\, {\it 16.69},\, {\it 21.38},\lb{\it 21.44},\, {\it 21.53}

Nnuga Mongo I.\,S.~
{\it 18.62},\, {\it 18.63},\, {\it 18.65},\lb {\it 20.120}

N\"obeling G.~ \textbf{A:}~{\it 2.36}

Noce M.~ 21.104

Norton S.~ \textbf{A:}~{\it 7.30}

 Noskov G.\,A.~ 6.31,\, 10.35--40,\, 10.42,\lb{\it 20.104},\, \textbf{A:}~{\it 1.70},\,
{\it 4.32},\, {\it
4.49},\, {\it 4.52}, \lb{\it 6.31}, \, 10.37,\, 10.41

Novikov P.\,S.~ \textbf{A:}~{\it 1.24}

Nunke R.~ 11.89

 Nuzhin Ya.\,N.~ 14.69,\, {\it 14.69},\, 15.67,\lb16.76,\, 19.61--63,\, {\it 19.63},\, {\it 21.76},\, 21.101,\lb\textbf{A:}~{\it 7.30},\, {\it 7.40},\,
{\it 14.49}

Nyberg-Brodda C.-F.~ {\it 17.99},\, 20.68--69

~

 \textbf{O}braztsov V.\,N.~ \textbf{A:}~{\it
6.46},\, {\it 8.32},\, {\it 9.79}, \linebreak
\makebox[3ex][r]{}{\it 9.81},\, {\it 10.14}, {\it 11.4},\, {\it
11.128},\, {\it 12.45}

 O'Brien E.\,A.~ 14.89,\, 18.90,\, 18.117,\lb\textbf{A:}~{\it 11.6},\, {\it 11.46},\, 11.112,\, {\it 14.92}

Occhipinti T.~ \textbf{A:}~{\it 11.2}

 O'Connor S.~ {\it 20.88}

Olive E.~ {\it 10.13}, {\it 20.88}

 Oliver A.~ {\it 16.3}

Oliveira J.\,R.~ {\it 15.19}

Oliynyk A.~ \textbf{A:}~{\it 16.85}

\newpage 
 Olshanskii A.\,Yu.~ 4.46,\, 4.48,\,
 5.42,\, 5.44,\lb8.52--55,\, 11.72--73,\, 12.64,\, 14.76,\lb15.68--70,\, 16.64,\, 17.99,\, 18.73--74,\lb19.111,\, 20.70,\, 21.102,\, \textbf{A:}~{\it
2.7},\, {\it 2.39},\lb{\it 2.41},\, {\it 2.45},\, {\it 2.61},\, {\it 2.66},\, {\it 2.73},\, 4.46--47,\lb{\it
4.61},\, {\it 4.73},\, {\it 5.22},\,
5.41,\, 5.43,\, {\it 5.69},\lb{\it 6.63},\, {\it 6.64},\, {\it 7.1},\, {\it 7.42},\,
{\it 8.7},\, 8.53,\, 8.56,\lb{\it 9.80},\, {\it 11.11},\, 11.74--75,\, {\it 11.74},\, {\it 11.75},\lb{\it
11.106},\, {\it 12.32},\, 14.104,\, {\it 16.101},\, {\it 17.79},\lb18.73,\, 18.75,\, {\it 19.75}

Olsson J.\,B.~ 8.51,\, 9.23,\, 13.43,\, {\it 21.59}

OpenAI~ {\it 12.20},\, {\it 15.42},\, {\it 17.94},\, {\it 18.87},\lb
 {\it 21.86}

 Ore O.~ 18.117

Orevkov S.~ {\it 14.102}

 Osin D.\,V.~ 21.80,\, \textbf{A:}~{\it 9.10},\, {\it 17.79}

 Ostylovski\u{\i} A.\,N.~ \textbf{A:}~7.42

Otmen N.~ {\it 14.53}

 Otto F.~ \textbf{A:}~{\it 9.30}

 Ould~Houcine A.~ {\it 5.47},\, {\it 10.49}

 Ozawa N.~ \textbf{A:}~{\it 15.7}

~

 \textbf{P}acifici E.~ {\it 20.78}

P\'alfy P.\,P.~ 12.62,\, 15.71--74,\, \textbf{A:}~{\it 10.63}

Palyutin E.\,A.~ 18.76,\, \textbf{A:}~{\it 3.22}

Pan J.~ {\it 21.8}

Pannone V.~ 17.100,\, 20.71

 Paras A.~ \textbf{A:}~{\it 11.26}

Parfenov V.\,A.~ 19.65

Park J.~ 20.54

Parker C.~ {\it 11.34},\, {\it 21.142}

Pasini A.~ \textbf{A:}~12.65

 Passi I.\,B.\,S.~ 16.66,\, \textbf{A:}~16.65,\, 16.104,\lb{\it 16.104}

Passman D.\,S.~ 12.29,\, 18.77--81,\, \textbf{A:}~{\it 4.28},\lb19.37

 Patne R.\,M.~ \textbf{A:}~18.82

 Pellegrini M.\,A.~ {\it 11.32},\, {\it 18.98},\, 18.117

Peng D.~ {\it 14.4},\, {\it 14.5},\, 21.103

 Pennington E.~ \textbf{A:}~{\it 5.17}

 Penttila T.~ \textbf{A:}~{\it 10.28}

P\'erez-Calabuig V.~ {\it 17.56}

Perkins P.~ \textbf{A:}~{\it 2.41}

 Perron B.~ {\it 8.11},\, \textbf{A:}~17.14

Petechuk V.\,M.~ 10.43--47,\, \textbf{A:}~{\it 8.46}

Petschick M.~ {\it 20.101},\, 21.104--106

Pettet A.~ \textbf{A:}~{\it 15.10}

Petyt H.~ {\it 19.92}

 Pfister G.~ \textbf{A:}~{\it 15.75}

 Phillips R.~ 11.107,\, 11.109,\, \textbf{A:}~{\it 5.1},
\lb11.106,\, 11.108

Pickel P.\,F.~ \textbf{A:}~4.53,\, {\it 4.53}

 Pink R.~  21.121,\, \textbf{A:}~{\it 1.75,\, 7.6}

Piskunov A.\,G.~ \textbf{A:}~{\it 8.63}

Piwek P.~ {\it 12.37}

 Platonov V.\,P.~ 1.74,\, 2.48,\, 2.53,\, 2.56,
\lb\textbf{A:}~{\it 1.32},\, 1.71--73,\,
{\it 1.71,\, 1.72},\lb1.75--76,\, {\it
2.16},\, {\it 2.46},\, 2.49--52,\lb2.54--55,\, 3.40--42,\, {\it 3.42}

Plesken W.~ 14.95

 Plotkin B.\,I.~ 3.43--49,\, 6.39,\, 14.70, \linebreak
\makebox[3ex][r]{}15.75--76,\, 16.66,\, 16.77,\, 17.101,\lb19.64--66,\, \textbf{A:}~1.69,\, 2.58--59,\, 2.61--63,\lb{\it
2.61},\, 2.64--65,\, 14.71,\, 15.75--76,\, 16.15

 Plotkin E.~ 15.76,\, 19.66,\, 20.72,\lb\textbf{A:}~{\it 15.75},\, 15.76

 Plunkett T.~ \textbf{A:}~{\it 11.98}

 Podufalov N.\,D.~ 9.43,\, 11.76--77,\, 12.66, \linebreak
\makebox[3ex][r]{}\textbf{A:}~{\it 4.23},\, 6.42,\, {\it 8.73},\,
9.43,\, {\it 9.43}, \lb10.48,\, {\it
10.48}

Poguntke W.~ 12.75

 Poizat B.~ 10.49,\, 16.78,\, 18.83

Poletskikh V.\,M.~ 8.58--60,\, 9.44,\, 12.67,\linebreak
\makebox[3ex][r]{}12.68

Ponomarenko I.~ 19.68,\, \textbf{A:}~19.67

Ponomar\"ev K.\,N.~ 11.78,\,
 12.69,\, 14.72,
\lb\textbf{A:}~11.79,\, {\it 11.79}

Popiel T.~ {\it 17.75}

Popov A.\,M.~ {\it 10.61}

Popov V.\,Yu.~ \textbf{A:}~{\it 2.40}

Popov Yu.\,D.~ 9.45

Pourgholi G.\,R.~ 19.25

Pouzet M.~ 19.18,\, \textbf{A:}~{\it 13.29}

Poznansky T.~ {\it 14.6},\, \textbf{A:}~{\it 15.7}

 Praeger C.\,E.~ 11.45,\, 11.70--71,\, 11.80, \linebreak
\makebox[3ex][r]{}{\it 11.80},\, 12.29,\, 14.46,\, 14.73,\, {\it 14.73},\lb15.47,\, 15.95,\, {\it 20.1},\, 20.1--3,\, 20.17,\lb21.23--24,\, \textbf{A:}~11.70,\, 18.23

Prajapati J.~ {\it 20.71}

Previtali A.~ {\it 11.71}

Pride S.\,J.~ \textbf{A:}~{\it 1.64},\, {\it 4.39},\, {\it 6.18},\, {\it 8.68}

 P{\v{r}}{\'{\i}}hoda P.~ \textbf{A:}~{\it 8.34}

Prishchepov M.\,I.~ \textbf{A:}~{\it 8.10}

Procesi C.~ 14.37

Protasov I.\,V.~ 8.59,\, 8.62,\, 9.45,
 9.47,\lb11.81,\, 13.25--25,\, 13.44,\, 13.48--49,\lb15.77--78,\, 15.80,\, {\it 15.80},\, 17.51,\lb17.102,\, 20.73,\, 21.107,\, \textbf{A:}~{\it
1.76},\, {\it 3.37},\lb8.61,\, 8.63,\, 9.46,\,
{\it 9.46},\, 9.48--49,\lb13.26,\, 13.45--47,\, {\it 13.45},\, 15.79

Pudl\'ak P.~ \textbf{A:}~{\it 5.3}

Puglisi O.~ {\it 12.27},\, {\it 12.28},\, {\it 18.110}

 Pyatetski\u{\i}-Shapiro I.\,I.~ 11.9

Pyber L.~ 12.56,\, 14.53--54, {\it 14.73},\,
\lb14.74--76,\, 17.41,\, {\it 19.11},\, 20.74--77,\lb{\it 20.77},\, \textbf{A:}~{\it 7.17}

 Pypka A.\,A.~ \textbf{A:}~{\it 18.91}

~

\textbf{Q}ian G.~ 20.78--80,\, {\it 20.78},\, {\it 20.80},\lb21.108--109,\, \, \textbf{A:}~{\it 15.62},\, {\it 17.111}

Qiu D.~ {\it 18.16}

Quick M.~ {\it 19.22}

Quillen D.~ \textbf{A:}~{\it 13.36}

~

\textbf{R}adulescu F.~ 14.9

Raimbault J.~ \textbf{A:}~{\it 15.49}

Rapaport Strasser E.~ 20.32

Raptis E.~ \textbf{A:}~{\it 10.33}

Razborov A.\,A.~ \textbf{A:}~{\it 1.18},\, {\it 6.25}

Razmyslov Yu.\,P.~ \textbf{A:}~{\it 5.56},\, 9.50

Rehmann U.~ 10.36

Reibol'd A.\,V.~ \textbf{A:}~9.62

Reid C.~ {\it 19.10},\, {\it 19.22},\, 19.69--73

Reinfeldt C.~ {\it 14.19}

 Remeslennikov V.\,N.~ 4.50,\, 5.47--48,\lb7.58,\, 11.65--66,\, 11.83--84,\, 13.39,\lb13.41--42,\, 14.18--22,\, 14.42,\, 16.70,\lb17.104,\, \textbf{A:}~{\it 2.16,\, 2.17},\, {\it 2.46},\, 4.49,\lb4.51--53,\, {\it 4.62},\, 5.46,\,
5.49,\, 11.82,\lb13.39,\, 13.40,\, {\it 14.71},\, 15.86

Reznichenko E.~ {\it 2.48}, {\it 18.71}

 Revin D.\,O.~ 17.43,\, {\it 17.43},\, 17.45,\, 17.103,\lb18.68,\, 18.84,\, 19.74,\, 20.54,\, 20.81--83,\lb21.57,\, 21.110--112,\, \textbf{A:}~{\it 3.62},\, {\it 5.65},\lb{\it 13.33},\, {\it 14.62},\, 17.44--45,\, {\it 17.44--45},\lb18.31--32,\, {\it 19.50}

 Rhemtulla A.\,H.~ {\it 2.82},\,
15.1,\, 17.31,\, 18.4,\, \lb\textbf{A:}~10.30

 Ribes L.~ \textbf{A:}~{\it 8.70},\, 12.70--71,\, {\it
12.70}

Ribnere E.~ {\it 15.75}

Riedl J.\,M.~ 15.72

 Riese U.~ {\it 9.23}

Rips E.~ {\it 14.70},\, 20.84,\, \textbf{A:}~{\it 2.59},\, {\it 11.52}

Ritter J.~ \textbf{A:}~{\it 12.1}

Rivas C.~ {\it 16.51}

Rizzoli A.~ {\it 7.15},\, {\it 21.3},\, {\it 21.29},\, {\it 21.68}

Robertson D.~ \textbf{A:}~{\it 19.40}

 Robertson E.\,F.~ 8.12,\, \textbf{A:}~8.12,\, {\it 8.12}

Robinson D.\,J.\,S.~ 9.51--53,\, \textbf{A:}~{\it 6.20},\lb9.54

Robinson G.\,R.~ 13.43,\, 21.92,\, 21.113,\lb\textbf{A:}~{\it 4.14},\, {\it 11.35}

 Roggenkamp K.\,W.\,\, 4.55--56,\, 9.55,\lb\textbf{A:}~4.54,\, 12.80

Roller M.\,A.~ {\it 20.118}

 Roman'kov V.\,A.~ 11.86--87,\, 16.92,\lb18.72,\, {\it 18.73},\, 18.85,\, 18.87,\, {\it 18.89},\lb\textbf{A:}~{\it 2.64},\, {\it 3.8},\,
{\it 4.49},\, {\it 8.71},\, {\it
11.54},\lb{\it 11.55},\, {\it 11.88},\, {\it 17.82},\, {\it 17.86},\, {\it 18.73},\lb18.86,\, {\it 18.86}

 Romanovski\u{\i} N.\,S.~ {\it 3.47},\,
11.85--86,\lb16.80,\, 17.105,\, 19.64,\,
\textbf{A:}~{\it 2.54},\, {\it 5.29},\lb{\it 6.36}

Roseblade J.\,E.~ 11.31

Rosenberger G.~ 14.14,\, 20.22,\, 20.41,\lb20.85,\, \textbf{A:}~{\it 9.2}

Rosendal C.~ {\it 18.36},\, {\it 21.15}

Ross K.~ \textbf{A:}~3.35

 Rothaus O.\,S.~ {\it 17.94},\, 18.87

Rowley P.\,J.~ 12.62

 Rozhkov A.\,V.~ 9.77,\, 13.55,\, 21.131,\lb\textbf{A:}~{\it 6.58},\, {\it 8.66},\, {\it 13.21},\,  16.79,\, {\it 16.79}

Rubin M.~ 12.13

Rud'ko V.\,P.~ 12.21

Rudloff Ch.~ 11.32

Rudvalis A.~ \textbf{A:}~4.77

Rukolaine A.\,V.~ 10.50

Rump W.~ 20.92

Rumyantsev A.\,K.~ 8.64,\, \textbf{A:}~{\it 6.43}

Rumynin D.~ {\it 4.55},\, 21.60

Rundstr\"om M.~ {\it 21.24}

Russo F.\,G.~ \textbf{A:}~{\it 9.32}

Ruzsa I.\,Z.~ \textbf{A:}~{\it 15.79}

Rychkov S.\,V.~ 10.51--55,\,
11.89--90

~

\textbf{S}abbagh G.~ \textbf{A:}~{\it 12.14}

Sabatini L.~  21.26,\, 21.114

 Sadovski\u{\i} L.\,E.~ 2.67--68,\, \textbf{A:}~2.66

 Safonov V.\,G.~ \textbf{A:}~{\it 14.80}

Safonova I.\,N.~ \textbf{A:}~{\it 19.87},\, {\it 19.88}

Sageev M~ {\it 14.6}

Sakai S.~ 14.9

Saksonov A.\,I.~ \textbf{A:}~3.50,\, 3.51

Salle, de la M.~ {\it 12.20},\, {\it 15.42},\, \textbf{A:}~{\it 19.98}

Salomon E.~ \textbf{A:}~{\it 12.96}

Sambale B.~ {\it 20.27},\, 21.115,\, {\it 21.115},\lb\textbf{A:}~{\it 12.80},\, 19.75

 Sangroniz J.~ 14.44,\, \textbf{A:}~{\it 11.43},\, 14.77,\lb{\it 17.77}

Sanov I.\,N.~ 20.36

Santos T.\,M.\,G.~ {\it 15.19},\, 21.41--42

Sanus L.~ {\it 20.78}

 Sapir M.\,V.~ 17.107,\, 18.88,\, 19.16,\, 19.76,\lb21.144,\, \textbf{A:}~{\it 5.22},\, {\it 8.7},\, 17.108,\, {\it 17.108}

Sarkisyan O.\,A.~ \textbf{A:}~{\it 4.4}

Sater R.~ {\it 10.34},\, {\it 19.57},\, {\it 19.58},\, {\it 21.87}

Sauer R.~ {\it 21.117}

 Savel'eva N.\,V.~ \textbf{A:}~{\it 13.50}

Saxl J.~ 14.69,\, 21.110,\, \textbf{A:}~{\it 6.6},\, {\it 9.21},\, 9.56

Schacher M.~ \textbf{A:}~{\it 7.16}

Scharf T.~ 21.113

Scheerer H.~ \textbf{A:}~9.33

 Scheiderer C.~ \textbf{A:}~{\it 9.33},\, {\it 12.70}

Schesler E.~ {\it 15.14},\, 20.43,\, 21.116--119\lb{\it 21.117}

 Schick T.~
\textbf{A:}~{\it 12.47},\, {\it 14.102},\, 16.52

 Schleimer S.~ \textbf{A:}~{\it 16.109}

 Schlitt G.~ \textbf{A:}~17.24

 Schmid P.~ {\it 9.23},\, 18.114,\, 18.116,\, 20.86,\lb\textbf{A:}~{\it 4.12},\, 17.2,\, 18.115

Schmidt R.\,A.~ \textbf{A:}~{\it 7.40}

Sch\"onheim J.~ 21.93

 Schupp P.\,E.~ 16.110,\, 19.77,\, \textbf{A:}~{\it 4.45},\lb{\it 4.60},\, 5.68

Schwerdtfeger H.~ \textbf{A:}~4.70

 Scoppola C.\,M.~ 17.47,\, 18.34,\, \textbf{A:}~{\it 11.6}

Scott F.~ {\it 20.31}

 Scott L.\,L.~ \textbf{A:}~{\it 7.17},\, {\it 16.55}

Scott P.~ \textbf{A:}~5.53

Seal D.\,J.~ \textbf{A:}~{\it 8.10}

 Segal D.~ {\it 11.5},\, 12.56,\, 17.109--110,\, {\it 20.77},\lb21.105,\, 21.120,\, \, \textbf{A:}~3.2,\, {\it
4.52,\, 4.53},\lb{\it 5.6},\, {\it 7.37,\,
8.70},\, 11.57,\, {\it 11.57},\, {\it 12.90},\lb{\it 18.75}

 Segev Y.~ \textbf{A:}~{\it 11.52}

Seitz G.~ \textbf{A:}~{\it 7.48}

Seksenbaev K.~ \textbf{A:}~{\it 1.84}

Sela Z.~ {\it 9.66},\, {\it 13.12},\, {\it 14.19},\, {\it 14.20},\lb\textbf{A:}~9.67

Selberg A.~ \textbf{A:}~2.49--50

 Sel'kin M.\,V.~ \textbf{A:}~16.82

 Semenchuk V.\,N.~ 12.72,\, \textbf{A:}~{\it 8.87},\, 11.91

Sem\"enov Yu.\,S.~ \textbf{A:}~{\it 6.22}

Semidetnov A.~ {\it 14.85},\, {\it 16.11},\, {\it 17.32},\lb{\it 17.47},\, {\it 19.94}

Senashov V.\,I.~ 10.75

 Serbin D.~ {\it 16.70}

 Sereda V.\,A.~ \textbf{A:}~{\it 11.101}

Ser\"ezhkin V.\,N.~ \textbf{A:}~{\it 4.29}

Sergiyenko V.\,I.~ \textbf{A:}~7.44

Serre J.-P.~ 15.57,\, 18.116,\, 19.78--79

 Sesekin N.\,F.~ 5.1,\, 6.1--2,\, \textbf{A:}~1.78,
{\it 1.78},\lb2.69--70,\, 4.57--58,\,
5.1--2,\, 7.6

Schaeffer Fry A.\,A.~ {\it 9.23}

Shafarevich I.\,R.~ 15.6,\, \textbf{A:}~3.40

Shahabi Shojaei M.\,A.~ \textbf{A:}~{\it 3.14}

Shahshahani M.~ \textbf{A:}~{\it 19.50}

 Shahverdi H.~ {\it 16.1},\, \textbf{A:}~{\it 16.1}

Shakhova S.\,A.~ \textbf{A:}~{\it 10.9}

 Shalev A.~ 12.95,\, {\it 12.95},\,
14.89,\, 18.117,\lb19.93,\, 20.102,\, 21.27,\, \textbf{A:}~{\it 7.17},\,  12.94,\lb13.36,\, 13.56

 Shaptala V.\,P.~ 11.115,\, \textbf{A:}~13.62

 Sharifi H.~ {\it 16.37}

Sheinfeld O.~ {\it 20.74}

 Shelah S.~ {\it 9.39},\, {\it 9.41},\, {\it 9.42},\, {\it 15.70},\lb17.121--122,\, ,\, {\it 20.40},\,  \textbf{A:}~{\it 1.66},\,
{\it 3.19},\lb  {\it 12.44},\, {\it 14.71}

 Shemetkov L.\,A.~ 3.60,\,
6.51,\, 9.75,\, 10.73,\lb11.117--121,\, 14.81,\, 14.99,\, 15.98,\lb17.38--39,\, \textbf{A:}~{\it 2.1},\, 3.59,\, {\it 3.59},\lb3.61--62,\,
4.71,\, {\it 4.71},\, 5.64--65,\, 6.52,\lb8.87,\, 9.58--60,\, 9.73--74,\, {\it
10.56},\, 10.72,\lb{\it 10.72},\, 11.24,\, {\it 11.24},\, 12.96--97,\, 14.80,\lb14.99

Shemetkova O.\,L.~ {\it 19.38}

Shen R.~ \textbf{A:}~{\it 18.112}

Shepelev V.\,D.~ {\it 17.43}

 Shepherd R.~ \textbf{A:}~16.103

Shevrin L.\,N.~ 2.81,\, 11.116,\, \textbf{A:}~1.81

Shi Shengming~ \textbf{A:}~12.26

 Shi W.\,J.~ 12.37,\, {\it 12.38},\, 13.64--65,\lb15.99,\, 17.123,\, 20.87,\,  {\it 21.134},\lb\textbf{A:}~12.38,\, {\it 12.38},\, 12.39,\, {\it 12.39},\, 13.63,\lb{\it 16.1},\, 16.106--107,\, 18.112,\, {\it 18.112},\lb19.80

Shilov A.\,V.~ 21.77

Shimura G.~ 8.83

Shirjian F.~ {\it 20.38},\, 21.59

Shirshov A.\,I.~ 2.82

Shirvani M.~ 8.2

 Shkuratski\u{\i} A.\,I.~ \textbf{A:}~{\it 4.41}

Shl\"epkin A.\,A.~ {\it 20.98}

Shl\"epkin A.\,K.~ 14.100--101,\, 15.100--101,\lb18.113,\, \textbf{A:}~{\it 12.38}

 Shmel'kin A.\,L.~ 1.86--87,\, 4.72,\,
11.122, \lb\textbf{A:}~1.82--85,\, 1.88--89,\, {\it
1.89},\, {\it 2.43},\lb4.73

Short M.~ 11.123,\, \textbf{A:}~11.47

Shpektorov S.\,V.~ \textbf{A:}~{\it 8.37}

\newpage 

 Shpilrain V.~ 11.124,\, 14.85,\, 14.102,\lb15.63,\, 15.102--103,\,
16.9--10,\, 16.108,\lb16.110,\, 17.124,\, 20.88,\, {\it 20.88},\lb\textbf{A:}\,{\it
5.49},\, {\it 10.24},\, 12.98,\, 13.66,\, 14.88,\lb14.102,\, 15.63,\, 16.109,\, 19.81

Shramov C.~ 21.121

Shult E.~ 5.30

 Shumyatski\u{\i} P.\,V.\, {\it 10.77},\, ,\, {\it 11.116},\lb11.125,\,
12.100,\, 15.104,\,
17.125--127,\lb{\it 17.125},\, 18.66--67,\, 18.110,\lb18.117--119,\, 19.82--83,\, 20.89--90,\lb21.122,\, \textbf{A:}~11.126,\, {\it 13.58},\, {\it 17.72}

 Shunkov V.\,P.~ 2.84,\, 4.74--75,\, 5.67,
6.55,\lb6.56,\, 6.59--62,\,
9.76--78,\, 9.83--84, \lb10.74--75,\,
10.77--78,\, 11.56,\, 11.127, \lb12.101,\, 13.67,\, 14.83,\, 14.100,\, 16.111,
\lb\textbf{A:}~{\it 1.23},\, {\it
1.81},\, 1.90, \,{\it 1.90},\, 2.83,\,
{\it 3.56}, \lb3.64,\, 4.74,\, 4.76, \,6.53--54,\, 6.57--58,
\lb6.63--64,\, {\it 6.63},\, 7.42,\,
8.66,\, 9.79--82, \lb10.66,\, 10.76,\,
12.102,\, 14.103--104

Shute G.~ 11.109,\, \textbf{A:}~{\it 1.75},\, {\it 7.6}

 Sibley D.~ 12.92

 Sidki S.~ 15.19,\, 16.74,\, 18.4,\, 18.90,\, 19.12,\lb21.41--42,\, \textbf{A:}~{\it 15.14}

Siemons J.~ {\it 16.29}

Sieradski A.\,J.~ \textbf{A:}~{\it 4.54}

Silberger D.\,M.~ 10.32

Sim H.-S.~ 12.71,\, 21.87

Simonyan L.\,A.~ 19.65,\, \textbf{A:}~{\it 2.58}

Sims C.~ 11.19,\, 15.66

Siniora D.~ {\it 17.76}

Sisto A.~ {\it 18.111}

Sitnikov V.\,M.~ \textbf{A:}~{\it 3.6}

 Skiba A.\,N.~ 9.57,\, 10.57--58,\, 11.92,\lb12.73,\, 13.51, \,14.47,\, 14.78--79,\, 14.81,\lb17.37,\, 17.112,\, 18.92--93,\, 19.86,\lb\textbf{A:}~{\it
6.52},\, 9.58--60,\, 10.56,\, {\it 10.56},\lb{\it 10.72},\, {\it 11.25},\, {\it
11.91},\, 12.74,\, {\it 12.74},\lb13.50,\, 14.80,\, 15.81,\, 16.82,\, 17.111,\lb18.91,\, 19.84--85,\, 19.87--88,\, {\it 19.87},\lb{\it 19.88}

Sklinos R.~ 19.64,\, 20.91

Skresanov S.\,V.~ {\it 11.96},\, {\it 15.53},\, {\it 17.76},\lb{\it 20.74},\, \textbf{A:}~{\it 12.79},\, {\it 14.60},\, {\it 15.79},\, {\it 17.73},\lb{\it 19.67}

Skutin A.~ \textbf{A:}~{\it 4.69}

 Slobodianiuk S.~ \textbf{A:}~{\it 13.45}

 Slobodsko\u{\i} A.\,M.~ \textbf{A:}~{\it 1.25}

Slonim Z.~ {\it 20.79}

 Smirnov D.\,M.~ 7.45,\, 8.64,\, 9.61,\, 12.75,\linebreak
\makebox[3ex][r]{}13.52,\, \textbf{A:}~{\it 1.44,\, 1.45},\, {\it
1.50},\, 2.71--72,\linebreak
\makebox[3ex][r]{}3.52--53,\, 4.59--60,\,
5.50--51, 6.43,\, 9.62,\linebreak
\makebox[3ex][r]{}{\it 9.62},\,
11.93,\, {\it 11.93}

Smith D.\,B.~ \textbf{A:}~{\it 3.18}

 Smith H.~ 15.85,\,
\textbf{A:}~12.45,\, {\it 18.33}

Smith M.\,K.~ 12.29

Smith S.\,M.~ {\it 19.22},\, 19.73,\, 19.89

Smoktunowicz A.~ 19.90,\, 20.92,\, {\it 20.92},\lb21.123,\, \textbf{A:}~19.90

Snopce I.~ \textbf{A:}~{\it 17.19}

Soifer G.\,A.~ \textbf{A:}~{\it 2.55}

Sokolov V.\,G.~ 7.58,\, \textbf{A:}~{\it 2.16},\, {\it 2.18},\lb{\it 2.79}

Solecki S.~ {\it 18.36}

 Solitar D.~ \textbf{A:}~{\it 3.15},\, 16.2

 Solomon R.\,M.~ 19.91--92,\, \textbf{A:}~{\it 16.1}

 Sonkin D.~ \textbf{A:}~{\it 16.101}

 Soroko I.~ 21.124--129

 Sosnovski\u{\i} Yu.\,V.~ 15.83--85,\, 16.83,\linebreak
\makebox[3ex][r]{}\textbf{A:}~{\it 8.47},\, 16.83

 Soul\'e C.~ 10.36

 Sozutov A.\,I.~ 8.67,\, 10.59--62,\, {\it 10.61},\lb13.53--54,\, 14.83,\, 15.82,\, \, 20.93--96,\lb\textbf{A:}~{\it 1.39},\, {\it
6.54}, \,{\it 6.63},\, 8.66,\, {\it 10.76},\lb{\it 11.11},\, {\it 11.13},\, {\it 11.101},\, 12.76,\, 13.54,\lb14.82,\, {\it 14.103,\, 14.104},\, {\it 16.15},\, {\it 17.3},\lb{\it 17.71},\, 18.94--95

Spaltenstein N.~ \textbf{A:}~{\it 8.13}

Sp\"ath B.~ {\it 8.51}

 Specht W.~ \textbf{A:}~2.39

 Spellman D.~ 11.19,\, 12.8--9,\, 13.39,\,
\lb18.35,\, \textbf{A:}~11.20,\, 13.18,\, 13.39,\, 14.32,\lb15.86,\, {\it 15.86}

 Spiezia F.~ {\it 14.46}

 Spiga P.~ {\it 9.69},\, {\it 11.71},\,
 {\it 14.73},\, 19.93,\lb20.17,\, \textbf{A:}~18.23

Springer T.~ 19.79

 Sprowl T.~ \textbf{A:}~{\it 7.17}

 Stafford J.\,T.~ 10.42

 Staiger M.~ {\it 19.8}

Stallings J.~ \textbf{A:}~1.13

Stammbach U.~ {\it 12.41}

 Staroletov A.\,M.~ {\it 20.58},\, \textbf{A:}~{\it 17.73}

Starostin A.\,I.~ 2.74,\, 5.14,\, 12.15,\lb\textbf{A:}~1.80,\, 2.1,\, 2.30,\,
2.73,\, 12.77,\, 12.78

Steinberg R.~ 7.28,\, \textbf{A:}~2.51

Stellmacher B.~ \textbf{A:}~{\it 6.7}

Stepanov A.\,V.~ \textbf{A:}~{\it 7.40}

St\"epin A.\,M.~ 9.7

St\"ohr R.~ \textbf{A:}~{\it 11.47}

Stonehewer S.~ \textbf{A:}~{\it 5.17}

 Storozhenko D.\,Yu.~ \textbf{A:}~{\it 18.91}

 Storozhev A.~ 16.64,\,
\textbf{A:}~{\it 2.72},\, {\it 4.46},\, {\it 4.73}, \linebreak
\makebox[3ex][r]{}{\it 6.16},\,
{\it 8.18},\, {\it 8.20},\, {\it 8.22},\, {\it 12.31}

Strebel R.~ \textbf{A:}~6.35

 Str\"ungmann L.~ \textbf{A:}~{\it 1.66}

 Strunkov S.\,P.~ 3.55,\, 11.95--96,\lb11.98--100,\, \textbf{A:}~2.75,\, 3.56,\, 11.94,\lb11.97--98,\, 12.79--80

 Suchkov N.\,M.~ 15.82,\, 15.87,\, 18.96--97,\lb20.97--98,\, {\it 20.98},\, \textbf{A:}~{\it 10.76}

Sun B.~ {\it 15.45}

Sun Y.~ {\it 20.79}

Sun Z.-W.~ 20.99--100,\, 21.130

\v{S}uni\'{c} Z.~ 15.12--14,\, 15.16,\, {\it 16.34},\lb  \textbf{A:}~15.15,\, 15.17--18

Suprunenko D.\,A.~ \textbf{A:}~2.76--77

Suprunenko I.\,D.~ 11.32

Sureaux S.~ {\it 21.97}

 Sushchanski\u{\i} V.\,I.~ 12.81,\, 15.19, \,16.34,\linebreak
\makebox[3ex][r]{}16.84,\, \textbf{A:}~{\it 1.37},\, {\it 2.83},\, 12.82,\, 16.85,\linebreak
\makebox[3ex][r]{}16.86

Suslin A.\,A.~ 10.36

 Suzuki M.~ 19.13,\, 20.93,\, \textbf{A:}~11.11,\, {\it 11.11},\linebreak
\makebox[3ex][r]{}15.33

 Swan J.~ 18.89

Swan R.\,G.~ \textbf{A:}~{\it 1.34}

 Sysak Ya.\,P.~ 8.67,\, 9.65,\, 10.64,\,
\textbf{A:}~{\it 1.36},\lb{\it 2.70},\, {\it
9.54},\, 9.63--64,\, 10.63,\, {\it 18.95}

 Syskin S.\,A.~ 12.85--86,\, 15.98,\, \textbf{A:}~{\it 3.6},
\lb{\it 3.59},\, {\it 4.16},\, {\it
4.22},\, 4.61,\, {\it 5.41},\, 7.48, \linebreak
\makebox[3ex][r]{}12.84

 Szab\'{o} E.~ {\it 14.73}

~

\textbf{T}ae Young L.~ {\it 21.43}

Taimanov A.\,D.~ 9.66,\, \textbf{A:}~9.67

Taketa K.~ 19.51

 Tamburini M.\,C.~ 15.67,\, 18.98,\, {\it 18.98},\lb\textbf{A:}~{\it 14.49},\, {\it 17.116}

Tang C.\,Y.~ 8.72,\, \textbf{A:}~8.68--71

Tanushevski S.~ \textbf{A:}~{\it 17.19}

Tao, J.~ \textbf{A:}~{\it 19.81}

 Tao T.~ 18.108

T\u{a}rn\u{a}uceanu M.~ 21.97,\, 21.115

 Tarski A.~ 10.12,\, 19.59,\, \textbf{A:}~1.18,\, 1.29,\lb1.68,\, 4.62--64,\, 9.67

Tavgen' O.\,I.~ 20.77

Taysnyov D.\,A.~ {\it 20.98}

Tent J.\,F.~ {\it 20.4}

 Tent K.~ 21.120,\, \textbf{A:}~{\it 11.52}

Teragaito M.~ \textbf{A:}~{\it 16.49}

Tertooy S.~ {\it 14.54}

Tessera R.~ \textbf{A:}~{\it 19.98}

Thillaisundaram A.~ 20.101--102

Thom A.~ {\it 16.69},\, {\it 16.70},\, 20.42,\, {\it 21.86},\lb 21.118

Thomas A.\,R.~ {\it 17.42},\, {\it 21.29}

Thomas S.~ {\it 9.39},\, 11.109,\, {\bf
A:}~{\it 1.75},\, {\it 7.6},\, \linebreak
\makebox[3ex][r]{}{\it
8.14},\, {\it 15.8}

Thomas V.\,Z.~ {\it 21.32}

 Thompson J.\,G.~ 4.65,\, 5.30,\, 6.10,\, 8.74, \linebreak
\makebox[3ex][r]{}8.75,\, 9.24,\, 12.37--38,\, 14.76,\,
16.33,\lb16.95,\, 20.65,\, 21.66,\, {\it 21.134},\, \textbf{A:}~3.27,\lb4.22,\, 8.73

 Thompson R.~ 12.20,\,
 15.42,\, 19.76,\lb\textbf{A:}~{\it 1.14},\, 19.40

Thumm A.~ {\it 16.9}

Thurston W.\,P.~ 19.42,\, \textbf{A:}~{\it 4.51}

Tiep P.\,H.~ {\it 9.23},\, 11.34,\, 18.117,\, {\it 20.74},\lb  21.27,\, 21.82

Tikhonenko T.\,V.~ {\it 10.34}

Timmesfeld F.~ 8.43

 Timofeenko A.\,V.~ 9.76,\, 13.55,\, {\it 15.67},\lb21.131--132,\, \, \textbf{A:}~{\it 6.58},\, {\it 7.30},\, 11.101

 Timoshenko E.\,I.~ {\it 4.33},\, 14.84--85,\,
15.88,\lb16.87,\, 17.104,\, 19.94--97,\, 20.103--104,\lb\textbf{A:}~{\it 2.16},\, {\it
2.20},\, {\it 13.66},\, {\it
14.88}

 Tits J.~ 12.95,\, 14.38,\, 19.79,\, 20.110,\lb\textbf{A:}~{\it 1.88},\, {\it
2.52},\, 3.42

 Tokarenko A.\,I.~ {\it 3.47},\, \textbf{A:}~{\it 1.69},\, {\it
2.46}

 Tolstykh V.\,A.~ 16.9--10,\, 16.88--94

Tomkinson M.\,J.~ 10.17--18,\, \textbf{A:}~{\it 6.8}

 Tong-Viet H.\,P.~ 18.99,\, {\it 19.11},\, {\it 20.79},\lb{\it 21.59},\, {\it 21.92},\, \textbf{A:}~{\it 11.8}

Topale A.\,G.~ 14.4--5,\, {\it 14.4}

Tornier S.~ {\it 19.22}

Toti T.~ {\it 14.53}

Touikan N.~ {\it 19.41}

 Trabelsi N.~ 18.100,\, 20.105

Tracey G.~  {\it 20.1},\, 20.1--3,\, 21.104

Trakhtman A.\,N.~ \textbf{A:}~{\it 2.41}

Trappeniers S.~ {\it 20.92}

 Traustason G.~ {\it 2.82},\,
16.96--97,\lb20.106--107,\, {\it 20.106},\, \textbf{A:}~{\it 9.50},\, {\it 17.12}

Trofimov P.\,I.~ 2.78,\, 3.57

 Trofimov V.\,I.~ 7.49,\, 12.87--89,\, 15.89--90,\linebreak
\makebox[3ex][r]{}{\it 15.90},\, \textbf{A:}~{\it 18.94},\, 19.98

Troitsky E.\,V.~ 19.28--29

Trombetti M.~ {\it 5.38},\, {\it 6.30},\, {\it 11.9},\, {\it 13.23},\, {\it 17.113},\, 21.133

Tsalenko M.\,S.~ \textbf{A:}~2.79

Tsang C.~ 20.108--109,\, \textbf{A:}~{\it 19.90}

Tsurkov A.~ {\it 15.76},\, {\it 20.72},\, \textbf{A:}~{\it 15.76}

Tsybenko Yu.\,V.~ \textbf{A:}~{\it 8.61},\, 9.49,\, {\it 9.49}

Tucker T.~ 19.89

Tulenbayev M.\,S.~ 10.36

 T\raisebox{-0.1ex}{$\stackrel{\circ}{{\rm u}}$}ma J.~
\textbf{A:}~{\it 5.3}

Turau V.~ 21.65

Tursunbaev A.~ {\it 11.124}

 Turull A.~ 5.30,\, 16.98--99,\, 18.102,\lb21.65--66

 Tushev A.\,V.~ 15.91,\, 18.103,\, \textbf{A:}~{\it 6.12},\lb{\it
9.12},\, {\it 12.46},
 {\it 12.90}

Tyrer Jones J.\,M.~ 21.47

 Tyutyanov V.\,N.~ {\it 10.34},\, 17.112,\, {\it 17.112},\lb19.57--58,\, {\it 19.100},\, \textbf{A:}~{\it 8.31},\, 18.30,\lb{\it 18.30},\, {\it 19.84},\, {\it 19.87},\, {\it 19.88},\, {\it 19.109}

~

 \textbf{U}lam S.\,M.~ \textbf{A:}~15.8

Umirbaev U.\,U.~ {\it 3.3}

Urban J.~ {\it 4.75},\, {\it 9.4},\, {\it 10.62},\, {\it 12.40},\lb{\it 11.116},\, {\it 13.42},\, {\it 14.72},\, {\it 15.92},\, {\it 16.20},\lb  {\it 16.38},\,  {\it 16.87},\, {\it 17.33},\, {\it  18.18},\, {\it 18.92},\lb {\it 19.61},\, {\it 19.62},\,  {\it 19.93},\, {\it 20.1},\, {\it 20.90},\lb {\it 20.108},\, {\it 20.122},\, {\it 21.60}

Uspenskii V.\,V.~ 21.103

~

\textbf{V}aillant A.\,G.~ \textbf{A:}~14.27

Valiev M.\,K.~ \textbf{A:}~4.20,\, {\it 1.17}

van Beek M.~ {\it 11.34}

 van der Kallen W.~ 10.38

 van~der~Waall R.\,W.~ \textbf{A:}~{\it 7.48}

van\,Doorn, W~ {\it 3.46},\, {\it 18.50},\, {\it 19.25},\lb{\it 20.125},\, {\it 21.8},\, {\it 21.24},\, {\it 21.147},\, {\it 21.150}

Vannacci M.~ {\it 14.53},\, {\it 19.10},\,  {\it 19.22},\lb 20.110--111

Vannier J.\,P.~ {\it 8.11},\, \textbf{A:}~17.14

Vapne Yu.\,E.~ \textbf{A:}~{\it 3.29}

Varsos D.~ \textbf{A:}~{\it 10.33}

 Vaserstein L.\,N.~ 14.15,\, 16.21--23

 Vasil'ev A.\,F.~ 14.28--29,\, 17.37--39,\lb{\it 17.112},\, 18.29,\, 19.99--100,\, 20.66,\lb20.112,\,  \textbf{A:}~11.21,\, {\it 11.24},\, 12.7,\lb15.38--39,\, 18.30

 Vasil'ev A.\,V.~ 16.26,\, 20.31,\, {\it 20.37},\lb20.58,\, {\it 20.58},\, 21.134--135,\, \textbf{A:}~{\it 12.38},\lb{\it 12.39},\, 16.24--25,\, 16.27,\, 17.36,\, {\it 17.73},\lb19.101

 Vasil'eva T.\,I.~ 19.99--100,\, 20.112,\lb\textbf{A:}~18.30

 Vaughan-Lee M.\,R.~ 8.4,\,
15.31,\, 15.41,\linebreak
\makebox[3ex][r]{}20.30,\, 20.113--114,\, \textbf{A:}~{\it 2.39},\, {\it 2.41},\lb{\it
5.56},\, 8.5,\, {\it 9.50},\, {\it 11.6},\, {\it 11.112}

 Vdovin E.\,P.~ 15.40,\, 16.28--29,\, 17.41--43,\lb17.45,\, {\it 17.41},\, \textbf{A:}~{\it 3.62},\, {\it 4.27,\, {\it 5.65},\lb13.33},\, {\it 14.82},\, 17.40,\, 17.44--45,\lb{\it 17.44--45},\, 18.31--32

 Vedernikov V.\,A.~ 16.30,\, 16.32,\, \textbf{A:}~3.6,\linebreak
\makebox[3ex][r]{}{\it 9.58},\, 16.31

 Vendramin L.~ 19.90,\, \textbf{A:}~19.49,\, 19.90

Venkataraman G.~ 15.66

 Ventura E.~ {\it 16.34},\, 19.102--105,\lb \textbf{A:}~{\it 10.26}

Vera-Lopez A.~ \textbf{A:}~{\it 17.77}

Veretennikov B.\,M.~ \textbf{A:}~17.46,\, 18.33

 Verma D.\,N.~ 16.39

Verret G.~ 20.21

Vershik A.\,M.~ 12.28

Vilyatser V.\,G.~ 2.5--6,\, 11.23,\, \textbf{A:}~3.7--8

 Vir\'ag B.~ \textbf{A:}~{\it 16.74}

Viudez F.~ {\it 19.57},\, {\it 19.58}

Vogtmann K.~ 16.90

Volochkov A.\,A.~ \textbf{A:}~{\it 14.63}

 Vorob'yev N.\,T.~ {\it 8.30},\, 11.25,\,
\, 14.30--31,\lb\textbf{A:}~{\it 9.18},\, {\it
9.58},\, 11.24--25,\, {\it 11.25},\lb{\it 13.50}

 Voskresenski\u{\i} V.\,E.~ \textbf{A:}~5.7,\, {\it 5.7}

Vovsi S.~ 3.43,\, 13.17,\, \textbf{A:}~{\it 2.61}

 Vsemirnov M.\,A.~ {\it 14.69},\, 18.98

Vyas R.~ \textbf{A:}~19.23--24

~

 \textbf{W}aall van der R.\,W.~ \textbf{A:}~{\it 7.48}

 Wagner A.~ \textbf{A:}~{\it 4.29}

 Wagner F.~ {\it 12.48},\, 15.37,\, 19.106,\lb\textbf{A:}~{\it 11.103},\,
13.16

Wagon S.~ 10.12

Wales D.\,B.~ 17.16

 Walker A.~ 18.42

 Wall G.\,E.~ 11.105,\, 15.65,\, \textbf{A:}~{\it 7.17},\lb11.104--105

 Walter J.\,H.~ \textbf{A:}~{\it 3.6}

Wamsley, J.\,W.~ 20.36

Wan S.~ {\it 17.102}

 Wang J.~ {\it 11.80}

Wang L.~ \textbf{A:}~{\it 12.38},\, {\it 16.1}

 Wang Yanming~ \textbf{A:}~{\it 15.81}

 Wanless I.\,M.~ \textbf{A:}~17.35

Ward J.~ {\it 14.69}

Ward M.\,B.~ \textbf{A:}~{\it 7.12}

Warfield R.~ 10.52

 Wehrfritz B.\,A.\,F.~ 5.5,\, 8.1--3,\, 11.22,\linebreak
\makebox[3ex][r]{}{\it 18.76},\, \textbf{A:}~{\it 2.46},\, {\it 4.37},\, 5.6,\, {\it
5.19}

 Wei G.\,M.~ 15.99

 Wei Huaquan~ \textbf{A:}~{\it 15.81}

Weidmann R.~ {\it 14.19}

Weigel T.~ {\it 14.53},\, 14.69,\, {\it 19.10},\,  {\it 20.90}

Weiss A.~ \textbf{A:}~12.24

Weiss B.~ 21.86

Weiss R.~ 14.73

Wesolek P.~ 19.69--70,\, 19.107

White S.~ 14.8

Whitehead~ 10.70,\, 17.86,\, \textbf{A:}~17.86

 Wiegold J.~ 4.66,\, 5.52,\,
5.54--56,\, 6.45,\lb11.102,\, 14.95,\, 15.31,\, 15.92,
\, 16.100,\lb17.113--115,\, 21.81,\, \textbf{A:}~ {\it
1.64},\, 4.67--69,\lb5.53,\, 5.56,\, {\it 5.56},\, 6.44,\, 12.45,\, 14.86,\lb16.101,\, 17.116

 Wielandt H.~ 14.43,\, 15.52,\, 17.43,\linebreak
\makebox[3ex][r]{}\textbf{A:}~6.6,\,
6.37,\, 9.64,\, 19.84

Wienke L.~ 20.22,\, 20.85

Wilde T.~ {\it 16.60},\, 19.108,\, 20.115

 Wilkens B.~ \textbf{A:}~{\it 17.46},\, {\it 18.33}

Wilkerson C.\,W.~ 19.92

Wilkes G.~ 20.119

Wille R.\,J.~ 6.29

 Willems W.~ 21.91,\, \textbf{A:}~{\it 3.50}

 Williams G.~ \textbf{A:}~{\it 8.10}

 Williams J.\,S.~ \textbf{A:}~{\it 5.11}

Willis G.\,A.~ 19.71,\, \textbf{A:}~{\it 9.48}

 Wilson J.\,S.~ 7.41,\, {\it 8.16},\ 9.68,\, 11.18,\linebreak
\makebox[3ex][r]{}{\it 12.95},\, 18.104--105,\, {\it 19.22},\, 20.116,\lb  21.136,\, \textbf{A:}~{\it 2.23},\, {\it 4.68},\, 6.46,\, {\it
8.17},\lb 8.76,\, {\it 9.54},\, 11.103,\, 12.90,\, {\it 14.49},\lb{\it 15.14},\, {\it 15.75}

Wilson L.\,E.~ 21.137,\, \textbf{A:}~{\it 16.104}

 Wilson R.\,A.~ \textbf{A:}~{\it
8.35},\, {\it 8.37},\, {\it 11.53},\lb{\it 15.43},\, {\it 15.75},\,
{\it 16.55},\, {\it 18.115}

Wilton H.~ {\it 20.22},\, {\it 20.41},\, 21.138--139

Wintenberger J.-P.~ \textbf{A:}~12.24

Wise D.\,T.~ 20.60,\, 21.124,\, \textbf{A:}~{\it 8.68},\lb{\it 17.108}

Wisliceny J.~ 5.26

Witzel S.~ 21.140

 Woldar A.\,J.~ \textbf{A:}~{\it 16.1}

Wolf J.\,A.~ \textbf{A:}~{\it 7.4}

 Wolf T.\,R.~ 15.2,\, \textbf{A:}~{\it 15.62}

Woods S.\,M.~ {\it 4.55}

Wright C.\,R.\,B.~ \textbf{A:}~{\it 3.8}

Wua J.~ {\it 4.13}

~

\textbf{X}iao W.~ \textbf{A:}~{\it 5.32},\, {\it 6.42},\, {\it 7.44}

 Xu F.~ {\it 18.16}

Xu M.\,C.~ \textbf{A:}~{\it 12.38},\, {\it 12.39}

Xu X.~  {\it 20.87}

Xu Y.~   {\it 20.87}

~

\textbf{Y}abanzhi G.\,G.~ 7.58

Yadav M.\,K.~ 20.126

Yagzhev A.\,V.~ {\it 3.3},\, {\it 8.3},\, \textbf{A:}~{\it 2.76,\, 2.77},\lb{\it 3.7},\, {\it 8.84}

Yakovlev B.\,V.~ 9.36,\, \textbf{A:}~5.69

Yakushevich A.\,V.~ \textbf{A:}~{\it 7.40}

 Yamada H.~ {\it 15.55}

 Yamauchi H.~ {\it 15.55}

Yang J.~ {\it 15.96}

Yang N.~ {\it 17.37},\, {\it 18.38},\, 20.31,\, {\it 20.58},\lb21.111--112

Yang Y.~ {\it 10.34},\, {\it 20.2},\, {\it 20.79},\, {\it 20.80},\lb {\it 20.87},\, 21.108,\, {\it 21.142},\, \textbf{A:}~{\it 15.62}

 Yi X.~ {\it 9.75},\, \textbf{A:}~{\it 14.99},\, {\it 19.85},\, {\it 19.87},\lb{\it 19.88}

Yoshiara S.~ \textbf{A:}~{\it 8.37},\, {\it 12.65}

Yu S.~ \textbf{A:}~{\it 12.78}

~

\textbf{Z}agorin D.\,L.~ \textbf{A:}~{\it 4.27},\, {\it 9.63}

 Za\u{\i}tsev D.I.~ 9.11,\, 9.13--14,
11.30--31,\lb\textbf{A:}~{\it 2.70},\, {\it
4.57,\, 4.58},\, 6.12,\, {\it 6.12},\, 9.12,\linebreak
\makebox[3ex][r]{}10.14

 Za\u{\i}tsev S.\,A.~ \textbf{A:}~18.75

 Zalesski\u{\i}~A.\,E.\,\, 11.32--34,\, {\it 11.32},\lb12.27--29,\, {\it 14.46},\, {\it 16.29},\, {\bf
A:}~3.2,\, 3.14,\lb4.28,\, {\it 4.28},\, 4.29,\,
{\it 4.29},\, 11.33,\, 11.35

 Zalesski\u{\i} P.\,A.~ 16.40--41,\, 20.117--119,\lb21.141--142,\, \textbf{A:}~{\it 8.69},\, {\it 8.70},\, {\it 12.70},\lb{\it 17.19}

Zamyatin A.\,P.~ 7.23,\, \textbf{A:}~{\it 4.63}

Zapata T.~ {\it 20.90}

Zaremsky M.\,C.\,B.~ 21.140,\, 21.143--146

Zassenhaus H.~ \textbf{A:}~17.67

 Zavarnitsine A.\,V.~ {\it 15.53},\, 18.38,\, 20.83,\lb21.54,\,  \textbf{A:}~{\it 13.63},\, {\it 14.60},\, {\it 17.74}

 Zelenyuk E.\,G.~ 13.24--25,\, {\it 13.24,\, 13.25,\lb13.26},\, {\it 15.80},\, 16.44,\, 17.49,\, 17.51,\lb\textbf{A:}~{\it 10.6},\, {\it 12.2},\, 13.26,\, {\it 13.26},\, {\it 13.46},\lb{\it 13.47},\, 16.42--43,\, {\it 16.43},\, 17.50

 Zel'manov E.\,I.~ 12.95,\, 13.21,\, 14.42,\linebreak
\makebox[3ex][r]{}15.41,\, 20.120,\, 21.63,\, \textbf{A:}~{\it 2.63},\, {\it 3.41},\lb{\it 4.74},\, {\it 9.34}

Zeng M.\,R.~ {\it 20.2},\, {\it 21.142}

Zeng Y.~ 21.109

Zenkov A.\,V.~ {\it 16.51},\, 21.147--149

 Zenkov V.\,I.~ {\it 16.29},\, 20.121--123,\, {\it 20.123},\lb21.150,\,  \textbf{A:}~{\it 14.62},\, {\it 5.12},\, {\it 5.32},\, {\it 9.27},\lb{\it 11.64},\, {\it 12.25},\, {\it 12.26},\, {\it 16.13},\, {\it 17.40},\lb{\it 18.31},\, {\it 19.37},\, 19.109

 Zhang J.~ {\it 13.43},\, {\it 16.3},\,
 \textbf{A:}~12.25--26,\linebreak
\makebox[3ex][r]{}{\it 12.78}

Zharov S.\,V.~ \textbf{A:}~{\it 9.19}

Zhitomirski G.~ 19.64

Zhu Y.\,Z.~ {\it 7.31}

 Zhurtov A.\,Kh.~ {\it 10.60},\, 15.98,\, \textbf{A:}~{\it
14.58}

Zieschang H.~ 10.70,\, \textbf{A:}~10.69

Zil'ber B.\,I.~ {\it 18.18},\, \textbf{A:}~{\it 1.26}

Zippin L.~ 6.11

Zisser I.~ 15.3

 Zubkov A.\,N.~ 14.37--42,\, {\it 14.37},\,
15.44,\lb15.93,\, \textbf{A:}~{\it
11.21},\, 15.44

 Zucca P.~ 15.67,\, \textbf{A:}~{\it 14.49}

Zusmanovich P.~ 19.110--111

 Zvezdina M.\,A.~
\textbf{A:}~{\it 16.24},\, {\it 17.36}

Zyabrev I.~ {\it 2.48}, {\it 18.71}

 Zyubin S.\,A.~ \textbf{A:}~{\it 6.38},\, {\it 11.70}

Zyulyarkina N.~ 20.54

\end{document}